\pdfoutput=1
\documentclass[11pt,a4paper]{article}
\usepackage[T1]{fontenc}
\usepackage{lmodern,microtype}
\usepackage[margin=26mm]{geometry}
\usepackage{amsmath,amssymb,amsthm,mathtools}
\usepackage{booktabs,array,needspace,longtable,enumitem,etoolbox}
\usepackage{graphicx}
\graphicspath{{figures/}}
\usepackage[section]{placeins}
\usepackage{pgfplots}
\pgfplotsset{compat=1.18}
\usepgfplotslibrary{groupplots}
\usepackage[colorlinks=true,linkcolor=blue,citecolor=blue,urlcolor=blue]{hyperref}
\expandafter\def\expandafter\UrlBreaks\expandafter{\UrlBreaks\do\-}% allow line breaks after hyphens in file paths
\usepackage[capitalise,nameinlink]{cleveref}
\newcounter{procedure}[section]
\renewcommand{\theprocedure}{\thesection.\arabic{procedure}}
\crefname{procedure}{Algorithm}{Algorithms}
\Crefname{procedure}{Algorithm}{Algorithms}
\newcounter{procedurestep}
\makeatletter
\newenvironment{procedure}[2]{%
  \par\Needspace{12\baselineskip}\small
  \setcounter{procedurestep}{0}%
  \let\c@table\c@procedure

  \setlength{\LTpre}{5pt}\setlength{\LTpost}{8pt}%
  \begin{longtable}{@{}>{\raggedright\arraybackslash}p{.09\linewidth}@{\hspace{.025\linewidth}}p{.885\linewidth}@{}}
  \toprule
  \multicolumn{2}{@{}p{\linewidth}@{}}{\textbf{Algorithm \theprocedure\quad #1}\label[procedure]{#2}}\\*
  \midrule
  \endfirsthead
  \toprule
  \multicolumn{2}{@{}p{\linewidth}@{}}{\textbf{Algorithm \theprocedure\quad #1} (continued)}\\*
  \midrule
  \textbf{Step} & \textbf{Procedure}\\*
  \midrule
  \endhead
  \midrule
  \multicolumn{2}{r@{}}{\footnotesize Continued on next page}\\
  \endfoot
  \bottomrule
  \endlastfoot
}{\end{longtable}\par}
\makeatother
\newcommand{\alginput}[1]{\textbf{Input} & #1\\*[3pt]}
\newcommand{\algoutput}[1]{\textbf{Output} & #1\\*[3pt]\midrule
  \textbf{Step} & \textbf{Procedure}\\*\midrule}
\newcommand{\algstep}[2][0]{\stepcounter{procedurestep}\arabic{procedurestep} &
  {\leftskip=#1em\relax #2\par}\\[3pt]}
\newtheorem{proposition}{Proposition}[section]
\newtheorem{lemma}[proposition]{Lemma}
\theoremstyle{remark}
\newtheorem{remark}[proposition]{Remark}
\newcommand{\R}{\mathbb R}
\newcommand{\T}{\mathbb T}
\newcommand{\one}{\mathbf 1}
\newcommand{\dx}{\Delta x}
\newcommand{\ip}[2]{\langle #1,#2\rangle_M}
\newcommand{\avg}[1]{\langle #1,\one\rangle_M}
\newcommand{\normone}[1]{\lVert #1\rVert_{1,M}}
\newcommand{\A}{\mathcal A}
\newcommand{\E}{\mathcal E}
\newcommand{\C}{\mathcal C}
\newcommand{\K}{\mathcal K}
\newcommand{\prox}{\operatorname{prox}}
\newcommand{\argmin}{\operatorname*{argmin}}
\newcommand{\diag}{\operatorname{diag}}
\newcommand{\supp}{\operatorname{supp}}
\DeclareMathOperator{\KL}{KL}
\definecolor{plotblue}{RGB}{31,91,145}
\definecolor{plotred}{RGB}{174,66,51}
\definecolor{plotgreen}{RGB}{43,120,88}
\definecolor{plotpurple}{RGB}{113,76,145}
\pgfplotsset{each axis/.append style={font=\small,tick label style={font=\footnotesize},
  label style={font=\small},legend style={font=\footnotesize,draw=none},
  axis line style={black!65},grid style={black!12},grid=major}}
\title{Martingale and signed diffusive optimal transport:\\geodesics and gradient flows}
\author{%
  \normalsize Aymeric Baradat\textsuperscript{1}\qquad
  \normalsize Jean-David Benamou\textsuperscript{2}\qquad
  \normalsize Fran\c{c}ois-Xavier Vialard\textsuperscript{3}\\[6pt]
  \small\textsuperscript{1}Univ Lyon, CNRS, Universit\'e Claude Bernard Lyon 1,\\
  \small Institut Camille Jordan, UMR 5208, Villeurbanne, France\\[3pt]
  \small\textsuperscript{2}Inria Paris; CEREMADE, Universit\'e Paris Dauphine--PSL, CNRS, Paris, France\\[3pt]
  \small\textsuperscript{3}LIGM, Universit\'e Gustave Eiffel, Marne-la-Vall\'ee, France}
\date{14 November 2026}
\definecolor{geotime0}{rgb}{0.18621300,0.01880300,0.58722800}
\definecolor{geotime1}{rgb}{0.35635900,0.00379800,0.64781000}
\definecolor{geotime2}{rgb}{0.51793300,0.02156300,0.65410900}
\definecolor{geotime3}{rgb}{0.65558000,0.12972500,0.59231700}
\definecolor{geotime4}{rgb}{0.76809000,0.24481700,0.49846500}
\definecolor{geotime5}{rgb}{0.86292700,0.36511900,0.40351900}
\definecolor{geotime6}{rgb}{0.93563000,0.48771200,0.31595200}
\definecolor{geotime7}{rgb}{0.98304100,0.62413100,0.22793700}
\definecolor{geotime8}{rgb}{0.99250500,0.77796700,0.15285500}
\definecolor{geotimeTransition}{rgb}{0.50645400,0.01633300,0.65620200}
\definecolor{geotimePeakedRefine}{rgb}{0.99385100,0.75930400,0.15909200}
\definecolor{geotimeDiffuseRefine}{rgb}{0.68375800,0.15627800,0.57166000}
\definecolor{flowtime0}{RGB}{68,1,84}
\definecolor{flowtime1}{RGB}{68,7,88}
\definecolor{flowtime2}{RGB}{67,12,93}
\definecolor{flowtime3}{RGB}{67,18,97}
\definecolor{flowtime4}{RGB}{67,24,101}
\definecolor{flowtime5}{RGB}{67,30,106}
\definecolor{flowtime6}{RGB}{66,35,110}
\definecolor{flowtime7}{RGB}{66,41,114}
\definecolor{flowtime8}{RGB}{66,47,119}
\definecolor{flowtime9}{RGB}{66,52,123}
\definecolor{flowtime10}{RGB}{65,58,127}
\definecolor{flowtime11}{RGB}{65,64,132}
\definecolor{flowtime12}{RGB}{65,69,135}
\definecolor{flowtime13}{RGB}{63,74,136}
\definecolor{flowtime14}{RGB}{61,78,136}
\definecolor{flowtime15}{RGB}{59,82,137}
\definecolor{flowtime16}{RGB}{57,87,138}
\definecolor{flowtime17}{RGB}{55,91,138}
\definecolor{flowtime18}{RGB}{53,96,139}
\definecolor{flowtime19}{RGB}{51,100,139}
\definecolor{flowtime20}{RGB}{49,105,140}
\definecolor{flowtime21}{RGB}{47,109,141}
\definecolor{flowtime22}{RGB}{45,113,141}
\definecolor{flowtime23}{RGB}{43,118,142}
\definecolor{flowtime24}{RGB}{42,122,142}
\definecolor{flowtime25}{RGB}{41,126,141}
\definecolor{flowtime26}{RGB}{40,130,140}
\definecolor{flowtime27}{RGB}{40,134,139}
\definecolor{flowtime28}{RGB}{39,138,138}
\definecolor{flowtime29}{RGB}{38,142,137}
\definecolor{flowtime30}{RGB}{38,147,136}
\definecolor{flowtime31}{RGB}{37,151,136}
\definecolor{flowtime32}{RGB}{36,155,135}
\definecolor{flowtime33}{RGB}{36,159,134}
\definecolor{flowtime34}{RGB}{35,163,133}
\definecolor{flowtime35}{RGB}{34,167,132}
\definecolor{flowtime36}{RGB}{40,171,129}
\definecolor{flowtime37}{RGB}{47,174,124}
\definecolor{flowtime38}{RGB}{55,178,120}
\definecolor{flowtime39}{RGB}{62,181,116}
\definecolor{flowtime40}{RGB}{70,185,111}
\definecolor{flowtime41}{RGB}{77,188,107}
\definecolor{flowtime42}{RGB}{85,192,103}
\definecolor{flowtime43}{RGB}{92,195,98}
\definecolor{flowtime44}{RGB}{100,199,94}
\definecolor{flowtime45}{RGB}{107,202,90}
\definecolor{flowtime46}{RGB}{115,206,85}
\definecolor{flowtime47}{RGB}{122,209,81}
\newcommand{\BaselineDualExceedances}{184}
\newcommand{\BaselineDualSteps}{192}
\newcommand{\BaselineDualExcess}{1.4\cdot10^{-12}}
\newcommand{\BaselineInputContrast}{1064}
\newcommand{\BaselineIntervalGap}{2.4\cdot10^{-17}}
\newcommand{\BaselineDensityError}{7.0\cdot10^{-9}}
\newcommand{\DiffuseInitialMoment}{0.0373732}
\newcommand{\DiffuseFinalMoment}{0.0128148}
\newcommand{\BaselineRecordedDipTime}{3.40\cdot10^{-5}}
\newcommand{\FineRecordedMinimum}{0.2213}
\newcommand{\FineRecordedDipTime}{3.60\cdot10^{-5}}
\newcommand{\FineStepGap}{1.8\cdot10^{-15}}
\newcommand{\FineMassError}{4.5\cdot10^{-16}}
\newcommand{\FineIterationMaximum}{982}
\newcommand{\RoundedSymbolError}{1.5\cdot10^{-7}}
\newcommand{\RoundedTrajectoryDifference}{4.1\cdot10^{-8}}
\newcommand{\OrderConvexSignedCoarse}{[6.32534504,6.32534506]\cdot10^{-5}}
\newcommand{\OrderConvexPositiveCoarse}{[7.71774326,7.71774327]\cdot10^{-5}}
\newcommand{\OrderConvexRatioCoarse}{[1.22013,1.22014]}
\newcommand{\OrderConvexSignedFine}{[6.32446394,6.32446396]\cdot10^{-5}}
\newcommand{\OrderConvexPositiveFine}{[7.71724273,7.71724276]\cdot10^{-5}}
\newcommand{\OrderConvexRatioFine}{[1.22022,1.22023]}
\newcommand{\OrderConvexRatioShort}{1.22}
\newcommand{\OrderReversedSignedCoarse}{[6.32534504,6.32534506]\cdot10^{-5}}
\newcommand{\OrderReversedPositiveCoarse}{[1.10307594,1.10307595]\cdot10^{-3}}
\newcommand{\OrderReversedRatioCoarse}{[17.43898,17.43899]}
\newcommand{\OrderReversedSignedFine}{[6.32446394,6.32446396]\cdot10^{-5}}
\newcommand{\OrderReversedPositiveFine}{[1.10297795,1.10297809]\cdot10^{-3}}
\newcommand{\OrderReversedRatioFine}{[17.43986,17.43987]}
\newcommand{\OrderReversedRatioShort}{17.4}
\newcommand{\OrderConvexDRIterations}{18600}
\newcommand{\OrderConvexIPMIterations}{25}
\newcommand{\OrderConvexDRResidual}{7.27\cdot10^{-11}}
\newcommand{\OrderConvexDensityDifference}{1.65\cdot10^{-8}}
\newcommand{\OrderConvexDiffusivityDifference}{6.35\cdot10^{-8}}
\newcommand{\OrderConvexKappaMax}{0.0136}
\newcommand{\OrderConvexSignedKappaMin}{-0.0510}
\newcommand{\OrderConvexMixture}{10^{-10}}
\newcommand{\OrderConvexSignedDRIterations}{12400}
\newcommand{\OrderConvexSignedIPMIterations}{22}
\newcommand{\OrderConvexSignedDRResidual}{9.78\cdot10^{-12}}
\newcommand{\OrderConvexSignedDensityDifference}{3.38\cdot10^{-7}}

\theoremstyle{plain}
\newtheorem{theorem}[proposition]{Theorem}
\newtheorem{corollary}[proposition]{Corollary}
\theoremstyle{definition}
\newtheorem{definition}[proposition]{Definition}
\newcommand{\Ac}{\mathcal A_{\mathrm c}}
\newcommand{\Ec}{\mathcal E_{\mathrm c}}
\newcommand{\Lc}{\mathcal L}
\newcommand{\Pp}{\mathcal P}
\newcommand{\M}{\mathcal M}
\newcommand{\Pcal}{\mathcal P}
\newcommand{\Mcal}{\mathcal M}
\newcommand{\Torus}{\mathbb T^d}
\newcommand{\dd}{\mathrm d}
\newcommand{\fr}{\mathrm{fr}}
\newcommand{\weakto}{\rightharpoonup}
\newcommand{\DivDiv}{\operatorname{div}\operatorname{Div}}

\renewcommand{\thesection}{\Alph{part}\arabic{section}}
\renewcommand{\theHsection}{\Alph{part}.\arabic{section}}
\newcommand{\mainpart}[1]{\clearpage\setcounter{section}{0}\part{#1}}
\makeatletter
\patchcmd{\l@section}{1.5em}{2.5em}{}{\errmessage{TOC number-width patch failed}}
\renewcommand*\l@subsection{\@dottedtocline{2}{1.5em}{3.1em}}
\makeatother
\hypersetup{pdftitle={Martingale and signed diffusive optimal transport: geodesics and gradient flows},pdfauthor={Aymeric Baradat, Jean-David Benamou, Francois-Xavier Vialard}}

\newcommand{\tfin}{t_{\mathrm f}}
\DeclareMathOperator{\spr}{spr}
\newcommand{\gaprel}{\operatorname{gap}_{\mathrm{rel}}}
\newcommand{\ind}{\mathbb I}

\newcommand{\Am}{{\mathcal A^{\mathrm M}}}
\newcommand{\Em}{\mathcal E_{\mathrm c,M}}
\newcommand{\Vm}{V_{\rm ext}}
\newcommand{\Ptwo}{\mathcal P_2}
\newcommand{\Sd}{\mathbb S^d}
\newcommand{\Sdp}{\mathbb S^d_+}
\DeclareMathOperator{\tr}{tr}
\DeclareMathOperator{\Law}{Law}
\newcommand{\MB}{\mathsf{MB}}
\newcommand{\fA}{\mathfrak A}
\newcommand{\fB}{\mathfrak B}

\newcommand{\cG}{\mathcal G}
\newcommand{\cH}{\mathcal H}
\theoremstyle{definition}
\newtheorem{example}[proposition]{Example}
\theoremstyle{remark}
\newtheorem{openquestion}[proposition]{Open question}

\begin{document}
\maketitle
\section*{Extended abstract}\addcontentsline{toc}{section}{Extended abstract}

\subsection*{Diffusivity of either sign}

Diffusion models usually prescribe a nonnegative diffusivity.
With a negative diffusivity, the equation can be ill-posed: for a constant $\kappa<0$, $\partial_t\rho=\partial_{xx}(\rho\kappa)$ is the backward heat equation, which is ill-posed as an initial-value problem.
Negative diffusivities nevertheless occur in nonlinear equations.
The Derrida--Lebowitz--Speer--Spohn (DLSS) equation \cite{DLSS91}
\[
 \partial_t\rho=-\partial_{xx}\bigl(\rho\,\partial_{xx}\log\rho\bigr)
\]
has the form $\partial_t\rho=\partial_{xx}(\rho\kappa)$ with $\kappa=-\partial_{xx}\log\rho$, and this diffusivity is negative near every positive nondegenerate minimum of $\rho$.

In signed diffusive transport \cite{MRSS,MRSSdlss}, the diffusivity is an unknown of a variational problem, and its sign is not constrained.
A path of probability densities $\rho_t$, $0\le t\le1$, on the circle $\T=\R/\mathbb Z$ and a momentum $m_t$ satisfy the conservation law $\partial_t\rho=\partial_{xx}m$, and the cost of the path is the action
\[
 \int_0^1\!\!\int_{\T}\frac{m^2}{2\rho}\,\dd x\,\dd t .
\]
The minimum of the action over the paths between two probability measures $\mu^0$ and $\mu^1$ on $\T$ is denoted $\mathsf G(\mu^0,\mu^1)$; a minimizing path is a geodesic, and $\mathsf D=(2\mathsf G)^{1/2}$ is the signed diffusive transport distance.
This is the dynamic formulation of optimal transport of Benamou and Brenier \cite{BB2000}, with the continuity equation replaced by a second-order conservation law.
Writing $m=\rho\kappa$ gives $\partial_t\rho=\partial_{xx}(\rho\kappa)$: the minimizer determines the diffusivity, with either sign, and no backward heat equation is solved as an initial-value problem.
On $\T$, the mass is the only conserved quantity of the form $\int_{\T}\varphi\rho_t\,\dd x$, because its time derivative $\int_{\T}\partial_{xx}\varphi\,m_t\,\dd x$ vanishes for every momentum only if $\varphi$ is constant. The minimum action $\mathsf G(\mu^0,\mu^1)$ is finite for all probability measures $\mu^0,\mu^1$, singular ones included (\cref{geo:prop:heat-semigroup-hypotheses}).

On the real line, a nonnegative diffusivity connects the model with probability theory.
With $m_t=\frac12a_t\rho_t$ and $a_t\ge0$, the conservation law is the Fokker--Planck equation of a continuous martingale with variance rate $a_t$, and the action equals $\frac18\int_0^1\!\int_{\R}a_t^2\rho_t\,\dd x\,\dd t$ (\cref{mart:prop:positive-eulerian}).
This is the martingale Benamou--Brenier problem of Huesmann and Trevisan \cite{HT19}; related convex formulations with nonnegative momentum are used to calibrate local volatility models from prescribed marginals \cite{GLW19}.
Its cost is finite only if the endpoints $\mu^0,\mu^1$ on $\R$ are in convex order, that is, if every convex function has a mean under $\mu^0$ at most equal to its mean under $\mu^1$.

The manuscript treats both models and two variational problems: geodesics, with two prescribed endpoints, and minimizing movements, which discretize gradient flows of an energy in time.
It studies three limits: the time step tends to zero, the spatial mesh size tends to zero, and the iterates of an optimization method converge on a fixed grid.
The main contributions are the time limit and the joint space--time limit of signed geodesics, the time limit of a modified Jordan--Kinderlehrer--Otto (JKO) scheme under a positivity assumption, a common limit of two discretizations of martingale transport, and computations of geodesics and modified JKO trajectories, with a certified optimization error for selected geodesics and for each step of selected trajectories, given the computed input density of the step.

\subsection*{Discretization}

A geodesic is approximated on $N$ time cells.
In each cell, the conservation law holds exactly between the two endpoint densities of the cell, and the action of the momentum is weighted by their average, the midpoint weight.
For the entropic evolution, a step of length $\tau$ of the JKO scheme \cite{JKO98} minimizes, over the new density $\rho$, the entropy of $\rho$ plus $\mathsf D(\rho^{\rm old},\rho)^2/(2\tau)$, where $\rho^{\rm old}$ is the density at the beginning of the step.
Evaluating $\mathsf D^2$ requires a minimization over all paths from $\rho^{\rm old}$ to $\rho$.
The modified JKO scheme replaces these paths by a single time cell: its step minimizes the entropy of $\rho$ plus $\tau^{-1}$ times the action of one cell, weighted by the density $\rho$ at the end of the step.
Each modified step is therefore a single convex problem, whose unknowns are the new density and one momentum.
In space, we use $M$ nodes on a periodic grid and the Fourier second derivative.
Both problems are then finite-dimensional convex programs, and each entropic step has a unique, strictly positive minimizer (\cref{prop:positive}).
We solve them by Douglas--Rachford (DR) splitting and by the interior-point method (IPM) Clarabel \cite{Clarabel}.
For martingale transport, we compare two time discretizations.
The first is the discretization with midpoint weights described above, on $\R$ and with a nonnegative momentum.
The second uses a martingale coupling between each pair of successive marginals; the cost of a coupling is proportional to the fourth moment of its increment (\cref{mart:sec:setting}).

\subsection*{Main results}

\begin{itemize}[leftmargin=*]
\item \textbf{Time limit of signed geodesics.}
For all probability endpoints on $\T$, the time-discrete problems $\Gamma$-converge to the continuous problem as $N\to\infty$, for the uniform-in-time weak convergence of measures.
The minimum actions converge to $\mathsf G(\mu^0,\mu^1)$, and the limits of convergent subsequences of discrete minimizers are geodesics (\cref{geo:thm:time-discrete-convergence,geo:prop:heat-semigroup-hypotheses}).
No positive lower bound on the densities is needed: the proof averages a continuous path in time, this averaging commutes with the conservation law, and the convexity of the action bounds the cost of the averaged path.
\item \textbf{Joint space--time limit of signed geodesics.}
The fully discrete problem solved in the computations $\Gamma$-converges as $M$ and $N$ tend to infinity, without a relation between the two meshes: its minimum actions converge to $\mathsf G(\mu^0,\mu^1)$, and subsequences of its minimizers converge to geodesics (\cref{geo:thm:joint-grid}).
The endpoints can be arbitrary probability measures, averaged on the grid by a positive Fej\'er-type kernel, or prescribed grid densities with a common positive lower bound that converge weakly to $\mu^0$ and $\mu^1$, such as the sampled profiles of the experiments (\cref{geo:cor:grid-endpoints}); the intermediate densities may vanish.
The averaging kernel commutes with the second derivative and does not increase the action, so it turns continuous paths into admissible grid paths with no larger action, without a positive lower bound on the densities.
\item \textbf{Time limit of the modified JKO scheme.}
On the torus $\Torus=(\R/\mathbb Z)^d$, with a matrix-valued momentum and an energy that sums a potential energy, an interaction energy and an internal energy, the scheme converges along subsequences, strongly in $L^1$, to a weak solution of a fourth-order equation; for the entropy in dimension one, this equation is the DLSS equation (\cref{jko:thm:main,jko:cor:dlss}).
The result is conditional: it assumes that the densities of the scheme stay above a positive constant on every finite time interval, for all small time steps (assumption (P) in \cref{sec:merged-jko}), and this assumption is not proved for any nonconstant initial density.
\item \textbf{Martingale transport.}
In dimension one, the two time discretizations of martingale transport have different variables and costs, but their normalized minimum values converge to the same martingale Benamou--Brenier cost (\cref{mart:cor:common-limit}); for the discretization by couplings, this convergence holds in every dimension (\cref{mart:c:cor:bb-convergence}).
For the potential energy $\mu\mapsto\int\Vm\,\dd\mu$ of probability measures on $\R$, with bounded second and third derivatives of $\Vm$, the minimizing movement with nonnegative momentum and the action weighted by $\tau^{-1}$ converges, without a density lower bound, to the solution $\mu_t$ of $\partial_t\mu=\partial_{xx}\bigl((-\Vm'')_+\mu\bigr)$ (\cref{mart:thm:positive-potential}).
For a uniformly concave potential, the minimizing movement with the fourth-moment coupling cost has the same limit after normalization (\cref{mart:cor:positive-potential-normalization}).
\end{itemize}

\subsection*{Computation and certification}

Neither DR nor IPM is faster in all tests (\cref{tab:geotarget}, \cref{sec:flowcost}).
We study the local convergence of DR through its linearization at a minimizer.
For the geodesic with uniform endpoints, the convergence factor of this linearization has an exact formula at one value of the proximal parameter of DR and a lower bound at every value; at fixed DR parameters, the lower bound, and therefore the factor, tends to one under grid refinement (\cref{app:spectrum}).
For the entropic step, a model of the linearization at the minimizer of the step gives a rule for the DR parameters (\cref{sec:drflow}).
The optimization error of a computed solution is bounded by dual lower bounds evaluated in interval arithmetic at this solution.
The dual objective reported by the solver is not used: it exceeds a certified upper bound on the optimal value in $\BaselineDualExceedances$ of the $\BaselineDualSteps$ steps of the baseline entropic trajectories (\cref{tab:reporteddual}).
For an entropic step with computed density $\nu$, normalized to unit mass, and exact minimizer $\nu^*$ of the certified problem, a certified objective gap $\eta$ also gives $\frac1M\sum_i|\nu_i-\nu_i^*|\le\sqrt{2\eta}$ (\cref{prop:gap}).
These certificates are evaluated for the IPM geodesics, for the entropic steps of \cref{tab:baseline-validation,tab:common-entropic,tab:flowcontrast} and for the computations of \cref{part:time-study}; other reported gaps are floating-point estimates. The certificates bound the optimization error of each discrete problem, not the discretization error.

\subsection*{Numerical observations}

\begin{itemize}[leftmargin=*]
\item \textbf{Near vacuum along geodesics.}
Between positive endpoints, the computed geodesics approach zero density inside the time interval.
For three pairs of endpoints at $(M,N)=(63,64)$, the smallest interior densities of the IPM solutions lie between $1.6\cdot10^{-11}$ and $5.2\cdot10^{-10}$ (\cref{tab:ipmvalidation}); IPM keeps its densities positive, but the proximal outputs of DR contain exact zeros (\cref{tab:zeros}).
For the peaked pair, these minima persist under joint refinement of space and time up to $(M,N)=(511,256)$ (\cref{fig:geofine}).
The diffusivity changes sign next to these minima: on the grids $M=63,127,255$ with $N=64$, each of the 104 sign changes per grid at times $0.1\le t\le0.9$ lies within $0.52/M$ of a local minimum of the midpoint weight in the same time cell (\cref{fig:geocolormaps}).
At the transition examined in \cref{sec:vacuumobservation}, near $x=0.55$ at $t=31/128$, the diffusivity goes from values above $0.045$ to values below $-0.045$ between two adjacent nodes, and this transition narrows under spatial refinement.
At two nodes where the diffusivity is at least $c>0$ and at most $-c$, the sum of the two midpoint weights is at most $1/c$ times the absolute difference of the two momenta divided by the time step.
The computations are therefore consistent with a limit density that vanishes where the diffusivity changes sign; they do not prove it (\cref{lem:vacuumcriterion}).
\item \textbf{Decrease of the density minimum along the modified JKO scheme.}
The entropy does not increase along the scheme, but the density minimum falls below its initial value for the four baseline initial densities (\cref{fig:minimum}); we call the time at which it is smallest the dip time.
At $M=128$, for the initial densities proportional to $\delta+\sin^{2r}(\pi(x-1/2))$, whose minimum has order $2r$ and offset $\delta>0$, quadratic minima ($2r=2$) do not decrease, minima of order $4$, $8$ and $16$ do, and for the orders $8$ and $16$ a smaller offset gives an earlier and deeper decrease (\cref{tab:isolatedsurvey}).
Rescaling space by $\delta^{1/(2r)}$, time by $\delta^{2/r}$ and the density by its initial minimum leaves the DLSS equation unchanged near the minimum; if the rescaled dip time stays bounded above and away from zero, the dip time therefore lies between two positive multiples of $\delta^{2/r}$, and the measured exponents differ from $2/r$ by at most $0.03$ (\cref{sec:inner,tab:dipscaling}).
For a peaked initial density at $M=128$, the density minima of the scheme with $256$ steps and of the spatially discretized DLSS equation, solved without time discretization, both fall to $0.21$ times the initial minimum (\cref{tab:merged-minima}).
The minimum remains positive, so this decrease does not contradict assumption (P); it shows, however, that the scheme on these grids violates a discrete minimum principle. A minimum principle for the scheme in continuous space would give (P) with the initial minimum as lower bound; these grid computations do not decide whether it holds.
\item \textbf{Cost of a nonnegative diffusivity.}
The constraint $\kappa\ge0$ keeps the geodesic problem convex.
For three pairs of endpoints with Gaussian profiles, it multiplies the minimum action by $\OrderConvexRatioShort$ when the endpoints are in convex order and by more than $17$ in the two other cases (\cref{tab:positive-cost}).
The signed geodesic of the pair in convex order still has negative diffusivities: for a positive density path, the momentum of least action has diffusivities of zero sum in every time cell (\cref{sec:nonnegative-diffusivity}).
\item \textbf{Time refinement on fixed spatial grids.}
For geodesics at $M=127$, the certified upper bounds on the minimum actions on $N$ and $2N$ time cells differ by amounts that decrease at observed orders between $1.82$ and $2.42$ in the time step (\cref{tab:geo-time-convergence}).
For the modified JKO scheme at $M=128$, the differences at the final time from the solution of the spatially discretized DLSS equation decrease at observed orders between $0.92$ and $1.01$ (\cref{tab:merged-convergence}).
\end{itemize}

\subsection*{Open questions}

The joint space--time limit of the modified JKO scheme is open: it requires estimates and a positive density lower bound that do not depend on the two meshes, and the convergence of the discrete nonlinear term to $\rho\,\partial_{xx}\log\rho$ (see the Conclusions).
It is also open whether the limit density vanishes where the diffusivity changes sign, and whether the diffusivity has distinct one-sided limits there.
For nonlinear energies, the limit of the minimizing movement with nonnegative momentum is not identified (\cref{jko:prop:positive-compactness}).
In dimension $d\ge2$, the convergence of the momentum formulation of martingale transport and its numerical methods are future work (\cref{mart:rem:higher-bb,mart:rem:computational-scope}), and so is its comparison with the coupling scheme (\cref{mart:open:tensor-bb}).
\Needspace{9\baselineskip}
\paragraph{Declaration of AI assistance.}
OpenAI Codex and Anthropic Claude assisted with drafting and editing this manuscript and with developing and checking mathematical arguments.
They also assisted with writing and revising numerical code, preparing computational checks, and producing tables and figures from saved numerical data.
AI assistance does not constitute independent verification of the mathematical arguments or numerical results.
The authors retain responsibility for the final manuscript, including its accuracy, originality, and conclusions.

\noindent\textbf{Keywords:} signed diffusive transport; geodesics; $\Gamma$-convergence; entropic evolution; Douglas--Rachford splitting; primal--dual validation.\\
\textbf{Mathematics Subject Classification (2020):} 49J45; 65K10; 65M12; 90C25.

\clearpage
\section*{Introduction}\addcontentsline{toc}{section}{Introduction}

We study variational models that describe changes in a probability
distribution through a second-order conservation equation. For a smooth
positive density $\rho_t$ on a one-dimensional domain $\Omega$, the
momentum $m_t$ and its kinetic action take the form
\[
 \partial_t\rho_t=\partial_{xx}m_t,\qquad
 \frac12\int_0^1\!\!\int_\Omega\frac{m_t^2}{\rho_t}\,\dd x\,\dd t.
\]
The equation specifies how momentum changes the distribution, while
the action penalizes its size relative to the density.
In signed diffusive transport on the torus $\T=\R/\mathbb Z$, the
momentum may have either sign \cite{MRSS,MRSSdlss}. This model extends
the action-based approach of Benamou and Brenier \cite{BB2000}.
The measure formulation below also allows singular marginals and
zero densities.

Nonnegative momentum on the real line gives a second model.
Writing $m_t=a_t\rho_t/2$, with $a_t\ge0$, turns the conservation
equation into the Fokker--Planck equation of a continuous martingale
with variance rate $a_t$. The momentum action is then one eighth of
the integral of $a_t^2$ against the time marginals. This identifies
the positive Eulerian model with the martingale Benamou--Brenier
formulation of Huesmann and Trevisan \cite{HT19}. Its admissible
endpoint marginals must satisfy the convex-order constraint.
Related convex Eulerian formulations are used for local-volatility
calibration: Guo, Loeper, and Wang \cite{GLW19} recover a local
volatility function from prescribed marginal distributions by
optimizing a diffusion with nonnegative momentum.

The martingale problem can also be discretized through successive
couplings. Each transition $(X,Y)$ satisfies
$\mathbb E[Y\mid X]=X$, and its cost is based on the fourth moment of
$Y-X$. The coupling variables retain information about transitions,
while the Eulerian variables describe marginals and momentum.
Their discrete costs differ, but their normalized geodesic values
converge to the same continuous cost in one dimension.

We distinguish two variational problems. The geodesic problem
prescribes both endpoint marginals and minimizes the action of a
connecting path. For gradient flows, a prescribed energy drives
successive choices of the next marginal through a minimizing-movement
scheme \cite{JKO98,AGS}. Our Eulerian steps use a single-cell action
weighted by the measure at the end of the step; we call them modified JKO steps.
A JKO step for the transport distance would instead use the minimum
action over all connecting paths inside each step.

The analysis and computations involve three limits: the limit as the
time step tends to zero, the limit as the spatial mesh size $\Delta x$
tends to zero, and convergence of the optimization iterates on a fixed
space--time grid.
Our goal is to explore efficient numerical solvers for these
variational problems. This first requires establishing that the
discrete problems are well defined and proving convergence towards
the corresponding continuous formulations under the stated assumptions.

Part~\ref{part:convergence} studies convergence in time and, for signed geodesics, joint refinement of space and time.
For signed Eulerian geodesics, \cref{geo:thm:time-discrete-convergence} proves, under hypotheses \ref{geo:hyp:real}--\ref{geo:hyp:order} on the spatial operator and for endpoints joined by a path of finite action, $\Gamma$-convergence of the time-discrete action, convergence of the minimum actions and convergence of subsequences of minimizing paths; for $\partial_{xx}$, all probability endpoints satisfy this condition (\cref{geo:prop:heat-semigroup-hypotheses}).

For the Fourier-grid geodesic problem, \cref{geo:thm:joint-grid} proves
joint $\Gamma$-convergence without a relation between the spatial and
temporal meshes. The proof constructs admissible grid paths that
approximate continuous paths and their actions. Positive Fourier
averaging of both marginals and momenta makes this construction
possible even for singular endpoint distributions.
\Cref{geo:cor:grid-endpoints} also covers endpoint densities specified
directly by grid values, such as the sampled positive profiles used in
the experiments. We assume a common positive lower bound at the two
endpoints and convergence of the distributions represented by these
grid values to the prescribed endpoint distributions. The intermediate
densities may vanish.

For the signed modified JKO scheme, \cref{jko:thm:main} proves subsequential convergence to a fourth-order equation under the positivity floor \eqref{jko:eq:floor}, a positive density lower bound uniform in the time step.

The remaining sections of Part~\ref{part:convergence} start from the continuous martingale
problem and distinguish two approaches to time discretization
(\cref{mart:sec:setting}): a momentum formulation imposing nonnegative
diffusion through a local equation, and a coupling formulation using
$M_4$ as a short-time static proxy. Section~\ref{mart:sec:limits} proves geodesic convergence of the coupling scheme in every dimension
and of the midpoint scheme with nonnegative momentum in dimension one. The latter result is a
corollary of the same time-averaging argument as signed geodesic
convergence, with convex-order tightness on $\R$; the fourth-moment
proof also gives continuous-path compactness for the coupled processes.
In one dimension, the momentum formulation avoids optimization over
successive couplings; the computations of Part~\ref{part:numerical}
with nonnegative momentum concern the periodic problem, which is not a
martingale problem (\cref{rem:martingale}).
Section~\ref{mart:g:sec:m4-theorem} then treats minimizing movements.
It identifies the constrained momentum flow for one-dimensional potentials
with bounded second and third derivatives, without a density floor,
and proves that it agrees with the fourth-moment flow in the latter's
uniformly concave setting after normalization. The corresponding
proofs are in Appendix~\ref{mart:app:proofs}.

For nonlinear energies satisfying the assumptions (E1)--(E2) of
\cref{sec:merged-jko}, \cref{jko:prop:positive-compactness} proves, on
the torus $\Torus=(\R/\mathbb Z)^d$ with nonnegative tensor momentum,
existence at each positive-momentum step and subsequential weak
convergence of the marginal curves and the space--time momentum
measures; no result for nonlinear energies is proved on $\R$. We have not identified the
limiting momentum as a function of the density and the first variation
of the energy, so the limiting gradient-flow equation remains unknown.
Future work includes extending the one-dimensional positive-momentum
geodesic and potential-flow results and the momentum implementation
to higher dimensions, identifying higher-dimensional coupling flows,
and comparing the tensor and coupling models.

Part~\ref{part:numerical} treats numerical optimization on periodic grids.
We solve the signed geodesic and modified JKO problems by Douglas--Rachford (DR) splitting (\cref{sec:dr}) and, after reformulation as conic programs, by the interior-point solver Clarabel \cite{Clarabel} (\cref{sec:ipm}).
Neither method is uniformly preferable: IPM is faster in the tested
entropic-step comparisons, while DR is competitive for peaked geodesics
(\cref{sec:geocompare,sec:flowcost}).
A linearization of DR at a minimizer gives, on the uniform geodesic, the exact convergence factor of the linearized iteration at one value of the DR proximal parameter and a lower bound on this factor at the other values (\cref{prop:constant}). A model of this linearization gives a rule for the DR parameters in the modified JKO scheme, and dual lower bounds evaluated in interval arithmetic bound the optimization error (\cref{sec:certificate}).
Interestingly, in the computed geodesics, the DR proximal outputs have exact zero density entries (\cref{sec:proximal-zeros}), and the diffusivity changes sign near density minima (\cref{sec:vacuumobservation}); along the computed trajectories of the modified JKO scheme, the density minimum decreases below its initial value (\crefrange{sec:minimum}{sec:inner}).

Part~\ref{part:time-study} measures the effect of time refinement on prescribed spatial grids, for geodesics in \cref{sec:geo-time-study} and for the modified JKO scheme in \cref{sec:merged-refinement}.

\begin{table}[htbp]
\centering\small
\begin{tabular}{@{}p{.23\linewidth}p{.69\linewidth}@{}}
\toprule
Notation & Meaning throughout the paper\\
\midrule
$M$, $\Delta x=1/M$ & Spatial node count and mesh size in the numerical problems.\\
$N$, $h=1/N$ & Number and length of geodesic time cells on $[0,1]$.\\
$\tau$, $\tfin$ & Minimizing-movement time step and final time. The coupling geodesic estimates also use $\tau=h=1/N$.\\
$L$, $\Lc$ & Periodic Fourier second-derivative operator on $\R^M$ and continuous spatial operator; $\Lc=\partial_{xx}$ for signed diffusive transport.\\
$\A$, $\Ac$ & Grid and measure actions, both with the normalization $1/2$.\\
$\E$, $\Ec$ & Grid entropy and continuous energy of \eqref{jko:eq:energy}, with potential $V_{\rm ext}$, interaction kernel $W$ and internal-energy density $f$; $\Ec$ is the entropy when $V_{\rm ext}=W=0$ and $f(r)=r\log r-r+1$.\\
$\nu$; $\mu,\sigma$; $\rho$ & Grid density vector; probability measures; density with respect to Lebesgue measure.\\
\bottomrule
\end{tabular}
\caption{Notation linking the numerical and continuous formulations. Local algorithm parameters are defined in Part~\ref{part:numerical}.}
\label{tab:merged-notation}
\end{table}

\clearpage
\setcounter{tocdepth}{2}
{\small\tableofcontents}
\clearpage
\mainpart{Time and space convergence}\label{part:convergence}
In this part, we study the time discretization of geodesic and
gradient-flow variational problems in two settings: the signed
diffusive metric on the torus and the martingale Benamou--Brenier
formulation on the real line introduced by Huesmann and Trevisan
\cite{HT19}. We first study time limits with space continuous and also
prove joint space--time convergence for the signed Fourier-grid
geodesic problem.

In the martingale setting, we consider two discretization methods
working at different levels. The first uses the Eulerian representation
and constrains the momentum to be nonnegative in one dimension.
The second works with martingale couplings and minimizes the fourth
moment of the increments.
Table~\ref{tab:time-convergence-results} summarizes our contributions.

\begin{table}[htbp]
\centering\small
\begin{tabular}{@{}>{\raggedright\arraybackslash}p{.17\linewidth}>{\raggedright\arraybackslash}p{.25\linewidth}>{\raggedright\arraybackslash}p{.49\linewidth}@{}}
\toprule
Problem & Variables and domain & Result\\
\midrule
Geodesic problem & Signed momentum on $\T$ & Time and joint Fourier-grid $\Gamma$-convergence; minimizing paths (\cref{geo:thm:time-discrete-convergence,geo:thm:joint-grid}).\\
& Nonnegative momentum on $\R$ & Midpoint action: normalized value convergence to the martingale cost (\cref{mart:cor:positive-midpoint}).\\
& Martingale couplings on $\R^d$ & Fourth-moment action: continuous-path compactness and convergence of values (\cref{mart:c:thm:compactness,mart:c:cor:bb-convergence}).\\
\midrule
Minimizing movements & Signed tensor momentum on $\Torus$ & Nonlinear energy: identified limit under the density floor \eqref{jko:eq:floor} (\cref{jko:thm:main}).\\
& Nonnegative tensor momentum on $\Torus$ & Nonlinear energy: existence and weak compactness (\cref{jko:prop:positive-compactness}).\\
& Nonnegative momentum on $\R$ & Potential energy: identified limit for bounded $V_{\rm ext}''$ and $V_{\rm ext}'''$ (\cref{mart:thm:positive-potential}).\\
& Martingale couplings on $\R$ & Potential energy: identified $M_4$ limit under uniform concavity (\cref{mart:thm:oned}).\\
\bottomrule
\end{tabular}
\caption{Convergence results in Part~\ref{part:convergence}. The assumptions, normalizations and convergence topologies are those of the cited statements. Space remains continuous except in the joint Fourier-grid geodesic result.}
\label{tab:time-convergence-results}
\end{table}

\section{Convergence of signed Eulerian geodesics}\label{sec:merged-geodesics}
We first discretize the signed Eulerian geodesic problem in time, with space kept continuous.
We study the limit $N\to\infty$ of the minimum value and of the minimizing paths of the midpoint action for the geodesic problem.
The endpoint measures may be singular and the paths may have zero densities, so the argument does not assume positive densities.
The theorem below identifies the limiting minimum action and the cluster points of minimizing paths. Section~\ref{geo:sec:joint-grid} then adds spatial Fourier discretization and proves the joint limit.

\subsection{Continuous action and midpoint discretization}
Let \(\T=\R/\mathbb Z\), and let \(\Pp(\T)\) denote the probability
measures on \(\T\), and \(\M(\T)\) the finite signed measures, with total variation denoted by $|m|$. For a finite positive measure \(\mu\) and a finite
signed measure \(m\), set
\begin{equation}\label{geo:eq:action-density}
  \Ac(\mu,m)
  =
  \begin{cases}
    \displaystyle
    \frac12\int_\T
      \left|\frac{\dd m}{\dd\mu}\right|^2\,\dd\mu,
      &m\ll\mu,\\[1ex]
    +\infty,&\text{otherwise}.
  \end{cases}
\end{equation}
We call \(\Ac\) the quadratic perspective action. The subscript
\(\mathrm c\) distinguishes this measure action from the grid action
\(\A\) introduced in Part~\ref{part:numerical}.
Let \(\Lc\) be a real Fourier multiplier with symbol \(\mathfrak p(k)\), normalized by
\(\mathfrak p(0)=0\).
Let \(\Lc^*\) denote the adjoint multiplier, with symbol \(\mathfrak p(-k)\), so that \(\int_\T\varphi\,\dd(\Lc m)=\int_\T\Lc^*\varphi\,\dd m\) for \(\varphi\in C^\infty(\T)\) and \(m\in\M(\T)\).
\begin{definition}[Admissible pair]\label{geo:def:admissible-pair}
An admissible pair is a weakly continuous curve
\(\mu\in C_{\mathrm w}([0,1];\Pp(\T))\) and a weakly measurable family
\(m_t\in\M(\T)\), defined for almost every \(t\), such that
\[
  \int_0^1 |m_t|(\T)\,\dd t<\infty,
  \qquad
  \partial_t\mu_t=\Lc m_t
\]
in distributions.  Equivalently, for every
\(\varphi\in C_c^\infty((0,1)\times\T)\),
\begin{equation}\label{geo:eq:continuity-weak}
  \int_0^1\!\!\int_\T \partial_t\varphi\,\dd\mu_t\,\dd t
  +
  \int_0^1\!\!\int_\T \Lc^*\varphi\,\dd m_t\,\dd t=0.
\end{equation}
\end{definition}
The associated space--time measures are
\begin{equation}\label{geo:eq:space-time-measures}
  R=\dd t\otimes\mu_t,
  \qquad
  \mathsf W(B)=\int_0^1m_t(B_t)\,\dd t .
\end{equation}
Here $B_t=\{x:(t,x)\in B\}$. Thus \(R\) is a positive measure whose time marginal is Lebesgue measure,
whereas \(\mathsf W\) is a signed measure and has no prescribed time marginal.
We use the same notation for the quadratic perspective action on
space--time measures: for a positive measure \(R\) and a signed
measure \(\mathsf W\) on \([0,1]\times\T\), set
\begin{equation}\label{geo:eq:space-time-action}
  \Ac(R,\mathsf W)
  =
  \begin{cases}
    \displaystyle
    \frac12\int
      \left|\frac{\dd \mathsf W}{\dd R}\right|^2\,\dd R,
      &\mathsf W\ll R,\\[1ex]
    +\infty,&\text{otherwise}.
  \end{cases}
\end{equation}
For the measures associated with an admissible pair, this is exactly
the time integral of the spatial action (\cref{geo:lem:space-time-action}):
\[
 \Ac(\dd t\otimes\mu_t,\dd t\otimes m_t)
 =\int_0^1\Ac(\mu_t,m_t)\,\dd t.
\]
Thus the two expressions describe the same action, either as an
integral in time or directly on space--time measures.

For endpoints \(\mu^0,\mu^1\in\Pp(\T)\), let
\begin{equation}\label{geo:eq:continuum-value}
  \mathsf G_{\Lc}(\mu^0,\mu^1)
  =
  \inf\left\{
    \int_0^1\Ac(\mu_t,m_t)\,\dd t:
    (\mu,m)\text{ admissible},\
    \mu_0=\mu^0,\ \mu_1=\mu^1
  \right\}.
\end{equation}
For $\Lc=\partial_{xx}$, set
\begin{equation}\label{geo:eq:distance}
 \mathsf D(\mu^0,\mu^1)^2=2\mathsf G_{\Lc}(\mu^0,\mu^1).
\end{equation}
The normalization gives the action density $m^2/(2\rho)$ for densities. The grid action in Part~\ref{part:numerical} uses the same convention.

Let \(N\ge1\), \(h=N^{-1}\), and \(t_j=jh\).  A discrete path consists of
\(\mu_h^0,\ldots,\mu_h^N\in\Pp(\T)\) and
\(m^0,\ldots,m^{N-1}\in\M(\T)\), with
\[
  \Lc m^j=\mu_h^{j+1}-\mu_h^j.
\]
Writing \(\bar\mu_h^j=(\mu_h^j+\mu_h^{j+1})/2\), define
\begin{equation}\label{geo:eq:time-discrete-value}
  \mathsf G_{\Lc,N}(\mu^0,\mu^1)
  =
  \inf
  \left\{
    \sum_{j=0}^{N-1}\frac1h\Ac(\bar\mu_h^j,m^j):
    \mu_h^0=\mu^0,\ \mu_h^N=\mu^1,\
    \Lc m^j=\mu_h^{j+1}-\mu_h^j
  \right\}.
\end{equation}
We call each time interval \([t_j,t_{j+1}]\), of length \(h\), a
\emph{time cell}; its midpoint is \(t_j+h/2\).
Here \(m^j\) is the momentum integrated over this cell, whereas \(m_t\) in \eqref{geo:eq:continuum-value} is a momentum per unit time.  If
the nodes are interpolated affinely and the momentum is set equal to
\(m^j/h\) on \([t_j,t_{j+1}]\), the constraint holds exactly.  For this pair, the continuum action of the \(j\)-th cell is \(h^{-2}\int_{t_j}^{t_{j+1}}\Ac(\mu_t,m^j)\,\dd t\).
The \(j\)-th term of \eqref{geo:eq:time-discrete-value} is
\(h^{-1}\Ac(\bar\mu_h^j,m^j)\): it evaluates the density at the
midpoint of the cell. Convexity of the perspective action gives
\[
 \frac1h\Ac(\bar\mu_h^j,m^j)
 \le\frac1{h^2}\int_{t_j}^{t_{j+1}}\Ac(\mu_t,m^j)\,\dd t.
\]
The discrete problem therefore both restricts the paths to affine
interpolations and replaces their continuum action by a midpoint
approximation. Inclusion of the admissible paths alone does not imply
\(\mathsf G_{\Lc}\le\mathsf G_{\Lc,N}\). The theorem below instead
establishes \(\mathsf G_{\Lc,N}\to\mathsf G_{\Lc}\) as
\(N\to\infty\).
\subsection{Assumptions on the operator}
\label{geo:sec:operator-hypotheses}
Although the motivating operator is the Laplacian
\(\partial_{xx}\), the time-convergence argument applies to more
general Fourier multipliers. We impose the following hypotheses on
their symbols.
\begin{enumerate}[label=\textup{(G\arabic*)},ref=G\arabic*]
\item\label{geo:hyp:real}
  \(\mathfrak p(-k)=\overline{\mathfrak p(k)}\) and \(\mathfrak p(0)=0\).
\item\label{geo:hyp:order}
  There are \(C>0\) and an integer \(r\ge1\) such that
  \[
    |\mathfrak p(k)|\le C(1+|k|)^r
    \qquad(k\in\mathbb Z).
  \]
\end{enumerate}
No semigroup hypothesis is needed for time convergence. For the upper
bound in the convergence proof, we approximate each continuous
finite-action admissible pair by feasible discrete pairs with the same
endpoints, converging to that pair and with no increase of action in
the limit. We call these approximations a \emph{recovery sequence},
and the procedure producing them a \emph{recovery construction}.
The construction in \cref{geo:lem:time-average-recovery} averages
only in time. It applies to measure-valued curves, including singular
marginals and vacuum, without a spatial density lower bound.
Let \(\mathcal X=C_{\mathrm w}([0,1];\Pp(\T))\).  The topology of uniform weak
convergence on \(\mathcal X\) is induced by
\[
  d_{\mathcal X}(\mu,\mu')
  =
  \sup_{t\in[0,1]}d_{-(r+2)}(\mu_t,\mu'_t),
  \qquad
  d_{-(r+2)}(\lambda,\lambda')
  =
  \sup_{\|\varphi\|_{C^{r+2}}\le1}
  \int_\T\varphi\,\dd(\lambda-\lambda').
\]
The unit ball of \(C^{r+2}(\T)\) is relatively compact in
\(C(\T)\) (Arzel\`a--Ascoli), so weak convergence in \(\Pp(\T)\) implies convergence in \(d_{-(r+2)}\).
Conversely, convergence in \(d_{-(r+2)}\) implies weak convergence, because \(C^{r+2}(\T)\) is dense in \(C(\T)\) and all the measures have mass one.
By \ref{geo:hyp:order}, \(\|\Lc^*\varphi\|_\infty\le C_{\Lc}\|\varphi\|_{C^{r+2}}\) for \(\varphi\in C^{r+2}(\T)\), with the constant \(C_{\Lc}=C\sum_{k\ne0}(1+|k|)^r(2\pi|k|)^{-(r+2)}<\infty\) of \eqref{geo:eq:Lstar-bound}.
Hence every admissible pair satisfies \(d_{-(r+2)}(\mu_t,\mu_s)\le C_{\Lc}\int_s^t|m_u|(\T)\,\dd u\) for \(0\le s\le t\le1\).
Since \(\mu_u\) is a probability measure, \(|m_u|(\T)\le(2\Ac(\mu_u,m_u))^{1/2}\), so the Cauchy--Schwarz inequality bounds the right-hand side by \(C_{\Lc}\bigl(2\int_0^1\Ac(\mu_u,m_u)\,\dd u\bigr)^{1/2}(t-s)^{1/2}\). Curves of bounded action are therefore equicontinuous for \(d_{-(r+2)}\); the discrete counterpart \eqref{geo:eq:holder-estimate} of this bound gives the equicoercivity in \cref{geo:thm:time-discrete-convergence}.
To compare the cost of individual paths, rather than only the minimum
values, we introduce two functionals on the same curve space
\(\mathcal X\). For \(\mu\in\mathcal X\), set \(\mathcal J_{\mathrm g}(\mu)=\inf\{\int_0^1\Ac(\mu_t,m_t)\,\dd t:\ m\text{ such that }(\mu,m)\text{ is admissible}\}\) if \(\mu_0=\mu^0\) and \(\mu_1=\mu^1\), and \(\mathcal J_{\mathrm g}(\mu)=+\infty\) otherwise.
Set
\[
  \mathcal J_{\mathrm g,N}(\mu)
  =
  \inf\left\{
    \sum_{j=0}^{N-1}\frac1h\Ac(\bar\mu_h^j,m^j):
    m^j\in\M(\T),\ \Lc m^j=\mu_{t_{j+1}}-\mu_{t_j}
  \right\},
  \qquad
  \bar\mu_h^j=\frac{\mu_{t_j}+\mu_{t_{j+1}}}2,
\]
if \(\mu_0=\mu^0\), \(\mu_1=\mu^1\) and \(t\mapsto\mu_t\) is affine on every \([t_j,t_{j+1}]\), and \(\mathcal J_{\mathrm g,N}(\mu)=+\infty\) otherwise.
The index \(\mathrm g\) refers to the geodesic problem.  Thus
\[
  \inf_{\mathcal X}\mathcal J_{\mathrm g}=\mathsf G_{\Lc}(\mu^0,\mu^1),
  \qquad
  \inf_{\mathcal X}\mathcal J_{\mathrm g,N}=\mathsf G_{\Lc,N}(\mu^0,\mu^1).
\]
In this notation, a recovery sequence for a path \(\mu\) satisfies
\[
 \mu_N\longrightarrow\mu\quad\text{in }\mathcal X,
 \qquad
 \limsup_{N\to\infty}\mathcal J_{\mathrm g,N}(\mu_N)
 \le\mathcal J_{\mathrm g}(\mu).
\]
Together with the lower bound for arbitrary convergent discrete paths,
this is the \(\Gamma\)-convergence asserted below. Its role is to
justify the discrete optimization problem: combined with compactness,
it yields both convergence of the minimum actions and convergence of
subsequences of minimizing paths to continuum geodesics.

\subsection{Convergence of actions and minimizing paths}
\begin{theorem}[Geodesic time limit]
\label{geo:thm:time-discrete-convergence}
Assume \ref{geo:hyp:real}--\ref{geo:hyp:order}, and suppose that
\(\mathsf G_{\Lc}(\mu^0,\mu^1)<\infty\).  Then
the family \(\mathcal J_{\mathrm g,N}\) is equicoercive and \(\Gamma\)-converges to \(\mathcal J_{\mathrm g}\) on
\(\mathcal X\). Its minimum values satisfy
\begin{equation}\label{geo:eq:value-convergence}
  \lim_{N\to\infty}\mathsf G_{\Lc,N}(\mu^0,\mu^1)
  =
  \mathsf G_{\Lc}(\mu^0,\mu^1).
\end{equation}
The continuum infimum and, for all sufficiently large \(N\), the discrete
infima are attained.  Every cluster point of an asymptotically minimizing
sequence of discrete paths is a continuum minimizer.
\end{theorem}
\begin{proposition}[The diffusive operator]
\label{geo:prop:heat-semigroup-hypotheses}
For \(\Lc=\partial_{xx}\), the symbol \(\mathfrak p(k)=-4\pi^2k^2\) satisfies \ref{geo:hyp:real} and \ref{geo:hyp:order} with \(r=2\).
Moreover \(\mathsf G_{\Lc}(\mu^0,\mu^1)<\infty\) for all \(\mu^0,\mu^1\in\Pp(\T)\).
Consequently \cref{geo:thm:time-discrete-convergence} applies to the
diffusive action for arbitrary probability endpoints.
\end{proposition}
The detailed proofs are in \cref{app:merged-geodesics}. Heat paths are
used only to establish finite signed action for arbitrary endpoints in
\cref{geo:prop:heat-semigroup-hypotheses}, not to construct the recovery
sequence for \cref{geo:thm:time-discrete-convergence}.

\begin{remark}[Nonnegative momentum and the real line]
\label{geo:rem:positive-extension}
The recovery construction preserves nonnegative momentum, and weak
limits preserve this constraint. The same convergence argument
therefore applies with $m\ge0$ on $\T$, under the hypothesis that the
minimum action with $m\ge0$ is finite.
\Cref{geo:prop:heat-semigroup-hypotheses} does not give this
hypothesis, because its proof uses momentum of both signs. On $\R$ with
$\Lc=\partial_{xx}$, one must additionally prove tightness of the
marginals and momenta and retain finite second moments. For prescribed geodesic
endpoint marginals in $\Ptwo(\R)$, convex-order domination supplies these
bounds; \cref{mart:cor:positive-midpoint} gives the resulting
martingale corollary and its proof verifies these noncompactness checks.
\end{remark}

\Cref{geo:thm:time-discrete-convergence} assumes no positive lower bound on the densities of \(\mu_t\).
Its conclusions are the convergence \eqref{geo:eq:value-convergence} of the minimum values and the uniform weak convergence in \(t\) of subsequences of asymptotically minimizing discrete paths.
It gives no rate.
It does not establish convergence of the signed diffusivity
\(\kappa_t=\dd m_t/\dd\mu_t\). In the numerical scheme of
Part~\ref{part:numerical}, its grid approximation is the quotient
\(\kappa_i^j=m_i^j/(h q_i^j)\) in \eqref{eq:kappa}, where
\(m_i^j/h\) is the momentum rate and \(q_i^j\) is the midpoint
density. Weak convergence of the marginals does not control this
quotient, particularly where the density becomes small.

\subsection{Joint convergence of the Fourier-grid geodesic problem}
\label{geo:sec:joint-grid}
We now discretize space as well as time, for $\Lc=\partial_{xx}$ on
$\T$. The result concerns the signed Fourier scheme used in
Part~\ref{part:numerical}. The spatial and temporal meshes may tend
to zero independently.

\begin{definition}[Fourier-grid operators]
\label{geo:def:fourier-grid-operators}
For $M\ge3$, set $\Delta x=1/M$, $x_i=i/M$ and
$\mathcal I_M=\{-\lfloor M/2\rfloor,\ldots,\lceil M/2\rceil-1\}$.
The atomic reconstruction of a grid vector $u$ is
\begin{equation}\label{geo:eq:grid-reconstruction}
 \iota_Mu:=\Delta x\sum_{i=0}^{M-1}u_i\delta_{x_i}.
\end{equation}
The self-adjoint discrete Laplacian $L_M$ has Fourier symbol
\[
 \widehat L_M(k)=-4\pi^2k^2,\qquad k\in\mathcal I_M.
\]
The Fej\'er kernel $F_K$, the positive averaging kernel $H_M$ and
its grid map $P_M$ are defined by
\begin{equation}\label{geo:eq:fejer-kernel}
 \begin{aligned}
 F_K(x)&=\frac1{K+1}
       \left|\sum_{k=0}^{K}e^{2\pi\mathrm i kx}\right|^2
       =\sum_{|k|\le K}\left(1-\frac{|k|}{K+1}\right)e^{2\pi\mathrm i kx},\\
 K_M&=\left\lfloor\frac{M-1}{4}\right\rfloor,\qquad
 H_M=(1-M^{-1})F_{K_M}+M^{-1},\\
 (P_M\lambda)_i&=\int_\T H_M(x_i-y)\,\dd\lambda(y)
 \quad\text{for a finite signed measure }\lambda.
 \end{aligned}
\end{equation}
\end{definition}
Here $\delta_{x_i}$ is the unit point mass at $x_i$, so $\iota_Mu$
places mass $\Delta x\,u_i$ at each grid point. A nonnegative vector
of grid mass $\Delta x\sum_i u_i=1$ reconstructs a probability measure.
The map $P_M$ evaluates the smooth convolution $H_M*\lambda$ at the
grid points; it does not require $\lambda$ to have a density.

Its grid perspective action is
\begin{equation}\label{geo:eq:grid-action}
 a(q,p)=\begin{cases}p^2/(2q),&q>0,\\0,&q=p=0,\\
 +\infty,&\text{otherwise},\end{cases}\qquad
 \A_M(q,p)=\Delta x\sum_i a(q_i,p_i).
\end{equation}
For $N\ge2$ and prescribed positive grid endpoints of mass one, let
$\mathcal J_{\mathrm g,M,N}(\lambda)$ be the infimum of
\begin{equation}\label{geo:eq:joint-grid-functional}
 \begin{aligned}
 &\sum_{j=0}^{N-1}\frac1h\A_M(q^j,m^j),\qquad
 q^j=\frac{\nu^j+\nu^{j+1}}2,\\
 &L_Mm^j=\nu^{j+1}-\nu^j,\qquad h=1/N,
 \end{aligned}
\end{equation}
over momenta, if $\lambda$ is the affine interpolation of
$\lambda(jh)=\iota_M\nu^j$, all $\nu^j$ are nonnegative of grid mass
one, and the endpoints agree with those prescribed. Set the functional
to $+\infty$ otherwise. Its minimum on $\mathcal X$ is denoted by
$E_{M,N}$. This is exactly the problem \eqref{eq:geodesic}.
Explicitly, the reconstructed path on the $j$-th time cell is
\[
 \lambda(t)=(1-\theta)\iota_M\nu^j+\theta\iota_M\nu^{j+1},
 \qquad \theta=\frac{t-jh}{h},\quad t\in[jh,(j+1)h].
\]
We call this the atomic affine reconstruction of the grid path.

The Fej\'er kernel $F_K$ is nonnegative and has integral one.
The added constant makes $H_M$ strictly positive, and its Fourier
support remains below the grid cutoff. The averaging attenuates the
retained Fourier modes. Its useful exact properties are mass
preservation, commutation with the conservation operator, and
contraction of the action; see \cref{geo:lem:positive-fourier-map}.

\begin{theorem}[Joint Fourier-grid geodesic limit]
\label{geo:thm:joint-grid}
Let $\mu^0,\mu^1\in\Pp(\T)$, and prescribe the grid endpoints
$\nu^0=P_M\mu^0$, $\nu^N=P_M\mu^1$ in
\eqref{geo:eq:joint-grid-functional}. For every sequence
$M_n,N_n\to\infty$, the functionals
$\mathcal J_{\mathrm g,M_n,N_n}$ are equicoercive and
$\Gamma$-converge on $\mathcal X$ to $\mathcal J_{\mathrm g}$ for
$\Lc=\partial_{xx}$. In particular,
\begin{equation}\label{geo:eq:joint-grid-value-limit}
 \lim_{n\to\infty}E_{M_n,N_n}
 =\mathsf G_{\partial_{xx}}(\mu^0,\mu^1).
\end{equation}
Each $\mathcal J_{\mathrm g,M_n,N_n}$ attains its minimum.
Every sequence of minimizers $\lambda_n$ has a subsequence converging
uniformly weakly in time, and every such subsequential limit is a
minimizer of $\mathcal J_{\mathrm g}$. The same conclusion holds for
asymptotically minimizing paths.
\end{theorem}

No density floor or relation between $\Delta x_n$ and $h_n$ is assumed.
The projection of the endpoint marginals with $P_M$ also applies to
singular measures. In computations, one often specifies the grid
endpoints directly, for example by sampling and normalizing positive
density profiles as in Part~\ref{part:numerical}. These vectors need
not equal $P_M\mu^0$ and $P_M\mu^1$. The following corollary ensures
that the same convergence result applies to such endpoint choices.
\begin{corollary}[Positive prescribed endpoint approximations]
\label{geo:cor:grid-endpoints}
Let $\mu^0,\mu^1\in\Pp(\T)$, and prescribe, for each $M$, initial and
final grid vectors $\nu_M^0$ and $\nu_M^N$ of grid mass one; they do not
depend on $N$. Their atomic
reconstructions use $\iota_Mu=\Delta x\sum_i u_i\delta_{x_i}$,
with $\Delta x=1/M$ and $x_i=i/M$, as in
\eqref{geo:eq:grid-reconstruction}. Suppose
\begin{equation}\label{geo:eq:grid-endpoint-assumptions}
 \nu_M^0,\nu_M^N\ge c\one\quad\text{for some }c>0\text{ independent of }M,
 \qquad
 \iota_M\nu_M^0\rightharpoonup\mu^0,
 \quad \iota_M\nu_M^N\rightharpoonup\mu^1.
\end{equation}
Then all conclusions of \cref{geo:thm:joint-grid} hold with these
vectors in place of the projections $P_M\mu^0$ and $P_M\mu^1$ of the
endpoint marginals.
\end{corollary}
The common positive lower bound is required only at the two prescribed
endpoints. The intermediate densities may vanish.
Since $\iota_M\nu_M^0\ge c\,\Delta x\sum_i\delta_{x_i}$ and
$c\,\Delta x\sum_i\delta_{x_i}\rightharpoonup c\,\dd x$, the hypotheses
imply $\mu^0\ge c\,\dd x$; similarly $\mu^1\ge c\,\dd x$. No other
regularity of the endpoint measures is assumed.
The proof in \cref{geo:app:joint-grid} combines atomic compactness and
perspective lower semicontinuity with the time recovery of
\cref{geo:lem:time-average-recovery}. Positive Fourier averaging
preserves the grid constraint exactly and, by convexity of the
perspective action, does not increase the action.
For \cref{geo:cor:grid-endpoints}, short signed-momentum bridges join
the prescribed endpoint vectors to the projections of the marginals
with $P_M$ at vanishing cost.

\section{Convergence of the signed modified JKO scheme}\label{sec:merged-jko}
We next study a minimizing movement scheme for the general energy
\eqref{jko:eq:energy}. We choose to weight the action at each step by
the unknown measure being optimized. Our main result for this
discretization proves convergence under a positive lower bound on the
optimized densities, uniform in the time step. We call this assumption
the density floor \eqref{jko:eq:floor}. We then collect results that
also apply when momentum is constrained to be nonnegative.

\subsection{Perspective action and frozen cost}
Let $d\ge1$ and $\Torus=(\R/\mathbb Z)^d$, with normalized Lebesgue measure. We identify a density $\rho$ with the probability measure $\rho\,\dd x$ when it appears in a transport cost.
For a matrix-valued distribution $\mathsf m=(\mathsf m_{ij})_{1\leq i,j\leq d}$, set
\[
  \DivDiv \mathsf m:=\sum_{i,j=1}^d \partial_{ij}\mathsf m_{ij}.
\]
The norm $|\mathsf m|$ and the product $\mathsf m:\mathsf n=\sum_{i,j}\mathsf m_{ij}\mathsf n_{ij}$ of two matrices are the Frobenius norm and product.

For a nonnegative finite measure $\sigma$ and a finite matrix-valued measure
$\mathsf m$, define the quadratic perspective action
\begin{equation}\label{jko:eq:perspective-action}
  \Ac(\sigma,\mathsf m):=
  \begin{cases}
    \displaystyle \frac12\int_{\Torus}|\mathsf v|^2\dd\sigma,
      & \mathsf m=\mathsf v\sigma,\quad \mathsf v\in L^2(\sigma;\R^{d\times d}),\\[4pt]
    +\infty,&\text{otherwise}.
  \end{cases}
\end{equation}
For \(d=1\), \eqref{jko:eq:perspective-action} coincides with \eqref{geo:eq:action-density}, with \(\mathsf v=\dd\mathsf m/\dd\sigma\).
For probability measures $\mu,\sigma\in\Pcal(\Torus)$, introduce the
frozen action, with the same normalization as
$\mathsf G_{\Lc}$ in \eqref{geo:eq:continuum-value},
\begin{equation}\label{jko:eq:frozen-cost}
  \mathsf G_{\fr}(\mu,\sigma)
  :=\inf\bigl\{\Ac(\sigma,\mathsf m):
       \mathsf m\in\Mcal(\Torus;\R^{d\times d}),\ 
       \DivDiv \mathsf m=\sigma-\mu\bigr\}.
\end{equation}
The equality in the constraint is distributional. The subscript
$\fr$ abbreviates ``frozen'': the action is weighted by the final
measure $\sigma$ throughout the time cell, rather than by a varying
marginal along a path. It uses a single integrated momentum and is
generally asymmetric. The geodesic action
minimizes over a whole connecting path. No triangle inequality is used
below.

\subsection{Energy, time step and positivity floor}
We consider energies of the form
\begin{equation}\label{jko:eq:energy}
  \Ec(\rho)
  :=\int_{\Torus}V_{\rm ext}\rho\dd x
    +\frac12\int_{\Torus}(W*\rho)\rho\dd x
    +\int_{\Torus}f(\rho)\dd x,
\end{equation}
where $\Ec(\rho)=+\infty$ if the probability measure is not absolutely
continuous or the last integral is not finite.  We make the following
assumptions.

\begin{enumerate}
\item[\textup{(E1)}] $V_{\rm ext},W\in C^2(\Torus)$ and $W$ is even.  Thus the first variation
  of the interaction energy is $W*\rho$.
\item[\textup{(E2)}] There are $A_0,A_1\in\R$ such that $f+A_0+A_1r$ maps $[0,+\infty)$ into $[0,+\infty]$,
  is lower semicontinuous, strictly convex,
  superlinear, belongs to $C^2((0,+\infty))$, and has second derivative
  $f''(r)>0$ for $r>0$.
  Since $\int_{\Torus}\rho\dd x=1$, replacing $f$ by $f+A_0+A_1r$ adds the constant $A_0+A_1$ to $\Ec$ on $\Pcal(\Torus)$ and does not change the minimizers in \eqref{jko:eq:scheme}.
  We assume this replacement is made, and \textup{(E3)} refers to the new $f$.
\item[\textup{(E3)}] For every $c>0$, there are $q_c>1$ and $C_c<+\infty$ such that
  \begin{equation}\label{jko:eq:fprime-growth}
    |f'(r)|^{q_c}\leq C_c\bigl(1+f(r)\bigr)
    \qquad\text{for all }r\geq c.
  \end{equation}
\end{enumerate}
The shifted Boltzmann entropy
\begin{equation}\label{jko:eq:boltzmann}
  f(r)=r\log r-r+1,\qquad f'(r)=\log r,
\end{equation}
with the usual convention at $r=0$, satisfies these assumptions.  

Given $\tau>0$ and $\rho_0\in\Pcal(\Torus)$ with
$\Ec(\rho_0)<+\infty$, choose recursively
\begin{equation}\label{jko:eq:scheme}
  \rho_{k+1}^\tau\in\operatorname*{argmin}_{\rho\in\Pcal(\Torus)}
  \left\{
    \frac{1}{\tau}\mathsf G_{\fr}(\rho_k^\tau,\rho)
    +\Ec(\rho)
  \right\},
  \qquad \rho_0^\tau=\rho_0.
\end{equation}
The minimum in \eqref{jko:eq:scheme} is attained whenever \(\Ec(\rho_k^\tau)<+\infty\) (\cref{jko:prop:existence}).
The competitor \(\rho=\rho_k^\tau\) has cost \(\mathsf G_{\fr}(\rho_k^\tau,\rho_k^\tau)=0\), hence \(\Ec(\rho_{k+1}^\tau)\le\Ec(\rho_k^\tau)\) and the recursion is defined for every \(k\).

We assume the following local-in-time positivity floor:
\begin{equation*}\tag{P}\label{jko:eq:floor}
  \begin{aligned}
  &\text{for every $\tfin>0$ there are $c_{\tfin}>0$ and $\tau_{\tfin}>0$ such that}\\
  &\rho_k^\tau\geq c_{\tfin}\ \text{a.e.\ on }\Torus\quad\text{for all }\tau\in(0,\tau_{\tfin}]\text{ and all $k$ with }k\tau\leq\tfin.
  \end{aligned}
\end{equation*}
For $k=0$, \eqref{jko:eq:floor} requires $\rho_0\geq c_{\tfin}$.
This is a strong assumption that must be verified case by case.
A usual strategy in JKO convergence arguments is first to establish
the corresponding lower bound for the continuous evolution, then to
prove an analogous estimate for the discrete minimizers, uniformly
in the time step. The second step requires a separate argument:
positivity of the continuous solution alone does not establish
\eqref{jko:eq:floor}. We leave verification of this assumption open
for nonconstant initial data, so \cref{jko:thm:main} is conditional
on \eqref{jko:eq:floor}.

Define the backward piecewise-constant interpolation
\begin{equation}\label{jko:eq:interpolation}
  \bar\rho_\tau(t):=\rho_{k+1}^\tau
  \quad\text{for }t\in(k\tau,(k+1)\tau],
  \qquad \bar\rho_\tau(0):=\rho_0.
\end{equation}
With \eqref{jko:eq:floor}, the growth condition \eqref{jko:eq:fprime-growth} makes \(f'(\rho_{k+1}^\tau)\) integrable and allows differentiation under the integral in the optimality conditions of each step (\cref{jko:prop:euler-system}).
For \(\tfin>0\), set \(\hat c_{\tfin}:=c_{\tfin+1}\) and \(\hat\tau_{\tfin}:=\min\{1,\tau_{\tfin+1}\}\), where \(c_{\tfin+1}\) and \(\tau_{\tfin+1}\) are the constants of \eqref{jko:eq:floor} on \([0,\tfin+1]\).
If \(\tau\in(0,\hat\tau_{\tfin}]\) and \(t\in(k\tau,(k+1)\tau]\) with \(t<\tfin\), then \((k+1)\tau<\tfin+\tau\leq\tfin+1\); hence \(\bar\rho_\tau\geq\hat c_{\tfin}\) a.e.\ on \((0,\tfin)\times\Torus\).
With \eqref{jko:eq:fprime-growth} for \(c=\hat c_{\tfin}\) and \(\Ec(\bar\rho_\tau(t))\leq\Ec(\rho_0)\), this floor bounds \(f'(\bar\rho_\tau)\) in \(L^{q_{\hat c_{\tfin}}}((0,\tfin)\times\Torus)\) uniformly in \(\tau\in(0,\hat\tau_{\tfin}]\), which identifies \(f'(\rho)\) in the limit.

\subsection{Limit equation and convergence}
\begin{theorem}[Continuous-space time-step limit]\label{jko:thm:main}
Assume \emph{(E1)--(E3)}, $\Ec(\rho_0)<+\infty$ and \eqref{jko:eq:floor}.  Let
$\tau_n\downarrow0$.  There are a subsequence, not relabelled, and a density
$\rho$ such that, for every $\tfin>0$,
\begin{align}
  \bar\rho_{\tau_n}&\longrightarrow\rho
    &&\text{strongly in }L^1((0,\tfin)\times\Torus),\label{jko:eq:strong-main}\\
  \bar\rho_{\tau_n}(t)&\longrightarrow\rho(t)
    &&\text{uniformly in $t\in[0,\tfin]$ for the weak topology of measures}.
    \label{jko:eq:weak-uniform-main}
\end{align}
Moreover, for every $\tfin>0$, $f'(\rho)\in L^1((0,\tfin)\times\Torus)$ and, with
\begin{equation}\label{jko:eq:chemical-potential}
  \xi:=V_{\rm ext}+W*\rho+f'(\rho),
\end{equation}
one has
\begin{equation}\label{jko:eq:limit-regularity}
  \sqrt\rho\,\nabla^2\xi\in
  L^2((0,\tfin)\times\Torus;\R^{d\times d}),
\end{equation}
and $\rho$ is a weak solution of
\begin{equation}\label{jko:eq:limit-pde}
  \partial_t\rho
  +\DivDiv\bigl(\rho\nabla^2\xi\bigr)=0,
  \qquad \rho|_{t=0}=\rho_0.
\end{equation}
Precisely, for every
$\zeta\in C_c^\infty([0,+\infty)\times\Torus)$,
\begin{equation}\label{jko:eq:weak-formulation}
  \int_0^\infty\!\!\int_{\Torus}
       \rho\,\partial_t\zeta\dd x\dd t
  +\int_{\Torus}\rho_0\zeta(0,\cdot)\dd x
  =\int_0^\infty\!\!\int_{\Torus}
       \rho\nabla^2\xi:\nabla^2\zeta\dd x\dd t.
\end{equation}
For almost every $t>0$, the limit also satisfies
\begin{equation}\label{jko:eq:limit-energy-inequality}
 \Ec(\rho(t))
 +\frac12\int_0^t\!\!\int_{\Torus}
     \rho|\nabla^2\xi|^2\dd x\dd s
 \leq\Ec(\rho_0).
\end{equation}
If \eqref{jko:eq:limit-pde} has at most one weak solution, in the sense of \eqref{jko:eq:weak-formulation}, that satisfies \(f'(\rho)\in L^1((0,\tfin)\times\Torus)\) and \eqref{jko:eq:limit-regularity} for every \(\tfin>0\), and \eqref{jko:eq:limit-energy-inequality}, then \(\bar\rho_\tau\to\rho\) as \(\tau\downarrow0\) in the senses \eqref{jko:eq:strong-main} and \eqref{jko:eq:weak-uniform-main}.
\end{theorem}

\begin{corollary}[Entropy and the DLSS equation]\label{jko:cor:dlss}
Under the hypotheses of \cref{jko:thm:main}, including the positivity floor
\eqref{jko:eq:floor}, the choice \eqref{jko:eq:boltzmann} gives the chemical potential
\[
  \xi=V_{\rm ext}+W*\rho+\log\rho.
\]
In particular, if $V_{\rm ext}=W=0$, the limit equation is
\begin{equation}\label{jko:eq:dlss}
  \partial_t\rho
  +\DivDiv\bigl(\rho\nabla^2\log\rho\bigr)=0.
\end{equation}
In one dimension this reads $\partial_t\rho=-\partial_{xx}(\rho\,\partial_{xx}\log\rho)$.
\end{corollary}

Under \eqref{jko:eq:floor}, let \(\varphi_{k+1}^\tau\) be the mean-zero solution of \(\DivDiv(\rho_{k+1}^\tau\nabla^2\varphi_{k+1}^\tau)=(\rho_{k+1}^\tau-\rho_k^\tau)/\tau\), let \(\varphi_\tau(t)=\varphi_{k+1}^\tau\) for \(t\in(k\tau,(k+1)\tau]\), and let \(r_\tau=\frac\tau2|\nabla^2\varphi_\tau|^2\) denote the defect.
The proof in \cref{app:merged-jko} first obtains an action estimate and shows that \(r_\tau\to0\) in \(L^1((0,\tfin)\times\Torus)\) as \(\tau\to0\).
The positivity floor \eqref{jko:eq:floor} converts the action estimate into an \(L^2(0,\tfin;H^2(\Torus))\) bound for \(\varphi_\tau\).
The Euler--Lagrange system \eqref{jko:eq:euler-system} then gives strong compactness of the densities, which identifies the limiting evolution equation.

\begin{remark}[Modified JKO scheme and metric JKO scheme]
The geodesic action of the tensor model is
formally
\[
 \mathsf G_{\DivDiv}(\mu^0,\mu^1)
 :=\inf\left\{
   \int_0^1\Ac(\mu_s,\mathsf m_s)\dd s:
   \partial_s\mu_s=\DivDiv \mathsf m_s,
   \ \mu_{0}=\mu^0,\ \mu_1=\mu^1
 \right\}.
\]
For $d=1$, this is $\mathsf G_{\partial_{xx}}$ of
\eqref{geo:eq:continuum-value}, with squared distance
$\mathsf D^2=2\mathsf G_{\partial_{xx}}$.
Replacing \(\mathsf G_{\fr}\) by \(\mathsf G_{\DivDiv}\) in \eqref{jko:eq:scheme} gives the metric JKO scheme.
The cost \(\mathsf G_{\fr}\) of \eqref{jko:eq:frozen-cost} differs from \(\mathsf G_{\DivDiv}\) in two ways: the action is weighted by the final measure \(\sigma\), and a single time cell is used.
\Cref{jko:thm:main} is proved for \eqref{jko:eq:scheme} only.
One possible route to transfer \eqref{jko:eq:strong-main}--\eqref{jko:eq:limit-energy-inequality} to the metric JKO scheme is a comparison of \(\mathsf G_{\DivDiv}\) with \(\mathsf G_{\fr}\), and of their first variations in the final measure, at consecutive iterates of the metric JKO scheme.
A comparison of the values of the two costs alone does not give for the metric JKO scheme an analogue of the Euler--Lagrange system \eqref{jko:eq:euler-system}, which the proof of \cref{jko:thm:main} uses.
Whether such a comparison holds is open.
\end{remark}

\begin{remark}[Tensor versus scalar multidimensional models]
In one dimension, the tensor field is scalar and $\DivDiv=\partial_{xx}$.
In dimension \(d>1\), the constraint \(\DivDiv\mathsf m=\sigma-\mu\) leads to the operator \(\DivDiv(\rho\nabla^2\xi)\) in \eqref{jko:eq:limit-pde}; with \(\xi=\log\rho\), this operator gives the tensor form \eqref{jko:eq:dlss} of the DLSS equation \cite{DLSS91,MRSSdlss}.
A scalar momentum \(m\) with the constraint \(\Delta m=\sigma-\mu\) leads instead to \(\Delta(\rho\Delta\xi)\).
The two operators coincide for \(d=1\) and differ in general for \(d\ge2\) and nonconstant \(\rho\).
The one-dimensional formulas therefore extend to \(d>1\) only for a chosen model.
\end{remark}

\begin{remark}[Role of the positivity floor \eqref{jko:eq:floor}]
The positivity floor \eqref{jko:eq:floor} is used at six points: to define the potential of \cref{jko:lem:duality}, to justify two-sided density variations (\cref{jko:lem:first-variation,jko:prop:euler-system}), to turn the weighted action estimate
into an \(L^2(0,\tfin;H^2(\Torus))\) bound for \(\varphi_\tau\), to make the defect \(r_\tau\)
vanish, to transfer
compactness from $f'(\bar\rho_\tau)$ to $\bar\rho_\tau$, and to bound \(f'(\bar\rho_\tau)\) in \(L^{q_{\hat c_{\tfin}}}((0,\tfin)\times\Torus)\) uniformly in \(\tau\in(0,\hat\tau_{\tfin}]\).
This bound makes \(f'(\bar\rho_{\tau_n})\) uniformly integrable, so that Vitali's theorem gives \(f'(\bar\rho_{\tau_n})\to f'(\rho)\) in \(L^1((0,\tfin)\times\Torus)\).  Strict positivity of each
individual step does not by itself imply \eqref{jko:eq:floor}.
Convergence without \eqref{jko:eq:floor}, in particular when the limit density vanishes on a set of positive measure, is not proved here.
In addition, for the entropy \eqref{jko:eq:boltzmann}, the chemical potential \(\xi\) contains \(\log\rho\), which is not defined where \(\rho=0\).
For this entropy, a positive lower bound on the density is not needed to solve the DLSS equation \eqref{jko:eq:dlss} on a bounded convex domain of \(\R^d\) with Neumann boundary conditions.
Matthes, Savar\'e, and Schlichting \cite{MSS26} construct there, for every nonnegative integrable initial density, a unique solution in a weak sense that allows vacuum: conservation of mass, local integrability in time of \(\int(\Delta\sqrt\rho)^2/\sqrt\rho\,\dd x\), and one-sided inequalities tested against smooth, uniformly positive functions \cite[Theorem~4.14]{MSS26}.
Whether the limits of the modified JKO scheme on \(\Torus\) satisfy the periodic analogue of these conditions without \eqref{jko:eq:floor} is open.
\end{remark}

\subsection{Nonnegative momentum: existence and compactness}
\label{jko:sec:positive}
The existence and compactness techniques developed for signed
momentum also apply under the additional nonnegativity constraint on
the momentum. On $\Torus$, fix $\kappa>0$ and consider the
modified JKO scheme
\begin{equation}\label{jko:eq:positive-scheme}
 (\rho_{k+1}^{\tau,+},\mathsf m_{k+1}^{\tau,+})
 \in\operatorname*{argmin}_{\substack{\rho\in\Pcal(\Torus),\ \mathsf m\succeq0\\
                   \DivDiv\mathsf m=\rho\dd x-\rho_k^{\tau,+}\dd x}}
 \left\{\frac\kappa\tau\Ac(\rho\dd x,\mathsf m)+\Ec(\rho)\right\},
 \qquad \rho_0^{\tau,+}=\rho_0.
\end{equation}
Here $\mathsf m\succeq0$ means a symmetric positive-semidefinite
matrix-valued measure. The choice $\kappa=1$ has the normalization of
\eqref{jko:eq:scheme}. Let $\bar\rho_\tau^+$ be the backward constant
interpolation, and set $Q_\tau^+=\mathsf m_{k+1}^{\tau,+}/\tau$ on
$(k\tau,(k+1)\tau]$.

\begin{proposition}[Positive-momentum existence and compactness]
\label{jko:prop:positive-compactness}
Assume \emph{(E1)--(E2)} and $\Ec(\rho_0)<\infty$.
Every step of \eqref{jko:eq:positive-scheme} is attained. If
$\mathsf m_{k+1}^{\tau,+}=\tau B_{k+1}^\tau\rho_{k+1}^{\tau,+}\dd x$,
then $B_{k+1}^\tau\succeq0$ and
\begin{equation}\label{jko:eq:positive-energy-bound}
 \Ec(\rho_n^{\tau,+})+
 \frac\kappa2\sum_{k<n}\tau\int_{\Torus}
              |B_{k+1}^\tau|_F^2\rho_{k+1}^{\tau,+}\dd x
 \le\Ec(\rho_0).
\end{equation}
For every sequence $\tau\downarrow0$, a subsequence has a narrowly
continuous limit $\rho_t\dd x$, uniformly weakly on each bounded time
interval, and
\[
 \dd t\otimes Q_\tau^+\rightharpoonup
 B_t(x)\rho_t(x)\dd x\dd t,\qquad B\succeq0.
\]
The limit satisfies $\rho|_{t=0}=\rho_0$,
$\partial_t\rho=\DivDiv(\rho B)$, and, for every $t\ge0$,
\begin{equation}\label{jko:eq:positive-limit-energy}
 \Ec(\rho_t)+\frac\kappa2\int_0^t\!\!\int_{\Torus}
                  |B_s|_F^2\rho_s\dd x\dd s\le\Ec(\rho_0).
\end{equation}
Neither \emph{(E3)} nor the density floor \eqref{jko:eq:floor} is
needed for these conclusions. The proposition does not identify $B$
as a function of the chemical potential.
\end{proposition}

The proof reuses the direct method, the action estimate, the time
modulus and the perspective lower semicontinuity from the signed
argument; it is given in \cref{jko:app:positive-compactness}.
The missing estimate concerns the transport dual potential, the
Lagrange multiplier of the conservation constraint. In the signed
proof, the action controls its full Hessian. With nonnegative
momentum, only the nonnegative part of this Hessian contributes to
the dual quadratic penalty. In one dimension, this means that only
the positive part of its second derivative is controlled; the negative
part may contain point masses. Thus the full square-integrable
second-derivative bound used to identify the limiting equation in the
signed proof is unavailable, even for smooth positive densities; see
\cref{jko:rem:positive-dual-obstruction}.

For the one-dimensional potential energy
$\Ec(\mu)=\int_\R V_{\rm ext}\dd\mu$, the limiting equation is
identified. If $V_{\rm ext}\in C^3(\R)$ has bounded second and third
derivatives, \cref{mart:thm:positive-potential} proves convergence
of the positive-momentum scheme to the unique solution of
\[
 \partial_t\mu=\partial_{xx}(D\mu),\qquad
 D(x)=\frac{(-V_{\rm ext}''(x))_+}{\kappa},\qquad
 \mu|_{t=0}=\mu_0.
\]
Here $r_+=\max\{r,0\}$. This result holds for any initial
$\mu_0\in\Ptwo(\R)$, without a density floor.

\section{Two discretizations of the martingale geodesic problem}
\label{mart:sec:setting}
We return to the geodesic problem of
\cref{sec:merged-geodesics}, now with nonnegative momentum on $\R$,
and compare it with martingale couplings. The momentum scheme minimizes
a quadratic
perspective action under a local conservation equation and nonnegative
momentum. The coupling scheme minimizes fourth moments over successive
martingale transitions. They use different variables and have different
finite-step costs. In dimension one, both converge, after normalization,
to the same continuous martingale cost; see
\cref{mart:cor:positive-midpoint,mart:cor:common-limit}.

Huesmann and Trevisan \cite[Section~6]{HT19} proposed optimized
short-time coupling costs and asked about their convergence to a
continuous martingale transport problem. For the fourth-moment cost,
we determine the optimized geodesic limit with the normalization
in \eqref{mart:c:eq:bb-convergence}. The momentum proof extends the
signed time-averaging argument to nonnegative momentum on $\R$,
with the additional moment and tightness estimates stated below.
\Cref{mart:rem:static-proxy} explains the distinction between optimizing
increments and sampling a fixed diffusion, and shows that the limit
proposed in \cite[Section~6]{HT19} fails for this cost.

The momentum formulation avoids solving a sequence of optimal
coupling problems when only the time marginals are sought.
This is a practical advantage in one dimension. However, we report no
computation of the martingale problem on $\R$ with either scheme. The
computations with $m\ge0$ in Part~\ref{part:numerical} concern the
periodic problem, which is not a martingale problem
(\cref{rem:martingale}).
Efficient inversion and constrained optimization for
higher-dimensional tensor momentum require further work; see
\cref{mart:rem:computational-scope}. The gradient-flow problem is
treated separately in \cref{mart:g:sec:m4-theorem}.
All proofs are collected in \cref{mart:app:proofs}.

Write $\Ptwo(\R^d)$ for probability measures with finite second moment,
$\Sd$ for symmetric matrices, and $\Sdp$ for positive semidefinite matrices.
The symbols $\Am$ and $\Em$ denote the martingale action and energy;
$\mathbb E$ denotes expectation and $\Law$ denotes distribution.
\subsection{The continuous martingale problem}
\label{mart:g:sec:continuous}
Fix endpoint marginals $\mu,\sigma\in\Ptwo(\R^d)$ and a time horizon $\tfin>0$.
We seek a continuous square-integrable martingale $X$, on a filtered
probability space, with these endpoint marginals and absolutely continuous
quadratic variation
\begin{equation}\label{mart:eq:continuous-bracket}
 \Law(X_0)=\mu,\qquad \Law(X_{\tfin})=\sigma,\qquad
 \dd\langle X\rangle_t=\alpha_t\dd t,\qquad \alpha_t\in\Sdp.
\end{equation}
Here $\alpha_t$ is a predictable covariance rate; it need not be a
function of $X_t$.
\begin{definition}[Martingale Benamou--Brenier problem]
The continuous endpoint cost is
\begin{equation}\label{mart:g:eq:mb-cost}
 \MB^2_{\tfin}(\mu,\sigma)
 :=\inf_X\mathbb E\int_0^{\tfin}(\tr\alpha_t)^2\dd t,
\end{equation}
where the infimum runs over the martingales in
\eqref{mart:eq:continuous-bracket}. The value is $+\infty$ if no
finite-action competitor exists.
\end{definition}
In dimension one, $\alpha_t$ is a nonnegative scalar and the integrand
is $\alpha_t^2$. Finite action requires the endpoint marginals to
satisfy the convex-order constraint
\[
 \mu\le_c\sigma
 \quad\Longleftrightarrow\quad
 \int_{\R^d}\phi\,\dd\mu\le\int_{\R^d}\phi\,\dd\sigma
 \quad\text{for every convex test function }\phi
\]
for which the two integrals are defined. A deterministic time change gives
\begin{equation}\label{mart:g:eq:mb-scaling}
 \MB^2_{\tfin}(\mu,\sigma)=\frac1{\tfin}\MB^2_1(\mu,\sigma).
\end{equation}

\subsection{Eulerian formulation and momentum discretization}
\label{mart:sec:momentum-approach}
The Eulerian representation replaces the path distribution by its marginals
and a local covariance field. It connects the martingale problem to
the perspective action used for signed diffusion.
For an admissible pair on $[0,\tfin]$, the action and equation are
\begin{align}
 \Am_{\tfin}(\mu,A)&:=\int_0^{\tfin}\!\int_{\R^d}
      (\tr A_t)^2\dd\mu_t\dd t,
      \label{mart:g:eq:intro-action}\\
 \partial_t\mu&=\frac12\partial_{ij}(A_{ij}\mu).
      \label{mart:g:eq:martingale-fp}
\end{align}
Repeated matrix indices are summed.
The coefficient $1/2$ is the usual covariance normalization in It\^o's formula.
The admissibility conditions and weak formulation are specified below.
\begin{definition}[Admissible martingale curve]
Let $\tfin>0$. An admissible martingale curve is a pair $(\mu,A)$ such that $(\mu_t)_{t\in[0,\tfin]}\subset\Ptwo(\R^d)$ is narrowly continuous, $A:(0,\tfin)\times\R^d\to\Sdp$ is Borel, and
\begin{equation}
  \int_0^{\tfin}\!\int_{\R^d}(\tr A_t)^2\dd\mu_t\dd t<\infty.
  \label{mart:g:eq:finite-action}
\end{equation}
For every $\varphi\in C_c^{1,2}([0,\tfin]\times\R^d)$, we require
\begin{equation}\label{mart:g:eq:weak-fp}
 \int\varphi(\tfin,x)\dd\mu_{\tfin}
 -\int\varphi(0,x)\dd\mu_0
 =\int_0^{\tfin}\!\int
   \left(\partial_t\varphi+\frac12A_t:D^2\varphi\right)
   \dd\mu_t\dd t.
\end{equation}
\end{definition}

The superposition principle \cite[Theorem~2.5]{Tre16} and conditional
projection give
\begin{equation}\label{mart:eq:martingale-eulerian-value}
 \MB^2_{\tfin}(\mu,\sigma)
 =\inf\{\Am_{\tfin}(\mu_t,A_t):
       (\mu,A)\text{ admissible with endpoints }\mu,\sigma\}.
\end{equation}
Indeed, the superposition representation has bracket density
$A_t(X_t)$. Conversely, projecting an adapted covariance rate gives
$A_t(x)=\mathbb E[\alpha_t\mid X_t=x]$. Since
$A\mapsto(\tr A)^2$ is convex, conditional Jensen shows that this
projection cannot increase the trace-square action. Thus optimization of the marginal
curve can be carried out entirely in Eulerian variables.
In dimension one, the common momentum variable is
$m_t=\frac12a_t\mu_t$. The nonnegative sign of the covariance then
becomes the additional constraint $m_t\ge0$ in the signed Eulerian
action.

\begin{definition}[One-dimensional positive Eulerian problem]
\label{mart:def:positive-eulerian}
For $\mu,\sigma\in\Ptwo(\R)$ and $\tfin>0$, set
\begin{equation}\label{mart:eq:positive-eulerian}
 \begin{aligned}
 \mathsf G^+_{[0,\tfin]}(\mu,\sigma)
 :=\inf\Bigl\{&\int_0^{\tfin}\Ac(\mu_t,m_t)\dd t:\\
 &\partial_t\mu_t=\partial_{xx}m_t,\quad m_t\ge0,\\
 &\mu_0=\mu,\quad\mu_{\tfin}=\sigma\Bigr\}.
 \end{aligned}
\end{equation}
Here $\mu_t\in\Ptwo(\R)$ is narrowly continuous, $m_t$ is a weakly
measurable family of finite nonnegative measures with integrable mass,
and the equation holds distributionally. We use the same perspective
normalization as \eqref{geo:eq:action-density}, now on $\R$:
\begin{equation}\label{mart:eq:positive-perspective}
 \Ac(\mu_t,m_t)=\frac12\int_\R
       \left|\frac{\dd m_t}{\dd\mu_t}\right|^2\dd\mu_t
 \quad\text{if }m_t\ll\mu_t,
\end{equation}
and $+\infty$ otherwise. Thus, for densities, the integrand is
$m_t^2/(2\rho_t)$. This is the signed quadratic Eulerian problem on the
same state space with the additional closed constraint $m_t\ge0$.
\end{definition}

\begin{proposition}[The common one-dimensional continuous action]
\label{mart:prop:positive-eulerian}
In dimension one, the substitution
\begin{equation}\label{mart:eq:momentum-covariance}
 m_t=\frac12a_t\mu_t,\qquad a_t\ge0,
 \qquad
 \int_0^{\tfin}\Ac(\mu_t,m_t)\dd t
 =\frac18\int_0^{\tfin}\!\int a_t^2\dd\mu_t\dd t
\end{equation}
identifies the finite-action Eulerian pairs with the admissible pairs
of \eqref{mart:g:eq:weak-fp}. In particular,
\begin{equation}\label{mart:eq:common-continuous-value}
 8\mathsf G^+_{[0,\tfin]}(\mu,\sigma)
 =\MB^2_{\tfin}(\mu,\sigma)
 =\inf_X\mathbb E\int_0^{\tfin}\alpha_t^2\dd t,
 \qquad \dd\langle X\rangle_t=\alpha_t\dd t,
\end{equation}
where the last infimum runs over continuous square-integrable
martingales with endpoint marginals $\mu,\sigma$. Every finite-action
Eulerian pair has such a representation with
$\dd\langle X\rangle_t=a_t(X_t)\dd t$.
\end{proposition}

The coefficient $a_t$ is a variance rate, whereas $m_t$ is a momentum
measure. The equality of values in
\eqref{mart:eq:common-continuous-value} does not identify a unique path distribution. Again, by convexity, for an arbitrary martingale representation, conditional projection
of $\alpha_t$ given $X_t=x$ can decrease the quadratic action.

\subsubsection{Direct time discretization with nonnegative momentum}
The midpoint scheme with nonnegative momentum corresponding to
\eqref{geo:eq:time-discrete-value}, now on $\R$, has variables
$(\mu^k,m^k)$ and value
\begin{equation}\label{mart:eq:positive-midpoint}
 \begin{aligned}
 \mathsf G^+_N(\mu,\sigma):=
 \inf\Bigl\{&\sum_{k=0}^{N-1}\frac1h\Ac(\bar\mu^k,m^k):\\
 &\partial_{xx}m^k=\mu^{k+1}-\mu^k,\quad m^k\ge0,\\
 &\mu^0=\mu,\quad\mu^N=\sigma\Bigr\},
 \end{aligned}
\end{equation}
where $h=1/N$ and $\bar\mu^k=(\mu^k+\mu^{k+1})/2$.
All nodes belong to $\Ptwo(\R)$, and each $m^k$ is a finite
nonnegative measure.
The measure $m^k$ is integrated over a time cell, so the momentum
per unit time is $m^k/h$. The local second-derivative equation and the
sign condition are imposed within the discrete optimization itself.
This is a discretization of diffusion in the Eulerian variables; no
transition kernel or jump mechanism is optimized. It should be
distinguished from asserting that every finite-step interpolation
already has a finite-action stochastic representation.

The underlying Eulerian model therefore builds in the local diffusion
property. Corollary~\ref{mart:cor:positive-midpoint} proves convergence of
this midpoint scheme, including recovery with the exact endpoints and
nonnegative momentum. The proof averages only in time and does not
require positive densities or spatial smoothing.
The preceding signed theory and the numerical computations are on a
torus, whereas \eqref{mart:eq:positive-midpoint} is on $\R$.
Some adaptations of the arguments are required on this unbounded
domain, in particular to control moments and prove tightness.

\subsection{Coupling discretization and the fourth-moment proxy}
\label{mart:g:sec:discrete}
Here the variables are the marginals and martingale couplings
between successive times.
The set $\Pi_M(\mu,\sigma)$ consists of the probability measures $\pi$ on $\R^d\times\R^d$ with marginals $\mu,\sigma$ and barycentric condition
\begin{equation}
  \int y\,\pi_x(\dd y)=x
  \quad\text{for $\mu$-almost every }x.
  \label{mart:g:eq:martingale-condition}
\end{equation}
The existence of such a coupling is equivalent to convex order by Strassen's theorem \cite{Str65}.

\begin{definition}[Radial fourth-moment coupling cost]
\label{mart:def:m4}
For $\mu,\sigma\in\Ptwo(\R^d)$, set
\begin{equation}\label{mart:g:eq:m4-def}
 M_4(\mu,\sigma):=\frac13
 \inf_{\pi\in\Pi_M(\mu,\sigma)}
       \int|y-x|^4\pi(\dd x,\dd y).
\end{equation}
Here $|y-x|$ is the Euclidean norm. We call the cost radial because it
uses only the length of the increment. The value is $+\infty$ when no
finite-cost martingale coupling exists.
\end{definition}
At each step the constraint fixes a conditional mean, not a local
diffusion coefficient. A coupling can assign mass to large increments;
the discrete martingale is not required to have diffusive transitions.
For chains of bounded rescaled fourth-moment action, the estimates in
\cref{mart:c:thm:compactness} control these increments and produce
continuous-path limits. This conclusion holds in every dimension.

For $\mu,\sigma\in\Ptwo(\R^d)$, $N\ge1$ and $\tau=1/N$, the
coupling discretization of the geodesic problem optimizes the intermediate marginals:
\begin{equation}
  \fB_N(\mu,\sigma):=\inf\left\{\frac1\tau\sum_{k=0}^{N-1}M_4(\mu_k,\mu_{k+1}):
  \mu_0=\mu,\ \mu_N=\sigma\right\}.
  \label{mart:c:eq:discrete-bb}
\end{equation}
\begin{remark}[Static proxy for the geodesic problem]
\label{mart:rem:static-proxy}
The endpoint cost \eqref{mart:g:eq:mb-cost} is dynamic. An exact static
representation would express it through a single endpoint coupling;
the existence of such a static formulation is an open question.
The explicit static cost $M_4$ is used as a short-time proxy,
motivated by Huesmann and Trevisan \cite[Section~6]{HT19}. They
propose a time discretization based on martingale couplings, without
imposing Gaussian transitions; its specialization to a fourth-moment cost gives
\eqref{mart:c:eq:discrete-bb}. This choice draws on their
Gaussian-averaging formula for a sufficiently regular continuous
martingale \cite[Theorem~3.1]{HT19}. In dimension one, for
$Z\sim\mathcal N(0,1)$, the identity
$\mathbb E|\sqrt{ha}\,Z|^4=3h^2a^2$ gives access to the squared
variance rate $a^2$. If volatility denotes the Brownian coefficient
$\sqrt a$, this quantity is its fourth power.

Although this Gaussian calculation motivates the fourth-moment cost,
we show that optimization over martingale couplings leads to a
different normalization from that proposed in
\cite[Section~6]{HT19}: \eqref{mart:c:eq:bb-convergence} gives
$3\fB_N(\mu,\sigma)\to\MB_1^2(\mu,\sigma)$.
Indeed, the centered two-point increment $\pm\sqrt{ha}$ has the same
variance $ha$ as the Gaussian increment but fourth moment $h^2a^2$,
one third of the Gaussian value. Sampling a fixed diffusion and
optimizing its discrete transitions therefore need not have the same
fourth-moment normalization.

In dimension one, \cref{mart:c:cor:bb-convergence} answers a question of
\cite[Section~6]{HT19} for the increment cost $|y-x|^4$.
For a partition $0=t_0<\dots<t_n=1$, they consider
\[
 \inf\Bigl\{\sum_{i=1}^n\frac1{t_i-t_{i-1}}
 \inf_{\pi\in\Pi_M(\mu_{i-1},\mu_i)}\int|y-x|^4\,\pi(\dd x,\dd y):
 \mu_0=\mu,\ \mu_n=\sigma\Bigr\},
\]
and ask whether its lower limit, as $\max_i(t_i-t_{i-1})\to0$, equals
the value of the continuous problem whose rate cost is the Gaussian
average $\mathbb E|\sqrt a\,Z|^4=3a^2$, that is,
$3\MB^2_1(\mu,\sigma)$ by \eqref{mart:eq:common-continuous-value}.
For the uniform partition with $N$ intervals, this infimum is
$3\fB_N(\mu,\sigma)$ by \eqref{mart:g:eq:m4-def} and
\eqref{mart:c:eq:discrete-bb}. Hence, by
\eqref{mart:c:eq:bb-convergence}, the lower limit is at most
$\MB^2_1(\mu,\sigma)$, and the proposed identity fails whenever
$0<\MB^2_1(\mu,\sigma)<\infty$. For example,
$\MB^2_1(\delta_0,\mathcal N(0,1))=1$: every admissible martingale
satisfies $\mathbb E\int_0^1\alpha_t^2\dd t\ge
(\mathbb E\int_0^1\alpha_t\dd t)^2=(\mathbb E X_1^2)^2=1$, with
equality for Brownian motion.

Both the unrestricted $M_4$ scheme and the midpoint scheme with
nonnegative momentum retain the prescribed endpoint constraints; their convergence
results include recovery with these exact endpoints whenever the
continuous action is finite.
\end{remark}

\subsection{Comparison of variables and computation of marginals}
\label{mart:sec:marginal-comparison}
In one dimension there is an exact projection from a coupling to an
integrated nonnegative momentum.
\begin{proposition}[Projection of a one-dimensional martingale coupling]
\label{mart:prop:coupling-momentum}
For $\pi\in\Pi_M(\mu,\sigma)$, with $\mu,\sigma\in\Ptwo(\R)$, let
$\xi=Y-X$ and define
\begin{equation}\label{mart:eq:segment-momentum}
 m_\pi(\dd z):=
 \mathbb E_\pi\left[\xi^2\int_0^1(1-s)
                     \delta_{X+s\xi}(\dd z)\dd s\right].
\end{equation}
Then $m_\pi$ is a finite nonnegative measure,
$\partial_{xx}m_\pi=\sigma-\mu$, and
$m_\pi(\R)=\mathbb E_\pi\xi^2/2$. Among finite measures on $\R$, this
momentum is uniquely determined by the endpoints. Existence of such a
nonnegative momentum is equivalent to $\mu\le_c\sigma$, hence to
existence of a martingale coupling.
\end{proposition}

The fourth moment also controls a quadratic action of the projected
momentum, with a distribution supported along the coupling segments.
Define the probability measure
\[
 \lambda_\pi:=2\mathbb E_\pi\int_0^1(1-s)
                         \delta_{X+s\xi}\,\dd s.
\]
Convexity of the square gives
\[
 \Ac(\lambda_\pi,m_\pi)\le\frac18\mathbb E_\pi|Y-X|^4.
\]
For an optimal finite-cost $M_4$ coupling, the right-hand side is
$3M_4(\mu,\sigma)/8$. This reference distribution $\lambda_\pi$
depends on the coupling and generally differs from the midpoint
measure $(\mu+\sigma)/2$. The bound therefore does not compare
$M_4(\mu,\sigma)$ with the midpoint-weighted momentum action
$\Ac((\mu+\sigma)/2,m_\pi)$. Its proof
is included with \cref{mart:prop:coupling-momentum} in
\cref{mart:app:proofs}.

\begin{remark}[Different discrete costs, the same geodesic limit]
\label{mart:rem:discrete-variables}
The $M_4$ coupling scheme and the midpoint scheme with nonnegative momentum have
different costs at a given time mesh, although their normalized minimum
actions converge to the same geodesic cost
(\cref{mart:cor:common-limit}). The reason is that the endpoint marginals
determine $m_\pi$ uniquely, while a coupling also records the conditional
distribution of the destination. Different couplings can therefore have
the same momentum but different fourth moments.

For example, the step from $\delta_0$ to
$(\delta_{-b}+\delta_b)/2$, $b>0$, has finite coupling cost $M_4=b^4/3$.
Its momentum has density $m_\pi(z)=(b-|z|)_+/2$, spread over $(-b,b)$,
whereas the midpoint measure is supported only at $0$ and $\pm b$.
Thus $m_\pi$ is not absolutely continuous with respect to the midpoint
measure, and the one-cell Eulerian action in
\eqref{mart:eq:positive-midpoint} is infinite.
\end{remark}

\begin{remark}[Choice of scheme and extension to higher dimensions]
\label{mart:rem:computational-scope}
In one dimension, momentum variables are more convenient when only
the time marginals are required. The optimization uses densities
and scalar momentum with a convex perspective cost, avoiding couplings
on the product state space. The one-dimensional periodic implementation
in Part~\ref{part:numerical} eliminates momentum up to one constant
per time cell. The $M_4$ scheme involves coupling variables but also
provides transition kernels and a discrete martingale. A continuous
martingale representation of a limiting Eulerian minimizer can instead
be obtained afterward by superposition.

The numerical approach for the momentum formulation in higher
dimensions is future work.
\end{remark}

\begin{remark}[Choice of moment]
\label{mart:rem:moment-choice}
Another legitimate cost minimizes the square of the conditional second
moment,
\[
 \mathbb E\!\left[\left(\mathbb E[|Y-X|^2\mid X]\right)^2\right],
\]
over martingale couplings. This choice does not rule out jumps.
For example, independent centered increments of size one occurring
with probability $h$ and zero otherwise, and symmetric increments
$\pm\sqrt h$, both have
conditional second moment $h$. Their fourth moments are $h$ and $h^2$,
respectively. Over $1/h$ steps, the former retain jumps, whereas the
fourth-moment penalty distinguishes the two behaviors.
\end{remark}

\section{Convergence of two different martingale geodesic schemes}\label{mart:sec:limits}
As for signed Eulerian geodesics, we first prove compactness and a
lower bound, then construct discrete paths that approximate continuous
admissible paths with the same endpoints and no increase of action in
the limit. For nonnegative momentum, this approximation uses the
shared time-averaging construction, supplemented by tightness on
$\R$. For coupling variables, fourth-moment estimates give
continuous-path compactness, and a separate construction provides the
matching upper bound. We compare the two values in
\cref{mart:cor:common-limit}; our study of their common limit is
restricted to the one-dimensional case.

\subsection{Convergence of the midpoint scheme with nonnegative momentum}
The common time-averaging construction preserves the sign of momentum.
Combining it with convex-order tightness on $\R$ gives the following
corollary of the geodesic convergence argument; the identification of
the limit is \cref{mart:prop:positive-eulerian}.

\begin{corollary}[Positive momentum and the martingale cost]
\label{mart:cor:positive-midpoint}
For every $\mu,\sigma\in\Ptwo(\R)$,
\begin{equation}\label{mart:eq:positive-value-convergence}
 \lim_{N\to\infty}\mathsf G_N^+(\mu,\sigma)
 =\mathsf G^+_{[0,1]}(\mu,\sigma)
 =\frac18\MB_1^2(\mu,\sigma).
\end{equation}
The equality includes infinite values. If the common value is finite,
the continuous infimum and all sufficiently fine discrete infima are
attained. Every sequence of asymptotically minimizing discrete pairs
has a subsequence whose affine marginal curves converge uniformly
narrowly and whose space--time momentum measures converge weakly.
Every resulting pair minimizes the continuous positive Eulerian
action and admits a continuous square-integrable martingale
representation with endpoint marginals $\mu,\sigma$.
\end{corollary}

The proof reuses \cref{geo:lem:time-average-recovery} and the
lower-bound argument of \cref{geo:thm:time-discrete-convergence}.
Its additional real-line compactness checks are given in
\cref{mart:app:positive-midpoint}; the periodic theorem alone does not
supply them.

\subsection{The fourth-moment coupling scheme}
We use the unit time interval below. The statements extend to a
finite interval by a deterministic time change.
We analyze the coupling discretization introduced in
\cref{mart:g:sec:discrete}. Continuous paths are not a finite-step
constraint; they follow from the fourth-moment estimates below.
The estimates below hold in every dimension, with constants independent of $d$.
The fourth-moment estimates used in the proofs are stated in
\cref{mart:app:moments}.
We abbreviate $\Am_1$ by $\Am$ on this interval.
For each admissible pair, choose the continuous martingale
representation of \eqref{mart:eq:martingale-eulerian-value} on a
filtered space $(\Omega,\cG,(\cG_t),\mathbb P)$, with
$\Law(X_t)=\mu_t$ and
$\langle X\rangle_t=\int_0^tA_s(X_s)\dd s$. We write
\begin{equation}
  R_t:=\tr A_t(X_t),
  \qquad
  \tr\langle X\rangle_t=\int_0^tR_s\dd s,
  \qquad
  \mathbb E\int_0^1R_t^2\dd t=\Am(\mu,A).
  \label{mart:c:eq:rate}
\end{equation}

Let $N\ge1$, $\tau:=1/N$ and $t_k:=k\tau$. A discrete martingale is an
$\R^d$-valued process $X=(X_k)_{0\le k\le N}$ adapted to a filtration
$(\cH_k)$, with $\mathbb E|X_0|^2<\infty$, $\mathbb E(X_{k+1}\mid\cH_k)=X_k$ and
$\xi_k:=X_{k+1}-X_k\in L^4$. Its discrete action is
\begin{equation}
  \fA_N(X):=\frac1\tau\sum_{k=0}^{N-1}\mathbb E|\xi_k|^4 .
  \label{mart:c:eq:discrete-action}
\end{equation}
Since $(X_k,X_{k+1})$ is a martingale coupling of its marginals,
\begin{equation}
  \frac3\tau\sum_{k=0}^{N-1}M_4\bigl(\Law(X_k),\Law(X_{k+1})\bigr)
  \le\fA_N(X),
  \label{mart:c:eq:action-cost}
\end{equation}
with equality when $X$ is the Markov gluing of optimal $M_4$ couplings. We
write $\widetilde X$ for the linear interpolation of the grid values,
$\mu_k:=\Law(X_k)$, $\bar\mu_t:=\mu_k$ for $t\in[t_k,t_{k+1})$,
$\bar\mu_1:=\mu_N$, and
\begin{equation}
  Q_k(x):=\mathbb E\bigl(\xi_k\otimes\xi_k\bigm|X_k=x\bigr),
  \qquad
  A^N_t(x):=\frac1\tau Q_k(x),
  \quad t\in[t_k,t_{k+1}).
  \label{mart:c:eq:discrete-diffusion}
\end{equation}
For $\omega\in C([0,1];\R^d)$ and $m\ge0$, set
\begin{equation}
  \Delta_m(\omega)
  :=\max_{0\le j<2^m}
  \left|\omega\Bigl(\frac{j+1}{2^m}\Bigr)-\omega\Bigl(\frac{j}{2^m}\Bigr)\right|.
  \label{mart:c:eq:dyadic-increment}
\end{equation}

\subsection{Compactness under a bound on the discrete action}
\label{mart:c:sec:compactness}

\begin{theorem}[Compactness]
\label{mart:c:thm:compactness}
Let $N_n\to\infty$, $\tau_n:=1/N_n$, and let $X^n$ be discrete martingales on
the grids $(k\tau_n)_{0\le k\le N_n}$ such that
\begin{equation}
  \sup_n\fA_{N_n}(X^n)\le\Lambda<\infty,
  \qquad
  \sup_n\mathbb E|X^n_0|^2<\infty .
  \label{mart:c:eq:uniform-action}
\end{equation}
Then the following hold.
\begin{enumerate}
\item[(i)] For every $\alpha\in(0,1/4)$,
\begin{equation}
  \sup_n\mathbb E\left[\Bigl(\sum_{m\ge0}2^{\alpha m}
  \Delta_m(\widetilde X^n)\Bigr)^4\right]<\infty,
  \label{mart:c:eq:holder}
\end{equation}
and the distributions of $\widetilde X^n$ are tight in $C([0,1];\R^d)$.
\item[(ii)] Every limit in distribution $X$ of a subsequence $(\widetilde X^{n_i})$ is a
continuous martingale, and $\bar\mu^{\,n_i}_t\to\mu_t:=\Law(X_t)$ narrowly
for every $t\in[0,1]$.
\item[(iii)] There is a Borel map $A:(0,1)\times\R^d\to\Sdp$ such that
$(\mu,A)$ is admissible and
\begin{equation}
  \Am(\mu,A)\le\liminf_{i\to\infty}\fA_{N_{n_i}}(X^{n_i}).
  \label{mart:c:eq:action-lsc}
\end{equation}
\end{enumerate}
\end{theorem}

\begin{remark}[Dyadic grids and the initial distribution]\label{mart:rem:dyadic}
For path tightness alone in \cref{mart:c:thm:compactness}(i), tightness of
the initial marginals suffices in place of their uniform second-moment bound.
The latter bound is used to pass the martingale identities and to identify
the Fokker--Planck limit in parts (ii)--(iii).
For $N_n=2^n$, $\tau_n=2^{-n}$, and Markov gluing of optimal $M_4$ couplings,
the bound $\sum_kM_4(\mu^n_k,\mu^n_{k+1})\le C\tau_n$ gives
$\sum_k\mathbb E|\xi^n_k|^4\le3C\tau_n$ and
$\int_0^1\!\int(\tr A^{N_n}_t)^2\dd\bar\mu^{\,n}_t\dd t\le3C$.
The dyadic estimates simplify to
\begin{equation}
 \mathbb E\Delta_m(\widetilde X^n)^4\le C'2^{-m},\qquad
 \Delta_m(\widetilde X^n)=2^{n-m}\Delta_n(\widetilde X^n)\quad(m\ge n).
 \label{mart:eq:dyadic-special}
\end{equation}
Every path limit is a continuous martingale with the limiting marginals.
The compactness result does not identify an explicit formula for the
limiting covariance field $A$.
\end{remark}

\subsection{Convergence of the coupling scheme for the geodesic problem}
\label{mart:c:sec:discrete-bb}

By \eqref{mart:c:eq:action-cost} and Markov gluing of optimal couplings,
\begin{equation}
  3\fB_N(\mu,\sigma)
  =\inf\bigl\{\fA_N(X):X\text{ discrete martingale},\
  \Law(X_0)=\mu,\ \Law(X_N)=\sigma\bigr\}.
  \label{mart:c:eq:discrete-bb-martingale}
\end{equation}

\begin{proposition}[Upper bound]
\label{mart:c:prop:upper-bound}
For every admissible $(\mu,A)$ with $\mu_0=\mu$, $\mu_1=\sigma$ and every
$N\ge1$,
\begin{equation}
  3\fB_N(\mu,\sigma)\le C_4\,\Am(\mu,A).
  \label{mart:c:eq:upper-bound}
\end{equation}
Here $C_4=(32/3)^2$ is the constant of \cref{mart:c:lem:continuous-bdg}.
In particular $3\fB_N(\mu,\sigma)\le C_4\MB^2_1(\mu,\sigma)$.
\end{proposition}

\begin{theorem}[Compactness and lower bound]
\label{mart:c:thm:bb-compactness}
Let $\mu,\sigma\in\Ptwo(\R^d)$ with $\MB^2_1(\mu,\sigma)<\infty$, and let $X^N$ be
discrete martingales with $\Law(X^N_0)=\mu$, $\Law(X^N_N)=\sigma$ and
$\fA_N(X^N)\le3\fB_N(\mu,\sigma)+1/N$. Then the conclusions of
\cref{mart:c:thm:compactness} hold for $(X^N)$, every limit curve satisfies
$\mu_0=\mu$ and $\mu_1=\sigma$, and, along the corresponding subsequence,
\begin{equation}
  \MB^2_1(\mu,\sigma)\le\Am(\mu,A)\le\liminf_{i\to\infty}3\fB_{N_i}(\mu,\sigma).
  \label{mart:c:eq:lower-bound}
\end{equation}
Moreover, $\MB^2_1(\mu,\sigma)\le\liminf_{N\to\infty}3\fB_N(\mu,\sigma)$ for all
$\mu,\sigma\in\Ptwo(\R^d)$.
\end{theorem}

\begin{theorem}[Sharp upper bound]
\label{mart:c:thm:sharp-upper-bound}
For every admissible $(\mu,A)$ with $\mu_0=\mu$ and $\mu_1=\sigma$,
\begin{equation}
  \limsup_{N\to\infty}3\fB_N(\mu,\sigma)\le\Am(\mu,A).
  \label{mart:c:eq:sharp-upper-bound}
\end{equation}
\end{theorem}

The proof of \cref{mart:c:thm:sharp-upper-bound} is valid in every
dimension. For a fixed $\varepsilon\in(0,1)$, it samples the martingale
$X$ of \eqref{mart:c:eq:rate} at successive exit times from balls
centered at the last sample.
The squared radius of each ball is $(1+\varepsilon)\tau$ times a
continuous adapted approximation of the rate $R$ at the last sampling
time; $R$ itself is defined only $\dd t\otimes\mathbb P$-almost
everywhere. The last $\lceil\varepsilon N/4\rceil$ steps reach time $1$,
hence the terminal marginal $\sigma$. For fixed $\varepsilon$, the
choice of the approximation of $R$ makes the limit superior as
$N\to\infty$ of their contribution to the action arbitrarily small.
In dimension one, the Dambis--Dubins--Schwarz theorem, on an enlarged
probability space, writes $X-X_0$ as a Brownian motion evaluated at the
random clock $\int_0^tR_s\dd s$. Each completed exit of $X$ from an interval is then
an exit of this Brownian motion from the same interval.

\begin{corollary}[Convergence]
\label{mart:c:cor:bb-convergence}
For all $\mu,\sigma\in\Ptwo(\R^d)$,
\begin{equation}
  \lim_{N\to\infty}3\fB_N(\mu,\sigma)=\MB^2_1(\mu,\sigma).
  \label{mart:c:eq:bb-convergence}
\end{equation}
If $\MB^2_1(\mu,\sigma)<\infty$, every limit curve in \cref{mart:c:thm:bb-compactness}
satisfies $\Am(\mu,A)=\MB^2_1(\mu,\sigma)$.
\end{corollary}

\subsection{Common limit and higher-dimensional questions}

\begin{corollary}[The common one-dimensional geodesic limit]
\label{mart:cor:common-limit}
For $\mu,\sigma\in\Ptwo(\R)$,
\begin{equation}\label{mart:eq:common-discrete-limit}
 \begin{aligned}
 \lim_{N\to\infty}8\mathsf G_N^+(\mu,\sigma)
 &=\lim_{N\to\infty}3\fB_N(\mu,\sigma)=\MB_1^2(\mu,\sigma),\\
 \lim_{N\to\infty}\mathsf G_N^+(\mu,\sigma)
 &=\lim_{N\to\infty}\frac38\fB_N(\mu,\sigma)
   =\mathsf G^+_{[0,1]}(\mu,\sigma).
 \end{aligned}
\end{equation}
If the common value is finite, the positive Eulerian infimum is attained.
For asymptotically minimizing $M_4$ chains, every path limit supplied by
\cref{mart:c:thm:bb-compactness} is a continuous martingale, and the
limiting Eulerian pair $(\mu_t,m_t)$, with $m_t=a_t\mu_t/2$, minimizes
\eqref{mart:eq:positive-eulerian}.
\end{corollary}

The factors in \eqref{mart:eq:common-discrete-limit} correspond to the
paper's Eulerian normalization $\int m^2/(2\rho)$ and the definition
$M_4=(1/3)\inf\mathbb E|Y-X|^4$. For the action $\int m^2/\rho$,
multiply the second line by two.

Both time discretizations therefore approximate the same geodesic
problem in dimension one. At finite $N$ the couplings and momenta
are different variables with different objectives, as explained in
\cref{mart:prop:coupling-momentum,mart:rem:discrete-variables}.
By \cref{mart:prop:positive-eulerian}, the continuous positive Eulerian
problem and the continuous-martingale problem have values related by
$8\mathsf G^+_{[0,1]}=\MB^2_1$. The substitution
\eqref{mart:eq:momentum-covariance} maps the minimizing pairs of
\eqref{mart:eq:positive-eulerian} to those of
\eqref{mart:eq:martingale-eulerian-value}.

\begin{remark}[Higher-dimensional choices of Eulerian action]
\label{mart:rem:higher-bb}
For $\mathsf m=\frac12A\mu\succeq0$, the Frobenius-square and
trace-square momentum actions are, respectively,
\begin{equation}\label{mart:eq:positive-tensor-action}
 \frac12\int\left|\frac{\dd\mathsf m_t}{\dd\mu_t}\right|_F^2
                    \dd\mu_t\dd t
 =\frac18\int|A_t|_F^2\dd\mu_t\dd t,
\end{equation}
\begin{equation}\label{mart:eq:trace-momentum-action}
 \frac12\int\left(\tr\frac{\dd\mathsf m_t}{\dd\mu_t}\right)^2
                    \dd\mu_t\dd t
 =\frac18\int(\tr A_t)^2\dd\mu_t\dd t.
\end{equation}
The trace action matches the radial $M_4$ limit
(\cref{mart:def:m4}); the Frobenius action extends the signed tensor
model. Both permit anisotropic covariance. Scalar momentum with
$\partial_t\mu=\Delta m$ instead imposes
$A=2(\dd m/\dd\mu)I$, hence isotropic diffusion. These choices agree
in one dimension but differ in higher dimensions.

Time averaging preserves positive-semidefinite momentum, but positive
Eulerian convergence is proved here only in one dimension.
Higher-dimensional convergence and efficient tensor optimization remain
open; see \cref{mart:rem:computational-scope}.
\end{remark}

\begin{openquestion}[A coupling approximation of the positive tensor action]
\label{mart:open:tensor-bb}
Does the cost
\begin{equation}\label{mart:eq:tensor-candidate}
 \mathsf C_h(\mu,\sigma):=\inf_{\pi\in\Pi_M(\mu,\sigma)}
 \int\left\{|Q_\pi(x)|_F^2
       +\varepsilon_h\int|y-x|^4\pi_x(\dd y)\right\}\mu(\dd x),
\end{equation}
with $Q_\pi(x):=\int(y-x)\otimes(y-x)\pi_x(\dd y)$,
$\varepsilon_h\to0$ and $h/\varepsilon_h\to0$, yield geodesic
convergence to \eqref{mart:eq:positive-tensor-action} after
normalization by $1/8$? Establish compactness, the lower bound and
recovery preserving the endpoint marginals, with continuous martingale
limits.
\end{openquestion}

\section{Convergence of martingale minimizing movements}
\label{mart:g:sec:m4-theorem}
We are now interested in proving convergence of the gradient-flow
discretizations in the martingale setting. We explore two
discretizations: the coupling scheme based on $M_4$ and the modified
JKO scheme with nonnegative momentum. Both optimize the next
marginal by balancing the step cost against an energy, as in the signed
modified JKO scheme of \cref{sec:merged-jko}. Geodesic convergence
alone does not identify the limiting flows; each scheme requires an
additional optimality argument.

For general nonlinear energies, the constrained momentum result in
\cref{jko:prop:positive-compactness} supplies existence and weak
compactness on the torus. We identify the limiting evolution equation
below for one-dimensional potential energies on $\R$.
The $M_4$ theorem assumes uniform concavity; the constrained momentum
theorem allows curvature of either sign when the second and third
derivatives are bounded. Under their common assumptions, the two
limiting curves coincide up to a constant rescaling of time determined
by the kinetic coefficients.

\subsection{Formal gradient flow and minimizing movements}
If $U_\mu=\delta\Em/\delta\mu$ is smooth, the formal first variation is
\begin{equation}
 \frac{\dd}{\dd t}\Em(\mu_t)
 =\frac12\int D^2U_{\mu_t}:A_t\dd\mu_t.
 \label{mart:g:eq:intro-first-variation}
\end{equation}
In one dimension, minimizing the corresponding quadratic dissipation gives
\begin{equation}
 a_\mu=\frac12(-U_\mu'')_+,
 \qquad \partial_t\mu=\frac14\partial_{xx}\bigl(\mu(-U_\mu'')_+\bigr).
 \label{mart:g:eq:intro-formal-gf}
\end{equation}
These formulas motivate the following tensor calculation.

\begin{remark}[Higher-dimensional trace dissipation]
\label{mart:prop:tensor}
Let $H\in\Sd$ and let $\lambda_{\min}(H)$ be its smallest eigenvalue. The minimizers of
\begin{equation}
  A\longmapsto\frac12(\tr A)^2+\frac12H:A,
  \qquad A\in\Sdp,
  \label{mart:g:eq:tensor-min}
\end{equation}
are characterized as follows. If $\lambda_{\min}(H)\ge0$, then $A=0$. If $\lambda_{\min}(H)<0$, then
\begin{equation}
  \tr A=-\frac12\lambda_{\min}(H),
  \qquad
  \operatorname{Ran}(A)\subseteq E_{\min}(H),
  \label{mart:g:eq:tensor-min-characterization}
\end{equation}
where $E_{\min}(H)$ is the eigenspace associated with $\lambda_{\min}(H)$.
\end{remark}

In particular, if $U_\mu=\delta\Em/\delta\mu$ is smooth, the formal gradient flow is
\begin{equation}
  \partial_t\mu=\frac12\partial_{ij}(A_{\mu,ij}\mu),
  \qquad
  A_\mu(x)\in\argmin_{A\ge0}
  \left\{\frac12(\tr A)^2+\frac12D^2U_\mu(x):A\right\}.
  \label{mart:g:eq:formal-gradient-flow-d}
\end{equation}
In one dimension the minimizer is unique and \eqref{mart:g:eq:formal-gradient-flow-d} reduces to \eqref{mart:g:eq:intro-formal-gf}.

\begin{remark}
The formula is a formal calculation. For a nonlinear energy, one must justify the first variation, compactness of the discrete subdifferentials, and the passage to the positive-part operation. For one-dimensional potential energies, \cref{mart:thm:oned,mart:thm:positive-potential} prove convergence to \eqref{mart:g:eq:intro-formal-gf}: the first for uniformly concave potentials and the cost $3M_4$, the second for the constrained momentum scheme \eqref{mart:eq:positive-potential-step} with $\kappa=4$ (\cref{mart:cor:normalization,mart:cor:positive-potential-normalization}). The tensor formula \eqref{mart:g:eq:formal-gradient-flow-d} is not justified for $d\ge2$.
\end{remark}

\subsubsection{Discrete minimizing movements}
The variational schemes below follow \cite{JKO98,AGS}. The two
approximations compared here are the fourth-moment coupling scheme and
the modified JKO scheme with nonnegative momentum.
For a time step $\tau>0$ and an energy $\Em$, the coupling step on
$\R^d$ is
\begin{equation}
  \mu^{M}_{k+1}
  \in\argmin_{\mu^M_k\le_c\mu}
  \left\{\frac1{2\tau}M_4(\mu^M_k,\mu)+\Em(\mu)\right\}.
  \label{mart:g:eq:m4-jko}
\end{equation}
For the potential energy $\Em(\mu)=\int_\R\Vm\dd\mu$ of a potential
$\Vm:\R\to\R$, defined in \eqref{mart:g:eq:linear-energy} below,
the constrained momentum step starts from $\mu_0\in\Ptwo(\R)$ and, for $\kappa>0$,
chooses
\begin{equation}\label{mart:eq:positive-potential-step}
 (\mu_{k+1}^\tau,m_{k+1}^\tau)
 \in\operatorname*{argmin}_{\substack{\nu\in\Ptwo(\R),\ m\ge0\\
                                     \partial_{xx}m=\nu-\mu_k^\tau}}
 \left\{\frac\kappa\tau\Ac(\nu,m)+\int_\R\Vm\dd\nu\right\}.
\end{equation}
The momentum $m$ is a finite nonnegative measure integrated over the
step and the action is
weighted by the new marginal $\nu$, just as in
\eqref{jko:eq:positive-scheme}. The choice $\kappa=1$ retains the
signed modified JKO normalization. If $m=\tau b\nu$, then $b$ is a
discrete momentum rate. In the continuous Fokker--Planck normalization,
a momentum rate $b$ corresponds to covariance rate $a=2b$; this
accounts for the different kinetic coefficients when comparing the
limits, without identifying finite-step momenta with coupling variances.
Let $\bar\mu_\tau(t)=\mu_{k+1}^\tau$ on $(k\tau,(k+1)\tau]$,
with $\bar\mu_\tau(0)=\mu_0$.

For a general energy, the coupling estimates below assume that the
chosen $M_4$ steps exist. Under the respective one-dimensional
potential hypotheses, existence follows from
\cref{mart:prop:pointwise,mart:thm:positive-potential} for the
fourth-moment coupling and constrained momentum schemes, respectively.
Whenever $\Em$ is bounded below along the fourth-moment scheme,
comparison with the unchanged measure gives
\begin{equation}\label{mart:g:eq:m4-energy}
 \sum_kM_4(\mu^M_k,\mu^M_{k+1})\le C\tau.
\end{equation}

\subsection{Compactness of the fourth-moment minimizing movement}
\label{mart:c:sec:gradient-flow}

Let $\Em:\Ptwo(\R^d)\to(-\infty,+\infty]$ and $\mu_{\mathrm{in}}\in\Ptwo(\R^d)$
with $\Em(\mu_{\mathrm{in}})<\infty$. For $N\ge1$ and $\tau=1/N$, start from
$\mu^N_0:=\mu_{\mathrm{in}}$ and use the recursion
\eqref{mart:g:eq:m4-jko} for $0\le k<N$.
Let $X^N$ be the Markov gluing of optimal $M_4$ couplings of
$(\mu^N_k)_{0\le k\le N}$.

\begin{corollary}[Minimizing movement]
\label{mart:c:cor:gradient-flow}
Assume that every step of \eqref{mart:g:eq:m4-jko} has a minimizer and that
\begin{equation}
  \Phi:=\sup_{N\ge1}\max_{0\le k\le N}
  \bigl(\Em(\mu_{\mathrm{in}})-\Em(\mu^N_k)\bigr)<\infty .
  \label{mart:c:eq:bounded-below}
\end{equation}
Then $\fA_N(X^N)\le6\Phi$ for every $N$. Consequently the conclusions of
\cref{mart:c:thm:compactness} hold for $(X^N)$, and every limit curve satisfies
$\mu_0=\mu_{\mathrm{in}}$ and $\Am(\mu,A)\le6\Phi$.
\end{corollary}

\subsection{Identification of the fourth-moment potential flow}

We first identify the $M_4$ scheme in dimension one for the potential
energy
\begin{equation}
  \Em(\mu)=\int_\R \Vm(x)\mu(\dd x).
  \label{mart:g:eq:linear-energy}
\end{equation}
The functional is linear in the measure; the potential $\Vm$ itself
need not be affine.
It is useful to keep a general kinetic normalization. Let $\eta>0$ and consider, conditionally on a starting point $x$,
\begin{equation}
  \mathcal J_{\tau,x}(K)
  :=\int\left(\frac\eta\tau|u|^4+\Vm(x+u)\right)K(\dd u),
  \qquad \int uK(\dd u)=0.
  \label{mart:g:eq:cell-problem}
\end{equation}
The scheme \eqref{mart:g:eq:m4-jko} has $\eta=1/6$, while replacing $M_4$ by $3M_4$ gives $\eta=1/2$.

\begin{proposition}[Pointwise reduction and extreme kernels]
\label{mart:prop:pointwise}
Suppose that $\Vm:\R\to\R$ is lower semicontinuous, $\Vm(z)\ge-C(1+|z|^2)$, and $\mu_k\in\Ptwo(\R)$.
Write $\Pcal_4(\R)$ for probability measures with finite fourth moment. The one-step problem associated with \eqref{mart:g:eq:cell-problem} separates in $x$:
\begin{equation}
  \inf_{\substack{\pi_1=\mu_k\\ \mathbb E_\pi[Y\mid X]=X}}
  \int\left(\frac\eta\tau|y-x|^4+\Vm(y)\right)\pi(\dd x,\dd y)
  =\int\inf_{\substack{K\in\Pcal_4(\R)\\ \int u\dd K=0}}
  \mathcal J_{\tau,x}(K)\,\mu_k(\dd x).
  \label{mart:g:eq:pointwise-separation}
\end{equation}
For each $x$, the value on the right is the convex envelope at zero of
$u\mapsto\eta u^4/\tau+\Vm(x+u)$. A minimizer can be chosen with at most two atoms; it is an extreme point of the set of centered probability measures. The infimum on the left is attained.
\end{proposition}

The case of a nonlinear energy is left open for future work.

\begin{theorem}[One-dimensional $M_4$ limit]
\label{mart:thm:oned}
Assume
\begin{equation}
  \Vm\in C^3(\R),
  \qquad
  -L\le \Vm''\le-L_0<0,
  \qquad
  \|\Vm'''\|_\infty<\infty.
  \label{mart:g:eq:F-assumptions}
\end{equation}
For each $x$, let $K_{\tau,x}$ be a minimizer of \eqref{mart:g:eq:cell-problem}, and set
\begin{equation}
  v_\tau(x):=\int u^2K_{\tau,x}(\dd u),
  \qquad
  w_\tau(x):=\int u^4K_{\tau,x}(\dd u).
\end{equation}
Then, uniformly in $x$ and in the choice of the minimizers,
\begin{equation}
  \frac{v_\tau(x)}\tau\longrightarrow-\frac{\Vm''(x)}{4\eta},
  \qquad
  w_\tau(x)=O(\tau^2).
  \label{mart:g:eq:variance-asymptotic}
\end{equation}
Let $\mu_0\in\Ptwo(\R)$. The Markov chain starts from $X_0\sim\mu_0$, and at step $k$ it moves from $X_k=x$ to $x+u$ with $u\sim K_{\tau,x}$; the family $x\mapsto K_{\tau,x}$ may change with $k$ and is measurable for the completion of the law of $X_k$. The linearly interpolated chains converge in distribution on $C([0,\tfin];\R)$ to the unique weak solution of
\begin{equation}
  \dd X_t=\sqrt{-\frac{\Vm''(X_t)}{4\eta}}\,\dd B_t,
  \qquad X_0\sim\mu_0.
  \label{mart:g:eq:limit-sde-beta}
\end{equation}
Consequently, for every $t\in[0,\tfin]$, their marginals converge narrowly to $\mu_t$, where $\mu$ is the unique narrowly continuous probability solution of
\begin{equation}
  \partial_t\mu
  =-\frac1{8\eta}\partial_{xx}(\Vm''\mu),
  \qquad \mu|_{t=0}=\mu_0.
  \label{mart:g:eq:limit-pde-beta}
\end{equation}
\end{theorem}

\begin{corollary}[Normalization of the limiting flow]
\label{mart:cor:normalization}
For the definition \eqref{mart:g:eq:m4-def}, one has $\eta=1/6$ and
\begin{equation}
  \partial_t\mu=-\frac34\partial_{xx}(\Vm''\mu).
  \label{mart:g:eq:m4-original-pde}
\end{equation}
For the cost $3M_4$, one has $\eta=1/2$ and
\begin{equation}
  \partial_t\mu=-\frac14\partial_{xx}(\Vm''\mu),
  \label{mart:g:eq:m4-rescaled-pde}
\end{equation}
which agrees with the formal gradient flow \eqref{mart:g:eq:intro-formal-gf} in the uniformly concave regime.
\end{corollary}

The effect of the kinetic normalization is already visible in the following exact example.

\begin{example}[Quadratic potential]\label{mart:ex:quadratic}
Let $\Vm(x)=-ax^2/2$, $a>0$. The unique optimizer is
\begin{equation}
  K_{\tau,x}
  =\frac12\delta_{-\sqrt{a\tau/(4\eta)}}
   +\frac12\delta_{\sqrt{a\tau/(4\eta)}}.
  \label{mart:g:eq:quadratic-bernoulli}
\end{equation}
The theorem is therefore exact in this case: the discrete process is a symmetric random walk converging to Brownian motion with variance rate $a/(4\eta)$.
\end{example}

\subsection{Identification of the constrained momentum potential flow}
\label{mart:sec:positive-potential}
We now treat the discretization with the constrained momentum
formulation \eqref{mart:eq:positive-potential-step}, for the same
one-dimensional potential energy. Completing the square and comparing
with a diffusion resolvent identify the limit without a density floor.
The conservation and lower-semicontinuity arguments remain those of
\cref{jko:prop:positive-compactness}. Since a pure potential energy
has no superlinear internal term and the domain is $\R$, existence and
compactness require separate moment and tail estimates.

\begin{theorem}[Constrained momentum limit in one dimension]
\label{mart:thm:positive-potential}
Assume
\begin{equation}\label{mart:eq:positive-potential-assumptions}
 \Vm\in C^3(\R),\qquad
 \|\Vm''\|_\infty+\|\Vm'''\|_\infty<\infty,
 \qquad D(x):=\frac{(-\Vm''(x))_+}{\kappa}.
\end{equation}
Every step of \eqref{mart:eq:positive-potential-step} is attained.
For every $T<\infty$, any choice of minimizing pairs satisfies
\begin{equation}\label{mart:eq:positive-potential-convergence}
 \sup_{0\le t\le T}W_2(\bar\mu_\tau(t),\mu_t)\longrightarrow0,
\end{equation}
where $\mu$ is the unique narrowly continuous probability solution of
\begin{equation}\label{mart:eq:positive-potential-pde}
 \partial_t\mu=\partial_{xx}(D\mu),\qquad \mu|_{t=0}=\mu_0.
\end{equation}
Equivalently, $\mu_t=\Law(X_t)$, where $X$ is the unique solution in distribution of
\begin{equation}\label{mart:eq:positive-potential-sde}
 \dd X_t=\sqrt{2D(X_t)}\dd B_t,\qquad X_0\sim\mu_0.
\end{equation}
If $m_{k+1}^\tau=\tau b_{k+1}^\tau\mu_{k+1}^\tau$, set
$b_\tau(t)=b_{k+1}^\tau$ on each time cell. Then, for a constant depending only on $\Vm,\kappa,T$,
\begin{equation}\label{mart:eq:positive-potential-flux-error}
 \int_0^T\!\!\int
       |b_\tau-D|\dd\bar\mu_\tau\dd t\le C_T\tau^{1/4}.
\end{equation}
In particular the space--time momentum rates converge weakly to
$\dd t\otimes D\mu_t$. The limit satisfies
\begin{equation}\label{mart:eq:positive-potential-energy}
 \int\Vm\dd\mu_t+
 \kappa\int_0^t\!\!\int D(x)^2\,\mu_s(\dd x)\dd s
 =\int\Vm\dd\mu_0.
\end{equation}
\end{theorem}

No density floor or absolute continuity of $\mu_0$ is assumed. The
convergence is a statement about time marginals; no transition
kernel is chosen in \eqref{mart:eq:positive-potential-step}.
The proof is in \cref{mart:app:positive-potential}. It shows that the
limiting momentum per unit time is $D\mu_t$, even where $D$ vanishes,
and hence gives the diffusion equation \eqref{mart:eq:positive-potential-pde}.

\begin{corollary}[Matching the two potential flows]
\label{mart:cor:positive-potential-normalization}
Under the uniform-concavity assumptions
\eqref{mart:g:eq:F-assumptions}, the constrained momentum scheme
\eqref{mart:eq:positive-potential-step} with $\kappa=8\eta$ and the
fourth-moment scheme \eqref{mart:g:eq:cell-problem} have the same
limiting marginal curve. In particular $\kappa=4$ matches $3M_4$,
and $\kappa=4/3$ matches the unscaled $M_4$ scheme.
For a convex $\Vm$ satisfying
\eqref{mart:eq:positive-potential-assumptions}, the constrained momentum
scheme is stationary.
\end{corollary}

For $\Vm(x)=-ax^2/2$, $a>0$, the constrained momentum step is exactly
the diffusion resolvent
$\mu_{k+1}^\tau-\tau(a/\kappa)\partial_{xx}\mu_{k+1}^\tau=\mu_k^\tau$.
It convolves $\mu_k^\tau$ with the density
\[
 z\longmapsto\frac1{2\sqrt{a\tau/\kappa}}
                  e^{-|z|/\sqrt{a\tau/\kappa}}.
\]
The fourth-moment step in \cref{mart:ex:quadratic} uses a two-point
increment instead. The two schemes have the same diffusion
limit after matching $\kappa=8\eta$.

\subsection{Higher-dimensional remarks and future work}

\begin{remark}[Diffusion in the most concave directions]
\label{mart:rem:high-d-cell}

The one-dimensional theorem does not extend to the equation $\partial_t\mu=-\frac1{8\eta}\partial_{ij}(\partial_{ij}\Vm\,\mu)$, obtained by replacing $\Vm''$ with $D^2\Vm$ in \eqref{mart:g:eq:limit-pde-beta}. For a time step $\tau$, the second-order Taylor expansion of $\Vm$ at $x$ turns the $d$-dimensional analogue of \eqref{mart:g:eq:cell-problem}, minus $\Vm(x)$, into
\begin{equation}
  \inf_{\mathbb EV=0}
  \left\{\frac\eta\tau\mathbb E|V|^4+\frac12\mathbb E[V^TD^2\Vm(x)V]\right\}.
  \label{mart:g:eq:high-d-cell}
\end{equation}
For a quadratic $\Vm$, this expansion is exact.
Let $\mathcal B=-D^2\Vm(x)$ and let $\lambda_{\max}(\mathcal B)>0$. If $r=\mathbb E|V|^2$, then
\begin{equation}
  \frac\eta\tau\mathbb E|V|^4-\frac12\mathbb E[V^T\mathcal BV]
  \ge\frac\eta\tau r^2-\frac12\lambda_{\max}(\mathcal B)r.
  \label{mart:g:eq:high-d-cell-bound}
\end{equation}
Equality is attained by a symmetric two-point distribution in a top eigendirection, with
$r=\tau\lambda_{\max}(\mathcal B)/(4\eta)$. More generally, every minimizing distribution has constant radius and is supported in the top eigenspace. Its covariance is therefore rank one when that eigenspace is one-dimensional. If the top eigenvalue has higher multiplicity, the covariance is not uniquely selected and may have higher rank. In particular, an isotropic quadratic potential permits an isotropic minimizing covariance but does not canonically select it.
For a quadratic potential, these minimizers are exactly the optimal increment distributions of the one-step problem, because the separation argument in the proof of \cref{mart:prop:pointwise} does not use the dimension. For a general potential, they solve only the second-order problem \eqref{mart:g:eq:high-d-cell}, and the higher-dimensional limit is not proved.
For a strongly concave anisotropic quadratic potential with a simple largest eigenvalue of $-D^2\Vm$, the trace action selects diffusion only in that eigendirection, despite concavity in every coordinate. Formally, the positive Frobenius action would instead produce the full tensor $(-D^2\Vm)_+$, up to the kinetic coefficient.
\end{remark}

\paragraph{Future work for the $M_4$ scheme.}
Identification of the higher-dimensional minimizing-movement limit
and convergence for nonlinear energies are left open for future work.

\paragraph{Future work for the constrained momentum scheme.}
For nonlinear energies on $\Torus$ satisfying \emph{(E1)--(E2)}, existence
and weak compactness follow from
\cref{jko:prop:positive-compactness}, but the limiting evolution
equation remains to be identified. Extending the one-dimensional
potential-flow result to higher dimensions is also left for future work.

\mainpart{Numerical methods and investigations}\label{part:numerical}
We now discretize the signed variational problems of Part~\ref{part:convergence} on periodic grids.
This part develops the numerical solvers and bounds their optimization error at fixed spatial grid and time step.

\section{Two discrete convex problems}\label{sec:problems}

We consider two convex problems for densities on a periodic grid.
The geodesic problem connects two prescribed densities by a path of least action.
The entropic problem computes one time step from a prescribed density $\nu^{\rm old}$: it minimizes the entropy of the density after the step plus $\tau^{-1}$ times the action of the transition from $\nu^{\rm old}$.
Both problems use a momentum variable, linked to the change in density by a discrete conservation equation.

\subsection{Grid and notation}\label{sec:grid-notation}
Let $M\ge3$ be the number of spatial grid points.
We use the periodic grid $x_i=i/M$, $0\le i<M$, on $\T=\R/\mathbb Z$.
Set $\dx=1/M$ and $\one=(1,\ldots,1)\in\R^M$.
For $u,v\in\R^M$, write
\[
 \ip uv=\dx\sum_{i=0}^{M-1}u_iv_i,\qquad
 \avg u=\dx\sum_{i=0}^{M-1}u_i,\qquad
 \normone u=\dx\sum_{i=0}^{M-1}|u_i|.
\]

For arrays, $\|\cdot\|_2$ is the unweighted Euclidean norm and $\|\cdot\|_\infty$ is the largest absolute entry.
For matrices, $\|\cdot\|_2$ is the induced Euclidean operator norm.
The inequalities $u\ge0$ and $u>0$ apply to each entry.
Products, quotients, powers, and scalar functions of vectors also apply to each entry.
We use the probability simplex and its positive part:
\[
 \mathcal P_M=\{\nu\in\R^M:\nu\ge0,\ \avg\nu=1\},\qquad
 \mathcal P_M^+=\{\nu\in\mathcal P_M:\nu>0\}.
\]
Thus $\nu_i$ is a density value and $\dx\nu_i$ is the probability mass at node $i$.

Let $\mathcal I_M=\{-\lfloor M/2\rfloor,\ldots,\lceil M/2\rceil-1\}$ be the Fourier index set.
Our Fourier convention is
\[
 \widehat u_k=\dx\sum_{i=0}^{M-1}u_i e^{-2\pi\mathrm i ki/M},\qquad
 u_i=\sum_{k\in\mathcal I_M}\widehat u_k e^{2\pi\mathrm i ki/M}.
\]

The operator $L:\R^M\to\R^M$ has Fourier symbol $\widehat L(k)=-4\pi^2k^2$.
Thus $\widehat{Lu}_k=\widehat L(k)\widehat u_k$.
The operator is self-adjoint for $\ip\cdot\cdot$.
Its kernel is $\R\one$, and its range is $\R^M_0=\{u:\avg u=0\}$.

Define $\C$ by the symbol $1/\widehat L(k)$ for $k\ne0$, and by zero for $k=0$.
Then
\[
 L\C u=\C Lu=u-\avg u\one.
\]

\subsection{Action and discrete geodesics}
For a density weight $q$ and a momentum $m$, define the action density and its grid integral:
\begin{equation}\label{eq:perspective}
 a(q,m)=
 \begin{cases}
 m^2/(2q),&q>0,\\
 0,&q=m=0,\\
 +\infty,&\text{otherwise}.
 \end{cases}
 \qquad
 \A(q,m)=\dx\sum_{i=0}^{M-1}a(q_i,m_i).
\end{equation}
The function $a$ is lower semicontinuous and convex.
Finite action requires $m_i=0$ at each entry where $q_i=0$.
We call a zero density entry a vacuum entry.

Fix $N\ge2$, set $h=1/N$, and prescribe $\nu^0,\nu^N\in\mathcal P_M^+$.
The free densities are $\nu^1,\ldots,\nu^{N-1}$.
Each time cell $j$ has an integrated momentum $m^j\in\R^M$, for $0\le j<N$.
Set $q^j=(\nu^j+\nu^{j+1})/2$.
We call $q^j$ the midpoint weight.

The density $\nu^j$ approximates the path at time $t_j=jh\in[0,1]$.
The momentum $m^j$ approximates $\int_{t_j}^{t_{j+1}}m_t\,\dd t$, where $m_t$ is the continuous momentum at time $t$, and the momentum rate of cell $j$ is $w^j=m^j/h$.
Integrating $\partial_t\rho=\partial_{xx}m_t$ over $[t_j,t_{j+1}]$ gives $\rho(t_{j+1})-\rho(t_j)=\partial_{xx}\int_{t_j}^{t_{j+1}}m_t\,\dd t$; the conservation equation $Lm^j=\nu^{j+1}-\nu^j$ in \eqref{eq:geodesic} is its grid version.
Since $\A(q,\cdot)$ is quadratic, the action of cell $j$ is $h\,\A(q^j,w^j)=h^{-1}\A(q^j,m^j)$.

The discrete geodesic problem is
\begin{equation}\label{eq:geodesic}
 \begin{aligned}
 E_{M,N}=\min_{\nu^1,\ldots,\nu^{N-1},\,m^0,\ldots,m^{N-1}}\quad
 &\sum_{j=0}^{N-1}h^{-1}\A(q^j,m^j),\\
 \text{subject to}\quad
 &Lm^j=\nu^{j+1}-\nu^j\quad(0\le j<N),\\
 &\nu^j\ge0\quad(1\le j<N).
 \end{aligned}
\end{equation}

We call the objective of \eqref{eq:geodesic} the midpoint action.
Since $\avg{Lu}=0$ for every $u$, the conservation equation gives $\avg{\nu^{j+1}}=\avg{\nu^j}$, so each feasible density and each midpoint weight has unit mass, $\avg{\nu^j}=\avg{q^j}=1$.
Because $\ker L=\R\one$, the solutions of the conservation equation are $m^j=\C(\nu^{j+1}-\nu^j)+c_\star^j\one$ with $c_\star^j\in\R$; for a fixed density path, minimization of the action fixes $c_\star^j$.

The signed diffusivity is
\begin{equation}\label{eq:kappa}
 \kappa_i^j=\frac{m_i^j}{h q_i^j}\qquad(q_i^j>0).
\end{equation}
We do not assign a value to this quotient at vacuum.
Since $a(q,-m)=a(q,m)$ and no sign constraint is imposed on $m$, the diffusivity $\kappa$ in \eqref{eq:kappa} can take either sign.
Its continuous counterpart occurs in $\partial_t\rho=\partial_{xx}(\rho\kappa)$.
For constant $\kappa$, $\partial_t\rho=\kappa\,\partial_{xx}\rho$ is the forward heat equation if $\kappa>0$ and the backward heat equation if $\kappa<0$.
For variable $\kappa$, $\partial_{xx}(\rho\kappa)=\kappa\,\partial_{xx}\rho+2\,\partial_x\kappa\,\partial_x\rho+\rho\,\partial_{xx}\kappa$; the sign of $\kappa(x)$ gives the sign of the coefficient of $\partial_{xx}\rho$ at $x$, not the sign of $\partial_t\rho(x)$.

\begin{proposition}[Existence of a discrete geodesic]\label{prop:geodesic-existence}
For $M\ge3$, $N\ge2$, and endpoints in $\mathcal P_M^+$, problem \eqref{eq:geodesic} has a minimizer of finite action.
For a specified feasible density path of finite action, the momentum of least action is unique.
\end{proposition}
The proof in \cref{app:existence} uses a positive feasible path, compactness of $\mathcal P_M$, and the bound $\normone m\le\sqrt{2\A(q,m)}$ for $q\in\mathcal P_M$ with finite action, which follows from the Cauchy--Schwarz inequality.

\begin{remark}[Uniqueness of the density path]
The proposition establishes existence of a minimizing density path.
Uniqueness is not proved here.
The action is not strictly convex, but this alone does not imply nonuniqueness, because admissible paths must also satisfy the conservation equations and prescribed endpoints.
For a minimizer with positive densities and midpoint weights, uniqueness of the minimizing density path is equivalent to the condition \eqref{eq:constrained-hessian} of the local DR analysis (\cref{app:lfp}).
\end{remark}

\subsection{Entropic time steps}\label{sec:entropic-step}
Define the closed entropy density by
\[
 e(r)=r\log r-r+1\quad(r>0),\qquad e(0)=1,\qquad e(r)=+\infty\quad(r<0).
\]
We use natural logarithms and the convention $0\log0=0$.
Set $\E(\nu)=\dx\sum_i e(\nu_i)$.

Given $\nu^{\rm old}\in\mathcal P_M^+$ and $\tau>0$, calculate
\begin{equation}\label{eq:step}
 \begin{aligned}
 H_\tau^*=\min_{\nu,m\in\R^M}\quad
 &H_\tau(\nu,m):=\tau^{-1}\A(\nu,m)+\E(\nu),\\
 \text{subject to}\quad &Lm=\nu-\nu^{\rm old}.
 \end{aligned}
\end{equation}
In \eqref{eq:step}, the action is weighted by the unknown density $\nu$, the density at the end of the step, and is evaluated on a single time cell. We call the minimizing density $\nu^*$ (\cref{prop:positive}) the successor density.
A JKO step for the signed diffusive transport distance $\mathsf D$ of \eqref{geo:eq:distance} would minimize $\E(\nu)+\mathsf D(\nu^{\rm old},\nu)^2/(2\tau)=\E(\nu)+\tau^{-1}\mathsf G_{\Lc}(\nu^{\rm old},\nu)$, where $\mathsf G_{\Lc}$, defined in \eqref{geo:eq:continuum-value}, is the minimum action over all paths from $\nu^{\rm old}$ to $\nu$.
Problem \eqref{eq:step} replaces $\mathsf G_{\Lc}(\nu^{\rm old},\nu)$ by the single-cell action $\A(\nu,m)$ with $Lm=\nu-\nu^{\rm old}$.
We call \eqref{eq:step} the entropic step, and its iteration in time the modified JKO scheme.

For smooth positive densities and with $L$ replaced by $\partial_{xx}$, one step of \eqref{eq:step} is consistent with the Derrida--Lebowitz--Speer--Spohn (DLSS) equation \cite{DLSS91,MRSSdlss} (\cref{prop:step-consistency}).
\Cref{jko:cor:dlss} identifies this equation as the limit $\tau\to0$ of the scheme \eqref{jko:eq:scheme}, the modified JKO scheme posed in continuous space, under the positivity floor \eqref{jko:eq:floor}:
\begin{equation}\label{eq:dlss}
 \partial_t\rho=-\partial_{xx}\bigl(\rho\,\partial_{xx}\log\rho\bigr).
\end{equation}
\begin{remark}[Related variational formulations]
Related Wasserstein schemes include the variational finite-volume method of Canc\`es, Gallou\"et, and Todeschi \cite{CancesGallouetTodeschi20} and the geodesic-extrapolation time discretization of Gallou\"et, Natale, and Todeschi \cite{GallouetNataleTodeschi24}.
Equation~\eqref{eq:dlss} also has a Wasserstein gradient-flow formulation based on Fisher information \cite{GST09}.
On bounded convex domains with Neumann boundary conditions, Matthes, Savar\'e, and Schlichting \cite{MSS26} obtain the solutions of the DLSS equation from a semigroup that is a contraction for the Hellinger distance, the $L^2$ distance between the square roots of two densities.
Their implicit Euler scheme, in the variable $\sqrt\rho$, solves a variational inequality at each step, whereas the entropic step \eqref{eq:step} is a convex minimization problem.
\end{remark}

\begin{proposition}[Existence, uniqueness, and positivity of an entropic step]\label{prop:positive}
For $M\ge3$, $\tau>0$, and $\nu^{\rm old}\in\mathcal P_M^+$, problem \eqref{eq:step} has a unique minimizer $(\nu^*,m^*)$.
Its density satisfies $\nu^*\in\mathcal P_M^+$.
Moreover,
\[
 \E(\nu^*)+\tau^{-1}\A(\nu^*,m^*)\le\E(\nu^{\rm old}).
\]
\end{proposition}
The proof in \cref{app:entropic-existence} uses the singular derivative $e'(r)=\log r$ at zero to exclude a vacuum entry by a feasible perturbation.
Strict convexity of entropy gives density uniqueness.
Positivity holds for each fixed problem; the proposition gives no lower bound on $\min_i\nu_i^*$ uniform in $M$, $\tau$, $\nu^{\rm old}$ or the number of steps.
\Cref{sec:minimum,sec:inner} record the relative minimum $d_{M,\mathcal T}$ of the computed trajectory, defined in \eqref{eq:dip}.

\begin{proposition}[Time consistency on a fixed grid]\label{prop:step-consistency}
Fix $M\ge3$, its operator $L$, and $\nu^{\rm old}\in\mathcal P_M^+$.
As $\tau\downarrow0$, the successor density $\nu^*_\tau$ of \eqref{eq:step}, that is, $\nu^*$ of \cref{prop:positive} for the time step $\tau$, satisfies
\[
 \nu^*_\tau=\nu^{\rm old}
 -\tau L\bigl(\nu^{\rm old}L\log\nu^{\rm old}\bigr)+O(\tau^2).
\]
The remainder is in the Euclidean norm; its constant can depend on $M$ and $\nu^{\rm old}$.
\end{proposition}
The proof in \cref{app:step-consistency} identifies the discrete vector field approximated by one step.
Replacing $L$ by $\partial_{xx}$ gives \eqref{eq:dlss} for smooth positive profiles.
Since the constant in $O(\tau^2)$ can depend on $M$ and $\nu^{\rm old}$, \cref{prop:step-consistency} gives no estimate as $M\to\infty$ or as $\min_i\nu_i^{\rm old}\to0$.

\begin{table}[htbp]
\centering\small
\begin{tabular}{@{}p{.22\linewidth}p{.70\linewidth}@{}}
\toprule
Symbol & Meaning\\
\midrule
$M$, $N$, $h=1/N$ & Spatial nodes, geodesic time cells, and cell length.\\
$\nu^j$, $q^j$ & Nodal density and midpoint weight $(\nu^j+\nu^{j+1})/2$.\\
$m^j$, $\kappa^j$ & Integrated geodesic momentum and signed diffusivity $m^j/(hq^j)$ where $q^j>0$.\\
$\nu^{\rm old}$, $\nu^*$, $\tau$ & Prescribed density, successor density, and time increment in an entropic step.\\
$L$, $\C$ & Periodic Fourier second-derivative operator on $\R^M$ and its inverse on mean-zero vectors.\\
$a$, $\A$, $\E$ & Action density, grid-integrated action, and entropy.\\
$\beta$, $\gamma$ & DR momentum weight and proximal parameter in \cref{sec:dr-algorithm}.\\
$b$, $\omega$ & Local proximal scales: $b=\gamma\dx/(h\beta)$ for geodesics; $\omega=\gamma\dx$ and $b=\omega/(\tau\beta)$ for an entropic step.\\
\bottomrule
\end{tabular}
\caption{Main unknowns, operators, and parameters. Geodesic momentum is integrated over a cell of length $h$; entropic momentum is integrated over the step of length $\tau$.}
\label{tab:notation}
\end{table}

\section{Douglas--Rachford splitting}\label{sec:dr}

We solve \eqref{eq:geodesic} and \eqref{eq:step} by Douglas--Rachford (DR) splitting.
Each DR iteration first applies the proximal map of the objective, which includes the constraint $\nu\ge0$ and splits into independent problems at each node and, for geodesics, each time cell.
It then applies the affine projection, the orthogonal projection onto the set defined by the affine constraints, to the reflected point $2x^\ell-z^\ell$ of \eqref{eq:dr}; we compute this projection with Fourier transforms in space and, for geodesics, sine transforms in time (\cref{app:projection}).
The proximal output and the projected output need not coincide at a finite iteration; under the assumptions of \cref{prop:drglobal}, both converge to the same minimizer.

\begin{remark}[Related proximal methods]
Proximal splitting for dynamical optimal transport was developed by Papadakis, Peyr\'e, and Oudet \cite{PPO2014} and extended to mixed finite elements by Natale and Todeschi \cite{NataleTodeschi22}.
\end{remark}

\subsection{The DR iteration and its outputs}\label{sec:dr-algorithm}

For geodesics, treat the midpoint weight $q$ as an independent variable and write
\[
 \mathbb H=(\R^M)^{N-1}\times(\R^M)^N\times(\R^M)^N,
 \qquad U=(\nu,m,q).
\]
For a set $S$, its indicator $\iota_S$ is zero on $S$ and $+\infty$ outside $S$.
Define
\[
 \begin{aligned}
 F(U)&=\sum_{j=0}^{N-1}h^{-1}\A(q^j,m^j)+\iota_{\{\nu\ge0\}}(\nu),\\
 \K&=\{U\in\mathbb H:Lm^j=\nu^{j+1}-\nu^j,\quad
                  2q^j=\nu^j+\nu^{j+1},\quad0\le j<N\}.
 \end{aligned}
\]
The endpoints are fixed in $\K$.
Thus the geodesic problem is $\min(F+\iota_\K)$.
With $q$ as an independent variable, $F$ is a sum of functions of single entries $\nu_i^j$ and of pairs $(q_i^j,m_i^j)$, so its proximal map is computed separately for each entry and each pair; the projection onto $\K$ imposes $2q^j=\nu^j+\nu^{j+1}$.

For $\beta>0$, use the norm
\begin{equation}\label{eq:metric}
 \|U\|_\beta^2=
 \sum_{j=1}^{N-1}\sum_{i=0}^{M-1}|\nu_i^j|^2+
 \sum_{j=0}^{N-1}\sum_{i=0}^{M-1}
       \bigl(\beta|m_i^j|^2+|q_i^j|^2\bigr).
\end{equation}
Let $Q_\beta$ replace $m$ by $v=\sqrt\beta m$ and leave the other components unchanged.
The norm is Euclidean in these scaled coordinates.
For geodesics, set $J=F\circ Q_\beta^{-1}$ and $\mathcal S_\beta=Q_\beta\K$.
For an entropic step, use the same momentum scaling and set
\[
 J=H_\tau\circ Q_\beta^{-1},\qquad
 \mathcal S_\beta=\{(\nu,v):Lv=\sqrt\beta(\nu-\nu^{\rm old})\},
 \qquad v=\sqrt\beta m.
\]

For $\gamma>0$, the proximal map and the affine projection are
\begin{equation}\label{eq:prox-projection}
 \prox_{\gamma J}(z)=\argmin_x\{\gamma J(x)+\tfrac12\|x-z\|_2^2\},
 \qquad
 P_{\mathcal S_\beta}(z)=\argmin_{x\in\mathcal S_\beta}\tfrac12\|x-z\|_2^2.
\end{equation}
Both maps are single valued.
Starting from any auxiliary iterate $z^0$, DR repeats
\begin{equation}\label{eq:dr}
 \begin{aligned}
 x^\ell&=\prox_{\gamma J}(z^\ell),\\
 y^\ell&=P_{\mathcal S_\beta}(2x^\ell-z^\ell),\\
 z^{\ell+1}&=z^\ell+y^\ell-x^\ell.
 \end{aligned}
\end{equation}
The physical outputs are $X^\ell=Q_\beta^{-1}x^\ell$ and $Y^\ell=Q_\beta^{-1}y^\ell$.
For geodesics, $X^\ell$ has nonnegative densities and midpoint weights, with zero momentum wherever its weight is zero.
The projected output $Y^\ell$ satisfies conservation and the midpoint equations, but need not be nonnegative.
The residual in \cref{sec:geodata} measures the difference between these two outputs.

\begin{proposition}[Global convergence of DR]\label{prop:drglobal}
Fix $M\ge3$ and $\beta,\gamma>0$.
Consider either the geodesic problem \eqref{eq:geodesic} with $N\ge2$ and $\nu^0,\nu^N\in\mathcal P_M^+$, or the entropic problem \eqref{eq:step} with $\tau>0$ and $\nu^{\rm old}\in\mathcal P_M^+$.
Keep these data and parameters fixed and evaluate the proximal map and affine projection exactly.
For every initial vector $z^0$, there exists a DR fixed point $z_*$ such that
\[
 z^\ell\longrightarrow z_*,\qquad
 x^\ell\longrightarrow x_*,\qquad
 y^\ell\longrightarrow x_*,\qquad
 x_*:=\prox_{\gamma J}(z_*)\in\argmin_{x\in\mathcal S_\beta}J(x).
\]
\end{proposition}
The proof is in \cref{app:drglobal}.
The parameters change the iteration, not the minimization problem.
The proposition gives no iteration count for a prescribed accuracy.
The local proximal maps are computed by solving scalar equations to finite accuracy.
By \cref{prop:prox-error}, if a bound $\varepsilon_\ell$ on the Euclidean norm of the local scalar residuals at iteration $\ell$ satisfies $\sum_\ell\varepsilon_\ell<\infty$, the conclusion of \cref{prop:drglobal} persists.
The scalar solves of \cref{alg:geo-prox,alg:jko-prox} stop after a fixed number of Newton iterations or when the Newton step falls below $10^{-13}$; neither stopping rule implies $\sum_\ell\varepsilon_\ell<\infty$, so \cref{prop:drglobal} is not established for the computed iterates.

For geodesics, $a(q,\beta^{-1/2}v)=\beta^{-1}a(q,v)$ and $\A$ contains the factor $\dx$, so the action term of cell $j$ in $\gamma J$ is
\begin{equation}\label{eq:b}
 \gamma h^{-1}\A(q^j,\beta^{-1/2}v^j)=\sum_{i=0}^{M-1}b\,a(q_i^j,v_i^j),\qquad
 b=\frac{\gamma\dx}{h\beta}.
\end{equation}
The proximal map of the action term therefore splits into $MN$ two-dimensional maps $\prox_{b\,a}$, which depend on $\gamma$, $\beta$, $h$ and $\dx$ only through the local proximal scale $b$.
The nodal-density update is $\nu_{\rm in}\mapsto\max(\nu_{\rm in},0)$.
For a scaled action input $(q_{\rm in},v_{\rm in})$, \eqref{eq:prox} returns $(0,0)$ when $2bq_{\rm in}+v_{\rm in}^2\le0$ and otherwise determines a positive weight from one scalar equation.
Except in the comparison of \cref{rem:checked-prox}, the geodesic DR runs of \cref{sec:geotests} solve this equation by 18 Newton iterations (\cref{alg:geo-prox}); the comparison uses a solver that accepts a positive weight only if the local scalar residual (\cref{app:prox-error}) is below a tolerance.
The affine projection uses spatial Fourier and temporal sine transforms, with a total cost of $O(NM(\log M+\log N))$ per DR iteration, independently of $\beta$ (\cref{app:projection}).

\subsection{Slow DR convergence}\label{sec:conditioning}\label{sec:dr-local}

Let $z_*$ be a DR fixed point, $x_*=\prox_{\gamma J}(z_*)$, and $\delta z=z-z_*$ the auxiliary error.
When the proximal map is differentiable at $z_*$, let $R=D\prox_{\gamma J}(z_*)$ be its Jacobian:
\[
 \prox_{\gamma J}(z_*+\delta z)=x_*+R\,\delta z+o(\|\delta z\|_2).
\]
Let $P$ be the orthogonal projector onto the direction space $\mathcal S_\beta^0=\mathcal S_\beta-x_*$, the space of perturbations that preserve the affine constraints.
Linearizing \eqref{eq:dr} gives
\begin{equation}\label{eq:T}
 \delta z^+=T\,\delta z,\qquad T=I-R+P(2R-I),
\end{equation}
where $I$ is the identity and $\delta z^+$ is the error after one linearized iteration.
If the spectral radius $\spr(T)$, the largest modulus of an eigenvalue of $T$, is less than one, the local error has a geometric bound with any factor $c\in(\spr(T),1)$.
If $T$ has a fixed direction, then $1$ is an eigenvalue of $T$, and the linearization gives no proved rate here.
\Cref{app:lfp} states the precise local theorem, including the additional assumptions at zero outputs.

If all nodal densities of the minimizing path are positive, the Jacobian $R$ leaves nodal-density perturbations unchanged.
At an action pair $(q,v)$ with $q>0$, put $s=v/q$.
Since $a(\lambda q,\lambda v)=\lambda\,a(q,v)$ for $\lambda>0$, the gradient satisfies $\nabla a(\lambda q,\lambda v)=\nabla a(q,v)$, and differentiating at $\lambda=1$ gives $\nabla^2a(q,v)\,(q,v)^{\mathsf T}=0$.
Hence the Hessian $\nabla^2a(q,v)$ has the eigenvalue zero in direction $(1,s)$, which we call the flat direction of the pair.
The corresponding block of $R$ leaves this direction unchanged and multiplies the perpendicular direction $(s,-1)$ by
\[
 \frac{q}{q+b(1+s^2)}.
\]

The flat subspace is $\ker(I-R)$; it contains the nodal-density perturbations and the flat directions $(1,s)$ of all action pairs.
A vector $\delta x\in\mathcal S_\beta^0\cap\ker(I-R)$ satisfies $T\delta x=(I-R)\delta x+P(2R-I)\delta x=P\delta x=\delta x$, so it is fixed by $T$.
For the uniform path at $b=1$, \eqref{eq:rate} gives $\spr(T)=\cos^2\theta_{\min}$ when this exceeds $1/2$, which tends to one as the smallest angle $\theta_{\min}$ between $\mathcal S_\beta^0$ and $\ker(I-R)$ tends to zero.
These angles are measured in the scaled coordinates $(\nu,v,q)$ and therefore depend on $\beta$.

\subsection{Choosing \texorpdfstring{$\beta$ and $b$}{beta and b}}\label{sec:dr-parameters}\label{sec:dr-choices}

At a fixed minimizing path, changing $\beta$ at fixed $b$ changes $\mathcal S_\beta^0$ and the flat directions $(1,s)$, with $s=\sqrt\beta\,h\kappa$ in scaled coordinates.
Changing $b$ at fixed $\beta$ changes only the factor $q/(q+b(1+s^2))$ and leaves $P$ and the flat directions unchanged.
By \eqref{eq:b}, varying $\beta$ at fixed $b$ requires $\gamma=bh\beta/\dx$.
\begin{remark}[Scope of global rate bounds]
Giselsson and Boyd \cite{GiselssonBoyd17} analyze metric and parameter selection under strong convexity and smoothness assumptions.
The flat directions of the geodesic action prevent a direct application of those global rate bounds; we use the local linearization instead.
\end{remark}

\paragraph{Illustration on a uniform geodesic: an exact angle calculation.}\label{sec:dr-uniform}
For endpoints $\nu^0=\nu^N=\one$, the minimizing path has $\nu^j=q^j=\one$ and $m^j=0$.
At this uniform path all local blocks of $R$ are equal:
\[
 R(\delta\nu,\delta v,\delta q)
   =(\delta\nu,(1+b)^{-1}\delta v,\delta q).
\]
Spatial Fourier and temporal sine transforms therefore reduce $T$ to small invariant blocks, and $\spr(T)$ depends only on $M$, $N$, $\beta$ and $b$.
At a nonzero spatial Fourier frequency $k$, conservation gives
\[
 \widehat{\delta v}_k^j
 =\frac{\sqrt\beta}{\widehat L(k)}
   (\widehat{\delta\nu}_k^{j+1}-\widehat{\delta\nu}_k^j),
 \qquad \widehat L(k)=-4\pi^2k^2.
\]
Hence $|\widehat{\delta v}_k^j|=\sqrt\beta\,|\widehat{\delta\nu}_k^{j+1}-\widehat{\delta\nu}_k^j|/(4\pi^2k^2)$.
At fixed $\beta$, the scaled momentum of a perturbation in $\mathcal S_\beta^0$ decreases as $k^{-2}$ relative to its density increment.
Increasing $\beta$ increases this ratio, and therefore the angle between $\mathcal S_\beta^0$ and the flat subspace $\delta v=0$.

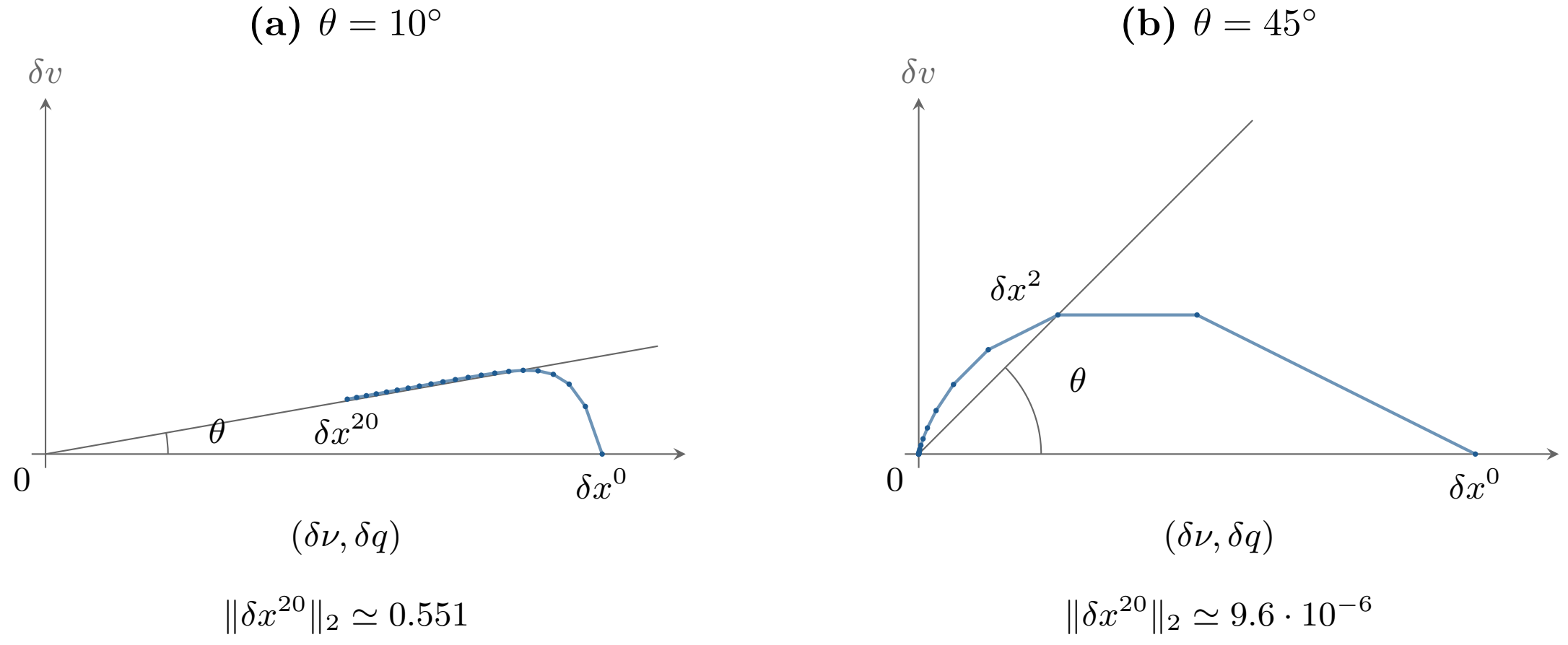
\begin{figure}[htbp]
\centering
\begin{tikzpicture}[font=\footnotesize,>=stealth]
\begin{scope}[xshift=0cm,x=5.1cm,y=5.1cm]
\node[font=\small\bfseries] at (.54,.77) {(a) $\theta=10^\circ$};
\draw[->,black!60] (-.025,0)--(1.15,0);
\draw[->,black!60] (0,-.025)--(0,.64) node[above] {$\delta v$};
\node at (.54,-.15) {$(\delta\nu,\delta q)$};
\node[below left] at (0,0) {$0$};
\draw[black!60] (0,0)--(1.10000000,0.19395968);
\draw[black!60] (.22,0) arc[start angle=0,end angle=10,radius=.22];
\node at (0.30882036,0.04201828) {$\theta$};
\draw[plotblue!65,thick] plot coordinates {(1.00000000,0.00000000) (0.96984631,0.08550504) (0.94060187,0.12567926) (0.91223925,0.14326583) (0.88473187,0.14963396) (0.85805394,0.15046601) (0.83218045,0.14860094) (0.80708714,0.14545609) (0.78275048,0.14173806) (0.75914767,0.13779814) (0.73625656,0.13381002) (0.71405571,0.12985865) (0.69252430,0.12598469) (0.67164213,0.12220666) (0.65138965,0.11853211) (0.63174784,0.11496315) (0.61269832,0.11149920) (0.59422320,0.10813839) (0.57630518,0.10487827) (0.55892745,0.10171613) (0.54207373,0.09864918)};
\fill[plotblue] (1.00000000,0.00000000) circle[radius=.7pt];
\fill[plotblue] (0.96984631,0.08550504) circle[radius=.7pt];
\fill[plotblue] (0.94060187,0.12567926) circle[radius=.7pt];
\fill[plotblue] (0.91223925,0.14326583) circle[radius=.7pt];
\fill[plotblue] (0.88473187,0.14963396) circle[radius=.7pt];
\fill[plotblue] (0.85805394,0.15046601) circle[radius=.7pt];
\fill[plotblue] (0.83218045,0.14860094) circle[radius=.7pt];
\fill[plotblue] (0.80708714,0.14545609) circle[radius=.7pt];
\fill[plotblue] (0.78275048,0.14173806) circle[radius=.7pt];
\fill[plotblue] (0.75914767,0.13779814) circle[radius=.7pt];
\fill[plotblue] (0.73625656,0.13381002) circle[radius=.7pt];
\fill[plotblue] (0.71405571,0.12985865) circle[radius=.7pt];
\fill[plotblue] (0.69252430,0.12598469) circle[radius=.7pt];
\fill[plotblue] (0.67164213,0.12220666) circle[radius=.7pt];
\fill[plotblue] (0.65138965,0.11853211) circle[radius=.7pt];
\fill[plotblue] (0.63174784,0.11496315) circle[radius=.7pt];
\fill[plotblue] (0.61269832,0.11149920) circle[radius=.7pt];
\fill[plotblue] (0.59422320,0.10813839) circle[radius=.7pt];
\fill[plotblue] (0.57630518,0.10487827) circle[radius=.7pt];
\fill[plotblue] (0.55892745,0.10171613) circle[radius=.7pt];
\fill[plotblue] (0.54207373,0.09864918) circle[radius=.7pt];
\node[below] at (1,0) {$\delta x^0$};
\node[below] at (0.54207373,0.09864918) {$\delta x^{20}$};
\node at (.54,-.29) {$\|\delta x^{20}\|_2\simeq 0.551$};
\end{scope}
\begin{scope}[xshift=8cm,x=5.1cm,y=5.1cm]
\node[font=\small\bfseries] at (.54,.77) {(b) $\theta=45^\circ$};
\draw[->,black!60] (-.025,0)--(1.15,0);
\draw[->,black!60] (0,-.025)--(0,.64) node[above] {$\delta v$};
\node at (.54,-.15) {$(\delta\nu,\delta q)$};
\node[below left] at (0,0) {$0$};
\draw[black!60] (0,0)--(0.60000000,0.60000000);
\draw[black!60] (.22,0) arc[start angle=0,end angle=45,radius=.22];
\node at (0.28640266,0.13363186) {$\theta$};
\draw[plotblue!65,thick] plot coordinates {(1.00000000,0.00000000) (0.50000000,0.25000000) (0.25000000,0.25000000) (0.12500000,0.18750000) (0.06250000,0.12500000) (0.03125000,0.07812500) (0.01562500,0.04687500) (0.00781250,0.02734375) (0.00390625,0.01562500) (0.00195312,0.00878906) (0.00097656,0.00488281) (0.00048828,0.00268555) (0.00024414,0.00146484) (0.00012207,0.00079346) (0.00006104,0.00042725) (0.00003052,0.00022888) (0.00001526,0.00012207) (0.00000763,0.00006485) (0.00000381,0.00003433) (0.00000191,0.00001812) (0.00000095,0.00000954)};
\fill[plotblue] (1.00000000,0.00000000) circle[radius=.7pt];
\fill[plotblue] (0.50000000,0.25000000) circle[radius=.7pt];
\fill[plotblue] (0.25000000,0.25000000) circle[radius=.7pt];
\fill[plotblue] (0.12500000,0.18750000) circle[radius=.7pt];
\fill[plotblue] (0.06250000,0.12500000) circle[radius=.7pt];
\fill[plotblue] (0.03125000,0.07812500) circle[radius=.7pt];
\fill[plotblue] (0.01562500,0.04687500) circle[radius=.7pt];
\fill[plotblue] (0.00781250,0.02734375) circle[radius=.7pt];
\fill[plotblue] (0.00390625,0.01562500) circle[radius=.7pt];
\fill[plotblue] (0.00195312,0.00878906) circle[radius=.7pt];
\fill[plotblue] (0.00097656,0.00488281) circle[radius=.7pt];
\fill[plotblue] (0.00048828,0.00268555) circle[radius=.7pt];
\fill[plotblue] (0.00024414,0.00146484) circle[radius=.7pt];
\fill[plotblue] (0.00012207,0.00079346) circle[radius=.7pt];
\fill[plotblue] (0.00006104,0.00042725) circle[radius=.7pt];
\fill[plotblue] (0.00003052,0.00022888) circle[radius=.7pt];
\fill[plotblue] (0.00001526,0.00012207) circle[radius=.7pt];
\fill[plotblue] (0.00000763,0.00006485) circle[radius=.7pt];
\fill[plotblue] (0.00000381,0.00003433) circle[radius=.7pt];
\fill[plotblue] (0.00000191,0.00001812) circle[radius=.7pt];
\fill[plotblue] (0.00000095,0.00000954) circle[radius=.7pt];
\node[below] at (1,0) {$\delta x^0$};
\node[above left] at (0.25000000,0.25000000) {$\delta x^2$};
\node at (.54,-.29) {$\|\delta x^{20}\|_2\simeq 9.6\cdot10^{-6}$};
\end{scope}
\end{tikzpicture}
\caption{Linearized DR iteration in one Fourier--sine mode of the uniform path, with $b=1$.
Horizontal axis: unit direction in the $(\delta\nu,\delta q)$ coordinates, on which $R$ is the identity; vertical axis: scaled momentum $\delta v$.
Gray line: intersection of the direction space $\mathcal S_\beta^0$ with the plane of the mode, at a prescribed angle $\theta$ to the horizontal axis; $\theta=10^\circ$ in (a) and $\theta=45^\circ$ in (b).
Both panels start from $\delta z^0=(1,0)$ and use the same scale.
Blue dots: proximal-output errors $\delta x^\ell=RT^\ell\delta z^0$, $0\le\ell\le20$, joined in iteration order; $\|\delta x^{20}\|_2$ is printed below each panel.
These are linearized mode trajectories, not nonlinear geodesic solves.}
\label{fig:drangle}
\end{figure}

Let $\theta_{\min}$ be the smallest angle between a nonzero vector in $\mathcal S_\beta^0$ and that flat subspace.
For every $M\ge3$, the exact uniform-path calculation at $b=1$ gives
\begin{equation}\label{eq:rate}
 \spr(T)=\max\{1/2,\cos^2\theta_{\min}\}.
\end{equation}
\Cref{app:spectrum} treats each Fourier index $k\in\mathcal I_M$ separately; for even $M$, the unpaired index $k=-M/2$ corresponds to the real eigenvector $((-1)^i)_{0\le i<M}$ of $L$.
Since $\theta_{\min}$ increases with $\beta$ (\cref{app:spectrum}), $\spr(T)$ is nonincreasing in $\beta$ and equals $1/2$ once $\cos^2\theta_{\min}\le1/2$.
The value $\spr(T)$ is an asymptotic factor, not an error reduction guaranteed at each iteration.
\Cref{fig:drangle} shows twenty linearized iterations at $\theta=10^\circ$, where $\cos^2\theta\approx0.97$ and $\|\delta x^{20}\|_2\simeq0.551$, and at $\theta=45^\circ$, where $\cos^2\theta=1/2$ and $\|\delta x^{20}\|_2\simeq9.6\cdot10^{-6}$.
\Cref{app:spectrum} gives the mode angles, the dependence on the grid, and the calculation for general $b$.

\paragraph{Nonuniform geodesics: rotation of the flat directions with $\beta$.}
At nonzero momentum, increasing $\beta$ also rotates each local flat direction $(1,s)$, since
\begin{equation}\label{eq:scaledkappa}
 s=\sqrt\beta\,h\kappa,
\end{equation}
with $\kappa$ from \eqref{eq:kappa}.
The angle between $(1,s)$ and the momentum axis is $\arctan(1/|s|)$, which decreases as $\beta$ increases.
Perturbations in $\mathcal S_\beta^0$ satisfy $\delta v^j-\sqrt\beta\,\C(\delta\nu^{j+1}-\delta\nu^j)\in\R\one$, so their scaled momentum also grows with $\beta$.
In cells with large $|s|$, the flat direction and these perturbations both approach the momentum axis, so the angle between $\mathcal S_\beta^0$ and $\ker(I-R)$ may decrease.
This is a heuristic explanation; we have no bound on $\spr(T)$ for nonuniform paths.
Spatially varying proximal blocks couple Fourier frequencies, so the uniform calculation no longer applies.
In \cref{sec:betasweep}, among the tested values of $\beta$, the relative action gap \eqref{eq:comparison-gap} after 3000 iterations is smallest at $\beta=10^6$ or $\beta=3\cdot10^6$ for each endpoint pair (\cref{fig:betasweep}).

Increasing $b$ decreases the factor $q/(q+b(1+s^2))$ perpendicular to each flat direction.
In the uniform mode at $\theta=45^\circ$, the DR trajectories oscillate for $b>1$, and $b=1$ leaves less error after twenty iterations than $b=3$ or $b=10$ (\cref{app:drillustrations}).
In the nonuniform tests of \cref{sec:drtests}, at $\beta=10^6$, $b=3$ gives the fewest iterations per factor-ten decrease of the DR residual $r_\ell$ of \cref{sec:geodata} over $10^{-5}\ge r_\ell\ge10^{-9}$, among $b\in\{0.3,1,3,10\}$, for each of the three endpoint pairs (\cref{tab:drstep}).
Unless stated otherwise, the geodesic DR tests of \cref{sec:geotests} use $\beta=10^6$ and $b=1$.

\subsection{DR convergence for entropic steps}\label{sec:drflow}

For the entropic problem \eqref{eq:step}, keep $\nu^{\rm old}>0$ and $\tau>0$ fixed during each DR solve.

Set $v=\sqrt\beta m$, $\omega=\gamma\dx$, and $b=\omega/(\tau\beta)$.
For input $(\nu_{{\rm in},i},v_{{\rm in},i})$, the local update is
\begin{equation}\label{eq:entropic-prox}
 (\nu_i,v_i)=\prox_\phi(\nu_{{\rm in},i},v_{{\rm in},i}),\qquad
 \phi(\nu,v)=b\,a(\nu,v)+\omega e(\nu).
\end{equation}
The proximal density output is strictly positive and is found from one scalar equation (\cref{app:localprox}).
The implementation solves for log-density to preserve positivity (\cref{alg:jko-prox}).

\Needspace{7\baselineskip}
The Hessian $H_{\rm JKO}=\nabla^2\phi(\nu_i,v_i)$ is positive definite at every positive output (\cref{app:lfp}).
The corresponding block of the Jacobian $R$ is $(I+H_{\rm JKO})^{-1}$, a symmetric matrix with eigenvalues in $(0,1)$; at the minimizer, therefore, $\|2R-I\|_2<1$.
Since $T=\tfrac12[I+(2P-I)(2R-I)]$ and $\|2P-I\|_2=1$, this gives $\|T\|_2\le\tfrac12(1+\|2R-I\|_2)<1$, hence $\spr(T)\le\|T\|_2<1$.
Thus exact DR with fixed parameters converges locally at a geometric rate for every fixed entropic problem (\cref{app:lfp}).
The bound $\tfrac12(1+\|2R-I\|_2)$ tends to one as $\min_i\nu_i^*\to0$: at node $i$, the largest eigenvalue of $H_{\rm JKO}$ is at least $\omega/\nu_i^*$, so the smallest eigenvalue of the corresponding block of $R$ is at most $\nu_i^*/(\nu_i^*+\omega)$.

To choose $\beta$ for the entropic step, first consider uniform density and zero momentum, where $H_{\rm JKO}=\diag(\omega,b)$.
Setting $b=\omega$, that is $\beta=1/\tau$, gives $H_{\rm JKO}=\omega I$ and $v=m/\sqrt\tau$.
For uniform input, $\omega=1$ then gives $R=T=I/2$.
We retain $\beta=1/\tau$ for nonuniform input and choose $\omega$ from an approximate model.
Replace the successor density $\nu^*$ by $\nu^{\rm old}$ and its momentum-to-density ratio by zero.
Define the model Jacobian $R_{\rm model}$ by local blocks $\varrho_iI_2$, where $\varrho_i=\nu_i^{\rm old}/(\nu_i^{\rm old}+\omega)$ and $I_2$ is the two-dimensional identity.
With the projector $P$ unchanged, its DR matrix is $T_{\rm model}=I-R_{\rm model}+P(2R_{\rm model}-I)$.
Writing $\nu_{\min}=\min_i\nu_i^{\rm old}$, $\nu_{\max}=\max_i\nu_i^{\rm old}$, and $C=\nu_{\max}/\nu_{\min}$, the parameter rule is
\[
 \omega=\sqrt{\nu_{\min}\nu_{\max}},\qquad
 \|T_{\rm model}\|_2\le\frac{\sqrt C}{1+\sqrt C}.
\]
\Cref{app:entropicbalance} obtains $\omega=\sqrt{\nu_{\min}\nu_{\max}}$ by minimizing over $\omega>0$ the upper bound on $\|T_{\rm model}\|_2$ derived there.
With this $\omega$, $\lceil(\sqrt C+1)\log10\rceil$ iterations of $\delta z\mapsto T_{\rm model}\delta z$ reduce $\|\delta z\|_2$ by at least a factor ten.
The bound $\sqrt C/(1+\sqrt C)$ concerns $T_{\rm model}$, which uses $\nu^{\rm old}$ and zero momentum; it is not a bound on $\|T\|_2$ at the minimizer $(\nu^*,m^*)$ of the entropic step.

\section{Interior-point methods and validation}\label{sec:ipm}

We also write both problems as conic programs and solve them by an interior-point method (IPM).
The IPM solutions are computed independently of DR, and the IPM dual multipliers are used in the objective bounds of \cref{sec:certificate}.
We represent the action epigraph by second-order cones and the entropy epigraph by exponential cones.
The resulting programs have linear objectives and affine constraints in addition to these cones.

\begin{remark}[Related interior-point methods]
Barrier methods for dynamical optimal transport include the finite-volume solver of Natale and Todeschi \cite{NataleTodeschi21} and the preconditioned interior-point method of Facca, Todeschi, Natale, and Benzi \cite{FaccaTodeschiNataleBenzi24}.
\end{remark}

\subsection{Equivalent conic problem}
We use the second-order cone and the exponential cone:
\[
 \mathcal Q=\{(y_1,y_2,y_3)\in\R^3:y_1\ge\sqrt{y_2^2+y_3^2}\},\qquad
 K_{\exp}=\overline{\{(y_1,y_2,y_3)\in\R^3:y_2>0,\ y_2e^{y_1/y_2}\le y_3\}}.
\]
The overline denotes closure in $\R^3$.
For a fixed $z>0$, their constraints give
\[
 \begin{aligned}
 (\nu+t,\nu-t,\sqrt2\,v)\in\mathcal Q
 &\ \Longleftrightarrow\ \nu,t\ge0,\quad 2\nu t\ge v^2,\\
 (-u,\nu,z)\in K_{\exp}
 &\ \Longleftrightarrow\ \nu\ge0,\quad u\ge\nu\log(\nu/z).
 \end{aligned}
\]
The convention at $\nu=0$ is $\nu\log(\nu/z)=0$.
\Cref{app:conic} proves both equivalences, including the case $\nu=0$.

Set $m=\sqrt\tau\,v$, so that $\tau^{-1}\A(\nu,m)=\A(\nu,v)$, and $z_i=\nu_i^{\rm old}$.
The entropy identity is
\[
 \E(\nu)=\avg{\nu\log(\nu/\nu^{\rm old})}
             +\ip{\log\nu^{\rm old}-1}{\nu}+1.
\]
The entropic step has the equivalent conic form
\begin{equation}\label{eq:conic}
 \begin{aligned}
 \min_{\nu,v,t,u\in\R^M}\quad
 &\sum_{i=0}^{M-1}[t_i+u_i+(\log\nu_i^{\rm old}-1)\nu_i],\\
 \text{subject to}\quad&\sqrt\tau\,Lv-\nu=-\nu^{\rm old},\\
 &(\nu_i+t_i,\nu_i-t_i,\sqrt2\,v_i)\in\mathcal Q,\\
 &(-u_i,\nu_i,\nu_i^{\rm old})\in K_{\exp}\quad(0\le i<M).
 \end{aligned}
\end{equation}

The minimum in \eqref{eq:conic} equals $M(H_\tau^*-1)$.
The problem has $4M$ scalar variables, $M$ affine equations, and $2M$ three-dimensional cones.
At the optimum, $t_i=a(\nu_i,v_i)$ and $u_i=\nu_i\log(\nu_i/\nu_i^{\rm old})$, with $0\log0=0$.
The cones impose only $\nu\ge0$; the positivity of the successor density $\nu^*$ follows from \cref{prop:positive}, not from an added constraint.

The recorded implementation uses Clarabel 0.11.1 \cite{Clarabel}.
It also uses a positive scaling $d_i=\sqrt{\nu_i^{\rm old}}$ at each node.
Write each of the four variables as $d_i$ times a new variable.
Divide the cone entries and conservation row $i$ by $d_i$.
This substitution gives an equivalent program because the two cones are positively homogeneous.
The matrix in the scaled conservation equation is $\diag(d)^{-1}L\diag(d)$.
The dual evaluation must undo this row scaling; \cref{app:records} gives the multiplier conversion.

Each IPM iteration factorizes a system containing the dense $M\times M$ matrix $L$, at cost $O(M^3)$.
In \cref{tab:flowgridsize}, from $M=128$ to $M=2048$, the solver-reported solve time of one entropic step grows from $0.033$\,s to $374$\,s, by factors between $7.6$ and $12.9$ at each doubling of $M$; per IPM iteration, the factors lie between $6$ and $10$.

\subsection{Objective upper bound and dual lower bound}\label{sec:certificate}

For the entropic problem \eqref{eq:step}, let $\nu\in\mathcal P_M$ be an approximate minimizing density.
The reduced objective $J_\tau(\nu)$ defined below is an upper bound on $H_\tau^*$, and the dual functional $D_\tau(\chi)$ in \eqref{eq:dual} is a lower bound on $H_\tau^*$ for every $\chi$.
\Cref{prop:gap} converts the computable difference $J_\tau(\nu)-D_\tau(\chi)$ into bounds on $J_\tau(\nu)-H_\tau^*$ and on $\normone{\nu-\nu^*}$.

For $\nu\in\mathcal P_M$, define the reduced objective by
\[
 J_\tau(\nu)=\inf\{H_\tau(\nu,m):Lm=\nu-\nu^{\rm old}\}.
\]
The infimum is $+\infty$ if no momentum has finite action.
For $\nu>0$, set $m_\perp=\C(\nu-\nu^{\rm old})$.
If $\avg\nu=\avg{\nu^{\rm old}}$, then $Lm_\perp=\nu-\nu^{\rm old}$; since $\ker L=\R\one$, the feasible momenta are $m_\perp+c\one$, $c\in\R$.
The map $c\mapsto\A(\nu,m_\perp+c\one)$ is a strictly convex quadratic, and its minimizer gives the unique momentum of least action
\begin{equation}\label{eq:elimination}
 m=m_\perp+c\one,\qquad
 c=-\frac{\avg{m_\perp/\nu}}{\avg{1/\nu}}.
\end{equation}
Then $J_\tau(\nu)=H_\tau(\nu,m)$.
The momentum $m$ is feasible in exact arithmetic when $\avg\nu=\avg{\nu^{\rm old}}$.

For a multiplier $\chi\in\R^M$ of the conservation equation $Lm=\nu-\nu^{\rm old}$, define
\begin{equation}\label{eq:dual}
 D_\tau(\chi)=\ip\chi{\nu^{\rm old}}
 -\log\avg{\exp\bigl(\chi+\tfrac\tau2(L\chi)^2\bigr)}.
\end{equation}
Since $\avg{\nu^{\rm old}}=1$ and $L\one=0$, $D_\tau(\chi+c\one)=D_\tau(\chi)$ for every $c\in\R$.
Hence $\sup_\chi D_\tau(\chi)=\sup_{\avg\chi=0}D_\tau(\chi)$.

\begin{proposition}[Dual bound and density error]\label{prop:gap}
Let $(\nu^*,m^*)$ be the minimizer of \eqref{eq:step}.
For each $\nu\in\mathcal P_M$ with $J_\tau(\nu)<\infty$, and each $\chi\in\R^M$,
\begin{equation}\label{eq:gap}
 \begin{gathered}
 D_\tau(\chi)\le H_\tau^*\le J_\tau(\nu),\\
 \tfrac12\normone{\nu-\nu^*}^2
 \le J_\tau(\nu)-H_\tau^*
 \le J_\tau(\nu)-D_\tau(\chi).
 \end{gathered}
\end{equation}
An optimal multiplier exists, and $\max_\chi D_\tau(\chi)=H_\tau^*$.
\end{proposition}
\Cref{app:gap,app:pinsker} prove \cref{prop:gap} by square completion and the grid-weighted Pinsker inequality.

The density estimate relies on the entropy term and therefore does not apply to the geodesic action.

Small conic feasibility residuals do not make the solver-reported dual objective a rigorous lower bound.
The baseline check in \cref{tab:reporteddual} detects violations at the precision considered here.
We therefore evaluate \eqref{eq:dual}, with the mass correction \eqref{eq:dualmass} of \cref{app:mass} when the stored input does not have mass one; in exact arithmetic, this dual value is a lower bound for every multiplier.

\begin{table}[htbp]
\centering\small
\begin{tabular}{@{}rrr@{}}
\toprule
Steps & Steps with reported dual $>$ upper bound & Largest excess\\
\midrule
\BaselineDualSteps & \BaselineDualExceedances & $\BaselineDualExcess$\\
\bottomrule
\end{tabular}
\caption{Check of the solver-reported dual objectives for the baseline entropic trajectories in \cref{sec:flowdata}. The columns count the baseline steps, count the steps whose solver-reported dual objective, converted to the scale of $H_\tau$, exceeds an interval upper bound on $J_{\tau,\varsigma}(\nu)$, where $\nu$ is the density returned by the solver normalized to the mass $\varsigma=\avg{\nu^{\rm old}}$ of the stored input (\cref{app:mass}), and give the largest difference (reported dual minus upper bound). Since $J_{\tau,\varsigma}(\nu)\ge J_{\tau,\varsigma}(\nu_\varsigma^*)$, the minimum of \eqref{eq:step} for the stored input, a positive difference shows that the reported value is not a lower bound on this minimum.}
\label{tab:reporteddual}
\end{table}

\subsection{Interval validation and its scope}
A direct subtraction of two floating-point objective values does not certify a gap near machine precision.
Let $\nu_{\rm cand}>0$ be the density returned by the solver and $\chi$ a multiplier computed from the solver output (for IPM, from its equality multipliers, \cref{app:records}).
The validation treats the stored arrays of the prescribed density $\nu^{\rm old}$ and of $\nu_{\rm cand}$, the time increment $\tau$ and the Fourier symbol supplied to the solver as exact data.
For the certificates of \cref{tab:baseline-validation,tab:flowcontrast} this symbol is the rounded symbol $\widehat L_{\rm round}$ of \cref{tab:operator-rounding}, not the exact symbol $\widehat L(k)=-4\pi^2k^2$; these certificates therefore concern \eqref{eq:step} with $L$ replaced by the operator of symbol $\widehat L_{\rm round}$, and \cref{app:records} gives the symbol of the other entropic certificates.
The proofs of \cref{prop:positive,prop:gap} and of \cref{app:mass} use only three properties of $L$: it is self-adjoint for $\ip\cdot\cdot$, $\ker L=\R\one$, and its range is $\R^M_0$.
These properties hold for every real symbol that is even in $k$, zero at $k=0$ and nonzero for $k\ne0$, in particular for $\widehat L_{\rm round}$.
Rounding of the stored input can make its exact mass $\varsigma=\avg{\nu^{\rm old}}$ differ from one.
In 133-bit interval arithmetic with outward rounding, set $\nu=\varsigma\,\nu_{\rm cand}/\avg{\nu_{\rm cand}}$, so that $\avg\nu=\varsigma$ and the momentum \eqref{eq:elimination} is feasible.
Enclose the reduced objective $J_{\tau,\varsigma}(\nu)$ and the dual value $D_{\tau,\varsigma}(\chi)$ of \eqref{eq:dualmass} for an input of mass $\varsigma$ (\cref{app:mass}); for $\varsigma=1$, they are $J_\tau(\nu)$ and $D_\tau(\chi)$.
Let $\eta$ be the upper endpoint of the enclosure of $J_{\tau,\varsigma}(\nu)-D_{\tau,\varsigma}(\chi)$.
\Cref{app:mass} gives $J_{\tau,\varsigma}(\nu)-J_{\tau,\varsigma}(\nu_\varsigma^*)\le\eta$ and $\normone{\nu-\nu_\varsigma^*}\le\sqrt{2\varsigma\eta}$, where $\nu_\varsigma^*$ is the successor density for the stored input.

If $\varsigma=1$, then $\nu_\varsigma^*=\nu^*$, and these bounds are those of \eqref{eq:gap}.
Each certificate concerns one discrete entropic problem and does not include discretization error.
A bound for the full trajectory also requires a stability argument for the propagation of errors between steps.

To assess the accuracy of individual density entries, we use the pointwise consequence of the discrete $L^1$ bound for $\varsigma=1$:
\begin{equation}\label{eq:pointwise}
 \|\nu-\nu^*\|_\infty
 \le M\normone{\nu-\nu^*}\le M\sqrt{2\eta}.
\end{equation}
The factor $M$ follows from the quadrature weight $1/M$.
For $\varsigma\ne1$, \eqref{eq:pointwise} holds with $\nu^*$ replaced by $\nu_\varsigma^*$ and $\eta$ by $\varsigma\eta$.
\begin{remark}[Absolute and relative accuracy]
For example, $M=128$ and $\eta=10^{-16}$, larger than the gap bound $\BaselineIntervalGap$ of \cref{tab:baseline-validation}, give $\|\nu-\nu^*\|_\infty\le1.8\cdot10^{-6}$.
For an entry $\nu_i^*=1$, this is a relative error of at most $1.8\cdot10^{-6}$.
For an entry $\nu_i^*=10^{-10}$, the bound exceeds the entry by a factor $1.8\cdot10^4$ and determines neither its value nor its order of magnitude.
\end{remark}

\subsection{Conic form for the geodesic comparison}\label{sec:geocone}
The geodesic problem \eqref{eq:geodesic} requires only second-order cones.
The conic form uses the momentum rate $w^j=m^j/h$.
For each cell and node, introduce an epigraph variable $t_i^j\ge a(q_i^j,w_i^j)$.
For $q^j=(\nu^j+\nu^{j+1})/2$, impose
\begin{equation}\label{eq:geocone}
 \left(\frac{q_i^j+t_i^j}{\sqrt2},
       \frac{q_i^j-t_i^j}{\sqrt2},w_i^j\right)\in\mathcal Q.
\end{equation}
This is equivalent to $q_i^j,t_i^j\ge0$ and $2q_i^jt_i^j\ge(w_i^j)^2$.
Minimize $h\dx\sum_{j,i}t_i^j$ subject to these constraints and
\[
 \nu^j\in\mathcal P_M\quad(1\le j<N),\qquad
 w^j=h^{-1}\C(\nu^{j+1}-\nu^j)+c^j\one\quad(0\le j<N),
\]
where $c^j\in\R$ is a free scalar for each time cell.
The endpoint densities stay fixed.
For $\nu^j,\nu^{j+1}\in\mathcal P_M$, $\avg{\nu^{j+1}-\nu^j}=0$, so $L\C(\nu^{j+1}-\nu^j)=\nu^{j+1}-\nu^j$.
Since $\ker L=\R\one$, the solutions of $Lm^j=\nu^{j+1}-\nu^j$ are $m^j=\C(\nu^{j+1}-\nu^j)+hc^j\one$, $c^j\in\R$, that is, the displayed constraint for $w^j=m^j/h$.
For fixed $(q^j,w^j)$, the minimum of $h\dx\sum_it_i^j$ under \eqref{eq:geocone} is $h\A(q^j,w^j)$, which equals $h^{-1}\A(q^j,m^j)$ because $w^j=m^j/h$.
At $q_i^j=0$, \eqref{eq:geocone} forces $w_i^j=0$ and $t_i^j\ge0$, which is the closed action $a(0,0)=0$.
Thus the conic problem has the same value as \eqref{eq:geodesic}.

\Cref{alg:geo-ipm} solves this conic problem with Clarabel.
It uses the constraints on $w^j$ displayed above, which involve $\C$ instead of $L$, minimizes $\sum_{j,i}t_i^j$, and multiplies the optimal conic value by $h\dx$ to obtain the action.
Finite interior-point iterates do not identify an exact zero set.
For the action comparison, we evaluate the original geodesic action on paths that satisfy its constraints, as specified in \cref{sec:geocompare}.

\subsection{Action evaluation for a prescribed path}\label{sec:geovalue}
For a positive unit-mass density path, conservation determines momentum up to one constant in each time cell.
The momentum of least action is
\begin{equation}\label{eq:geovalue}
 m_\perp^j=\C(\nu^{j+1}-\nu^j),\qquad
 c_\star^j=-\frac{\avg{m_\perp^j/q^j}}{\avg{1/q^j}},\qquad
 m^j=m_\perp^j+c_\star^j\one.
\end{equation}
This follows by minimizing the strictly convex quadratic action of $m_\perp^j+c\one$ over $c$.
Since \eqref{eq:geovalue} depends only on the density path, \cref{sec:geocompare} applies it to the paths built from the DR and IPM outputs.
The primal value $V=\sum_j h^{-1}\A(q^j,m^j)$ and the dual value $D_G(\chi)$ from \cref{app:geodual} satisfy
\[
 D_G(\chi)\le E_{M,N}\le V.
\]
These action bounds give no density-path error bound, so agreement of the density arrays is checked separately in \cref{sec:ipmgeovalidation}.

\section{Signed diffusive geodesics: accuracy and profiles}\label{sec:geotests}

We compare DR and IPM for three endpoint pairs, using action bounds to assess optimization accuracy and grid refinement to assess path sensitivity.
We then examine density minima, diffusivity sign changes, and the action cost of imposing nonnegative diffusivity.

\subsection{Data, residual and action evaluation}\label{sec:geodata}
We use three pairs of positive endpoints on the periodic interval.
For each profile $\mathfrak f$, the endpoint density is $(\mathfrak f+10^{-2})/\avg{\mathfrak f+10^{-2}}$.
For the peaked pair, $\nu^0$ and $\nu^N$ use $\mathfrak f(x)=\exp(-d_{\T}(x,x_c)^2/(2\cdot0.025^2))$ with the centers $x_c=1/4$ and $x_c=3/4$, respectively.
Here $d_{\T}(x,y)=\min_{n\in\mathbb Z}|x-y+n|$.
For the box pair, $\nu^0$ and $\nu^N$ use the indicator functions of $[0.20,0.40]$ and $[0.60,0.80]$, respectively.

For the diffuse pair, $\nu^0$ and $\nu^N$ use $\mathfrak f(x)=\exp(-(x-\tfrac12)^2/(2\cdot0.20^2))$ and $\mathfrak f(x)=\exp(-(x-\tfrac12)^2/(2\cdot0.10^2))$, respectively, for $x\in[0,1)$, extended periodically.
Unless stated otherwise, the tests use $M=63$, $N=64$, binary64 arithmetic, and linear interpolation as the initial density path.
The DR weight is $\beta=10^6$.

We record the DR residual
\[
 r_\ell=\frac{\|Y^\ell-X^\ell\|_\beta}{\|Y^0-X^0\|_\beta}.
\]
Here $(X^0,Y^0)$ are the outputs of the first DR iteration from $z^0$, which \cref{alg:geo-dr} numbers 1.
The denominator is nonzero in the reported runs.
The residual $r_\ell$ is not an objective gap or a density-error bound.

A finite-iteration solver output need not satisfy nonnegativity, unit mass and the conservation equations.
The upper bound $V\ge E_{M,N}$ of \cref{sec:geovalue}, valid in exact arithmetic, is therefore evaluated on a positive unit-mass path $\bar\nu$ with the prescribed endpoints, called the action-evaluation path; \cref{sec:geocompare} builds it from each solver output.

To study zero densities, we count the entries equal to zero in the density component of the DR proximal output $X^\ell$.
The projection onto $\mathcal P_M$ and the mixture with $\varepsilon_{\rm mix}\one$ of \cref{sec:geocompare} are not applied to this array; the action-evaluation path $\bar\nu$ has all interior entries at least $\varepsilon_{\rm mix}$ and is not used for this count.

\subsection{Convergence in iterations and elapsed time}\label{sec:geocompare}

The DR residual $r_\ell$ is the normalized norm $\|Y^\ell-X^\ell\|_\beta$ of the difference between the two DR outputs, whereas the IPM stopping test uses the feasibility residuals and primal--dual gaps of the conic problem of \cref{sec:geocone}.
These quantities have different definitions and scales.
To compare objective accuracy, we evaluate the action on paths satisfying the geodesic constraints and compare it with a common dual lower bound for each endpoint pair.

\begin{figure}[p]
\centering
\begin{tikzpicture}
\begin{groupplot}[group style={group size=2 by 3,horizontal sep=1.3cm,vertical sep=1.8cm},
 width=.43\textwidth,height=4.9cm,xmode=log,ymode=log,
 ymin=1e-11,ymax=10,ytick={1,1e-3,1e-6,1e-9},
 legend columns=2,legend style={draw=none,/tikz/column sep=6pt},unbounded coords=jump]
\nextgroupplot[title={Peaked},xlabel={},ylabel={$\gaprel(V)$},legend to name=legend:geocompare]
\addplot[plotgreen,thick,mark=*,mark size=1.2pt] coordinates {(1,1.991249049) (100,6.283928761e-05) (200,5.028030992e-06) (300,2.114965888e-06) (400,2.727584622e-08) (500,8.350806615e-09) (600,3.158654167e-09) (700,1.366304017e-09) (800,6.520559032e-10) (900,3.389598752e-10) (1000,1.898931802e-10) (1100,1.148771237e-10) (1200,7.520648243e-11) (1300,5.328743351e-11) (1400,4.075382087e-11) (1500,3.340030066e-11) (1600,2.893441125e-11) (1700,2.609922885e-11) (1800,2.427343075e-11) (1900,2.304817392e-11) (2000,2.221437158e-11) (2100,2.162471584e-11) (2200,2.119820149e-11) (2300,2.08839352e-11) (2400,2.06484122e-11) (2500,2.047000281e-11) (2600,2.033244946e-11) (2700,2.02254321e-11) (2800,2.014131674e-11) (2900,2.007501405e-11) (3000,2.002171743e-11)};
\addlegendentry{DR evaluations}
\addplot[plotblue,only marks,mark=square*,mark size=2.2pt] coordinates {(11,0.0002851673256) (14,1.566176892e-06) (17,1.607358986e-09) (18,1.175081673e-10)};
\addlegendentry{IPM completed runs}
\nextgroupplot[title={Peaked},xlabel={},ylabel={}]
\addplot[plotgreen,thick,mark=*,mark size=1.2pt] coordinates {(0.03704649,1.991249049) (0.109724697,6.283928761e-05) (0.182608017,5.028030992e-06) (0.255956221,2.114965888e-06) (0.329099248,2.727584622e-08) (0.402109304,8.350806615e-09) (0.475965819,3.158654167e-09) (0.549280499,1.366304017e-09) (0.6222788,6.520559032e-10) (0.695193552,3.389598752e-10) (0.768138632,1.898931802e-10) (0.841252155,1.148771237e-10) (0.913985286,7.520648243e-11) (0.987056965,5.328743351e-11) (1.06038959,4.075382087e-11) (1.133678612,3.340030066e-11) (1.20680573,2.893441125e-11) (1.279897049,2.609922885e-11) (1.352877724,2.427343075e-11) (1.425987949,2.304817392e-11) (1.498990558,2.221437158e-11) (1.5723985,2.162471584e-11) (1.645449268,2.119820149e-11) (1.71806931,2.08839352e-11) (1.791012779,2.06484122e-11) (1.865191359,2.047000281e-11) (1.938136706,2.033244946e-11) (2.010805917,2.02254321e-11) (2.085289251,2.014131674e-11) (2.158359837,2.007501405e-11) (2.231211654,2.002171743e-11)};
\addplot[plotblue,only marks,mark=square*,mark size=2.2pt] coordinates {(0.528827095,0.0002851673256) (0.64811851,1.566176892e-06) (0.759004495,1.607358986e-09) (0.795191876,1.175081673e-10)};
\nextgroupplot[title={Box},xlabel={},ylabel={$\gaprel(V)$}]
\addplot[plotgreen,thick,mark=*,mark size=1.2pt] coordinates {(1,5.781204667) (100,0.06863054828) (200,5.142836241e-05) (300,1.296094317e-05) (400,7.96355854e-06) (500,6.500833543e-06) (600,5.648295535e-07) (700,3.262030713e-07) (800,1.904050371e-07) (900,1.117490804e-07) (1000,6.615947369e-08) (1100,3.949345924e-08) (1200,2.377591377e-08) (1300,1.445530052e-08) (1400,8.900246666e-09) (1500,5.572290173e-09) (1600,3.563805476e-09) (1700,2.337323669e-09) (1800,1.579288209e-09) (1900,1.579288209e-09) (2000,1.440968756e-09) (2100,1.06908364e-09) (2200,8.206862749e-10) (2300,6.358902154e-10) (2400,4.96598723e-10) (2500,3.915224337e-10) (2600,3.122901672e-10) (2700,2.525278567e-10) (2800,2.074235639e-10) (2900,1.732109415e-10) (3000,1.471533617e-10)};
\addplot[plotblue,only marks,mark=square*,mark size=2.2pt] coordinates {(11,0.0002107244632) (14,1.144517074e-06) (17,2.097503556e-08) (19,1.905454556e-10)};
\nextgroupplot[title={Box},xlabel={},ylabel={}]
\addplot[plotgreen,thick,mark=*,mark size=1.2pt] coordinates {(0.001404326002,5.781204667) (0.073809978,0.06863054828) (0.146769973,5.142836241e-05) (0.219237125,1.296094317e-05) (0.292437051,7.96355854e-06) (0.365189555,6.500833543e-06) (0.437655906,5.648295535e-07) (0.509980692,3.262030713e-07) (0.582308122,1.904050371e-07) (0.654590693,1.117490804e-07) (0.727111436,6.615947369e-08) (0.799360152,3.949345924e-08) (0.871746404,2.377591377e-08) (0.944344502,1.445530052e-08) (1.016468835,8.900246666e-09) (1.088591896,5.572290173e-09) (1.161107665,3.563805476e-09) (1.233795902,2.337323669e-09) (1.306397164,1.579288209e-09) (1.378914892,1.579288209e-09) (1.451344467,1.440968756e-09) (1.523589378,1.06908364e-09) (1.596238897,8.206862749e-10) (1.668518221,6.358902154e-10) (1.741092062,4.96598723e-10) (1.813163177,3.915224337e-10) (1.885514727,3.122901672e-10) (1.95840001,2.525278567e-10) (2.030943162,2.074235639e-10) (2.103199945,1.732109415e-10) (2.175341606,1.471533617e-10)};
\addplot[plotblue,only marks,mark=square*,mark size=2.2pt] coordinates {(0.523833458,0.0002107244632) (0.640535902,1.144517074e-06) (0.743232248,2.097503556e-08) (0.823965224,1.905454556e-10)};
\nextgroupplot[title={Diffuse},xlabel={Iterations},ylabel={$\gaprel(V)$}]
\addplot[plotgreen,thick,mark=*,mark size=1.2pt] coordinates {(1,0.05714715876) (100,0.001110930649) (200,0.0001521573304) (300,4.238565046e-05) (400,1.737002794e-05) (500,9.068235381e-06) (600,5.472751689e-06) (700,3.58606506e-06) (800,2.462171831e-06) (900,1.738232058e-06) (1000,1.249105536e-06) (1100,9.085093995e-07) (1200,6.665437672e-07) (1300,4.922280882e-07) (1400,3.6536168e-07) (1500,2.723131115e-07) (1600,2.036549588e-07) (1700,1.527483056e-07) (1800,1.148537254e-07) (1900,8.655195204e-08) (2000,6.535551743e-08) (2100,4.944262451e-08) (2200,3.747162158e-08) (2300,2.844989292e-08) (2400,2.164018427e-08) (2500,1.649307514e-08) (2600,1.259793836e-08) (2700,9.647106046e-09) (2800,7.409544734e-09) (2900,5.711431092e-09) (3000,4.421762247e-09)};
\addplot[plotblue,only marks,mark=square*,mark size=2.2pt] coordinates {(11,0.001248060383) (14,9.26629318e-06) (18,1.032463866e-08) (20,4.797118859e-10)};
\nextgroupplot[title={Diffuse},xlabel={Elapsed time (s)},ylabel={}]
\addplot[plotgreen,thick,mark=*,mark size=1.2pt] coordinates {(0.001427792013,0.05714715876) (0.07953567401,0.001110930649) (0.152398581,0.0001521573304) (0.225076182,4.238565046e-05) (0.297852476,1.737002794e-05) (0.370511211,9.068235381e-06) (0.44346685,5.472751689e-06) (0.516169394,3.58606506e-06) (0.589895807,2.462171831e-06) (0.664237027,1.738232058e-06) (0.736966563,1.249105536e-06) (0.809734347,9.085093995e-07) (0.882376076,6.665437672e-07) (0.955321552,4.922280882e-07) (1.028150917,3.6536168e-07) (1.10071536,2.723131115e-07) (1.17345278,2.036549588e-07) (1.245946195,1.527483056e-07) (1.318771515,1.148537254e-07) (1.391551852,8.655195204e-08) (1.464074864,6.535551743e-08) (1.537102802,4.944262451e-08) (1.609728187,3.747162158e-08) (1.682355896,2.844989292e-08) (1.754932824,2.164018427e-08) (1.828637592,1.649307514e-08) (1.901394584,1.259793836e-08) (1.974435619,9.647106046e-09) (2.046973339,7.409544734e-09) (2.119504774,5.711431092e-09) (2.192125788,4.421762247e-09)};
\addplot[plotblue,only marks,mark=square*,mark size=2.2pt] coordinates {(0.52311086,0.001248060383) (0.638209184,9.26629318e-06) (0.787506618,1.032463866e-08) (0.856209338,4.797118859e-10)};
\end{groupplot}
\end{tikzpicture}

\pgfplotslegendfromname{legend:geocompare}
\caption{Geodesic objective comparison at $M=63$ and $N=64$.
Left: abscissa, iteration count (DR iterations for the green curves, interior-point iterations for the blue squares).
Right: abscissa, elapsed time (DR including the checkpoint evaluations of $V$; IPM median of three runs).
Ordinate in both columns: $\gaprel(V)$ of \eqref{eq:comparison-gap}.
In each row, DR is evaluated after iteration 1 and then every 100 iterations during one solve.
Each green point shows $\gaprel(V)$ for the lowest action $V$ evaluated up to that iteration.
The green segments join these recorded evaluations.
The four blue squares are final results from separate IPM solves with stopping tolerances $10^{-4}$, $10^{-6}$, $10^{-8}$, and $10^{-10}$.
The benchmark saved no intermediate action evaluations for IPM, so its final results are shown as unconnected points.
The common reference in \eqref{eq:comparison-gap} is evaluated in floating-point arithmetic.}
\label{fig:geocompare}
\end{figure}
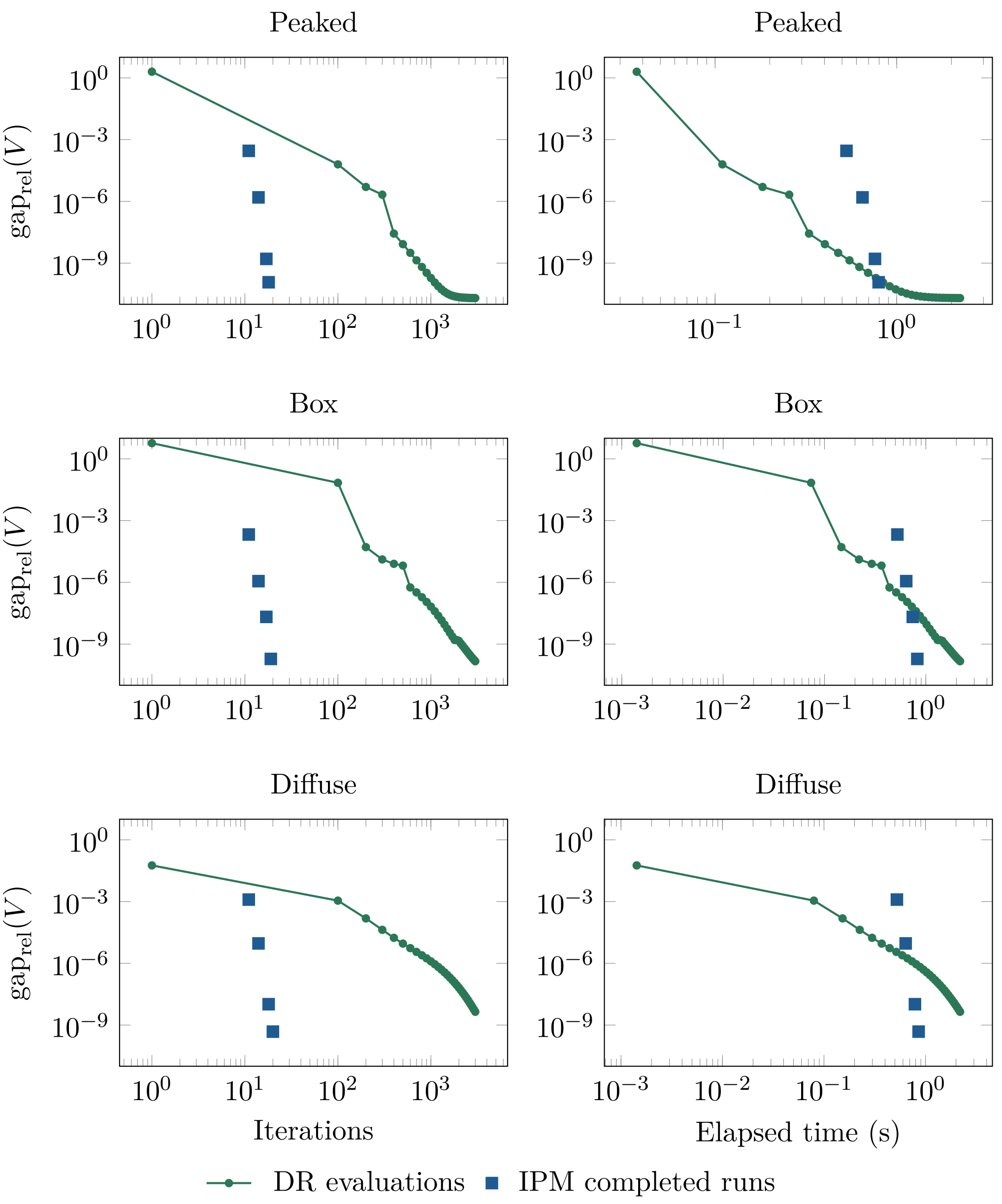

IPM uses the conic formulation in \cref{sec:geocone}.
The timings of \cref{tab:geotarget,fig:geocompare} were measured on one performance core of an Intel i9-13950HX processor, with one thread per solver; \cref{tab:geofine,tab:checked-prox} state their own settings.
DR uses $\beta=10^6$, $b=1$, and 3000 iterations from the linear density interpolation.

Clarabel computes its search directions with the \href{https://clarabel.org/stable/user_guide_linsolvers/}{QDLDL} direct $LDL^{\mathsf T}$ solver.
The four IPM runs use the stopping tolerances $10^{-4}$, $10^{-6}$, $10^{-8}$ and $10^{-10}$, each assigned to Clarabel's feasibility tolerance and to its absolute and relative primal--dual gap tolerances.

For each endpoint pair, fix a positive dual reference value $D_{\rm ref}$ from a separate IPM run with tolerance $10^{-12}$.
The dual functional is the one in \cref{app:geodual}.
For the action $V$ of a path that satisfies the geodesic constraints, define
\begin{equation}\label{eq:comparison-gap}
 \gaprel(V)=\frac{V-D_{\rm ref}}{D_{\rm ref}}.
\end{equation}
\Cref{fig:geocompare} plots $\gaprel(V)$ using the same $D_{\rm ref}$ for all methods within each endpoint pair.
The IPM stopping test uses feasibility residuals and primal--dual gaps computed from its current conic iterates, rather than the gap to this fixed reference.
With exact evaluations, $0<D_{\rm ref}\le E_{M,N}\le V$ gives
\[
 0\le\frac{V-E_{M,N}}{E_{M,N}}\le\gaprel(V).
\]
The evaluations of $V$ and $D_{\rm ref}$ use floating-point arithmetic.
Thus the plotted values are numerical gap estimates, without interval certificates.

DR evaluates the action of the initial path and at the checkpoints $\ell=1,100,200,\ldots,3000$.
At a checkpoint $\ell$, let $y_\nu^{\ell,j}$ be the density component of the projected output $Y^\ell$ at time level $j$.
The affine projection enforces conservation and unit mass in exact arithmetic, but it can leave negative density entries.
At each interior time level, set $\nu_{\rm proj}^j=\argmin_{p\in\mathcal P_M}\|p-y_\nu^{\ell,j}\|_2^2$, $1\le j<N$.
To obtain positive midpoint weights for \eqref{eq:geovalue}, set
\[
 \bar\nu^j=(1-\varepsilon_{\rm mix})\nu_{\rm proj}^j+\varepsilon_{\rm mix}\one,
 \qquad \varepsilon_{\rm mix}=10^{-12},\qquad 1\le j<N.
\]
This mixture retains unit mass and makes every interior entry at least $\varepsilon_{\rm mix}$.
Keep $\bar\nu^0=\nu^0$ and $\bar\nu^N=\nu^N$.
Set $\bar q^j=(\bar\nu^j+\bar\nu^{j+1})/2$.
Each cell $0\le j<N$ contains an interior level $\bar\nu^j$ or $\bar\nu^{j+1}$ with entries at least $\varepsilon_{\rm mix}$, so $\bar q^j\ge(\varepsilon_{\rm mix}/2)\one>0$ and \eqref{eq:geovalue} is well defined.
The mixture is used only for action evaluation; the geodesic problem still permits zero densities.
\Cref{sec:proximal-zeros} gives the conservation condition required at zero midpoint weights.

Compute $\bar m^j$ by \eqref{eq:geovalue} with $(\nu,q)=(\bar\nu,\bar q)$.
This momentum satisfies $L\bar m^j=\bar\nu^{j+1}-\bar\nu^j$.
The action $V$ of \cref{sec:geovalue} is evaluated with $(q^j,m^j)=(\bar q^j,\bar m^j)$.
The projection, mixture, and new momentum calculation apply only to the evaluation of $V$.
Each evaluated $V$ is an upper bound on $E_{M,N}$; at checkpoint $\ell$, the DR curve therefore plots $\gaprel(V_{\min}^\ell)$, where $V_{\min}^\ell$ is the smallest action evaluated on the initial path and at the checkpoints $\ell'\le\ell$.
For IPM, \cref{alg:geo-candidates} clips the final density iterate at zero and renormalizes each interior row; it then forms mixtures with the uniform density path and with the linear interpolation of the endpoints, with weights $10^{-k/2}$, $8\le k\le32$, adds the row-normalized path of the slacks of the constraints $\nu^j\ge0$ when all these slacks are positive, and keeps the candidate of least action, computed with \eqref{eq:geovalue}.
The DR path uses the single weight $\varepsilon_{\rm mix}=10^{-12}$, so the two values of $V$ come from different path constructions.

DR times cover the iterations and the evaluations of $V$ at the checkpoints, but not the evaluation of the dual value $D_G(\chi)$ of \cref{alg:geo-dr}.
IPM times are medians of three runs; they cover solver setup, factorization and iterations, but not matrix assembly or the evaluations of $V$ and $D_G(\chi)$.
Neither includes the IPM run that computes $D_{\rm ref}$.

\begin{table}[htbp]
\centering\small
\begin{tabular}{@{}lrr@{}}
\toprule
Pair & DR seconds & IPM seconds\\
\midrule
Peaked & 0.40 & 0.76\\
Box & 1.02 & 0.82\\
Diffuse & 1.97 & 0.86\\
\bottomrule
\end{tabular}
\caption{Elapsed time to reach $\gaprel(V)\le10^{-8}$ in the geodesic benchmark at $M=63$, $N=64$. DR ($\beta=10^6$, $b=1$): time at the first checkpoint (iteration 1 or a multiple of 100) with $\gaprel(V)\le10^{-8}$, an upper estimate of the first crossing; iterations and checkpoint evaluations of $V$ included. IPM: the run with the largest stopping tolerance among $10^{-4},10^{-6},10^{-8},10^{-10}$ whose final $\gaprel(V)\le10^{-8}$ ($10^{-8}$ for the peaked pair, $10^{-10}$ for the box and diffuse pairs), with the median time of three repetitions; setup and factorization included, matrix assembly and evaluation of $V$ excluded. The computation of $D_{\rm ref}$ is excluded for both methods.}
\label{tab:geotarget}
\end{table}

At the threshold $\gaprel(V)\le10^{-8}$, the ratio of DR time to IPM time is $0.53$ for the peaked pair, $1.24$ for the box pair and $2.29$ for the diffuse pair (\cref{tab:geotarget}).

\paragraph{Finer spatial grids.}
To measure how the cost of each method grows with $M$, \cref{tab:geofine} reports the peaked and diffuse pairs at $M\in\{127,255\}$, and IPM for the peaked pair at $M=511$, with $N=64$.
DR uses $\beta=10^6$, $b=1$ and 3000 iterations; IPM uses QDLDL and stopping tolerance $10^{-10}$.

For these grids, \eqref{eq:comparison-gap} uses in place of $D_{\rm ref}$ the reference $D_{\rm ref}^{\rm fine}=D_G(\chi)$ of \eqref{eq:geodual}, with $\chi$ obtained by the multiplier conversion of \cref{alg:geo-ipm} from the final iterate of the listed IPM solve; the same value is used for both methods on the same endpoint pair and grid.
The IPM entry of $\gaprel(V)$ is therefore the relative primal--dual gap of that solve, and the difference between the DR and IPM entries equals $(V_{\rm DR}-V_{\rm IPM})/D_{\rm ref}^{\rm fine}$.
Since the runs are not stopped at a common accuracy, \cref{tab:geofine} gives for each run its final $\gaprel(V)$ and its time.
DR retains its lowest evaluated action over 3000 iterations; IPM reports its final action.
The DR times include checkpoint action evaluations.
The IPM times include solver setup and factorization but exclude matrix assembly and subsequent action and dual evaluations.

\begin{table}[htbp]
\centering\small
\begin{tabular}{@{}lrlrrr@{}}
\toprule
Pair & $M$ & Method & Iterations & Seconds & $\gaprel(V)$\\
\midrule
Peaked & 127 & DR & 3000 & 3.37 & $5.22\cdot10^{-9}$\\
 &  & IPM & 18 & 4.73 & $5.19\cdot10^{-9}$\\
\addlinespace
Peaked & 255 & DR & 3000 & 5.04 & $8.47\cdot10^{-10}$\\
 &  & IPM & 20 & 35.60 & $5.73\cdot10^{-10}$\\
\addlinespace
Peaked & 511 & IPM & 21 & 286.62 & $3.46\cdot10^{-9}$\\
\addlinespace
Diffuse & 127 & DR & 3000 & 3.34 & $1.80\cdot10^{-7}$\\
 &  & IPM & 20 & 5.25 & $1.73\cdot10^{-8}$\\
\addlinespace
Diffuse & 255 & DR & 3000 & 5.10 & $4.28\cdot10^{-7}$\\
 &  & IPM & 23 & 41.68 & $2.08\cdot10^{-9}$\\
\bottomrule
\end{tabular}
\caption{Stored geodesic computations on finer spatial grids, with $N=64$.
DR uses four threads and IPM one.
Each DR entry comes from one solve.
IPM times are medians of three solves at $M=127,255$ and one recorded time at $M=511$.
Iterations are DR splitting iterations or IPM iterations.
DR reports the lowest evaluated action over 3000 iterations, and its time includes the checkpoint action evaluations.
$\gaprel(V)$ is evaluated in floating-point arithmetic with $D_{\rm ref}^{\rm fine}$, the dual value $D_G(\chi)$ of the listed IPM solve on the same pair and grid, in place of $D_{\rm ref}$.}
\label{tab:geofine}
\end{table}

For the peaked pair, the DR and IPM values of $\gaprel(V)$ differ by at most $3\cdot10^{-10}$, and DR takes less time ($3.37$ s versus $4.73$ s at $M=127$, $5.04$ s versus $35.60$ s at $M=255$).
For the diffuse pair, the IPM value of $\gaprel(V)$ is smaller than the DR value by a factor $10.4$ at $M=127$ and $206$ at $M=255$, and the DR value increases with $M$ at the fixed budget of 3000 iterations.
When $M$ doubles, the IPM time grows by a factor $7.5$ to $8.1$, whereas the DR time grows by a factor $1.5$ (\cref{tab:geofine}).
The differing thread counts and attained gaps prevent a comparison at equal computational resources and accuracy.
DR uses Fourier and sine transforms with cost $O(MN(\log M+\log N))$ per iteration; its total cost also depends on the required iteration count.
\FloatBarrier

\begin{remark}[Residual-checked proximal solver]\label{rem:checked-prox}
\Cref{alg:geo-prox} performs 18 Newton iterations without a residual test, and its zero outputs result from the clipping at zero, not from the test $2bq_{\rm in}+v_{\rm in}^2\le0$.
To check the DR outputs obtained with it, we also ran DR for the three endpoint pairs at $(M,N)=(63,64)$ with the residual-checked solver of \cref{alg:geo-prox-checked}.
This solver applies the zero test, keeps each Newton iterate in a bracket that contains the root, and accepts a positive weight $q$ only if the scalar residual $g_0$ of \cref{app:prox-error} satisfies $|g_0(q)|\le2\cdot10^{-14}(1+\max\{|q_{\rm in}|,|v_{\rm in}|\})$; it reports a failure if an entry fails this test after 160 iterations.
After 3000 DR iterations, for each pair, the final interior densities of the two runs have the same zero entries and differ by at most $1.64\cdot10^{-12}$ in maximum norm, and the DR loop with \cref{alg:geo-prox-checked} takes $1.3$ to $2.0$ times as long as with \cref{alg:geo-prox} (\cref{tab:checked-prox}).
This comparison concerns these three runs; it does not bound the error of each proximal call in the other runs.
Residual-checked solvers are also used in the tests of \cref{tab:common-entropic,tab:spectralchecks} and in the notebooks (\cref{app:notebooks}).
\end{remark}

\begin{table}[htbp]
\centering\small
\begin{tabular}{@{}lrrrrrr@{}}
\toprule
Pair & Zeros & Changed & $\|\Delta\nu\|_\infty$ & $\max\|g_0\|_2$ & Alg.~\ref{alg:geo-prox} (s) & Alg.~\ref{alg:geo-prox-checked} (s)\\
\midrule
Peaked & 61 & 0 & $5.76\cdot10^{-14}$ & $1.00\cdot10^{-12}$ & 2.31 & 3.82\\
Box & 59 & 0 & $9.79\cdot10^{-14}$ & $7.97\cdot10^{-13}$ & 2.34 & 4.78\\
Diffuse & 26 & 0 & $1.64\cdot10^{-12}$ & $3.15\cdot10^{-13}$ & 2.33 & 3.01\\
\bottomrule
\end{tabular}
\caption{DR with the proximal solver of \cref{alg:geo-prox} (18 Newton iterations) and with that of \cref{alg:geo-prox-checked} (residual-checked), after 3000 DR iterations, with $M=63$, $N=64$, $\beta=10^6$, and $b=1$. Both runs start from the same linear path and use the same Fourier operator. Zeros: number of entries equal to zero in the interior density component of the final proximal output (identical for both solvers). Changed: number of interior density entries that are zero in one run and positive in the other. $\|\Delta\nu\|_\infty$: largest absolute difference between the interior density components of the two final proximal outputs. $\max\|g_0\|_2$: largest Euclidean norm of the scalar residual $g_0$ of \cref{app:prox-error} ($\omega=0$ for geodesics) over all calls of \cref{alg:geo-prox-checked}. Alg.~\ref{alg:geo-prox}, Alg.~\ref{alg:geo-prox-checked}: time of one local execution of the DR loop with each solver, including residual recording but not the final action evaluation, on logical CPU 0 with one thread. No call of \cref{alg:geo-prox-checked} raised its failure. The comparison concerns these finite runs, not the zero set of the exact minimizer.}
\label{tab:checked-prox}
\end{table}

\subsection{Effect of the momentum weight \texorpdfstring{$\beta$}{beta}}\label{sec:betasweep}
\Cref{sec:dr-choices} shows that, for the uniform geodesic at $b=1$, increasing $\beta$ decreases the asymptotic factor $\spr(T)$ of \eqref{eq:rate} until it reaches $1/2$ and does not increase it afterwards.
We test whether the three endpoint pairs follow this behavior, and which tested $\beta$ gives the smallest final gap after 3000 iterations, by varying $\beta$ in the DR norm \eqref{eq:metric}, keeping $b=1$, $M=63$, and $N=64$ fixed.
Each of the three endpoint pairs uses the same linear density interpolation for every run.
Each run performs 3000 iterations on one computational thread.
Since $b=\gamma\dx/(h\beta)$, the proximal parameter $\gamma$ changes with $\beta$.
The minimized action and the reference $D_{\rm ref}$ do not change.

\begin{figure}[p]
\centering
\begin{tikzpicture}
\begin{groupplot}[group style={group size=2 by 3,horizontal sep=1.25cm,vertical sep=1.9cm},width=.43\textwidth,height=4.8cm,ymode=log,ymin=1e-12,ymax=10,ytick={1,1e-4,1e-8,1e-12},legend columns=4,legend style={draw=none,font=\scriptsize}]
\nextgroupplot[title={(a) Peaked},xlabel={},ylabel={$\gaprel(V)$},xmin=0,xmax=3000,xtick={0,1000,2000,3000},legend to name=legend:betasweep]
\addplot[plotpurple,dashed,thick] coordinates {(0,2.904330866) (1,2.873545625) (100,0.6191765808) (200,0.418321223) (300,0.3202819792) (400,0.2588372976) (500,0.215620435) (600,0.1831775934) (700,0.1577610923) (800,0.1372412114) (900,0.1202962986) (1000,0.1060609335) (1100,0.09393607218) (1200,0.08349902213) (1300,0.07443616032) (1400,0.06650777355) (1500,0.05953233908) (1600,0.05336572817) (1700,0.04789092267) (1800,0.04301241211) (1900,0.0386521329) (2000,0.03474661625) (2100,0.03124261075) (2200,0.0280936517) (2300,0.02525995024) (2400,0.02270732198) (2500,0.0204062233) (2600,0.01834602455) (2700,0.0164732119) (2800,0.01478974102) (2900,0.01325177428) (3000,0.01187661864)};
\addlegendentry{$\beta=10^{4}$}
\addplot[plotgreen,thick] coordinates {(0,2.904330866) (1,1.991249049) (100,6.283928761e-05) (200,5.028030992e-06) (300,2.114965888e-06) (400,2.727584622e-08) (500,8.350806615e-09) (600,3.158654167e-09) (700,1.366304017e-09) (800,6.520559032e-10) (900,3.389598752e-10) (1000,1.898931802e-10) (1100,1.148771237e-10) (1200,7.520648243e-11) (1300,5.328743351e-11) (1400,4.075382087e-11) (1500,3.340030066e-11) (1600,2.893441125e-11) (1700,2.609922885e-11) (1800,2.427343075e-11) (1900,2.304817392e-11) (2000,2.221437158e-11) (2100,2.162471584e-11) (2200,2.119820149e-11) (2300,2.08839352e-11) (2400,2.06484122e-11) (2500,2.047000281e-11) (2600,2.033244946e-11) (2700,2.02254321e-11) (2800,2.014131674e-11) (2900,2.007501405e-11) (3000,2.002171743e-11)};
\addlegendentry{$\beta=10^{6}$}
\addplot[plotblue,dashdotted,thick] coordinates {(0,2.904330866) (1,2.904330866) (100,2.904330866) (200,2.904330866) (300,2.904330866) (400,2.904330866) (500,2.904330866) (600,2.904330866) (700,2.904330866) (800,2.904330866) (900,2.904330866) (1000,2.904330866) (1100,0.002837560487) (1200,0.0002808112782) (1300,0.0002808112782) (1400,0.0002808112782) (1500,0.0002808112782) (1600,7.915580625e-05) (1700,7.449600914e-05) (1800,7.249883868e-05) (1900,6.590237919e-05) (2000,3.64551061e-05) (2100,3.64551061e-05) (2200,3.00518868e-05) (2300,2.594215694e-05) (2400,1.766737611e-05) (2500,1.550789135e-05) (2600,1.389313376e-05) (2700,1.249964073e-05) (2800,1.200338315e-05) (2900,1.200338315e-05) (3000,1.200338315e-05)};
\addlegendentry{$\beta=10^{8}$}
\addplot[plotred,dotted,thick] coordinates {(0,2.904330866) (1,2.904330866) (100,2.904330866) (200,2.904330866) (300,2.904330866) (400,2.904330866) (500,2.904330866) (600,2.904330866) (700,2.904330866) (800,2.904330866) (900,2.904330866) (1000,2.904330866) (1100,2.904330866) (1200,2.904330866) (1300,2.904330866) (1400,2.904330866) (1500,2.904330866) (1600,2.904330866) (1700,2.904330866) (1800,2.904330866) (1900,2.904330866) (2000,2.904330866) (2100,2.904330866) (2200,2.904330866) (2300,2.904330866) (2400,2.904330866) (2500,2.904330866) (2600,2.904330866) (2700,2.904330866) (2800,2.904330866) (2900,2.904330866) (3000,2.904330866)};
\addlegendentry{$\beta=10^{12}$}
\nextgroupplot[title={(b) After 3000 iterations},xlabel={},xmode=log,xmin=1,xmax=1e12,xtick={1,1e4,1e8,1e12}]
\addplot[black,mark=*,mark size=1.5pt,thick] coordinates {(1,2.647585309) (100,1.026754015) (10000,0.01187661864) (100000,6.775540013e-08) (300000,2.748970521e-10) (1000000,2.002171743e-11) (3000000,1.976032365e-11) (100000000,1.200338315e-05) (1e+10,2.904330866) (1e+12,2.904330866)};
\nextgroupplot[title={(c) Box},xlabel={},ylabel={$\gaprel(V)$},xmin=0,xmax=3000,xtick={0,1000,2000,3000}]
\addplot[plotpurple,dashed,thick] coordinates {(0,5.781204667) (1,5.510228022) (100,0.1614883925) (200,0.07857180289) (300,0.0459264517) (400,0.02875212166) (500,0.01893537859) (600,0.01329737214) (700,0.009435124171) (800,0.007037383338) (900,0.004997126874) (1000,0.003664278366) (1100,0.002739381801) (1200,0.001960244738) (1300,0.001497310578) (1400,0.001219348887) (1500,0.001005614147) (1600,0.0007869354955) (1700,0.0006255855281) (1800,0.0006255855281) (1900,0.0004977847421) (2000,0.0003999011433) (2100,0.0003333895025) (2200,0.0002682322708) (2300,0.0002191897678) (2400,0.0001826134758) (2500,0.0001530094849) (2600,0.0001292245449) (2700,0.0001101015638) (2800,9.452146432e-05) (2900,8.188837062e-05) (3000,7.160958804e-05)};
\addplot[plotgreen,thick] coordinates {(0,5.781204667) (1,5.781204667) (100,0.06863054828) (200,5.142836241e-05) (300,1.296094317e-05) (400,7.96355854e-06) (500,6.500833543e-06) (600,5.648295535e-07) (700,3.262030713e-07) (800,1.904050371e-07) (900,1.117490804e-07) (1000,6.615947369e-08) (1100,3.949345924e-08) (1200,2.377591377e-08) (1300,1.445530052e-08) (1400,8.900246666e-09) (1500,5.572290173e-09) (1600,3.563805476e-09) (1700,2.337323669e-09) (1800,1.579288209e-09) (1900,1.579288209e-09) (2000,1.440968756e-09) (2100,1.06908364e-09) (2200,8.206862749e-10) (2300,6.358902154e-10) (2400,4.96598723e-10) (2500,3.915224337e-10) (2600,3.122901672e-10) (2700,2.525278567e-10) (2800,2.074235639e-10) (2900,1.732109415e-10) (3000,1.471533617e-10)};
\addplot[plotblue,dashdotted,thick] coordinates {(0,5.781204667) (1,5.781204667) (100,5.781204667) (200,5.781204667) (300,5.781204667) (400,5.781204667) (500,5.781204667) (600,5.781204667) (700,5.781204667) (800,5.781204667) (900,5.781204667) (1000,5.781204667) (1100,5.781204667) (1200,5.781204667) (1300,5.781204667) (1400,5.781204667) (1500,5.781204667) (1600,5.781204667) (1700,5.781204667) (1800,5.781204667) (1900,5.781204667) (2000,5.781204667) (2100,5.781204667) (2200,5.781204667) (2300,5.781204667) (2400,5.781204667) (2500,5.781204667) (2600,5.781204667) (2700,5.781204667) (2800,5.781204667) (2900,5.781204667) (3000,5.781204667)};
\addplot[plotred,dotted,thick] coordinates {(0,5.781204667) (1,5.781204667) (100,5.781204667) (200,5.781204667) (300,5.781204667) (400,5.781204667) (500,5.781204667) (600,5.781204667) (700,5.781204667) (800,5.781204667) (900,5.781204667) (1000,5.781204667) (1100,5.781204667) (1200,5.781204667) (1300,5.781204667) (1400,5.781204667) (1500,5.781204667) (1600,5.781204667) (1700,5.781204667) (1800,5.781204667) (1900,5.781204667) (2000,5.781204667) (2100,5.781204667) (2200,5.781204667) (2300,5.781204667) (2400,5.781204667) (2500,5.781204667) (2600,5.781204667) (2700,5.781204667) (2800,5.781204667) (2900,5.781204667) (3000,5.781204667)};
\nextgroupplot[title={(d) After 3000 iterations},xlabel={},xmode=log,xmin=1,xmax=1e12,xtick={1,1e4,1e8,1e12}]
\addplot[black,mark=*,mark size=1.5pt,thick] coordinates {(1,2.710916924) (100,0.4182661918) (10000,7.160958804e-05) (100000,8.619593722e-07) (300000,4.274532707e-08) (1000000,1.471533617e-10) (3000000,4.969237943e-09) (100000000,5.781204667) (1e+10,5.781204667) (1e+12,5.781204667)};
\nextgroupplot[title={(e) Diffuse},xlabel={Iterations},ylabel={$\gaprel(V)$},xmin=0,xmax=3000,xtick={0,1000,2000,3000}]
\addplot[plotpurple,dashed,thick] coordinates {(0,0.0577575311) (1,0.05774652897) (100,0.04290446733) (200,0.03604969056) (300,0.03167806736) (400,0.02851134395) (500,0.02606028046) (600,0.02407671369) (700,0.02241664046) (800,0.02099083464) (900,0.0197408618) (1000,0.018627083) (1100,0.01762225858) (1200,0.01670722719) (1300,0.01586778715) (1400,0.01509460442) (1500,0.01437897705) (1600,0.01371397948) (1700,0.01309338069) (1800,0.01251259244) (1900,0.01196738975) (2000,0.01145468862) (2100,0.01097134421) (2200,0.01051488916) (2300,0.01008323327) (2400,0.009674552644) (2500,0.009287110236) (2600,0.008919418496) (2700,0.008570150553) (2800,0.008238111435) (2900,0.007922217515) (3000,0.007621480656)};
\addplot[plotgreen,thick] coordinates {(0,0.0577575311) (1,0.05714715876) (100,0.001110930649) (200,0.0001521573304) (300,4.238565046e-05) (400,1.737002794e-05) (500,9.068235381e-06) (600,5.472751689e-06) (700,3.58606506e-06) (800,2.462171831e-06) (900,1.738232058e-06) (1000,1.249105536e-06) (1100,9.085093995e-07) (1200,6.665437672e-07) (1300,4.922280882e-07) (1400,3.6536168e-07) (1500,2.723131115e-07) (1600,2.036549588e-07) (1700,1.527483056e-07) (1800,1.148537254e-07) (1900,8.655195204e-08) (2000,6.535551743e-08) (2100,4.944262451e-08) (2200,3.747162158e-08) (2300,2.844989292e-08) (2400,2.164018427e-08) (2500,1.649307514e-08) (2600,1.259793836e-08) (2700,9.647106046e-09) (2800,7.409544734e-09) (2900,5.711431092e-09) (3000,4.421762247e-09)};
\addplot[plotblue,dashdotted,thick] coordinates {(0,0.0577575311) (1,0.0433041126) (100,0.00312885895) (200,0.0005920087873) (300,0.0002195631386) (400,0.0001040083589) (500,4.519367014e-05) (600,2.776200782e-05) (700,1.605105837e-05) (800,1.038161763e-05) (900,6.914123649e-06) (1000,4.819037086e-06) (1100,3.46605522e-06) (1200,2.581467301e-06) (1300,1.968201781e-06) (1400,1.527648196e-06) (1500,1.198468055e-06) (1600,9.485944769e-07) (1700,7.57342057e-07) (1800,6.085167044e-07) (1900,4.560179073e-07) (2000,3.709992847e-07) (2100,2.436542396e-07) (2200,1.832786583e-07) (2300,1.410824735e-07) (2400,1.096510946e-07) (2500,8.593233731e-08) (2600,6.611271048e-08) (2700,5.140594876e-08) (2800,4.04726568e-08) (2900,3.221950605e-08) (3000,2.589327973e-08)};
\addplot[plotred,dotted,thick] coordinates {(0,0.0577575311) (1,0.0577575311) (100,0.0577575311) (200,0.0577575311) (300,0.0577575311) (400,0.0577575311) (500,0.0577575311) (600,0.0577575311) (700,0.0577575311) (800,0.0577575311) (900,0.0577575311) (1000,0.0577575311) (1100,0.0577575311) (1200,0.0577575311) (1300,0.0577575311) (1400,0.0577575311) (1500,0.0577575311) (1600,0.0577575311) (1700,0.0577575311) (1800,0.0577575311) (1900,0.0577575311) (2000,0.0577575311) (2100,0.0577575311) (2200,0.0577575311) (2300,0.0577575311) (2400,0.0577575311) (2500,0.0577575311) (2600,0.0577575311) (2700,0.0577575311) (2800,0.0577575311) (2900,0.0577575311) (3000,0.0577575311)};
\nextgroupplot[title={(f) After 3000 iterations},xlabel={$\beta$},xmode=log,xmin=1,xmax=1e12,xtick={1,1e4,1e8,1e12}]
\addplot[black,mark=*,mark size=1.5pt,thick] coordinates {(1,0.05759691538) (100,0.05044356994) (10000,0.007621480656) (100000,4.131514186e-05) (300000,1.721323979e-06) (1000000,4.421762247e-09) (3000000,3.153123069e-10) (100000000,2.589327973e-08) (1e+10,0.01222639452) (1e+12,0.0577575311)};
\end{groupplot}
\end{tikzpicture}
\pgfplotslegendfromname{legend:betasweep}
\caption{Effect of the momentum weight $\beta$ on DR, with $b=1$, $M=63$, and $N=64$.
Left: $\gaprel(V)$ of the running minimum of $V$ versus the DR iteration, for $\beta\in\{10^4,10^6,10^8,10^{12}\}$, starting from the same linear density interpolation.
Each curve joins the relative gaps of the lowest action obtained so far, evaluated at iteration 0, at iteration 1 and every 100 iterations.
Right: the same quantity after 3000 iterations versus $\beta$, at $\beta\in\{1,10^2,10^4,10^5,3\cdot10^5,10^6,3\cdot10^6,10^8,10^{10},10^{12}\}$; segments only connect these tested values.
Each row uses its own common IPM dual reference, independent of $\beta$.
The ordinate is $\gaprel(V)$ of \eqref{eq:comparison-gap}, evaluated in floating-point arithmetic.}
\label{fig:betasweep}
\end{figure}
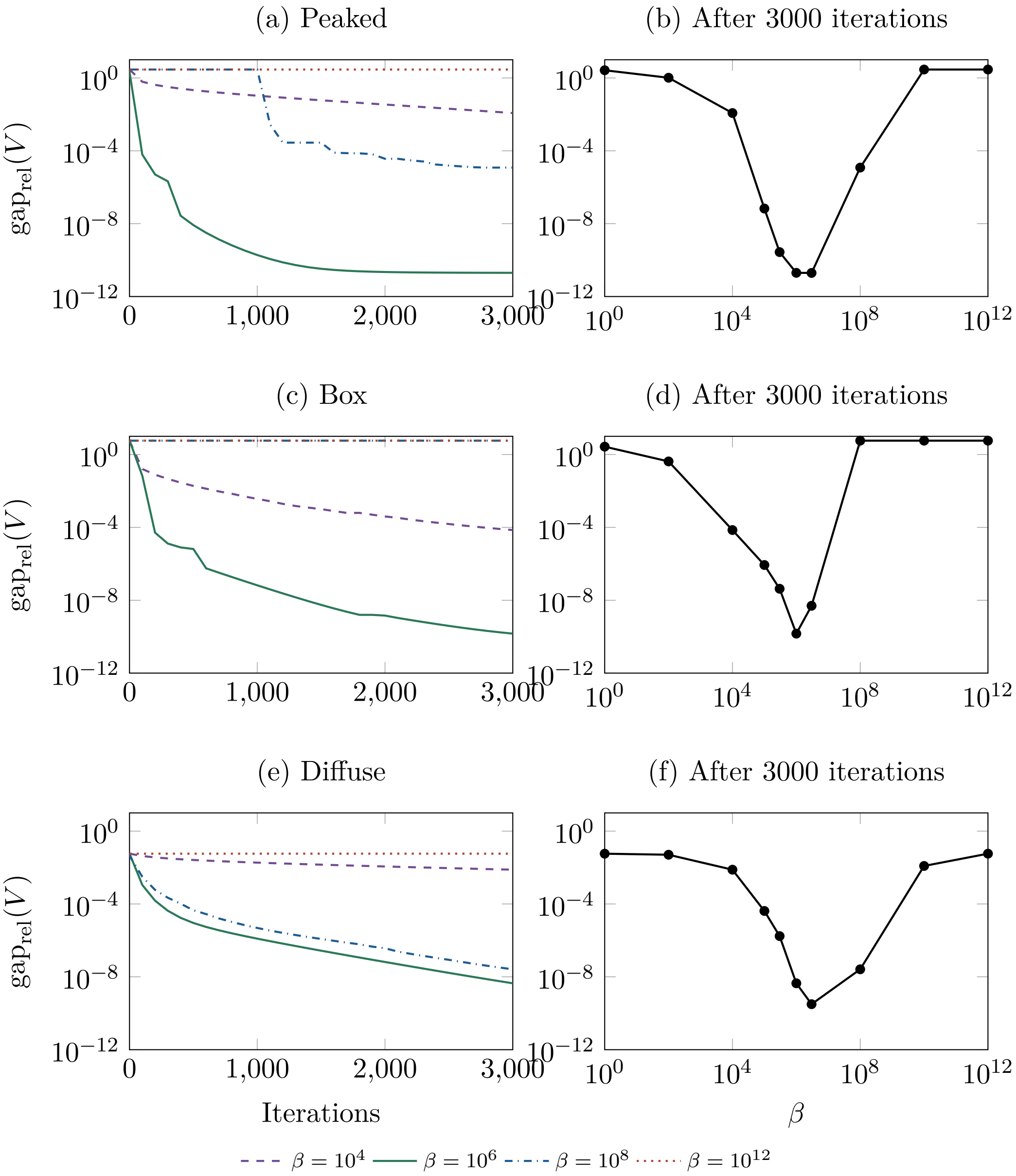

\Cref{fig:betasweep} compares the relative objective gaps from \eqref{eq:comparison-gap}.
This gap provides a common accuracy measure for different values of $\beta$; the norm defining the DR residual itself depends on $\beta$.
After 3000 iterations, the smallest final gaps occur at $\beta=10^6$ or $\beta=3\cdot10^6$ for all three pairs; every tested $\beta\le10^5$ or $\beta\ge10^8$ leaves a final gap above $10^{-8}$ (\cref{fig:betasweep}, right column).
For all three pairs, the final gap increases from $\beta=3\cdot10^6$ to $\beta=10^8$, so after 3000 iterations it does not follow the monotone dependence of $\spr(T)$ on $\beta$ found for the uniform geodesic.
For $\beta\in\{10^6,3\cdot10^6\}$, \cref{tab:betatarget} gives the first checkpoint with $\gaprel(V)\le10^{-8}$ and its time.
Among these two values, the one that reaches this threshold at the earlier checkpoint depends on the endpoint pair.
The final gaps of the peaked pair at $\beta=10^6$ and $3\cdot10^6$, $2.00\cdot10^{-11}$ and $1.98\cdot10^{-11}$, are binary64 evaluations without an error estimate; this test therefore does not determine which of the two actions is smaller.

\begin{table}[htbp]
\centering\small
\begin{tabular}{@{}llrr@{}}
\toprule
Pair & $\beta$ & First checkpoint & Seconds\\
\midrule
Peaked & $10^6$ & 500 & 0.41\\
Peaked & $3\cdot10^6$ & 700 & 0.57\\
Box & $10^6$ & 1400 & 1.07\\
Box & $3\cdot10^6$ & 2800 & 2.14\\
Diffuse & $10^6$ & 2700 & 2.04\\
Diffuse & $3\cdot10^6$ & 900 & 0.69\\
\bottomrule
\end{tabular}
\caption{First DR checkpoint with $\gaprel(V)\le10^{-8}$ in the parameter sweep.
Each entry is one run, with $b=1$, $M=63$, and $N=64$.
Times include iteration and checkpoint action evaluation, but exclude initialization and the reference IPM solve.
Checkpoints occur after iteration 1 and every 100 iterations, so the threshold may have been crossed earlier.
These runs are separate from those of \cref{tab:geotarget}; for $\beta=10^6$ the times differ from it by at most 5\%.}
\label{tab:betatarget}
\end{table}

\subsection{Effect of the local proximal scale \texorpdfstring{$b$}{b}}\label{sec:drtests}

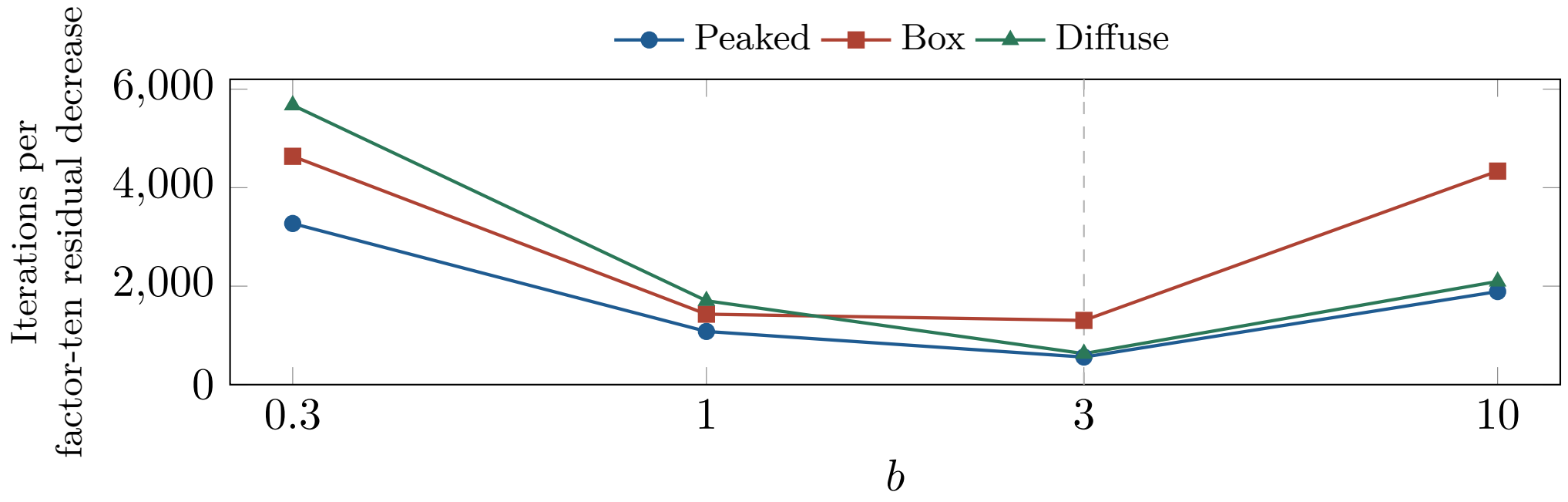
\begin{figure}[htbp]
\centering
\begin{tikzpicture}
\begin{axis}[width=13cm,height=4.2cm,xmode=log,xmin=.25,xmax=12,xtick={.3,1,3,10},xticklabels={$0.3$,$1$,$3$,$10$},xlabel={$b$},ymin=0,ymax=6200,ytick={0,2000,4000,6000},ylabel={Iterations per\\factor-ten residual decrease},ylabel style={align=center,font=\footnotesize},legend style={at={(.5,1.02)},anchor=south,legend columns=3,font=\footnotesize,draw=none}]
\addplot[black!30,dashed,forget plot] coordinates {(3,0) (3,6200)};
\addplot[plotblue,thick,mark=*,mark size=1.7pt] coordinates {(0.3,3272) (1,1080) (3,560) (10,1890)};
\addlegendentry{Peaked}
\addplot[plotred,thick,mark=square*,mark size=1.7pt] coordinates {(0.3,4637) (1,1431) (3,1304) (10,4336)};
\addlegendentry{Box}
\addplot[plotgreen,thick,mark=triangle*,mark size=1.7pt] coordinates {(0.3,5674) (1,1703) (3,630) (10,2096)};
\addlegendentry{Diffuse}
\end{axis}
\end{tikzpicture}
\caption{DR iterations per factor-ten decrease of $r_\ell$ (data of \cref{tab:drstep}) for the three endpoint pairs, at $\beta=10^6$, $M=63$, $N=64$.
Each marker is computed over $10^{-5}\ge r_\ell\ge10^{-9}$ as in \cref{tab:drstep}; segments connect the tested values of $b$, and the dashed vertical line marks $b=3$, which gives the smallest count for all three pairs.
These counts over a finite interval are not asymptotic convergence factors.}
\label{fig:drbcomparison}
\end{figure}

We vary $b$ at fixed $\beta=10^6$, grid, and endpoints.
The DR residual $r_\ell$ is recorded at every iteration.
For an interval $r_{\rm hi}\ge r_\ell\ge r_{\rm lo}$, let $\ell_{\rm hi}$ and $\ell_{\rm lo}$ be the first iterations with $r_\ell\le r_{\rm hi}$ and $r_\ell\le r_{\rm lo}$.
The number of iterations per factor-ten decrease of $r_\ell$ over this interval is $(\ell_{\rm lo}-\ell_{\rm hi})/\log_{10}(r_{\rm hi}/r_{\rm lo})$, rounded to an integer.
The linearization \eqref{eq:T} describes the DR error only near the fixed point, so \cref{tab:drstep,fig:drbcomparison} compare these counts over the late interval $10^{-5}\ge r_\ell\ge10^{-9}$.
Among the tested values, $b=3$ gives the smallest count for each pair.
For the peaked pair at $b=1$, the count is 87 over $10^{-2}\ge r_\ell\ge10^{-4}$ and 1629 over $10^{-9}\ge r_\ell\ge10^{-12}$ (\cref{tab:drstep}); a count computed on an early interval therefore underestimates the iterations needed to reach a small residual.
The ranking of $b=1$ before $b=3$ in \cref{app:drillustrations} concerns one Fourier--sine mode of the uniform path after twenty iterations and does not apply to these nonuniform paths.

\begin{table}[htbp]
\centering\small
\begin{tabular}{@{}llrrrr@{}}
\toprule
Pair & Residual interval & $b=0.3$ & $b=1$ & $b=3$ & $b=10$\\
\midrule
Peaked & $10^{-2}\ge r_\ell\ge10^{-4}$ & -- & 87 & -- & --\\
Peaked & $10^{-5}\ge r_\ell\ge10^{-9}$ & 3272 & 1080 & 560 & 1890\\
Peaked & $10^{-9}\ge r_\ell\ge10^{-12}$ & -- & 1629 & -- & --\\
Box & $10^{-5}\ge r_\ell\ge10^{-9}$ & 4637 & 1431 & 1304 & 4336\\
Diffuse & $10^{-5}\ge r_\ell\ge10^{-9}$ & 5674 & 1703 & 630 & 2096\\
\bottomrule
\end{tabular}
\caption{DR iterations per factor-ten decrease of $r_\ell$ over the stated residual intervals, at $\beta=10^6$, $M=63$, $N=64$. Each entry is $(\ell_{\rm lo}-\ell_{\rm hi})/\log_{10}(r_{\rm hi}/r_{\rm lo})$, rounded to an integer, where $r_\ell$ is recorded at every iteration and $\ell_{\rm hi}$, $\ell_{\rm lo}$ are the first iterations with $r_\ell\le r_{\rm hi}$ and $r_\ell\le r_{\rm lo}$ (\cref{sec:drtests}). The first and third rows show the dependence on the interval for the peaked pair at $b=1$; a dash marks a combination not computed. These counts over finite intervals are not asymptotic convergence factors.}
\label{tab:drstep}
\end{table}

A separate small-grid test compares the measured residual decay with the linearized DR iteration in \cref{app:numericalspectrum}.

\FloatBarrier

\subsection{Interval action bounds and DR--IPM agreement for geodesics}\label{sec:ipmgeovalidation}
For each IPM run we compute two constraint defects of the raw solver arrays and enclose by interval arithmetic the action of the IPM candidate and a dual lower bound; at $(M,N)=(63,64)$ we also compute the discrete $L^1$ distance between this candidate and the DR proximal density on the same grid.
The interval action bounds give an upper bound on the relative optimization error of the candidate action; the $L^1$ distance measures DR--IPM agreement.

For each saved IPM solution, let $\nu_{\rm raw}$ be its density and $m_{\rm raw}^j=hw^j$ its integrated momentum, where $w^j$ is the momentum-rate variable of \cref{alg:geo-ipm}.
Let $\nu_{\rm IPM}$ be the positive unit-mass density path that \cref{alg:geo-candidates} selects from the raw IPM output (clipped path, positive mixture, or slack path); we call it the IPM candidate.
We evaluate the constraint defects
\[
 \delta_{\rm mass}=\max_{0\le j\le N}|\avg{\nu_{\rm raw}^j}-1|,\qquad
 \delta_{\rm cons}=\max_{0\le j<N,\,0\le i<M}|(Lm_{\rm raw}^j)_i-\nu_{{\rm raw},i}^{j+1}+\nu_{{\rm raw},i}^{j}|.
\]
They measure how accurately the raw arrays satisfy unit mass and the conservation equations; the interval action bounds do not use them.
The interval action bounds use the IPM candidate, normalized to unit mass, and the momenta in \eqref{eq:geovalue}.

Each stored binary64 density array $\nu_{\rm IPM}^j$, $0\le j\le N$, is treated as exact real data and replaced by $\bar\nu^j=\nu_{\rm IPM}^j/\avg{\nu_{\rm IPM}^j}$, which has unit mass in exact arithmetic (\cref{app:geointerval}).
In the display below, $E_{M,N}$ is the minimum of \eqref{eq:geodesic} for the endpoints $\bar\nu^0,\bar\nu^N$.
The operator has the exact Fourier symbol $\widehat L(k)=-4\pi^2k^2$.
Let $V$ be the action of $\bar\nu$ with the momenta \eqref{eq:geovalue}, and let $D_G(\chi)$ be the dual value \eqref{eq:geodual} at $\chi^j=-\C\psi^j$, where $\psi^j$ is the IPM equality multiplier of time cell $j$ with its spatial mean subtracted.
Interval arithmetic with 40 decimal digits and outward rounding encloses $V$ and $D_G(\chi)$; $V_{\rm up}$ is the upper endpoint of the enclosure of $V$, and $D_{\rm lo}$ is the lower endpoint of the enclosure of $D_G(\chi)$.
Since $D_G(\chi)\le E_{M,N}\le V$ (\cref{sec:geovalue}), we obtain the interval action bounds and define their relative gap:
\[
 D_{\rm lo}\le E_{M,N}\le V_{\rm up},\qquad
 \gaprel^{\rm up}=\frac{V_{\rm up}-D_{\rm lo}}{D_{\rm lo}}.
\]
All computed values satisfy $D_{\rm lo}>0$.
Hence the relative optimization error of the IPM candidate satisfies $0\le(V-E_{M,N})/E_{M,N}\le\gaprel^{\rm up}$, as in \eqref{eq:comparison-gap} with $D_{\rm lo}$ in place of $D_{\rm ref}$.
\Cref{app:geointerval} specifies the normalization and interval evaluation.
The three runs on the base grid $(M,N)=(63,64)$ use IPM tolerance $10^{-12}$; the eight runs with $M\in\{127,255\}$ or $N\in\{128,256\}$ use $10^{-10}$.
The floating-point dual values of the three base-grid runs are the reference values $D_{\rm ref}$ of the three endpoint pairs in \eqref{eq:comparison-gap}.
\Cref{tab:ipmvalidation} reports $\gaprel^{\rm up}$ and the constraint defects of the raw solver arrays.
For all eleven runs, $\gaprel^{\rm up}\le4.7\cdot10^{-8}$, so the action of each IPM candidate is within this relative distance of the minimum $E_{M,N}$ for the normalized endpoints; the raw arrays satisfy unit mass and conservation to $\delta_{\rm mass}\le1.7\cdot10^{-15}$ and $\delta_{\rm cons}\le4.5\cdot10^{-11}$.

\begin{table}[htbp]
\centering\small
\begin{tabular}{@{}lrrlrrrr@{}}
\toprule
Pair & $M$ & $N$ & Candidate & Smallest entry & $\gaprel^{\rm up}$ & $\delta_{\rm mass}$ & $\delta_{\rm cons}$\\
\midrule
Peaked & 63 & 64 & Unit mix. & $1.65\cdot10^{-11}$ & $2.02\cdot10^{-11}$ & $7.56\cdot10^{-16}$ & $2.43\cdot10^{-11}$\\
Peaked & 127 & 64 & Clipped & $4.74\cdot10^{-9}$ & $5.20\cdot10^{-9}$ & $1.15\cdot10^{-15}$ & $6.26\cdot10^{-12}$\\
Peaked & 255 & 64 & Clipped & $8.75\cdot10^{-10}$ & $5.73\cdot10^{-10}$ & $1.65\cdot10^{-15}$ & $4.04\cdot10^{-11}$\\
Peaked & 63 & 128 & Linear mix. & $6.88\cdot10^{-9}$ & $6.01\cdot10^{-9}$ & $7.49\cdot10^{-16}$ & $9.94\cdot10^{-13}$\\
Peaked & 63 & 256 & Clipped & $5.22\cdot10^{-8}$ & $4.66\cdot10^{-8}$ & $5.98\cdot10^{-16}$ & $8.83\cdot10^{-13}$\\
Box & 63 & 64 & Clipped & $2.46\cdot10^{-11}$ & $5.36\cdot10^{-11}$ & $4.16\cdot10^{-16}$ & $1.67\cdot10^{-11}$\\
Diffuse & 63 & 64 & Unit mix. & $5.19\cdot10^{-10}$ & $3.19\cdot10^{-10}$ & $7.72\cdot10^{-16}$ & $4.42\cdot10^{-11}$\\
Diffuse & 127 & 64 & Unit mix. & $2.49\cdot10^{-8}$ & $1.73\cdot10^{-8}$ & $8.91\cdot10^{-16}$ & $4.47\cdot10^{-12}$\\
Diffuse & 255 & 64 & Clipped & $3.04\cdot10^{-9}$ & $2.08\cdot10^{-9}$ & $9.68\cdot10^{-16}$ & $2.43\cdot10^{-11}$\\
Diffuse & 63 & 128 & Unit mix. & $4.81\cdot10^{-9}$ & $3.00\cdot10^{-9}$ & $8.74\cdot10^{-16}$ & $9.82\cdot10^{-13}$\\
Diffuse & 63 & 256 & Unit mix. & $3.66\cdot10^{-8}$ & $2.29\cdot10^{-8}$ & $8.67\cdot10^{-16}$ & $7.04\cdot10^{-13}$\\
\bottomrule
\end{tabular}
\caption{Interval action bounds and constraint defects of eleven IPM geodesics.
Candidate: path selected by \cref{alg:geo-candidates}, either the clipped path (Clipped) or its mixture with the unit-density base (Unit mix.) or with the linear-interpolation base (Linear mix.).
Smallest entry: $\min_{0<j<N,\,0\le i<M}\nu_{{\rm IPM},i}^j$.
$\gaprel^{\rm up}=(V_{\rm up}-D_{\rm lo})/D_{\rm lo}$, rounded upward from the interval evaluations.
$\delta_{\rm mass}$ and $\delta_{\rm cons}$ are the constraint defects of the raw arrays defined in \cref{sec:ipmgeovalidation}, evaluated in NumPy \texttt{longdouble} arithmetic (x87 extended precision, 64-bit significand) without interval enclosure.
IPM tolerance $10^{-12}$ at $(M,N)=(63,64)$ and $10^{-10}$ otherwise.
Every run ended with Clarabel status \texttt{Solved}; the interval bounds do not use this status.
After exact normalization, the IPM candidate differs from the raw density by less than $9.1\cdot10^{-14}$ at every entry and by less than $1.6\cdot10^{-15}$ at the endpoints.
The interval action bounds apply to the normalized endpoint data and fixed grids specified in the text; they do not enclose a continuous geodesic.}
\label{tab:ipmvalidation}
\end{table}

\begin{figure}[p]
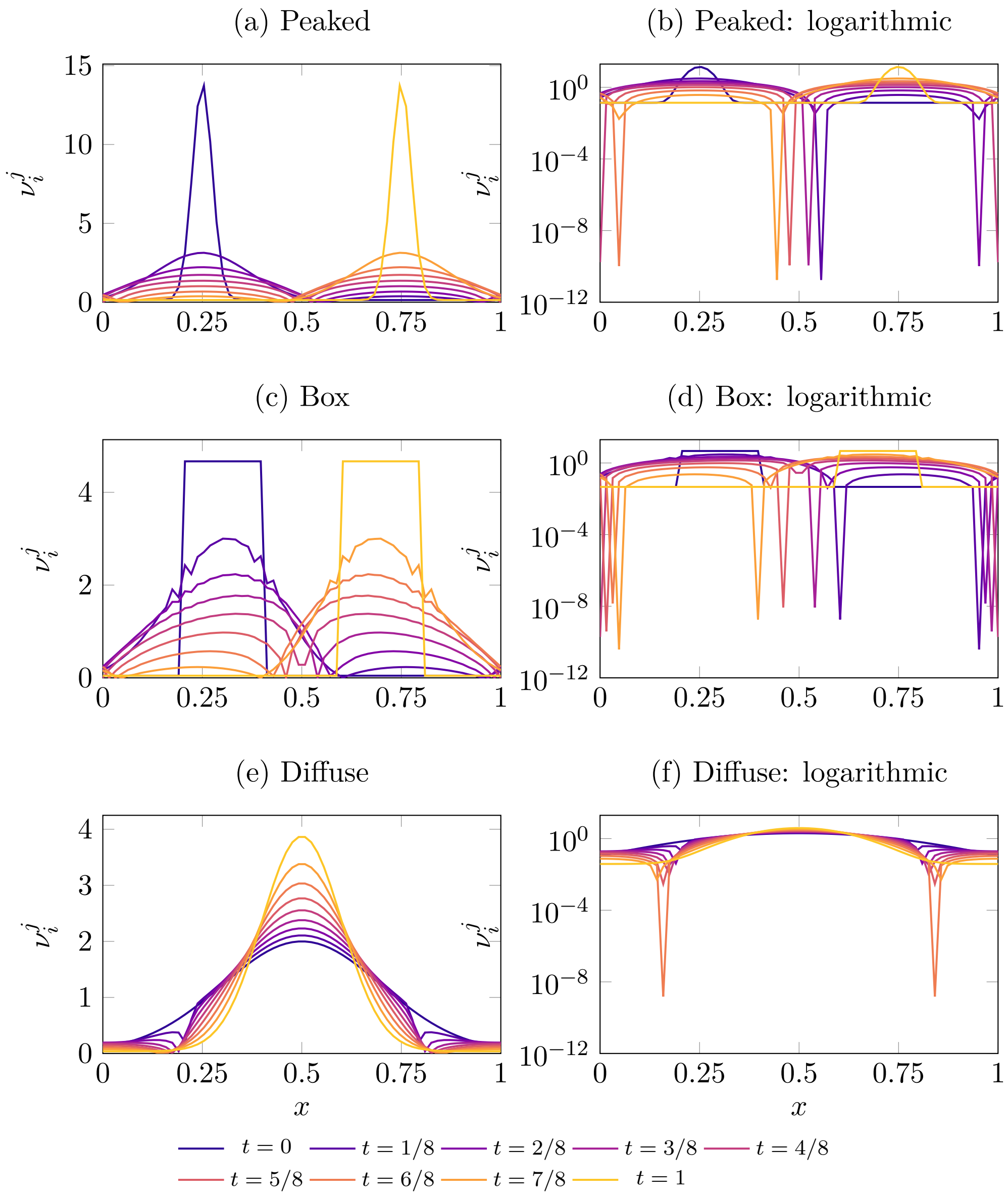

\centering

\pgfplotslegendfromname{legend:ipmgeoprofiles}
\caption{IPM geodesic densities at $M=63$, $N=64$, and solver tolerance $10^{-12}$.
The plotted arrays are the IPM candidates used for the interval action bounds in \cref{tab:ipmvalidation}.
Colors identify the nine times $j/N=0,1/8,\ldots,1$.
Lines connect every spatial node and close the periodic interval; no spatial smoothing or positive display floor is applied.
The linear and logarithmic panels show the same values.
The interval action bounds of \cref{tab:ipmvalidation} give no pointwise bound on these densities.}
\label{fig:ipmgeoprofiles}
\end{figure}

\Cref{fig:ipmgeoprofiles} shows the IPM candidates at the nine times $t=j/N$, $j=0,N/8,\ldots,N$.
At interior times their smallest entries are $1.65\cdot10^{-11}$ (peaked), $2.46\cdot10^{-11}$ (box) and $5.19\cdot10^{-10}$ (diffuse) (\cref{tab:ipmvalidation}), although both endpoints are positive.
Interior-point iterates keep the density slacks strictly positive, so these values do not show whether a minimizer of \eqref{eq:geodesic} vanishes at these nodes.

For the path comparison on the base grid, let $\nu_{\rm DR}$ be the unmodified density component of the DR proximal output $X^\ell$ from the saved run with $\beta=10^6$, $b=1$, stopped when $r_\ell<10^{-12}$.
We compute
\[
 \Delta_{\rm path}=\max_{1\le j<N}\normone{\nu_{\rm IPM}^j-\nu_{\rm DR}^j}.
\]
\Cref{tab:pathchecks} reports $\Delta_{\rm path}$ and the largest entry of the IPM candidate at the nodes where the DR density is exactly zero.
For the three pairs, $\Delta_{\rm path}\le1.5\cdot10^{-7}$, and this entry is at most $7.8\cdot10^{-8}$.
On the base grid the two solvers therefore give the same path to this accuracy; neither quantity certifies either path or its zero set.

\subsection{IPM paths under grid refinement}\label{sec:ipmgeorefine}
\begin{figure}[p]
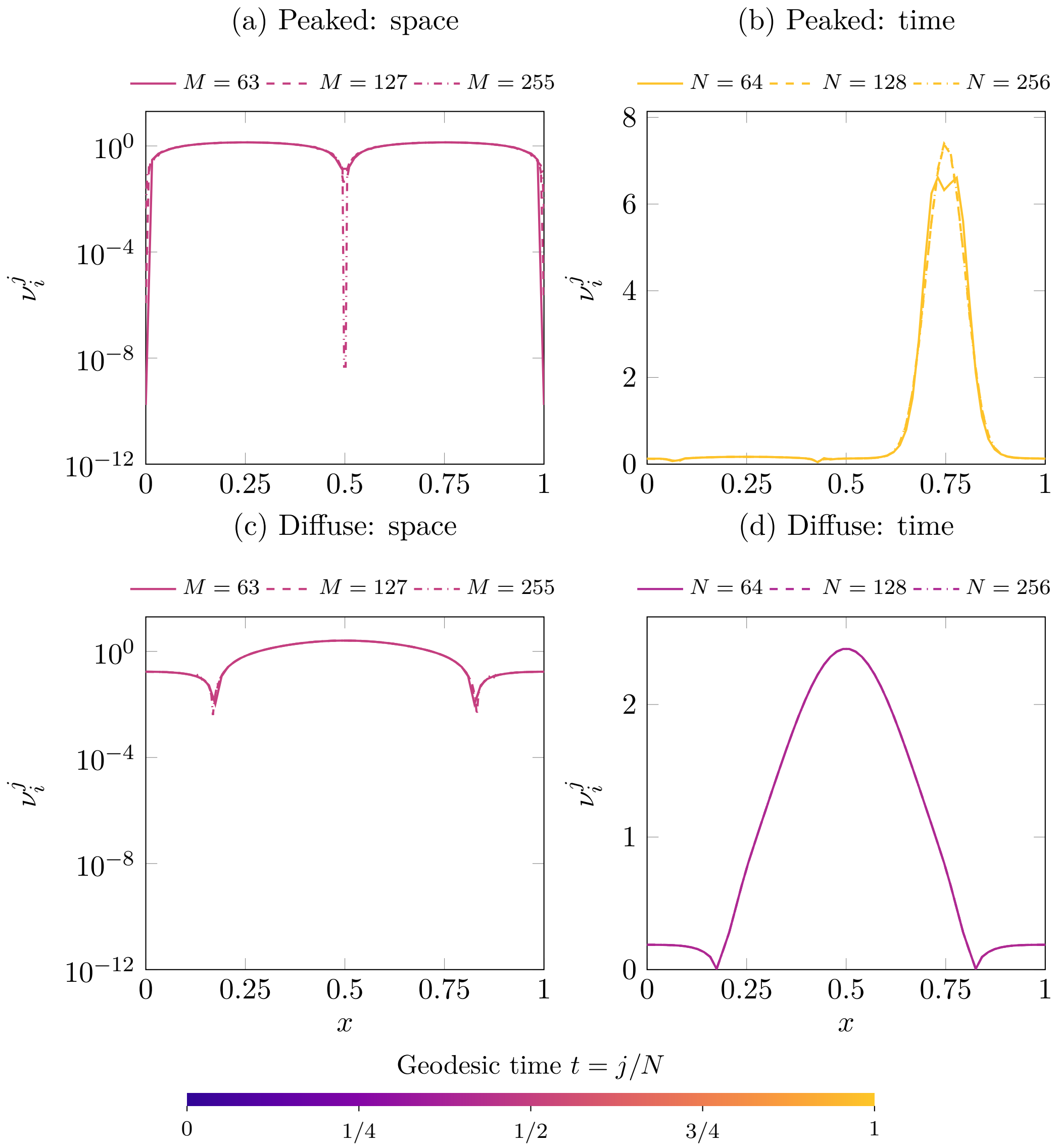

\centering

\caption{IPM densities under separate space and time refinements.
Left: $M=63,127,255$ at fixed $N=64$, shown at $t=1/2$ on a logarithmic scale.
Right: $N=64,128,256$ at fixed $M=63$, on a linear scale at $t=63/64$ for the peaked pair and $t=13/32$ for the diffuse pair.
These are the times $j/64$ that attain the maximum over $j$ defining $\Delta_{64}$ (\cref{tab:pathchecks}), the largest discrete $L^1$ difference between the $N=64$ and $N=128$ IPM candidates on common time levels.
Colors give the path time on the color bar below the panels; solid, dashed and dash-dotted lines distinguish the grids listed in each legend.
Each line connects all nodes of its own grid and closes the periodic interval.
The plotted densities are the IPM candidates of \cref{tab:ipmvalidation}.
The IPM tolerances are $10^{-12}$ on the base grid and $10^{-10}$ on the refined grids.
No positive display floor is used.
The curves are numerical path comparisons, not enclosures of the continuous solution.}
\label{fig:ipmgeorefine}
\end{figure}

We solve \eqref{eq:geodesic} by IPM for the peaked and diffuse pairs on refined spatial grids, $M=127,255$ at $N=64$, and on refined time grids, $N=128,256$ at $M=63$; \cref{fig:ipmgeorefine} compares the resulting IPM candidates.
The interval action bounds for these paths are in \cref{tab:ipmvalidation}.

At fixed $M=63$, let $\nu_{(N)}^j$ be the density at time $j/N$ of the IPM candidate computed with $N$ time steps.
On the common time levels, define
\[
 \Delta_N=\max_{0\le j\le N}\normone{\nu_{(N)}^j-\nu_{(2N)}^{2j}}.
\]
\Cref{tab:pathchecks} compares these time-refinement differences with the DR--IPM difference $\Delta_{\rm path}$ on the base grid.
For the peaked pair, $\Delta_{64}=8.14\cdot10^{-2}$ exceeds $\Delta_{\rm path}=1.28\cdot10^{-7}$ by a factor of about $6\cdot10^5$.
These computed-path differences assess grid sensitivity; they are not certified discretization errors.
Doubling $N$ from $64$ to $128$ reduces $\Delta_N$: the ratio $\Delta_{64}/\Delta_{128}$ is about $3.1$ for the peaked pair and $2.4$ for the diffuse pair.
No difference measure is computed under spatial refinement; \cref{fig:ipmgeorefine} compares the profiles of the three spatial grids at $t=1/2$ only.

\begin{table}[htbp]
\centering\small
\begin{tabular}{@{}lrrrr@{}}
\toprule
Pair & $\Delta_{\rm path}$ & IPM density at DR zeros & $\Delta_{64}$ & $\Delta_{128}$\\
\midrule
Peaked & $1.28\cdot10^{-7}$ & $1.41\cdot10^{-9}$ & $8.14\cdot10^{-2}$ & $2.62\cdot10^{-2}$\\
Box & $1.48\cdot10^{-7}$ & $7.79\cdot10^{-8}$ & -- & --\\
Diffuse & $1.22\cdot10^{-7}$ & $3.60\cdot10^{-9}$ & $2.70\cdot10^{-4}$ & $1.10\cdot10^{-4}$\\
\bottomrule
\end{tabular}
\caption{DR--IPM path difference and IPM time-refinement differences. At $M=63$, $N=64$, $\Delta_{\rm path}=\max_{1\le j<N}\normone{\nu_{\rm IPM}^j-\nu_{\rm DR}^j}$ compares the IPM candidate with the unmodified DR proximal density; ``IPM density at DR zeros'' is the largest entry of the IPM candidate at the nodes where the DR density is exactly zero. DR uses $\beta=10^6$, $b=1$, and $r_\ell<10^{-12}$. $\Delta_{64}$ and $\Delta_{128}$ are the differences $\Delta_N$ of \cref{sec:ipmgeorefine}, at $M=63$; the box pair is not refined. The interval action bounds of the IPM candidates are in \cref{tab:ipmvalidation}.}
\label{tab:pathchecks}
\end{table}

\subsection{Zero entries in proximal outputs}\label{sec:proximal-zeros}
We ask whether a minimizer of \eqref{eq:geodesic} between positive endpoints has vacuum entries.
A first indication is the presence of exact zero entries in DR proximal outputs stopped at $r_\ell<10^{-12}$; these finite-iteration outputs do not determine the zero set of a minimizer.
The diffusivity transitions in \cref{sec:vacuumobservation} give a second test based on independent DR and IPM computations.

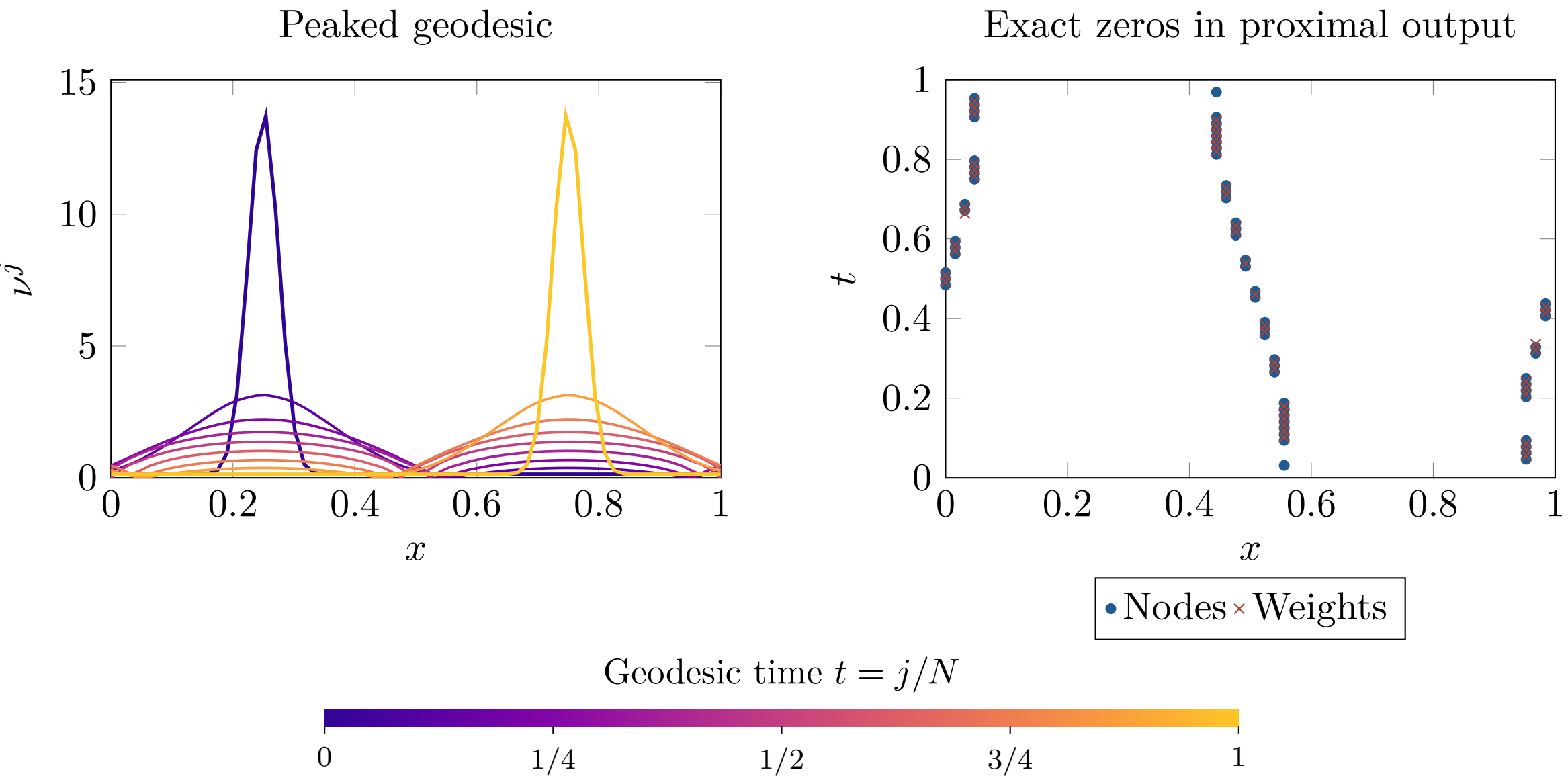
\begin{figure}[htbp]
\centering
\begin{tikzpicture}
\begin{axis}[width=.48\textwidth,height=5.5cm,xmin=0,xmax=1,ymin=0,
 xlabel={$x$},ylabel={$\nu^j$},title={Peaked geodesic},
 legend style={at={(.98,.98)},anchor=north east}]
\addplot[geotime0,line width=1pt] coordinates {(0,0.1376164966) (0.01587301587,0.1376164966) (0.03174603175,0.1376164966) (0.04761904762,0.1376164966) (0.06349206349,0.1376164966) (0.07936507937,0.1376164976) (0.09523809524,0.1376165622) (0.1111111111,0.1376192293) (0.126984127,0.1376925223) (0.1428571429,0.1390298533) (0.1587301587,0.1551742186) (0.1746031746,0.2833669547) (0.1904761905,0.9461109931) (0.2063492063,3.134498105) (0.2222222222,7.560753634) (0.2380952381,12.42418935) (0.253968254,13.72698938) (0.2698412698,10.18130997) (0.2857142857,5.097972683) (0.3015873016,1.774650328) (0.3174603175,0.4986336403) (0.3333333333,0.1908179351) (0.3492063492,0.1428554472) (0.3650793651,0.1379612357) (0.380952381,0.1376316553) (0.3968253968,0.137616942) (0.4126984127,0.1376165053) (0.4285714286,0.1376164967) (0.4444444444,0.1376164966) (0.4603174603,0.1376164966) (0.4761904762,0.1376164966) (0.4920634921,0.1376164966) (0.5079365079,0.1376164966) (0.5238095238,0.1376164966) (0.5396825397,0.1376164966) (0.5555555556,0.1376164966) (0.5714285714,0.1376164966) (0.5873015873,0.1376164966) (0.6031746032,0.1376164966) (0.619047619,0.1376164966) (0.6349206349,0.1376164966) (0.6507936508,0.1376164966) (0.6666666667,0.1376164966) (0.6825396825,0.1376164966) (0.6984126984,0.1376164966) (0.7142857143,0.1376164966) (0.7301587302,0.1376164966) (0.746031746,0.1376164966) (0.7619047619,0.1376164966) (0.7777777778,0.1376164966) (0.7936507937,0.1376164966) (0.8095238095,0.1376164966) (0.8253968254,0.1376164966) (0.8412698413,0.1376164966) (0.8571428571,0.1376164966) (0.873015873,0.1376164966) (0.8888888889,0.1376164966) (0.9047619048,0.1376164966) (0.9206349206,0.1376164966) (0.9365079365,0.1376164966) (0.9523809524,0.1376164966) (0.9682539683,0.1376164966) (0.9841269841,0.1376164966) (1,0.1376164966)};
\addplot[geotime1,line width=.7pt] coordinates {(0,0.2642677847) (0.01587301587,0.3824801482) (0.03174603175,0.5211671648) (0.04761904762,0.6852722192) (0.06349206349,0.8735478869) (0.07936507937,1.085901721) (0.09523809524,1.318682909) (0.1111111111,1.567717719) (0.126984127,1.82587731) (0.1428571429,2.083640654) (0.1587301587,2.329941559) (0.1746031746,2.560803695) (0.1904761905,2.768143901) (0.2063492063,2.928895298) (0.2222222222,3.036356442) (0.2380952381,3.114364762) (0.253968254,3.132579191) (0.2698412698,3.075282433) (0.2857142857,2.978484419) (0.3015873016,2.844316754) (0.3174603175,2.653506487) (0.3333333333,2.430490595) (0.3492063492,2.189937287) (0.3650793651,1.937598339) (0.380952381,1.679985416) (0.3968253968,1.429480476) (0.4126984127,1.192862527) (0.4285714286,0.9760917275) (0.4444444444,0.7809052019) (0.4603174603,0.6100266822) (0.4761904762,0.4617696154) (0.4920634921,0.3358549744) (0.5079365079,0.2293913635) (0.5238095238,0.1389136825) (0.5396825397,0.06515079932) (0.5555555556,0) (0.5714285714,0.05534972563) (0.5873015873,0.1157507991) (0.6031746032,0.1559644531) (0.619047619,0.2062114826) (0.6349206349,0.244034917) (0.6507936508,0.2798242081) (0.6666666667,0.3083080443) (0.6825396825,0.3323263032) (0.6984126984,0.3505342772) (0.7142857143,0.3641697835) (0.7301587302,0.3726341006) (0.746031746,0.376649258) (0.7619047619,0.3757648022) (0.7777777778,0.3704277116) (0.7936507937,0.360161953) (0.8095238095,0.345284851) (0.8253968254,0.3252891106) (0.8412698413,0.3005625424) (0.8571428571,0.2707525015) (0.873015873,0.236217155) (0.8888888889,0.1975573403) (0.9047619048,0.1524057462) (0.9206349206,0.09982865062) (0.9365079365,0.04570673843) (0.9523809524,0.01733725778) (0.9682539683,0.08345970166) (0.9841269841,0.1677974405) (1,0.2642677847)};
\addplot[geotime2,line width=.7pt] coordinates {(0,0.4679851983) (0.01587301587,0.6326703983) (0.03174603175,0.7935289149) (0.04761904762,0.9562769546) (0.06349206349,1.118805387) (0.07936507937,1.275357319) (0.09523809524,1.426636913) (0.1111111111,1.567601868) (0.126984127,1.699073887) (0.1428571429,1.818002417) (0.1587301587,1.924285845) (0.1746031746,2.015398615) (0.1904761905,2.091507857) (0.2063492063,2.150807416) (0.2222222222,2.191338335) (0.2380952381,2.214768844) (0.253968254,2.219287411) (0.2698412698,2.203867435) (0.2857142857,2.170785304) (0.3015873016,2.119831707) (0.3174603175,2.050373359) (0.3333333333,1.965582032) (0.3492063492,1.865478198) (0.3650793651,1.752005023) (0.380952381,1.626574838) (0.3968253968,1.491717062) (0.4126984127,1.348050618) (0.4285714286,1.19902304) (0.4444444444,1.042136626) (0.4603174603,0.8835717074) (0.4761904762,0.7134006901) (0.4920634921,0.5442460152) (0.5079365079,0.3661544044) (0.5238095238,0.1850209434) (0.5396825397,0.03359133426) (0.5555555556,0.1048267824) (0.5714285714,0.2149241257) (0.5873015873,0.3100110129) (0.6031746032,0.3796010569) (0.619047619,0.4492993257) (0.6349206349,0.5001793347) (0.6507936508,0.5476154777) (0.6666666667,0.5848747589) (0.6825396825,0.6163284322) (0.6984126984,0.6399579866) (0.7142857143,0.657818688) (0.7301587302,0.6688010353) (0.746031746,0.6742449072) (0.7619047619,0.6731855162) (0.7777777778,0.6666581241) (0.7936507937,0.6535608435) (0.8095238095,0.6347725752) (0.8253968254,0.608826694) (0.8412698413,0.5765057171) (0.8571428571,0.5357223681) (0.873015873,0.4865005632) (0.8888888889,0.4261662052) (0.9047619048,0.3549548091) (0.9206349206,0.2626651692) (0.9365079365,0.1511082535) (0.9523809524,0) (0.9682539683,0.168569692) (0.9841269841,0.327576627) (1,0.4679851983)};
\addplot[geotime3,line width=.7pt] coordinates {(0,0.3457806553) (0.01587301587,0.555379051) (0.03174603175,0.7511331044) (0.04761904762,0.9040351456) (0.06349206349,1.049196446) (0.07936507937,1.162669494) (0.09523809524,1.270647346) (0.1111111111,1.360219522) (0.126984127,1.443288804) (0.1428571429,1.512940582) (0.1587301587,1.574999037) (0.1746031746,1.625279711) (0.1904761905,1.667092032) (0.2063492063,1.698114768) (0.2222222222,1.719912929) (0.2380952381,1.731554991) (0.253968254,1.733721659) (0.2698412698,1.725932588) (0.2857142857,1.708469941) (0.3015873016,1.681506755) (0.3174603175,1.644260494) (0.3333333333,1.597764384) (0.3492063492,1.540333812) (0.3650793651,1.473176501) (0.380952381,1.393806437) (0.3968253968,1.303970763) (0.4126984127,1.199713328) (0.4285714286,1.085949267) (0.4444444444,0.9523440206) (0.4603174603,0.8169450304) (0.4761904762,0.64223057) (0.4920634921,0.4941838778) (0.5079365079,0.2679325376) (0.5238095238,0) (0.5396825397,0.2364831795) (0.5555555556,0.4044658305) (0.5714285714,0.5083365715) (0.5873015873,0.6170917615) (0.6031746032,0.6985418662) (0.619047619,0.7746330602) (0.6349206349,0.8314846271) (0.6507936508,0.8836176604) (0.6666666667,0.923019194) (0.6825396825,0.9565269487) (0.6984126984,0.9806795441) (0.7142857143,0.9993045066) (0.7301587302,1.01001057) (0.746031746,1.015656993) (0.7619047619,1.013902541) (0.7777777778,1.007290736) (0.7936507937,0.993129081) (0.8095238095,0.9740479965) (0.8253968254,0.9465993549) (0.8412698413,0.9137712229) (0.8571428571,0.8707240945) (0.873015873,0.8207379043) (0.8888888889,0.7556035006) (0.9047619048,0.6805097484) (0.9206349206,0.5787973752) (0.9365079365,0.4548698418) (0.9523809524,0.2766427644) (0.9682539683,0.1167883799) (0.9841269841,0.1222475607) (1,0.3457806553)};
\addplot[geotime4,line width=.7pt] coordinates {(0,0) (0.01587301587,0.2985881143) (0.03174603175,0.5133608207) (0.04761904762,0.6389008922) (0.06349206349,0.7910741048) (0.07936507937,0.8935490463) (0.09523809524,0.9959105763) (0.1111111111,1.073559937) (0.126984127,1.145151666) (0.1428571429,1.200001654) (0.1587301587,1.248820393) (0.1746031746,1.285324593) (0.1904761905,1.316196767) (0.2063492063,1.337292464) (0.2222222222,1.352957508) (0.2380952381,1.360152958) (0.253968254,1.361806622) (0.2698412698,1.355505982) (0.2857142857,1.343335192) (0.3015873016,1.323247918) (0.3174603175,1.296802083) (0.3333333333,1.261923318) (0.3492063492,1.220303334) (0.3650793651,1.168592253) (0.380952381,1.10953175) (0.3968253968,1.036984495) (0.4126984127,0.9545630325) (0.4285714286,0.8490558623) (0.4444444444,0.7222966187) (0.4603174603,0.5638614746) (0.4761904762,0.3439071753) (0.4920634921,0.137441394) (0.5079365079,0.137441394) (0.5238095238,0.3439071753) (0.5396825397,0.5638614746) (0.5555555556,0.7222966187) (0.5714285714,0.8490558623) (0.5873015873,0.9545630325) (0.6031746032,1.036984495) (0.619047619,1.10953175) (0.6349206349,1.168592253) (0.6507936508,1.220303334) (0.6666666667,1.261923318) (0.6825396825,1.296802083) (0.6984126984,1.323247918) (0.7142857143,1.343335192) (0.7301587302,1.355505982) (0.746031746,1.361806622) (0.7619047619,1.360152958) (0.7777777778,1.352957508) (0.7936507937,1.337292464) (0.8095238095,1.316196767) (0.8253968254,1.285324593) (0.8412698413,1.248820393) (0.8571428571,1.200001654) (0.873015873,1.145151666) (0.8888888889,1.073559937) (0.9047619048,0.9959105763) (0.9206349206,0.8935490463) (0.9365079365,0.7910741048) (0.9523809524,0.6389008922) (0.9682539683,0.5133608207) (0.9841269841,0.2985881143) (1,0)};
\addplot[geotime5,line width=.7pt] coordinates {(0,0.3457806553) (0.01587301587,0.1222475607) (0.03174603175,0.1167883799) (0.04761904762,0.2766427644) (0.06349206349,0.4548698418) (0.07936507937,0.5787973752) (0.09523809524,0.6805097484) (0.1111111111,0.7556035006) (0.126984127,0.8207379043) (0.1428571429,0.8707240945) (0.1587301587,0.9137712229) (0.1746031746,0.9465993549) (0.1904761905,0.9740479965) (0.2063492063,0.993129081) (0.2222222222,1.007290736) (0.2380952381,1.013902541) (0.253968254,1.015656993) (0.2698412698,1.01001057) (0.2857142857,0.9993045066) (0.3015873016,0.9806795441) (0.3174603175,0.9565269487) (0.3333333333,0.923019194) (0.3492063492,0.8836176604) (0.3650793651,0.8314846271) (0.380952381,0.7746330602) (0.3968253968,0.6985418662) (0.4126984127,0.6170917615) (0.4285714286,0.5083365715) (0.4444444444,0.4044658305) (0.4603174603,0.2364831795) (0.4761904762,0) (0.4920634921,0.2679325376) (0.5079365079,0.4941838778) (0.5238095238,0.64223057) (0.5396825397,0.8169450304) (0.5555555556,0.9523440206) (0.5714285714,1.085949267) (0.5873015873,1.199713328) (0.6031746032,1.303970763) (0.619047619,1.393806437) (0.6349206349,1.473176501) (0.6507936508,1.540333812) (0.6666666667,1.597764384) (0.6825396825,1.644260494) (0.6984126984,1.681506755) (0.7142857143,1.708469941) (0.7301587302,1.725932588) (0.746031746,1.733721659) (0.7619047619,1.731554991) (0.7777777778,1.719912929) (0.7936507937,1.698114768) (0.8095238095,1.667092032) (0.8253968254,1.625279711) (0.8412698413,1.574999037) (0.8571428571,1.512940582) (0.873015873,1.443288804) (0.8888888889,1.360219522) (0.9047619048,1.270647346) (0.9206349206,1.162669494) (0.9365079365,1.049196446) (0.9523809524,0.9040351456) (0.9682539683,0.7511331044) (0.9841269841,0.555379051) (1,0.3457806553)};
\addplot[geotime6,line width=.7pt] coordinates {(0,0.4679851983) (0.01587301587,0.327576627) (0.03174603175,0.168569692) (0.04761904762,0) (0.06349206349,0.1511082535) (0.07936507937,0.2626651692) (0.09523809524,0.3549548091) (0.1111111111,0.4261662052) (0.126984127,0.4865005632) (0.1428571429,0.5357223681) (0.1587301587,0.5765057171) (0.1746031746,0.608826694) (0.1904761905,0.6347725752) (0.2063492063,0.6535608435) (0.2222222222,0.6666581241) (0.2380952381,0.6731855162) (0.253968254,0.6742449072) (0.2698412698,0.6688010353) (0.2857142857,0.657818688) (0.3015873016,0.6399579866) (0.3174603175,0.6163284322) (0.3333333333,0.5848747589) (0.3492063492,0.5476154777) (0.3650793651,0.5001793347) (0.380952381,0.4492993257) (0.3968253968,0.3796010569) (0.4126984127,0.3100110129) (0.4285714286,0.2149241257) (0.4444444444,0.1048267824) (0.4603174603,0.03359133426) (0.4761904762,0.1850209434) (0.4920634921,0.3661544044) (0.5079365079,0.5442460152) (0.5238095238,0.7134006901) (0.5396825397,0.8835717074) (0.5555555556,1.042136626) (0.5714285714,1.19902304) (0.5873015873,1.348050618) (0.6031746032,1.491717062) (0.619047619,1.626574838) (0.6349206349,1.752005023) (0.6507936508,1.865478198) (0.6666666667,1.965582032) (0.6825396825,2.050373359) (0.6984126984,2.119831707) (0.7142857143,2.170785304) (0.7301587302,2.203867435) (0.746031746,2.219287411) (0.7619047619,2.214768844) (0.7777777778,2.191338335) (0.7936507937,2.150807416) (0.8095238095,2.091507857) (0.8253968254,2.015398615) (0.8412698413,1.924285845) (0.8571428571,1.818002417) (0.873015873,1.699073887) (0.8888888889,1.567601868) (0.9047619048,1.426636913) (0.9206349206,1.275357319) (0.9365079365,1.118805387) (0.9523809524,0.9562769546) (0.9682539683,0.7935289149) (0.9841269841,0.6326703983) (1,0.4679851983)};
\addplot[geotime7,line width=.7pt] coordinates {(0,0.2642677847) (0.01587301587,0.1677974405) (0.03174603175,0.08345970166) (0.04761904762,0.01733725778) (0.06349206349,0.04570673843) (0.07936507937,0.09982865062) (0.09523809524,0.1524057462) (0.1111111111,0.1975573403) (0.126984127,0.236217155) (0.1428571429,0.2707525015) (0.1587301587,0.3005625424) (0.1746031746,0.3252891106) (0.1904761905,0.345284851) (0.2063492063,0.360161953) (0.2222222222,0.3704277116) (0.2380952381,0.3757648022) (0.253968254,0.376649258) (0.2698412698,0.3726341006) (0.2857142857,0.3641697835) (0.3015873016,0.3505342772) (0.3174603175,0.3323263032) (0.3333333333,0.3083080443) (0.3492063492,0.2798242081) (0.3650793651,0.244034917) (0.380952381,0.2062114826) (0.3968253968,0.1559644531) (0.4126984127,0.1157507991) (0.4285714286,0.05534972563) (0.4444444444,0) (0.4603174603,0.06515079932) (0.4761904762,0.1389136825) (0.4920634921,0.2293913635) (0.5079365079,0.3358549744) (0.5238095238,0.4617696154) (0.5396825397,0.6100266822) (0.5555555556,0.7809052019) (0.5714285714,0.9760917275) (0.5873015873,1.192862527) (0.6031746032,1.429480476) (0.619047619,1.679985416) (0.6349206349,1.937598339) (0.6507936508,2.189937287) (0.6666666667,2.430490595) (0.6825396825,2.653506487) (0.6984126984,2.844316754) (0.7142857143,2.978484419) (0.7301587302,3.075282433) (0.746031746,3.132579191) (0.7619047619,3.114364762) (0.7777777778,3.036356442) (0.7936507937,2.928895298) (0.8095238095,2.768143901) (0.8253968254,2.560803695) (0.8412698413,2.329941559) (0.8571428571,2.083640654) (0.873015873,1.82587731) (0.8888888889,1.567717719) (0.9047619048,1.318682909) (0.9206349206,1.085901721) (0.9365079365,0.8735478869) (0.9523809524,0.6852722192) (0.9682539683,0.5211671648) (0.9841269841,0.3824801482) (1,0.2642677847)};
\addplot[geotime8,line width=1pt] coordinates {(0,0.1376164966) (0.01587301587,0.1376164966) (0.03174603175,0.1376164966) (0.04761904762,0.1376164966) (0.06349206349,0.1376164966) (0.07936507937,0.1376164966) (0.09523809524,0.1376164966) (0.1111111111,0.1376164966) (0.126984127,0.1376164966) (0.1428571429,0.1376164966) (0.1587301587,0.1376164966) (0.1746031746,0.1376164966) (0.1904761905,0.1376164966) (0.2063492063,0.1376164966) (0.2222222222,0.1376164966) (0.2380952381,0.1376164966) (0.253968254,0.1376164966) (0.2698412698,0.1376164966) (0.2857142857,0.1376164966) (0.3015873016,0.1376164966) (0.3174603175,0.1376164966) (0.3333333333,0.1376164966) (0.3492063492,0.1376164966) (0.3650793651,0.1376164966) (0.380952381,0.1376164966) (0.3968253968,0.1376164966) (0.4126984127,0.1376164966) (0.4285714286,0.1376164966) (0.4444444444,0.1376164966) (0.4603174603,0.1376164966) (0.4761904762,0.1376164966) (0.4920634921,0.1376164966) (0.5079365079,0.1376164966) (0.5238095238,0.1376164966) (0.5396825397,0.1376164966) (0.5555555556,0.1376164966) (0.5714285714,0.1376164967) (0.5873015873,0.1376165053) (0.6031746032,0.137616942) (0.619047619,0.1376316553) (0.6349206349,0.1379612357) (0.6507936508,0.1428554472) (0.6666666667,0.1908179351) (0.6825396825,0.4986336403) (0.6984126984,1.774650328) (0.7142857143,5.097972683) (0.7301587302,10.18130997) (0.746031746,13.72698938) (0.7619047619,12.42418935) (0.7777777778,7.560753634) (0.7936507937,3.134498105) (0.8095238095,0.9461109931) (0.8253968254,0.2833669547) (0.8412698413,0.1551742186) (0.8571428571,0.1390298533) (0.873015873,0.1376925223) (0.8888888889,0.1376192293) (0.9047619048,0.1376165622) (0.9206349206,0.1376164976) (0.9365079365,0.1376164966) (0.9523809524,0.1376164966) (0.9682539683,0.1376164966) (0.9841269841,0.1376164966) (1,0.1376164966)};
\end{axis}
\begin{axis}[at={(.52\textwidth,0)},anchor=south west,
 width=.48\textwidth,height=5.5cm,xmin=0,xmax=1,ymin=0,ymax=1,
 xlabel={$x$},ylabel={$t$},title={Exact zeros in proximal output},
 legend style={at={(.5,-.25)},anchor=north,legend columns=2}]
\addplot[only marks,mark=*,mark size=1.3pt,plotblue] coordinates {(0.555555556,0.03125) (0.952380952,0.046875) (0.952380952,0.0625) (0.952380952,0.078125) (0.555555556,0.09375) (0.952380952,0.09375) (0.555555556,0.109375) (0.555555556,0.125) (0.555555556,0.140625) (0.555555556,0.15625) (0.555555556,0.171875) (0.555555556,0.1875) (0.952380952,0.203125) (0.952380952,0.21875) (0.952380952,0.234375) (0.952380952,0.25) (0.53968254,0.265625) (0.53968254,0.28125) (0.53968254,0.296875) (0.968253968,0.3125) (0.968253968,0.328125) (0.523809524,0.359375) (0.523809524,0.375) (0.523809524,0.390625) (0.984126984,0.40625) (0.984126984,0.421875) (0.984126984,0.4375) (0.507936508,0.453125) (0.507936508,0.46875) (0,0.484375) (0,0.5) (0,0.515625) (0.492063492,0.53125) (0.492063492,0.546875) (0.0158730159,0.5625) (0.0158730159,0.578125) (0.0158730159,0.59375) (0.476190476,0.609375) (0.476190476,0.625) (0.476190476,0.640625) (0.0317460317,0.671875) (0.0317460317,0.6875) (0.46031746,0.703125) (0.46031746,0.71875) (0.46031746,0.734375) (0.0476190476,0.75) (0.0476190476,0.765625) (0.0476190476,0.78125) (0.0476190476,0.796875) (0.444444444,0.8125) (0.444444444,0.828125) (0.444444444,0.84375) (0.444444444,0.859375) (0.444444444,0.875) (0.444444444,0.890625) (0.0476190476,0.90625) (0.444444444,0.90625) (0.0476190476,0.921875) (0.0476190476,0.9375) (0.0476190476,0.953125) (0.444444444,0.96875)};
\addlegendentry{Nodes}
\addplot[only marks,mark=x,mark size=2pt,plotred] coordinates {(0.952380952,0.0546875) (0.952380952,0.0703125) (0.952380952,0.0859375) (0.555555556,0.1015625) (0.555555556,0.1171875) (0.555555556,0.1328125) (0.555555556,0.1484375) (0.555555556,0.1640625) (0.555555556,0.1796875) (0.952380952,0.2109375) (0.952380952,0.2265625) (0.952380952,0.2421875) (0.53968254,0.2734375) (0.53968254,0.2890625) (0.968253968,0.3203125) (0.968253968,0.3359375) (0.523809524,0.3671875) (0.523809524,0.3828125) (0.984126984,0.4140625) (0.984126984,0.4296875) (0.507936508,0.4609375) (0,0.4921875) (0,0.5078125) (0.492063492,0.5390625) (0.0158730159,0.5703125) (0.0158730159,0.5859375) (0.476190476,0.6171875) (0.476190476,0.6328125) (0.0317460317,0.6640625) (0.0317460317,0.6796875) (0.46031746,0.7109375) (0.46031746,0.7265625) (0.0476190476,0.7578125) (0.0476190476,0.7734375) (0.0476190476,0.7890625) (0.444444444,0.8203125) (0.444444444,0.8359375) (0.444444444,0.8515625) (0.444444444,0.8671875) (0.444444444,0.8828125) (0.444444444,0.8984375) (0.0476190476,0.9140625) (0.0476190476,0.9296875) (0.0476190476,0.9453125)};
\addlegendentry{Weights}
\end{axis}
\end{tikzpicture}

\begin{tikzpicture}[x=1cm,y=1cm]
\shade[left color=geotime0,right color=geotime1] (0.000000,0) rectangle (1.125000,.17);
\shade[left color=geotime1,right color=geotime2] (1.125000,0) rectangle (2.250000,.17);
\shade[left color=geotime2,right color=geotime3] (2.250000,0) rectangle (3.375000,.17);
\shade[left color=geotime3,right color=geotime4] (3.375000,0) rectangle (4.500000,.17);
\shade[left color=geotime4,right color=geotime5] (4.500000,0) rectangle (5.625000,.17);
\shade[left color=geotime5,right color=geotime6] (5.625000,0) rectangle (6.750000,.17);
\shade[left color=geotime6,right color=geotime7] (6.750000,0) rectangle (7.875000,.17);
\shade[left color=geotime7,right color=geotime8] (7.875000,0) rectangle (9.000000,.17);
\draw (0,0)--(0,-.07) node[below,font=\scriptsize] {$0$};
\draw (2.25,0)--(2.25,-.07) node[below,font=\scriptsize] {$1/4$};
\draw (4.5,0)--(4.5,-.07) node[below,font=\scriptsize] {$1/2$};
\draw (6.75,0)--(6.75,-.07) node[below,font=\scriptsize] {$3/4$};
\draw (9,0)--(9,-.07) node[below,font=\scriptsize] {$1$};
\node[font=\small,anchor=south] at (4.5,.20) {Geodesic time $t=j/N$};
\end{tikzpicture}
\caption{Peaked geodesic at $M=63$, $N=64$, $\beta=10^6$, and $b=1$.
Left: densities at $t=0,1/8,\ldots,1$, namely the prescribed geodesic endpoint profiles and the density components $\nu_X^j$ of the DR proximal output $X^\ell$, with time given by the color scale.
Right: points $(x_i,j/N)$ with $\nu_{X,i}^j=0$ (Nodes) and points $(x_i,(j+1/2)/N)$ with $q_{X,i}^j=0$ (Weights) in the same output.
Spatial line segments connect nodes and close the periodic interval.
The display uses no positive floor or zero threshold.}
\label{fig:geodesic}
\end{figure}

\begin{figure}[htbp]
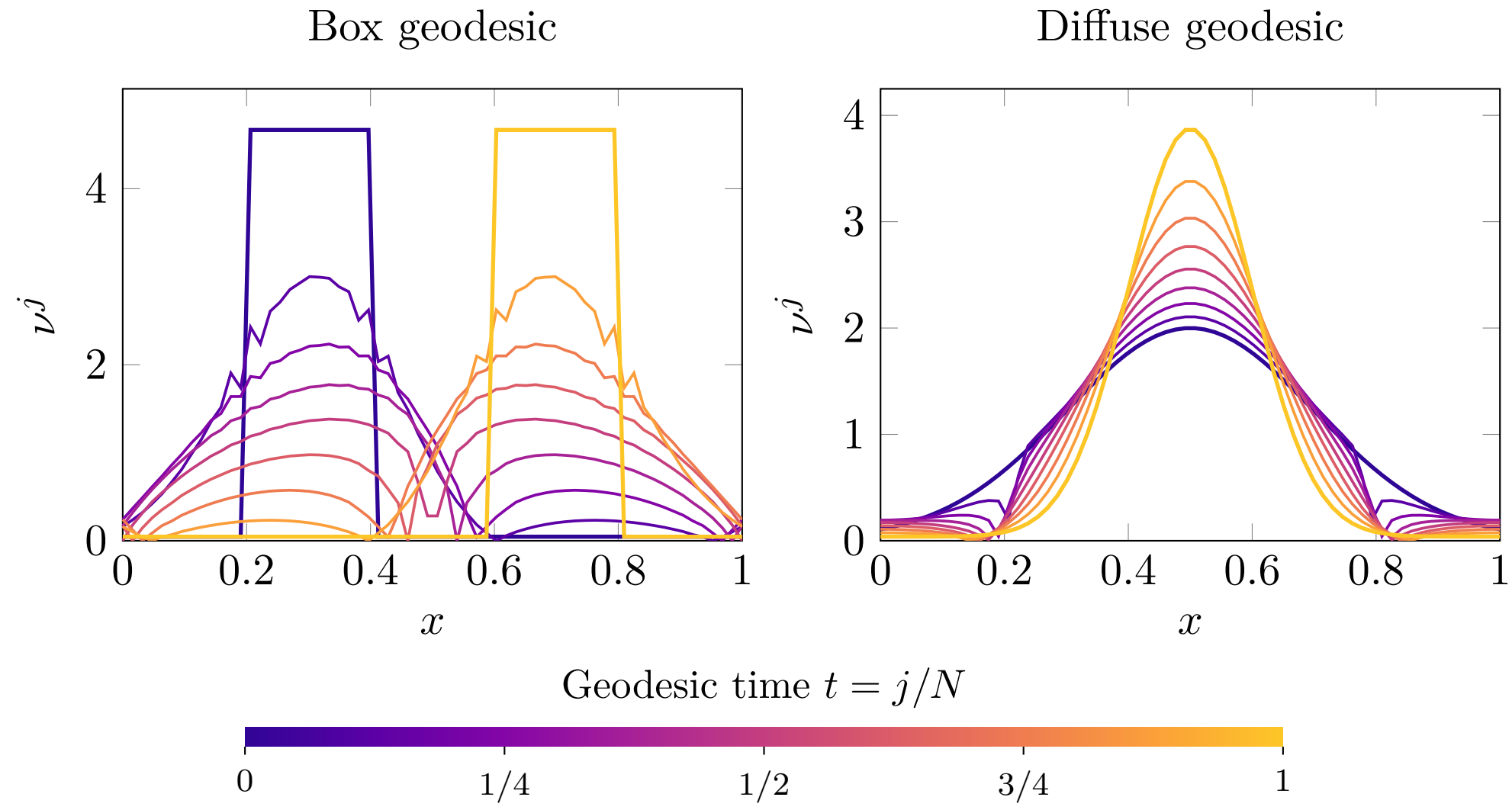

\centering

\caption{Box and diffuse geodesic profiles at nine times $t=j/N=0,1/8,\ldots,1$.
The color scale specifies the same times in both panels.
The grid and solver parameters are those of \cref{tab:zeros}.
Interior profiles are saved nonnegative proximal outputs, without a positive display floor.
Line segments connect the spatial nodes and close the periodic interval.}
\label{fig:geoprofiles}
\end{figure}

\Cref{tab:zeros,fig:geodesic,fig:geoprofiles} use the density and midpoint-weight components $\nu_X^j$ and $q_X^j$ of $X^\ell$ without modification (\cref{sec:geodata}).

\Cref{tab:zeros} uses the saved runs with $b=1$ and $r_\ell<10^{-12}$.
In these runs, the exact zero entries come from the nodal-density update $\nu_{\rm in}\mapsto\max(\nu_{\rm in},0)$ and from the clipping at zero in step 3 of \cref{alg:geo-prox}, which computes the action proximal map \eqref{eq:prox}.
In $X^\ell$ the midpoint weight $q_X^j$ is an independent variable.
The relation $2q^j=\nu^j+\nu^{j+1}$ and the conservation equations of $\K$ are imposed by the affine projection; they hold for $Y^\ell$, and for $X^\ell$ only up to the difference $Y^\ell-X^\ell$ measured by $r_\ell$.
The last column of \cref{tab:zeros} gives $\max_{1\le j\le N-2}\|q_X^j-(\nu_X^j+\nu_X^{j+1})/2\|_\infty$, at most $4.3\cdot10^{-14}$; the cells $1\le j\le N-2$ are those with two interior endpoints.

All three DR outputs contain exact zero density entries (\cref{tab:zeros}), and the IPM candidate is at most $7.8\cdot10^{-8}$ at these entries (\cref{tab:pathchecks}).
For the diffuse pair the smallest positive density entry is $2.3\cdot10^{-4}$, so the zero entries are separated from the positive entries.
For the peaked and box pairs, positive entries of order $10^{-14}$ occur, so at residual $r_\ell<10^{-12}$ an exact zero cannot be distinguished from a positive value of this order.
These observations are consistent with a discrete minimizer with vacuum entries; they do not prove it.

\begin{table}[htbp]
\centering\small
\begin{tabular}{@{}lrrrr@{}}
\toprule
Pair & \begin{tabular}[b]{@{}r@{}}Zero density\\ entries\end{tabular} & \begin{tabular}[b]{@{}r@{}}Zero midpoint\\ weights\end{tabular} & \begin{tabular}[b]{@{}r@{}}Smallest positive\\ density entry\end{tabular} & Midpoint defect\\
\midrule
Peaked & 61 & 44 & $2.6\cdot10^{-14}$ & $4.3\cdot10^{-14}$\\
Box & 59 & 36 & $7.1\cdot10^{-14}$ & $4.1\cdot10^{-14}$\\
Diffuse & 26 & 20 & $2.3\cdot10^{-4}$ & $5.2\cdot10^{-15}$\\
\bottomrule
\end{tabular}
\caption{Zero entries of the saved DR proximal outputs $X^\ell$ at $M=63$, $N=64$, $\beta=10^6$, $b=1$, and $r_\ell<10^{-12}$.
Density entries are the $63\times63$ entries of $\nu_X^1,\ldots,\nu_X^{N-1}$; midpoint weights are the $64\times63$ entries of $q_X^0,\ldots,q_X^{N-1}$.
The midpoint defect is $\max_{1\le j\le N-2}\|q_X^j-(\nu_X^j+\nu_X^{j+1})/2\|_\infty$, over the cells with two interior endpoints.
A positive entry of order $10^{-14}$ is counted as nonzero, so the zero counts depend on the stopping residual $r_\ell<10^{-12}$.}
\label{tab:zeros}
\end{table}

\Needspace{6\baselineskip}
At a zero midpoint weight the action density $a$ is not differentiable, so the proximal map \eqref{eq:prox} involves its convex subdifferential
\begin{equation}\label{eq:vacsubdiff}
 \partial a(0,0)=\{(\alpha_q,\alpha_m)\in\R^2:\alpha_q+\tfrac12\alpha_m^2\le0\},
\end{equation}
where $(\alpha_q,\alpha_m)$ is a subgradient when $a(q,m)\ge\alpha_q q+\alpha_m m$ for all $(q,m)$.
\Cref{app:drmaps} proves \eqref{eq:vacsubdiff}.
Since $\prox_{ba}(q_{\rm in},v_{\rm in})=(0,0)$ if and only if $(q_{\rm in},v_{\rm in})/b\in\partial a(0,0)$, this set gives the zero branch $2bq_{\rm in}+v_{\rm in}^2\le0$ of \eqref{eq:prox}, which produces the exact zero weights counted in \cref{tab:zeros}.
The local DR analysis at zero weights in \cref{app:lfp} also uses \eqref{eq:vacsubdiff}.

Zero midpoint weights of a feasible finite-action path also satisfy a condition imposed by conservation.
In cell $j$, conservation gives $m^j=m_\perp^j+c\one$ with $m_\perp^j=\C(\nu^{j+1}-\nu^j)$ and $c\in\R$, as in \eqref{eq:geovalue}.
Since $a(0,m)=+\infty$ for $m\ne0$ in \eqref{eq:perspective}, finite action requires $m_i^j=0$ wherever $q_i^j=0$; hence $(m_\perp^j)_i=-c$ on the set $\{i:q_i^j=0\}$, and $m_\perp^j$ must be constant on this set.
A zero set of $q^j$ obtained by replacing by zero all entries below a threshold need not satisfy this condition; when it does not, every momentum that satisfies conservation has infinite action.

\subsection{Density minima and diffusivity transitions}\label{sec:vacuumobservation}
\begin{figure}[p]
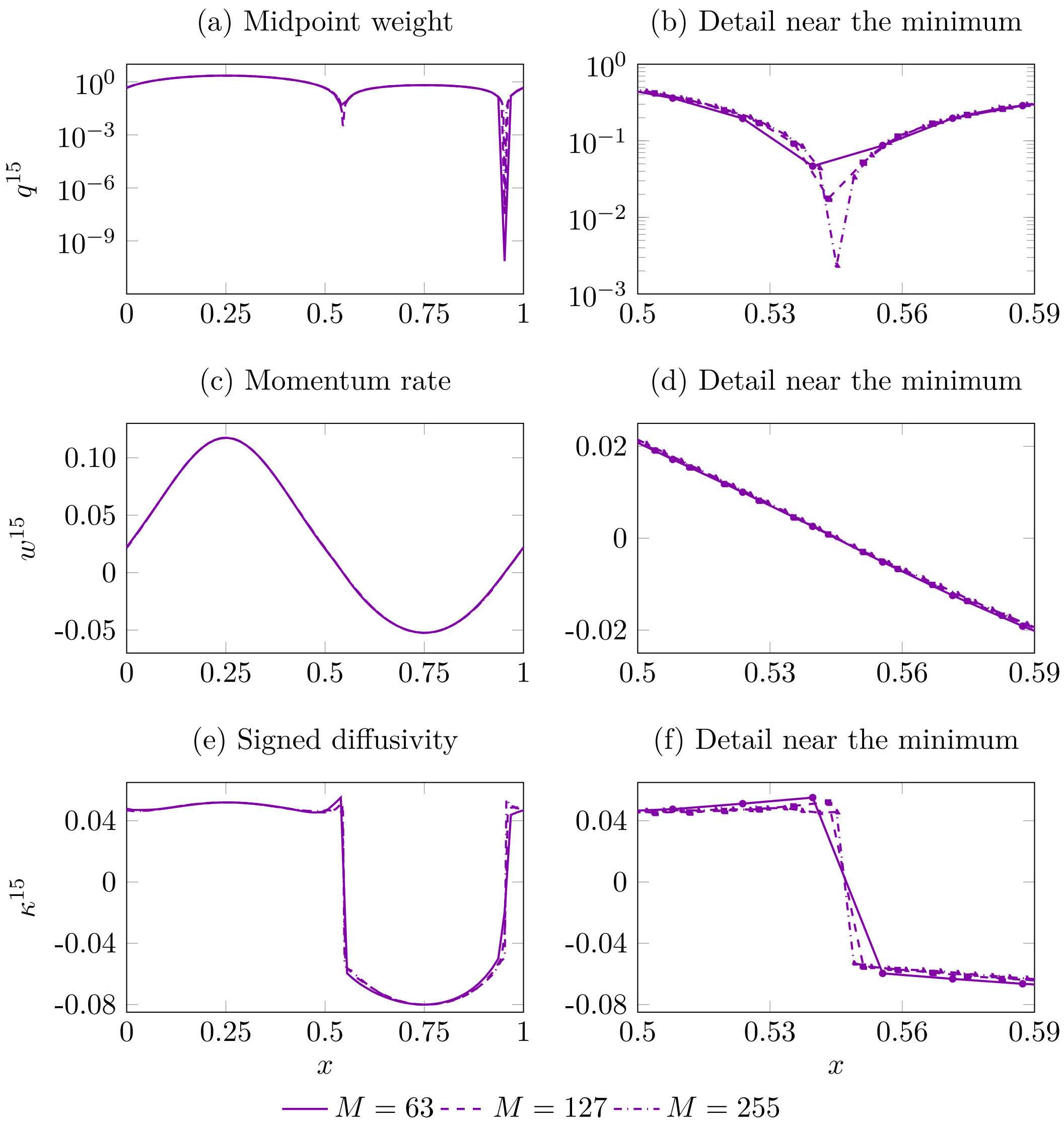

\centering

\pgfplotslegendfromname{legend:diffusivityobservation}
\caption{Density minimum and diffusivity sign change for the peaked geodesic.
The three IPM runs use $N=64$ and the same cell $j=15$, with midpoint time $31/128$.
All curves are at $t=31/128$; solid, dashed and dash-dotted lines denote $M=63,127,255$.
Left: the full periodic interval. Right: the same arrays on $[0.50,0.59]$, which contains the sign change; markers show grid nodes.
The midpoint weight $q^{15}$ is positive and is shown on a logarithmic scale.
The momentum is reconstructed by \eqref{eq:geovalue} from the IPM candidates of \cref{tab:ipmvalidation}, normalized to unit mass; $\kappa^{15}=m^{15}/(hq^{15})$.
The IPM tolerances are $10^{-12}$ for $M=63$ and $10^{-10}$ for $M=127,255$.
Lines connect the computed nodes, with no smoothing or density floor.}
\label{fig:diffusivityobservation}
\end{figure}

For the peaked geodesic, \cref{fig:diffusivityobservation} shows a sign change of the diffusivity between two adjacent nodes, near a density minimum and a zero of the momentum.
This fixed-time-grid study tests whether the profiles are consistent with vanishing density at the transition; it does not establish continuum vacuum.
The test is inequality \eqref{eq:signbound}, which bounds the sum of the midpoint weights at two nodes of opposite diffusivity sign by their momentum-rate difference divided by a lower bound $c_\kappa>0$ of the two diffusivity magnitudes.

We fix $N=64$ and examine cell $j=15$, whose midpoint time is $31/128$.
For each $M\in\{63,127,255\}$, we use the IPM candidate $\nu_{\rm IPM}$ checked in \cref{tab:ipmvalidation}.
We form its midpoint weight $q^{15}$, reconstruct the least-action momentum by \eqref{eq:geovalue}, and calculate $\kappa^{15}=m^{15}/(hq^{15})$.
The right panels of \cref{fig:diffusivityobservation} show the interval $[0.50,0.59]$.
Let $i_M$ be the unique index with $0.50\le x_{i_M}<x_{i_M+1}\le0.59$, $\kappa^{15}_{i_M}>0$ and $\kappa^{15}_{i_M+1}<0$.
Write $q_L,\kappa_L,w_L$ and $q_R,\kappa_R,w_R$ for the values of $q^{15}$, $\kappa^{15}$ and the momentum rate $w^{15}=m^{15}/h$ at $i_M$ and $i_M+1$, and $\delta\kappa_M=\kappa_L-\kappa_R$ for the transition size.
Let $\kappa^{15}_{\rm DR}$ be computed by \eqref{eq:geovalue} and \eqref{eq:kappa} from the action-evaluation path $\bar\nu$ (\cref{sec:geocompare}) of the saved DR run on the same grid, and set $e_\kappa=\max_{i\in\{i_M,i_M+1\}}|\kappa^{15}_{{\rm DR},i}-\kappa^{15}_i|$.
\Cref{tab:diffusivityobservation} records both midpoint weights, both diffusivities, $\delta\kappa_M$, the momentum-rate difference $|w_L-w_R|$, $e_\kappa$ and $e_\kappa/\delta\kappa_M$.
In $[0.50,0.59]$, the sign change of $\kappa^{15}$ lies between the adjacent nodes $i_M,i_M+1$, at distance $1/M$.
On the three grids, $\min(\kappa_L,-\kappa_R)\ge0.045$, $\delta\kappa_M$ decreases from $0.115$ to $0.099$, and $q_L$ decreases from $4.7\cdot10^{-2}$ to $2.4\cdot10^{-3}$ (\cref{tab:diffusivityobservation}).
These three grids do not determine whether $\min(\kappa_L,-\kappa_R)$ stays above a positive constant as $M\to\infty$.

\begin{table}[htbp]
\centering\small
\begin{tabular}{@{}rrrrrrrrr@{}}
\toprule
$M$ & $q_L$ & $q_R$ & $\kappa_L$ & $\kappa_R$ & $\delta\kappa_M$ & $|w_L-w_R|$ & $e_\kappa$ & $e_\kappa/\delta\kappa_M$\\
\midrule
63 & $4.67\cdot10^{-2}$ & $8.68\cdot10^{-2}$ & 0.05513 & -0.05965 & 0.11478 & $7.75\cdot10^{-3}$ & $7.50\cdot10^{-8}$ & $6.53\cdot10^{-7}$\\
127 & $1.78\cdot10^{-2}$ & $5.26\cdot10^{-2}$ & 0.05226 & -0.05490 & 0.10716 & $3.82\cdot10^{-3}$ & $4.70\cdot10^{-7}$ & $4.38\cdot10^{-6}$\\
255 & $2.36\cdot10^{-3}$ & $3.36\cdot10^{-2}$ & 0.04574 & -0.05303 & 0.09876 & $1.89\cdot10^{-3}$ & $2.67\cdot10^{-6}$ & $2.71\cdot10^{-5}$\\
\bottomrule
\end{tabular}
\caption{Two adjacent nodes at the positive-to-negative diffusivity transition of the saved peaked IPM candidate, with $N=64$, $j=15$, and $w^{15}=m^{15}/h=q^{15}\kappa^{15}$. The subscripts $L,R$ denote the nodes $i_M,i_M+1$ defined in the text, with $0.50\le x_{i_M}<x_{i_M+1}\le0.59$, and $\delta\kappa_M=\kappa_L-\kappa_R$. The column $e_\kappa$ is the maximum DR--IPM diffusivity difference at these two nodes, defined in the text. The DR runs use $\beta=10^6$ and $b=1$ and stop when $r_\ell\le10^{-12}$, after 9707, 22770 and 42607 iterations for $M=63,127,255$. The same saved candidates and their interval transition checks are used as in the original comparison; entries are displayed rounded values, not error bounds for the exact minimizing path. The complete positions and interval values remain in the source record.}
\label{tab:diffusivityobservation}
\end{table}

On each grid, $e_\kappa/\delta\kappa_M\le2.71\cdot10^{-5}$ (\cref{tab:diffusivityobservation}).
Since also $e_\kappa<\min(\kappa_L,-\kappa_R)$, the DR diffusivities at the two nodes have the same signs as the IPM diffusivities.

\begin{figure}[p]
\centering
\includegraphics[width=\linewidth]{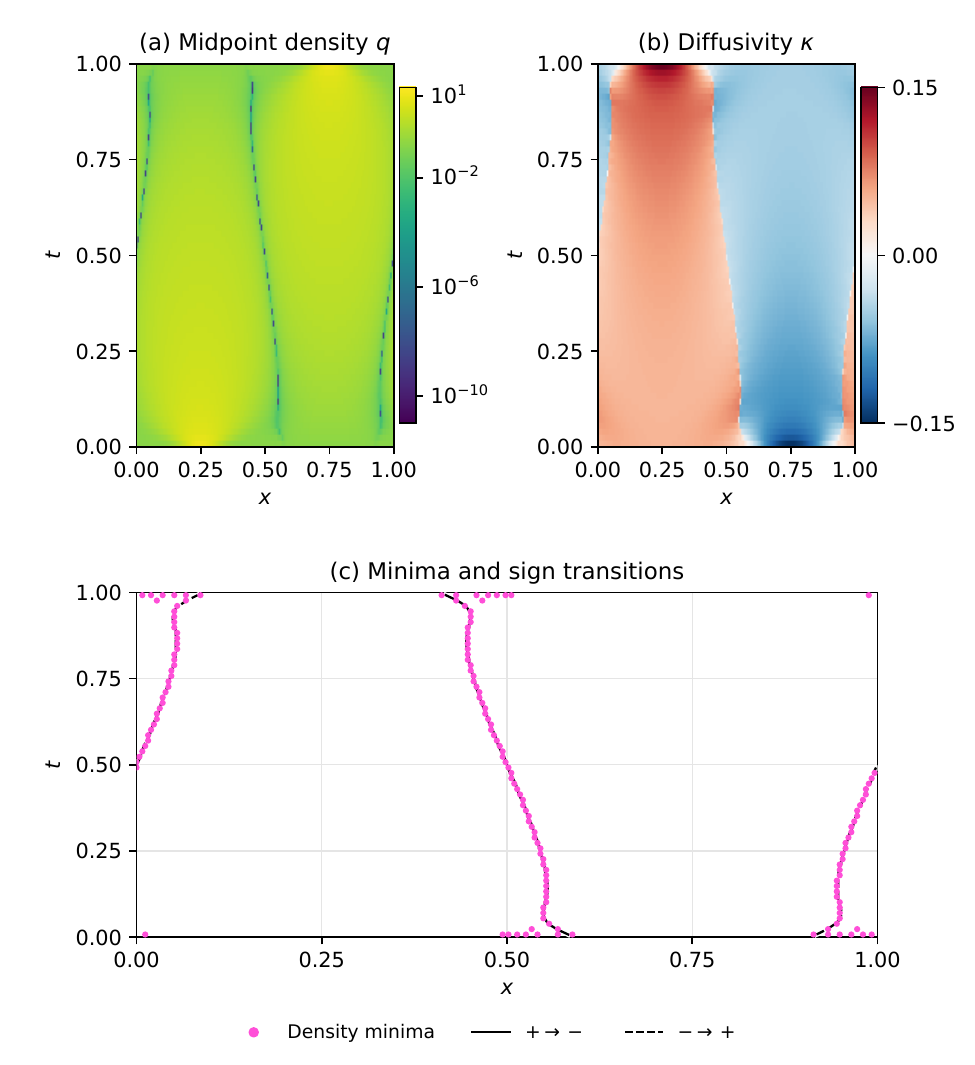}
\caption{Space--time maps for the peaked IPM candidate $\nu_{\rm IPM}$ at $M=255$, $N=64$, validated in \cref{tab:ipmvalidation}.
\textbf{(a)} Midpoint weight $q^j=(\nu^j+\nu^{j+1})/2$, with logarithmic colors from $10^{-11}$ to 20.
\textbf{(b)} Signed diffusivity $\kappa^j=m^j/(hq^j)$, with linear colors from $-0.15$ to $0.15$.
Both maps use cell midpoint times $t=(j+1/2)/N$ and display unsmoothed cell colors.
Each density row is normalized to unit mass; momentum is reconstructed by \eqref{eq:geovalue}. No additional positive mixture is used.
\textbf{(c)} Pink points are all discrete local minima of $q^j$, using periodic neighbors.
Black curves join the zeros of the piecewise linear interpolant of $\kappa^j$ between adjacent nodes where $\kappa^j$ changes sign: solid for positive to negative, dashed for negative to positive.
Curves are split at the periodic boundary.
The IPM densities are strictly positive; these points are minima, not exact density zeros. Only $M=255$ is displayed. On each grid $M=63,127,255$, consider the sign transitions at midpoint times $0.1\le t\le0.9$, 104 per grid. The distance $d_{\T}$ from each transition to the nearest discrete local minimum of $q^j$ in the same cell is at most $0.52/M$.}
\label{fig:geocolormaps}
\end{figure}

\Needspace{8\baselineskip}
\paragraph{Space--time location of the transitions.}
\Cref{fig:geocolormaps} compares the midpoint weight and diffusivity over the full computed path at $M=255$, $N=64$.
The third panel compares all discrete local density minima with the interpolated sign transitions.
On all three grids, the sampled sign transitions lie within $0.52/M$ of a discrete local minimum of $q^j$ (\cref{fig:geocolormaps}).
The interpolated curves locate transitions between opposite signs; they do not determine whether the limiting transition is continuous or discontinuous.

\paragraph{Continuous momentum and discontinuous diffusivity.}
In a fixed cell $j$, write $q,w,\kappa$ for $q^j$, the momentum rate $w^j=m^j/h$ and $\kappa^j$; then $w=q\kappa$ where $q>0$.
For continuous spatial profiles at a fixed time, use $w=\rho\kappa$ where $\rho>0$.
If $\rho(x)\to0$ at a transition $x_0$ and $\kappa$ has finite limits on both sides, then $w(x)\to0$ from both sides, even when the two diffusivity limits differ.
The momentum therefore extends continuously through $x_0$ with value zero.

\paragraph{A criterion for vanishing density.}
For the discrete profiles defined above, suppose the diffusivity values remain separated from zero across a transition.
Suppose two nodes $x_+,x_-$ with $q(x_+),q(x_-)>0$ satisfy $\kappa(x_+)\ge c_\kappa>0$ and $\kappa(x_-)\le-c_\kappa$.
Then
\begin{equation}\label{eq:signbound}
 q(x_+)+q(x_-)\le\frac{|w(x_+)-w(x_-)|}{c_\kappa}.
\end{equation}
Thus both midpoint weights $q(x_+),q(x_-)$ tend to zero if the nodes approach one point, the opposite diffusivity values stay uniformly separated from zero, and the momentum difference tends to zero.
If the density profiles also converge locally uniformly to a continuous limit $\rho$, that limit vanishes at the transition point.
\Cref{lem:vacuumcriterion} states the assumptions precisely.

At a fixed time step $h$, conservation and nonnegative unit-mass densities provide a spatial continuity estimate for $w$ that is uniform in $M$; \cref{app:vacuumcriterion} proves it.
With \eqref{eq:signbound}, this estimate gives the adjacent-node bound $q(x_+)+q(x_-)\le1/(c_\kappa h\sqrt{3M})$ for two nodes at distance $1/M$.
For $h=1/64$, $M=255$ and $c_\kappa=\min(\kappa_L,-\kappa_R)=0.0457$ (\cref{tab:diffusivityobservation}), the adjacent-node bound equals $50.6$, whereas $q_L+q_R\le0.134$ on the three grids; the bound does not explain the tabulated weights.
The three grids with $h=1/64$ do not determine whether the hypothesis of \eqref{eq:signbound} holds at adjacent nodes with a constant $c_\kappa>0$ independent of $M$ as $M\to\infty$.
Identifying a zero of the continuous geodesic requires local uniform
convergence of the density profiles, stronger than the weak convergence
provided by \cref{geo:thm:joint-grid}.

\Needspace{4\baselineskip}
\paragraph{Least action at prescribed momentum.}
For $\rho>0$, the local action per unit time is
\[
 a(\rho,w)=\frac{w^2}{2\rho}=\frac12\rho\kappa^2,
 \qquad w=\rho\kappa.
\]
At fixed nonzero momentum $w$, decreasing $\rho$ increases the cost.
Fix $w$ and minimize $\int_\T w^2/(2\rho)\,dx$ over densities $\rho\ge0$ with $\int_\T\rho\,dx=1$; the conservation and endpoint constraints of the geodesic problem are not imposed.
Assume $0<c_w:=\int_\T|w|\,dx<\infty$.
For every such density (\cref{app:fixed-momentum}),
\[
 \int_\T\frac{w^2}{2\rho}\,dx\ge\frac12c_w^2.
\]
With $a(0,0)=0$ as in \eqref{eq:perspective}, equality holds for $\rho=|w|/c_w$; then $\kappa=c_w\operatorname{sgn}(w)$ wherever $w\ne0$.
If $w$ is continuous and changes sign at an isolated zero $x_0$, the minimizer $\rho=|w|/c_w$ vanishes at $x_0$, and $\kappa=c_w\operatorname{sgn}(w)$ jumps between $-c_w$ and $c_w$ there.
The action $\frac12\int_\T\rho\kappa^2\,dx$ contains no derivative of $\kappa$, so this jump adds no cost.
In the geodesic problem, conservation and the endpoint conditions also couple $\rho$ and $w$.
The calculation therefore only interprets the computed profiles; it proves neither a jump of $\kappa$ for a limiting minimizing geodesic nor the number of vacuum points.

\begin{figure}[htbp]
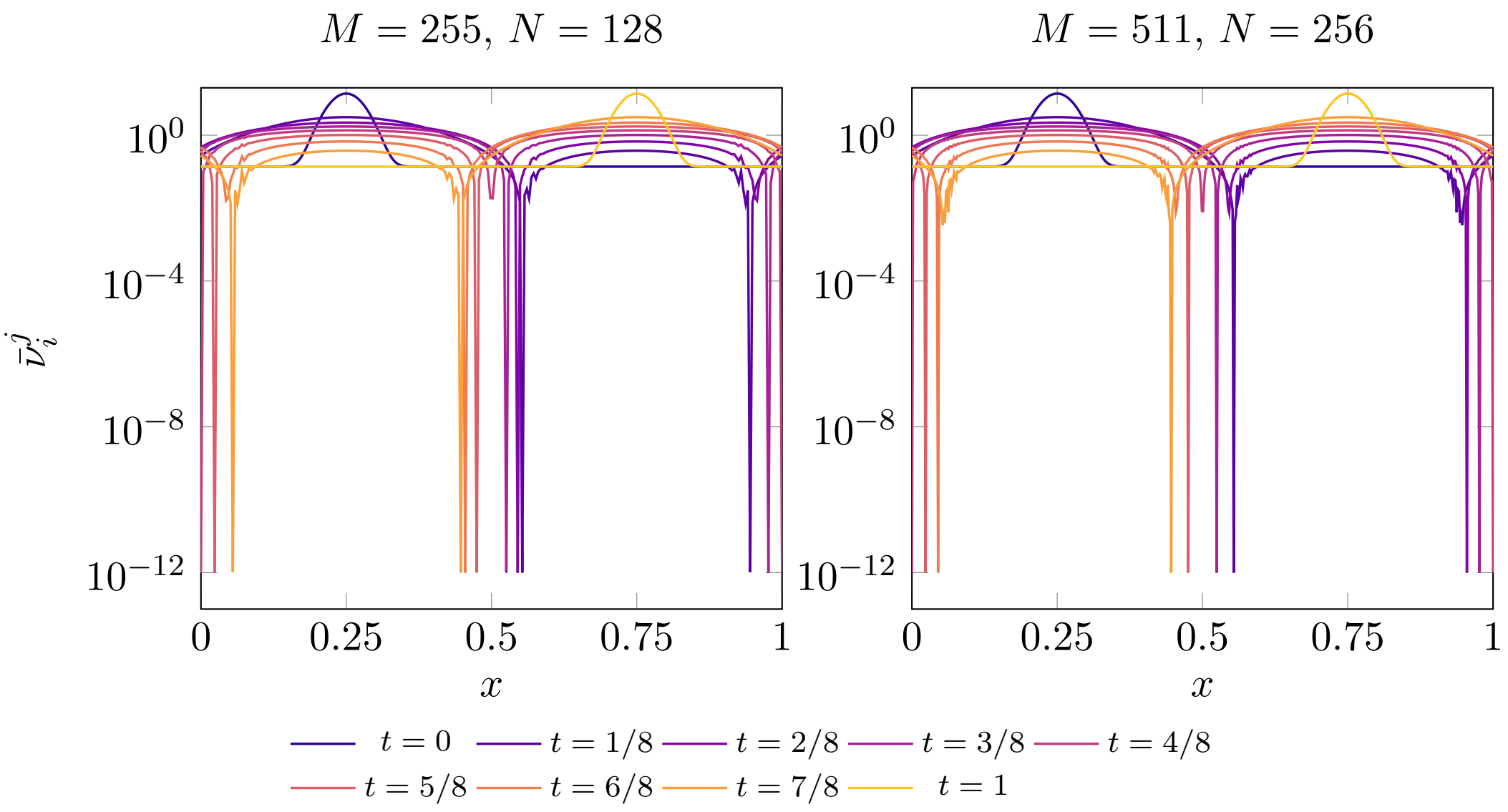

\centering


\pgfplotslegendfromname{legend:geofine}
\caption{Peaked geodesic densities under joint space--time refinement, on a logarithmic scale at nine computed time levels.
Both space and time are refined from the base grid $(M,N)=(63,64)$.
The same colors identify the same path times in both panels.
The plotted densities $\bar\nu$ include the uniform mixture of weight $10^{-12}$ used for action evaluation.
\Cref{tab:geofineruns} lists the two DR runs.
All spatial nodes are plotted, with periodic closure and no spatial smoothing.}
\label{fig:geofine}
\end{figure}

\paragraph{Joint space--time refinement of the peaked geodesic.}
\Cref{fig:geofine} shows the peaked densities at $(M,N)=(255,128)$ and $(511,256)$ and tests whether the interior density minima observed at $(M,N)=(63,64)$ persist.
\begin{table}[htbp]
\centering\small
\begin{tabular}{@{}rrrrr@{}}
\toprule
$M$ & $N$ & Iterations & $r_\ell$ & $(V-D)/D$\\
\midrule
255 & 128 & 60000 & $2.2\cdot10^{-11}$ & $1.2\cdot10^{-10}$\\
511 & 256 & 60000 & $5.9\cdot10^{-9}$ & $5.5\cdot10^{-8}$\\
\bottomrule
\end{tabular}
\caption{DR runs used for the peaked profiles in \cref{fig:geofine}. Both use $b=1$ and $\beta=10^6(N/64)^2$; this scaling gives $\sqrt\beta\,h=10^3/64$ for every $N$, so the ratio $v/q=\sqrt\beta\,h\kappa$ of scaled momentum to midpoint weight does not depend on $N$ at a given diffusivity. Both runs stop at the iteration cap. The column $r_\ell$ is the DR residual of \cref{sec:geodata}. $V$ is the action of the action-evaluation path of \cref{sec:geocompare}, and $D=D_G(\chi)$ is the dual value \eqref{eq:geodual} at a multiplier $\chi$ obtained by numerical maximization started from the final DR iterate; this gap uses no IPM reference value. Both values are computed in floating point, so $(V-D)/D$ is an estimate, not an interval certificate or a density-error bound.}
\label{tab:geofineruns}
\end{table}

The adjacent-node bound $1/(c_\kappa h\sqrt{3M})$ tends to zero as $M\to\infty$ only if $h\sqrt M\to\infty$.
Both runs have $M+1=2N$ and $h=1/N$, so $h\sqrt M=\sqrt{2N-1}/N$ decreases from $0.125$ ($N=128$) to $0.088$ ($N=256$) and tends to zero along this sequence.
The adjacent-node bound therefore does not tend to zero along these runs.
The plotted arrays are the unit-mass action-evaluation densities $\bar\nu$, including the $10^{-12}$ uniform mixture from \cref{sec:geocompare}.
They therefore cannot identify proximal zero entries.
At each of the seven interior plotted times $t=1/8,\ldots,7/8$ and on both grids, $\min_i\bar\nu_i^j\le1.01\cdot10^{-12}$: the minima reach the mixture weight $10^{-12}$.
The interior minima of the peaked path at $(M,N)=(63,64)$ therefore persist at these two refinements; this does not verify the hypotheses of \cref{lem:vacuumcriterion}.

\FloatBarrier
\subsection{Nonnegative diffusivity: feasibility and action cost}\label{sec:nonnegative-diffusivity}

We now add $m_i^j\ge0$ to \eqref{eq:geodesic} and denote the minimum by $E^+_{M,N}$.
This is equivalent to $\kappa_i^j\ge0$ wherever $q_i^j>0$; finite action still requires $m_i^j=0$ at $q_i^j=0$.
The feasible set is smaller, so $E_{M,N}\le E^+_{M,N}$.

Any two strictly positive unit-mass densities on $\T$ are feasible endpoints for the periodic problem.
\Cref{rem:martingale} distinguishes this periodic problem from a martingale problem on $\R$.
For a positive density path, set $m_\perp^j=\C(\nu^{j+1}-\nu^j)$ as in \eqref{eq:geovalue}.
Since $L\one=0$, the least-action nonnegative momentum is
\begin{equation}\label{eq:positive-momentum}
 m_+^j=m_\perp^j-\min_i(m_\perp^j)_i\one.
\end{equation}
The positive linear interpolation of the endpoints with \eqref{eq:positive-momentum} gives a finite-action competitor.
For a given density path, $m_+^j$ has a zero entry in each cell.
Under the open constraint $m>0$ the infimum is therefore not attained, so we impose the closed constraint $m\ge0$.
\Cref{app:positive-dual} proves these statements and gives the corresponding dual bound.

Both solvers require only a local change.
For DR, define $a_+(q,v)=a(q,v)+\iota_{\{v\ge0\}}(q,v)$ for the scaled momentum $v=\sqrt\beta m$.
Then, in the notation of \eqref{eq:prox},
\begin{equation}\label{eq:positive-prox}
 \prox_{b a_+}(q_{\rm in},v_{\rm in})
 =\prox_{b a}(q_{\rm in},\max\{v_{\rm in},0\}).
\end{equation}
\Cref{app:positive-prox} proves \eqref{eq:positive-prox}; the affine projection is unchanged.
For IPM, add $w_i^j\ge0$ to the conic formulation in \cref{sec:geocone}.
The implementation stores the array $-w^j$, so the implemented constraint is $-w_i^j\le0$.

We use the diffuse endpoints defined in \cref{sec:geodata}.
A nodal density $\nu$ defines the law on $\R$ with mass $\nu_i/M$ at $x_i\in[0,1)$.
For two probability laws $\nu,\mu$ on $\R$, $\nu\le_{\rm cx}\mu$ (convex order) means that every convex function on $\R$ has a mean under $\nu$ at most equal to its mean under $\mu$.
For the diffuse pair, $\nu^0\not\le_{\rm cx}\nu^N$ (\cref{rem:martingale}), whereas the periodic problem with $m\ge0$ is feasible.
\Cref{tab:positive-cost} gives interval bounds for both optimal actions on two grids.
The ratio $E^+_{M,N}/E_{M,N}$ equals $17.39$ for $(M,N)=(63,64)$ and $17.43$ for $(127,128)$ (\cref{tab:positive-cost}).

\begin{table}[htbp]
\centering\small
\begin{tabular}{@{}rrrrr@{}}
\toprule
$M$ & $N$ & $E_{M,N}$ & $E^+_{M,N}$ & $E^+_{M,N}/E_{M,N}$\\
\midrule
\multicolumn{5}{@{}l}{Diffuse pair}\\
63 & 64 & $[6.32113248,6.32113249]\cdot10^{-5}$
 & $[1.09945745,1.09945746]\cdot10^{-3}$ & $[17.39336,17.39337]$\\
127 & 128 & $[6.32343401,6.32343402]\cdot10^{-5}$
 & $[1.10208203,1.10208204]\cdot10^{-3}$ & $[17.42853,17.42854]$\\
\multicolumn{5}{@{}l}{Convex-order pair}\\
63 & 64 & $\OrderConvexSignedCoarse$ & $\OrderConvexPositiveCoarse$ & $\OrderConvexRatioCoarse$\\
127 & 128 & $\OrderConvexSignedFine$ & $\OrderConvexPositiveFine$ & $\OrderConvexRatioFine$\\
\multicolumn{5}{@{}l}{Reversed pair}\\
63 & 64 & $\OrderReversedSignedCoarse$ & $\OrderReversedPositiveCoarse$ & $\OrderReversedRatioCoarse$\\
127 & 128 & $\OrderReversedSignedFine$ & $\OrderReversedPositiveFine$ & $\OrderReversedRatioFine$\\
\bottomrule
\end{tabular}
\caption{Interval enclosures of the discrete optimal actions with and without nonnegative momentum, for the diffuse pair of \cref{sec:geodata}, the convex-order pair (Gaussian profiles of standard deviations $0.10$ for $\nu^0$ and $0.20$ for $\nu^N$, centered at $(M-1)/(2M)$) and the reversed pair (the same endpoints exchanged).
In the two solves of each row, the endpoints are the exact unit-mass normalizations of the same stored binary64 arrays.
All displayed interval endpoints are rounded outward.
The bounds use 40-decimal-digit interval arithmetic and the exact Fourier symbol; they do not bound discretization error.}
\label{tab:positive-cost}
\end{table}

\begin{table}[htbp]
\centering\small
\begin{tabular}{@{}lrrrrr@{}}
\toprule
Momentum & DR $\beta$ & DR iterations & IPM iterations & DR $r_\ell$ & Density difference\\
\midrule
\multicolumn{6}{@{}l}{Diffuse pair}\\
Signed & $10^6$ & 12800 & 21 & $9.41\cdot10^{-12}$ & $1.22\cdot10^{-7}$\\
Nonnegative & $10^7$ & 105200 & 28 & $9.87\cdot10^{-11}$ & $2.29\cdot10^{-7}$\\
\multicolumn{6}{@{}l}{Convex-order pair}\\
Signed & $10^6$ & \OrderConvexSignedDRIterations & \OrderConvexSignedIPMIterations & $\OrderConvexSignedDRResidual$ & $\OrderConvexSignedDensityDifference$\\
Nonnegative & $10^7$ & \OrderConvexDRIterations & \OrderConvexIPMIterations & $\OrderConvexDRResidual$ & $\OrderConvexDensityDifference$\\
\bottomrule
\end{tabular}
\caption{Solver agreement for the diffuse pair and the convex-order pair at $M=63$, $N=64$. DR uses $b=1$ and stops when $r_\ell\le10^{-11}$ (signed) or $r_\ell\le10^{-10}$ (nonnegative); IPM tolerances are $10^{-11}$; IPM uses the baseline settings of \cref{alg:positive-experiment} for the signed run of the diffuse pair and the modified settings for the other three rows. Both signed rows use $(\beta,b)=(10^6,1)$, the values of the geodesic DR tests (\cref{sec:dr-choices}). The density difference is $\max_j\|\bar\nu_{\rm DR}^j-\bar\nu_{\rm IPM}^j\|_{1,M}$ for the action-evaluation paths defined in the text. For nonnegative momentum, $\max_{j,i}|\kappa_{{\rm DR},i}^j-\kappa_{{\rm IPM},i}^j|$ is at most $2.66\cdot10^{-8}$ for the diffuse pair and $\OrderConvexDiffusivityDifference$ for the convex-order pair. Values of $r_\ell$ are rounded upward; density differences are rounded to three significant digits.}
\label{tab:positive-agreement}
\end{table}

\begin{figure}[p]
\centering
\includegraphics[width=.98\textwidth]{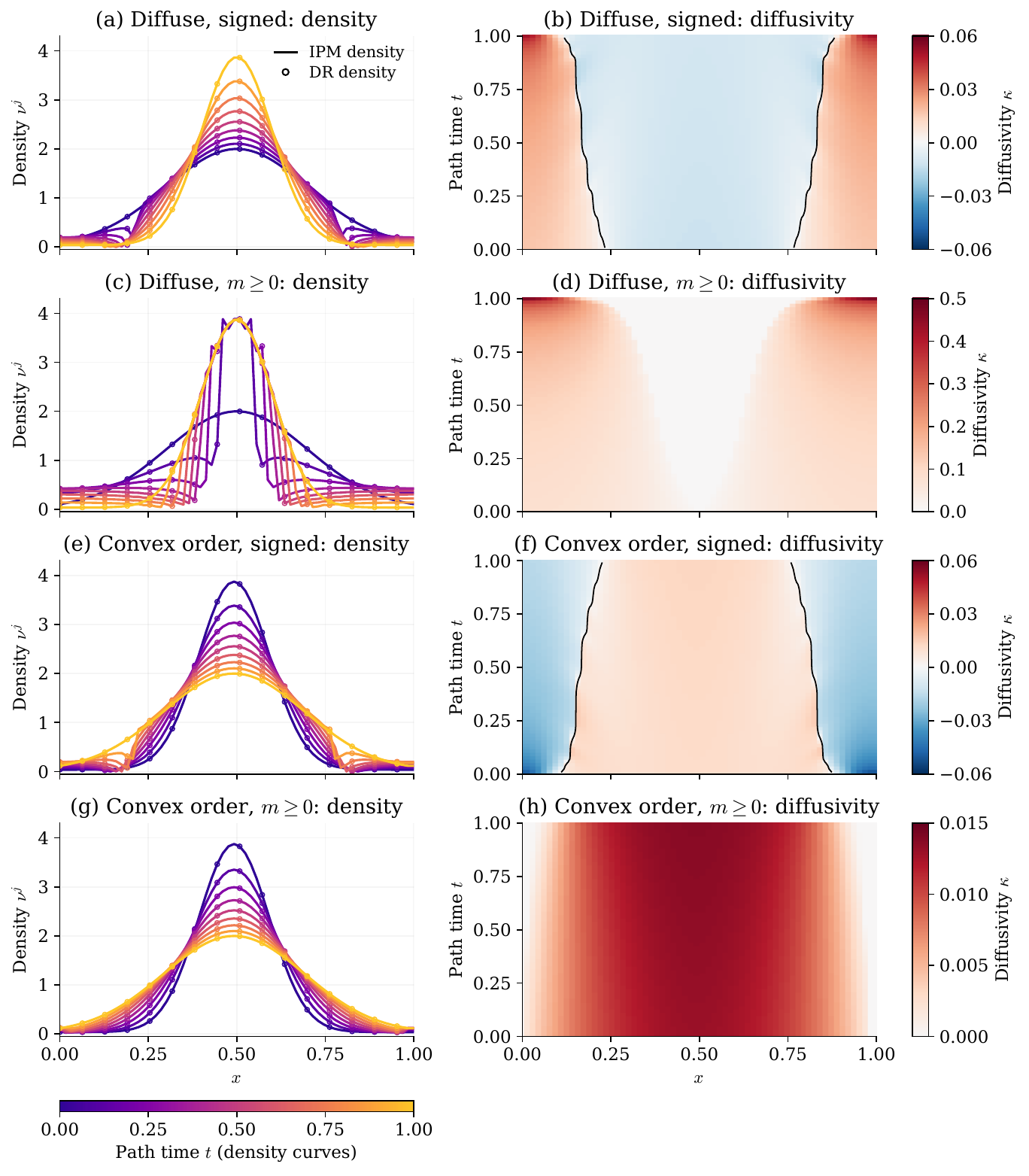}
\caption{Solutions at $M=63$, $N=64$: diffuse pair with signed momentum (first row) and with nonnegative momentum (second row), and convex-order pair with signed momentum (third row) and with nonnegative momentum (fourth row).
Left: IPM density curves at $t=j/N$, $j=0,8,\ldots,64$, with independent DR values marked at every fourth spatial node; the color bar below the left column records path time.
Right: space--time maps of the IPM diffusivity $\kappa^j=m^j/(h q^j)$ at all 64 cell times and all 63 spatial nodes.
Blue denotes negative diffusivity, red positive diffusivity, and white zero; each row has its own color bar with a linear value scale, and the two signed maps use the same scale.
Black curves in the two signed maps interpolate zero crossings between opposite-sign samples.
Map colors are constant on display cells centered at the spatial nodes and cell midpoint times, with periodic closure and no smoothing.
The density panels share a vertical scale.
All densities and momenta are the action-evaluation arrays specified in the text; DR stops when $r_\ell\le10^{-11}$ (signed momentum) or $r_\ell\le10^{-10}$ ($m\ge0$), IPM tolerances are $10^{-11}$, and solver agreement is quantified in \cref{tab:positive-agreement}.}
\label{fig:positive-diffusivity}
\end{figure}

For action evaluation, DR uses the simplex-projected density with a $10^{-12}$ uniform mixture.
In this section, the IPM candidate, which is also the IPM action-evaluation path, is selected by \cref{alg:positive-experiment} instead of \cref{alg:geo-candidates}: each raw interior row is clipped at zero and normalized to unit mass; among the mixtures of this path with the uniform density, with mixture weights $0,10^{-14},10^{-12},10^{-10},10^{-8}$, those with positive midpoint weights are evaluated and the one of least action is kept.
At $(M,N)=(63,64)$, the mixture weight selected for the IPM candidate is $\OrderConvexMixture$ for the nonnegative solution of the convex-order pair and $0$ for the three other solutions of \cref{tab:positive-agreement}.
The action of each path is evaluated with the momentum recomputed by \eqref{eq:geovalue} for signed momentum and by \eqref{eq:positive-momentum} for $m\ge0$.
\Cref{tab:positive-agreement} quantifies the DR--IPM density differences with and without nonnegative momentum.
\Cref{fig:positive-diffusivity} shows density profiles and space--time maps of the diffusivity.
For the diffuse pair, the IPM diffusivities of the signed solution satisfy $\max_{j,i}|\kappa_i^j|=0.0513$ (first row of \cref{fig:positive-diffusivity}).
Those of the nonnegative solution satisfy $\max_{j,i}\kappa_i^j=0.468$, attained at $x=0$ in the last cell, and $\max_j\kappa_i^j>0.0513$ exactly at the nodes $x_i\notin[1/3,2/3]$ (second row of \cref{fig:positive-diffusivity}).
Thus, for the diffuse pair at $(M,N)=(63,64)$, the constraint $m\ge0$ replaces the negative diffusivities of the signed solution by nonnegative values, produces diffusivities up to $0.468$ outside $[1/3,2/3]$, and multiplies the minimum action by $17.39$ (\cref{tab:positive-cost}).
This comparison concerns the fixed discrete problems; agreement of the solvers does not establish convergence to a continuous geodesic.

To test whether $E^+_{M,N}/E_{M,N}$ depends on the direction of the convex order between the endpoints, we add two pairs.
For the convex-order pair, $\nu^0$ and $\nu^N$ use $\mathfrak f(x)=\exp(-(x-\bar x)^2/(2\cdot0.10^2))$ and $\mathfrak f(x)=\exp(-(x-\bar x)^2/(2\cdot0.20^2))$, respectively, normalized as in \cref{sec:geodata}, with $\bar x=(M-1)/(2M)$ the center of the nodes $x_i$.
The reversed pair exchanges these two endpoints; it differs from the diffuse pair only by the center, $\bar x$ instead of $1/2$.
For the convex-order pair, $\nu^0\le_{\rm cx}\nu^N$.
Indeed, both endpoints are symmetric about $\bar x$ and therefore have the same mean, and on the nodes a convex function coincides with an affine function plus a nonnegative combination of the functions $\max\{x-x_l,0\}$, $0\le l<M$.
With equal masses and means, $\nu^0\le_{\rm cx}\nu^N$ therefore reduces to the $M$ inequalities $\avg{\max\{x-x_l,0\}\nu^0}\le\avg{\max\{x-x_l,0\}\nu^N}$.
We verified these inequalities and the equality of the means in exact rational arithmetic, for the exact unit-mass normalizations of the stored binary64 arrays on both grids.
For the reversed pair, $\nu^N\le_{\rm cx}\nu^0$, and $\nu^0\not\le_{\rm cx}\nu^N$ since $\avg{(x-\bar x)^2\nu^N}<\avg{(x-\bar x)^2\nu^0}$.
If $(\nu^j,m^j)$ is feasible in \eqref{eq:geodesic} for one of these two pairs, then the path with density $\nu^{N-j}$ and momentum $-m^{N-1-j}$ at index $j$ is feasible for the other pair and has the same action, since $L(-m^{N-1-j})=\nu^{N-1-j}-\nu^{N-j}$, its midpoint weight at index $j$ is $q^{N-1-j}$, and $a(q,-m)=a(q,m)$.
Hence $E_{M,N}$ is the same for both pairs, whereas $E^+_{M,N}$ can differ, since this map exchanges the constraints $m\ge0$ and $m\le0$.
At $(M,N)=(127,128)$, the ratio $E^+_{M,N}/E_{M,N}$ equals $\OrderConvexRatioShort$ for the convex-order pair and $\OrderReversedRatioShort$ for the reversed pair (\cref{tab:positive-cost}); the latter agrees with the ratio of the diffuse pair to three significant digits.
For these two laws, $E^+_{M,N}/E_{M,N}$ is therefore smaller when $\nu^0\le_{\rm cx}\nu^N$ than when $\nu^N\le_{\rm cx}\nu^0$.
These two pairs do not establish a general relation between the direction of the convex order and $E^+_{M,N}/E_{M,N}$.

When $\nu^N\ne\nu^0$, every minimizer of \eqref{eq:geodesic} with a positive density path satisfies $\min_{j,i}\kappa_i^j<0$, even if $\nu^0\le_{\rm cx}\nu^N$.
Indeed, the momentum of such a minimizer is the signed momentum of least action \eqref{eq:geovalue} of its density path, which satisfies $\avg{m^j/q^j}=0$ by the choice of $c_\star^j$, that is $\sum_i\kappa_i^j=0$ by \eqref{eq:kappa}.
If $\nu^{j+1}\ne\nu^j$, which holds for some $j$ since $\nu^N\ne\nu^0$, then $m^j\ne0$ because $Lm^j=\nu^{j+1}-\nu^j$; hence $\min_i\kappa_i^j<0<\max_i\kappa_i^j$.
For the convex-order pair at $(M,N)=(63,64)$, the IPM diffusivities of the signed solution satisfy $\min_{j,i}\kappa_i^j=\OrderConvexSignedKappaMin$, so this solution is not feasible for $m\ge0$.
This minimum is attained in the first cell at $x=0$ and $x=(M-1)/M$, the two nodes farthest from $\bar x$.
In every cell $j$, no entry of $\kappa^j$ vanishes, and $|x_i-\bar x|>|x_l-\bar x|$ whenever $\kappa_i^j<0<\kappa_l^j$ (third row of \cref{fig:positive-diffusivity}).
Neither the interval action enclosures of \cref{tab:positive-cost} nor the density difference of \cref{tab:positive-agreement} bounds the error of $\kappa$, since an error $\delta m$ in $m$ changes $\kappa$ by $\delta m/(hq)$, which has no bound as $q\to0$ (\cref{app:algorithm-positive}).
The signs of the IPM diffusivities of the signed solution are therefore checked against the independent DR run: the DR and IPM diffusivities of this solution have the same sign at every node and in every cell.
The IPM diffusivities of the nonnegative solution satisfy $\max_{j,i}\kappa_i^j=\OrderConvexKappaMax$, attained at $x=\bar x$ in the last cell (fourth row of \cref{fig:positive-diffusivity}).
Conversely, the periodic problem with $m\ge0$ does not require $\nu^0\le_{\rm cx}\nu^N$: it is feasible for any positive endpoints by \eqref{eq:positive-momentum}.
On $\R$, by contrast, under the assumptions of \cref{rem:martingale}, a nonnegative diffusivity requires $\nu^0\le_{\rm cx}\nu^N$.

\begin{remark}[Relation to a martingale constraint]\label{rem:martingale}
If $\kappa\ge0$ and $\rho$ solves $\partial_t\rho=\partial_{xx}(\rho\kappa)$ on $\R$ under the integrability assumptions of \cite{HT19}, then $\rho(t,\cdot)$ is the law of $X_t$ for a continuous martingale $(X_t)_{0\le t\le1}$ with quadratic variation $[X]_t=\int_0^t2\kappa(s,X_s)\,ds$.
Then $\rho(0,\cdot)\le_{\rm cx}\rho(1,\cdot)$.
Here the endpoints are prescribed on $\T$.
A law on $\T$ does not determine a law on $\R$: paths on $\T$ with the same endpoint laws can give different laws of the unwrapped position on $\R$.
The second moment about $1/2$, $\avg{(x-1/2)^2\nu}$, decreases from $\DiffuseInitialMoment$ to $\DiffuseFinalMoment$ between the endpoints (\cref{tab:diffuse-moments}).
Since $x\mapsto(x-1/2)^2$ is convex, this moment would be nondecreasing along a martingale.
Thus these same nodal laws, viewed on $[0,1)\subset\R$, cannot be joined by a real-valued martingale.
The periodic nonnegative-diffusivity problem is nevertheless feasible by \eqref{eq:positive-momentum}.
Moreover, $L$ is not a Markov generator, so $m\ge0$ alone does not define a discrete martingale.
The ratio $E^+_{M,N}/E_{M,N}$ in \cref{tab:positive-cost} therefore measures the cost of the constraint $m\ge0$ in the periodic discrete problem, not the cost of a martingale constraint.
\end{remark}

\begin{table}[htbp]
\centering\small
\begin{tabular}{@{}rrr@{}}
\toprule
$M$ & $\avg{(x-1/2)^2\nu^0}$ & $\avg{(x-1/2)^2\nu^N}$\\
\midrule
63 & $\DiffuseInitialMoment$ & $\DiffuseFinalMoment$\\
\bottomrule
\end{tabular}
\caption{Second moment about $1/2$ of the diffuse endpoints of \cref{sec:geodata}, represented at $x_i=i/M$ in $[0,1)\subset\R$. The quadrature moment decreases between the endpoints, excluding a real-valued martingale with these nodal laws. This test does not obstruct the periodic nonnegative-diffusivity problem.}
\label{tab:diffuse-moments}
\end{table}

\section{Entropic evolution: accuracy and density minima}\label{sec:flowtests}

We compare DR and IPM for repeated entropic steps, using primal--dual bounds to assess optimization accuracy and recorded runtimes to assess computational cost.
We then ask how far the minimum $\min_i\nu_i^k$ decreases along the entropic trajectory, how this decrease depends on the initial profile and on the time grid, and whether the error bound \eqref{eq:pointwise} or the computed dual multiplier controls the smallest computed entries.

\subsection{Baseline data and validation}\label{sec:flowdata}
We use the family from \cite[Figure~1]{MRSSdlss}:
\begin{equation}\label{eq:flowdata}
 \nu_i^0=Z^{-1}\left(\sqrt\varepsilon+
       \left[\frac{1+\cos(2\pi x_i)}2\right]^p\right)^2,
 \qquad p\in\{2,4,8,16\},\quad\varepsilon=10^{-3}.
\end{equation}
The constant $Z$ gives $\avg{\nu^0}=1$.
For a positive vector $\nu$, its density contrast is $\max_i\nu_i/\min_i\nu_i$.
The baseline uses $M=128$; \cref{tab:baseline-validation} records its initial contrast and validation bounds.
The time grid has 48 steps:
\[
 t_0=0,\qquad t_k=10^{-6}10^{4(k-1)/47}\quad(1\le k\le48),\qquad
 \tau_k=t_k-t_{k-1}.
\]
Here $\tau_1=10^{-6}$ and, for $k\ge2$, $\tau_k=(1-10^{-4/47})\,t_k$, so $\tau_k$ increases by the factor $10^{4/47}\approx1.22$ per step, from $\tau_2\approx2.2\cdot10^{-7}$ to $\tau_{48}\approx1.8\cdot10^{-3}$.
The early steps are small because the density minimum decreases at times of order $10^{-5}$ (\cref{fig:minimum}); \cref{tab:dip} compares the recorded minima under time refinement.

For each $p$, IPM computes $\nu^1,\dots,\nu^{48}$ on this time grid.
A step is accepted by the test of \cref{alg:jko-robust}: the floating-point primal--dual gap, replaced by its interval upper bound when it lies in $(10^{-14},10^{-12}]$, is at most $10^{-14}$.
The accepted output is the prescribed density $\nu^{\rm old}$ of the next step.
The interval bounds in \cref{tab:baseline-validation} certify the optimization accuracy of every accepted step.
They do not include propagation of input errors through the trajectory or space and time discretization errors.

\begin{table}[htbp]
\centering\small
\begin{tabular}{@{}lr@{}}
\toprule
Quantity & Value or upper bound\\
\midrule
Initial contrast, each $p$ & $\simeq\BaselineInputContrast$\\
Primal--dual gap, every accepted step & $\le\BaselineIntervalGap$\\
Discrete $L^1$ density error, every accepted step & $\le\BaselineDensityError$\\
\bottomrule
\end{tabular}
\caption{Baseline IPM validation for \eqref{eq:flowdata} at $M=128$, over all 48 steps of each of the four trajectories. The contrast is $\max_i\nu_i^0/\min_i\nu_i^0$. The gap bound uses outward-rounded interval arithmetic; the density-error bound is the bound $\sqrt{2\varsigma\eta}+|\avg{\nu_{\rm cand}}-\varsigma|$ on $\normone{\nu_{\rm cand}-\nu_\varsigma^*}$ of \cref{app:mass}, with $\varsigma=\avg{\nu^{\rm old}}$. Each bound concerns the successor density $\nu_\varsigma^*$, the minimizer of \eqref{eq:step} with $\nu^{\rm old}$ equal to the density passed to that step, $\tau=\tau_k$ and the operator $L$; it does not include errors in the passed density.}
\label{tab:baseline-validation}
\end{table}

\subsection{Density profiles}
\begin{figure}[p]
\centering
\input{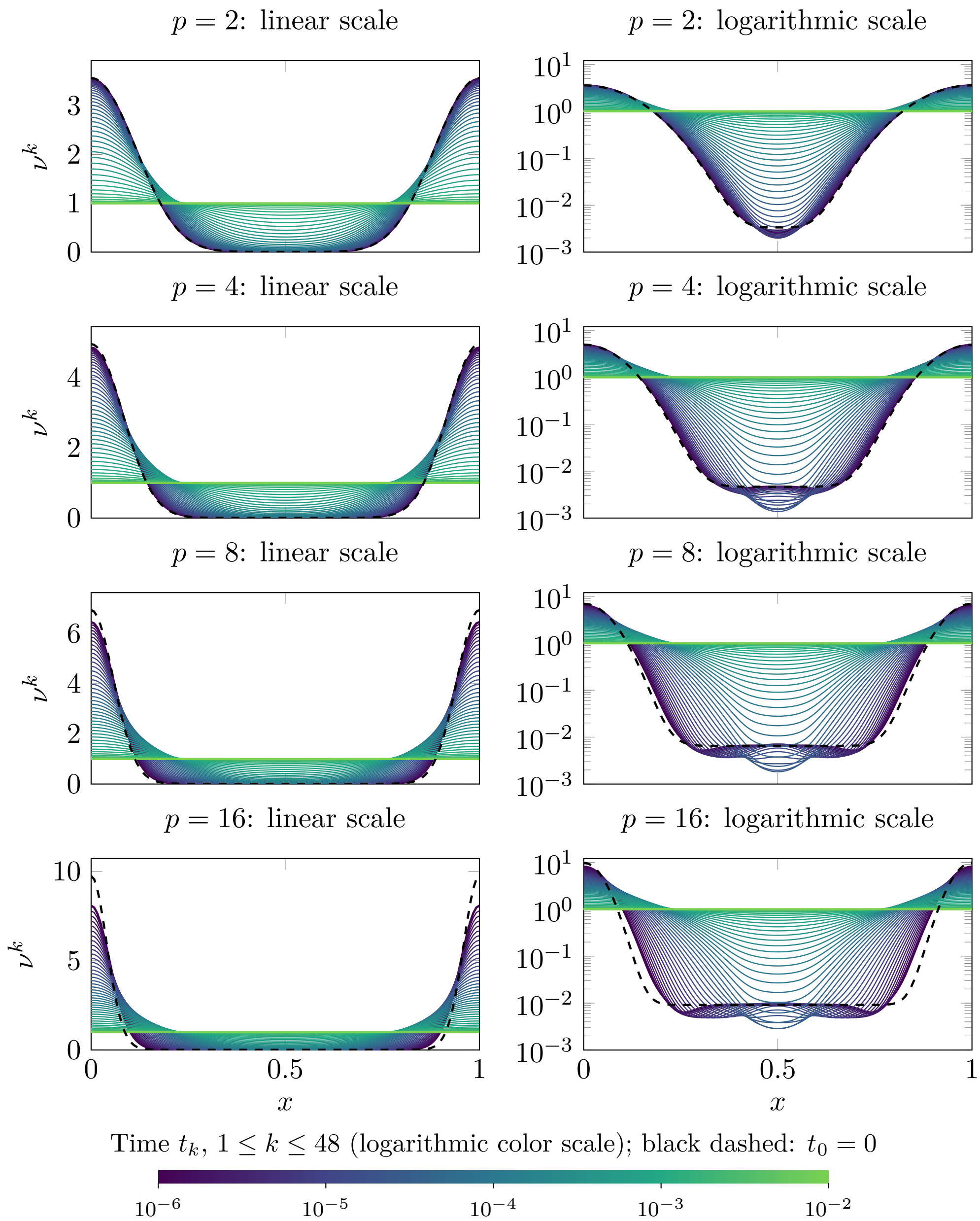}

\begin{tikzpicture}[x=1cm,y=1cm]
\shade[left color=flowtime0,right color=flowtime1] (0.000000,0) rectangle (0.191489,.17);
\shade[left color=flowtime1,right color=flowtime2] (0.191489,0) rectangle (0.382979,.17);
\shade[left color=flowtime2,right color=flowtime3] (0.382979,0) rectangle (0.574468,.17);
\shade[left color=flowtime3,right color=flowtime4] (0.574468,0) rectangle (0.765957,.17);
\shade[left color=flowtime4,right color=flowtime5] (0.765957,0) rectangle (0.957447,.17);
\shade[left color=flowtime5,right color=flowtime6] (0.957447,0) rectangle (1.148936,.17);
\shade[left color=flowtime6,right color=flowtime7] (1.148936,0) rectangle (1.340426,.17);
\shade[left color=flowtime7,right color=flowtime8] (1.340426,0) rectangle (1.531915,.17);
\shade[left color=flowtime8,right color=flowtime9] (1.531915,0) rectangle (1.723404,.17);
\shade[left color=flowtime9,right color=flowtime10] (1.723404,0) rectangle (1.914894,.17);
\shade[left color=flowtime10,right color=flowtime11] (1.914894,0) rectangle (2.106383,.17);
\shade[left color=flowtime11,right color=flowtime12] (2.106383,0) rectangle (2.297872,.17);
\shade[left color=flowtime12,right color=flowtime13] (2.297872,0) rectangle (2.489362,.17);
\shade[left color=flowtime13,right color=flowtime14] (2.489362,0) rectangle (2.680851,.17);
\shade[left color=flowtime14,right color=flowtime15] (2.680851,0) rectangle (2.872340,.17);
\shade[left color=flowtime15,right color=flowtime16] (2.872340,0) rectangle (3.063830,.17);
\shade[left color=flowtime16,right color=flowtime17] (3.063830,0) rectangle (3.255319,.17);
\shade[left color=flowtime17,right color=flowtime18] (3.255319,0) rectangle (3.446809,.17);
\shade[left color=flowtime18,right color=flowtime19] (3.446809,0) rectangle (3.638298,.17);
\shade[left color=flowtime19,right color=flowtime20] (3.638298,0) rectangle (3.829787,.17);
\shade[left color=flowtime20,right color=flowtime21] (3.829787,0) rectangle (4.021277,.17);
\shade[left color=flowtime21,right color=flowtime22] (4.021277,0) rectangle (4.212766,.17);
\shade[left color=flowtime22,right color=flowtime23] (4.212766,0) rectangle (4.404255,.17);
\shade[left color=flowtime23,right color=flowtime24] (4.404255,0) rectangle (4.595745,.17);
\shade[left color=flowtime24,right color=flowtime25] (4.595745,0) rectangle (4.787234,.17);
\shade[left color=flowtime25,right color=flowtime26] (4.787234,0) rectangle (4.978723,.17);
\shade[left color=flowtime26,right color=flowtime27] (4.978723,0) rectangle (5.170213,.17);
\shade[left color=flowtime27,right color=flowtime28] (5.170213,0) rectangle (5.361702,.17);
\shade[left color=flowtime28,right color=flowtime29] (5.361702,0) rectangle (5.553191,.17);
\shade[left color=flowtime29,right color=flowtime30] (5.553191,0) rectangle (5.744681,.17);
\shade[left color=flowtime30,right color=flowtime31] (5.744681,0) rectangle (5.936170,.17);
\shade[left color=flowtime31,right color=flowtime32] (5.936170,0) rectangle (6.127660,.17);
\shade[left color=flowtime32,right color=flowtime33] (6.127660,0) rectangle (6.319149,.17);
\shade[left color=flowtime33,right color=flowtime34] (6.319149,0) rectangle (6.510638,.17);
\shade[left color=flowtime34,right color=flowtime35] (6.510638,0) rectangle (6.702128,.17);
\shade[left color=flowtime35,right color=flowtime36] (6.702128,0) rectangle (6.893617,.17);
\shade[left color=flowtime36,right color=flowtime37] (6.893617,0) rectangle (7.085106,.17);
\shade[left color=flowtime37,right color=flowtime38] (7.085106,0) rectangle (7.276596,.17);
\shade[left color=flowtime38,right color=flowtime39] (7.276596,0) rectangle (7.468085,.17);
\shade[left color=flowtime39,right color=flowtime40] (7.468085,0) rectangle (7.659574,.17);
\shade[left color=flowtime40,right color=flowtime41] (7.659574,0) rectangle (7.851064,.17);
\shade[left color=flowtime41,right color=flowtime42] (7.851064,0) rectangle (8.042553,.17);
\shade[left color=flowtime42,right color=flowtime43] (8.042553,0) rectangle (8.234043,.17);
\shade[left color=flowtime43,right color=flowtime44] (8.234043,0) rectangle (8.425532,.17);
\shade[left color=flowtime44,right color=flowtime45] (8.425532,0) rectangle (8.617021,.17);
\shade[left color=flowtime45,right color=flowtime46] (8.617021,0) rectangle (8.808511,.17);
\shade[left color=flowtime46,right color=flowtime47] (8.808511,0) rectangle (9.000000,.17);
\draw (0,0)--(0,-.07) node[below,font=\scriptsize] {$10^{-6}$};
\draw (2.25,0)--(2.25,-.07) node[below,font=\scriptsize] {$10^{-5}$};
\draw (4.5,0)--(4.5,-.07) node[below,font=\scriptsize] {$10^{-4}$};
\draw (6.75,0)--(6.75,-.07) node[below,font=\scriptsize] {$10^{-3}$};
\draw (9,0)--(9,-.07) node[below,font=\scriptsize] {$10^{-2}$};
\node[font=\small,anchor=south] at (4.5,.20) {Time $t_k$, $1\le k\le48$ (logarithmic color scale); black dashed: $t_0=0$};
\end{tikzpicture}
\caption{Densities $\nu^k$, $0\le k\le48$, of the four baseline IPM trajectories for \eqref{eq:flowdata}, with $M=128$ and $\varepsilon=10^{-3}$; $\nu^0$ is the initial density and, for $1\le k\le48$, $\nu^k$ is the accepted density at time $t_k$ of \cref{sec:flowdata}.
Rows: $p=2,4,8,16$, from top to bottom.
Each row shows the same densities on a linear (left) and a logarithmic (right) vertical scale.
The black dashed curve is $\nu^0$.
Colors encode $t_k$, $1\le k\le48$, on the logarithmic scale of the color bar; the two thicker colored curves are $\nu^1$ ($t_1=10^{-6}$) and $\nu^{48}$ ($t_{48}=10^{-2}$).
Only the recorded times $t_k$ are plotted; there is no interpolation in time.
For $p=16$, the dip time $t_{\rm dip}$ of \eqref{eq:dip}, at which $\min_i\nu_i^k$ is smallest over the recorded times, is $t_{\rm dip}=t_{19}\simeq\BaselineRecordedDipTime$.
Spatial line segments close the periodic interval, without a positive display floor.}
\label{fig:flowprofiles}
\end{figure}

\Cref{fig:flowprofiles} shows the densities $\nu^k$ of the four baseline IPM trajectories; the decrease of their minima below $\min_i\nu_i^0$ is quantified in \cref{sec:minimum}.
The initial peaks decrease and spread, and at late times the densities approach the uniform density.
The logarithmic panels show the values of $\nu^k$ near $\min_i\nu_i^k$.
Each curve joins the nodal values $\nu_i^k$ by line segments.
The bounds of \cref{tab:baseline-validation} concern the optimization accuracy at the nodes and time levels $t_k$ only; they give no bound between nodes or between time levels.

\subsection{Computational cost}\label{sec:flowcost}
We compare the time that DR and IPM need to reach the same validated accuracy on one entropic step from $\nu^{\rm old}=\nu^0$ of \eqref{eq:flowdata}, with $(M,p,\varepsilon)$ as listed in \cref{tab:common-entropic} and $\tau=10^{-6}$.
Each solver stops at the first candidate whose interval upper bound on the primal--dual gap $J_\tau(\nu)-D_\tau(\chi)$ is at most $10^{-9}$.
For DR, the multiplier inserted in \eqref{eq:dual} is $\chi=(z_\nu-x_\nu)/\omega$, where $z_\nu$ is the density block of the auxiliary iterate and $x_\nu$ the density block of its proximal output (\cref{alg:jko-dr}).
The bound $D_\tau(\chi)\le H_\tau^*$ in \eqref{eq:gap} holds for every $\chi\in\R^M$ (\cref{prop:gap}), so this choice requires no optimality property of $\chi$.
This choice approximates an optimal multiplier when the auxiliary iterate is close to the DR fixed point $z_*$ of \cref{prop:drglobal}.
At $z_*$, the proximal output is $x_*=(\nu^*,\sqrt\beta m^*)$, and the first-order condition of the proximal map, with $\omega=\gamma\dx$, gives the entrywise identity $(z_{*,\nu}-x_{*,\nu})/\omega=\log\nu^*-\tfrac\tau2\bigl(m^*/(\tau\nu^*)\bigr)^2$.
By the first-order conditions in \cref{app:gap}, the right side is an optimal multiplier minus a constant vector, and constants do not change $D_\tau$.
At each checkpoint (every 25 DR updates), the proximal density is normalized, its least-action momentum is reconstructed by \eqref{eq:elimination}, and the interval bound $\eta$ is computed when the floating-point gap is at most $10^{-9}$.
The IPM candidates are subject to the same test $\eta\le10^{-9}$, without the floating-point screening.
These tests differ from the stopping rules of \cref{alg:jko-dr,alg:jko-robust} used for the trajectories of \cref{tab:timing}.
All attempts and their times are retained, including IPM attempts rejected by the gap test.

\begin{table}[htbp]
\centering\small
\begin{tabular}{@{}rrrlrrrr@{}}
\toprule
$M$ & $p$ & $\varepsilon$ & Solver ($c_\omega$) & Iter. & Solve (s) & Total (s) & Gap bound\\
\midrule
32 & 4 & $1\cdot10^{-3}$ & DR (0.3) & 775 & 0.179 & 0.207 & $8.1\cdot10^{-10}$\\
32 & 4 & $1\cdot10^{-3}$ & DR (1) & 300 & 0.071 & 0.082 & $9.5\cdot10^{-10}$\\
32 & 4 & $1\cdot10^{-3}$ & DR (3) & 900 & 0.201 & 0.213 & $7.4\cdot10^{-10}$\\
32 & 4 & $1\cdot10^{-3}$ & IPM & 11 & 0.001 & 0.011 & $6.3\cdot10^{-10}$\\
64 & 16 & $1\cdot10^{-3}$ & DR (0.3) & 725 & 0.172 & 0.223 & $8.3\cdot10^{-10}$\\
64 & 16 & $1\cdot10^{-3}$ & DR (1) & 350 & 0.085 & 0.118 & $8.0\cdot10^{-10}$\\
64 & 16 & $1\cdot10^{-3}$ & DR (3) & 1025 & 0.233 & 0.267 & $9.7\cdot10^{-10}$\\
64 & 16 & $1\cdot10^{-3}$ & IPM & 31 & 0.008 & 0.072 & $3.9\cdot10^{-11}$\\
64 & 16 & $1\cdot10^{-6}$ & DR (0.3) & 21750 & 4.659 & 4.769 & $1.0\cdot10^{-9}$\\
64 & 16 & $1\cdot10^{-6}$ & DR (1) & 6750 & 1.494 & 1.550 & $9.6\cdot10^{-10}$\\
64 & 16 & $1\cdot10^{-6}$ & DR (3) & 18900 & 3.824 & 3.923 & $1.0\cdot10^{-9}$\\
64 & 16 & $1\cdot10^{-6}$ & IPM & 30 & 0.008 & 0.072 & $7.5\cdot10^{-11}$\\
\bottomrule
\end{tabular}
\caption{Single entropic steps \eqref{eq:step} with $\tau=10^{-6}$ from the old density $\nu^{\rm old}=\nu^0$ of \eqref{eq:flowdata}, the same for all solvers. The gap bound is the interval upper bound $\eta$ on $J_\tau(\nu)-D_\tau(\chi)$ (\cref{sec:certificate}), with the mass corrections of \cref{app:mass}; it is at most $10^{-9}$ in every row, and displayed values are rounded upward. The operator $L$ has the symbol $\widehat L(k)=-(2\pi k)^2$ with $2\pi k$ evaluated in binary64, not the rounded symbol of \cref{tab:operator-rounding}. DR uses the residual-checked proximal solver of \cref{alg:jko-prox-checked} (\cref{rem:checked-prox}), $\beta=1/\tau$, and $\omega=c_\omega\sqrt{\min\nu^{\rm old}\max\nu^{\rm old}}$, with all three tested factors $c_\omega$ reported. DR checks the gap every 25 updates. Solve time records its iteration loop, or the sum of IPM solver times over all attempts. Total time also includes initialization, assembly, candidate normalization, floating-point screening, interval validation, and candidate recording. IPM tries tolerances $10^{-7},10^{-9},10^{-11}$ until acceptance; all attempted iterations and times are counted. Times are single local executions on logical CPU 0 with one thread, not medians. The times exclude the computation of the local spectra reported in \cref{tab:entropic-parameter-check}.}
\label{tab:common-entropic}
\end{table}

For the three inputs, the total time of DR with $c_\omega=1$ is 7.5 times ($M=32$, $p=4$, $\varepsilon=10^{-3}$), 1.6 times ($M=64$, $p=16$, $\varepsilon=10^{-3}$) and 22 times ($M=64$, $p=16$, $\varepsilon=10^{-6}$) the IPM total time.
For each input, the factor $c_\omega=1$, that is, the rule $\omega=\sqrt{\min\nu^{\rm old}\max\nu^{\rm old}}$, gives the fewest DR iterations and the smallest DR times among $c_\omega\in\{0.3,1,3\}$.
\Cref{tab:entropic-parameter-check} gives the total time of the three DR runs needed to compare the three factors, together with the final local spectra.
This single-step comparison does not establish a speed ratio for arbitrary trajectories, grids, or implementations.

\paragraph{Complete 48-step trajectories at $M=128$.}

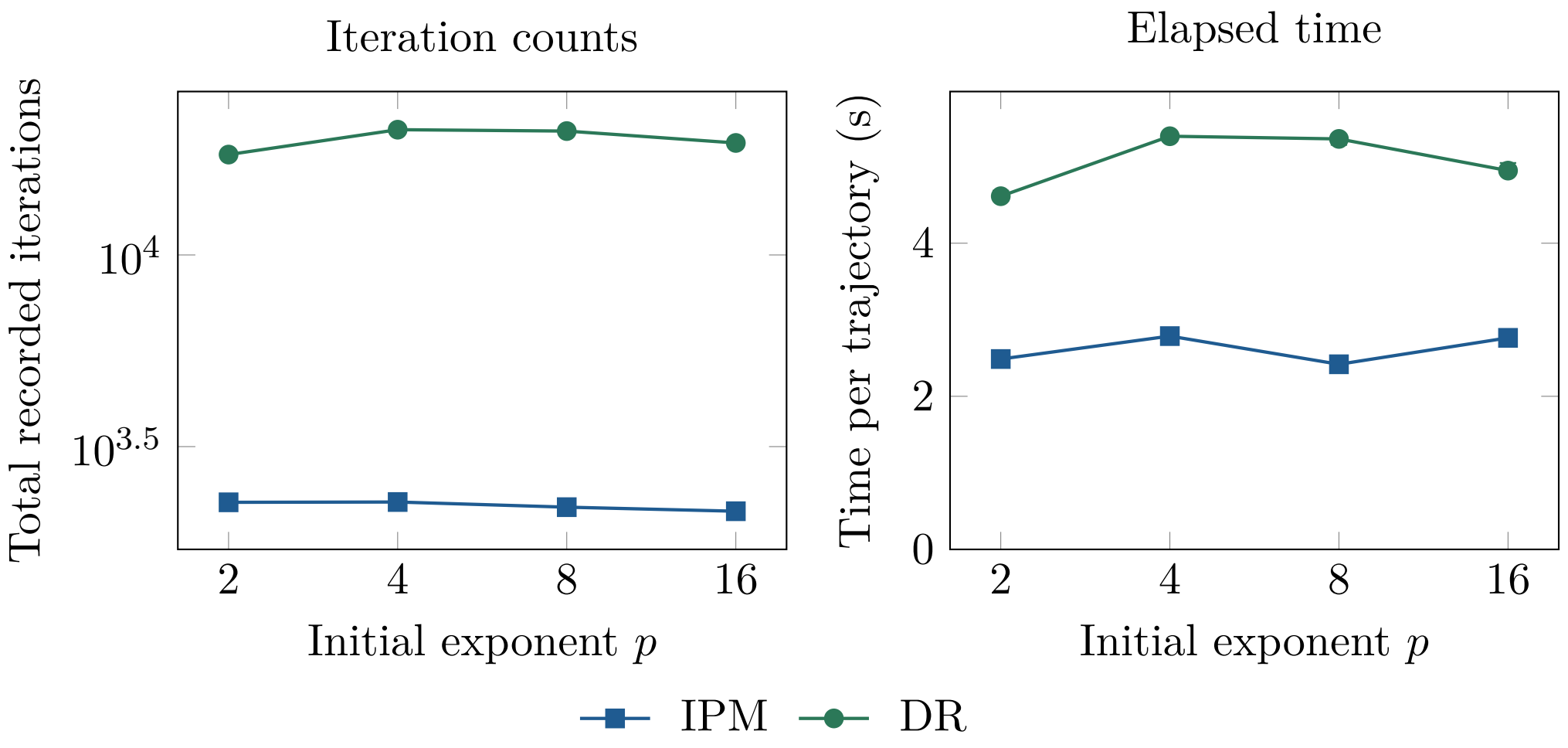
\begin{figure}[htbp]
\centering
\begin{tikzpicture}
\begin{groupplot}[group style={group size=2 by 1,horizontal sep=1.4cm},
 width=.43\textwidth,height=5.5cm,xmin=.7,xmax=4.3,xtick={1,2,3,4},xticklabels={2,4,8,16},
 xlabel={Initial exponent $p$},legend columns=2,legend style={draw=none,/tikz/column sep=5pt}]
\nextgroupplot[ymode=log,ylabel={Total recorded iterations},title={Iteration counts},legend to name=legend:flowcost]
\addplot[plotblue,mark=square*,thick] coordinates {(1,2263) (2,2267) (3,2198) (4,2145)};
\addlegendentry{IPM}
\addplot[plotgreen,mark=*,thick] coordinates {(1,18278) (2,21224) (3,21048) (4,19601)};
\addlegendentry{DR}
\nextgroupplot[ymin=0,ylabel={Time per trajectory (s)},title={Elapsed time}]
\addplot[plotblue,mark=square*,thick,error bars/.cd,y dir=both,y explicit] coordinates { (1,2.486951981) += (0,0.04953750499) -= (0,0.008655059006) (2,2.786488282) += (0,0.05910266101) -= (0,0.03861417502) (3,2.418339997) += (0,0.06425993601) -= (0,0.02751232198) (4,2.763447168) += (0,0.06105952698) -= (0,0.05114039298)};
\addplot[plotgreen,mark=*,thick,error bars/.cd,y dir=both,y explicit] coordinates { (1,4.613908703) += (0,0.07238030102) -= (0,0.06130980901) (2,5.397413717) += (0,0.039968631) -= (0,0.05984783801) (3,5.362321466) += (0,0.05893939399) -= (0,0.07823014501) (4,4.948262537) += (0,0.100521841) -= (0,0.02847555903)};
\end{groupplot}
\end{tikzpicture}

\pgfplotslegendfromname{legend:flowcost}
\caption{Work for one complete 48-step entropic trajectory at $M=128$ and $\varepsilon=10^{-3}$.
Left: recorded iteration totals, with a logarithmic scale.
The IPM count includes only the solve accepted at each time step; elapsed time includes both solves when a step is repeated.
DR counts splitting iterations.
Right: median elapsed time over five runs; bars give the smallest and largest times.
Lines connect the four tested exponents for display.
\Cref{tab:timing} specifies the stopping rules and timing exclusions.}
\label{fig:flowcost}
\end{figure}

\Cref{tab:timing,fig:flowcost} compare complete 48-step trajectories from the four baseline initial data \eqref{eq:flowdata}, computed by DR (\cref{alg:jko-dr} with \cref{alg:jko-prox}) and by the IPM procedure of \cref{alg:jko-robust}; the caption of \cref{tab:timing} lists both stopping rules.
One DR iteration applies the proximal map and affine projection; one IPM iteration contains a conic Newton factorization.
The iteration counts therefore measure different work.

\begin{table}[htbp]
\centering\small
\begin{tabular}{@{}crrrr@{}}
\toprule
 & \multicolumn{2}{c}{IPM} & \multicolumn{2}{c}{DR}\\
\cmidrule(lr){2-3}\cmidrule(l){4-5}
$p$ & Iterations & Seconds & Iterations & Seconds\\
\midrule
16 & 2145 & 2.76 & 19601 & 4.95\\
8 & 2198 & 2.42 & 21048 & 5.36\\
4 & 2267 & 2.79 & 21224 & 5.40\\
2 & 2263 & 2.49 & 18278 & 4.61\\
\bottomrule
\end{tabular}
\caption{Iteration counts and median seconds per 48-step trajectory.
Data are \eqref{eq:flowdata}, with $M=128$ and $\varepsilon=10^{-3}$.
Times are medians of five alternating local runs per method, using one thread on one Intel i9-13950HX performance core.
IPM iterations count only the accepted solve at each step; times include all 49 calls for $p=16,4$ and all 48 calls for $p=8,2$.
For $p=16$ and $p=4$, the 49th call is the additional attempt that \cref{alg:jko-robust} makes after an accepted attempt whose status is not \texttt{Solved}, on the step from $t_{12}$ to $t_{13}$ and from $t_8$ to $t_9$, respectively.
Times include initialization and, for IPM, matrix assembly, density normalization, momentum reconstruction by \eqref{eq:elimination}, and floating-point gap evaluation.
Separate dual optimization and interval validation are excluded; validation adds approximately 5.4 seconds per IPM trajectory.
DR stops when $\|y^\ell-x^\ell\|_2\le10^{-13}\|z^{\ell+1}\|_2$, with the scaled iterates of \eqref{eq:dr} (\cref{alg:jko-dr}).
IPM uses tolerances $10^{-12}$ and acceptance gap $10^{-14}$, with the fallback in \cref{alg:jko-robust}.
These are different stopping rules; \cref{sec:flowdata} gives the attained IPM accuracy.}\label{tab:timing}
\end{table}

For each $p\in\{2,4,8,16\}$, the median DR time is between 1.8 and 2.2 times the median IPM time (\cref{tab:timing}, with the stopping rules of its caption).
These complete-trajectory timings do not measure convergence within an individual solve.

To measure how the solver cost depends on the density contrast and on $M$, we vary $\varepsilon$ at $M=128$, $p=16$ on the 48-step time grid (\cref{tab:flowcontrast,tab:flowcontrastdr}), and we vary $M$ for the first step of the $p=16$, $\varepsilon=10^{-3}$ data at $\tau=10^{-6}$ (\cref{tab:flowgridsize}).

\begin{table}[htbp]
\centering\small
\begin{tabular}{@{}rrrrr@{}}
\toprule
$\varepsilon$ & Maximum input contrast & IPM iterations & Seconds & Maximum absolute gap bound $\eta$\\
\midrule
$10^{-5}$ & $3.03\cdot10^{5}$ & 2094 & 2.36 & $7.74\cdot10^{-18}$\\
$10^{-6}$ & $4.52\cdot10^{6}$ & 2211 & 2.42 & $1.61\cdot10^{-17}$\\
$10^{-7}$ & $4.86\cdot10^{7}$ & 2078 & 3.08 & $5.32\cdot10^{-17}$\\
$10^{-8}$ & $1.10\cdot10^{9}$ & 2249 & 2.51 & $2.81\cdot10^{-17}$\\
\bottomrule
\end{tabular}
\caption{IPM contrast test for \eqref{eq:flowdata}, at $p=16$, $M=128$, with the baseline 48-step time grid. Input contrast is $\max_{1\le k\le48}(\max_i\nu_i^{k-1}/\min_i\nu_i^{k-1})$. Each time is one recorded complete-trajectory run, excluding interval validation. The last column bounds the gap of every accepted step by outward-rounded interval evaluation; its displayed bounds are rounded upward. Counts are totals over the 48 accepted solves, one per step. At $\varepsilon=10^{-7}$, the accepted solves of two steps have the Clarabel status \texttt{AlmostSolved} (\cref{tab:flowgridsize}), and \cref{alg:jko-robust} makes one additional attempt at each of these two steps; the time of this run includes its 50 solver calls. All other accepted solves have status \texttt{Solved}.}
\label{tab:flowcontrast}
\end{table}

\begin{table}[htbp]
\centering\small
\begin{tabular}{@{}rrr@{}}
\toprule
$\varepsilon$ & DR iterations & Maximum DR--IPM $L^1$ difference\\
\midrule
$10^{-5}$ & 138264 & $1.61\cdot10^{-9}$\\
$10^{-6}$ & 386869 & $2.29\cdot10^{-9}$\\
\bottomrule
\end{tabular}
\caption{DR check of the two contrast-test trajectories at $p=16$, $M=128$, on the same 48-step time grid as \cref{tab:flowcontrast}. DR uses \cref{alg:jko-dr} with \cref{alg:jko-prox}, $\beta=1/\tau$, $\omega=\sqrt{\min\nu^{\rm old}\max\nu^{\rm old}}$, and the stopping rule of \cref{tab:timing}; DR iteration counts are totals over the 48 steps. The difference is $\max_k\|\nu_{\rm DR}^k-\nu_{\rm IPM}^k\|_{1,M}$ at saved time levels, where $\nu_{\rm DR}^k$ is the normalized proximal density. The IPM candidates match those in the contrast test; these independent DR runs are not part of the five-run baseline timing benchmark.}
\label{tab:flowcontrastdr}
\end{table}

\begin{table}[htbp]
\centering\small
\begin{tabular}{@{}rrrl@{}}
\toprule
$M$ & IPM iterations & Solver seconds & Status\\
\midrule
128 & 32 & 0.033 & Solved\\
256 & 53 & 0.33 & Solved\\
512 & 51 & 2.5 & Solved\\
1024 & 60 & 29 & Solved\\
2048 & 79 & 374 & AlmostSolved\\
\bottomrule
\end{tabular}
\caption{Spatial-size test for the first entropic step of the $p=16$, $\varepsilon=10^{-3}$ baseline problem, with $\tau=10^{-6}$. Each row is one saved IPM run; seconds are the solver-reported solve time, excluding setup, matrix assembly, and separate interval validation. Each run uses QDLDL on one thread. These are fixed-step costs, not trajectory times. \texttt{AlmostSolved} is the Clarabel status of a solve that meets only its reduced-accuracy tolerances. The $M=2048$ candidate has floating-point gap $6.7\cdot10^{-16}$ and no interval bound.}
\label{tab:flowgridsize}
\end{table}

At $M=128$, the IPM totals lie between 2078 and 2249 iterations and between 2.36 and 3.08 s for maximum input contrasts from $3.0\cdot10^5$ to $1.1\cdot10^9$ (\cref{tab:flowcontrast}), against 2145 iterations for $p=16$ at $\varepsilon=10^{-3}$ (\cref{tab:timing}).
DR needs $1.4\cdot10^5$ and $3.9\cdot10^5$ iterations at $\varepsilon=10^{-5}$ and $10^{-6}$ (\cref{tab:flowcontrastdr}), against $2.0\cdot10^4$ at $\varepsilon=10^{-3}$ (\cref{tab:timing}).
For the first step, the IPM solve time increases from 0.033 s at $M=128$ to 374 s at $M=2048$ (\cref{tab:flowgridsize}).
These results support using IPM for the moderate-grid minimum-density study; they do not establish a cost bound independent of grid size or contrast.

\begin{figure}[htbp]
\centering
\input{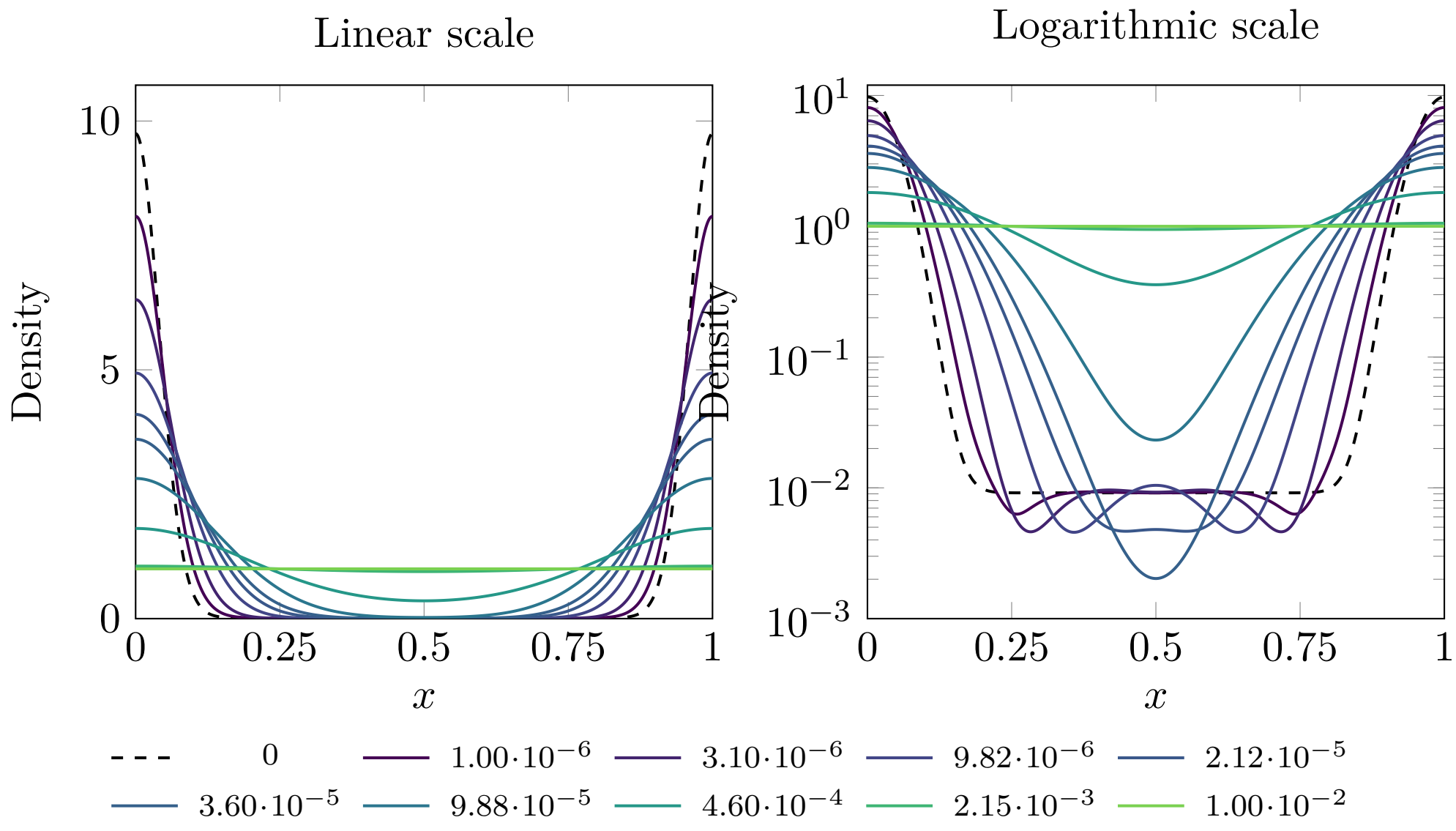}

\pgfplotslegendfromname{legend:flowfine}
\caption{Finer entropic profiles for $p=16$, $\varepsilon=10^{-3}$, $M=1024$, and 384 DR time steps. Positive times are geometrically spaced from $10^{-6}$ to $10^{-2}$; all 385 densities are saved. \Cref{tab:flowfine-diagnostics} records the minimum and solver diagnostics. The run uses the rounded symbol $\widehat L_{\rm round}$ of \cref{tab:operator-rounding} (\cref{app:records}).
The two panels show the same ten computed profiles, at $t_k$ for $k=0,1,48,96,128,192,256,320,384$ and at the dip time $t_{150}$; the legend gives their times.
The black dashed curve is the initial density.
All 1024 nodes of each selected profile are plotted.
There is no temporal interpolation or positive display floor.}
\label{fig:flowfine}
\end{figure}

\paragraph{Finer entropic profiles.}
For the finer space--time computation, DR avoids the dense factorization used by the IPM implementation.
\Cref{fig:flowfine} resolves the decrease and subsequent recovery of the density minimum on a finer space--time grid.
The relative minimum $\min_{k,i}\nu_i^k/\min_i\nu_i^0$ of this computation is $\FineRecordedMinimum$, at $t\simeq\FineRecordedDipTime$ (\cref{tab:flowfine-diagnostics}).
This value is below the value $0.312$ on the baseline grid and close to the value $0.216$ on the finest time grid, both in the row $\varepsilon=10^{-3}$ of \cref{tab:dip}, at $M=128$.
Since both the spatial and the time grids change, this comparison does not isolate the time-discretization error.

\begin{table}[htbp]
\centering\small
\begin{tabular}{@{}lr@{}}
\toprule
Quantity & Recorded value\\
\midrule
Relative minimum & $\FineRecordedMinimum$\\
Dip time & $\simeq\FineRecordedDipTime$\\
Maximum floating-point step gap & $\FineStepGap$\\
Maximum mass error & $\FineMassError$\\
Maximum iterations in one step & $\FineIterationMaximum$\\
\bottomrule
\end{tabular}
\caption{Diagnostics for the finer DR trajectory in \cref{fig:flowfine}, with $p=16$, $\varepsilon=10^{-3}$, $M=1024$, and 384 time steps. The relative minimum is $\min_{k,i}\nu_i^k/\min_i\nu_i^0$; the dip time is a saved time attaining it. The mass error is $|\avg{\nu^k}-1|$. Each step uses \cref{alg:jko-dr} with \cref{alg:jko-prox}, $\beta=1/\tau$ and $\omega=\sqrt{\min\nu^{\rm old}\max\nu^{\rm old}}$; it stops when $\|y^\ell-x^\ell\|_2\le10^{-13}\|z^{\ell+1}\|_2$ or after 50000 iterations. The step gap is $J_\tau(\nu)-D_\tau(\chi)$ for the normalized proximal density $\nu$ and the multiplier $\chi$ of \cref{alg:jko-dual}. Gap, mass-error, and iteration maxima are taken over the trajectory. The gap and mass-error values are floating-point checks, not interval bounds.}
\label{tab:flowfine-diagnostics}
\end{table}

\FloatBarrier
\subsection{Decrease of the minimum density}\label{sec:minimum}

We measure how far the density minimum decreases below its initial value.
The entropy inequality in \cref{prop:positive} does not imply that the density minimum is nondecreasing.

Let $\mathcal T=\{t_0,\ldots,t_K\}$ be the recorded time grid, with $t_0=0$.
For a positive initial density, define
\begin{equation}\label{eq:dip}
 d_{M,\mathcal T}(\nu^0)
 =\min_{0\le k\le K}\frac{\min_{0\le i<M}\nu_i^k}{\min_{0\le i<M}\nu_i^0},\qquad
 t_{\rm dip}\in\argmin_{t_k\in\mathcal T}\min_i\nu_i^k.
\end{equation}
We call $d_{M,\mathcal T}$ the relative minimum and $t_{\rm dip}$ a dip time.
Both quantities refer to the recorded nodes and times; $d_{M,\mathcal T}<1$ indicates a decrease below the initial minimum.

\begin{figure}[htbp]
\centering
\begin{tikzpicture}
\begin{semilogxaxis}[width=.84\textwidth,height=6.5cm,
 xmin=1e-6,xmax=2e-4,ymin=0,ymax=2,
 xlabel={$t$},ylabel={$\min_i\nu_i^k/\min_i\nu_i^0$},
 legend style={at={(.02,.98)},anchor=north west,legend columns=2}]
\addplot[plotblue,solid,thick,mark=*,mark size=1.2pt] coordinates {(1e-06,0.853816875) (1.21648395e-06,0.827086341) (1.4798332e-06,0.798891764) (1.80019333e-06,0.766720792) (2.18990629e-06,0.731262317) (2.66398586e-06,0.69419876) (3.24069604e-06,0.657949897) (3.94225471e-06,0.625774643) (4.79568958e-06,0.601878889) (5.8338794e-06,0.591449284) (7.09682065e-06,0.600652512) (8.63316841e-06,0.636794759) (1.05021108e-05,0.708927232) (1.27756492e-05,0.82907609) (1.55413722e-05,1.01409379) (1.89058298e-05,1.2881121) (2.29986385e-05,1.68568814) (2.79774746e-05,2.2557733) (3.40341488e-05,3.06650636) (4.14019958e-05,4.21055038) (5.03648633e-05,5.81031388) (6.12680478e-05,8.02200334) (7.45315967e-05,11.0371681) (9.06664911e-05,15.0803662) (0.000110294331,20.4019158) (0.000134171284,27.2653729) (0.000163217213,35.9301556) (0.00019855112,46.6302099) (0.00024153425,59.5494924) (0.000293822538,74.7944417)};
\addlegendentry{$p=2$}
\addplot[plotred,dashed,thick,mark=*,mark size=1.2pt] coordinates {(1e-06,0.898415626) (1.21648395e-06,0.884902334) (1.4798332e-06,0.865595457) (1.80019333e-06,0.841932556) (2.18990629e-06,0.816335424) (2.66398586e-06,0.789526261) (3.24069604e-06,0.763160218) (3.94225471e-06,0.7384571) (4.79568958e-06,0.71490431) (5.8338794e-06,0.688487896) (7.09682065e-06,0.653666434) (8.63316841e-06,0.600000118) (1.05021108e-05,0.513921302) (1.27756492e-05,0.40599807) (1.55413722e-05,0.325006016) (1.89058298e-05,0.296125232) (2.29986385e-05,0.335966747) (2.79774746e-05,0.461480388) (3.40341488e-05,0.708000135) (4.14019958e-05,1.13685429) (5.03648633e-05,1.83945839) (6.12680478e-05,2.94069213) (7.45315967e-05,4.59973783) (9.06664911e-05,7.00669648) (0.000110294331,10.3741494) (0.000134171284,14.9239174) (0.000163217213,20.8702718) (0.00019855112,28.4013012) (0.00024153425,37.6597309) (0.000293822538,48.7234071)};
\addlegendentry{$p=4$}
\addplot[plotgreen,densely dotted,thick,mark=*,mark size=1.2pt] coordinates {(1e-06,0.792681986) (1.21648395e-06,0.782329694) (1.4798332e-06,0.749297863) (1.80019333e-06,0.710425829) (2.18990629e-06,0.675729617) (2.66398586e-06,0.642428573) (3.24069604e-06,0.61415872) (3.94225471e-06,0.591449784) (4.79568958e-06,0.571927501) (5.8338794e-06,0.560300403) (7.09682065e-06,0.558623054) (8.63316841e-06,0.565545777) (1.05021108e-05,0.575081379) (1.27756492e-05,0.571609856) (1.55413722e-05,0.527124553) (1.89058298e-05,0.415696629) (2.29986385e-05,0.301427442) (2.79774746e-05,0.280334243) (3.40341488e-05,0.364701315) (4.14019958e-05,0.575355597) (5.03648633e-05,0.969030397) (6.12680478e-05,1.63445497) (7.45315967e-05,2.68959167) (9.06664911e-05,4.27761111) (0.000110294331,6.55932914) (0.000134171284,9.70218126) (0.000163217213,13.8668427) (0.00019855112,19.1929641) (0.00024153425,25.7850987) (0.000293822538,33.6989448)};
\addlegendentry{$p=8$}
\addplot[plotpurple,dashdotted,thick,mark=*,mark size=1.2pt] coordinates {(1e-06,0.687796188) (1.21648395e-06,0.683583625) (1.4798332e-06,0.640205831) (1.80019333e-06,0.597789495) (2.18990629e-06,0.567937507) (2.66398586e-06,0.549483939) (3.24069604e-06,0.535299732) (3.94225471e-06,0.53006404) (4.79568958e-06,0.526125983) (5.8338794e-06,0.526001308) (7.09682065e-06,0.529299104) (8.63316841e-06,0.538430912) (1.05021108e-05,0.557150438) (1.27756492e-05,0.580030392) (1.55413722e-05,0.584658986) (1.89058298e-05,0.541900313) (2.29986385e-05,0.420260462) (2.79774746e-05,0.312787498) (3.40341488e-05,0.311584962) (4.14019958e-05,0.423886506) (5.03648633e-05,0.673897346) (6.12680478e-05,1.11916546) (7.45315967e-05,1.84242353) (9.06664911e-05,2.94612878) (0.000110294331,4.54603973) (0.000134171284,6.76275216) (0.000163217213,9.71180099) (0.00019855112,13.493319) (0.00024153425,18.1819797) (0.000293822538,23.8172685)};
\addlegendentry{$p=16$}
\addplot[black!60,dashed,forget plot] coordinates {(1e-6,1) (2e-4,1)};
\end{semilogxaxis}
\end{tikzpicture}
\caption{Ratio $\min_i\nu_i^k/\min_i\nu_i^0$, whose minimum over $k$ is the relative minimum \eqref{eq:dip}, for the baseline IPM trajectories of \eqref{eq:flowdata}, with $\varepsilon=10^{-3}$ and $M=128$, at the saved times $t_k$ in the displayed window $10^{-6}\le t\le2\cdot10^{-4}$.
Markers identify the saved time levels; segments connect them for display only.
The horizontal dashed line is the value 1.}
\label{fig:minimum}
\end{figure}
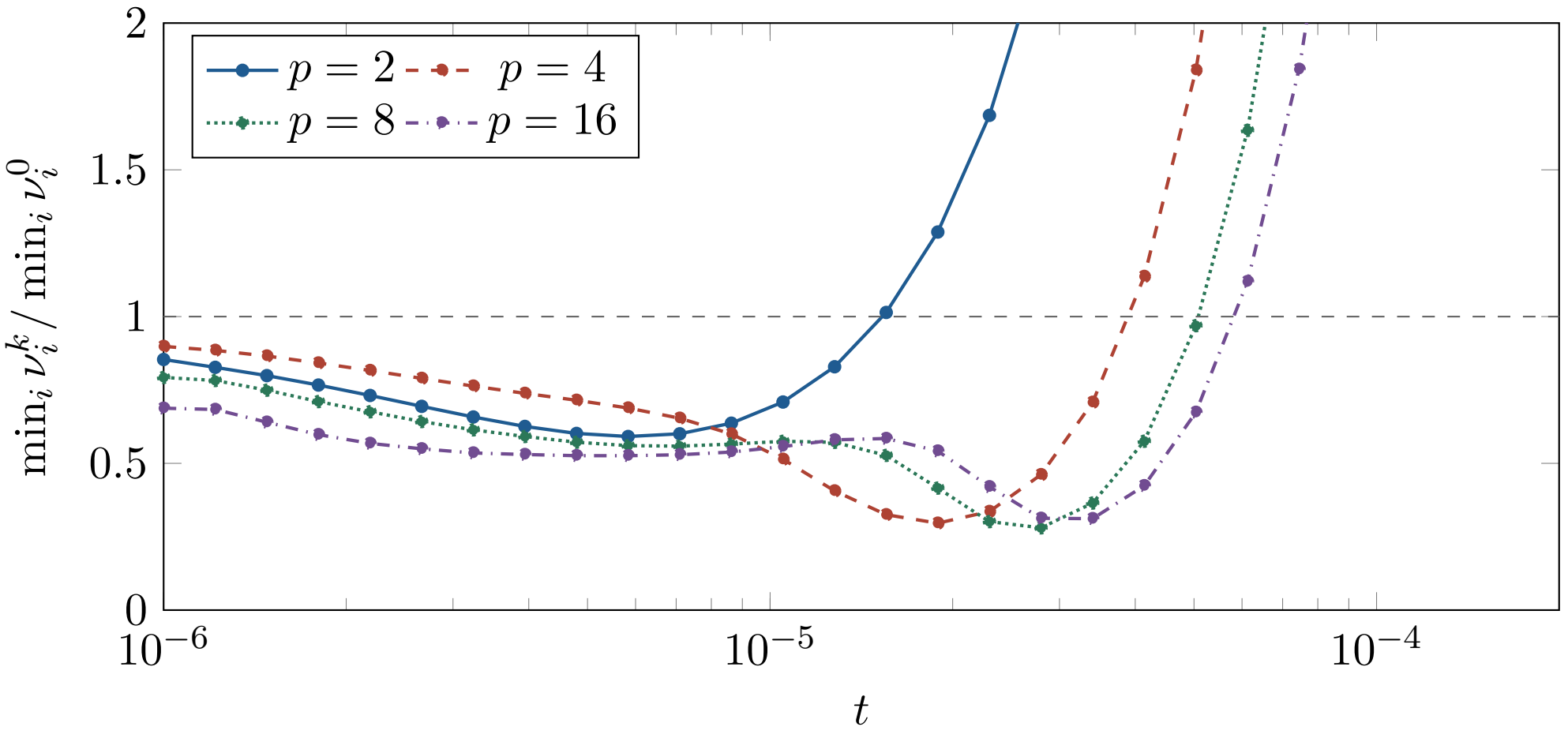

For each $p\in\{2,4,8,16\}$, the ratio $\min_i\nu_i^k/\min_i\nu_i^0$ of the baseline trajectory falls below 1 and then increases (\cref{fig:minimum}); \cref{fig:flowprofiles} shows all saved profiles up to $t=10^{-2}$.
\Cref{tab:dip} shows that time refinement gives smaller relative minima.

\subsection{Time refinement}

If $\min_i\nu_i$ decreases and then increases within a time interval shorter than a step of the baseline grid, the baseline grid can miss its smallest value.
Refinement changes both the discrete evolution and the times at which its minimum is sampled.
The comparison below tests their combined effect while keeping the initial datum and spatial grid fixed within each row.

We fix $p=16$, $M=128$, and compare the baseline grid with the refined grids $\mathcal T_1,\mathcal T_2,\mathcal T_3$ of \cref{tab:dipgrids}, up to $t=2\cdot10^{-4}$.
This end time exceeds the $\mathcal T_3$ dip times of \cref{tab:dip}, which are at most $3.6\cdot10^{-5}$.

\begin{table}[htbp]
\centering\small
\begin{tabular}{@{}lrrr@{}}
\toprule
Grid & First positive time & Consecutive-time ratio & Steps\\
\midrule
$\mathcal T_1$ & $10^{-8}$ & $10^{1/47}$ & 204\\
$\mathcal T_2$ & $10^{-9}$ & $10^{1/94}$ & 500\\
$\mathcal T_3$ & $10^{-10}$ & $10^{1/188}$ & 1186\\
\bottomrule
\end{tabular}
\caption{Refined time grids $\mathcal T_1,\mathcal T_2,\mathcal T_3$ for the minimum-density refinement test. Each grid starts at zero, then uses geometric positive times with the listed first positive time and consecutive-time ratio. The step counts are those through the end $t=2\cdot10^{-4}$ of the test interval of \cref{tab:dip}; \cref{tab:isolatedsurvey} uses the same grids up to $t=10^{-3}$. The baseline uses the time grid in \cref{sec:flowdata}, restricted to this interval. Both the first step and the later steps change under refinement.}
\label{tab:dipgrids}
\end{table}

\begin{table}[htbp]
\centering\small
\begin{tabular}{@{}crrrrr@{}}
\toprule
$\varepsilon$ & Baseline & $\mathcal T_1$ & $\mathcal T_2$ & $\mathcal T_3$ & $t_{\rm dip}$ on $\mathcal T_3$\\
\midrule
$10^{-3}$ & .312 & .239 & .224 & .216 & $3.6\cdot10^{-5}$\\
$10^{-4}$ & .208 & .147 & .129 & .119 & $2.2\cdot10^{-5}$\\
$10^{-5}$ & .157 & .098 & .081 & .071 & $1.5\cdot10^{-5}$\\
$10^{-6}$ & .115 & .071 & .056 & .046 & $1.1\cdot10^{-5}$\\
$10^{-7}$ & .117 & .053 & .040 & .031 & $7.9\cdot10^{-6}$\\
$10^{-8}$ & .054 & .034 & .024 & .018 & $6.0\cdot10^{-6}$\\
\bottomrule
\end{tabular}
\caption{Recorded relative minimum $d_{M,\mathcal T}$ defined in \eqref{eq:dip}, for the initial-data family \eqref{eq:flowdata} with $p=16$ and $M=128$.
The columns use the baseline grid and $\mathcal T_1,\mathcal T_2,\mathcal T_3$ of \cref{tab:dipgrids}, all restricted to $0\le t\le2\cdot10^{-4}$.
The last column is the dip time on $\mathcal T_3$.
All trajectories are computed by IPM (\cref{alg:jko-robust}) with solver tolerance $10^{-12}$ and step-gap acceptance threshold $10^{-14}$; the minima are not interval-certified.
A separate baseline check doubles $M$ from 128 to 256 and leaves the three displayed digits unchanged.}
\label{tab:dip}
\end{table}

In every row of \cref{tab:dip}, $d_{M,\mathcal T}$ on $\mathcal T_3$ is smaller than on the baseline grid, and on $\mathcal T_1$, $\mathcal T_2$ and $\mathcal T_3$ it decreases as $\varepsilon$ decreases, unlike on the baseline grid.
The relative change from $\mathcal T_2$ to $\mathcal T_3$ is still between 4\% and 25\%, so the $\mathcal T_3$ values are not converged in time.
On the baseline grid, doubling $M$ leaves the three displayed digits unchanged, whereas $\mathcal T_1$ lowers $d_{M,\mathcal T}$ by at least 23\%; the time grid, not $M$, therefore dominates the baseline error in $d_{M,\mathcal T}$.
The $M$ check does not apply to the refined grids.

\Needspace{7\baselineskip}
\begin{remark}[Initial-data family and vacuum]
Varying $\varepsilon$ changes the initial datum.
At $p=16$ and $\varepsilon=0$, the profile is proportional to $\cos^{64}(\pi x)$ and has an isolated zero of order 64 at $x=1/2$.
\Cref{tab:dip} therefore concerns isolated minima; it establishes neither the filling of a vacuum interval nor finite-time vanishing from a fixed positive initial datum.
\end{remark}

\subsection{Order and scale of the initial minimum}\label{sec:inner}
We study how the decrease of the density minimum depends on the value and shape of the initial minimum. Consider
\begin{equation}\label{eq:isolated}
 \rho^0(x)=Z^{-1}\bigl(\delta+\sin^{2r}(\pi(x-1/2))\bigr),
 \qquad \delta>0,\quad r\in\{1,2,\ldots\}.
\end{equation}
Here $Z$ gives unit continuous mass, and we call $\delta$ the offset of \eqref{eq:isolated}.
The minimum is $\rho^0(1/2)=\delta/Z$.
Near $x=1/2$,
\[
 \rho^0(x)-\rho^0(1/2)\sim\frac{\pi^{2r}}{Z}|x-1/2|^{2r}.
\]
For $r=1$ the minimum is quadratic; for $r>1$ the derivatives of $\rho^0$ of orders $1,\ldots,2r-1$ vanish at $x=1/2$.
Only when $\delta=0$ does this minimum become a zero, of order $2r$.
For discrete computations, grid quadrature supplies the normalization.

\Cref{tab:isolatedsurvey} reports the relative minima and dip times for \eqref{eq:isolated} on three time grids.
Quadratic minima do not decrease at the recorded nodes and times.
For the two highest orders tested, decreasing the offset $\delta$ gives earlier and deeper minima.
For $2r=4$ and $\delta=10^{-6}$, the dip time on $\mathcal T_3$ is $3.05\cdot10^{-10}$, below the first positive times $10^{-8}$ and $10^{-9}$ of $\mathcal T_1$ and $\mathcal T_2$; on these two grids $d_{128,\mathcal T_1}=d_{128,\mathcal T_2}=1$, so they miss the decrease (\cref{tab:isolatedsurvey}).

\begin{table}[htbp]
\centering\small
\begin{tabular}{@{}rrrrrr@{}}
\toprule
Order $2r$ & $\delta$ & $d_1$ & $d_2$ & $d_3$ & $t_{\rm dip}$ on $\mathcal T_3$\\
\midrule
2 & $10^{-2}$ & 1.0000 & 1.0000 & 1.0000 & $0$\\
2 & $10^{-4}$ & 1.0000 & 1.0000 & 1.0000 & $0$\\
2 & $10^{-6}$ & 1.0000 & 1.0000 & 1.0000 & $0$\\
\addlinespace
4 & $10^{-2}$ & 0.8339 & 0.8328 & 0.8323 & $2.66\cdot10^{-6}$\\
4 & $10^{-4}$ & 0.8091 & 0.8020 & 0.8012 & $3.24\cdot10^{-8}$\\
4 & $10^{-6}$ & 1.0000 & 1.0000 & 0.8227 & $3.05\cdot10^{-10}$\\
\addlinespace
8 & $10^{-2}$ & 0.5281 & 0.5220 & 0.5188 & $2.63\cdot10^{-5}$\\
8 & $10^{-4}$ & 0.3491 & 0.3390 & 0.3334 & $2.97\cdot10^{-6}$\\
8 & $10^{-6}$ & 0.3016 & 0.2905 & 0.2844 & $2.97\cdot10^{-7}$\\
\addlinespace
16 & $10^{-2}$ & 0.4005 & 0.3908 & 0.3858 & $4.62\cdot10^{-5}$\\
16 & $10^{-4}$ & 0.1446 & 0.1297 & 0.1215 & $1.26\cdot10^{-5}$\\
16 & $10^{-6}$ & 0.0848 & 0.0721 & 0.0649 & $3.66\cdot10^{-6}$\\
\bottomrule
\end{tabular}
\caption{Initial-minimum experiment for \eqref{eq:isolated}: 36 IPM trajectories with $M=128$.
Here $d_j=d_{128,\mathcal T_j}$ is the relative minimum from \eqref{eq:dip}, evaluated at every recorded time, including $t=0$.
The time grids are $\mathcal T_1,\mathcal T_2,\mathcal T_3$ of \cref{tab:dipgrids}, continued to $t=10^{-3}$, which they reach after 236, 565 and 1317 steps.
The last column gives the dip time on the finest grid; zero means that the recorded minimum is attained initially.
The stored densities are obtained by replacing entries below $10^{-300}$ by $10^{-300}$ and normalizing to unit mass.
The runs use the rounded symbol $\widehat L_{\rm round}$ of \cref{tab:operator-rounding}, solver tolerance $10^{-12}$ and step-gap acceptance threshold $10^{-14}$.
No step failed the acceptance test; the displayed minima are not interval-certified.
The source records are listed in \cref{app:source-records}.}
\label{tab:isolatedsurvey}
\end{table}

To compare the earlier dips in \cref{tab:isolatedsurvey} as $\delta$ decreases, we rescale space and time near the initial minimum.
The two terms $\delta$ and $\sin^{2r}(\pi(x-1/2))$ in \eqref{eq:isolated} are equal when $|\sin(\pi(x-1/2))|=\delta^{1/(2r)}$, that is, for $|x-1/2|=\pi^{-1}\arcsin(\delta^{1/(2r)})\sim\delta^{1/(2r)}/\pi$ as $\delta\downarrow0$.
For a positive smooth solution $\rho$ of \eqref{eq:dlss}, write $\rho^0_{\min}=\min_x\rho^0(x)>0$ and set
\begin{equation}\label{eq:inner}
 \ell_\delta=\delta^{1/(2r)},\quad
 \mathsf y=\frac{x-1/2}{\ell_\delta},\quad \mathsf t=\frac{t}{\ell_\delta^4},\quad
 \mathsf u(\mathsf t,\mathsf y)=\frac{\rho(\ell_\delta^4\mathsf t,1/2+\ell_\delta\mathsf y)}{\rho^0_{\min}}.
\end{equation}
Substitution into \eqref{eq:dlss} gives
\[
 \rho^0_{\min}\ell_\delta^{-4}\partial_{\mathsf t}\mathsf u
 =-\rho^0_{\min}\ell_\delta^{-4}\partial_{\mathsf y\mathsf y}(\mathsf u\,\partial_{\mathsf y\mathsf y}\log\mathsf u).
\]
The four spatial derivatives contribute the factor $\ell_\delta^{-4}$, which is matched by the time scale
$\ell_\delta^4=(\delta^{1/(2r)})^4=\delta^{2/r}$.
The density amplitude $\rho^0_{\min}$ cancels, so $\mathsf u$ satisfies the same DLSS equation in $(\mathsf t,\mathsf y)$.
The rescaled circle has length $\mathsf L_\delta=1/\ell_\delta$.
The initial datum of $\mathsf u$ is
\begin{equation}\label{eq:innerdata}
 \mathsf u(0,\mathsf y)=1+\bigl[\mathsf L_\delta\sin(\pi\mathsf y/\mathsf L_\delta)\bigr]^{2r}
 \longrightarrow 1+(\pi\mathsf y)^{2r}
 \quad\text{locally as }\mathsf L_\delta\longrightarrow\infty.
\end{equation}
\Cref{app:inner} derives the space and time derivatives under this change of variables and states the regularity assumptions.
The limiting datum has infinite mass on $\R$, so this comparison concerns local density behavior near the minimum, not convergence of probability laws on $\R$.

For fixed $r$, the limit $\delta\downarrow0$ corresponds to $\mathsf L_\delta\to\infty$.
We test two properties of the computed $\mathsf u$ as $\mathsf L_\delta$ increases at fixed $r$: whether the rescaled dip time $\mathsf t_{\rm dip}$, defined below, stays bounded above and away from zero (\cref{tab:dipscaling}), and whether the minimum $c_r(\mathsf L_\delta)$ of $\mathsf u$ over a fixed rescaled time interval, also defined below, has a positive limit (\cref{tab:inner}).

Expressing the recorded dip time in the variable $\mathsf t$ of \eqref{eq:inner} gives
\[
 \mathsf t_{\rm dip}:=\frac{t_{\rm dip}}{\ell_\delta^4},
 \qquad t_{\rm dip}=\delta^{2/r}\mathsf t_{\rm dip}.
\]
The predicted exponent $2/r$ assumes that $\mathsf t_{\rm dip}$ remains bounded above and away from zero as $\delta\downarrow0$ at fixed $r$.
The change of variables alone does not establish this assumption.
\begin{table}[htbp]
\centering\small
\begin{tabular}{@{}rrr@{}}
\toprule
Order $2r$ & Measured exponent & Predicted $2/r$\\
\midrule
4 & .99 & 1\\
8 & .49 & $1/2$\\
16 & .28 & $1/4$\\
\bottomrule
\end{tabular}
\caption{Power-law exponent of the recorded dip time as the offset $\delta$ decreases in the initial data \eqref{eq:isolated}. The measured exponents are least-squares slopes of $\log t_{\rm dip}$ against $\log\delta$, using the $\mathcal T_3$ times in \cref{tab:isolatedsurvey}; the predicted exponent assumes that $\mathsf t_{\rm dip}=t_{\rm dip}/\delta^{2/r}$ remains bounded above and away from zero as $\delta\downarrow0$ at fixed $r$.}
\label{tab:dipscaling}
\end{table}

Each measured exponent, fitted to three values of $\delta$, differs from $2/r$ by at most 0.03 (\cref{tab:dipscaling}); convergence of the rescaled solutions is not established.

We use the minimum of $\mathsf u$ to test whether the density remains a positive fraction of its initial minimum on a fixed rescaled time interval as $\delta\downarrow0$.
For each finite $\mathsf L_\delta$, let $c_r(\mathsf L_\delta)$ be the minimum of the recorded values of $\mathsf u$ on $0\le\mathsf t\le0.003$.
For $2r\in\{4,8,16\}$, the rescaled dip times $t_{\rm dip}/\delta^{2/r}$ of \cref{tab:isolatedsurvey} lie between $1.1\cdot10^{-4}$ and $3.3\cdot10^{-4}$, so this interval contains these dips.
The notation suppresses dependence on the spatial and temporal grids.

Each rescaled computation solves the entropic steps \eqref{eq:step} for the data \eqref{eq:isolated} with $\delta=\mathsf L_\delta^{-2r}$ on the unit circle, by IPM (\cref{alg:jko-robust}) with the rounded symbol $\widehat L_{\rm round}$ of \cref{tab:operator-rounding}, solver tolerance $10^{-12}$ and step-gap acceptance threshold $10^{-14}$.
The grid has $M=32\mathsf L_\delta$ nodes, that is, 32 nodes per unit of $\mathsf y$, except $M=128$ for $\mathsf L_\delta=3$ and $M=256$ for $\mathsf L_\delta=16$; \cref{tab:innertime} uses $M=128$ at $\mathsf L_\delta=8$.
The 328-step grid has $\mathsf t_0=0$ and $\mathsf t_k=10^{-6}\cdot10^{(k-1)/94}$ for $1\le k\le328$; the 843-step grid has $\mathsf t_0=0$ and $\mathsf t_k=10^{-7}\cdot10^{(k-1)/188}$ for $1\le k\le843$.
Both grids end at the first $\mathsf t_k\ge3\cdot10^{-3}$, and the computation uses the physical times $t_k=\ell_\delta^4\mathsf t_k$.
The entropic step \eqref{eq:step} is homogeneous in the mass: for $\lambda>0$ and $\avg\mu=1$, $\A(\lambda\mu,\lambda m)=\lambda\A(\mu,m)$ and $\E(\lambda\mu)=\lambda\E(\mu)+\lambda\log\lambda-\lambda+1$, and the constraint $Lm=\nu-\nu^{\rm old}$ is linear.
Hence, if $(\nu^*,m^*)$ solves \eqref{eq:step} for an old density $\nu^{\rm old}$ of unit mass, then $(\lambda\nu^*,\lambda m^*)$ solves the same step for the old density $\lambda\nu^{\rm old}$.
We therefore normalize the initial profile to unit mass, solve, and multiply the computed densities by the initial mass to recover $\mathsf u$; \cref{app:inner} gives the corresponding invariance of \eqref{eq:dlss}.

\begin{table}[htbp]
\centering\small
\begin{tabular}{@{}crrrrr@{}}
\toprule
Order $2r$ & $\mathsf L_\delta=3$ & $\mathsf L_\delta=4$ & $\mathsf L_\delta=6$ & $\mathsf L_\delta=8$ & $\mathsf L_\delta=16$\\
\midrule
4  & -- & .8206 & -- & .8040 & .7997\\
6  & -- & .5360 & -- & .5027 & .4944\\
8  & -- & .3117 & -- & .2800 & .2722\\
10 & -- & .1773 & .1598 & .1525 & --\\
12 & -- & .1051 & .0830 & -- & --\\
16 & .0534 & -- & -- & -- & --\\
\bottomrule
\end{tabular}
\caption{Recorded relative minima $c_r(\mathsf L_\delta)$ of the rescaled computations over $0\le\mathsf t\le0.003$, on the 328-step grid (first positive rescaled time $\mathsf t=10^{-6}$). The grids have 32 nodes per unit of $\mathsf y$, except 16 at $\mathsf L_\delta=16$ and $M=128$ nodes at $\mathsf L_\delta=3$.
The omitted runs fail the step-gap test of \cref{alg:jko-robust} or have log differences $|\log(\nu_i/(\nu_\chi)_i)|$ of $0.02$ and $0.03$ between the computed density and the Gibbs density \eqref{eq:gibbs} at a minimum-density node, against at most $3\cdot10^{-4}$ for the retained runs (\cref{tab:innerchecks}).
Values are rounded to four decimals; they are computed values, not error bounds or limiting constants.}\label{tab:inner}
\end{table}

\begin{table}[htbp]
\centering\small
\begin{tabular}{@{}rrr@{}}
\toprule
Order $2r$ & 328-step grid & 843-step grid\\
\midrule
4 & .8040 & .8033\\
8 & .2800 & .2737\\
\bottomrule
\end{tabular}
\caption{Time sensitivity of the rescaled relative minimum $c_r(8)$ at $\mathsf L_\delta=8$, with $M=128$ nodes. The columns use the 328-step grid (first positive rescaled time $\mathsf t=10^{-6}$) and the 843-step grid (first positive rescaled time $\mathsf t=10^{-7}$). Doubling $M$ on the 328-step grid leaves the displayed values unchanged. All entries are computed values, not certified limiting constants.}
\label{tab:innertime}
\end{table}

For orders 4 and 6, $c_r$ decreases by $0.0209$ and $0.0416$ from $\mathsf L_\delta=4$ to $\mathsf L_\delta=16$, where the $\mathsf L_\delta=16$ runs have half the spatial resolution of the others (\cref{tab:inner}).
For order 10, $c_r$ decreases from $0.1773$ at $\mathsf L_\delta=4$ to $0.1525$ at $\mathsf L_\delta=8$, the largest $\mathsf L_\delta$ computed for this order.
For order 12, $c_r$ decreases from $0.1051$ at $\mathsf L_\delta=4$ to $0.0830$ at $\mathsf L_\delta=6$, and the $\mathsf L_\delta=8$ run is omitted; for order 16, only $\mathsf L_\delta=3$ is retained (\cref{tab:innerchecks}).
For order 8, time refinement lowers $c_4(8)$ by $0.0063$ (\cref{tab:innertime}), and increasing $\mathsf L_\delta$ from 8 to 16 lowers it by $0.0078$ (\cref{tab:inner}).
The latter comparison also halves the number of spatial nodes per rescaled unit.
These data therefore do not separate increasing-domain behavior from discretization error and do not support an extrapolation to a limiting constant.
The study establishes a decrease of the computed minimum on the tested grids; a limit would require common spatial resolution per rescaled unit, time refinement for each $\mathsf L_\delta$, and pointwise control of the computed minima.

\subsection{Accuracy at the smallest entries}

The absolute error bound \eqref{eq:pointwise} can exceed the smallest computed density entries.
We therefore compare the primal candidate $\nu$ with the Gibbs density obtained from its calculated multiplier $\chi$:
\begin{equation}\label{eq:gibbs}
 s(\chi)=\chi+\tfrac\tau2(L\chi)^2,\qquad
 \nu_\chi=\frac{\exp s(\chi)}{\avg{\exp s(\chi)}}.
\end{equation}
For a maximizer $\chi^*$ of $D_\tau$ in \eqref{eq:dual}, $\nu_{\chi^*}=\nu^*$ (\cref{app:gap}).
At node $i$, the log difference $|\log(\nu_i/(\nu_\chi)_i)|$ is the logarithm of the larger of the ratios $\nu_i/(\nu_\chi)_i$ and $(\nu_\chi)_i/\nu_i$.
This comparison between $\nu$ and $\nu_\chi$ does not bound either density relative to the exact minimizer without further control of $\chi$.
The candidates on the refined grids of \cref{tab:dip} and those of \cref{tab:isolatedsurvey,tab:inner,tab:innertime} are accepted by the step-gap test of \cref{alg:jko-robust} (threshold $10^{-14}$); they are not interval-certified.

\begin{table}[htbp]
\centering\small
\begin{tabular}{@{}lrrrr@{}}
\toprule
Runs & $\mathsf L_\delta$ & Relative minimum & Log difference & Failed gap tests\\
\midrule
Retained checks & -- & -- & $10^{-8}$ to $3\cdot10^{-4}$ & 0\\
Order 12 & 8 & .0132 & .03 & --\\
Order 16 & 4 & .0038 & .02 & --\\
Order 16 & 5 & -- & -- & 89 of 328\\
\bottomrule
\end{tabular}
\caption{Acceptance diagnostics for the rescaled computations. The log difference is $|\log(\nu_i/(\nu_\chi)_i)|$ at the node and time attaining $c_r(\mathsf L_\delta)$, using \eqref{eq:gibbs}; the retained-check row gives the observed range. The last three runs are excluded from \cref{tab:inner}. The two displayed rejected minima have mass-normalized entries of order $10^{-12}$. Dashes mark diagnostics not reported here.}
\label{tab:innerchecks}
\end{table}

The rejected runs of orders 12 and 16 have log differences $0.03$ and $0.02$ between $\nu$ and $\nu_\chi$, against at most $3\cdot10^{-4}$ for the retained runs, or fail the gap test (\cref{tab:innerchecks}).
The limits of $c_r(\mathsf L_\delta)$ as $\mathsf L_\delta\to\infty$ for $2r\ge12$ therefore remain unresolved.
Per-step optimization certificates concern the supplied old density; they do not control accumulated trajectory error without a stability estimate.

The absolute-error estimate \eqref{eq:pointwise} certifies a positive minimum of the exact grid minimizer when $\min_i\nu_i>M\sqrt{2\eta}$.
The bound is $\min_i\nu_i^*\ge\min_i\nu_i-M\sqrt{2\eta}$.
For $\min_i\nu_i\le M\sqrt{2\eta}$, the relative bound below requires a bound $\varepsilon_s$ on $\|s(\chi)-s(\chi^*)\|_\infty$ for an optimal multiplier $\chi^*$.
The optimal multipliers are the vectors $\chi^*+c\one$, $c\in\R$ (\cref{app:gap}), and $s(\chi^*+c\one)=s(\chi^*)+c\one$; the bound below holds for each of them.
If $\|s(\chi)-s(\chi^*)\|_\infty\le\varepsilon_s$, then
\[
 e^{-2\varepsilon_s}(\nu_\chi)_i\le\nu_i^*\le e^{2\varepsilon_s}(\nu_\chi)_i.
\]
\Cref{app:dual-exponents} proves the two-sided density bound.

The computations of \cref{tab:isolatedsurvey,tab:inner} provide no certified bound on $\|s(\chi)-s(\chi^*)\|_\infty$; the relative estimate therefore does not validate their smallest density values.

\section{Numerical conclusions}
For the uniform geodesic at $b=1$, \eqref{eq:rate} gives $\spr(T)=\max\{1/2,\cos^2\theta_{\min}\}$; at fixed $\beta$, $\theta_{\min}\to0$ as $M,N\to\infty$, so $\spr(T)\to1$ under grid refinement (\cref{app:spectrum}).
For nonuniform paths, the local action Hessians \eqref{eq:rankone} vary across cells and couple Fourier modes, so this uniform-path formula does not apply (\cref{sec:dr-choices}).
\Cref{prop:prox-error} bounds the error of the local proximal output by the scalar residual $|g_\omega(q)|$.
Interval primal--dual bounds certify the optimization accuracy of the IPM geodesics (\cref{tab:ipmvalidation,tab:positive-cost}) and of the entropic steps of \cref{tab:baseline-validation,tab:common-entropic,tab:flowcontrast}; the other reported gaps are floating-point estimates.
In the single-step tests of \cref{tab:common-entropic}, both solvers stop at an interval gap bound at most $10^{-9}$; for each of the three inputs, IPM attains the smaller bound and the smaller total time.

The DR proximal outputs contain exact zero density entries (\cref{tab:zeros}); these counts describe the stored iterates and do not certify the zero set of the minimizer.
The smallest final action gaps in the parameter sweep occur at $\beta=10^6$ or $3\cdot10^6$ (\cref{fig:betasweep}), and the tested $\beta$ requiring the fewest recorded iterations at the common action-gap threshold depends on the endpoint pair (\cref{tab:betatarget}).
At $\beta=10^6$, $b=3$ gives the fewest iterations per factor-ten residual decrease over $10^{-5}\ge r_\ell\ge10^{-9}$ for the three pairs among $b\in\{0.3,1,3,10\}$ (\cref{tab:drstep}).

The interval action bounds in \cref{tab:ipmvalidation} establish optimization accuracy on each checked grid.
Density-path accuracy remains a separate question: doubling the number of time cells from $N=64$ to $128$ changes the peaked IPM path by $\Delta_{64}=8.14\cdot10^{-2}$, whereas changing the solver on the same grid changes it by $\Delta_{\rm path}=1.28\cdot10^{-7}$ (\cref{tab:pathchecks}).
Thus agreement between DR and IPM alone is insufficient to assess path accuracy.

The zero entries of the DR proximal outputs (\cref{tab:zeros}) and the diffusivity transitions of \cref{tab:diffusivityobservation}, which narrow under spatial refinement with $|\kappa|\ge0.045$ on both sides, are consistent with the hypothesis that the limiting geodesic density vanishes near $x=0.55$, where the diffusivity changes sign; they do not prove it.
\Cref{lem:vacuumcriterion} gives a sufficient condition: if nodes $x_M^+,x_M^-\to x_0$ satisfy \eqref{eq:vacuumcriterionhyp}, then $q_M(x_M^+)+q_M(x_M^-)\to0$.
Whether the computed paths satisfy these hypotheses under joint space--time refinement, and whether their densities converge locally uniformly near the transition, remains open.

Imposing nonnegative diffusivity increases the minimum grid action of the diffuse pair by the factor $E^+_{M,N}/E_{M,N}\in[17.39,17.43]$ on the two tested grids (\cref{tab:positive-cost}).
For the convex-order pair of \cref{sec:nonnegative-diffusivity}, which satisfies $\nu^0\le_{\rm cx}\nu^N$, the factor $E^+_{M,N}/E_{M,N}$ is $\OrderConvexRatioShort$ at $(M,N)=(127,128)$; for the reversed pair, which satisfies $\nu^N\le_{\rm cx}\nu^0$, it is $\OrderReversedRatioShort$ (\cref{tab:positive-cost}).

The interval gap bounds of \cref{tab:baseline-validation,tab:flowcontrast} certify the optimization accuracy of every accepted IPM step of the baseline and contrast-test trajectories at $M=128$, for the supplied old densities.
In the trajectory timings of \cref{tab:timing}, the DR time is 1.8 to 2.2 times the IPM time under the stated stopping rules.
In \cref{tab:flowgridsize}, the IPM solve time of the first step grows from $0.033$ s at $M=128$ to $374$ s at $M=2048$, a factor $1.1\cdot10^4$ for a 16-fold increase of $M$; the test does not isolate the cost of the dense factorization, which is the presumed cause.

The ratio $\min_i\nu_i^k/\min_i\nu_i^0$ falls below 1 for the four baseline data (\cref{fig:minimum}), and on $\mathcal T_3$ the relative minimum $d_{128,\mathcal T_3}$ is below 1 for \eqref{eq:isolated} with $2r\in\{4,8,16\}$ and equal to 1 for $2r=2$, at each tested offset $\delta$ (\cref{tab:isolatedsurvey}).
Refining the time grid never increases the relative minimum in \cref{tab:dip,tab:isolatedsurvey}.
The measured dip-time exponents differ by at most $0.03$ from $2/r$, the exponent derived from \eqref{eq:inner} under the assumption that $\mathsf t_{\rm dip}$ stays bounded above and away from zero (\cref{tab:dipscaling}).
Whether the rescaled minima approach a strictly positive limit remains unresolved.
The computed values in \cref{tab:inner} require stronger pointwise validation and controlled space--time refinement before they can support such a conclusion.

\paragraph{Code and data availability.}
The repository \nolinkurl{signed-diffusive-optimal-transport} contains the solvers, the saved runs, the scripts that generate the tables and figures, and the two notebooks of \cref{app:notebooks}; it is available from the authors.
For Part~\ref{part:numerical}, its directory \nolinkurl{paper/dr-ipm-evidence/revision-20261026/} holds the saved candidates and multipliers, the dependency records, and the commands, listed in its file \nolinkurl{README.txt}, that revalidate the saved candidates and rerun the geodesic, spectral, and common-gap experiments. The records for Part~\ref{part:time-study} are described in \cref{sec:merged-records}.

\mainpart{Numerical time convergence of geodesics and entropic evolution}\label{part:time-study}
The two problems have different time discretizations.
For the geodesic problem, the path interval \([0,1]\) between the prescribed endpoint laws is divided into $N$ cells of length \(h=1/N\).
For the modified JKO scheme, $n$ steps of length \(\tau=\tfin/n\) lead from the initial density to the final time $\tfin$.
Each study fixes the spatial grid within a refinement series and reports optimization checks separately from time-mesh comparisons.
\section{Geodesic time convergence}\label{sec:geo-time-study}
\Cref{geo:thm:joint-grid,geo:cor:grid-endpoints} give convergence of the minimum actions and subsequential weak convergence of the minimizing paths under joint space--time refinement for the stated endpoint approximations.
Here the time mesh is refined at each prescribed spatial resolution.
Minimum actions and density paths are compared separately, because a small difference between two minimum actions does not bound the distance between the corresponding minimizing paths.

\subsection{Endpoints, time meshes and interval action bounds}
We use the peaked and diffuse endpoint pairs defined in \cref{sec:geodata}.
The endpoints remain fixed throughout each time-refinement series.
For each $M$, $E_{M,N}$ is the minimum value of \eqref{eq:geodesic} with $N$ time cells.
The interval evaluation of \cref{app:geointerval} gives
\begin{equation}\label{eq:geo-time-bracket}
 D_{{\rm lo},N}\le E_{M,N}\le V_{{\rm up},N},\qquad
 {\gaprel^{\rm up}}_N:=\frac{V_{{\rm up},N}-D_{{\rm lo},N}}{D_{{\rm lo},N}}.
\end{equation}
The spatial index $M$ is suppressed in $D_{{\rm lo},N}$, $V_{{\rm up},N}$ and ${\gaprel^{\rm up}}_N$ within one series.

\begin{table}[!htbp]
\centering\small
\begin{tabular}{@{}p{.24\linewidth}p{.68\linewidth}@{}}
\toprule
Quantity & Prescribed value or procedure\\
\midrule
Endpoint pairs & Endpoint densities \((\mathfrak f+10^{-2})/\avg{\mathfrak f+10^{-2}}\) with the peaked and diffuse profiles \(\mathfrak f\) of \cref{sec:geodata}.\\
Spatial grids & $M\in\{63,127\}$; periodic Fourier operator with symbol $-(2\pi k)^2$. Each time-refinement series holds $M$ fixed.\\
Time meshes & $N\in\{16,32,64,128,256,512\}$, $h=1/N$, path interval $[0,1]$.\\
Conic solve & Clarabel 0.11.1, QDLDL, one thread, tolerance $10^{-12}$, at most $400$ iterations and $180$ seconds per solve; equilibration off. No warm start.\\
Candidate selection & Clipping, positive mixtures and positive slacks, as in \cref{alg:geo-candidates}; choose the finite candidate with the least action computed with the momentum \eqref{eq:geovalue}.\\
Acceptance & $D_{{\rm lo},N}>0$ and ${\gaprel^{\rm up}}_N\le10^{-8}$, evaluated with 40-decimal interval arithmetic using exact $\pi$ and the saved density rows normalized exactly to unit mass $\avg{\nu}=1$.\\
Tighter repeat & Repeat $N=512$ at tolerance $10^{-13}$ for each pair and spatial grid; use it only to check the selected path's sensitivity to the solver tolerance.\\
\bottomrule
\end{tabular}
\caption{Geodesic refinement settings. There are six time meshes for each of two endpoint pairs on each of two spatial grids, plus four tighter repeats. Binary64 samples define the endpoints; their exact normalizations to unit mass $\avg{\nu}=1$ define the problems in the interval checks. The conic assembly uses the Fourier symbol $\widehat L$ evaluated in binary64, while interval evaluation uses the exact symbol $\widehat L$.}
\label{tab:geo-time-protocol}
\end{table}

\begin{procedure}{Geodesic time refinement}{alg:geo-time-study}
\alginput{the endpoint pairs, spatial grids, time meshes, solver tolerances and interval threshold in \cref{tab:geo-time-protocol}.}
\algoutput{raw conic arrays, selected density paths, action intervals, successive-mesh comparisons and tighter-repeat differences.}
\algstep{\textbf{For each} endpoint pair and $M$, sample the two endpoints and divide each by its mass, so that it has unit mass $\avg{\nu}=1$. Fix these arrays for the complete time series.}
\algstep[1]{\textbf{For each} $N$, assemble the inverse-operator formulation in \cref{alg:geo-ipm}. Solve at tolerance $10^{-12}$ without warm start.}
\algstep[2]{Save the primal variables, multipliers, cone slacks and diagnostics. Apply \cref{alg:geo-candidates}; retain every tested mixture weight and action. \textbf{Stop} if no candidate has finite action.}
\algstep[2]{Replace each selected row \(\nu_{\rm IPM}^j\) by \(\nu_{\rm IPM}^j/\avg{\nu_{\rm IPM}^j}\) in exact arithmetic and evaluate \eqref{eq:geo-time-bracket} by \cref{app:geointerval}. Reconstruct the least-action momentum during this evaluation.}
\algstep[2]{\textbf{If} interval acceptance fails, repeat once with tolerance $10^{-13}$, static regularization constant $10^{-9}$, and iterative-refinement absolute and relative tolerances $10^{-13}$. Save both attempts. \textbf{Stop} if this repeat also fails.}
\algstep[2]{Save the accepted candidate, its interval action bounds and its change from the raw density. A stopping status alone does not imply acceptance.}
\algstep[1]{\textbf{End for.} Compare $N$ and $2N$ at their common time nodes using \eqref{eq:geo-time-differences}.}
\algstep[1]{Repeat $N=512$ at tolerance $10^{-13}$, with the same acceptance and fallback. Record \(R_{512}=\max_{0\le j\le512}\|\bar\nu_{(512)}^j-\bar\nu_{(512),\rm rep}^j\|_{1,M}\), where \(\bar\nu_{(512),\rm rep}\) is the selected path of the repeat at tolerance \(10^{-13}\).}
\algstep{\textbf{End for.} Generate the tables and plots from the saved arrays and interval bounds. Keep all primary runs in the refinement series; the tighter repeat does not replace its finest member.}
\end{procedure}

\subsection{Action and density-path refinement}
Let $\bar\nu_{(N)}$ be the selected density path after row normalization.
Define the successive-mesh differences
\begin{equation}\label{eq:geo-time-differences}
 \Delta^V_N=|V_{{\rm up},2N}-V_{{\rm up},N}|,\qquad
 \Delta_N=\max_{0\le j\le N}
 \|\bar\nu_{(N)}^j-\bar\nu_{(2N)}^{2j}\|_{1,M}.
\end{equation}
No reference geodesic is assumed known.
The interval bounds distinguish a change of the discrete minimum from optimization error:
\begin{equation}\label{eq:geo-time-action-change}
 D_{{\rm lo},2N}-V_{{\rm up},N}\le E_{M,2N}-E_{M,N}\le V_{{\rm up},2N}-D_{{\rm lo},N}.
\end{equation}
For the density comparison, we also record the modification made after the conic solve,
\begin{equation}\label{eq:geo-time-correction}
 c_N=\max_{0\le j\le N}
 \|\bar\nu_{(N)}^j-\nu_{(N),\mathrm{raw}}^j\|_{1,M}.
\end{equation}

\begin{table}[!htbp]
\centering\small
\begin{tabular}{@{}llrrrrrr@{}}
\toprule
Pair & $N$ & $\Delta^V_N/V_{{\rm up},2N}$ & $D_{{\rm lo},2N}-V_{{\rm up},N}$ & $r_A$ & $\Delta_N$ & $j_\Delta$ & $r_\nu$\\
\midrule
Diffuse & 16 & $6.66\cdot10^{-5}$ & $4.21\cdot10^{-9}$ & --- & $1.24\cdot10^{-2}$ & 2 & ---\\
Diffuse & 32 & $1.25\cdot10^{-5}$ & $7.89\cdot10^{-10}$ & 2.42 & $3.62\cdot10^{-3}$ & 4 & 1.77\\
Diffuse & 64 & $2.46\cdot10^{-6}$ & $1.56\cdot10^{-10}$ & 2.34 & $5.75\cdot10^{-4}$ & 25 & 2.65\\
Diffuse & 128 & $6.15\cdot10^{-7}$ & $3.89\cdot10^{-11}$ & 2.00 & $2.24\cdot10^{-4}$ & 28 & 1.36\\
Diffuse & 256 & $1.51\cdot10^{-7}$ & $9.50\cdot10^{-12}$ & 2.02 & $9.97\cdot10^{-5}$ & 72 & 1.17\\
\midrule
Peaked & 16 & $7.47\cdot10^{-4}$ & $1.15\cdot10^{-6}$ & --- & $3.81\cdot10^{-1}$ & 15 & ---\\
Peaked & 32 & $2.12\cdot10^{-4}$ & $3.25\cdot10^{-7}$ & 1.82 & $1.95\cdot10^{-1}$ & 31 & 0.97\\
Peaked & 64 & $4.83\cdot10^{-5}$ & $7.41\cdot10^{-8}$ & 2.14 & $8.17\cdot10^{-2}$ & 1 & 1.25\\
Peaked & 128 & $1.20\cdot10^{-5}$ & $1.84\cdot10^{-8}$ & 2.01 & $2.79\cdot10^{-2}$ & 1 & 1.55\\
Peaked & 256 & $2.88\cdot10^{-6}$ & $4.42\cdot10^{-9}$ & 2.05 & $7.05\cdot10^{-3}$ & 255 & 1.98\\
\bottomrule
\end{tabular}
\caption{Geodesic time refinement at $M=127$. Each row compares $N$ and $2N$ time cells with the same endpoints and spatial operator. The action change $\Delta^V_N$ and common-node path difference $\Delta_N$ are defined in \eqref{eq:geo-time-differences}. The column $D_{{\rm lo},2N}-V_{{\rm up},N}$ is the lower bound on $E_{M,2N}-E_{M,N}$ in \eqref{eq:geo-time-action-change}, and $j_\Delta$ is the index $j$ at which the maximum in $\Delta_N$ is attained. The empirical orders are $r_A=\log_2(\Delta^V_{N/2}/\Delta^V_N)$ and $r_\nu=\log_2(\Delta_{N/2}/\Delta_N)$; their first values require an uncomputed coarser mesh. The final row of each pair compares $N=256$ with $N=512$. These are successive-mesh differences, not errors relative to an exact geodesic.}
\label{tab:geo-time-convergence}
\end{table}

\begin{table}[!htbp]
\centering\small
\begin{tabular}{@{}llrrrr@{}}
\toprule
Pair & $M$ & $V_{{\rm up},512}$ & $\max_N{\gaprel^{\rm up}}_N$ & $\max_N c_N$ & $R_{512}$\\
\midrule
Diffuse & 63 & $6.3211547\cdot10^{-5}$ & $5.56\cdot10^{-10}$ & $2.13\cdot10^{-15}$ & $7.76\cdot10^{-8}$\\
Diffuse & 127 & $6.3234389\cdot10^{-5}$ & $9.61\cdot10^{-10}$ & $1.31\cdot10^{-15}$ & $9.59\cdot10^{-8}$\\
Peaked & 63 & $1.5339424\cdot10^{-3}$ & $1.49\cdot10^{-10}$ & $1.07\cdot10^{-15}$ & $0$\\
Peaked & 127 & $1.5341426\cdot10^{-3}$ & $5.40\cdot10^{-10}$ & $4.81\cdot10^{-15}$ & $9.29\cdot10^{-8}$\\
\bottomrule
\end{tabular}
\caption{Optimization checks for the geodesic refinement series. $V_{{\rm up},512}$ is the interval upper action bound on the finest time mesh. The maximum relative interval gap ${\gaprel^{\rm up}}_N=(V_{{\rm up},N}-D_{{\rm lo},N})/D_{{\rm lo},N}$ and the maximum path modification $c_N$ in \eqref{eq:geo-time-correction} include all six time meshes and the tighter repeat. $R_{512}$ is the maximum common-time grid $L^1$ difference between the selected paths at conic tolerances $10^{-12}$ and $10^{-13}$, with $N=512$. For the peaked pair at $M=63$, the solves at tolerances $10^{-12}$ and $10^{-13}$ both stop after 22 iterations, so the selected paths coincide and $R_{512}=0$. None of the 28 solves (24 primary solves and 4 tighter repeats) used the fallback of step~5 of \cref{alg:geo-time-study}. The gap bounds certify action accuracy, and the repeat tests numerical stability of the selected density path.}
\label{tab:geo-time-validation}
\end{table}

\begin{figure}[!htbp]
\centering
\includegraphics[width=\linewidth]{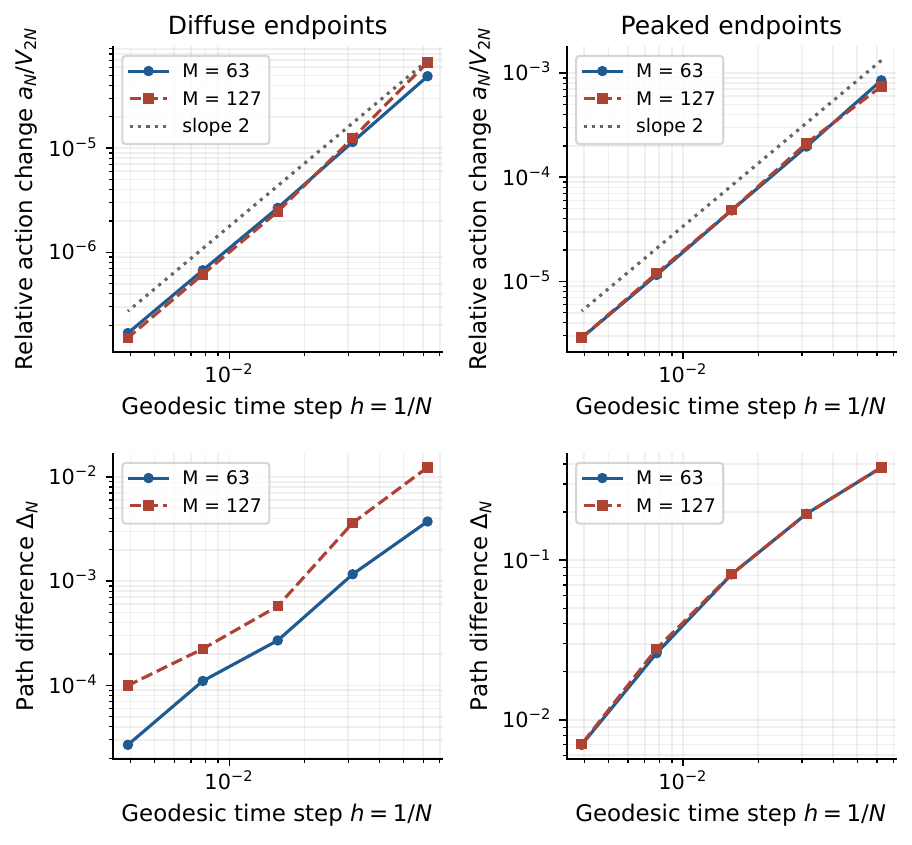}
\caption{Time refinement for the geodesic problem with prescribed endpoint laws. Top row: relative change of the upper action values, $\Delta^V_N/V_{{\rm up},2N}$. Bottom row: maximum grid $L^1$ difference $\Delta_N$ at common time nodes. Each curve fixes the spatial grid and compares $N$ with $2N$, for $N=16,32,64,128,256$. The dotted lines have slope two, scaled to $1.8$ times the finest plotted $M=127$ action change; they are comparison slopes and supply no theoretical rate. The interval action bounds and tolerance-repeat checks are in \cref{tab:geo-time-validation}.}
\label{fig:geo-time-convergence}
\end{figure}

For \(N\ge128\), \(2.00\le r_A\le2.05\) for both pairs in \cref{tab:geo-time-convergence}, consistent with \(\Delta^V_N=O(h^2)\) on these meshes (\cref{fig:geo-time-convergence}).
In every comparison, at \(M=63\) and at \(M=127\), \(D_{{\rm lo},2N}-V_{{\rm up},N}>0\) (the values at \(M=127\) are in \cref{tab:geo-time-convergence}); hence \(E_{M,2N}>E_{M,N}\) by \eqref{eq:geo-time-action-change}.
The minimum action therefore increases under each refinement, and this ordering holds for the exact minima whatever the optimization error within the interval bounds.

\begin{figure}[!htbp]
\centering
\includegraphics[width=\linewidth]{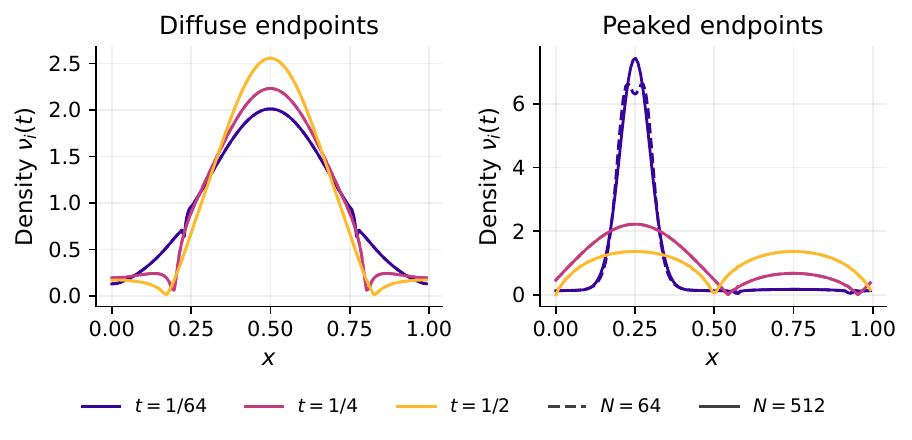}
\caption{Selected geodesic density paths at $M=127$, sampled at $t=1/64,1/4,1/2$. Colors identify the path time; dashed curves use $N=64$ and solid curves use $N=512$. The earliest slice is the first interior node of the coarser mesh. These are the candidates used in the action evaluation, after row normalization, rather than raw conic densities. The maximum differences over all common time nodes are measured in \cref{tab:geo-time-convergence}.}
\label{fig:geo-time-profiles}
\end{figure}

The common-time density differences \(\Delta_N\) also decrease for both pairs.
In \cref{tab:geo-time-convergence}, \(r_\nu\) ranges from 0.97 to 2.65, whereas \(r_A\) ranges from 1.82 to 2.42.
At \(M=127\) and each \(N\), \(\Delta_N\) of the peaked pair is 31 to 142 times \(\Delta_N\) of the diffuse pair.
For the peaked pair at $M=127$, the largest differences between paths on the meshes $N$ and $2N$ occur at $j_\Delta=1$ or $j_\Delta=N-1$, the time nodes adjacent to $t=0$ and $t=1$ (\cref{tab:geo-time-convergence}); \cref{fig:geo-time-profiles} includes an early slice to show this part of the path.
At \(M=127\), \(\max_Nc_N\le4.81\cdot10^{-15}\) and \(R_{512}\le9.59\cdot10^{-8}\) in \cref{tab:geo-time-validation}, whereas \(\Delta_{256}\ge9.97\cdot10^{-5}\) in \cref{tab:geo-time-convergence}.
The post-solve modification and the effect of the solver tolerance are therefore below \(10^{-3}\Delta_{256}\), so \(\Delta_N\) records the change of the selected path between time meshes.

Thus, at \(M=127\), the interval bounds show that the action changes \(\Delta^V_N\) exceed the optimization error, and both \(\Delta^V_N\) and \(\Delta_N\) decrease under each refinement in \cref{tab:geo-time-convergence}.
\Cref{geo:thm:time-discrete-convergence} gives no rate, so the observed orders \(r_A\) and \(r_\nu\) are numerical observations only.

\begin{remark}[Scope of the geodesic experiment]
The interval bounds in \cref{tab:geo-time-validation} certify the optimization accuracy of the action on each grid.
The comparisons \(\Delta_N\) and \(R_{512}\) concern the selected paths: they do not prove uniqueness of the minimizing density path of \eqref{eq:geodesic} and do not bound the distance from the selected path to a minimizing path.
The series at \(M=63\) and \(M=127\) are separate time refinements and do not establish a joint limit \((h,\Delta x)\to(0,0)\), the convergence of the signed diffusivity \(\kappa\) of \eqref{eq:kappa}, or pointwise convergence of the density where it tends to zero.
\end{remark}

\newcommand{\MergedRunCount}{25}
\newcommand{\MergedStepCount}{2480}
\newcommand{\MergedFallbackCount}{11}
\newcommand{\MergedSecondLookCount}{37}

\section{Time convergence of the modified JKO scheme}\label{sec:merged-refinement}
\subsection{Initial data, time meshes and fixed-grid ODE}
The conditional theorem in \cref{sec:merged-jko} concerns continuous-space minimizers.
The present study refines the modified JKO scheme, whose steps are \eqref{eq:step}, with the spatial operator, initial density and final time fixed in each series.
For the smooth positive profile, the computed trajectories are compared with the solution of the fixed-grid ODE \eqref{eq:merged-ode} below, whose right-hand side is the first-order term in \cref{prop:step-consistency}.
For the peaked profile, the study records how the relative minimum \(d_{M,\mathcal T}\) of \eqref{eq:dip} changes under time refinement.

For each grid, the fixed-grid ODE is
\begin{equation}\label{eq:merged-ode}
 \dot u_M=-L\bigl(u_M L\log u_M\bigr),\qquad u_M(0)=\nu^0.
\end{equation}
Products and logarithms are componentwise.
We use the binary64 Fourier symbol $\widehat L(k)=-(2\pi k)^2$, including the real Nyquist mode for even $M$.
Thus the ODE reference and the entropic steps use the same specified grid operator.
To obtain the ODE reference independently of the entropic steps, we solve \eqref{eq:merged-ode} with Radau IIA, a fifth-order implicit Runge--Kutta method, using SciPy's \texttt{Radau} solver \cite{SciPyRadau}.
We use both tolerance pairs in \cref{tab:merged-protocol} and retain the tighter solve as the ODE reference.
The difference between these two ODE solutions measures sensitivity to the solver tolerances; it does not give a rigorous bound on the error of the ODE reference.

\begin{table}[htbp]
\centering\small
\begin{tabular}{@{}p{.25\linewidth}p{.67\linewidth}@{}}
\toprule
Quantity & Prescribed value or procedure\\
\midrule
Smooth initial profile & $\nu_i^0=1+\tfrac12\cos(2\pi x_i)$; $M\in\{32,64,128\}$; $\tfin=10^{-4}$.\\
Peaked initial profile & \(\nu^0\) from \eqref{eq:flowdata} with \(p=16\) and \(\varepsilon=10^{-3}\); \(M\in\{64,128\}\); \(\tfin=8\cdot10^{-5}\).\\
Time meshes & $n\in\{16,32,64,128,256\}$ equal steps; $t_k=k\tfin/n$ and $\tau=\tfin/n$. Stored binary64 time differences are supplied to each solve.\\
Entropic solver & Clarabel 0.11.1; one thread; tolerance $10^{-12}$; the seven-variant fallback in \cref{alg:jko-robust}.\\
Acceptance & Each accepted step must have a 133-bit interval upper gap at most $10^{-14}$. A failure stops the trajectory.\\
ODE reference & SciPy 1.15.3 \texttt{Radau} (Radau IIA), analytic Jacobian; relative/absolute tolerances $(10^{-10},10^{-12})$ and $(10^{-12},10^{-14})$; the tighter solve defines the ODE reference.\\
Recorded arrays & All density nodes and multipliers, the operator symbol, time nodes, interval bounds, and attempted solver variants.\\
\bottomrule
\end{tabular}
\caption{Settings of the time-refinement study. Initial vectors and each entropic output are normalized to unit mass $\avg{\nu}=1$. The entropic solver replaces raw density entries below $10^{-300}$ by $10^{-300}$ before normalization. Every raw density entry of an accepted step was at least $10^{-300}$, so this replacement did not change the density of any accepted step. The interval check accounts for the exact mass of the stored input.}
\label{tab:merged-protocol}
\end{table}

\begin{procedure}{Time refinement and ODE comparison}{alg:merged-refinement}
\alginput{the profiles, spatial grids, final times, time meshes and tolerances in \cref{tab:merged-protocol}.}
\algoutput{saved trajectories, step validation records, ODE references, and the differences in \eqref{eq:merged-errors}.}
\algstep{\textbf{For each} profile and spatial grid, sample and normalize $\nu^0$; construct the Fourier symbol and the corresponding matrix $L$.}
\algstep[1]{Solve \eqref{eq:merged-ode} with Radau IIA for each tolerance pair in \cref{tab:merged-protocol}. Supply the analytic Jacobian
$-L\,\diag(L\log u)-L\,\diag(u)L\,\diag(u^{-1})$.
\textbf{Stop} if the integration fails or a stage has a nonpositive density; do not clip ODE stages. Save both Radau solutions at the times \(k\tfin/256\), \(0\le k\le256\), which contain the nodes of every mesh \(n\in\{16,\ldots,256\}\).}
\algstep[1]{\textbf{For each} $n\in\{16,32,64,128,256\}$, set $t_k=k\tfin/n$ and initialize the density to $\nu^0$.}
\algstep[2]{\textbf{For} $k=0,\ldots,n-1$, call \cref{alg:jko-robust} with the current density and $t_{k+1}-t_k$. Retain every attempted variant.}
\algstep[3]{Call \cref{alg:jko-interval} even if the floating-point acceptance test passed. \textbf{Stop} on a nonfinite interval upper gap or a gap exceeding $10^{-14}$.}
\algstep[3]{Save the normalized density, multiplier, interval check and attempt log; use this density as the next input.}
\algstep[2]{\textbf{End for.} Compare the density nodes with the ODE reference at the same times; save the minimum and entropy at all computed nodes.}
\algstep[1]{\textbf{End for.} Compare consecutive time meshes at their common nodes.}
\algstep{\textbf{End for.} Compare nested spatial grids at the same time nodes by restricting the fine-grid vectors to the coarse-grid nodes. Generate all reported tables and figures from the saved arrays.}
\end{procedure}

\subsection{Trajectory differences and observed orders}
For the grid with $M$ nodes and $\tau=\tfin/n$, let $\nu_{M,\tau}^k$ be the density after $k$ entropic steps \eqref{eq:step} started from $\nu^0$; the index $M$ or $\tau$ is omitted when it is fixed.
At the common time nodes, define
\begin{equation}\label{eq:merged-errors}
\begin{aligned}
 \mathrm{err}_{\rm f}(M,\tau)&=\|\nu_{M,\tau}^n-u_M^{\rm ref}(\tfin)\|_{1,M},\\
 \mathrm{err}_{\rm path}(M,\tau)&=\max_{0\le k\le n}\|\nu_{M,\tau}^k-u_M^{\rm ref}(t_k)\|_{1,M},\\
 \Delta_\tau(M)&=\max_{0\le k\le n}\|\nu_{M,\tau}^{k}-\nu_{M,\tau/2}^{2k}\|_{1,M}.
\end{aligned}
\end{equation}
Here $u_M^{\rm ref}$ denotes the Radau solution of \eqref{eq:merged-ode} with the tighter tolerances of \cref{tab:merged-protocol}, and $\|z\|_{1,M}$ is the norm of \cref{sec:grid-notation}.
The difference $\Delta_\tau(M)$ compares two entropic trajectories only, so it does not depend on the ODE integrator.

\begin{table}[htbp]
\centering\small
\begin{tabular}{@{}lrrrrrrr@{}}
\toprule
Profile & $n$ & $\tau$ & $\mathrm{err}_{\rm f}$ & $r_{\rm f}$ & $\mathrm{err}_{\rm path}$ & $r_{\rm path}$ & $\Delta_{\tau}$\\
\midrule
Smooth & 16 & $6.25\cdot10^{-6}$ & $2.03\cdot10^{-4}$ & --- & $4.59\cdot10^{-4}$ & --- & $2.21\cdot10^{-4}$\\
Smooth & 32 & $3.13\cdot10^{-6}$ & $1.01\cdot10^{-4}$ & 1.01 & $2.38\cdot10^{-4}$ & 0.95 & $1.17\cdot10^{-4}$\\
Smooth & 64 & $1.56\cdot10^{-6}$ & $5.03\cdot10^{-5}$ & 1.00 & $1.21\cdot10^{-4}$ & 0.97 & $6.00\cdot10^{-5}$\\
Smooth & 128 & $7.81\cdot10^{-7}$ & $2.51\cdot10^{-5}$ & 1.00 & $6.12\cdot10^{-5}$ & 0.99 & $3.05\cdot10^{-5}$\\
Smooth & 256 & $3.91\cdot10^{-7}$ & $1.25\cdot10^{-5}$ & 1.00 & $3.08\cdot10^{-5}$ & 0.99 & ---\\
\midrule
Peaked & 16 & $5.00\cdot10^{-6}$ & $1.76\cdot10^{-2}$ & --- & $1.08\cdot10^{-1}$ & --- & $4.86\cdot10^{-2}$\\
Peaked & 32 & $2.50\cdot10^{-6}$ & $9.32\cdot10^{-3}$ & 0.92 & $7.59\cdot10^{-2}$ & 0.51 & $3.39\cdot10^{-2}$\\
Peaked & 64 & $1.25\cdot10^{-6}$ & $4.85\cdot10^{-3}$ & 0.94 & $4.70\cdot10^{-2}$ & 0.69 & $2.09\cdot10^{-2}$\\
Peaked & 128 & $6.25\cdot10^{-7}$ & $2.49\cdot10^{-3}$ & 0.96 & $2.61\cdot10^{-2}$ & 0.85 & $1.24\cdot10^{-2}$\\
Peaked & 256 & $3.13\cdot10^{-7}$ & $1.26\cdot10^{-3}$ & 0.98 & $1.43\cdot10^{-2}$ & 0.87 & ---\\
\bottomrule
\end{tabular}
\caption{Time refinement at $M=128$. The final-time difference $\mathrm{err}_{\rm f}$, the maximum difference at computed times $\mathrm{err}_{\rm path}$, and the difference $\Delta_\tau$ between consecutive time meshes are defined in \eqref{eq:merged-errors}. The observed orders are $r_{\rm f}=\log_2(\mathrm{err}_{\rm f}(2\tau)/\mathrm{err}_{\rm f}(\tau))$ and $r_{\rm path}=\log_2(\mathrm{err}_{\rm path}(2\tau)/\mathrm{err}_{\rm path}(\tau))$. The ODE reference $u_M^{\rm ref}$ is computed with the tighter tolerances in \cref{tab:merged-protocol}. No reference at a finer entropic mesh is available in the last row of each profile.}
\label{tab:merged-convergence}
\end{table}

\begin{figure}[htbp]
\centering
\includegraphics[width=\linewidth]{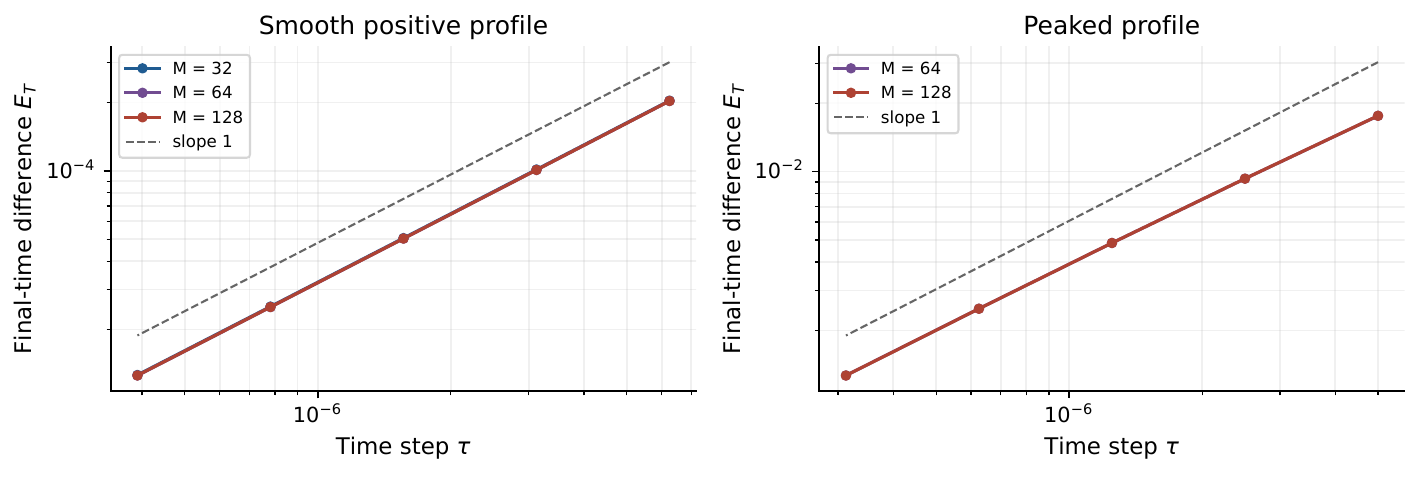}
\caption{Final-time differences $\mathrm{err}_{\rm f}$ under time refinement for the two profiles in \cref{tab:merged-protocol}. Each colored curve holds $M$ fixed; at each $n$, the plotted values of $\mathrm{err}_{\rm f}$ for different $M$ differ by less than $0.6\%$ of the smallest of them. The gray dashed lines have slope one and are scaled to $1.5$ times the finest-mesh difference at $M=128$; they are comparison slopes, not fitted error bounds. Each entropic step is accepted only if its 133-bit interval upper gap is at most $10^{-14}$. The ODE reference is the Radau solve with relative and absolute tolerances $(10^{-12},10^{-14})$; its difference from the $(10^{-10},10^{-12})$ solve is in \cref{tab:merged-validation}.}
\label{fig:merged-convergence}
\end{figure}

\Needspace{8\baselineskip}
In \cref{tab:merged-convergence}, the observed orders of $\mathrm{err}_{\rm f}$, $\mathrm{err}_{\rm path}$ and $\Delta_\tau$ for the smooth profile lie between 0.92 and 1.01.
For the peaked profile, $r_{\rm f}$ increases from 0.92 to 0.98, whereas $r_{\rm path}$ increases from 0.51 to 0.87.
In the right panel of \cref{fig:merged-minima}, the largest peaked-profile difference $\|\nu_{M,\tau}^k-u_M^{\rm ref}(t_k)\|_{1,M}$ occurs at a time $t_k\le\tfin/16$ for every $n$.
For this profile, the ratio $\mathrm{err}_{\rm path}/\mathrm{err}_{\rm f}$ increases from 6.1 at $n=16$ to 11.3 at $n=256$ (\cref{tab:merged-convergence}), so $\mathrm{err}_{\rm f}$ alone understates the maximum difference along the trajectory.

\begin{table}[htbp]
\centering\small
\begin{tabular}{@{}lrrrr@{}}
\toprule
Profile & $M$ & Maximum interval gap & Maximum step $L^1$ bound & ODE tolerance difference\\
\midrule
Smooth & 32 & $1.09\cdot10^{-16}$ & $1.48\cdot10^{-8}$ & $9.28\cdot10^{-12}$\\
Smooth & 64 & $1.26\cdot10^{-16}$ & $1.59\cdot10^{-8}$ & $9.32\cdot10^{-12}$\\
Smooth & 128 & $1.59\cdot10^{-16}$ & $1.78\cdot10^{-8}$ & $9.31\cdot10^{-12}$\\
Peaked & 64 & $1.13\cdot10^{-17}$ & $4.76\cdot10^{-9}$ & $7.49\cdot10^{-12}$\\
Peaked & 128 & $3.37\cdot10^{-17}$ & $8.22\cdot10^{-9}$ & $5.50\cdot10^{-12}$\\
\bottomrule
\end{tabular}
\caption{Accuracy checks over all five time meshes for each profile and spatial grid. For each step, $\varsigma=\avg{\nu^{\rm old}}$ is the exact mass of the recorded input, $\nu_{\rm cand}$ is the positive candidate density before mass normalization, $V$ is the objective of \eqref{eq:step} at $\varsigma\nu_{\rm cand}/\avg{\nu_{\rm cand}}$, and $D_{\tau,\varsigma}$ is the dual lower bound \eqref{eq:dualmass}; $V$ and $D_{\tau,\varsigma}$ are evaluated in interval arithmetic by \cref{alg:jko-interval}. Maximum interval gap: the largest $\eta=\sup(V-D_{\tau,\varsigma})$ over all steps. Maximum step $L^1$ bound: the largest $\sqrt{2\varsigma\eta}+|\avg{\nu_{\rm cand}}-\varsigma|$, which bounds $\normone{\nu_{\rm cand}-\nu^*}$ for the recorded input (\cref{app:mass}). Both quantities concern one step with its recorded input density; they are not trajectory error bounds. The last column is the maximum grid $L^1$ difference between the two Radau reference integrations at the 257 common times. All \MergedStepCount\ entropic steps satisfy the prescribed interval acceptance test. A fallback event is a step for which the first variant of \cref{alg:jko-robust} was rejected. The logs record \MergedFallbackCount\ fallback events and \MergedSecondLookCount\ cases in which a candidate passed the gap test without solver status \texttt{Solved}, so that \cref{alg:jko-robust} attempted the next variant once.}
\label{tab:merged-validation}
\end{table}

\Needspace{6\baselineskip}
The step $L^1$ bounds in \cref{tab:merged-validation} are at most $1.78\cdot10^{-8}$, whereas the smallest $\Delta_\tau$ in \cref{tab:merged-convergence} is $3.05\cdot10^{-5}$.
The step bounds concern each optimization problem with its recorded input; $\mathrm{err}_{\rm f}$, $\mathrm{err}_{\rm path}$ and $\Delta_\tau$ compare complete trajectories.
A trajectory bound would require a stability estimate of the entropic step with respect to its input density.
No such estimate is proved here, so no trajectory bound is derived.

\subsection{Density minima and spatial sensitivity}
For a fixed time mesh, define the nested-grid comparison
\begin{equation}\label{eq:merged-spatial}
 S_M=\max_{0\le k\le n}\frac1M\sum_{i=0}^{M-1}
 \left|\nu_{M,i}^{k}-\nu_{2M,2i}^{k}\right|.
\end{equation}
Node $2i$ of the grid with $2M$ nodes is the point $x_i=i/M$, so $S_M$ compares the two solutions at the same points and the same times.
The two grids use the same analytic initial profile and separate grid normalizations.

\begin{table}[htbp]
\centering\small
\begin{tabular}{@{}lrrrr@{}}
\toprule
Profile & $M$ & Relative minimum, computed & Relative minimum, ODE reference & $S_M$\\
\midrule
Smooth & 32 & 1.000000 & 1.000000 & $7.81\cdot10^{-9}$\\
Smooth & 64 & 1.000000 & 1.000000 & $5.95\cdot10^{-9}$\\
Smooth & 128 & 1.000000 & 1.000000 & ---\\
Peaked & 64 & 0.209587 & 0.206598 & $3.04\cdot10^{-7}$\\
Peaked & 128 & 0.209587 & 0.206598 & ---\\
\bottomrule
\end{tabular}
\caption{Spatial comparison and density minima at $n=256$. Each relative minimum divides the smallest sampled density, including the initial time, by the initial minimum on that grid. The ODE reference is sampled at the same times. The nested-grid difference $S_M$ is defined in \eqref{eq:merged-spatial}; no grid with $M=256$ was computed. This comparison tests spatial sensitivity at a prescribed time step and does not bound continuum discretization error.}
\label{tab:merged-minima}
\end{table}

\begin{figure}[htbp]
\centering
\includegraphics[width=\linewidth]{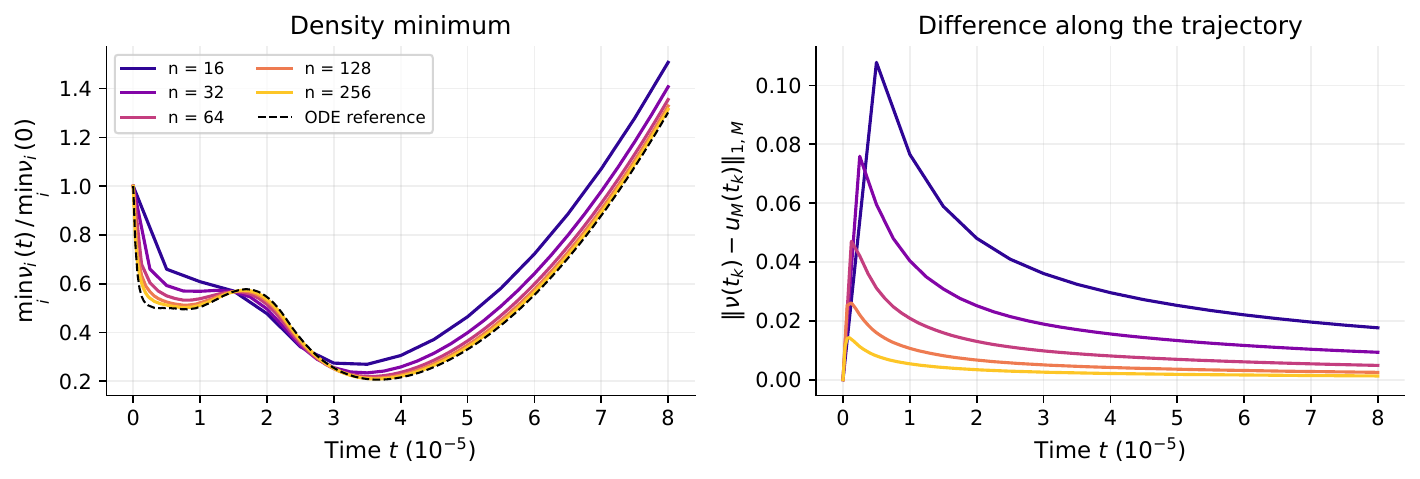}
\caption{Peaked-profile evolution at $M=128$. Left: the minimum at every computed time, divided by the initial minimum; the dashed curve is the tighter ODE reference sampled at $257$ times. Right: the grid $L^1$ difference from that reference at each entropic time node. Colors identify the number $n$ of equal steps in both panels. Lines join recorded values and do not certify extrema between samples.}
\label{fig:merged-minima}
\end{figure}

For the peaked profile at $M=128$, the relative minimum $d_{M,\mathcal T}$ of \eqref{eq:dip} on the time grid $\mathcal T=\{t_0,\ldots,t_n\}$ decreases as $n$ increases (left panel of \cref{fig:merged-minima}).
At $n=256$, $d_{M,\mathcal T}$ exceeds the relative minimum of the ODE reference by $3.0\cdot10^{-3}$ for both $M=64$ and $M=128$ (\cref{tab:merged-minima}).
For each profile, $\max_M S_M$ at $n=256$ (\cref{tab:merged-minima}) is smaller than $\mathrm{err}_{\rm path}$ at $M=128$, $n=256$ (\cref{tab:merged-convergence}) by a factor larger than $3\cdot10^3$.
On these grids, the difference between the grids $M$ and $2M$ at fixed $\tau$ is therefore smaller than the difference from the ODE reference at fixed $M$; this comparison does not bound the spatial discretization error.

On the tested grids and final times, $\mathrm{err}_{\rm f}$ decreases with observed order between 0.92 and 1.01 for both profiles, in agreement with the first-order expansion of \cref{prop:step-consistency}.
For the peaked profile, $\mathrm{err}_{\rm path}$ decreases with observed order between 0.51 and 0.87.
At $n=256$, the relative minimum of the peaked profile exceeds that of the ODE reference by $3.0\cdot10^{-3}$, and, for both profiles, $\max_M S_M$ is smaller than $\mathrm{err}_{\rm path}$ at $M=128$ by more than three orders of magnitude.
These are observations on finite grids; they prove neither the positivity floor \eqref{jko:eq:floor} nor a space--time limit.

\Needspace{8\baselineskip}
\begin{remark}[Relation to the positivity floor]
The observed trajectories remain positive on the tested time interval.
This does not establish the positivity floor \eqref{jko:eq:floor}, which must hold for all sufficiently small time steps in continuous space.
In \cref{tab:merged-minima}, the relative minimum at $n=256$ is $0.209587$ for the peaked profile, and that of the ODE reference is $0.206598$; \cref{fig:merged-minima} shows the ratio $\min_i\nu_i^k/\min_i\nu_i^0$ at every computed time.
These values are consistent with \eqref{jko:eq:floor} on the tested grids, but they neither prove a lower bound uniform in $\tau$ and $M$ nor exclude vacuum for other data or longer times.
\end{remark}

\begin{remark}[Reproduction records]\label{sec:merged-records}
The directory \nolinkurl{paper/merged-evidence} contains both time-refinement studies, their solver snapshots, saved arrays, interval checks, and environment records.
The driver \nolinkurl{run_geodesic_refinement.py} records every geodesic candidate modification and tolerance repeat; \nolinkurl{run_time_refinement.py} records every solver attempt and fallback of the modified JKO scheme.
The corresponding figure scripts generate the tables and plots from these records.
The records for Part~\ref{part:numerical} are described in \cref{app:records,app:notebooks}.
\end{remark}

\clearpage
\section*{Conclusions and space--time convergence}
\addcontentsline{toc}{section}{Conclusions and space--time convergence}
\label{sec:merged-conclusion}
Part~\ref{part:convergence} proves time convergence for signed geodesics
without a density lower bound (\cref{geo:thm:time-discrete-convergence}),
joint space--time convergence of the signed Fourier-grid geodesic
problem (\cref{geo:thm:joint-grid,geo:cor:grid-endpoints}), and the
continuous-space modified JKO limit under the positivity floor
\eqref{jko:eq:floor} (\cref{jko:thm:main}).
Part~\ref{part:numerical} bounds the optimization error of selected grid computations by interval primal--dual bounds (\cref{tab:ipmvalidation,tab:baseline-validation,tab:common-entropic,tab:flowcontrast}); the other reported gaps are floating-point estimates.
In Part~\ref{part:time-study}, on fixed spatial grids, the action change between $N$ and $2N$ time cells decreases at observed orders between 1.82 and 2.42 (\cref{tab:geo-time-convergence}), and the final-time difference of the modified JKO scheme from the ODE reference $u_M^{\rm ref}$ of \eqref{eq:merged-ode} decreases at observed orders between 0.92 and 1.01 (\cref{tab:merged-convergence}).

For geodesics, the joint theorem allows independent spatial and temporal
refinement. Atomic reconstructions give the lower bound, while positive
Fourier averaging and time averaging give recovery without an interior
density floor. The endpoint prescriptions cover the projections with
$P_M$ of arbitrary marginals and weakly consistent grid endpoints with a common positive
lower bound. The theorem gives convergence of the minimum actions and
uniform weak convergence in time of subsequences of minimizing paths
to geodesics, but no finite-grid discretization-error bound.

The joint space--time limit of the modified JKO scheme remains open.
The following issues require additional arguments.
\begin{enumerate}[leftmargin=*,label=\arabic*.]
\item \textbf{Compactness and positivity for the modified JKO scheme.}
The proof of \cref{jko:thm:main} uses the positivity floor \eqref{jko:eq:floor} to turn weighted action estimates into spatial regularity of the chemical potential $\xi$ and then strong compactness of the densities.
For a joint limit, the estimates and positive lower bound must be
independent of both $\tau$ and $M$. Strict positivity of each
finite-grid minimizer supplies no such uniform bound. Removing the
floor would require estimates that remain valid at vacuum.

\item \textbf{Identification of the limiting equation.}
The grid proof must identify the limit of the nodal products and
logarithms in the discrete optimality equations. Positive
piecewise-constant reconstructions and discrete translation estimates
must justify convergence of the nonlinear term
$\rho\,\partial_{xx}\log\rho$ in the continuum equation.
The fixed-grid expansion of \cref{prop:step-consistency} has no
remainder bound uniform in $M$; by itself it proves no joint limit.

\item \textbf{Optimization errors along refinement.}
For geodesics, feasible grid paths whose certified action gaps vanish
are covered by the joint theorem. Fixed stopping tolerances or the
selected finite-grid certificates do not alone establish this
asymptotic condition. For the modified JKO scheme, a convergence result
for computed trajectories must also control accumulated step errors;
a small gap for each supplied input density is not a bound on the
error of the complete trajectory.
\end{enumerate}

Pointwise questions still require stronger local estimates. Weak
convergence of density and momentum does not control their quotient
near vacuum and therefore does not establish a zero density or distinct
one-sided diffusivity limits at a sign transition.

\paragraph{Two discretizations of the martingale geodesic problem.}
The midpoint scheme with nonnegative momentum and the $M_4$ coupling scheme
have different variables and finite-step objectives. In dimension one,
both converge to the same martingale geodesic cost, with the
normalizations in \cref{mart:cor:common-limit}. The Eulerian recovery
sequence preserves the sign and the endpoints by time averaging.
In one dimension, momentum
variables avoid solving a succession of optimal coupling problems;
superposition supplies a continuous martingale representation of a
limiting Eulerian minimizer afterward. The computations of
Part~\ref{part:numerical} with nonnegative momentum concern the
periodic problem, which is not a martingale problem
(\cref{rem:martingale}).
Numerical methods for the momentum formulation in higher dimensions
are future work.
For gradient flows, \cref{mart:thm:positive-potential} proves
convergence of the modified JKO scheme with nonnegative momentum for energies given
by one-dimensional potentials with bounded second and third derivatives,
without a density floor. Its limiting marginals agree with those of the $M_4$
scheme under the latter's uniform-concavity assumptions and matched
normalization. For nonlinear energies \eqref{jko:eq:energy} satisfying
\textup{(E1)--(E2)}, the constrained momentum scheme
\eqref{jko:eq:positive-scheme} on $\Torus$ has the existence and weak
compactness of \cref{jko:prop:positive-compactness}; on $\R$, and for
the identification of the limiting evolution equation, an additional
argument is required.

\clearpage\part*{Appendices}\addcontentsline{toc}{part}{Appendices}
\appendix
\renewcommand{\theHsection}{appendix.\Alph{section}}

\section{Existence and positivity}\label{app:existence}

On a fixed grid, mass conservation and an action bound make minimizing sequences bounded.
Lower semicontinuity then gives existence.
For the entropic problem, a feasible perturbation excludes zero entries, and strict convexity gives uniqueness.

\subsection{Closed action and momentum bound}
The identity
\[
 a(q,m)=\sup_{z\in\R}\{zm-\tfrac12qz^2\}
\]
holds for each $(q,m)\in\R^2$.
For $q>0$, the maximizer is $z=m/q$.
For $q=0$, the supremum is zero if $m=0$.
It is infinite if $m\ne0$.
For $q<0$, the quadratic expression is unbounded above.
As a supremum of continuous affine functions of $(q,m)$, the action is convex and lower semicontinuous.

Suppose $q\in\mathcal P_M$ and $\A(q,m)<\infty$.
At each zero entry of $q$, finite action gives $m_i=0$.
The Cauchy--Schwarz inequality on the positive entries gives
\[
 \begin{aligned}
 \normone m^2
 &=\left(\dx\sum_{q_i>0}\sqrt{q_i}\frac{|m_i|}{\sqrt{q_i}}\right)^2\\
 &\le\left(\dx\sum_{q_i>0}q_i\right)
       \left(\dx\sum_{q_i>0}\frac{m_i^2}{q_i}\right)
 =2\A(q,m).
 \end{aligned}
\]
The last equality uses $\dx\sum_{q_i>0}q_i=\avg q=1$ for $q\in\mathcal P_M$, and $\dx\sum_{q_i>0}m_i^2/q_i=2\A(q,m)$ since $m_i=0$ where $q_i=0$.
On a fixed grid, $\dx|m_i|\le\normone m$ bounds every entry of the momentum.
The constant in this entrywise estimate depends on the grid size; only fixed-grid compactness is needed here.

\subsection{Proof for discrete geodesics}
Set
\[
 \nu^j=(1-jh)\nu^0+jh\nu^N,\qquad
 m^j=h\C(\nu^N-\nu^0).
\]
The endpoint difference has mean zero.
Thus $Lm^j=h(\nu^N-\nu^0)=\nu^{j+1}-\nu^j$.
Since $\nu^0,\nu^N\in\mathcal P_M^+$, each $\nu^j$ is positive, hence every midpoint weight $q^j=(\nu^j+\nu^{j+1})/2$ of this path is positive.
The action of this path is therefore finite and is an upper bound for $E_{M,N}$ in \eqref{eq:geodesic}.

Choose a minimizing sequence whose objective values are bounded by the action of the feasible path above plus one.
Mass conservation places all its density vectors and midpoint weights in the compact set $\mathcal P_M$.
The estimate $\normone m^2\le2\A(q,m)$ bounds each momentum block.
A subsequence thus converges in the full finite-dimensional array space.
The limit satisfies the affine equations and nonnegativity.
Lower semicontinuity of the action shows that the limit attains the infimum.

Now fix a feasible density path that admits finite action.
At zero weights, each admissible finite-action momentum is zero.
On the remaining entries, the momentum objective is a positive definite quadratic form.
Its restriction to the affine conservation equations is strictly convex.
Coercivity gives existence, and strict convexity gives uniqueness of the momentum of least action.
This argument does not give uniqueness of the density path.

\subsection{Proof for an entropic step}\label{app:entropic-existence}
The pair $(\nu^{\rm old},0)$ is feasible and has finite objective.
The entropy density satisfies $e(r)\ge0$ for $r\ge0$.
Thus a bounded objective gives a bounded action.
Let $(\nu_{(n)},m_{(n)})$ be a minimizing sequence with finite objective.
Since $e(r)=+\infty$ for $r<0$, finiteness of $\E(\nu_{(n)})$ gives $\nu_{(n)}\ge0$, and $\avg{Lm_{(n)}}=0$ gives $\avg{\nu_{(n)}}=\avg{\nu^{\rm old}}=1$, so $\nu_{(n)}\in\mathcal P_M$.
The estimate $\normone{m_{(n)}}^2\le2\A(\nu_{(n)},m_{(n)})$ of \cref{app:existence} bounds $m_{(n)}$, and compactness of $\mathcal P_M$ gives a convergent subsequence.
Lower semicontinuity gives a minimizer $(\nu^*,m^*)$.

Suppose its zero set $Z_0=\{i:\nu_i^*=0\}$ is nonempty.
For $0<\epsilon<1$, form the feasible pair
\[
 \nu_\epsilon=(1-\epsilon)\nu^*+\epsilon\nu^{\rm old},\qquad
 m_\epsilon=(1-\epsilon)m^*.
\]
Convexity of $\A$ and $\A(\nu^{\rm old},0)=0$ give $\A(\nu_\epsilon,m_\epsilon)\le(1-\epsilon)\A(\nu^*,m^*)\le\A(\nu^*,m^*)$.
At a zero entry,
\[
 e(\epsilon\nu_i^{\rm old})-e(0)
 =\epsilon\nu_i^{\rm old}\log\epsilon
   +\epsilon\nu_i^{\rm old}(\log\nu_i^{\rm old}-1).
\]
At each positive entry of $\nu^*$, the entropy change is $O(\epsilon)$.
The grid is finite, so
\[
 \E(\nu_\epsilon)-\E(\nu^*)
 =c_0\epsilon\log\epsilon+O(\epsilon),\qquad
 c_0=\dx\sum_{i\in Z_0}\nu_i^{\rm old}>0.
\]
Since $\log\epsilon\to-\infty$, the negative term $c_0\epsilon\log\epsilon$ dominates all terms bounded in magnitude by a constant times $\epsilon$.
By the convexity bound on $\A(\nu_\epsilon,m_\epsilon)$, the action term does not increase, so the change of the objective is at most $\E(\nu_\epsilon)-\E(\nu^*)$.
For sufficiently small $\epsilon$, the total objective is therefore smaller than at $(\nu^*,m^*)$, which contradicts minimality.
Hence $\nu^*>0$.

If two minimizers have different densities, strict convexity of the entropy gives a smaller objective at their midpoint.
Thus the successor density $\nu^*$ is unique.
At this positive density, the action is strictly convex in momentum.
Hence the minimizing momentum is also unique.
Comparison with the feasible pair $(\nu^{\rm old},0)$ gives the entropy inequality in \cref{prop:positive}.

\subsection{Proof of fixed-grid time consistency}\label{app:step-consistency}
Write $w=m/\tau$.
For a positive stationary point, the momentum and density equations of \eqref{eq:step} give
\[
 w=-\nu L\chi,\qquad
 \chi=\log\nu-\frac\tau2(w/\nu)^2,
 \qquad \nu=\nu^{\rm old}+\tau Lw.
\]
A constant added to $\chi$ has no effect on the first equation.
Eliminate $\nu$ and $\chi$ by defining, where $\nu^{\rm old}+\tau Lw>0$,
\[
 \mathcal F(w,\tau)=w+(\nu^{\rm old}+\tau Lw)L\!\left[
 \log(\nu^{\rm old}+\tau Lw)-\frac\tau2\left(\frac{w}{\nu^{\rm old}+\tau Lw}\right)^2\right].
\]
At $\tau=0$, $\mathcal F(w,0)=w+\nu^{\rm old}L\log\nu^{\rm old}$ and $D_w\mathcal F(w,0)=I$.
The implicit function theorem therefore gives a smooth branch
$w(\tau)=-\nu^{\rm old}L\log\nu^{\rm old}+O(\tau)$ with $\nu^{\rm old}+\tau Lw(\tau)>0$ for sufficiently small $\tau$.
The pair $(\nu^*_\tau,\tau w(\tau))$ with $\nu^*_\tau=\nu^{\rm old}+\tau Lw(\tau)$ satisfies the stationarity conditions, which are sufficient for optimality because \eqref{eq:step} is convex.
Uniqueness in \cref{prop:positive} identifies this pair with the minimizer.
Substitution of $w(\tau)$ into $\nu^*_\tau=\nu^{\rm old}+\tau Lw(\tau)$ gives the expansion of \cref{prop:step-consistency}.
The neighborhood and remainder depend on the fixed operator and on $\min_i\nu^{\rm old}_i>0$; this argument supplies no uniform estimate as either changes.

\section{Proximal maps and DR convergence}\label{app:drmaps}

\subsection{Local proximal maps}\label{app:localprox}
The action $a$ is the perspective of the scalar quadratic; its proximal map is a special case of the calculus in \cite{CombettesMuller18}.
Fix $b>0$ and an input $(q_{\rm in},v_{\rm in})\in\R^2$.
The local proximal problem is
\[
 \min_{\substack{q\ge0\\v\in\R}}
 \left\{b\,a(q,v)+\tfrac12(q-q_{\rm in})^2+\tfrac12(v-v_{\rm in})^2\right\}.
\]
Its minimizing pair $(q_*,v_*)$ is
\begin{equation}\label{eq:prox}
 \begin{cases}
 (q_*,v_*)=(0,0),&2bq_{\rm in}+v_{\rm in}^2\le0,\\[2pt]
 (q_*-q_{\rm in})(q_*+b)^2=\dfrac b2v_{\rm in}^2,
 \quad v_*=\dfrac{q_*}{q_*+b}v_{\rm in},
 &2bq_{\rm in}+v_{\rm in}^2>0.
 \end{cases}
\end{equation}

We first fix $q>0$ and minimize this objective with respect to $v$.
Since $a(q,v)=v^2/(2q)$, the objective is a strictly convex quadratic function of $v$.
Its derivative with respect to $v$ vanishes exactly when
\[
 \frac{b}{q}v+v-v_{\rm in}=0,
 \qquad\text{that is,}\qquad
 v=\frac{q}{q+b}v_{\rm in}.
\]
Substitution of this value of $v$ gives the reduced objective
\[
 f(q)=\tfrac12(q-q_{\rm in})^2+
             \frac{b\,v_{\rm in}^2}{2(q+b)}.
\]
At $q=0$, finite action requires $v=0$, and $a(0,0)=0$.
The value of the proximal objective at this pair is
\[
 b\,a(0,0)+\tfrac12(0-q_{\rm in})^2+\tfrac12(0-v_{\rm in})^2
 =\tfrac12q_{\rm in}^2+\tfrac12v_{\rm in}^2.
\]
This value equals $f(0)=\tfrac12q_{\rm in}^2+\dfrac{b\,v_{\rm in}^2}{2b}$.
Hence $f$ is the reduced objective on the closed half-line $q\ge0$.
For $q\ge0$,
\[
 f'(q)=q-q_{\rm in}-\frac{b\,v_{\rm in}^2}{2(q+b)^2},\qquad
 f''(q)=1+\frac{b\,v_{\rm in}^2}{(q+b)^3}>0.
\]
Since $f''>0$, $f'$ is increasing.
The minimizer over $q\ge0$ is therefore $q_*=0$ exactly when the right derivative $f'(0)=-q_{\rm in}-v_{\rm in}^2/(2b)$ is nonnegative, that is, when $2bq_{\rm in}+v_{\rm in}^2\le0$.
Otherwise $f'(0)<0$ and $f'(q)\to+\infty$ as $q\to\infty$, so $f'$ has a unique zero in $(0,\infty)$, which is the minimizer; the equation $f'(q_*)=0$ is the second case of \eqref{eq:prox}.

\paragraph{Entropic proximal map.}
For the map \eqref{eq:entropic-prox} at node $i$, with density input $q_{\rm in}=\nu_{{\rm in},i}$ and momentum input $v_{\rm in}=v_{{\rm in},i}$, elimination of $v$ gives
\[
 f_{\rm ent}(q)=\tfrac12(q-q_{\rm in})^2
       +\frac{b\,v_{\rm in}^2}{2(q+b)}+\omega e(q),\qquad q\ge0.
\]
For $q>0$, differentiation with respect to $q$ gives
\[
 g(q):=f_{\rm ent}'(q)=q-q_{\rm in}+\omega\log q
             -\frac{b\,v_{\rm in}^2}{2(q+b)^2}.
\]
The density output is the root $g(q)=0$, and the scaled momentum is $v=q\,v_{\rm in}/(q+b)$.
\Cref{alg:jko-prox} solves $g(q)=0$ for the unknown $\log q$ with a bracketed Newton iteration.
The derivative of $g$ satisfies
\[
 g'(q)=1+\frac{\omega}{q}+\frac{b\,v_{\rm in}^2}{(q+b)^3}>0.
\]
Moreover, $g(q)\to-\infty$ as $q\downarrow0$ and $g(q)\to+\infty$ as $q\to\infty$.
Hence there is a unique positive root.
Let $q_*$ be the zero of $g$.
Since $f_{\rm ent}'=g<0$ on $(0,q_*)$ and $f_{\rm ent}$ is continuous at $0$, $f_{\rm ent}(q)<f_{\rm ent}(0)$ for $0<q\le q_*$.
Hence zero is not a minimizer, and the minimizer is $q_*$ with $v=q_*v_{\rm in}/(q_*+b)$.
This proves positivity of every entropic proximal output, without a uniform lower bound.

\paragraph{Nonnegative momentum.}\label{app:positive-prox}
For $v_{\rm in}<0$ and any fixed $q\ge0$, the action and squared-distance terms are minimized over $v\ge0$ at $v=0$.
Minimizing the remaining density term gives $q=\max(q_{\rm in},0)$, which is also $\prox_{ba}(q_{\rm in},0)$.
For $v_{\rm in}\ge0$, replacing a negative candidate momentum by its absolute value preserves the action and cannot increase its squared distance to $v_{\rm in}$.
The unconstrained proximal minimizer is therefore nonnegative.
These two cases prove \eqref{eq:positive-prox}.

\subsection{Proximal residual and accumulated error}\label{app:prox-error}
Let $b>0$, $\omega\ge0$, and use the positive branch of the local proximal problem with entropy coefficient $\omega$.
For $q>0$, set
\[
 v=\frac{q\,v_{\rm in}}{q+b},\qquad
 g_\omega(q)=q-q_{\rm in}+\omega\log q
                -\frac{b\,v_{\rm in}^2}{2(q+b)^2}.
\]
For $\omega=0$, omit the logarithmic term.
For $\omega>0$, $g_\omega$ is the function $g$ of \cref{app:localprox}.
We call $g_\omega(q)$ the scalar residual of the candidate $q$; it gives a stopping test for each local proximal solve, independent of the DR update norm $\|z^{\ell+1}-z^\ell\|_2$.
The zero branch of \eqref{eq:prox} is evaluated separately.

\begin{proposition}[Control of the proximal error]\label{prop:prox-error}
Let $(q_*,v_*)$ be the exact local proximal output.
For the positive candidate above,
\[
 \|(q,v)-(q_*,v_*)\|_2\le |g_\omega(q)|.
\]
For a DR iteration with exact affine projection, let $\varepsilon_\ell$ bound the Euclidean norm of all local scalar residuals at iteration $\ell$; evaluate the nodal-density update $\nu_{\rm in}\mapsto\max(\nu_{\rm in},0)$ and the zero branch of \eqref{eq:prox} exactly.
If $\sum_{\ell=0}^\infty\varepsilon_\ell<\infty$, the auxiliary iterates $z^\ell$ converge to a DR fixed point $z_*$, and the proximal outputs $x^\ell$ and the projected outputs $y^\ell=P_{\mathcal S_\beta}(2x^\ell-z^\ell)$ converge to $x_*=\prox_{\gamma J}(z_*)$, a minimizer of $J$ over $\mathcal S_\beta$ under the assumptions of \cref{prop:drglobal}.
\end{proposition}
\begin{proof}
The local objective is 1-strongly convex because of its squared-distance term.
At the candidate, its first-order residual is $(g_\omega(q),0)$.
Strong monotonicity of its subdifferential, followed by Cauchy--Schwarz, gives the pair-error estimate.
The block estimates imply $\|\delta x^\ell\|_2\le\varepsilon_\ell$, where $\delta x^\ell$ is the error in the full proximal output $x^\ell$.

Let $\mathcal D$ be the exact DR update and $P$ the orthogonal projector onto the direction space $\mathcal S_\beta^0$ of $\mathcal S_\beta$.
The approximate update is
\[
 z^{\ell+1}=\mathcal D z^\ell+(2P-I)\delta x^\ell,
 \qquad \|(2P-I)\delta x^\ell\|_2=\|\delta x^\ell\|_2.
\]
For any fixed point $z_*$, nonexpansiveness and summability first imply boundedness of $z^\ell$.
Firm nonexpansiveness then gives
\[
 \|z^{\ell+1}-z_*\|_2^2
 \le \|z^\ell-z_*\|_2^2
     -\|z^\ell-\mathcal D z^\ell\|_2^2+c\varepsilon_\ell+\varepsilon_\ell^2
\]
for a finite constant $c$.
The distances to each fixed point converge, and the squared update residuals are summable.
In finite dimension, a cluster point exists; continuity of $\mathcal D$ makes it a fixed point.
Convergence of the distances to this cluster point gives convergence of the full sequence.
Continuity of the exact proximal map and $\delta x^\ell\to0$ give convergence of the proximal outputs.
At the fixed point, $P_{\mathcal S_\beta}(2x_*-z_*)=x_*$, so continuity of $P_{\mathcal S_\beta}$ gives $y^\ell\to x_*$.
\end{proof}

For $n_{\rm b}$ scalar blocks and constants $\varepsilon_0>0$, $\alpha>0$, the tolerance $|g_{\omega,i}|\le\varepsilon_0(\ell+1)^{-1-\alpha}/\sqrt{n_{\rm b}}$ gives $\varepsilon_\ell=\varepsilon_0(\ell+1)^{-1-\alpha}$ and $\sum_\ell\varepsilon_\ell<\infty$, a sufficient schedule in exact arithmetic.
The residual-checked solvers (\cref{alg:geo-prox-checked,alg:jko-prox-checked}) instead accept a computed $q$ only if
\[
 |g_{\omega,i}|\le 2\cdot10^{-14}
       \bigl(1+\max\{|q_{{\rm in},i}|,|v_{{\rm in},i}|\}\bigr)
\]
for every scalar block; otherwise the solve stops with a failure status.
The computed residual is not an interval bound on rounding error.
Fixed tolerances and rounded affine projections therefore do not establish convergence to an exact fixed point as the iteration count tends to infinity.
The local rate statements in \cref{prop:lfp} remain statements about the exact iteration.

\subsection{Subdifferential at vacuum}
A pair $(\alpha_q,\alpha_m)$ belongs to $\partial a(0,0)$ if
\[
 \alpha_q q+\alpha_m m\le a(q,m)\quad\text{for all }(q,m).
\]

For $q>0$, minimize $m^2/(2q)-\alpha_m m-\alpha_q q$ over $m$.
The minimum is $-q(\alpha_q+\alpha_m^2/2)$.
It is nonnegative for each $q>0$ exactly when $\alpha_q+\alpha_m^2/2\le0$.
The inequality also holds at $(0,0)$.
At the remaining points, where $q<0$ or $q=0\ne m$, $a(q,m)=+\infty$, so the subgradient inequality holds.
This proves \eqref{eq:vacsubdiff}.

\subsection{Scaled splitting and global convergence}\label{app:drglobal}

For the geodesic variables, the weighted maps are
\[
 \prox_{\gamma F}^\beta(U_{\rm in})=\argmin_U\{\gamma F(U)+\tfrac12\|U-U_{\rm in}\|_\beta^2\},
 \qquad P_\K^\beta(U_{\rm in})=\argmin_{U\in\K}\tfrac12\|U-U_{\rm in}\|_\beta^2.
\]
Recall $J=F\circ Q_\beta^{-1}$ and $\mathcal S_\beta=Q_\beta\K$.
Since $\|Q_\beta U\|_2=\|U\|_\beta$, a change of variables in the proximal and projection problems gives
\[
 \prox_{\gamma J}=Q_\beta\prox_{\gamma F}^\beta Q_\beta^{-1},
 \qquad
 P_{\mathcal S_\beta}=Q_\beta P_\K^\beta Q_\beta^{-1}.
\]
The same identities hold for the entropic objective and its affine set.
Hence DR \eqref{eq:dr} in scaled coordinates is the image under $Q_\beta$ of DR for $(F,\K)$ in the norm $\|\cdot\|_\beta$, and similarly for the entropic step.

Both $J$ and $\iota_{\mathcal S_\beta}$ are proper, lower semicontinuous, and convex.
The positive geodesic interpolation from \cref{app:existence} is feasible and belongs to the interior of the objective domain.
For the entropic step, $(\nu^{\rm old},0)$ has the same properties.
Hence $\mathcal S_\beta\cap\operatorname{ri}(\operatorname{dom}J)\ne\varnothing$ in both cases.
Here $\operatorname{dom}J=\{x:J(x)<+\infty\}$, and $\operatorname{ri}$ denotes interior within the affine hull of a convex set.
Minimizers exist by \cref{app:existence}.

For a convex function $f$ and a point $x$ with $f(x)<+\infty$, the subdifferential is $\partial f(x)=\{g:f(y)\ge f(x)+\langle g,y-x\rangle\text{ for all }y\}$, with the Euclidean inner product in these coordinates.
Let $x_*$ be a minimizer of $J+\iota_{\mathcal S_\beta}$, and recall that $\mathcal S_\beta^0=\mathcal S_\beta-x_*$.
Since $\mathcal S_\beta$ is affine, $\partial\iota_{\mathcal S_\beta}(x_*)=(\mathcal S_\beta^0)^\perp$.
The intersection $\mathcal S_\beta\cap\operatorname{ri}(\operatorname{dom}J)\ne\varnothing$ permits the subdifferential sum rule, so optimality gives
\[
 0\in\partial(J+\iota_{\mathcal S_\beta})(x_*)
   =\partial J(x_*)+(\mathcal S_\beta^0)^\perp.
\]
By this inclusion, some element of $\partial J(x_*)$ has its opposite in $(\mathcal S_\beta^0)^\perp$; let $z_*$ be $x_*$ plus $\gamma$ times this element. Since $(\mathcal S_\beta^0)^\perp$ is a linear subspace and $\gamma>0$, $z_*$ satisfies $\gamma^{-1}(z_*-x_*)\in\partial J(x_*)$ and $x_*-z_*\in(\mathcal S_\beta^0)^\perp$.

The proximal optimality condition from \eqref{eq:prox-projection} holds at $x_*$ because
\[
 0=(z_*-x_*)+(x_*-z_*)\in\gamma\partial J(x_*)+x_*-z_*.
\]
The proximal objective is strongly convex, so $\prox_{\gamma J}(z_*)=x_*$.

Since $x_*\in\mathcal S_\beta$, the orthogonal projection of a vector $u$ onto $\mathcal S_\beta$ equals $x_*$ exactly when $u-x_*\in(\mathcal S_\beta^0)^\perp$.
For the reflected proximal point,
\[
 (2x_*-z_*)-x_*=x_*-z_*\in(\mathcal S_\beta^0)^\perp,
 \qquad P_{\mathcal S_\beta}(2x_*-z_*)=x_*.
\]
The DR update at $z_*$ is therefore
\[
 z_*+P_{\mathcal S_\beta}\bigl(2\prox_{\gamma J}(z_*)-z_*\bigr)
       -\prox_{\gamma J}(z_*)
 =z_*+x_*-x_*=z_*.
\]
Since a fixed point exists, the convergence theorem for DR of Lions and Mercier \cite{LM79,BC11} applies: for every $z^0$, the auxiliary iterates $z^\ell$ converge to a fixed point, which in general depends on $z^0$ and need not be the fixed point constructed above. With $z_*$ now denoting this limit and $x_*=\prox_{\gamma J}(z_*)$, the same theorem gives $x^\ell\to x_*$ and $y^\ell\to x_*$, and $x_*$ is a minimizer of $J$ over $\mathcal S_\beta$.
This proves \cref{prop:drglobal}.

\subsection{Local linear convergence of DR}\label{app:lfp}\label{app:drlinearization}

Let $z_*$ be the limit given by \cref{prop:drglobal}, let $x_*=\prox_{\gamma J}(z_*)$, and set $\delta z^\ell=z^\ell-z_*$.
This subsection applies the local convergence theory of Liang, Fadili and Peyr\'e \cite{LFP17}, abbreviated LFP, to the DR iteration \eqref{eq:dr}.
Let $R$ be the Jacobian of $\prox_{\gamma J}$ at $z_*$ defined in \cref{sec:dr-local}, when it exists.
For geodesics, $R$ is block diagonal because the proximal map acts separately on nodal-density entries and action pairs.
At a positive limiting nodal-density input, the scalar map $u\mapsto\max(u,0)$ is locally the identity, so its derivative is one.

Each action block depends on $b$ defined in \eqref{eq:b}.
At the limiting output of an action pair with $q>0$, put $s=v/q$.
The Hessian of $a$ with respect to the output coordinates $(q,v)$ and its nonzero eigenvalue are
\begin{equation}\label{eq:rankone}
 H:=\nabla^2_{(q,v)}a(q,v)
 =\frac1q\begin{pmatrix}s^2&-s\\-s&1\end{pmatrix},\qquad
 \lambda=\frac{1+s^2}q.
\end{equation}
The corresponding $2\times2$ block of $R$, written in $(q,v)$ coordinates, is $(I+bH)^{-1}$.
It leaves the direction $(1,s)$ unchanged and multiplies the component along $(s,-1)$ by $(1+b\lambda)^{-1}$.

Recall the direction space of the affine constraints from \cref{sec:dr-local},
$\mathcal S_\beta^0=\mathcal S_\beta-x_*=\{y-x_*:y\in\mathcal S_\beta\}$.
Its vectors preserve the affine equations; nonnegativity is imposed separately by $J$.
For geodesics, $\delta x=(\delta\nu,\delta v,\delta q)$ belongs to $\mathcal S_\beta^0$ exactly when
\[
 L\delta v^j=\sqrt\beta\,(\delta\nu^{j+1}-\delta\nu^j),
 \qquad 2\delta q^j=\delta\nu^j+\delta\nu^{j+1},\qquad 0\le j<N,
\]
with $\delta\nu^0=\delta\nu^N=0$ because the endpoints are fixed.
For an entropic step, the equation is $L\delta v=\sqrt\beta\,\delta\nu$.
If $P$ is the orthogonal projector onto $\mathcal S_\beta^0$, then
$P_{\mathcal S_\beta}(u+\delta x)-P_{\mathcal S_\beta}(u)=P\delta x$ for every input $u$ and every perturbation $\delta x$.
This identity is exact because the constraint set is affine.

The derivative, at input $z_*$, of the auxiliary update $z\mapsto z+P_{\mathcal S_\beta}(2\prox_{\gamma J}(z)-z)-\prox_{\gamma J}(z)$ is $T$ in \eqref{eq:T}; this is proved in the paragraph on the linearization of the DR update below.
Thus
\[
 \delta z^{\ell+1}=T\,\delta z^\ell+o(\|\delta z^\ell\|_2).
\]
To cover the case where $T$ has eigenvalue one, let $T_\infty$ be the orthogonal projector onto $\ker(I-T)=\{\delta x:T\delta x=\delta x\}$.
\Cref{prop:lfp} gives a geometric rate only when $T_\infty=0$, that is, when $T$ has no fixed direction.
For $T_\infty\ne0$, the linearization gives no proved rate here, because the remainder $o(\|\delta z^\ell\|_2)$ need not be small compared with $\|(I-T_\infty)\delta z^\ell\|_2$.
In both cases, $\delta z^\ell\to0$ by \cref{prop:drglobal}.

\Needspace{16\baselineskip}
\begin{proposition}[Local convergence estimate for both problems]\label{prop:lfp}
Assume the hypotheses of \cref{prop:drglobal}: $M\ge3$, $\beta,\gamma>0$, and either a geodesic with $N\ge2$ and endpoints $\nu^0,\nu^N\in\mathcal P_M^+$, or an entropic step with $\tau>0$ and $\nu^{\rm old}\in\mathcal P_M^+$.
Let $(x^\ell,y^\ell,z^\ell)$ be generated by \eqref{eq:dr} with exact proximal maps and projections.
For a geodesic, write the minimizing output as $x_*=(\nu_*,v_*,q_*)$
and the limiting proximal input as $z_*=(\nu_{\rm in},v_{\rm in},q_{\rm in})$.
Assume
\begin{equation}\label{eq:dr-strict-input}
 \nu_{{\rm in},i}^j<0\quad\text{if }\nu_{*,i}^j=0,
 \qquad
 2b\,q_{{\rm in},i}^j+(v_{{\rm in},i}^j)^2<0
 \quad\text{if }q_{*,i}^j=0.
\end{equation}
No such condition is needed for an entropic step.
Then $R$ exists.
For a geodesic, write the proximal output as $x^\ell=(\nu^\ell,v^\ell,q^\ell)$.
There exists an iteration $\ell_0$ such that, for every $\ell\ge\ell_0$,
\[
 (\nu^\ell)_i^j=0\ \Longleftrightarrow\ \nu_{*,i}^j=0,
 \qquad
 (q^\ell)_i^j=0\ \Longleftrightarrow\ q_{*,i}^j=0,
\]
for all entries of the respective arrays.
Moreover, $\lim_{n\to\infty}T^n=T_\infty$ and $\spr(T-T_\infty)<1$.
If $T_\infty=0$, then for every $\spr(T)<c<1$ there exist a constant $C_{\rm loc}>0$ and an iteration $\ell_{\rm loc}$ such that, for all $\ell\ge\ell_{\rm loc}$,
\[
 \begin{aligned}
 \|z^\ell-z_*\|_2&\le C_{\rm loc}c^\ell,\\
 \|x^\ell-x_*\|_2&\le C_{\rm loc}c^\ell,\qquad
 \|y^\ell-x_*\|_2\le C_{\rm loc}c^\ell.
 \end{aligned}
\]
For every entropic step, $T_\infty=0$.
\end{proposition}

Consider a positive geodesic, that is, a geodesic minimizer $x_*$ whose densities and midpoint weights are all positive.
Under the hypotheses of \cref{prop:lfp}, $T_\infty=0$ if and only if
\begin{equation}\label{eq:constrained-hessian}
 \langle\delta x,\nabla^2J(x_*)\delta x\rangle>0
 \quad\text{for every }\delta x\in\mathcal S_\beta^0\setminus\{0\},
\end{equation}
because $\ker(I-T)=\mathcal S_\beta^0\cap\ker\mathcal H$ (paragraph on the directions unchanged by the linearized DR map below) and $\mathcal H=\gamma\nabla^2J(x_*)$ is positive semidefinite.
Under \eqref{eq:constrained-hessian}, \cref{prop:lfp} gives, for every $\spr(T)<c<1$, the bounds $\|z^\ell-z_*\|_2,\|x^\ell-x_*\|_2,\|y^\ell-x_*\|_2\le C_{\rm loc}c^\ell$ for $\ell\ge\ell_{\rm loc}$.
For a positive geodesic, \eqref{eq:constrained-hessian} holds if and only if \eqref{eq:geodesic} has a unique minimizing density path.
Indeed, $\nabla^2J(x_*)$ is positive semidefinite, so \eqref{eq:constrained-hessian} fails exactly when some $\delta x\in\mathcal S_\beta^0\setminus\{0\}$ satisfies $\nabla^2J(x_*)\delta x=0$.
In each action pair this means $(\delta q,\delta v)=\delta q\,(1,s)$ with $s=v_*/q_*$, so that $(q_*,v_*)+t(\delta q,\delta v)=(1+t\,\delta q/q_*)(q_*,v_*)$.
For $|t|$ small, $x_*+t\delta x$ therefore lies in $\mathcal S_\beta$, has positive densities and midpoint weights, and, by the homogeneity $a(\lambda q,\lambda v)=\lambda\,a(q,v)$, satisfies $J(x_*+t\delta x)=J(x_*)+t\,c_J$ for a constant $c_J$.
Since $x_*$ is a minimizer, $c_J=0$, and every such $x_*+t\delta x$ is a minimizer.
Its density path differs from that of $x_*$ for $t\ne0$: if $\delta\nu=0$, the midpoint equations and the fixed endpoints give $\delta q=0$, hence $\delta v=0$ and $\delta x=0$.
Conversely, under \eqref{eq:constrained-hessian}, $J$ is twice differentiable near $x_*$, so $x_*$ is a strict local minimizer of $J$ on $\mathcal S_\beta$; by convexity of $J$, it is the unique minimizer.
Uniqueness of the minimizing density path is not proved (remark after \cref{prop:geodesic-existence}), so it is open whether positivity of all densities and midpoint weights implies \eqref{eq:constrained-hessian}.
For positive geodesics, \cref{app:geocontrast} shows that \eqref{eq:constrained-hessian} is equivalent to $\mu_s>0$.
For an entropic step, the objective Hessian is positive definite and always gives $T_\infty=0$.

Two hypotheses of \cref{prop:lfp} are not verified for the computations.
First, \eqref{eq:dr-strict-input} has not been checked for the geodesic runs of \cref{tab:zeros}, whose proximal outputs contain zero entries; a zero output does not imply \eqref{eq:dr-strict-input}.
Second, \cref{prop:lfp} assumes exact proximal maps and fixed parameters, whereas the proximal solvers (\cref{alg:geo-prox,alg:geo-prox-checked,alg:jko-prox,alg:jko-prox-checked}) stop the Newton iterations for the scalar proximal equations after finitely many steps; \cref{prop:lfp} does not bound the effect of the resulting errors in the proximal outputs on the DR iterates.

\begin{proof}[Proof of \cref{prop:lfp}]
The proof establishes the optimality relations at the limit, linearizes the DR update, identifies the zero entries of geodesic outputs, treats the entropic step, applies \cite{LFP17}, and determines the fixed directions of $T$.

\paragraph{Optimality at the limit.}
The fixed-point relations in \cref{app:drglobal} imply
\[
 z_*-x_*\in\gamma\partial J(x_*),\qquad
 x_*-z_*\perp\mathcal S_\beta^0,\qquad \mathcal S_\beta^0=\mathcal S_\beta-x_*.
\]
Hence $0\in\partial J(x_*)+(\mathcal S_\beta^0)^\perp$, as in \cref{app:drglobal}.

\paragraph{Linearization of the DR update.}
For $u_0\in\mathcal S_\beta$, write $P_{\mathcal S_\beta}(z)=u_0+P(z-u_0)$.
If $\prox_{\gamma J}(z_*+\delta z)=x_*+R\,\delta z+o(\|\delta z\|_2)$, the three DR steps give
\[
 x=x_*+R\,\delta z+o(\|\delta z\|_2),\qquad
 y=x_*+P(2R-I)\delta z+o(\|\delta z\|_2),
\]
and hence $z+y-x=z_*+[I-R+P(2R-I)]\delta z+o(\|\delta z\|_2)$.
This proves \eqref{eq:T}.
For $\delta x\in\mathcal S_\beta^0$, use $P\delta x=\delta x$ to obtain
\[
 T\delta x-\delta x=(2P-I)(R-I)\delta x.
\]
Since $2P-I$ is an orthogonal reflection, $\|(T-I)\delta x\|_2=\|(R-I)\delta x\|_2$ for every $\delta x\in\mathcal S_\beta^0$: the linearized DR update changes a perturbation in $\mathcal S_\beta^0$ exactly as much as the Jacobian $R$ does, and $T\delta x=\delta x$ if and only if $R\delta x=\delta x$.
In the physical variables $(\nu,m,q)$, the linearized iteration has matrix $Q_\beta^{-1}TQ_\beta$, which has the same eigenvalues as $T$; hence $\spr(T)$ does not depend on the coordinates used to measure the error.
If $J$ is smooth at $x_*$, differentiation of $z=x+\gamma\nabla J(x)$ gives
\[
 R=(I+\mathcal H)^{-1},\qquad \mathcal H=\gamma\nabla^2J(x_*).
\]
The Hessian is evaluated at the proximal output, not its input.
For a positive geodesic cell, $\gamma J$ contributes $b\,a(q,v)$.

\paragraph{Geodesic entries that vanish at the minimizer.}
Keep each zero density coordinate of $x_*$ equal to zero and each zero action pair equal to $(0,0)$.
Allow all remaining coordinates to vary, with the remaining densities and midpoint weights positive.
The resulting set $\mathfrak M$ is open in the linear space $\mathfrak E$ obtained by fixing these zero coordinates.
The restriction of $J$ to $\mathfrak M$ is smooth.
For a nodal-density coordinate, the objective term is $\iota_{[0,\infty)}$; its subdifferential at zero is $(-\infty,0]$.
At a zero action pair, \eqref{eq:vacsubdiff} gives a subdifferential whose affine hull is $\R^2$.
At a positive entry, the subdifferential contains only the ordinary gradient, which varies smoothly.
The linear space generated by differences of subgradients is therefore exactly $\mathfrak E^\perp$.
Moreover, the subdifferentials vary continuously as the positive coordinates vary within $\mathfrak M$.
By the definition of partial smoothness in \cite{LFP17}, $J$ is partly smooth at $x_*$ relative to $\mathfrak M$ if (i) $J$ is $C^2$ on $\mathfrak M$ near $x_*$, (ii) the span of the differences of subgradients of $J$ at $x_*$ equals $\mathfrak E^\perp$, and (iii) $\partial J$ is continuous along $\mathfrak M$.
The three properties above are (i)--(iii).

The LFP non-degeneracy condition is $z_*-x_*\in\operatorname{ri}\partial(\gamma J)(x_*)$.
At zero coordinates it reads
\[
 \nu_{\rm in}<0,\qquad
 (q_{\rm in},v_{\rm in})\in
 \operatorname{ri}\partial(ba)(0,0)
 =\{(\alpha_q,\alpha_m):2b\alpha_q+\alpha_m^2<0\}.
\]
These are the assumptions \eqref{eq:dr-strict-input}.
Because the limiting inequalities are strict, they continue to hold for all sufficiently close inputs.
For sufficiently close inputs, $u\mapsto\max(u,0)$ returns zero on the limiting density zero set, and \eqref{eq:prox} returns $(0,0)$ on the limiting action zero set.
The corresponding blocks of the Jacobian $R$ are therefore zero.
At positive nodal-density outputs, the scalar block of $R$ is one; at positive midpoint weights, the action block of $R$ is $(I+bH)^{-1}$.
By \cref{prop:drglobal}, $z^\ell\to z_*$ and $x^\ell\to x_*$.
Hence there is $\ell_0$ such that, for $\ell\ge\ell_0$, the entries of $z^\ell$ satisfy the strict inequalities of \eqref{eq:dr-strict-input} at every zero entry of $x_*$, so $(\nu^\ell)_i^j=0$ and $(q^\ell)_i^j=0$ there.
After enlarging $\ell_0$, $(\nu^\ell)_i^j>0$ wherever $\nu_{*,i}^j>0$ and $(q^\ell)_i^j>0$ wherever $q_{*,i}^j>0$, because $x^\ell\to x_*$.
This proves the identification statement of \cref{prop:lfp}.
\Cref{prop:lfp} does not cover limiting inputs with $\nu_{{\rm in},i}^j=0$ or $2b\,q_{{\rm in},i}^j+(v_{{\rm in},i}^j)^2=0$, where differentiability of $\prox_{\gamma J}$ is not claimed.

\paragraph{Entropic step.}
The minimizing density is positive by \cref{prop:positive}.
Thus $J$ is smooth near $x_*$, so $\partial(\gamma J)(x_*)=\{\gamma\nabla J(x_*)\}$.
A single point is its own relative interior, and $z_*-x_*\in\gamma\partial J(x_*)$ by the optimality relation at the limit, so $z_*-x_*\in\operatorname{ri}\partial(\gamma J)(x_*)$.
Write $\nu$ for one density entry of a positive proximal output and put $s=v/\nu$.
The Hessian of the action and entropy terms, with respect to $(\nu,v)$, is
\[
 H_{\rm JKO}=\frac1\nu
 \begin{pmatrix}\omega+bs^2&-bs\\-bs&b\end{pmatrix},
 \qquad \det H_{\rm JKO}=\frac{b\omega}{\nu^2}>0.
\]
Here $H_{\rm JKO}=\nabla^2\phi(\nu,v)$, with $\phi$ from \eqref{eq:entropic-prox}.
For an input perturbation $(\delta\nu_{\rm in},\delta v_{\rm in})$ and the corresponding first-order output perturbation $(\delta\nu,\delta v)$, differentiation of the proximal optimality condition gives
\[
 (I+H_{\rm JKO})
 \begin{pmatrix}\delta\nu\\\delta v\end{pmatrix}
 =\begin{pmatrix}\delta\nu_{\rm in}\\\delta v_{\rm in}\end{pmatrix}.
\]
Thus the corresponding block of the Jacobian $R$ is $(I+H_{\rm JKO})^{-1}$, with $H_{\rm JKO}$ evaluated at the output.

At the minimizer, each block $H_{\rm JKO}$ of $\mathcal H$ is positive definite because $b/\nu>0$ and $\det H_{\rm JKO}>0$.
If $0<\lambda_{\min}\le\lambda_{\max}$ are the extreme eigenvalues of $\mathcal H$, then
\[
 \|T\|_2
 \le\tfrac12(1+\|2R-I\|_2)
 =\max\left\{\frac1{1+\lambda_{\min}},\frac{\lambda_{\max}}{1+\lambda_{\max}}\right\}<1.
\]
Here we used $T=\tfrac12[I+(2P-I)(2R-I)]$ and the orthogonality of $2P-I$.
The inequality $\|T\|_2<1$ excludes a nonzero vector satisfying $T\delta x=\delta x$.
Hence $T_\infty=0$ for every entropic step.
This bound needs only that the minimizing densities are positive.
It holds for every minimizing momentum, without the substitutions that define $R_{\rm model}$ in \cref{app:entropicbalance}.

\paragraph{Application of LFP.}
The reference \cite{LFP17} treats $\min_xJ(x)+G(x)$; here $G=\iota_{\mathcal S_\beta}$, whose domain $\mathcal S_\beta$ is affine with direction space $\mathcal S_\beta^0$.
At each point of $\mathcal S_\beta$, $\partial G=(\mathcal S_\beta^0)^\perp$.
A linear space is equal to its relative interior, so the fixed-point relation $x_*-z_*\in(\mathcal S_\beta^0)^\perp$ gives the required non-degeneracy condition for $G$.
The linearizations of $\prox_{\gamma J}$ and $\prox_{\gamma G}=P_{\mathcal S_\beta}$ at the fixed point are $R$ and $P$.
Equation (5.4) of \cite{LFP17} gives $T$ in \eqref{eq:T}, and Lemma 5.4(i) of that reference gives $\lim_{n\to\infty}T^n=T_\infty$ and $\spr(T-T_\infty)<1$.
When $T_\infty=0$, Theorem 5.6(i) of that reference gives the geometric bound for $\|z^\ell-z_*\|_2$ in \cref{prop:lfp}.
The parameter-error term in Theorem 5.6(i) of \cite{LFP17} is zero here because the DR map is the same at every iteration: the scaled objective, the proximal parameter, and the update rule are fixed.
The proximal map, its reflection, and the affine projection do not increase distances.
Since $x_*=\prox_{\gamma J}(z_*)=P_{\mathcal S_\beta}(2x_*-z_*)$, this gives
\[
 \|x^\ell-x_*\|_2\le\|z^\ell-z_*\|_2,\qquad
 \|y^\ell-x_*\|_2\le\|z^\ell-z_*\|_2.
\]
When $T_\infty=0$, the geometric bound for the full auxiliary error therefore gives both output bounds.
One can use the same constant $C_{\rm loc}$ in all three inequalities.

\paragraph{Directions unchanged by the linearized DR map.}
Since $R$ is symmetric with spectrum in $[0,1]$,
\[
 \ker(I-T)
 =(\mathcal S_\beta^0\cap\ker(I-R))\oplus((\mathcal S_\beta^0)^\perp\cap\ker R).
\]
To verify the identity, first suppose $T\delta x=\delta x$.
Equation \eqref{eq:T} gives $P(2R-I)\delta x=R\delta x$, so $R\delta x$ belongs to $\mathcal S_\beta^0$.
Hence $PR\delta x=R\delta x$, and $P(2R-I)\delta x=R\delta x$ becomes $2R\delta x-P\delta x=R\delta x$, that is, $P\delta x=R\delta x$.
Since $R$ is symmetric with spectrum in $[0,1]$, the matrix $R(I-R)$ is symmetric positive semidefinite.
Consequently,
\[
 \langle\delta x,R(I-R)\delta x\rangle
 =\langle R\delta x,(I-R)\delta x\rangle=\langle P\delta x,\delta x-P\delta x\rangle=0,
\]
so $R(I-R)\delta x=0$, that is, $\delta x\in\ker R\oplus\ker(I-R)$.
The eigenvalue-one component is $R\delta x=P\delta x\in\mathcal S_\beta^0$; the eigenvalue-zero component is $\delta x-P\delta x\in(\mathcal S_\beta^0)^\perp$.
Conversely, substitution of either type of component in \eqref{eq:T} gives $T\delta x=\delta x$.
For a positive geodesic, $R=(I+\mathcal H)^{-1}$ is invertible, so
$\ker(I-T)=\mathcal S_\beta^0\cap\ker\mathcal H$.
At zero outputs, $(\mathcal S_\beta^0)^\perp\cap\ker R$ can be nonzero.
For $\delta x$ in this space, $R\delta x=0$ and $P(2R-I)\delta x=0$, so replacing $z_*$ by $z_*+t\,\delta x$ changes neither output to first order in $t$.
Such a fixed direction of $T$ concerns the auxiliary iterate $z$, not the minimizer $x_*$.
This proves \cref{prop:lfp}.
\end{proof}

The following identity and estimate quantify the local geometry described after \eqref{eq:scaledkappa}.
By the definition \eqref{eq:kappa} of $\kappa$, $s=v/q=\sqrt\beta\,m/q=\sqrt\beta\,h\kappa$ for $q>0$; this proves \eqref{eq:scaledkappa}.
At $q=0$, finite action forces $m=0$, but neither ratio is defined.
For the action block $R_a=(I+bH)^{-1}$, let $u_m=(0,1)$ denote a unit momentum perturbation in $(q,v)$ coordinates.
The rank-one Hessian formula gives
\[
 \|(I-R_a)u_m\|_2
 =\frac{b\sqrt{1+s^2}}{q+b(1+s^2)}
 \le\frac1{\sqrt{1+s^2}}.
\]
Thus $R_a$ changes a unit momentum perturbation by at most $(1+s^2)^{-1/2}$.
This quantifies the statement after \eqref{eq:scaledkappa} for large $|s|$; it concerns the action block $R_a$ of $R$, not $T$.

\begin{remark}[Rotation of the flat direction with $\beta$]
By \eqref{eq:scaledkappa}, for a cell with $\kappa\ne0$, $|s|=\sqrt\beta\,h|\kappa|>1$ exactly when $\beta>1/(h^2\kappa^2)$; the flat direction $(1,s)$ then makes an angle smaller than $45^\circ$ with the momentum axis.
For $h=1/64$, as in \cref{fig:geocolormaps}, this threshold is $1.8\cdot10^5$ at $|\kappa|=0.15$, the limit of the color scale of \cref{fig:geocolormaps}(b), and $1.6\cdot10^6$ at $|\kappa|=0.05$.
The threshold describes the local geometry; it is not a proved convergence threshold.
\end{remark}

\subsection{Midpoint-weight contrast in the constrained geodesic Hessian}\label{app:geocontrast}\label{par:geocontrast}

This subsection rewrites \eqref{eq:constrained-hessian} as $\mu_s>0$ and bounds the condition number of the constrained Hessian by $C_q\Lambda_s/\mu_s$, where the midpoint-weight contrast $C_q$ and the constants $\mu_s$ and $\Lambda_s$ are defined below.
Assume that all density entries and midpoint weights of $x_*$ are positive.
Use unweighted Euclidean array norms in scaled coordinates.
For a perturbation $\delta x=(\delta\nu,\delta v,\delta q)$, define the linear map
\[
 (B_s\delta x)_i^j=\delta v_i^j-s_i^j\delta q_i^j,
 \qquad s_i^j=v_{*,i}^j/q_{*,i}^j.
\]
The cell Hessian in \eqref{eq:rankone} gives the exact identity
\[
 \langle\delta x,\mathcal H\delta x\rangle
 =b\sum_{i,j}\frac{|(B_s\delta x)_i^j|^2}{q_{*,i}^j},
 \qquad \mathcal H=\gamma\nabla^2J(x_*).
\]
The part of $J$ depending directly on the nodal density is the nonnegativity indicator.
At a positive density it is locally zero, so this part contributes no second derivative.
Define $q_{\min}=\min_{i,j}q_{*,i}^j$, $q_{\max}=\max_{i,j}q_{*,i}^j$, and $C_q=q_{\max}/q_{\min}$.
In cell $(i,j)$, the local Hessian $(b/q_{*,i}^j)\begin{pmatrix}(s_i^j)^2&-s_i^j\\-s_i^j&1\end{pmatrix}$ in the variables $(\delta q,\delta v)$ has eigenvalue $0$ in the direction $(1,s_i^j)$ and eigenvalue $b(1+(s_i^j)^2)/q_{*,i}^j$ in the direction $(s_i^j,-1)$.
If $s_i^j$ does not depend on $(i,j)$, the largest positive eigenvalue divided by the smallest is $C_q$.
For every perturbation $\delta x$,
\[
 \frac b{q_{\max}}\|B_s\delta x\|_2^2
 \le\langle\delta x,\mathcal H\delta x\rangle
 \le\frac b{q_{\min}}\|B_s\delta x\|_2^2.
\]
Restrict to perturbations in $\mathcal S_\beta^0$ and set
\[
 \mu_s=\min_{\substack{\delta x\in\mathcal S_\beta^0\\\|\delta x\|_2=1}}\|B_s\delta x\|_2^2,
 \qquad
 \Lambda_s=\max_{\substack{\delta x\in\mathcal S_\beta^0\\\|\delta x\|_2=1}}\|B_s\delta x\|_2^2.
\]
These extrema exist because the unit sphere of $\mathcal S_\beta^0$ is compact.
Suppose $\mu_s>0$.
Let $\mathcal H_{\mathcal S}=P\mathcal H|_{\mathcal S_\beta^0}$ be the Hessian of $\gamma J$ restricted to $\mathcal S_\beta$, evaluated at $x_*$.
Its smallest and largest eigenvalues satisfy
\[
 \lambda_{\min}(\mathcal H_{\mathcal S})\ge\frac{b\mu_s}{q_{\max}},
 \qquad
 \lambda_{\max}(\mathcal H_{\mathcal S})\le\frac{b\Lambda_s}{q_{\min}}.
\]
Consequently,
\[
 \operatorname{cond}_2(\mathcal H_{\mathcal S})
 :=\frac{\lambda_{\max}(\mathcal H_{\mathcal S})}{\lambda_{\min}(\mathcal H_{\mathcal S})}
 \le C_q\frac{\Lambda_s}{\mu_s}.
\]
The factor $\Lambda_s/\mu_s$ depends on the ratios $s$, on the affine constraints, and on the scaling parameter $\beta$.
Thus a condition-number bound based on $C_q$ also requires control of $\Lambda_s/\mu_s$.
Since $q_*>0$ and $\mathcal H$ is positive semidefinite, the identity for $\langle\delta x,\mathcal H\delta x\rangle$ gives $\ker\mathcal H=\ker B_s$.
Hence $\mu_s>0$ if and only if $\mathcal S_\beta^0\cap\ker\mathcal H=\{0\}$, which is equivalent to \eqref{eq:constrained-hessian}; then $T_\infty=0$ by the identity $\ker(I-T)=\mathcal S_\beta^0\cap\ker\mathcal H$ of \cref{app:lfp}.
The bound on $\operatorname{cond}_2(\mathcal H_{\mathcal S})$ does not bound $\spr(T)$, which requires the matrix $T$.
At zero midpoint weights the Hessian formula is unavailable, so the local DR analysis must instead use the blocks of the Jacobian $R$ at zero outputs described in \cref{app:lfp}.

\subsection{Parameter balance for the entropic model}\label{app:entropicbalance}
Fix $\nu_i^{\rm old}>0$ and $\omega>0$.
In the model of \cref{sec:drflow}, the model Jacobian $R_{\rm model}$ has the two diagonal entries $\varrho_i=\nu_i^{\rm old}/(\nu_i^{\rm old}+\omega)$ at each node.
Write $\nu_{\min}$ and $\nu_{\max}$ for the extrema of $\nu^{\rm old}$ and $C=\nu_{\max}/\nu_{\min}$, as in \cref{sec:drflow}.
Let $P$ be any Euclidean orthogonal projector and set
\[
 T_{\rm model}=I-R_{\rm model}+P(2R_{\rm model}-I)
 =\tfrac12\bigl[I+(2P-I)(2R_{\rm model}-I)\bigr].
\]
Since $2P-I$ is orthogonal, its operator norm is one.
The triangle inequality gives
\[
 \begin{aligned}
 \|T_{\rm model}\|_2
 &\le\tfrac12\bigl(1+\|2R_{\rm model}-I\|_2\bigr)\\
 &=\max_i\max\{\varrho_i,1-\varrho_i\}\\
 &=\max\left\{\frac{\nu_{\max}}{\nu_{\max}+\omega},
              \frac{\omega}{\nu_{\min}+\omega}\right\}.
 \end{aligned}
\]
This proves the bound without assuming that $P$ and $R_{\rm model}$ commute.
In the entropic problem, $P$ is diagonal by spatial Fourier blocks because $L$ has constant coefficients.
Multiplication by a spatially varying $\varrho_i$ becomes convolution in Fourier coordinates, so $R_{\rm model}$ can mix frequencies even though $P$ does not.
The norm bound uses only the extrema of $\varrho_i$ and discards this spatial information.
Its two terms are respectively decreasing and increasing functions of $\omega$.
Their unique intersection satisfies
\[
 \nu_{\max}(\nu_{\min}+\omega)
 =\omega(\nu_{\max}+\omega)
 \quad\Longleftrightarrow\quad \omega^2=\nu_{\min}\nu_{\max}.
\]
For smaller $\omega$, the first term is larger than its value at the intersection; for larger $\omega$, the second term is larger.
The maximum is therefore minimized at $\omega=\sqrt{\nu_{\min}\nu_{\max}}$.
Substitution gives $\|T_{\rm model}\|_2\le c_{\rm model}$ with $c_{\rm model}=\sqrt C/(1+\sqrt C)$.
For the linear recurrence $\delta z^{\ell+1}=T_{\rm model}\delta z^\ell$, the bound is $\|\delta z^\ell\|_2\le c_{\rm model}^\ell\|\delta z^0\|_2$.
To guarantee $\|\delta z^\ell\|_2\le\varepsilon_{\rm red}\|\delta z^0\|_2$ for a prescribed reduction factor $\varepsilon_{\rm red}\in(0,1)$, it is sufficient to take
\[
 \ell\ge\frac{\log(1/\varepsilon_{\rm red})}{-\log c_{\rm model}}.
\]
Since $-\log c_{\rm model}=\log(1+C^{-1/2})\sim C^{-1/2}$ as $C\to\infty$, this sufficient iteration count grows as $O(\sqrt C)$ when $\varepsilon_{\rm red}$ is fixed.
The bound $\|T_{\rm model}\|_2\le c_{\rm model}$ concerns the matrix $T_{\rm model}$ only.
The matrix $T_{\rm model}$ differs from the DR matrix $T$ at the minimizer because $R_{\rm model}$ replaces the minimizing density by $\nu^{\rm old}$ and the momentum-to-density ratio by zero.
No bound on $\|T-T_{\rm model}\|_2$ is proved here.

For the three inputs of \cref{tab:common-entropic}, \cref{tab:entropic-parameter-check} gives $\spr(T)$ at the last DR iterate for $\omega=c_\omega\sqrt{\nu_{\min}\nu_{\max}}$, $c_\omega=0.3,1,3$.
In every row, $c_\omega=1$ has the smallest $\spr(T)$ and the smallest iteration count in \cref{tab:common-entropic}.
For $M=64$, $\varepsilon=10^{-6}$, $\spr(T)$ is smaller at $c_\omega=0.3$ than at $c_\omega=3$, but \cref{tab:common-entropic} records more iterations at $c_\omega=0.3$ (21750) than at $c_\omega=3$ (18900).
Hence the final $\spr(T)$ does not order the three factors by iteration count.

\begin{table}[htbp]
\centering\small
\begin{tabular}{@{}rrrlr@{}}
\toprule
$M$ & $p$ & $\varepsilon$ & $\spr(T)$ for $c_\omega=0.3,1,3$ & Scan time (s)\\
\midrule
32 & 4 & $1\cdot10^{-3}$ & 0.989787 / 0.973159 / 0.990891 & 0.502\\
64 & 16 & $1\cdot10^{-3}$ & 0.986445 / 0.978457 / 0.992715 & 0.608\\
64 & 16 & $1\cdot10^{-6}$ & 0.999538 / 0.999467 / 0.999822 & 10.243\\
\bottomrule
\end{tabular}
\caption{Spectral radius $\spr(T)$ of \eqref{eq:T} for the three entropic inputs of \cref{tab:common-entropic}, with $\beta=1/\tau$, $\omega=c_\omega\sqrt{\nu_{\min}\nu_{\max}}$ and $c_\omega=0.3,1,3$. The matrix $T$ is formed from the Jacobian $R$ at the last DR iterate of each run and from the projector $P$. Scan time is the sum of the three DR Total times of \cref{tab:common-entropic} ($c_\omega=0.3,1,3$). The iteration counts, in which $c_\omega=1$ is smallest, are in \cref{tab:common-entropic}; a smaller final spectral radius need not predict the total work from the initial iterate.}
\label{tab:entropic-parameter-check}
\end{table}

\subsection{Numerical comparison of local contraction factors}\label{app:numericalspectrum}
This test compares $\spr(T)$ with the late decay of the DR residual $r_\ell$ for the peaked pair at $M=15$, $N=16$ and the four pairs $(\beta,b)$ of \cref{tab:spectralchecks}.
Let $z^{\ell_{\rm end}}$ be the last auxiliary iterate of a run, with components $(\nu_{\rm in},v_{\rm in},q_{\rm in})$.
Form $R$ blockwise at $z^{\ell_{\rm end}}$: $(I+bH)^{-1}$ with $H$ from \eqref{eq:rankone} on cells where the action output is positive, $0$ on cells where $2b\,q_{\rm in}+v_{\rm in}^2<0$, and $1$ or $0$ for each nodal density according to the sign of its input.
Apply the linear part $P$ of the affine projection to coordinate vectors and form $T=I-R+P(2R-I)$.
The branch margin is the smallest of the numbers $|\nu_{{\rm in},i}^j|$ and $|q_{{\rm in},i}^j+(v_{{\rm in},i}^j)^2/(2b)|$ over all entries of $z^{\ell_{\rm end}}$; their signs select the branches.
In the four runs this margin is positive, so $\prox_{\gamma J}$ is differentiable at $z^{\ell_{\rm end}}$; this does not verify \eqref{eq:dr-strict-input} at $z_*$.

To separate fixed and contracting directions, use the decomposition of $\ker(I-T)$ proved in \cref{app:lfp} (paragraph on the directions unchanged by the linearized DR map).
Construct orthonormal bases $\Xi_1,\Xi_0$ for $\ker(I-R),\ker R$ from the local blocks.
The null spaces of $(I-P)\Xi_1$ and $P\Xi_0$ determine the intersections.
Singular-value thresholds are varied as specified in \cref{tab:spectralchecks}; a dimension that changes with the threshold would be marked unresolved.
Let the columns of $\Xi_{\rm fix}$ be an orthonormal basis of $\ker(I-T)$ and the columns of $\Xi_\perp$ an orthonormal basis of its orthogonal complement.
We report $\spr_{\rm an}=\spr(\Xi_\perp^\top T\Xi_\perp)$, the spectral radius of the linearized map $T$ on the orthogonal complement of its fixed directions; it equals $\spr(T-T_\infty)$ and is not a proved rate for the nonlinear iteration.
No fixed direction is detected in these four tests.

For the observed decay, let $r_\ell$ be the DR residual of \cref{sec:geodata}.
Keep the indices $\ell$ in the last two thirds of the run with $r_\ell>10^3\min_{\ell'}r_{\ell'}$; call these residuals eligible.
For a window length $\Delta\ell\in\{500,1000,2000\}$, or the whole eligible segment if it is shorter, compute $(r_{\ell+\Delta\ell}/r_\ell)^{1/\Delta\ell}$ on the last eligible window, and report the range over the three windows.
To check $R$, perturb all density entries, all momentum entries or all weight entries of $z^{\ell_{\rm end}}$ by $\pm\epsilon_{\rm FD}$, for each step $\epsilon_{\rm FD}\in\{10^{-5},10^{-7},10^{-9}\}$.
Since $\prox_{\gamma J}$ acts separately on nodal-density entries and action pairs, these three centered differences of $\prox_{\gamma J}$, divided by $2\epsilon_{\rm FD}$, give all entries of a block-diagonal matrix $R^{\rm FD}$.
With $T^{\rm FD}=I-R^{\rm FD}+P(2R^{\rm FD}-I)$, the column FD difference of \cref{tab:spectralchecks} is the largest value of $|\spr(\Xi_\perp^\top T^{\rm FD}\Xi_\perp)-\spr_{\rm an}|$ over the three steps.
In the four runs, $\spr_{\rm an}$ differs from the observed factors by at most $6.7\cdot10^{-5}$, except for $(\beta,b)=(10^6,0.1)$, where the difference is $6.4\cdot10^{-4}$; it differs from the finite-difference values by less than $10^{-7}$ (\cref{tab:spectralchecks}).
These are evaluations at a final iterate: they are not a bound on the convergence rate, and they do not show that the final local rate predicts the whole iteration history.

\begin{table}[htbp]
\centering\small
\begin{tabular}{@{}rrrrrr@{}}
\toprule
$\beta$ & $b$ & $\spr_{\rm an}$ & Observed range & FD difference & Fixed dim.\\
\midrule
$1\cdot10^{4}$ & 1 & 0.9941444 & 0.9940783--0.9940812 & $7.23\cdot10^{-8}$ & 0\\
$1\cdot10^{4}$ & 10 & 0.9919838 & 0.9919739--0.9919759 & $5.95\cdot10^{-10}$ & 0\\
$1\cdot10^{6}$ & 0.1 & 0.9818062 & 0.9811645--0.9811645 & $3.19\cdot10^{-8}$ & 0\\
$1\cdot10^{6}$ & 1 & 0.9962664 & 0.9962551--0.9962638 & $2.92\cdot10^{-9}$ & 0\\
\bottomrule
\end{tabular}
\caption{Analytic local DR spectra for the peaked pair, $M=15$, $N=16$, using the residual-checked proximal solver of \cref{alg:geo-prox-checked} and stopping when the DR residual $r_\ell$ of \cref{sec:geodata} satisfies $r_\ell\le10^{-12}$ (cap 16000 iterations). No run reached the cap. Here $\spr_{\rm an}=\spr(\Xi_\perp^\top T\Xi_\perp)$ (\cref{app:numericalspectrum}). Fixed-space dimensions agree at singular-value tolerances $10^{-8},10^{-10},10^{-12}$. Observed ranges use the window lengths $\Delta\ell=500,1000,2000$ of \cref{app:numericalspectrum}, ending at the last eligible residual, or the whole eligible segment if it is shorter. The numbers of eligible residuals are 1457, 1003, 495 and 1983, in the order of the rows; for $(\beta,b)=(10^6,0.1)$ the three windows are the whole segment of 495 residuals, so the observed range reduces to one value. FD difference is $\max|\spr(\Xi_\perp^\top T^{\rm FD}\Xi_\perp)-\spr_{\rm an}|$ over the centered-difference steps $10^{-5},10^{-7},10^{-9}$ (\cref{app:numericalspectrum}). These are floating-point diagnostics at a final iterate, not bounds at the exact fixed point.}
\label{tab:spectralchecks}
\end{table}

\section{The affine projection in Fourier coordinates}\label{app:projection}

At each nonzero spatial frequency, the affine constraints express momentum and midpoint weights in terms of the free density values.
Substitution leaves, at each such frequency, a strictly convex quadratic problem in the $N-1$ interior density coefficients.

Fix a Fourier frequency $k$.
Write $\boldsymbol\nu_k=(\widehat\nu_k^1,\ldots,\widehat\nu_k^{N-1})$ for the unknown interior density coefficients.
For a vector $u=(u^1,\ldots,u^{N-1})\in\mathbb C^{N-1}$, extended by $u^0=u^N=0$, define
\[
 (Du)^j=u^{j+1}-u^j,\qquad
 (Bu)^j=\tfrac12(u^j+u^{j+1}),\qquad0\le j<N.
\]

These definitions give matrices $D,B:\mathbb C^{N-1}\to\mathbb C^N$.
Adjoints in this appendix use the Euclidean Hermitian inner product.

The specified endpoint terms are
\[
 d_\partial^j=-\widehat\nu_k^0\,\ind_{\{j=0\}}
          +\widehat\nu_k^N\,\ind_{\{j=N-1\}},\qquad
 q_\partial^j=\tfrac12\widehat\nu_k^0\,\ind_{\{j=0\}}
          +\tfrac12\widehat\nu_k^N\,\ind_{\{j=N-1\}}.
\]
The indicator $\ind_{\{j=j_0\}}$ is one if $j=j_0$ and zero if $j\ne j_0$.

Let $(\widehat\nu_{\rm in},\widehat m_{\rm in},\widehat q_{\rm in})$ be the three Fourier blocks of the projection input.
For $k\ne0$, feasibility gives
\begin{equation}\label{eq:projection-recovery}
 \widehat m=\widehat L(k)^{-1}(D\boldsymbol\nu_k+d_\partial),\qquad
 \widehat q=B\boldsymbol\nu_k+q_\partial.
\end{equation}

By Parseval's identity, the projection minimizes this quadratic function of $\boldsymbol\nu_k$:
\[
 \frac12\|\boldsymbol\nu_k-\widehat\nu_{\rm in}\|_2^2
 +\frac\beta2\|\widehat L(k)^{-1}(D\boldsymbol\nu_k+d_\partial)-\widehat m_{\rm in}\|_2^2
 +\frac12\|B\boldsymbol\nu_k+q_\partial-\widehat q_{\rm in}\|_2^2.
\]
Since $\overline{\widehat L(k)}^{\,-1}=\widehat L(k)/|\widehat L(k)|^2$, setting the gradient of this function to zero gives the first-order condition for $\boldsymbol\nu_k$:
\[
 \begin{aligned}
 \left(I+\frac\beta{|\widehat L(k)|^2}D^*D+B^*B\right)\boldsymbol\nu_k
 =\widehat\nu_{\rm in}
 +\frac\beta{|\widehat L(k)|^2}D^*(\widehat L(k)\widehat m_{\rm in}-d_\partial)
 +B^*(\widehat q_{\rm in}-q_\partial).
 \end{aligned}
\]
The matrix is tridiagonal and positive definite, since it is $I$ plus two positive semidefinite matrices.

The normalized sine vectors are
\[
 \varphi_n^j=\sqrt{\frac2N}\sin\frac{\pi nj}{N},\qquad
 1\le j,n<N.
\]
Direct multiplication gives
\[
 D^*D\varphi_n=d_n\varphi_n,\quad d_n=4\sin^2\frac{\pi n}{2N},\qquad
 B^*B=I-\tfrac14D^*D.
\]
The eigenvalue of the projection system on the $n$th sine vector is
\begin{equation}\label{eq:projection}
 2-\frac{d_n}{4}+\frac{\beta d_n}{|\widehat L(k)|^2},\qquad1\le n<N.
\end{equation}
The eigenvalues in \eqref{eq:projection} reduce the solve for $\boldsymbol\nu_k$ to scalar divisions after a temporal sine transform, as implemented in \cref{alg:geo-project}.
Then recover $\widehat m$ and $\widehat q$ by \eqref{eq:projection-recovery}.

For $k=0$, the endpoint masses are one and $\widehat L(0)=0$.
The constraints give $\widehat\nu_0^j=1$ and $\widehat q_0^j=1$.
They impose no condition on $\widehat m_0^j$.
The nearest momentum thus has $\widehat m_0^j=\widehat m_{{\rm in},0}^j$.

\section{DR spectrum at the uniform path}\label{app:spectrum}

\subsection{Mode angles and convergence factor}
The calculation below gives the rate in \eqref{eq:rate} and determines how grid refinement changes it at fixed $\beta$.
For a nonzero spatial Fourier frequency $k$, use $\widehat L(k)=-4\pi^2k^2$.
For $k\ne0$ and $1\le n<N$, the perturbations in $\mathcal S_\beta^0$ of spatial frequency $k$ and temporal sine mode $n$ are the multiples of the mode
\[
 \delta x_{k,n}=(\varphi_n,c_kD\varphi_n,B\varphi_n),\qquad c_k=\frac{\sqrt\beta}{\widehat L(k)},
\]
with $D$, $B$ and $\varphi_n$ from \cref{app:projection}.
Let $\theta_{k,n}$ be the angle between $\delta x_{k,n}$ and the subspace $\delta v=0$ in scaled Euclidean coordinates, and define $\sigma_{k,n}>0$ by
\[
 \sin^2\theta_{k,n}
 =\frac{\|\delta v\|_2^2}
 {\|\delta\nu\|_2^2+\|\delta v\|_2^2+\|\delta q\|_2^2}
 =\sigma_{k,n}^2,\qquad(\delta\nu,\delta v,\delta q)=\delta x_{k,n}.
\]
Let $K=\lfloor M/2\rfloor$, the largest value of $|k|$ for $k\in\mathcal I_M$.
The calculation below treats each index $k\in\mathcal I_M$ separately, so it does not require $\mathcal I_M$ to be symmetric.
For even $M$, the unpaired index $k=-M/2$ corresponds to the real eigenvector $((-1)^i)_{0\le i<M}$ of $L$, with eigenvalue $\widehat L(-M/2)=-\pi^2M^2$.
The conservation and midpoint equations give
\begin{equation}\label{eq:sigma}
 t_{k,n}=\frac{\beta d_n}{|\widehat L(k)|^2},\qquad
 \sigma_{k,n}^2=\frac{t_{k,n}}{2-d_n/4+t_{k,n}},
 \qquad d_n=4\sin^2\frac{\pi n}{2N}.
\end{equation}
Since $\sigma_{k,n}^2=1/\bigl(1+(2-d_n/4)/t_{k,n}\bigr)$ increases with $t_{k,n}$ and with $d_n$, and since $d_n$ and $t_{k,n}$ increase with $n$ while $t_{k,n}$ decreases with $|k|$, the minimum of $\sigma_{k,n}$ over $k\ne0$, $1\le n<N$ is attained at $|k|=K$, $n=1$.

Thus the smallest angle in \cref{sec:dr-uniform} is $\theta_{\min}=\theta_{K,1}$.

\begin{proposition}[Convergence factor for the uniform path]\label{prop:constant}
Take uniform endpoints $\nu^0=\nu^N=\one$, $M\ge3$, $N\ge2$, and fixed $\beta,b>0$.
Use $K=\lfloor M/2\rfloor$ and $\sigma_{k,n}$ from \eqref{eq:sigma}.
At the uniform path, the spectral radius of $T$ is given by \eqref{eq:rate} when $b=1$.
For each $b>0$,
\[
 1-\spr(T)\le\sqrt{\frac b{1+b}}\,\sigma_{K,1}\le\sigma_{K,1}.
\]
If $b=1$ and $\beta\ge|\widehat L(K)|^2/d_1$, then $\spr(T)\le2/3$.
\end{proposition}

For the mode $\delta x_{k,n}=(\delta\nu,\delta v,\delta q)$, conservation and the midpoint equations give
\[
 \|\delta v\|_2^2=t_{k,n}\|\delta\nu\|_2^2,\qquad
 \|\delta q\|_2^2=(1-d_n/4)\|\delta\nu\|_2^2.
\]
Consequently, the angle between $\delta x_{k,n}$ and the subspace $\delta v=0$ satisfies
\[
 \sin^2\theta_{k,n}
 =\frac{\|\delta v\|_2^2}{\|\delta x_{k,n}\|_2^2}
 =\frac{t_{k,n}}{2-d_n/4+t_{k,n}}=\sigma_{k,n}^2.
\]
\begin{remark}[Grid refinement at fixed $\beta$]
At the uniform path, $R\,\delta x_{k,n}=(\delta\nu,r\,\delta v,\delta q)$ with $r=(1+b)^{-1}$, so one application of $R$ removes the fraction
\[
 \frac{\|\delta x_{k,n}\|_2^2-\|R\,\delta x_{k,n}\|_2^2}{\|\delta x_{k,n}\|_2^2}
 =(1-r^2)\sigma_{k,n}^2
\]
of the squared norm of the mode $\delta x_{k,n}$; this fraction is smallest at $(k,n)=(K,1)$.
As $N\to\infty$ and $M\to\infty$ at fixed $\beta$, $t_{K,1}=\beta d_1/|\widehat L(K)|^2\to0$ and
\[
 \sigma_{K,1}^2\sim\frac{\beta}{32\pi^2N^2K^4}.
\]
Hence, at $b=1$, \eqref{eq:rate} gives $\spr(T)=\max\{1/2,1-\sigma_{K,1}^2\}\to1$, although all midpoint weights equal one, so that $C_q=1$.
No bound $\spr(T)\le g(C_q)<1$ with $g$ independent of $M$ and $N$ can hold at fixed $\beta$.
\end{remark}

The following argument shows that the uniform path satisfies \eqref{eq:constrained-hessian}, so that $T_\infty=0$ in \cref{prop:lfp}.
At this path, the null space of $\mathcal H=\gamma\nabla^2J(x_*)$ consists of perturbations with $\delta v=0$.
Conservation then gives $\delta\nu^{j+1}=\delta\nu^j$.
The perturbations at the prescribed geodesic endpoints are zero, so all density perturbations vanish.
The midpoint relation gives $\delta q=0$.
Hence $\mathcal S_\beta^0\cap\ker\mathcal H=\{0\}$ and $T_\infty=0$.

\begin{proof}[Proof of \cref{prop:constant}]
At the uniform path, the action blocks of $R$ are identical at every cell.
Fourier and sine coordinates therefore reduce $\delta z^+=T\delta z$ to invariant subspaces of dimension at most two, which we now construct.
For $k\ne0$, the modes $\delta x_{k,n}$, $1\le n<N$, span the perturbations in $\mathcal S_\beta^0$ at frequency $k$.
Their momentum components are $\delta v_{k,n}=(0,c_kD\varphi_n,0)$.
We have
\[
 \|\delta x_{k,n}\|_2^2=2-d_n/4+\beta d_n/|\widehat L(k)|^2,\qquad
 \|\delta v_{k,n}\|_2^2=\beta d_n/|\widehat L(k)|^2.
\]
Thus $\|\delta v_{k,n}\|_2^2/\|\delta x_{k,n}\|_2^2=\sigma_{k,n}^2$ from \eqref{eq:sigma}.

The planes $\operatorname{span}(\delta x_{k,n},\delta v_{k,n})$ are mutually orthogonal.
Set $r=(1+b)^{-1}$.
The linear maps $P$ and $R$ preserve each plane, because
\[
 \begin{aligned}
 &P\delta x_{k,n}=\delta x_{k,n},&&P\delta v_{k,n}=\sigma_{k,n}^2\delta x_{k,n},\\
 &R\delta x_{k,n}=\delta x_{k,n}+(r-1)\delta v_{k,n},&&R\delta v_{k,n}=r\,\delta v_{k,n}.
 \end{aligned}
\]
Choose the orthonormal basis
\[
 e_1=\frac{\delta v_{k,n}}{\|\delta v_{k,n}\|_2},\qquad
 e_2=\frac{\delta x_{k,n}-\delta v_{k,n}}{\|\delta x_{k,n}-\delta v_{k,n}\|_2}.
\]
In this basis, the normalized mode $\delta x_{k,n}/\|\delta x_{k,n}\|_2$ has coordinates $(\sigma,\sqrt{1-\sigma^2})$, where $\sigma=\sigma_{k,n}$.
The matrices are
\[
 R=\begin{pmatrix}r&0\\0&1\end{pmatrix},\qquad
 P=\begin{pmatrix}
 \sigma^2&\sigma\sqrt{1-\sigma^2}\\
 \sigma\sqrt{1-\sigma^2}&1-\sigma^2
 \end{pmatrix}.
\]
Let $T_{k,n}$ be the restriction of $T$ to $\operatorname{span}(\delta x_{k,n},\delta v_{k,n})$.
Substitution into \eqref{eq:T} gives its characteristic polynomial:
\begin{equation}\label{eq:poly}
 \det(\mu I-T_{k,n})=\mu^2-
 (1+\bar r-2\bar r\sigma_{k,n}^2)\mu
 +\bar r(1-\sigma_{k,n}^2),\qquad \bar r=1-r=\frac b{1+b}.
\end{equation}

The remaining momentum direction is perpendicular to every $D\varphi_n$ and gives eigenvalue $\bar r$.
The $N$ remaining density and midpoint-weight directions are unchanged by $R$ and annihilated by $P$, so \eqref{eq:T} gives eigenvalue zero on them.
This accounts for the orthogonal complement of the two-dimensional planes.
For $k=0$, only the momentum means are free directions of $\mathcal S_\beta^0$.
This mode contributes $r$ with multiplicity $N$ and zero with multiplicity $2N-1$.
Thus the plane decomposition accounts for all $3N-1$ variables per spatial frequency.

If $b=1$, the two roots of \eqref{eq:poly} are $1/2$ and $1-\sigma_{k,n}^2$.
Since $\sigma_{k,n}$ is smallest at $(k,n)=(K,1)$, as shown after \eqref{eq:sigma}, and $\sin\theta_{K,1}=\sigma_{K,1}$, this proves \eqref{eq:rate}.

For general $b>0$, write
\[
 T=\tfrac12\bigl(I+(2P-I)(2R-I)\bigr).
\]
The operator $2P-I$ is an isometry, and $\|2R-I\|_2\le1$.
Hence $\spr(T)\le1$.

Let $\mu_1,\mu_2$ be the roots in one plane.
Since $|\mu_j|\le\spr(T)$, the reverse triangle inequality gives
\[
 (1-\spr(T))^2\le|1-\mu_1|\,|1-\mu_2|
 =\det(I-T_{k,n})=\bar r\sigma_{k,n}^2.
\]

Use $(k,n)=(K,1)$ to get the bound in \cref{prop:constant}.
Finally, $\beta d_1/|\widehat L(K)|^2\ge1$ gives $t_{k,n}\ge1$ for each nonzero mode.
Thus $\sigma_{k,n}^2\ge t_{k,n}/(2+t_{k,n})\ge1/3$.
For $b=1$, this gives $\spr(T)\le2/3$.
This proves \cref{prop:constant}.
\end{proof}

\begin{remark}[Numerical values for the uniform path]
For example, a smallest angle of $5^\circ$ gives $\spr(T)\simeq0.9924$, whereas $45^\circ$ gives $\spr(T)=1/2$.
Increasing $\beta$ beyond the value that gives $45^\circ$ leaves this spectral radius at $1/2$.
By \cref{prop:constant}, at fixed $\beta$ every $b>0$ gives $1-\spr(T)\le\sigma_{K,1}$; hence no choice of $b$ makes $1-\spr(T)$ larger than $\sigma_{K,1}$.

At $M=63$, $N=64$, and $\beta=1$, the smallest $\sigma_{k,n}^2$ is approximately $8.4\cdot10^{-13}$.
At $b=1$, a sufficient weight $\beta$ for $\spr(T)\le2/3$ is about $6\cdot10^{11}$.
These are exact uniform-path predictions, not parameter recommendations for nonuniform geodesics.
\end{remark}

\subsection{Geometry of one linearized DR iteration}\label{app:drillustrations}
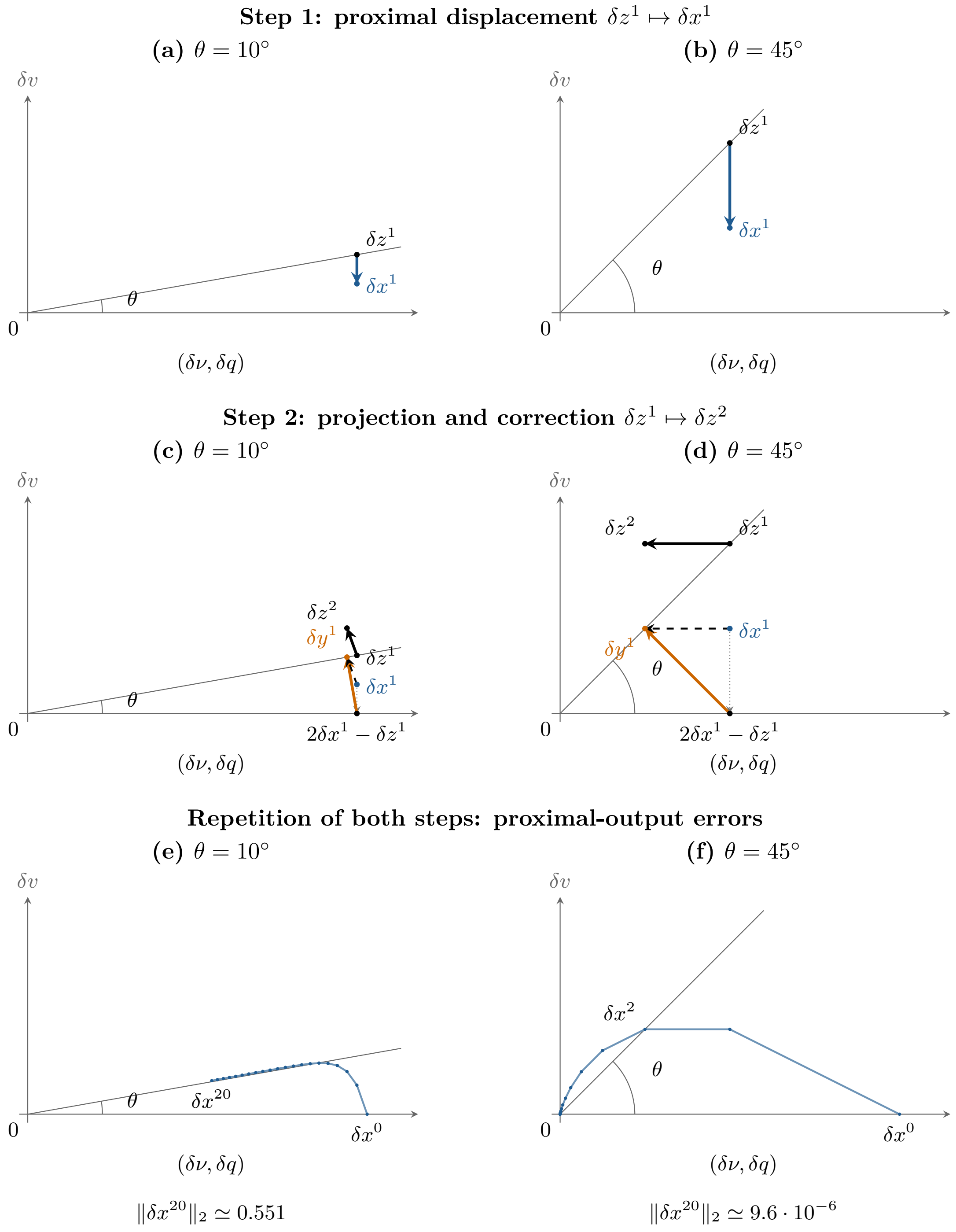
\begin{figure}[p]
\centering
{\small\bfseries Step 1: proximal displacement $\delta z^1\mapsto\delta x^1$\par}
\begin{tikzpicture}[font=\footnotesize,>=stealth]
\begin{scope}[xshift=0cm,x=5.1cm,y=5.1cm]
\node[font=\small\bfseries] at (.54,.77) {(a) $\theta=10^\circ$};
\draw[->,black!60] (-.025,0)--(1.15,0);
\draw[->,black!60] (0,-.025)--(0,.64) node[above] {$\delta v$};
\node at (.54,-.15) {$(\delta\nu,\delta q)$};
\node[below left] at (0,0) {$0$};
\draw[black!60] (0,0)--(1.10000000,0.19395968);
\draw[black!60] (.22,0) arc[start angle=0,end angle=10,radius=.22];
\node at (0.30882036,0.04201828) {$\theta$};
\draw[plotblue,very thick,->] (0.96984631,0.17101007)--(0.96984631,0.08550504);
\fill[black] (0.96984631,0.17101007) circle[radius=1.2pt] node[above right] {$\delta z^1$};
\fill[plotblue] (0.96984631,0.08550504) circle[radius=1.2pt] node[right] {$\delta x^1$};
\end{scope}
\begin{scope}[xshift=8cm,x=5.1cm,y=5.1cm]
\node[font=\small\bfseries] at (.54,.77) {(b) $\theta=45^\circ$};
\draw[->,black!60] (-.025,0)--(1.15,0);
\draw[->,black!60] (0,-.025)--(0,.64) node[above] {$\delta v$};
\node at (.54,-.15) {$(\delta\nu,\delta q)$};
\node[below left] at (0,0) {$0$};
\draw[black!60] (0,0)--(0.60000000,0.60000000);
\draw[black!60] (.22,0) arc[start angle=0,end angle=45,radius=.22];
\node at (0.28640266,0.13363186) {$\theta$};
\draw[plotblue,very thick,->] (0.50000000,0.50000000)--(0.50000000,0.25000000);
\fill[black] (0.50000000,0.50000000) circle[radius=1.2pt] node[above right] {$\delta z^1$};
\fill[plotblue] (0.50000000,0.25000000) circle[radius=1.2pt] node[right] {$\delta x^1$};
\end{scope}
\end{tikzpicture}

\medskip
{\small\bfseries Step 2: projection and correction $\delta z^1\mapsto\delta z^2$\par}
\begin{tikzpicture}[font=\footnotesize,>=stealth]
\begin{scope}[xshift=0cm,x=5.1cm,y=5.1cm]
\node[font=\small\bfseries] at (.54,.77) {(c) $\theta=10^\circ$};
\draw[->,black!60] (-.025,0)--(1.15,0);
\draw[->,black!60] (0,-.025)--(0,.64) node[above] {$\delta v$};
\node at (.54,-.15) {$(\delta\nu,\delta q)$};
\node[below left] at (0,0) {$0$};
\draw[black!60] (0,0)--(1.10000000,0.19395968);
\draw[black!60] (.22,0) arc[start angle=0,end angle=10,radius=.22];
\node at (0.30882036,0.04201828) {$\theta$};
\draw[black!55,densely dotted,->] (0.96984631,0.08550504)--(0.96984631,0.00000000);
\draw[orange!80!black,very thick,->] (0.96984631,0.00000000)--(0.94060187,0.16585349);
\draw[black,dashed,thick,->] (0.96984631,0.08550504)--(0.94060187,0.16585349);
\draw[black,very thick,->] (0.96984631,0.17101007)--(0.94060187,0.25135852);
\fill[black] (0.96984631,0.17101007) circle[radius=1.2pt] node[right] {$\delta z^1$};
\fill[plotblue] (0.96984631,0.08550504) circle[radius=1.2pt] node[right] {$\delta x^1$};
\fill[black] (0.96984631,0.00000000) circle[radius=1.2pt] node[below] {$2\delta x^1-\delta z^1$};
\fill[orange!80!black] (0.94060187,0.16585349) circle[radius=1.2pt] node[above left] {$\delta y^1$};
\fill[black] (0.94060187,0.25135852) circle[radius=1.2pt] node[above left] {$\delta z^2$};
\end{scope}
\begin{scope}[xshift=8cm,x=5.1cm,y=5.1cm]
\node[font=\small\bfseries] at (.54,.77) {(d) $\theta=45^\circ$};
\draw[->,black!60] (-.025,0)--(1.15,0);
\draw[->,black!60] (0,-.025)--(0,.64) node[above] {$\delta v$};
\node at (.54,-.15) {$(\delta\nu,\delta q)$};
\node[below left] at (0,0) {$0$};
\draw[black!60] (0,0)--(0.60000000,0.60000000);
\draw[black!60] (.22,0) arc[start angle=0,end angle=45,radius=.22];
\node at (0.28640266,0.13363186) {$\theta$};
\draw[black!55,densely dotted,->] (0.50000000,0.25000000)--(0.50000000,0.00000000);
\draw[orange!80!black,very thick,->] (0.50000000,0.00000000)--(0.25000000,0.25000000);
\draw[black,dashed,thick,->] (0.50000000,0.25000000)--(0.25000000,0.25000000);
\draw[black,very thick,->] (0.50000000,0.50000000)--(0.25000000,0.50000000);
\fill[black] (0.50000000,0.50000000) circle[radius=1.2pt] node[above right] {$\delta z^1$};
\fill[plotblue] (0.50000000,0.25000000) circle[radius=1.2pt] node[right] {$\delta x^1$};
\fill[black] (0.50000000,0.00000000) circle[radius=1.2pt] node[below] {$2\delta x^1-\delta z^1$};
\fill[orange!80!black] (0.25000000,0.25000000) circle[radius=1.2pt] node[below left] {$\delta y^1$};
\fill[black] (0.25000000,0.50000000) circle[radius=1.2pt] node[above left] {$\delta z^2$};
\end{scope}
\end{tikzpicture}

\medskip
{\small\bfseries Repetition of both steps: proximal-output errors\par}
\begin{tikzpicture}[font=\footnotesize,>=stealth]
\begin{scope}[xshift=0cm,x=5.1cm,y=5.1cm]
\node[font=\small\bfseries] at (.54,.77) {(e) $\theta=10^\circ$};
\draw[->,black!60] (-.025,0)--(1.15,0);
\draw[->,black!60] (0,-.025)--(0,.64) node[above] {$\delta v$};
\node at (.54,-.15) {$(\delta\nu,\delta q)$};
\node[below left] at (0,0) {$0$};
\draw[black!60] (0,0)--(1.10000000,0.19395968);
\draw[black!60] (.22,0) arc[start angle=0,end angle=10,radius=.22];
\node at (0.30882036,0.04201828) {$\theta$};
\draw[plotblue!65,thick] plot coordinates {(1.00000000,0.00000000) (0.96984631,0.08550504) (0.94060187,0.12567926) (0.91223925,0.14326583) (0.88473187,0.14963396) (0.85805394,0.15046601) (0.83218045,0.14860094) (0.80708714,0.14545609) (0.78275048,0.14173806) (0.75914767,0.13779814) (0.73625656,0.13381002) (0.71405571,0.12985865) (0.69252430,0.12598469) (0.67164213,0.12220666) (0.65138965,0.11853211) (0.63174784,0.11496315) (0.61269832,0.11149920) (0.59422320,0.10813839) (0.57630518,0.10487827) (0.55892745,0.10171613) (0.54207373,0.09864918)};
\fill[plotblue] (1.00000000,0.00000000) circle[radius=.7pt];
\fill[plotblue] (0.96984631,0.08550504) circle[radius=.7pt];
\fill[plotblue] (0.94060187,0.12567926) circle[radius=.7pt];
\fill[plotblue] (0.91223925,0.14326583) circle[radius=.7pt];
\fill[plotblue] (0.88473187,0.14963396) circle[radius=.7pt];
\fill[plotblue] (0.85805394,0.15046601) circle[radius=.7pt];
\fill[plotblue] (0.83218045,0.14860094) circle[radius=.7pt];
\fill[plotblue] (0.80708714,0.14545609) circle[radius=.7pt];
\fill[plotblue] (0.78275048,0.14173806) circle[radius=.7pt];
\fill[plotblue] (0.75914767,0.13779814) circle[radius=.7pt];
\fill[plotblue] (0.73625656,0.13381002) circle[radius=.7pt];
\fill[plotblue] (0.71405571,0.12985865) circle[radius=.7pt];
\fill[plotblue] (0.69252430,0.12598469) circle[radius=.7pt];
\fill[plotblue] (0.67164213,0.12220666) circle[radius=.7pt];
\fill[plotblue] (0.65138965,0.11853211) circle[radius=.7pt];
\fill[plotblue] (0.63174784,0.11496315) circle[radius=.7pt];
\fill[plotblue] (0.61269832,0.11149920) circle[radius=.7pt];
\fill[plotblue] (0.59422320,0.10813839) circle[radius=.7pt];
\fill[plotblue] (0.57630518,0.10487827) circle[radius=.7pt];
\fill[plotblue] (0.55892745,0.10171613) circle[radius=.7pt];
\fill[plotblue] (0.54207373,0.09864918) circle[radius=.7pt];
\node[below] at (1,0) {$\delta x^0$};
\node[below] at (0.54207373,0.09864918) {$\delta x^{20}$};
\node at (.54,-.29) {$\|\delta x^{20}\|_2\simeq 0.551$};
\end{scope}
\begin{scope}[xshift=8cm,x=5.1cm,y=5.1cm]
\node[font=\small\bfseries] at (.54,.77) {(f) $\theta=45^\circ$};
\draw[->,black!60] (-.025,0)--(1.15,0);
\draw[->,black!60] (0,-.025)--(0,.64) node[above] {$\delta v$};
\node at (.54,-.15) {$(\delta\nu,\delta q)$};
\node[below left] at (0,0) {$0$};
\draw[black!60] (0,0)--(0.60000000,0.60000000);
\draw[black!60] (.22,0) arc[start angle=0,end angle=45,radius=.22];
\node at (0.28640266,0.13363186) {$\theta$};
\draw[plotblue!65,thick] plot coordinates {(1.00000000,0.00000000) (0.50000000,0.25000000) (0.25000000,0.25000000) (0.12500000,0.18750000) (0.06250000,0.12500000) (0.03125000,0.07812500) (0.01562500,0.04687500) (0.00781250,0.02734375) (0.00390625,0.01562500) (0.00195312,0.00878906) (0.00097656,0.00488281) (0.00048828,0.00268555) (0.00024414,0.00146484) (0.00012207,0.00079346) (0.00006104,0.00042725) (0.00003052,0.00022888) (0.00001526,0.00012207) (0.00000763,0.00006485) (0.00000381,0.00003433) (0.00000191,0.00001812) (0.00000095,0.00000954)};
\fill[plotblue] (1.00000000,0.00000000) circle[radius=.7pt];
\fill[plotblue] (0.50000000,0.25000000) circle[radius=.7pt];
\fill[plotblue] (0.25000000,0.25000000) circle[radius=.7pt];
\fill[plotblue] (0.12500000,0.18750000) circle[radius=.7pt];
\fill[plotblue] (0.06250000,0.12500000) circle[radius=.7pt];
\fill[plotblue] (0.03125000,0.07812500) circle[radius=.7pt];
\fill[plotblue] (0.01562500,0.04687500) circle[radius=.7pt];
\fill[plotblue] (0.00781250,0.02734375) circle[radius=.7pt];
\fill[plotblue] (0.00390625,0.01562500) circle[radius=.7pt];
\fill[plotblue] (0.00195312,0.00878906) circle[radius=.7pt];
\fill[plotblue] (0.00097656,0.00488281) circle[radius=.7pt];
\fill[plotblue] (0.00048828,0.00268555) circle[radius=.7pt];
\fill[plotblue] (0.00024414,0.00146484) circle[radius=.7pt];
\fill[plotblue] (0.00012207,0.00079346) circle[radius=.7pt];
\fill[plotblue] (0.00006104,0.00042725) circle[radius=.7pt];
\fill[plotblue] (0.00003052,0.00022888) circle[radius=.7pt];
\fill[plotblue] (0.00001526,0.00012207) circle[radius=.7pt];
\fill[plotblue] (0.00000763,0.00006485) circle[radius=.7pt];
\fill[plotblue] (0.00000381,0.00003433) circle[radius=.7pt];
\fill[plotblue] (0.00000191,0.00001812) circle[radius=.7pt];
\fill[plotblue] (0.00000095,0.00000954) circle[radius=.7pt];
\node[below] at (1,0) {$\delta x^0$};
\node[above left] at (0.25000000,0.25000000) {$\delta x^2$};
\node at (.54,-.29) {$\|\delta x^{20}\|_2\simeq 9.6\cdot10^{-6}$};
\end{scope}
\end{tikzpicture}
\caption{Two steps of one linearized DR iteration \eqref{eq:T}, then twenty iterations, on one mode of the uniform path.
The angle $\theta=\theta_{k,n}$ of \eqref{eq:sigma} between the gray constraint line $\operatorname{span}(\delta x_{k,n})$ and the horizontal axis is $10^\circ$ in the left column and $45^\circ$ in the right column.
All panels use $b=1$, $\delta z^0=(1,0)$, and identical axis scales; the origin is zero error.
The horizontal axis is along $(\delta x_{k,n}-\delta v_{k,n})/\|\delta x_{k,n}-\delta v_{k,n}\|_2$ (density and midpoint weight, no momentum); the vertical axis is along $\delta v_{k,n}/\|\delta v_{k,n}\|_2$ (scaled momentum).
\textbf{(a,b)} Blue arrow: proximal displacement $\delta x^1-\delta z^1$, with $\delta z^1=P\delta z^0$.
\textbf{(c,d)} Dotted gray arrow: from $\delta x^1$ to the reflected point $2\delta x^1-\delta z^1$; orange arrow: orthogonal projection onto the constraint line; dashed and solid black arrows: DR correction $\delta y^1-\delta x^1$, based at $\delta x^1$ and at $\delta z^1$.
\textbf{(e,f)} Blue dots are the proximal-output errors $\delta x^\ell=RT^\ell\delta z^0$, $0\le\ell\le20$, joined in iteration order.
The final norms are shown below the plots.
These are trajectories of the linearized mode, not nonlinear geodesic solves.}
\label{fig:drgeometry}
\end{figure}

For one mode of the uniform path, this subsection draws the two steps of the linearized iteration \eqref{eq:T} and shows how the constraint angle $\theta_{k,n}$ of \eqref{eq:sigma}, whose minimum $\theta_{\min}$ appears in \eqref{eq:rate}, and the parameter $b$ in \cref{prop:constant} determine the error reduction.
\Cref{fig:drgeometry} follows one Fourier--sine mode $(k,n)$ of the linearized iteration, with $k\ne0$ and $1\le n<N$.
In the plane $\operatorname{span}(\delta x_{k,n},\delta v_{k,n})$ of \cref{app:spectrum}, the horizontal coordinate is along $(\delta x_{k,n}-\delta v_{k,n})/\|\delta x_{k,n}-\delta v_{k,n}\|_2$ (density and midpoint-weight components, zero momentum) and the vertical coordinate is along $\delta v_{k,n}/\|\delta v_{k,n}\|_2$ (scaled momentum).
These are the coordinates $(e_2,e_1)$ of \cref{app:spectrum}, so $R=\diag(1,r)$ with $r=(1+b)^{-1}$.
The gray constraint line spanned by $\delta x_{k,n}$ makes the angle $\theta=\theta_{k,n}$ with the horizontal axis, where $\sin\theta_{k,n}=\sigma_{k,n}$ in \eqref{eq:sigma}.
The panels treat $\theta$ as a free parameter and use $\theta=10^\circ$ and $\theta=45^\circ$.
We choose $b=1$, so $R$ halves the vertical coordinate and leaves the horizontal coordinate unchanged.
Write $\delta z^\ell=T^\ell\delta z^0$ for the auxiliary error, and $\delta x^\ell=R\,\delta z^\ell$ and $\delta y^\ell=P(2\delta x^\ell-\delta z^\ell)$ for the first-order parts of $x^\ell-x_*$ and $y^\ell-x_*$ in \eqref{eq:dr}.
We call $\delta x^\ell$ the proximal-output error.
Then $\delta z^{\ell+1}=\delta z^\ell+\delta y^\ell-\delta x^\ell$, which is \eqref{eq:T}.
Both columns start from the same auxiliary error $\delta z^0=(1,0)$ and use the same axis scales.
Since $R\,\delta z^0=\delta z^0$, the first DR iteration gives $\delta z^1=P\delta z^0$.
Panels (a--d) of \cref{fig:drgeometry} show the next iteration, from $\delta z^1$ to $\delta z^2$.

In panels (a,b), the blue arrow is the proximal displacement $\delta x^1-\delta z^1$, with $\delta x^1=R\,\delta z^1$.
In panels (c,d), the dotted gray arrow goes from $\delta x^1$ to the reflected point $2\delta x^1-\delta z^1=\delta x^1+(\delta x^1-\delta z^1)$.
The orange arrow goes from the reflected point to its orthogonal projection $\delta y^1=P(2\delta x^1-\delta z^1)$ onto the gray constraint line.
The DR correction $\delta y^1-\delta x^1$ is drawn dashed from $\delta x^1$ and solid from $\delta z^1$; the solid arrow ends at $\delta z^2=\delta z^1+(\delta y^1-\delta x^1)$.

With $\delta z^0=(1,0)$, $\delta z^1=\cos\theta\,(\cos\theta,\sin\theta)$, and at $b=1$, $\delta z^2=(\cos^4\theta,\ \cos\theta\sin\theta\,(\cos^2\theta+1/2))$.
The step from $\delta z^1$ to $\delta z^2$ thus multiplies the horizontal component of the auxiliary error by $\cos^2\theta$: by $0.970$ at $\theta=10^\circ$, and by $1/2$ at $\theta=45^\circ$, where the correction $\delta y^1-\delta x^1$ is horizontal.
Panels (e,f) show the proximal-output errors $\delta x^0,\ldots,\delta x^{20}$ obtained by repeating these two steps.

\begin{figure}[p]
\centering
{\small\bfseries 1. Proximal displacement from $\delta z^1=(1/2,1/2)$\par}
\begin{tikzpicture}[font=\footnotesize,>=stealth]
\begin{scope}[xshift=0cm,x=4.2cm,y=4.2cm]
\node[font=\small] at (.44,.88) {$b=0.3$};
\draw[->,black!60] (-.24,0)--(1.08,0);
\draw[->,black!60] (0,-0.02)--(0,.65) node[above] {$\delta v$};
\node at (.44,-0.15) {$(\delta\nu,\delta q)$};
\draw[black!60] (0,0)--(.60,.60);
\draw[black!60] (.24,0) arc[start angle=0,end angle=45,radius=.24];
\draw[plotblue,very thick,->] (0.50000000,0.50000000)--(0.50000000,0.38461538);
\fill[black] (0.50000000,0.50000000) circle[radius=1pt] node[above right] {$\delta z^1$};
\fill[plotblue] (0.50000000,0.38461538) circle[radius=1pt] node[right] {$\delta x^1$};
\end{scope}
\begin{scope}[xshift=8cm,x=4.2cm,y=4.2cm]
\node[font=\small] at (.44,.88) {$b=1$};
\draw[->,black!60] (-.24,0)--(1.08,0);
\draw[->,black!60] (0,-0.02)--(0,.65) node[above] {$\delta v$};
\node at (.44,-0.15) {$(\delta\nu,\delta q)$};
\draw[black!60] (0,0)--(.60,.60);
\draw[black!60] (.24,0) arc[start angle=0,end angle=45,radius=.24];
\draw[plotblue,very thick,->] (0.50000000,0.50000000)--(0.50000000,0.25000000);
\fill[black] (0.50000000,0.50000000) circle[radius=1pt] node[above right] {$\delta z^1$};
\fill[plotblue] (0.50000000,0.25000000) circle[radius=1pt] node[right] {$\delta x^1$};
\end{scope}
\end{tikzpicture}\par
\medskip
{\small\bfseries 2. Affine projection and DR correction\par}
\begin{tikzpicture}[font=\footnotesize,>=stealth]
\begin{scope}[xshift=0cm,x=4.2cm,y=4.2cm]
\node[font=\small] at (.44,.88) {$b=0.3$};
\draw[->,black!60] (-.24,0)--(1.08,0);
\draw[->,black!60] (0,-0.49)--(0,.65) node[above] {$\delta v$};
\node at (.44,-0.61) {$(\delta\nu,\delta q)$};
\draw[black!60] (-.12,-.12)--(.60,.60);
\draw[black!60] (.24,0) arc[start angle=0,end angle=45,radius=.24];
\draw[black!55,densely dotted,->] (0.50000000,0.38461538)--(0.50000000,0.26923077);
\draw[orange!80!black,thick,->] (0.50000000,0.26923077)--(0.38461538,0.38461538);
\draw[black,dashed,thick,->] (0.50000000,0.38461538)--(0.38461538,0.38461538);
\draw[black,thick,->] (0.50000000,0.50000000)--(0.38461538,0.50000000);
\fill[black] (0.50000000,0.50000000) circle[radius=1pt] node[above right] {$\delta z^1$};
\fill[plotblue] (0.50000000,0.38461538) circle[radius=1pt] node[right] {$\delta x^1$};
\fill[black] (0.50000000,0.26923077) circle[radius=1pt] node[below right] {$2\delta x^1-\delta z^1$};
\fill[orange!80!black] (0.38461538,0.38461538) circle[radius=1pt] node[below left] {$\delta y^1$};
\fill[black] (0.38461538,0.50000000) circle[radius=1pt] node[above left] {$\delta z^2$};
\end{scope}
\begin{scope}[xshift=8cm,x=4.2cm,y=4.2cm]
\node[font=\small] at (.44,.88) {$b=1$};
\draw[->,black!60] (-.24,0)--(1.08,0);
\draw[->,black!60] (0,-0.49)--(0,.65) node[above] {$\delta v$};
\node at (.44,-0.61) {$(\delta\nu,\delta q)$};
\draw[black!60] (-.12,-.12)--(.60,.60);
\draw[black!60] (.24,0) arc[start angle=0,end angle=45,radius=.24];
\draw[black!55,densely dotted,->] (0.50000000,0.25000000)--(0.50000000,0.00000000);
\draw[orange!80!black,thick,->] (0.50000000,0.00000000)--(0.25000000,0.25000000);
\draw[black,dashed,thick,->] (0.50000000,0.25000000)--(0.25000000,0.25000000);
\draw[black,thick,->] (0.50000000,0.50000000)--(0.25000000,0.50000000);
\fill[black] (0.50000000,0.50000000) circle[radius=1pt] node[above right] {$\delta z^1$};
\fill[plotblue] (0.50000000,0.25000000) circle[radius=1pt] node[right] {$\delta x^1$};
\fill[black] (0.50000000,0.00000000) circle[radius=1pt] node[below right] {$2\delta x^1-\delta z^1$};
\fill[orange!80!black] (0.25000000,0.25000000) circle[radius=1pt] node[below left] {$\delta y^1$};
\fill[black] (0.25000000,0.50000000) circle[radius=1pt] node[above left] {$\delta z^2$};
\end{scope}
\end{tikzpicture}\par
\medskip
{\small\bfseries 3. Proximal-output errors over twenty iterations\par}
\begin{tikzpicture}[font=\footnotesize,>=stealth]
\begin{scope}[xshift=0cm,x=4.2cm,y=4.2cm]
\node[font=\small] at (.44,.88) {$b=0.3$};
\draw[->,black!60] (-.24,0)--(1.08,0);
\draw[->,black!60] (0,-0.18)--(0,.65) node[above] {$\delta v$};
\node at (.44,-0.3) {$(\delta\nu,\delta q)$};
\draw[black!60] (-.12,-.12)--(.60,.60);
\draw[black!60] (.24,0) arc[start angle=0,end angle=45,radius=.24];
\draw[plotblue,thick] plot coordinates {(1.00000000,0.00000000) (0.50000000,0.38461538) (0.38461538,0.38461538) (0.32692308,0.34023669) (0.28254438,0.29585799) (0.24482249,0.25659991) (0.21222121,0.22246245) (0.18397246,0.19285477) (0.15948540,0.16718602) (0.13825781,0.14493355) (0.11985565,0.12564285) (0.10390282,0.10891975) (0.09007332,0.09442250) (0.07808454,0.08185484) (0.06769146,0.07095993) (0.05868171,0.06151514) (0.05087115,0.05332746) (0.04410019,0.04622956) (0.03823044,0.04007639) (0.03314196,0.03474221) (0.02873075,0.03011801)};
\fill[plotblue] (1.00000000,0.00000000) circle[radius=.65pt];
\fill[plotblue] (0.50000000,0.38461538) circle[radius=.65pt];
\fill[plotblue] (0.38461538,0.38461538) circle[radius=.65pt];
\fill[plotblue] (0.32692308,0.34023669) circle[radius=.65pt];
\fill[plotblue] (0.28254438,0.29585799) circle[radius=.65pt];
\fill[plotblue] (0.24482249,0.25659991) circle[radius=.65pt];
\fill[plotblue] (0.21222121,0.22246245) circle[radius=.65pt];
\fill[plotblue] (0.18397246,0.19285477) circle[radius=.65pt];
\fill[plotblue] (0.15948540,0.16718602) circle[radius=.65pt];
\fill[plotblue] (0.13825781,0.14493355) circle[radius=.65pt];
\fill[plotblue] (0.11985565,0.12564285) circle[radius=.65pt];
\fill[plotblue] (0.10390282,0.10891975) circle[radius=.65pt];
\fill[plotblue] (0.09007332,0.09442250) circle[radius=.65pt];
\fill[plotblue] (0.07808454,0.08185484) circle[radius=.65pt];
\fill[plotblue] (0.06769146,0.07095993) circle[radius=.65pt];
\fill[plotblue] (0.05868171,0.06151514) circle[radius=.65pt];
\fill[plotblue] (0.05087115,0.05332746) circle[radius=.65pt];
\fill[plotblue] (0.04410019,0.04622956) circle[radius=.65pt];
\fill[plotblue] (0.03823044,0.04007639) circle[radius=.65pt];
\fill[plotblue] (0.03314196,0.03474221) circle[radius=.65pt];
\fill[plotblue] (0.02873075,0.03011801) circle[radius=.65pt];
\node[below] at (1,0) {$\delta x^0$};
\node[right] at (0.50000000,0.38461538) {$\delta x^1$};
\node at (.44,-.52) {$\|\delta x^{20}\|_2\simeq 4.2\cdot10^{-2}$};
\end{scope}
\begin{scope}[xshift=8cm,x=4.2cm,y=4.2cm]
\node[font=\small] at (.44,.88) {$b=1$};
\draw[->,black!60] (-.24,0)--(1.08,0);
\draw[->,black!60] (0,-0.18)--(0,.65) node[above] {$\delta v$};
\node at (.44,-0.3) {$(\delta\nu,\delta q)$};
\draw[black!60] (-.12,-.12)--(.60,.60);
\draw[black!60] (.24,0) arc[start angle=0,end angle=45,radius=.24];
\draw[plotblue,thick] plot coordinates {(1.00000000,0.00000000) (0.50000000,0.25000000) (0.25000000,0.25000000) (0.12500000,0.18750000) (0.06250000,0.12500000) (0.03125000,0.07812500) (0.01562500,0.04687500) (0.00781250,0.02734375) (0.00390625,0.01562500) (0.00195312,0.00878906) (0.00097656,0.00488281) (0.00048828,0.00268555) (0.00024414,0.00146484) (0.00012207,0.00079346) (0.00006104,0.00042725) (0.00003052,0.00022888) (0.00001526,0.00012207) (0.00000763,0.00006485) (0.00000381,0.00003433) (0.00000191,0.00001812) (0.00000095,0.00000954)};
\fill[plotblue] (1.00000000,0.00000000) circle[radius=.65pt];
\fill[plotblue] (0.50000000,0.25000000) circle[radius=.65pt];
\fill[plotblue] (0.25000000,0.25000000) circle[radius=.65pt];
\fill[plotblue] (0.12500000,0.18750000) circle[radius=.65pt];
\fill[plotblue] (0.06250000,0.12500000) circle[radius=.65pt];
\fill[plotblue] (0.03125000,0.07812500) circle[radius=.65pt];
\fill[plotblue] (0.01562500,0.04687500) circle[radius=.65pt];
\fill[plotblue] (0.00781250,0.02734375) circle[radius=.65pt];
\fill[plotblue] (0.00390625,0.01562500) circle[radius=.65pt];
\fill[plotblue] (0.00195312,0.00878906) circle[radius=.65pt];
\fill[plotblue] (0.00097656,0.00488281) circle[radius=.65pt];
\fill[plotblue] (0.00048828,0.00268555) circle[radius=.65pt];
\fill[plotblue] (0.00024414,0.00146484) circle[radius=.65pt];
\fill[plotblue] (0.00012207,0.00079346) circle[radius=.65pt];
\fill[plotblue] (0.00006104,0.00042725) circle[radius=.65pt];
\fill[plotblue] (0.00003052,0.00022888) circle[radius=.65pt];
\fill[plotblue] (0.00001526,0.00012207) circle[radius=.65pt];
\fill[plotblue] (0.00000763,0.00006485) circle[radius=.65pt];
\fill[plotblue] (0.00000381,0.00003433) circle[radius=.65pt];
\fill[plotblue] (0.00000191,0.00001812) circle[radius=.65pt];
\fill[plotblue] (0.00000095,0.00000954) circle[radius=.65pt];
\node[below] at (1,0) {$\delta x^0$};
\node[right] at (0.50000000,0.25000000) {$\delta x^1$};
\node at (.44,-.52) {$\|\delta x^{20}\|_2\simeq 9.6\cdot10^{-6}$};
\end{scope}
\end{tikzpicture}\par
\medskip
\caption{Effect of $b$ on a uniform-path mode at $\theta=45^\circ$, for $b=0.3$ (left) and $b=1$ (right).
Here $R=\diag(1,(1+b)^{-1})$ and $\delta z^0=(1,0)$.
Read down each column: proximal displacement, projection of the reflected point and DR correction, then the proximal-output errors $\delta x^\ell=RT^\ell\delta z^0$, $0\le\ell\le20$.
The axes and colors follow \cref{fig:drgeometry}; all panels use the same coordinate scale.
Dots show every iterate, and the final norms are printed below the trajectories.
\Cref{fig:drbgeometry-large} shows the same calculation for $b=3,10$.}
\label{fig:drbgeometry}
\end{figure}
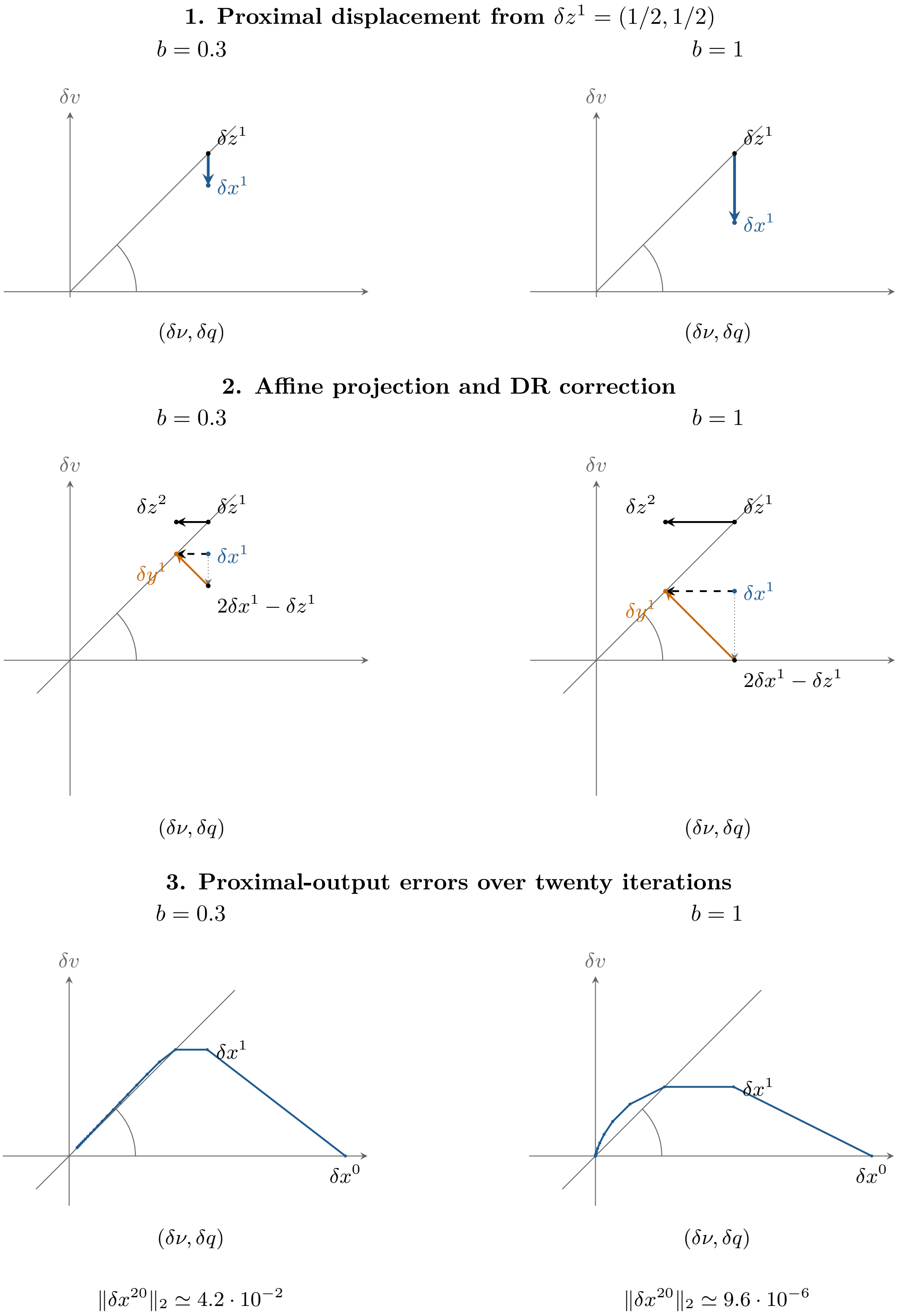

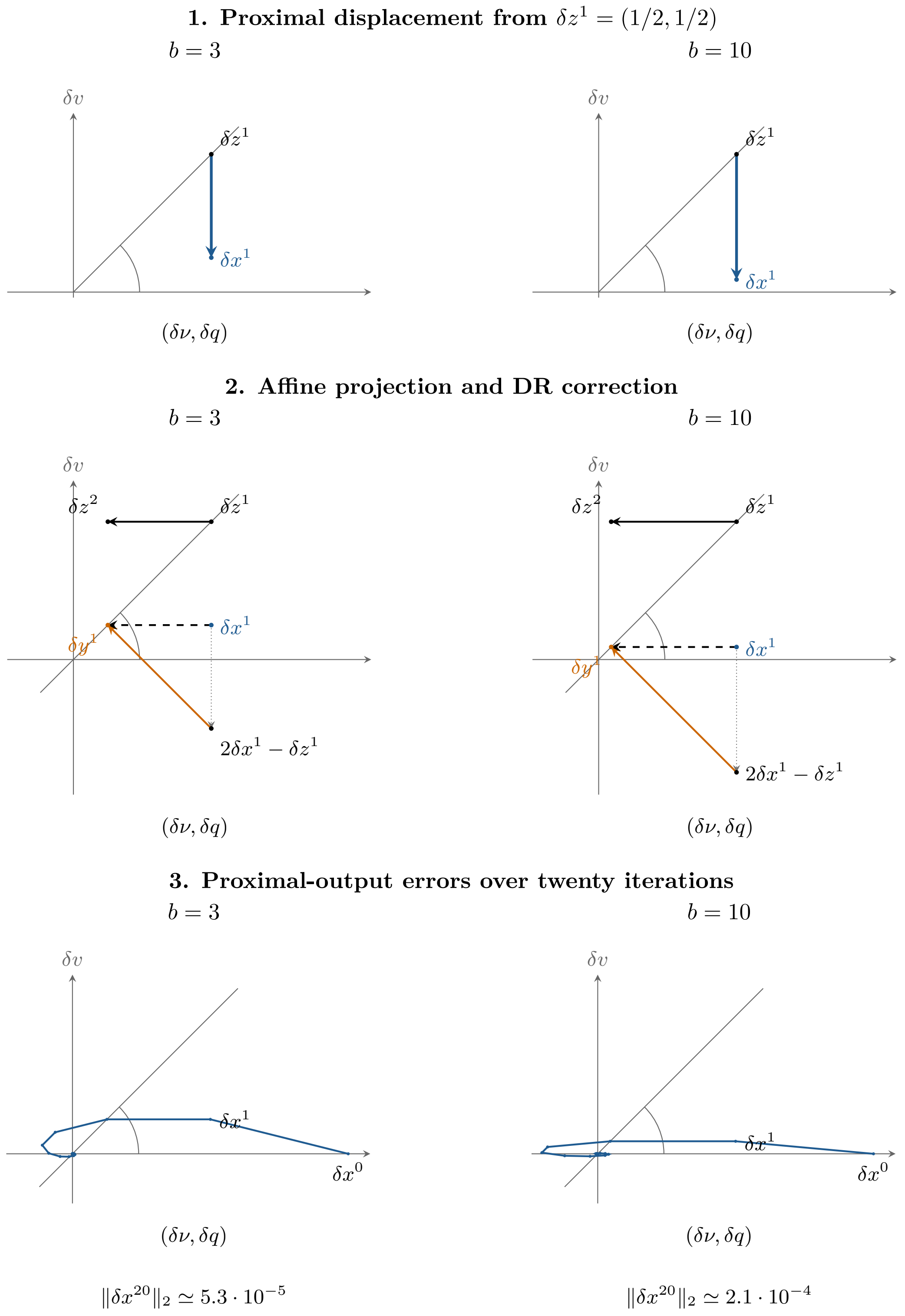
\begin{figure}[p]
\centering
{\small\bfseries 1. Proximal displacement from $\delta z^1=(1/2,1/2)$\par}
\begin{tikzpicture}[font=\footnotesize,>=stealth]
\begin{scope}[xshift=0cm,x=4.2cm,y=4.2cm]
\node[font=\small] at (.44,.88) {$b=3$};
\draw[->,black!60] (-.24,0)--(1.08,0);
\draw[->,black!60] (0,-0.02)--(0,.65) node[above] {$\delta v$};
\node at (.44,-0.15) {$(\delta\nu,\delta q)$};
\draw[black!60] (0,0)--(.60,.60);
\draw[black!60] (.24,0) arc[start angle=0,end angle=45,radius=.24];
\draw[plotblue,very thick,->] (0.50000000,0.50000000)--(0.50000000,0.12500000);
\fill[black] (0.50000000,0.50000000) circle[radius=1pt] node[above right] {$\delta z^1$};
\fill[plotblue] (0.50000000,0.12500000) circle[radius=1pt] node[right] {$\delta x^1$};
\end{scope}
\begin{scope}[xshift=8cm,x=4.2cm,y=4.2cm]
\node[font=\small] at (.44,.88) {$b=10$};
\draw[->,black!60] (-.24,0)--(1.08,0);
\draw[->,black!60] (0,-0.02)--(0,.65) node[above] {$\delta v$};
\node at (.44,-0.15) {$(\delta\nu,\delta q)$};
\draw[black!60] (0,0)--(.60,.60);
\draw[black!60] (.24,0) arc[start angle=0,end angle=45,radius=.24];
\draw[plotblue,very thick,->] (0.50000000,0.50000000)--(0.50000000,0.04545455);
\fill[black] (0.50000000,0.50000000) circle[radius=1pt] node[above right] {$\delta z^1$};
\fill[plotblue] (0.50000000,0.04545455) circle[radius=1pt] node[right] {$\delta x^1$};
\end{scope}
\end{tikzpicture}\par
\medskip
{\small\bfseries 2. Affine projection and DR correction\par}
\begin{tikzpicture}[font=\footnotesize,>=stealth]
\begin{scope}[xshift=0cm,x=4.2cm,y=4.2cm]
\node[font=\small] at (.44,.88) {$b=3$};
\draw[->,black!60] (-.24,0)--(1.08,0);
\draw[->,black!60] (0,-0.49)--(0,.65) node[above] {$\delta v$};
\node at (.44,-0.61) {$(\delta\nu,\delta q)$};
\draw[black!60] (-.12,-.12)--(.60,.60);
\draw[black!60] (.24,0) arc[start angle=0,end angle=45,radius=.24];
\draw[black!55,densely dotted,->] (0.50000000,0.12500000)--(0.50000000,-0.25000000);
\draw[orange!80!black,thick,->] (0.50000000,-0.25000000)--(0.12500000,0.12500000);
\draw[black,dashed,thick,->] (0.50000000,0.12500000)--(0.12500000,0.12500000);
\draw[black,thick,->] (0.50000000,0.50000000)--(0.12500000,0.50000000);
\fill[black] (0.50000000,0.50000000) circle[radius=1pt] node[above right] {$\delta z^1$};
\fill[plotblue] (0.50000000,0.12500000) circle[radius=1pt] node[right] {$\delta x^1$};
\fill[black] (0.50000000,-0.25000000) circle[radius=1pt] node[below right] {$2\delta x^1-\delta z^1$};
\fill[orange!80!black] (0.12500000,0.12500000) circle[radius=1pt] node[below left] {$\delta y^1$};
\fill[black] (0.12500000,0.50000000) circle[radius=1pt] node[above left] {$\delta z^2$};
\end{scope}
\begin{scope}[xshift=8cm,x=4.2cm,y=4.2cm]
\node[font=\small] at (.44,.88) {$b=10$};
\draw[->,black!60] (-.24,0)--(1.08,0);
\draw[->,black!60] (0,-0.49)--(0,.65) node[above] {$\delta v$};
\node at (.44,-0.61) {$(\delta\nu,\delta q)$};
\draw[black!60] (-.12,-.12)--(.60,.60);
\draw[black!60] (.24,0) arc[start angle=0,end angle=45,radius=.24];
\draw[black!55,densely dotted,->] (0.50000000,0.04545455)--(0.50000000,-0.40909091);
\draw[orange!80!black,thick,->] (0.50000000,-0.40909091)--(0.04545455,0.04545455);
\draw[black,dashed,thick,->] (0.50000000,0.04545455)--(0.04545455,0.04545455);
\draw[black,thick,->] (0.50000000,0.50000000)--(0.04545455,0.50000000);
\fill[black] (0.50000000,0.50000000) circle[radius=1pt] node[above right] {$\delta z^1$};
\fill[plotblue] (0.50000000,0.04545455) circle[radius=1pt] node[right] {$\delta x^1$};
\fill[black] (0.50000000,-0.40909091) circle[radius=1pt] node[right] {$2\delta x^1-\delta z^1$};
\fill[orange!80!black] (0.04545455,0.04545455) circle[radius=1pt] node[below left] {$\delta y^1$};
\fill[black] (0.04545455,0.50000000) circle[radius=1pt] node[above left] {$\delta z^2$};
\end{scope}
\end{tikzpicture}\par
\medskip
{\small\bfseries 3. Proximal-output errors over twenty iterations\par}
\begin{tikzpicture}[font=\footnotesize,>=stealth]
\begin{scope}[xshift=0cm,x=4.2cm,y=4.2cm]
\node[font=\small] at (.44,.88) {$b=3$};
\draw[->,black!60] (-.24,0)--(1.08,0);
\draw[->,black!60] (0,-0.18)--(0,.65) node[above] {$\delta v$};
\node at (.44,-0.3) {$(\delta\nu,\delta q)$};
\draw[black!60] (-.12,-.12)--(.60,.60);
\draw[black!60] (.24,0) arc[start angle=0,end angle=45,radius=.24];
\draw[plotblue,thick] plot coordinates {(1.00000000,0.00000000) (0.50000000,0.12500000) (0.12500000,0.12500000) (-0.06250000,0.07812500) (-0.10937500,0.03125000) (-0.08593750,0.00195312) (-0.04492188,-0.00976562) (-0.01269531,-0.01049805) (0.00415039,-0.00683594) (0.00891113,-0.00289917) (0.00735474,-0.00033569) (0.00401306,0.00075150) (0.00125504,0.00087738) (-0.00024986,0.00059557) (-0.00072050,0.00026655) (-0.00062680,0.00004321) (-0.00035661,-0.00005674) (-0.00012156,-0.00007295) (0.00001217,-0.00005167) (0.00005775,-0.00002431) (0.00005319,-0.00000494)};
\fill[plotblue] (1.00000000,0.00000000) circle[radius=.65pt];
\fill[plotblue] (0.50000000,0.12500000) circle[radius=.65pt];
\fill[plotblue] (0.12500000,0.12500000) circle[radius=.65pt];
\fill[plotblue] (-0.06250000,0.07812500) circle[radius=.65pt];
\fill[plotblue] (-0.10937500,0.03125000) circle[radius=.65pt];
\fill[plotblue] (-0.08593750,0.00195312) circle[radius=.65pt];
\fill[plotblue] (-0.04492188,-0.00976562) circle[radius=.65pt];
\fill[plotblue] (-0.01269531,-0.01049805) circle[radius=.65pt];
\fill[plotblue] (0.00415039,-0.00683594) circle[radius=.65pt];
\fill[plotblue] (0.00891113,-0.00289917) circle[radius=.65pt];
\fill[plotblue] (0.00735474,-0.00033569) circle[radius=.65pt];
\fill[plotblue] (0.00401306,0.00075150) circle[radius=.65pt];
\fill[plotblue] (0.00125504,0.00087738) circle[radius=.65pt];
\fill[plotblue] (-0.00024986,0.00059557) circle[radius=.65pt];
\fill[plotblue] (-0.00072050,0.00026655) circle[radius=.65pt];
\fill[plotblue] (-0.00062680,0.00004321) circle[radius=.65pt];
\fill[plotblue] (-0.00035661,-0.00005674) circle[radius=.65pt];
\fill[plotblue] (-0.00012156,-0.00007295) circle[radius=.65pt];
\fill[plotblue] (0.00001217,-0.00005167) circle[radius=.65pt];
\fill[plotblue] (0.00005775,-0.00002431) circle[radius=.65pt];
\fill[plotblue] (0.00005319,-0.00000494) circle[radius=.65pt];
\node[below] at (1,0) {$\delta x^0$};
\node[right] at (0.50000000,0.12500000) {$\delta x^1$};
\node at (.44,-.52) {$\|\delta x^{20}\|_2\simeq 5.3\cdot10^{-5}$};
\end{scope}
\begin{scope}[xshift=8cm,x=4.2cm,y=4.2cm]
\node[font=\small] at (.44,.88) {$b=10$};
\draw[->,black!60] (-.24,0)--(1.08,0);
\draw[->,black!60] (0,-0.18)--(0,.65) node[above] {$\delta v$};
\node at (.44,-0.3) {$(\delta\nu,\delta q)$};
\draw[black!60] (-.12,-.12)--(.60,.60);
\draw[black!60] (.24,0) arc[start angle=0,end angle=45,radius=.24];
\draw[plotblue,thick] plot coordinates {(1.00000000,0.00000000) (0.50000000,0.04545455) (0.04545455,0.04545455) (-0.18181818,0.02479339) (-0.20247934,0.00413223) (-0.11983471,-0.00713749) (-0.02779865,-0.00901578) (0.02667168,-0.00577146) (0.03930742,-0.00167338) (0.02718394,0.00095001) (0.00931692,0.00171064) (-0.00303941,0.00127882) (-0.00727438,0.00050125) (-0.00589283,-0.00008003) (-0.00258629,-0.00030787) (0.00009227,-0.00027149) (0.00126785,-0.00013155) (0.00122591,-0.00000815) (0.00064962,0.00005165) (0.00009238,0.00005535) (-0.00020290,0.00003188)};
\fill[plotblue] (1.00000000,0.00000000) circle[radius=.65pt];
\fill[plotblue] (0.50000000,0.04545455) circle[radius=.65pt];
\fill[plotblue] (0.04545455,0.04545455) circle[radius=.65pt];
\fill[plotblue] (-0.18181818,0.02479339) circle[radius=.65pt];
\fill[plotblue] (-0.20247934,0.00413223) circle[radius=.65pt];
\fill[plotblue] (-0.11983471,-0.00713749) circle[radius=.65pt];
\fill[plotblue] (-0.02779865,-0.00901578) circle[radius=.65pt];
\fill[plotblue] (0.02667168,-0.00577146) circle[radius=.65pt];
\fill[plotblue] (0.03930742,-0.00167338) circle[radius=.65pt];
\fill[plotblue] (0.02718394,0.00095001) circle[radius=.65pt];
\fill[plotblue] (0.00931692,0.00171064) circle[radius=.65pt];
\fill[plotblue] (-0.00303941,0.00127882) circle[radius=.65pt];
\fill[plotblue] (-0.00727438,0.00050125) circle[radius=.65pt];
\fill[plotblue] (-0.00589283,-0.00008003) circle[radius=.65pt];
\fill[plotblue] (-0.00258629,-0.00030787) circle[radius=.65pt];
\fill[plotblue] (0.00009227,-0.00027149) circle[radius=.65pt];
\fill[plotblue] (0.00126785,-0.00013155) circle[radius=.65pt];
\fill[plotblue] (0.00122591,-0.00000815) circle[radius=.65pt];
\fill[plotblue] (0.00064962,0.00005165) circle[radius=.65pt];
\fill[plotblue] (0.00009238,0.00005535) circle[radius=.65pt];
\fill[plotblue] (-0.00020290,0.00003188) circle[radius=.65pt];
\node[below] at (1,0) {$\delta x^0$};
\node[right] at (0.50000000,0.04545455) {$\delta x^1$};
\node at (.44,-.52) {$\|\delta x^{20}\|_2\simeq 2.1\cdot10^{-4}$};
\end{scope}
\end{tikzpicture}\par
\medskip
\caption{The uniform-path calculation of \cref{fig:drbgeometry} with $b=3$ (left) and $b=10$ (right).
The initial error, constraint angle, coordinate scales, and plotted quantities are unchanged.
For $b=3$ and $b=10$, the reflected point $2\delta x^1-\delta z^1=(1/2,(1-b)/(2(1+b)))$ has a negative vertical coordinate.
The eigenvalues of $T$ on this mode are complex, and the coordinates of $\delta x^\ell$ change sign along the trajectory.
Dots show the proximal-output errors $\delta x^\ell=RT^\ell\delta z^0$ for $0\le\ell\le20$; the printed final norms quantify the remaining error.
\Cref{fig:drbcomparison} reports the DR iterations per factor-ten residual decrease for $b=0.3,1,3,10$ on three nonuniform geodesics.}
\label{fig:drbgeometry-large}
\end{figure}

By \cref{sec:dr-choices}, at a fixed minimizing path and fixed $\beta$, $b$ changes neither $P$ nor the flat subspace $\ker(I-R)$.
\Cref{fig:drbgeometry,fig:drbgeometry-large} repeat the two steps and the twenty iterations of \cref{fig:drgeometry} for $b=0.3,1,3,10$, in the plane $\operatorname{span}(\delta x_{k,n},\delta v_{k,n})$ of a mode with $\theta_{k,n}=45^\circ$.
At $\theta=45^\circ$, $\delta z^1=(1/2,1/2)$ and $\delta x^1-\delta z^1=(0,-b/(2(1+b)))$, whose length increases with $b$.
The reflected point $2\delta x^1-\delta z^1=(1/2,(1-b)/(2(1+b)))$ has a negative vertical coordinate when $b>1$.
With $\sigma_{k,n}^2=1/2$, \eqref{eq:poly} becomes $\mu^2-\mu+b/(2(1+b))=0$.
For $b>1$, its roots are complex with modulus $\sqrt{b/(2(1+b))}$, so the coordinates of $\delta x^\ell$ change sign along the trajectory; for $b\le1$, they are real and the larger root is $(1+\sqrt{(1-b)/(1+b)})/2$.
For $b<1$ the larger real root exceeds $1/2$, and for $b>1$ the modulus $\sqrt{b/(2(1+b))}$ exceeds $1/2$; both equal $1/2$ at $b=1$.
The spectral radius on this mode is therefore smallest at $b=1$ (\cref{tab:drbmode}), so a larger proximal displacement does not give a smaller spectral radius.

\begin{table}[htbp]
\centering\small
\begin{tabular}{@{}lrrrr@{}}
\toprule
$b$ & $0.3$ & $1$ & $3$ & $10$\\
\midrule
$\spr(T_{k,n})$ & $0.867$ & $1/2$ & $0.612$ & $0.674$\\
$\|\delta x^{20}\|_2$ & $4.2\cdot10^{-2}$ & $9.6\cdot10^{-6}$ & $5.3\cdot10^{-5}$ & $2.1\cdot10^{-4}$\\
\bottomrule
\end{tabular}
\caption{The uniform-path mode of \cref{fig:drbgeometry,fig:drbgeometry-large}, with $\theta_{k,n}=45^\circ$ and $\delta z^0=(1,0)$.
The spectral radius $\spr(T_{k,n})$ of the restriction $T_{k,n}$ of $T$ to $\operatorname{span}(\delta x_{k,n},\delta v_{k,n})$ is the largest modulus of a root of \eqref{eq:poly} with $\sigma_{k,n}^2=1/2$, rounded to three digits.
The norm $\|\delta x^{20}\|_2$ of the proximal-output error after twenty linearized iterations is the value printed below the trajectories in the two figures.}
\label{tab:drbmode}
\end{table}

\FloatBarrier

\section{Conic equivalence and dual error bounds}\label{app:conic}

\subsection{Cone constraints}
The second-order cone condition implies $\nu+t\ge|\nu-t|$.
Thus $\nu\ge0$ and $t\ge0$.
The squared inequality is
\[
 (\nu+t)^2\ge(\nu-t)^2+2v^2
 \quad\Longleftrightarrow\quad 2\nu t\ge v^2.
\]
The converse implication follows from the same calculation and $\nu+t\ge0$.

For $\nu>0$, the smallest admissible $t$ is $v^2/(2\nu)$.
For $\nu=0$, feasibility requires $v=0$, and the smallest $t$ is zero.
These are exactly the values of $a(\nu,v)$.

For $\nu>0$, the exponential cone condition is $\nu e^{-u/\nu}\le z$.
The logarithm gives $u\ge\nu\log(\nu/z)$.
For $\nu=0$, closure of the exponential cone gives $u\ge0$ and $z\ge0$.
Since $z=\nu_i^{\rm old}>0$ is fixed data, only $u\ge0$ remains at $\nu=0$.
This is the bound $u\ge\nu\log(\nu/z)$ at $\nu=0$ with the convention $0\log0=0$ of \cref{sec:entropic-step}.

The objective of \eqref{eq:conic} increases with each $t_i$ and $u_i$.
Minimization over $(t,u)$ at fixed $(\nu,v)$ thus gives $t_i=a(\nu_i,v_i)$ and $u_i=\nu_i\log(\nu_i/\nu_i^{\rm old})$, with the zero-density values above.
With $m=\sqrt\tau\,v$, the equality constraint of \eqref{eq:conic} is $Lm=\nu-\nu^{\rm old}$ and $\A(\nu,m)=\tau\A(\nu,v)$; the entropy identity before \eqref{eq:conic} then gives the objective value $M(H_\tau(\nu,m)-1)$.
Hence \eqref{eq:conic} and \eqref{eq:step} have the same minimizers, with $m=\sqrt\tau\,v$.

\subsection{Square-completion identity}\label{app:gap}

\Cref{prop:gap} follows from an identity that writes $H_\tau(\nu,m)-D_\tau(\chi)$ as the sum of two nonnegative terms.

For unit-mass densities $\nu$ and $\nu'>0$, define the Kullback--Leibler divergence (relative entropy) with grid weights,
\[
 \KL_M(\nu\mid\nu')=\avg{\nu\log(\nu/\nu')},
\]
with $0\log0=0$.

Fix any feasible finite-action pair $(\nu,m)$ and any $\chi\in\R^M$.
Use $s(\chi)$ and the Gibbs density $\nu_\chi$ of \eqref{eq:gibbs}, and set $Z_\chi=\avg{e^{s(\chi)}}$, so that $\nu_\chi=e^{s(\chi)}/Z_\chi$.
Add the zero term $\ip\chi{Lm-\nu+\nu^{\rm old}}$ to $H_\tau(\nu,m)$.
Since $\avg\nu=1$, $\E(\nu)=\avg{\nu\log\nu}$.
Self-adjointness of $L$ gives $\ip\chi{Lm}=\ip{L\chi}m$, and
\[
 \begin{aligned}
 \frac1{2\tau}\avg{\frac{(m+\tau\nu L\chi)^2}{\nu}}
 &=\tau^{-1}\A(\nu,m)+\ip{L\chi}m+\frac\tau2\ip{\nu}{(L\chi)^2},\\
 \KL_M(\nu\mid\nu_\chi)&=\avg{\nu\log\nu}-\ip{s(\chi)}{\nu}+\log Z_\chi.
 \end{aligned}
\]
Adding these two equalities, and using the zero term and \eqref{eq:dual}, gives
\[
 \begin{aligned}
 H_\tau(\nu,m)-D_\tau(\chi)
 =\KL_M(\nu\mid\nu_\chi)
 +\frac1{2\tau}\avg{\frac{(m+\tau\nu L\chi)^2}{\nu}}.
 \end{aligned}
\]

At $\nu_i=m_i=0$, set the quotient in the last term to zero.
The two terms are nonnegative.
Therefore $D_\tau(\chi)\le H_\tau(\nu,m)$ for every feasible pair and every multiplier $\chi$.
Taking the infimum over feasible pairs gives $D_\tau(\chi)\le H_\tau^*$.

By \cref{prop:positive}, $\nu^*>0$, so $H_\tau$ is differentiable at $(\nu^*,m^*)$.
Since the constraints are affine, there exist multipliers $\chi^*\in\R^M$ for $Lm=\nu-\nu^{\rm old}$ and $\lambda\in\R$ for $\avg\nu=1$ such that
\[
 m^*=-\tau\nu^*L\chi^*,\qquad
 \log\nu^*=\chi^*+\tfrac\tau2(L\chi^*)^2-\lambda\one.
\]
The unit-mass constraint follows from conservation, since $\avg{Lm}=0$ and $\avg{\nu^{\rm old}}=1$; it is kept so that $\lambda$ is determined by normalization.
The condition $\avg{\nu^*}=1$ gives
$\lambda=\log\avg{\exp s(\chi^*)}$.
Thus $\nu_{\chi^*}=\nu^*$, and the two terms in the square-completion identity vanish at the minimizer.
This proves attainment and equality of the primal and dual values.
For any maximizer $\tilde\chi$ of $D_\tau$, the square-completion identity at $(\nu^*,m^*)$ gives $\KL_M(\nu^*\mid\nu_{\tilde\chi})=0$ and $m^*=-\tau\nu^*L\tilde\chi$.
The equality case of the Gibbs inequality gives $\nu_{\tilde\chi}=\nu^*$.
Since $m^*=-\tau\nu^*L\chi^*$ and $\nu^*>0$, also $L\tilde\chi=L\chi^*$, so $\tilde\chi-\chi^*\in\ker L=\R\one$.
Hence the maximizers of $D_\tau$ are the vectors $\chi^*+c\one$, $c\in\R$, and each of them gives the Gibbs density $\nu^*$.

Apply the square-completion identity with $\chi=\chi^*$.
Since $D_\tau(\chi^*)=H_\tau^*$ and $\nu_{\chi^*}=\nu^*$, every feasible finite-action pair $(\nu,m)$ satisfies
\[
 H_\tau(\nu,m)-H_\tau^*\ge\KL_M(\nu\mid\nu^*).
\]
Minimizing the left side over momenta satisfying conservation replaces $H_\tau(\nu,m)$ by $J_\tau(\nu)$.
Pinsker's inequality in \cref{app:pinsker} gives the density-error term in \eqref{eq:gap}.

\subsection{Pinsker's inequality with grid weights}\label{app:pinsker}
We prove the weighted Pinsker inequality
\[
 \KL_M(\nu\mid\nu^*)\ge\tfrac12\normone{\nu-\nu^*}^2.
\]
Set $\varpi_i=\dx\nu_i$ and $\varpi_i^*=\dx\nu_i^*$.
These are probability vectors, and $\varpi_i^*>0$ because $\nu^*>0$ (\cref{prop:positive}).
For $\varpi=\varpi^*$, equality holds.
Assume now $\varpi\ne\varpi^*$, and let $A=\{i:\varpi_i>\varpi_i^*\}$, $\varpi_A=\sum_{i\in A}\varpi_i$, and $\varpi_A^*=\sum_{i\in A}\varpi_i^*$.
Then $0<\varpi_A^*<\varpi_A\le1$.
Since $\varpi$ and $\varpi^*$ both have total mass one, the excess on $A$ equals the deficit on its complement:
\[
 \varpi_A-\varpi_A^*=\sum_{i\in A}(\varpi_i-\varpi_i^*)=\sum_{i\notin A}(\varpi_i^*-\varpi_i).
\]
Splitting the absolute-value sum over these two sets gives
\[
 2(\varpi_A-\varpi_A^*)=\sum_i|\varpi_i-\varpi_i^*|
 =\dx\sum_i|\nu_i-\nu_i^*|
 =\normone{\nu-\nu^*}.
\]

On $A$, apply convexity of $x\log x$ with weights $\varpi_i^*/\varpi_A^*$ to obtain
\[
 \sum_{i\in A}\varpi_i\log(\varpi_i/\varpi_i^*)\ge\varpi_A\log(\varpi_A/\varpi_A^*).
\]
Apply the same convexity inequality on the complement of $A$ and add to obtain
\[
 \sum_i\varpi_i\log(\varpi_i/\varpi_i^*)\ge f(\varpi_A),\qquad
 f(x)=x\log(x/\varpi_A^*)+(1-x)\log\bigl((1-x)/(1-\varpi_A^*)\bigr),
\]
where $0<x\le1$ and $f(1)$ is defined by continuity.

Then $f(\varpi_A^*)=f'(\varpi_A^*)=0$ and $f''(x)=1/x+1/(1-x)\ge4$ for $0<x<1$.
Taylor's formula with integral remainder at $\varpi_A^*$ gives $f(x)\ge2(x-\varpi_A^*)^2$ for $0<x<1$, and continuity of $f$ extends this bound to $x=1$.
Since $0<\varpi_A^*<\varpi_A\le1$, this bound applies at $x=\varpi_A$, and $2(\varpi_A-\varpi_A^*)^2=\tfrac12\normone{\nu-\nu^*}^2$ because $2(\varpi_A-\varpi_A^*)=\normone{\nu-\nu^*}$.
Hence $\KL_M(\nu\mid\nu^*)\ge\tfrac12\normone{\nu-\nu^*}^2$.

\subsection{Two-sided density bound}\label{app:dual-exponents}
For entries with $\nu_i\le M\sqrt{2\eta}$, \eqref{eq:pointwise} gives no relative accuracy.
A bound on $\|s(\chi)-s(\chi^*)\|_\infty$, with $s$ from \eqref{eq:gibbs}, gives the two-sided density bound stated after \eqref{eq:gibbs}.
Suppose $\|s(\chi)-s(\chi^*)\|_\infty\le\varepsilon_s$.
At each entry,
\[
 e^{-\varepsilon_s}e^{s_i(\chi)}\le e^{s_i(\chi^*)}\le e^{\varepsilon_s}e^{s_i(\chi)}.
\]

Summing with the positive quadrature weights gives
$e^{-\varepsilon_s}Z_\chi\le Z_{\chi^*}\le e^{\varepsilon_s}Z_\chi$.
Dividing the bound on $e^{s_i(\chi^*)}$ by the bound on $Z_{\chi^*}$ gives
\[
 e^{-2\varepsilon_s}(\nu_\chi)_i\le(\nu_{\chi^*})_i\le e^{2\varepsilon_s}(\nu_\chi)_i.
\]
Since $\nu_{\chi^*}=\nu^*$ at an optimal multiplier $\chi^*$, this is the two-sided density bound stated after \eqref{eq:gibbs}.

\section{Dual bounds, rounded masses and sign transitions}

\subsection{Geodesic dual bound}\label{app:geodual}
For arbitrary cell multipliers $\chi^0,\ldots,\chi^{N-1}\in\R^M$, set
\[
 \Theta^j=\chi^j-\chi^{j-1}
 -\frac h4\bigl[(L\chi^{j-1})^2+(L\chi^j)^2\bigr],\qquad1\le j<N.
\]
For every $\chi=(\chi^0,\ldots,\chi^{N-1})$, the function $D_G$ below satisfies $D_G(\chi)\le E_{M,N}$, the minimum value of \eqref{eq:geodesic}:
\begin{equation}\label{eq:geodual}
 \begin{split}
 D_G(\chi)={}&\ip{\chi^0-\frac h4(L\chi^0)^2}{\nu^0}
 -\ip{\chi^{N-1}+\frac h4(L\chi^{N-1})^2}{\nu^N}
 +\sum_{j=1}^{N-1}\min_i \Theta_i^j.
 \end{split}
\end{equation}

To prove this bound, retain the unit-mass density constraints and add
\[
 \sum_{j=0}^{N-1}\ip{\chi^j}{Lm^j-\nu^{j+1}+\nu^j}
\]
to the action; the result is the Lagrangian.
The added sum vanishes on every feasible path, so the infimum of the Lagrangian over all $m^j$ and all $\nu^j\in\mathcal P_M$, $1\le j<N$, is at most $E_{M,N}$.
By self-adjointness, $\ip{\chi^j}{Lm^j}=\ip{L\chi^j}{m^j}$.
For fixed $q^j$, the minimum over $m^j$ of the $m^j$-dependent terms is $-\tfrac h2\avg{q^j(L\chi^j)^2}$, attained at $m^j=-hq^jL\chi^j$.
Insert $q^j=(\nu^j+\nu^{j+1})/2$ and collect the coefficients of each density.
After minimization over the momenta, $\Theta^j$ is the coefficient of the interior density $\nu^j$ in the Lagrangian.
For each interior density, the minimum of $\ip{\Theta^j}{\nu^j}$ over $\mathcal P_M$ is $\min_i\Theta_i^j$.
The endpoint terms stay fixed and give \eqref{eq:geodual}.

An upper bound requires a feasible path of finite action, with the constraints satisfied in exact arithmetic.
For positive candidate densities, normalize each time level and calculate the momenta of least action under the conservation equations, as in \cref{sec:geovalue}.
At zero weights, the compatibility condition in \cref{sec:proximal-zeros} is also necessary.
The reported upper bounds use the positive density arrays $\bar\nu$ and momenta $\bar m$ specified in \cref{sec:geocompare}.
The relative action gaps $\gaprel(V)$ in \eqref{eq:comparison-gap} are floating-point estimates.
The IPM bounds in \cref{sec:ipmgeovalidation} use the interval procedure in \cref{app:geointerval}.

\subsection{Interval action bounds for the IPM geodesics}\label{app:geointerval}
The interval check starts from the IPM candidate $\nu_{\rm IPM}^j$ selected by \cref{alg:geo-candidates} for each IPM run, and from the equality multipliers $\psi^j$ of \cref{alg:geo-ipm}.
Every stored binary64 number is treated as an exact real input.
Set $\bar\nu^j=\nu_{\rm IPM}^j/\avg{\nu_{\rm IPM}^j}$ for every time level, including the endpoints.
The certified problem is \eqref{eq:geodesic} with the unit-mass endpoints $\bar\nu^0$ and $\bar\nu^N$, which can differ from the prescribed endpoints by the rounding of the stored arrays; this normalization gives unit mass in exact arithmetic.
All entries of $\nu_{\rm IPM}^j$ for the runs of \cref{tab:ipmvalidation} and of \cref{sec:geo-time-study} are positive, so no positive lower bound is imposed on $\bar\nu$.

For odd $M$, the operator $\C$ of \cref{sec:grid-notation}, with symbol $1/\widehat L(k)$ for $k\ne0$ and zero for $k=0$, is the circulant matrix
\[
 \C_{ii'}=-\frac2M\sum_{k=1}^{(M-1)/2}
 \frac{\cos(2\pi k(i-i')/M)}{4\pi^2k^2}.
\]
Evaluate this matrix with interval enclosures for $\pi$, cosine, and all arithmetic operations.
Apply \eqref{eq:geovalue} to $\bar\nu$ to enclose $V$.
The interval enclosures of the midpoint weights $\bar q_i^j=(\bar\nu_i^j+\bar\nu_i^{j+1})/2$ have positive lower endpoints, so every division by $\bar q_i^j$ in \eqref{eq:geovalue} is defined in interval arithmetic.

In the implementation, the equality rows of \cref{alg:geo-ipm} are written for the stored array $-w^j$, with $w^j=m^j/h$, as $-w^j+h^{-1}\C(\nu^{j+1}-\nu^j)+c^j\one=0$.
Their multiplier $\psi^j$ is therefore not directly a multiplier of $Lm^j=\nu^{j+1}-\nu^j$ in \eqref{eq:geodual}; the conversion below, $\chi^j=-\C\psi^j$ after centering, gives such a multiplier.
Center each $\psi^j$ by subtracting its spatial mean and set $\chi^j=-\C\psi^j$.
Then $L\chi^j=-\psi^j$, and \eqref{eq:geodual} gives the lower bound.
If $[\underline\Theta_i^j,\overline\Theta_i^j]$ encloses $\Theta_i^j$ for each $i$, then $[\min_i\underline\Theta_i^j,\min_i\overline\Theta_i^j]$ encloses $\min_i\Theta_i^j$ in \eqref{eq:geodual}.
Hence the term $\min_i\Theta_i^j$ of $D_G(\chi)$ in \eqref{eq:geodual} is enclosed as defined for the given $\chi$; it is not replaced by a smooth approximation, and $\chi$ is not optimized further.

The upper endpoint of the interval for $V$ and the lower endpoint of the interval for $D_G$ give $V_{\rm up}$ and $D_{\rm lo}$.
The computation uses the \texttt{mpmath.iv} interval context with 40 decimal digits.
Reported decimal upper bounds are rounded upward and lower bounds downward.
These enclosures bound the optimal action for the normalized endpoint data and exact grid operator.
The constraint defects $\delta_{\rm mass}$ and $\delta_{\rm cons}$ of the raw IPM arrays are separate floating-point calculations and are not used as certificates.

\subsection{Rounded input mass}\label{app:mass}
Suppose the supplied input $\nu^{\rm old}>0$ has the exact quadrature mass $\varsigma=\avg{\nu^{\rm old}}$.
Every feasible density of \eqref{eq:step} has mass $\varsigma$, because $\avg{Lm}=0$.
Because $\nu^{\rm old}$ is stored in binary64 arithmetic, its exact quadrature mass $\varsigma$ can differ from one by rounding error.
To obtain a dual lower bound for the supplied input, the minimization over densities must use mass $\varsigma$.
The resulting formula is
\begin{equation}\label{eq:dualmass}
 D_{\tau,\varsigma}(\chi)=\ip\chi{\nu^{\rm old}}+1-\varsigma+\varsigma\log\varsigma
                   -\varsigma\log\avg{e^{s(\chi)}}.
\end{equation}
Indeed, write $\nu=\varsigma\mu$ with $\avg\mu=1$.
Then
\[
 \E(\nu)-\ip{s(\chi)}\nu
 =1-\varsigma+\varsigma\log\varsigma+\varsigma\bigl(\avg{\mu\log\mu}-\ip{s(\chi)}\mu\bigr).
\]
Since $\min_{\mu\in\mathcal P_M}\{\avg{\mu\log\mu}-\ip{s(\chi)}\mu\}=-\log\avg{e^{s(\chi)}}$, minimization over $\mu$ gives \eqref{eq:dualmass}, which equals \eqref{eq:dual} for $\varsigma=1$.

For $\avg\mu=1$, $H_\tau(\varsigma\mu,\varsigma m)=\varsigma H_\tau(\mu,m)+1-\varsigma+\varsigma\log\varsigma$, and $L(\varsigma m)=\varsigma\mu-\nu^{\rm old}$ if and only if $Lm=\mu-\nu^{\rm old}/\varsigma$.
Thus minimizing $H_\tau$ subject to $Lm=\nu-\nu^{\rm old}$ is \eqref{eq:step} for the input $\nu^{\rm old}/\varsigma\in\mathcal P_M^+$, with the density and the momentum multiplied by $\varsigma$, and the objective multiplied by $\varsigma$ and increased by $1-\varsigma+\varsigma\log\varsigma$.
By \cref{prop:positive}, this problem has a unique minimizer $(\nu_\varsigma^*,m_\varsigma^*)$, and $\nu_\varsigma^*>0$.
For $\nu\ge0$ with $\avg\nu=\varsigma$, let $J_{\tau,\varsigma}(\nu)=\inf\{H_\tau(\nu,m):Lm=\nu-\nu^{\rm old}\}$ be the reduced objective.
For a positive candidate $\nu_{\rm cand}$, use $\nu=\varsigma\nu_{\rm cand}/\avg{\nu_{\rm cand}}$.
Compute its momentum by \eqref{eq:elimination}.
Suppose interval arithmetic gives $J_{\tau,\varsigma}(\nu)-D_{\tau,\varsigma}(\chi)\le\eta$.
By \eqref{eq:dualmass}, $D_{\tau,\varsigma}(\chi)=\varsigma\bigl[\ip\chi{\nu^{\rm old}/\varsigma}-\log\avg{e^{s(\chi)}}\bigr]+1-\varsigma+\varsigma\log\varsigma$, where the bracket is \eqref{eq:dual} for the input $\nu^{\rm old}/\varsigma$.
The bounds of \cref{app:gap} for this input therefore give $D_{\tau,\varsigma}(\chi)\le J_{\tau,\varsigma}(\nu_\varsigma^*)$ and $J_{\tau,\varsigma}(\nu)-J_{\tau,\varsigma}(\nu_\varsigma^*)\ge\varsigma\,\KL_M(\nu/\varsigma\mid\nu_\varsigma^*/\varsigma)$.
\Cref{app:pinsker} then gives
\[
 \normone{\nu-\nu_\varsigma^*}\le\sqrt{2\varsigma\eta}.
\]
Since $\nu_{\rm cand}-\nu=(1-\varsigma/\avg{\nu_{\rm cand}})\nu_{\rm cand}$ and $\nu_{\rm cand}>0$, we have
$\normone{\nu_{\rm cand}-\nu}=|1-\varsigma/\avg{\nu_{\rm cand}}|\,\avg{\nu_{\rm cand}}=|\avg{\nu_{\rm cand}}-\varsigma|$.
The triangle inequality gives
\[
 \normone{\nu_{\rm cand}-\nu_\varsigma^*}
 \le\sqrt{2\varsigma\eta}+|\avg{\nu_{\rm cand}}-\varsigma|.
\]
If $\varsigma\ne1$, a candidate of mass one is infeasible for $\nu^{\rm old}$ without the mass correction.
Moreover, $D_\tau(\chi+c\one)=D_\tau(\chi)+c(\varsigma-1)$ is then unbounded above in $c\in\R$, so $D_\tau(\chi)>J_{\tau,\varsigma}(\nu_\varsigma^*)$ for some $\chi$, whereas $D_{\tau,\varsigma}(\chi+c\one)=D_{\tau,\varsigma}(\chi)$.
For a fixed $\chi$,
\[
 D_{\tau,\varsigma}(\chi)-D_\tau(\chi)
 =1-\varsigma+\varsigma\log\varsigma-(\varsigma-1)\log\avg{e^{s(\chi)}},
\]
whose first-order term in $\varsigma-1$ is $-(\varsigma-1)\log\avg{e^{s(\chi)}}$.
The correction therefore matters whenever $\eta$ is not large compared with $|\varsigma-1|\,|\log\avg{e^{s(\chi)}}|$.

\subsection{Sign transitions and vanishing density}\label{app:vacuumcriterion}

\begin{lemma}[Vanishing density at separated sign transitions]\label{lem:vacuumcriterion}
For each spatial grid size $M$, consider a density profile $q_M$, a momentum rate $w_M$, and nodes $x_M^+,x_M^-$ with $q_M(x_M^+),q_M(x_M^-)>0$.
Set $\kappa_M=w_M/q_M$ at these nodes.
Assume $x_M^+,x_M^-\to x_0$ and that constants $c_\kappa,A_w,\alpha_w>0$, independent of $M$, satisfy
\begin{equation}\label{eq:vacuumcriterionhyp}
 \kappa_M(x_M^+)\ge c_\kappa,\qquad
 \kappa_M(x_M^-)\le-c_\kappa,\qquad
 |w_M(x_M^+)-w_M(x_M^-)|\le A_w d_{\T}(x_M^+,x_M^-)^{\alpha_w}.
\end{equation}
Then
\[
 q_M(x_M^+)+q_M(x_M^-)
 \le\frac{A_w}{c_\kappa}d_{\T}(x_M^+,x_M^-)^{\alpha_w}\longrightarrow0.
\]
If continuous interpolants of $q_M$ converge locally uniformly near $x_0$ to a continuous density $\rho$, then $\rho(x_0)=0$.
\end{lemma}

\begin{proof}
The two sign bounds in \eqref{eq:vacuumcriterionhyp} give
\[
 w_M(x_M^+)-w_M(x_M^-)
 \ge c_\kappa\bigl(q_M(x_M^+)+q_M(x_M^-)\bigr).
\]
Combine this inequality with the momentum estimate in \eqref{eq:vacuumcriterionhyp}.
Both density values tend to zero because they are nonnegative.
Local uniform convergence and continuity of $\rho$ give
$q_M(x_M^+)-\rho(x_0)\to0$, hence $\rho(x_0)=0$.
\end{proof}

For the paper's Fourier discretization, the momentum estimate follows from conservation at fixed $h$.
Assume $M$ is odd and that the two density rows $\nu^j,\nu^{j+1}$ are nonnegative and have quadrature mass one.
Let $w_M$ be the trigonometric interpolant of $m^j/h$.
With $K=(M-1)/2$ and $\Delta\nu=\nu^{j+1}-\nu^j$, conservation gives
\[
 \widehat w_{M,k}=\frac{\widehat{\Delta\nu}_k}{h\widehat L(k)}
 \quad(0<|k|\le K),\qquad \widehat L(k)=-4\pi^2k^2.
\]
Here $\widehat{\Delta\nu}_k=M^{-1}\sum_i\Delta\nu_i e^{-2\pi\mathrm i ki/M}$, so $|\widehat{\Delta\nu}_k|\le\normone{\Delta\nu}$ for every $k$; nonnegativity and unit mass of $\nu^j,\nu^{j+1}$ give $\normone{\Delta\nu}\le2$.
Parseval's identity therefore gives
\[
 \|w_M'\|_{L^2(\T)}^2
 =\sum_{0<|k|\le K}\frac{|\widehat{\Delta\nu}_k|^2}{4\pi^2h^2 k^2}
 \le\frac1{3h^2},
\]
since $\sum_{k\ne0}k^{-2}=\pi^2/3$.
The mean momentum does not affect this derivative.
Integration along a shortest periodic arc and the Cauchy--Schwarz inequality give
\[
 |w_M(x)-w_M(y)|
 \le\frac{\sqrt{d_{\T}(x,y)}}{\sqrt3\,h},
 \qquad d_{\T}(x,y)=\min_{n\in\mathbb Z}|x-y+n|.
\]
This supplies
$A_w=1/(\sqrt3 h)$ and $\alpha_w=1/2$, uniformly in $M$.
The inequality of \cref{lem:vacuumcriterion} follows at each $M$ from \eqref{eq:vacuumcriterionhyp} at that $M$, also with constants that depend on $M$.
For adjacent nodes, it consequently bounds their density sum by
$1/(c_\kappa\sqrt{3M}\,h)$ whenever the sign separation holds.
At fixed $M$, the bound $1/(c_\kappa\sqrt{3M}\,h)$ tends to infinity as $h\to0$.
Under simultaneous refinement with persistent sign separation, $A_w=1/(\sqrt3h)$ depends on $M$ through $h$, and the bound tends to zero if and only if $h\sqrt M\to\infty$.
A momentum estimate uniform in $h$ would remove this requirement.
Identifying a zero of a continuous geodesic still requires convergence of the density profiles at the selected times.

\subsection{Action minimization at fixed momentum}\label{app:fixed-momentum}
Let $w$ be a prescribed integrable momentum on $\T$, with $0<c_w=\int_\T|w|\,dx<\infty$.
For a nonnegative unit-mass density $\rho$ of finite action, $w=0$ almost everywhere on $\{\rho=0\}$.
Cauchy--Schwarz on $\{\rho>0\}$ gives
\[
 c_w^2
 =\left(\int_{\{\rho>0\}}\sqrt\rho\,\frac{|w|}{\sqrt\rho}\,dx\right)^2
 \le\left(\int_{\{\rho>0\}}\rho\,dx\right)
       \left(\int_{\{\rho>0\}}\frac{w^2}{\rho}\,dx\right)
 =2\int_\T a(\rho,w)\,dx.
\]
The density $\rho=|w|/c_w$ has unit mass and attains equality, with the convention $a(0,0)=0$.
At every point where $w\ne0$, its diffusivity is $\kappa=c_w\operatorname{sgn}(w)$.
This minimization fixes $w$ and imposes only nonnegativity and unit mass on $\rho$; it does not impose the geodesic conservation and endpoint constraints.

\subsection{Nonnegative momentum: feasibility and dual bound}\label{app:positive-dual}

For a positive density path, conservation gives $m^j=m_\perp^j+c\one$, where $m_\perp^j=\C(\nu^{j+1}-\nu^j)$.
The constraint $m^j\ge0$ is equivalent to $c\ge-\min_i(m_\perp^j)_i$.
On this half-line,
\[
 \frac{d}{dc}\A(q^j,m_\perp^j+c\one)
 =\avg{(m_\perp^j+c\one)/q^j}\ge0.
\]
Thus the minimum is attained at $c=-\min_i(m_\perp^j)_i$, proving \eqref{eq:positive-momentum}.
A positive linear interpolation supplies finite action and proves feasibility for any positive endpoints.
With $c=-\min_i(m_\perp^j)_i$, the entry of the least-action nonnegative momentum at an index attaining $\min_i(m_\perp^j)_i$ is zero, so every cell has a zero entry.

For $\mathsf s\in\R$, write $\mathsf s_-=\max\{-\mathsf s,0\}$.
For $q>0$, minimization over nonnegative momentum gives
\[
 \inf_{m\ge0}\left\{\frac{m^2}{2hq}+\mathsf s m\right\}
 =-\frac{hq}{2}(\mathsf s_-)^2.
\]
The identity also holds for the closed action at $q=0$.
Consequently, define $D_G^+(\chi)$ by replacing each $(L\chi^j)^2$ in \eqref{eq:geodual} and in $\Theta^j$ by $((L\chi^j)_-)^2$.
The proof in \cref{app:geodual} then gives
\[
 D_G^+(\chi)\le E^+_{M,N}\le V^+,
 \qquad V^+=\sum_j h^{-1}\A(q^j,m_+^j),
\]
for every multiplier $\chi$ and every positive unit-mass density path, with $m_+^j$ from \eqref{eq:positive-momentum}.

The interval evaluation follows \cref{app:geointerval}, with two changes: the constant in the primal momentum is $-\min_i(m_\perp^j)_i$, and the dual uses the squared negative part.
For intervals $[\underline y_i,\overline y_i]$, the interval $[\min_i\underline y_i,\min_i\overline y_i]$ encloses their minimum.
The map $y\mapsto y_-$ is enclosed by $[\max\{-\overline y,0\},\max\{-\underline y,0\}]$ on $[\underline y,\overline y]$.
The upper endpoint $V^+_{\rm up}$ of the enclosure of $V^+$ and the lower endpoint $D^+_{\rm lo}$ of the enclosure of $D_G^+(\chi)$ then satisfy $D^+_{\rm lo}\le E^+_{M,N}\le V^+_{\rm up}$, for the endpoint data normalized as in \cref{app:geointerval}.
Let $D_{\rm lo},V_{\rm up}$ be the interval endpoints for $E_{M,N}$ from \cref{app:geointerval}.
If $D_{\rm lo}>0$, then
\[
 \frac{D^+_{\rm lo}}{V_{\rm up}}\le\frac{E^+_{M,N}}{E_{M,N}}\le\frac{V^+_{\rm up}}{D_{\rm lo}}.
\]
These are the ratio enclosures in \cref{tab:positive-cost}.

The scripts and arrays are stored in
\path{paper/dr-ipm-evidence/nonnegative-diffusivity/}.
The scalar nonnegative proximal map was checked against an independent one-dimensional minimization, and the affine projection against conservation and midpoint identities.
The IPM runs use QDLDL, with Clarabel's feasibility, absolute-gap and relative-gap tolerances set to $10^{-11}$.
For the nonnegative-momentum runs, the signed run of the diffuse pair with $M=127$, $N=128$, and all signed runs of the convex-order and reversed pairs, equilibration is enabled, static regularization is $10^{-7}$, and iterative refinement uses absolute and relative tolerances $10^{-14}$ with at most 50 steps.
The nonnegative runs of the convex-order and reversed pairs with $M=127$, $N=128$ stop with status \texttt{AlmostSolved}; repeated without equilibration, they stop with status \texttt{Solved}.
For each solve, \cref{tab:positive-cost} uses the run with the narrower interval action enclosure.
The validation uses saved candidates and multipliers independently of the solver status.

\section{Rescaling near an isolated minimum}\label{app:inner}

Let $\rho$ be a positive classical solution of \eqref{eq:dlss}: $\rho$ is continuously differentiable in $t$ and four times continuously differentiable in $x$, so that every term of \eqref{eq:dlss} is defined and continuous.
Write $\rho^0_{\min}=\min_x\rho^0(x)>0$.
The change of variables in \eqref{eq:inner} is $\rho(t,x)=\rho^0_{\min}\mathsf u(\mathsf t,\mathsf y)$.
The chain rule gives
\[
 \partial_t\rho=\rho^0_{\min}\ell_\delta^{-4}\partial_{\mathsf t}\mathsf u,\qquad
 \partial_{xx}\log\rho=\ell_\delta^{-2}\partial_{\mathsf y\mathsf y}\log\mathsf u.
\]

The constant $\log\rho^0_{\min}$ disappears under spatial differentiation.
Hence
\[
 \partial_{xx}(\rho\,\partial_{xx}\log\rho)
 =\rho^0_{\min}\ell_\delta^{-4}\partial_{\mathsf y\mathsf y}(\mathsf u\,\partial_{\mathsf y\mathsf y}\log\mathsf u).
\]
Both sides of the DLSS equation have the same factor $\rho^0_{\min}\ell_\delta^{-4}$.
After division, $\mathsf u$ satisfies the same equation in $(\mathsf t,\mathsf y)$.
Multiplying a positive smooth solution $\rho$ by a constant $c>0$ therefore gives another solution of \eqref{eq:dlss}, with initial value $c\rho^0$.
The computed densities come from the entropic step \eqref{eq:step}, which has the same homogeneity in the mass (\cref{sec:inner}); this discrete property justifies normalizing the initial profile to unit mass before the solve and multiplying the computed densities by the initial mass afterwards.

For \eqref{eq:isolated}, $\rho^0_{\min}=\delta/Z$ and $\ell_\delta^{2r}=\delta$.
Thus
\[
 \mathsf u(0,\mathsf y)=1+\delta^{-1}\sin^{2r}(\pi\ell_\delta\mathsf y)
       =1+[\mathsf L_\delta\sin(\pi\mathsf y/\mathsf L_\delta)]^{2r}.
\]
On each bounded interval of $\mathsf y$, the sine expansion gives the limit in \eqref{eq:innerdata}.

This is convergence of initial data on bounded intervals.
It does not prove convergence of solutions as the circle length increases.
If a rescaled dip occurs at time $\mathsf t_{\rm dip}$, the exact change of variables gives $t_{\rm dip}=\delta^{2/r}\mathsf t_{\rm dip}$.
If the rescaled dip time has a finite nonzero limit, this identity gives the exponent $2/r$.

\section{Numerical conventions and source records}\label{app:records}

\subsection{Multiplier scaling and acceptance}
For the unscaled program \eqref{eq:conic}, the solver's equality multiplier is $\chi$ from \eqref{eq:dual}.
The objective factor $M$ cancels the quadrature factor $1/M$ in the Lagrangian.
After division of conservation row $i$ by $d_i=\sqrt{\nu_i^{\rm old}}$, a row multiplier $\psi_i$ corresponds to $\chi_i=\psi_i/d_i$.
Each dual calculation uses this conversion.

The recorded IPM settings use absolute, relative, and feasibility tolerances of $10^{-12}$, with a limit of 150 iterations.
Acceptance follows \cref{alg:jko-robust}.
With the floating-point gap $g=J_\tau(\nu)-D_\tau(\chi)$, a step is accepted if $g\le10^{-14}$, or if $10^{-14}<g\le10^{-12}$ and its interval upper gap bound is at most $10^{-14}$.
\Cref{tab:ipm-fallback} lists the solver variants in the order in which \cref{alg:jko-robust} tries them. \Cref{alg:jko-robust} tries the next variant after a rejected candidate, and once after an accepted candidate whose solver status is not \texttt{Solved}; if both candidates are then accepted, it keeps the one with the smaller interval upper gap bound, and otherwise the first.

Every accepted step in \cref{tab:baseline-validation,tab:flowcontrast} has an interval gap bound.
The minimum-density surveys of \cref{tab:dip,tab:isolatedsurvey,tab:inner,tab:innertime} use the floating-point acceptance of \cref{alg:jko-robust} unless an interval check is stated.
The comparison of $\nu$ with the Gibbs density $\nu_\chi$ from \eqref{eq:gibbs} checks relative agreement at the smallest entries.
Without a rigorous bound on the multiplier error, this agreement does not certify those entries.

\subsection{Specified Fourier operator}
The entropic trajectories computed by \cref{alg:flow} use the rounded symbol $\widehat L_{\rm round}(k)$, obtained by forming $2\pi k$ in binary32, converting this value to binary64, and evaluating $-(2\pi k)^2$ in binary64.
\Cref{tab:operator-rounding} reports the error of $\widehat L_{\rm round}$ relative to the exact symbol $\widehat L$ and its effect on the baseline entropic trajectory with exponent 16.
The interval certificates of the trajectories computed with $\widehat L_{\rm round}$ concern the grid problems with this symbol.
They do not include the error relative to the exact symbol $\widehat L$ or to the continuous equation.

The geodesic certificates of \cref{app:geointerval,app:positive-dual} use the exact symbol $\widehat L$.
The entropic certificates of \cref{tab:baseline-validation,tab:flowcontrast} use $\widehat L_{\rm round}$; those of \cref{tab:common-entropic} and \cref{sec:merged-refinement} use the binary64 symbol $-(2\pi k)^2$ of their archives.

\begin{table}[htbp]
\centering\small
\begin{tabular}{@{}p{.69\linewidth}r@{}}
\toprule
Check & Recorded maximum\\
\midrule
Relative symbol error at nonzero frequencies & $\RoundedSymbolError$\\
Discrete $L^1$ trajectory difference at saved times & $\RoundedTrajectoryDifference$\\
\bottomrule
\end{tabular}
\caption{Rounding check for the rounded symbol $\widehat L_{\rm round}$. The first row reports the maximum of $|\widehat L_{\rm round}(k)-\widehat L(k)|/|\widehat L(k)|$ over nonzero grid frequencies, with the exact symbol $\widehat L(k)=-4\pi^2k^2$. The second row reports the maximum, over the saved times $t$, of $\normone{\nu(t)-\tilde\nu(t)}$, where $\nu(t)$ is the density at time $t$ of the baseline entropic trajectory of \cref{sec:flowdata} with exponent 16 in \eqref{eq:flowdata}, computed with $\widehat L_{\rm round}$, and $\tilde\nu(t)$ is the density at time $t$ of the same trajectory computed with binary64 frequencies. These floating-point differences are not included in the interval certificates, which concern $\widehat L_{\rm round}$.}
\label{tab:operator-rounding}
\end{table}

\subsection{Source records}\label{app:source-records}
The paths in this appendix are relative to the root of the repository \nolinkurl{signed-diffusive-optimal-transport}: \nolinkurl{code/} holds the solvers and the saved runs, \nolinkurl{paper/} the manuscript and its evidence, and \nolinkurl{notebooks/} the two interactive notebooks.
The directory \nolinkurl{paper/dr-ipm-evidence/revision-20261026/} contains the code of \cref{alg:geo-prox,alg:geo-prox-checked,alg:jko-prox-checked} and all candidates used in \cref{tab:checked-prox,tab:spectralchecks,tab:common-entropic}.
Its \nolinkurl{results/} directory records inner residuals, failed stopping tests, complete IPM attempt logs, operator arrays, and timing environments.
The script \nolinkurl{make_revision_tables.py} of this directory generates the entries of \cref{tab:checked-prox,tab:spectralchecks,tab:common-entropic,tab:entropic-parameter-check} from these records, and those of \cref{tab:diffusivityobservation} from \nolinkurl{paper/dr-ipm-evidence/data/diffusivity_observation.json}.

The saved arrays of the other experiments are:
\begin{itemize}
 \item \nolinkurl{code/results/dr_conditioning/bsweep_DIFFUSIVE_*_b1.npz}:
 saved proximal outputs for the zero counts of \cref{tab:zeros} and for \cref{fig:geodesic}.
 \item \nolinkurl{code/results/drs_refine/geo_peaked_M*_N64.npz}:
 positive DR paths for the DR--IPM comparisons of \cref{tab:diffusivityobservation}.
 \item \nolinkurl{paper/dr-ipm-evidence/data/historical_fine_comparison.json}:
 DR actions and timings for \cref{tab:geofine}, extracted from the DR records.
 \item \nolinkurl{code/ipm_clarabel/results/jko_flows_M128.npz}:
 baseline IPM trajectories of \cref{fig:flowprofiles}.
 \item \nolinkurl{code/ipm_clarabel/results/jko_timing.json}:
 timing records for \cref{tab:timing,fig:flowcost}.
 \item \nolinkurl{code/ipm_clarabel/results/jko_dip_study.json}:
 time-refinement records for \cref{tab:dip}.
 \item \nolinkurl{code/ipm_clarabel/results/jko_vacuum_survey.json} and \nolinkurl{jko_vacuum_survey.npz} in the same directory:
 minimum histories and profiles for \cref{tab:isolatedsurvey,tab:dipscaling}.
 \item \nolinkurl{code/ipm_clarabel/results/jko_inner_problem.json}:
 records of the rescaled minimum-density study in \cref{tab:inner,tab:innertime}.
\end{itemize}

The additional geodesic records are under \nolinkurl{code/ipm_clarabel/results/}:
\begin{itemize}
 \item \nolinkurl{geodesic_references/diffusive_*_threads1.json}:
 matched geodesic timings, DR checkpoints, and final IPM values for \cref{fig:geocompare}.
 \item \nolinkurl{geodesic_runs/diffusive_*_M63_N64_qdldl1_C_tol1e-12.json}:
 the dual reference values $D_{\rm ref}$ of \eqref{eq:comparison-gap}.
 \item \nolinkurl{geodesic_scaling.json} and its source files under \nolinkurl{geodesic_runs/}:
 IPM refinement counts, timings, actions, and dual values for \cref{tab:geofine}.
\end{itemize}
The timings of \cref{tab:timing,fig:flowcost} are the five alternating runs per method selected from \nolinkurl{jko_timing.json}: the runs with one thread and the baseline data \eqref{eq:flowdata}, without the separate validation computations.
The DR conditioning arrays also supply the density profiles of \cref{fig:geoprofiles}.

The geodesic comparison driver is \nolinkurl{code/ipm_clarabel/run_geodesic_references.py}.
The IPM points of \cref{fig:geocompare} are completed solves at different tolerances, not samples from one solve.
The DR checkpoints contain objective evaluations and elapsed times.

The DR conditioning arrays \nolinkurl{bsweep_DIFFUSIVE_*_b1.npz} store the momentum array $\texttt{x\_m}=-\sqrt\beta\,m$: the DR code uses the opposite sign of $L$ (\cref{app:notebook-conventions}) and stores the scaled momentum $\sqrt\beta\,m$.
Convert by $m=-\texttt{x\_m}/\sqrt\beta$ (that is, $m=-\texttt{x\_m}/1000$ at $\beta=10^6$) before use in \eqref{eq:kappa}.
The saved density array contains the auxiliary DR density.
Its nonnegative part gives the proximal density used in the zero counts.

\Needspace{8\baselineskip}
The DR parameter sweep and independent IPM checks are stored in \nolinkurl{paper/dr-ipm-evidence/}.
The scripts \nolinkurl{run_beta_sweep.py} and \nolinkurl{validate_ipm.py} of this directory produce \nolinkurl{data/beta_sweep.json} and \nolinkurl{data/ipm_validation.json}.
The latter records source archive paths, SHA-256 hashes, the largest entrywise difference between the raw IPM density and the normalized path $\bar\nu$ of \cref{app:geointerval} (over all entries and over the endpoint entries), interval endpoints, and the constraint defects recomputed from the raw density and momentum arrays.
The checked IPM density arrays are saved in \nolinkurl{data/diffusive_*.npz}.
The script \nolinkurl{check_diffusivity.py} reconstructs the midpoint weights, momenta, and diffusivities in \cref{fig:diffusivityobservation}.
Its record \nolinkurl{data/diffusivity_observation.json} gives the source hashes, DR--IPM comparisons, and 40-digit interval enclosures of the diffusivity values at the transition nodes of \cref{fig:diffusivityobservation}, for the IPM candidates.
\nolinkurl{data/crosscheck.json} records the same-grid DR--IPM density differences.
The parameter sweep uses one thread and one selected CPU, recorded in each row, with a warm-up outside timing.
Each plotted parameter value corresponds to one run started from the linear interpolation of the endpoint densities.

\subsection{Technical terms}

\begin{center}
\small
\begin{tabular}{@{}p{.25\textwidth}p{.69\textwidth}@{}}
\toprule
Term & Meaning in this text\\
\midrule
Midpoint weight & The time average $q^j=(\nu^j+\nu^{j+1})/2$ used in the geodesic action.\\[3pt]
Integrated momentum & The cell variable $m^j$. Its momentum rate is $w^j=m^j/h$.\\[3pt]
Affine projection & $P_{\mathcal S_\beta}$ in \eqref{eq:prox-projection}; its weighted form is given in \cref{app:drglobal}.\\[3pt]
Proximal map & $\prox_{\gamma J}$ defined in \eqref{eq:prox-projection}.\\[3pt]
DR residual & The normalized residual $r_\ell$ of \cref{sec:geodata}.\\[3pt]
Proximal residual & $|g_\omega(q)|$ in \cref{app:prox-error} ($g_0$ in the geodesic case).\\[3pt]
Linear-solve residual & The norm of the residual of a linear system, specified where the system is used.\\[3pt]
Primal--dual gap & A feasible primal value minus a valid dual value for the same discrete problem.\\[3pt]
Interval certificate & An error bound from outward-rounded interval evaluations of the stated mathematical bounds.\\[3pt]
Relative minimum & $d_{M,\mathcal T}$ defined in \eqref{eq:dip}.\\[3pt]
Vacuum entry & An exactly zero density entry. Continuous vacuum requires a separate limiting statement.\\
\bottomrule
\end{tabular}
\end{center}

\section{Interactive notebooks and implementation algorithms}\label{app:notebooks}

The two notebooks call the DR and IPM solver modules and record the method and the parameters with each result.
Their defaults differ from some experimental settings in the main text; each experiment specifies the parameters used.

\subsection{Notebook interfaces and stored outputs}
The geodesic notebook \nolinkurl{notebooks/OT-GEODESIC.ipynb} and the entropic notebook \nolinkurl{notebooks/OT-JKO.ipynb} call the module \nolinkurl{code/notebook/paper_notebook_tools.py}.
Its functions \texttt{run\_geodesic} and \texttt{run\_jko} run the solvers with the defaults \texttt{GEO\_DEFAULTS} and \texttt{JKO\_DEFAULTS}, and \texttt{notebook\_app} builds the interactive panel.
The panel has four parts.
\begin{itemize}
\item \emph{Saved results.} The list contains the saved runs of the paper and the archives written by the panel.
A selection loads an archive with its recorded method and parameters, independently of the settings for a new solve; no solve starts.
\item \emph{New solve.} The user selects the solver, the grid and the data, and sets the parameters that the selected solver reads.
For a geodesic computed by DR or IPM, the parameters also select the constraint $m^j\ge0$ of \cref{sec:nonnegative-diffusivity}, and the data also include the convex-order pair of \cref{sec:nonnegative-diffusivity}.
For the selected grid, the code evaluates the $M$ inequalities of \cref{sec:nonnegative-diffusivity} in binary64 arithmetic and stops with an error if one fails by more than $10^{-12}$.
A progress bar follows the DR iterations of a geodesic and the steps of an entropic trajectory.
A Jupyter kernel interrupt stops a DR geodesic or an entropic trajectory and keeps the last iterate or the completed steps.
The text printed by the solver goes to a separate solver log.
A new result stays in memory until the user saves it.
Its archive receives a new name under \nolinkurl{code/results/interactive/}; no existing file is replaced.
\item \emph{View.} A time slider selects a computed time level.
The overview shows the density profile, the space--time density map and the minimum of the density over the nodes; the frame view shows one movie frame (\cref{app:notebook-movies}).
Logarithmic plots omit exact zeros and add no positive display floor.
\item \emph{Movie.} The panel writes the movie of the displayed result (\cref{app:notebook-movies}).
\end{itemize}

\begin{remark}[Runtime setup and batch execution]
The dependency list \nolinkurl{code/notebook/requirements.txt} includes NumPy, SciPy, PyTorch, Matplotlib, Jupyter, ipywidgets, mpmath, and Clarabel 0.11.1.
The notebooks do not install packages.
If Clarabel is absent from the Python path, the interface also checks \texttt{\~{}/.cache/pylib-clarabel}.
The setup variable \texttt{PUBLIC}, or the environment variable \texttt{ANTIDIFFUSIVE\_PUBLIC}, overrides the automatic detection of the directory \nolinkurl{code/}.
An mp4 movie requires the program \texttt{ffmpeg}; without it, movies are written as gif files.
The driver \nolinkurl{code/notebook/run_experiments.py} executes copies of the two notebooks and stores the figures of their cells without the panel; it does not edit the notebooks.
With the switches of the movie and of the solver comparison set to false, which is the stored setting, it loads saved examples and starts no solve.
\end{remark}

\subsection{Common conventions and reported quantities}\label{app:notebook-conventions}
Stored density arrays contain the density values $\nu_i$ of \cref{sec:grid-notation}; the probability mass at node $i$ is $\nu_i/M$.
The interface reports momenta and diffusivities for the paper convention $L=\partial_{xx}$.
The notebooks pass to the DR modules of \cref{app:algorithm-settings} the Fourier symbol of the inverse of $-\partial_{xx}$; with this symbol, the modules store $m_{\rm code}=-m$ (scaled: $-\sqrt\beta\,m$).
The geodesic momentum is integrated over a time cell; its diffusivity is $m/(hq)$.
For an entropic step, the diffusivity is $m/(\tau\nu)$.

The entropic notebook uses the same Fourier symbol for DR and IPM.
It uses the rounded symbol $\widehat L_{\rm round}$ of \cref{app:records} (\cref{tab:operator-rounding}): the frequencies and their multiplication by $2\pi$ are formed in binary32, then squared in binary64.
Thus the notebook comparison of DR and IPM concerns the single operator with symbol $\widehat L_{\rm round}$.
The geodesic interface constructs its symbol in \texttt{run\_newton\_benchmark\_suite.py}, where the frequencies $k$ are converted to binary64 before the multiplication by $2\pi$.
These conventions are recorded because the binary32 construction changes the symbol by a relative error of at most $\RoundedSymbolError$ (\cref{tab:operator-rounding}), which exceeds the DR and IPM tolerances $10^{-13}$ and $10^{-12}$ of the entropic defaults.

Each output contains the full density array, time levels, grid, and metadata.
The metadata records the method, parameters, elapsed time, stopping status, and accuracy checks.
New runs use one computational thread by default.
The reported elapsed time includes setup and postprocessing performed by the interface, including its dual calculations.
It is not the benchmark timing defined in \cref{sec:geocompare,sec:flowcost}.
A DR iteration is one update \eqref{eq:dr}; an IPM iteration is one interior-point step of Clarabel.
The counts of the two methods are therefore not comparable.

\subsection{Movies}\label{app:notebook-movies}
The function \texttt{write\_movie} writes the frames of a result to a new mp4 or gif file and never replaces an existing file.
Each frame has three panels: the density, its base-10 logarithm, and the signed diffusivity.
By default, all frames share their axes.
On a geodesic, frame $j$ is the time cell $j$, at $t=(j+1/2)/N$; it shows the midpoint weight $q^j$ and the diffusivity $\kappa^j$ of \cref{eq:kappa}.
In an entropic trajectory, frame $k$ is the level $t_k$; it shows the successor density $\nu^k$ and the diffusivity $m^k/(\tau_k\nu^k)$.
Frame $0$ shows the initial density and no diffusivity.
The momentum is the least-action momentum of the stored positive path, computed by \cref{alg:momentum-value}.
The diffusivity is not drawn at the nodes where the density is below $10^{-10}$.
On the saved DR geodesics, these are $0.08$ to $0.35$ per cent of the nodes; there the density lies between $\varepsilon_{\rm mix}/2$ and $10^{-10}$, where $\varepsilon_{\rm mix}=10^{-12}$ is the weight of the mixture $\bar\nu^j=(1-\varepsilon_{\rm mix})\nu_{\rm proj}^j+\varepsilon_{\rm mix}\one$ of \cref{sec:geocompare}, and the quotient depends on this weight.
By default, the diffusivity axis is $[-1.05\max|\kappa|,1.05\max|\kappa|]$, with the maximum over all frames; the options are the 98th percentile of $|\kappa|$ over all frames, a separate axis for each frame, and a fixed half-width.

\subsection{Shared routines and solver defaults}\label{app:algorithm-common}
The procedures below specify the DR and IPM computations, including all changes to solver outputs before objective evaluation.
All formulas use the paper convention for $L$ and $\C$ unless a stored code variable is explicitly identified.
With the symbol passed by the notebooks, the DR modules store $m_{\rm code}=-m$ (\cref{app:notebook-conventions}); negating that block converts their iteration to the formulas below.
Fourier transforms use the normalization in \cref{sec:grid-notation}.
\Cref{tab:notebook-defaults} lists the notebook defaults; the settings attached to each experiment take precedence.

\begin{table}[htbp]
\centering\small
\begin{tabular}{@{}lll@{}}
\toprule
Parameter & Geodesic notebook & Entropic notebook\\
\midrule
Grid & $M=63$, $N=64$ & $M=128$, 48 time steps\\
Initial data & Peaked, endpoint offset $10^{-2}$ & \eqref{eq:flowdata}, exponent 16, $\varepsilon=10^{-3}$\\
DR metric $\beta$ & $10^6(N/64)^2$ & $1/\tau$\\
DR proximal scale & $b=3$ & $\omega=\sqrt{\min\nu^{\rm old}\max\nu^{\rm old}}$, $b=\omega$\\
DR outer tolerance / cap & $10^{-10}$ / 20000 & $10^{-13}$ / 50000\\
IPM tolerances & $10^{-10}$ & $10^{-12}$\\
IPM iteration cap & 400 & 150 or 400; see \cref{tab:ipm-fallback}\\
IPM time limit per call & 550 seconds & No finite limit\\
Computational threads & 1 & 1\\
\bottomrule
\end{tabular}
\caption{Notebook defaults. Endpoint offset: the constant $10^{-2}$ added to the profile $\mathfrak f$ in \cref{sec:geodata}. DR outer tolerance and cap: $\varepsilon_{\rm DR}$ and $K$ in \cref{alg:geo-dr}, and the outer tolerance and cap $K$ in \cref{alg:jko-dr}. IPM tolerances: the feasibility, absolute-gap and relative-gap tolerances of Clarabel. IPM time limit: the limit passed to Clarabel for one call.}
\label{tab:notebook-defaults}
\end{table}
These three Clarabel tolerances do not replace the separately evaluated action or entropic gap.
The geodesic default $\beta=10^6(N/64)^2$ gives $\sqrt\beta\,h=10^3/64$ for every $N$, so that the slope $s=\sqrt\beta\,h\kappa$ of \eqref{eq:scaledkappa} takes the same value at fixed diffusivity $\kappa$.
The saved fine-grid DR geodesics use $b=1$, rather than the notebook default.

\begin{procedure}{Positive path for DR action evaluation}{alg:positive-path}
\alginput{the interior density rows $y_\nu^j$ of the projected output, fixed unit-mass endpoints, and $\varepsilon_{\rm mix}=10^{-12}$.}
\algoutput{an action-evaluation path $\bar\nu$; the proximal output is not changed.}
\algstep{For each interior row, compute the Euclidean projection of $\varpi=\dx\,y_\nu^j$ onto the probability simplex: sort its entries as $\varpi_{(1)}\ge\cdots\ge\varpi_{(M)}$ and compute
\[
 k_* = \max\left\{k:\varpi_{(k)}-\frac{\sum_{i=1}^k\varpi_{(i)}-1}{k}>0\right\},
 \qquad \lambda=\frac{\sum_{i=1}^{k_*}\varpi_{(i)}-1}{k_*}.
\]
The index $k=1$ satisfies the inequality, so $k_*\ge1$ in exact arithmetic; the code also replaces $k_*$ by $\max(k_*,1)$ to guard against rounding.}
\algstep{Set $\varpi_i^+=\max(\varpi_i-\lambda,0)$ and
$\bar\nu_i^j=((1-\varepsilon_{\rm mix})\varpi_i^++\varepsilon_{\rm mix}/M)/\dx$.
Keep $\bar\nu^0=\nu^0$ and $\bar\nu^N=\nu^N$.}
\algstep{Return $\bar\nu$ and record $\varepsilon_{\rm mix}$.
Use this path for action evaluation only; retain the unmodified proximal density for zero counts.}
\end{procedure}

\begin{procedure}{Momentum and objective for prescribed densities}{alg:momentum-value}
\alginput{a geodesic path or an entropic pair $(\nu^{\rm old},\nu)$, the specified operator, and the choice between the direct and the stabilized formula of step 2.}
\algoutput{the least-action momentum and the objective for those densities.}
\algstep{For a geodesic cell, set $m_\perp=\C(\nu^{j+1}-\nu^j)$ and $q=(\nu^{j+1}+\nu^j)/2$.
For an entropic step, set $m_\perp=\C(\nu-\nu^{\rm old})$ and $q=\nu$.
Compute $m_\perp$ by a real Fourier transform, multiplication by $\widehat L(k)^{-1}$ for $k\ne0$, zero at $k=0$, and inverse transform.}
\algstep{Require positive $q_i$ for the following division.
For the direct calculation, set
\[
 c=-\frac{\sum_i(m_\perp)_i/q_i}{\sum_i1/q_i},\qquad m=m_\perp+c\one.
\]
For the stabilized calculation, choose $i_*$ with $q_{i_*}=\min_iq_i$, set $r_i=(m_\perp)_i-(m_\perp)_{i_*}$, and return
\[
 m_i=r_i-\frac{\sum_k(q_{i_*}/q_k)\,r_k}{\sum_k q_{i_*}/q_k}.
\]
The two formulas are algebraically identical. The second avoids large reciprocal weights and reduces cancellation.
No positive cutoff replaces $q$ in either formula.}
\algstep{For geodesics, sum $\dx\sum_{j,i}(m_i^j)^2/(2h q_i^j)$.
For an entropic step, return
$\dx\sum_i[m_i^2/(2\tau\nu_i)+\nu_i\log\nu_i-\nu_i+1]$.}
\algstep{If the nonnegative-momentum constraint is imposed, replace step 2 by
$m_i=(m_\perp)_i-\min_k(m_\perp)_k$ and evaluate the same action.}
\end{procedure}
The notebook's final objective and diffusivity use the stabilized calculation.
The candidate selection of \cref{alg:geo-candidates} and the IPM backend of the entropic notebook use the direct calculation.
\Cref{alg:geo-candidates} assigns action $+\infty$ to a path with a nonpositive midpoint weight; it does not use the value $a(0,0)=0$ at entries with $q_i^j=m_i^j=0$.
This rejection is a rule for selecting action-evaluation paths, not a claim that every path with a zero weight has infinite action.

\subsection{Geodesic DR}\label{app:notebook-geodr}
Write $v=\sqrt\beta m$ and store the auxiliary iterate $z=(z_\nu,z_v,z_q)$ of \eqref{eq:dr}.
The proximal output $x=(x_\nu,x_v,x_q)$ and the projected output $y=(y_\nu,y_v,y_q)$ are the scaled iterates of \eqref{eq:dr}; the density block $y_\nu$ of $y$ is the density component of $Y^\ell=Q_\beta^{-1}y$ used in \cref{sec:geocompare}.
The stored density has $N-1$ interior rows; momentum and midpoint weight have $N$ rows.
The endpoint densities are fixed.
Set $\gamma=bh\beta/\dx$.

\begin{procedure}{Perspective proximal map by 18 Newton iterations}{alg:geo-prox}
\alginput{$b>0$ and arrays $(q_{\rm in},v_{\rm in})$.}
\algoutput{$(q,v)$.}
\algstep{Initialize $q=\max(q_{\rm in},10^{-12})+(bv_{\rm in}^2/2)^{1/3}$ entrywise; the cutoff $10^{-12}$ is an initialization value, not an output floor.}
\algstep{\textbf{for} $n=1,\ldots,18$ \textbf{do}}
\algstep[1]{Evaluate the cubic $\mathfrak c$ and its derivative,
\[
 \mathfrak c(q)=(q-q_{\rm in})(q+b)^2-bv_{\rm in}^2/2,\qquad
 \mathfrak c'(q)=(q+b)(3q+b-2q_{\rm in}),
\]
and set $q\leftarrow\max\{0,q-\mathfrak c(q)/\max(\mathfrak c'(q),10^{-300})\}$.}
\algstep{\textbf{end for}}
\algstep{Return $q$ and $v=v_{\rm in}q/(q+b)$.
There is no residual test and no failure flag; an output $q=0$ results from the clipping in step 3.}
\end{procedure}
\begin{procedure}{Residual-checked perspective proximal map}{alg:geo-prox-checked}
\alginput{$b>0$, arrays $(q_{\rm in},v_{\rm in})$, and positive entrywise residual tolerances $\varepsilon_i^{\rm res}$.}
\algoutput{$(q,v)$ and maximum residual and iteration count, or an explicit failure.}
\algstep{Set $(q_i,v_i)=(0,0)$ wherever $2bq_{{\rm in},i}+v_{{\rm in},i}^2\le0$.
For each remaining entry, initialize
\[
 q_{\rm lo}=0,\qquad q_{\rm hi}=\max\{q_{\rm in},0\}+(bv_{\rm in}^2/2)^{1/3},
 \qquad q=q_{\rm hi}.
\]
The monotone function $g_0$ from \cref{app:prox-error} is negative at the lower endpoint and nonnegative at the upper endpoint.}
\algstep{\textbf{for} $n=1,\ldots,160$ \textbf{do}}
\algstep[1]{Evaluate $g_0(q)$.
Stop updating the entries with $|g_0(q)|\le\varepsilon_i^{\rm res}$.
For each remaining entry, replace the lower endpoint by $q$ if $g_0(q)<0$, and the upper endpoint otherwise.
Set $q$ to $q-g_0(q)/g_0'(q)$, or to the midpoint of the updated bracket if this value is outside the open bracket or is nonfinite.}
\algstep[1]{\textbf{if} every entry passes its residual test, \textbf{break}.}
\algstep{\textbf{end for}}
\algstep{Return $v=v_{\rm in}q/(q+b)$ only after every entry passes the test.
Raise a failure on a nonfinite bracket or output, or on reaching the cap without meeting the tolerance.
Record $\max_i|g_0(q_i)|$, the Euclidean norm $\|(g_0(q_i))_i\|_2$, and the iteration count for every call.}
\end{procedure}
If $2bq_{\rm in}+v_{\rm in}^2>0$, the increasing function $g_0$ of \cref{app:prox-error} has one positive root.
Since $\mathfrak c(q)=(q+b)^2g_0(q)$, this root is also the only root of $\mathfrak c$ in $(-b,\infty)$.
\Cref{alg:geo-prox} applies Newton's method with clipping at zero to $\mathfrak c$, and \cref{alg:geo-prox-checked} applies bracketed Newton steps to $g_0$.

\begin{procedure}{Geodesic affine projection}{alg:geo-project}
\alginput{$(\nu_{\rm in},v_{\rm in},q_{\rm in})$ in scaled coordinates, endpoints, $\widehat L(k)$, and $\beta$.}
\algoutput{$y=P_{\mathcal S_\beta}(\nu_{\rm in},v_{\rm in},q_{\rm in})$.}
\algstep{Transform every spatial row. Set $c_k=\sqrt\beta/\widehat L(k)$ for $k\ne0$, and $c_0=0$.
Precompute, for $1\le n<N$,
\[
 \lambda_{nk}=\frac32+2|c_k|^2+
 2\left(\frac14-|c_k|^2\right)\cos(\pi n/N).
\]}
\algstep{At interior times $1\le j<N$, form
\[
 R_k^j=\widehat\nu_{{\rm in},k}^j-\overline{c_k}
 (\widehat v_{{\rm in},k}^{j}-\widehat v_{{\rm in},k}^{j-1})
 +\tfrac12(\widehat q_{{\rm in},k}^{j-1}+\widehat q_{{\rm in},k}^{j}).
\]
Add $(|c_k|^2-1/4)\widehat\nu_k^0$ to $R_k^1$ and
$(|c_k|^2-1/4)\widehat\nu_k^N$ to $R_k^{N-1}$.
When $N=2$, both additions enter the same row.}
\algstep{Apply the orthonormal type-I discrete sine transform in time, divide each coefficient by $\lambda_{nk}$, and apply the same transform again; the second application inverts the first because this transform is its own inverse.
This gives the interior Fourier density $\widehat y_\nu$.
The implementation evaluates the transform by the odd extension of length $2N$ and a fast Fourier transform.}
\algstep{Override its zero spatial mode by $\widehat y_{\nu,0}^j=\widehat\nu_0^0$.
Append the endpoint transforms and set
\[
 \widehat y_{v,k}^j=c_k(\widehat y_{\nu,k}^{j+1}-\widehat y_{\nu,k}^{j}),
 \qquad
 \widehat y_{q,k}^j=\tfrac12(\widehat y_{\nu,k}^{j+1}+\widehat y_{\nu,k}^{j}).
\]
Override $\widehat y_{v,0}^j=\widehat v_{{\rm in},0}^j$: the momentum mean is free.}
\algstep{Apply inverse spatial transforms and return the interior density, momentum, and midpoint arrays.}
\end{procedure}

\begin{procedure}{Geodesic DR iteration and reported outputs}{alg:geo-dr}
\alginput{positive unit-mass endpoints, $M,N,\beta,b$, tolerance $\varepsilon_{\rm DR}$, cap $K$, and a proximal solver: \cref{alg:geo-prox}, or \cref{alg:geo-prox-checked} with the tolerances $\varepsilon_i^{\rm res}$ of \cref{app:prox-error}.}
\algoutput{proximal arrays, a separate action-evaluation path, residual history, objective gap, and stopping status.}
\algstep{Initialize $\nu^j=(1-j/N)\nu^0+(j/N)\nu^N$.
Set $z_\nu=\nu^{1:N-1}$, $z_q^j=(\nu^j+\nu^{j+1})/2$, and
$z_v^j=\sqrt\beta m^j$, with $m^j$ from the direct formula in \cref{alg:momentum-value}.
Precompute the projection factors in \cref{alg:geo-project}.}
\algstep{\textbf{for} $\ell=1,\ldots,K$ \textbf{do}}
\algstep[1]{Compute
\[
 x_\nu=\max(z_\nu,0),\qquad
 (x_q,x_v)=\text{proximal solver}(z_q,z_v;\,b),\qquad
 y=\text{Algorithm~\ref{alg:geo-project}}(2x-z),
\]
Set $\Delta_{\rm DR}=y-x$ and $z\leftarrow z+\Delta_{\rm DR}$.
At the first iteration, set $r_{\rm ref}=\max(\|\Delta_{\rm DR}\|_2,10^{-300})$.
The normalized residual of \cref{sec:geodata} at iterate $\ell-1$ is $r_{\ell-1}=\|\Delta_{\rm DR}\|_2/r_{\rm ref}$.
The norm includes all three scaled blocks; do not multiply the momentum contribution by $\beta$ again.}
\algstep[1]{Record $r_{\ell-1}$ at iteration 1, every 100 iterations, and termination.
\textbf{if} $r_{\ell-1}\le\varepsilon_{\rm DR}$, \textbf{break}.
The notebook checks this condition at every iteration.}
\algstep{\textbf{end for}}
\algstep{Report convergence or termination at the cap.
Retain the last $x$ computed in the loop, without an additional proximal evaluation after the final update.
The display arrays are $X=Q_\beta^{-1}x$: append the endpoints to $x_\nu$ for display and exact-zero counts; store $x_q$ and $x_v/\sqrt\beta$ separately.
They satisfy the affine equations only to the attained residual accuracy.}
\algstep{Apply \cref{alg:positive-path} to $y_\nu$ and \cref{alg:momentum-value} to the resulting path.
Use this pair to report $V$ and diffusivity $m/(hq)$.}
\algstep{Since $D_G(\chi)$ is a lower bound for every $\chi$ (\cref{app:geodual}), the choice of $\chi$ affects only the tightness of the bound.
Use the multiplier estimate
\[
 \chi^j=\C\left(\frac{x_v^j-z_v^j}{b\sqrt\beta h}\right),
\]
where $z_v$ is the updated variable and $x_v$ the retained proximal output.
This is the paper form of the multiplier estimate $\psi=-(z_{m,\rm code}-x_{m,\rm code})/(b\sqrt\beta h)$ of the code, followed by $\chi=-\C\psi$.
Evaluate \eqref{eq:geodual} with $\min_i\Theta_i^j$ computed over all $M$ entries, without smoothing and without maximizing over $\chi$.
Return $V$, $D_G(\chi)$, and their signed floating-point difference.}
\end{procedure}
The DR runs of \cref{sec:geotests} use \cref{alg:geo-prox}, except the residual-checked runs of \cref{tab:checked-prox}; the geodesic notebook and the runs of \cref{tab:spectralchecks} use \cref{alg:geo-prox-checked}.
The DR runs of \cref{sec:geocompare} perform 3000 iterations without a stopping test.
They evaluate the action on the initial path, at iteration 1 and every 100 iterations, and \cref{fig:geocompare} plots $\gaprel(V)$ for the smallest action $V$ evaluated so far.
This running minimum differs from the last-iterate output of the notebook.

\subsection{Geodesic IPM}\label{app:algorithm-geoipm}
\Cref{alg:geo-ipm} solves the conic program of \cref{sec:geocone}, with the cone $\mathcal Q$, the epigraph variables $t_i^j$, and the free scalars $c^j$ defined there.
The implementation stores the array $-w^j$ instead of the momentum rate $w^j=m^j/h$.

\begin{procedure}{Geodesic conic solve}{alg:geo-ipm}
\alginput{endpoints, odd $M$, $N\ge2$, and the geodesic IPM settings of \cref{tab:notebook-defaults}.}
\algoutput{raw conic arrays, the selected action-evaluation path, a dual value, and solver diagnostics.}
\algstep{Use variables $(\nu^{1:N-1},w^{0:N-1},t^{0:N-1},c^{0:N-1})$ and minimize $\sum_{j,i}t_i^j$.
For each cell, impose
\[
 w^j-h^{-1}\C(\nu^{j+1}-\nu^j)-c^j\one=0,
 \qquad \avg{\nu^j}=1\quad(1\le j<N),
\]
and $\nu^j\ge0$ at every interior time.
For $q^j=(\nu^j+\nu^{j+1})/2$, impose
\[
 \left((q_i^j+t_i^j)/\sqrt2,
       (q_i^j-t_i^j)/\sqrt2,w_i^j\right)\in\mathcal Q.
\]
By \cref{sec:geocone}, $h\dx$ times the minimum of this program is the minimum of \eqref{eq:geodesic}.}
\algstep{Write the constraints in Clarabel's standard conic form: a linear function of the variables plus a slack vector equals a constant vector, and the slack vector lies in a product cone.
This product cone is the zero cone of the equalities, then the nonnegative orthant of the constraints $\nu^j\ge0$, then the $NM$ copies of $\mathcal Q$, in this order.
Clarabel minimizes the sum of a convex quadratic function and a linear function; take the quadratic term zero and the linear term $\sum_{j,i}t_i^j$.
Store the constraint matrix in compressed-column sparse format.
Call Clarabel with QDLDL, one thread, the three tolerances $10^{-10}$, at most 400 iterations, a time limit of 550 seconds, and equilibration disabled.
Each call starts a new solver; no preceding conic state is supplied.}
\algstep{Save the raw density, $m_{\rm raw}^j=hw^j$, cone slacks, equality multipliers, status, iterations, residuals, and timings.
Select the action-evaluation path by \cref{alg:geo-candidates}.}
\algstep{Let $\psi^j\in\R^M$ be the raw multipliers of the $M$ conservation rows of cell $j$ (the first $NM$ equality rows, ordered by cell).
Subtract the spatial mean of $\psi^j$ and set $\chi^j=-\C\psi^j$.
Evaluate $D_G(\chi)$ from \eqref{eq:geodual}.
No factor is applied to $\psi^j$.
With $m^j=hw^j$ and unit-mass densities,
\[
 \sum_j\ip{\chi^j}{Lm^j-\nu^{j+1}+\nu^j}
 =h\sum_j\ip{L\chi^j}{w^j-h^{-1}\C(\nu^{j+1}-\nu^j)-c^j\one},
\]
and the action is $h\dx\sum_{j,i}t_i^j$ at minimal $t$.
The Lagrangian of \cref{app:geodual} is therefore $h\dx$ times the Lagrangian of the conic program, whose multipliers for the rows written with the stored array $-w^j$ are $\psi^j=-L\chi^j$.}
\algstep{Recompute momentum and $V$ for the selected path with the stabilized formula in \cref{alg:momentum-value}.
Return the signed gap $V-D_G(\chi)$.
Keep Clarabel's termination statuses distinct: \texttt{Solved} (the three tolerances are met), \texttt{AlmostSolved} (only Clarabel's reduced tolerances are met), and the other statuses (iteration or time limit, infeasibility, numerical failure).
Only \texttt{Solved} sets the geodesic notebook's convergence flag; none of these statuses is an interval certificate.}
\end{procedure}

\begin{procedure}{Geodesic IPM candidate selection}{alg:geo-candidates}
\alginput{raw conic density $\nu_{\rm raw}$, slacks of the constraints $\nu^j\ge0$, and prescribed geodesic endpoint data.}
\algoutput{the selected candidate $\nu_{\rm IPM}$, which is the candidate with the smallest finite action, or the clipped path $\nu_{\rm clip}$ if every action is nonfinite; the tested weights and actions, the name of the selected candidate, and $\max_{j,i}|\nu_{{\rm IPM},i}^j-\nu_{{\rm raw},i}^j|$.}
\algstep{Clip negative interior entries of $\nu_{\rm raw}$ to zero and divide each interior row by its mean.
Denote this clipped path by $\nu_{\rm clip}$; it is a candidate, and the endpoints stay fixed throughout.
If clipping makes an interior row zero, the division by its mean is undefined, and $\nu_{\rm clip}$ and the mixtures of step 3 have nonfinite actions.}
\algstep{Construct two reference paths $B$: unit density at every interior time, and the linear interpolation of the endpoints.
Use the second reference path only if all its interior entries are positive; normalize its rows.}
\algstep{For each reference path $B$ and each $\varepsilon_{\rm mix}=10^{-k/2}$, $k=32,31,\ldots,8$, form
$\nu_{\rm mix}^j=(1-\varepsilon_{\rm mix})\nu_{\rm clip}^j+\varepsilon_{\rm mix}B^j$ and normalize each interior row again.
Evaluate each action by the direct version of \cref{alg:momentum-value}, with no density cutoff; retain the lowest-action mixture for each reference path.}
\algstep{Let $s_\nu^j$ be the slack of the constraint $\nu^j\ge0$; it equals $\nu^j$ at an exact solution.
If every interior entry of $s_\nu$ is positive, normalize each interior row of $s_\nu$, append the prescribed geodesic endpoint data, and add this path as a candidate; it is positive without clipping.}
\algstep{Evaluate the clipped path and all retained candidates, and select the candidate with the smallest finite action.
Record all tested mixture weights and actions, the name of the selected candidate, and its maximum change from the raw density.
If every action is nonfinite, return $\nu_{\rm clip}$ and report no upper bound on $V$; the gap of \cref{alg:geo-ipm} is then not defined.}
\end{procedure}
The gap $V-D_G(\chi)$ returned in step 5 of \cref{alg:geo-ipm} is evaluated in binary64 arithmetic and is therefore not an error bound; \cref{app:geointerval} encloses $V$ and $D_G(\chi)$ for the saved arrays in interval arithmetic.

\subsection{Entropic outer time loop}\label{app:algorithm-flow}
\begin{procedure}{Repeated entropic steps}{alg:flow}
\alginput{a positive initial density of unit mass, increasing times, solver choice, and an optional interval check.}
\algoutput{saved densities, diffusivities, per-step diagnostics, and a complete or partial trajectory status.}
\algstep{Construct the operator $L$ with the rounded symbol $\widehat L_{\rm round}$ of \cref{app:records}.
For the notebook defaults, use \eqref{eq:flowdata},
$t_0=0$, and $t_k=10^{-6}\cdot10^{4(k-1)/47}$ for $k=1,\ldots,48$, so that $t_1=10^{-6}$ and $t_{48}=10^{-2}$.
Store the initial density and set $\nu^{\rm old}=\nu^0$.}
\algstep{\textbf{for each} interval $[t_k,t_{k+1}]$ \textbf{do}}
\algstep[1]{Set $\tau=t_{k+1}-t_k$.
Call \cref{alg:jko-dr} for DR or \cref{alg:jko-robust} for IPM.
Do not reuse a previous DR auxiliary state or a previous conic interior-point state.}
\algstep[1]{For DR, accept the outer solve only when its recorded residual passes its tolerance, and calculate the dual value by \cref{alg:jko-dual}.
For IPM, use the acceptance flag from \cref{alg:jko-robust}.
If the interval check is selected, apply \cref{alg:jko-interval} to the selected IPM candidate and require its upper gap bound to be at most the acceptance threshold $\eta_{\rm acc}$ of \cref{alg:jko-robust}.
The notebook provides this check for IPM only.}
\algstep[1]{If an entry of $\nu^{\rm old}$ or $\nu$ is nonfinite or nonpositive, raise an error and stop the trajectory; the steps completed before this step are returned with a partial status that names the error.
Otherwise compute with the stabilized \cref{alg:momentum-value} the momentum $m$, the entropic objective $V$ of that procedure, and $\kappa=m/(\tau\nu)$.
Store $V-D_\tau(\chi)$, with $D_\tau(\chi)$ from \cref{alg:jko-dual} (DR) or step 6 of \cref{alg:jko-attempt} (IPM), without replacing negative values by zero.}
\algstep[1]{Save the candidate, time, objective values, status, iterations, acceptance flag, and elapsed time.
If the step is accepted, set $\nu^{\rm old}=\nu$ and continue.
Otherwise retain this final candidate for inspection, stop the trajectory, and report a partial result.
There is no automatic time-step reduction or switch to another method.}
\algstep{\textbf{end for}; \textbf{return} the saved trajectory and its status.}
\end{procedure}

\subsection{Entropic DR step and dual check}\label{app:algorithm-jkodr}
\begin{procedure}{Entropic proximal solve}{alg:jko-prox}
\alginput{$(q_{\rm in},v_{\rm in})$, $b,\omega>0$, and previous log-density $s$.}
\algoutput{$(q,v)$ and stored $s$.}
\algstep{This procedure computes $\prox_\phi$ of \eqref{eq:entropic-prox} entrywise by solving $g(e^s)=0$ for $s=\log q$, with $g$ from \cref{app:localprox}.
Set
\[
 s_{\rm lo}=\min\{0,(q_{\rm in}-1)/\omega\}-1,\qquad
 s_{\rm hi}=\log\bigl(\max\{1,q_{\rm in}\}+v_{\rm in}^2/(2b)\bigr).
\]
Clip $s$ to this bracket.}
\algstep{\textbf{for} $n=1,\ldots,100$ \textbf{do}}
\algstep[1]{Set $q=e^s$ and evaluate
\[
 g=q-q_{\rm in}+\omega s-\frac{bv_{\rm in}^2}{2(q+b)^2},\qquad
 \frac{dg}{ds}=q+\omega+\frac{bv_{\rm in}^2q}{(q+b)^3}.
\]}
\algstep[1]{\textbf{if} $\max_i|g_i/(dg_i/ds)|<10^{-13}$, \textbf{break}.}
\algstep[1]{Replace $s_{\rm lo}$ by $s$ where $g<0$ and $s_{\rm hi}$ by $s$ where $g>0$; $g$ is increasing in $s$.
Set $s$ to $s-g/(dg/ds)$, or to the midpoint of the updated bracket where this value lies outside the bracket.}
\algstep{\textbf{end for}}
\algstep{Return $q=e^s$ and $v=v_{\rm in}q/(q+b)$, and store $s$.
If the test of step 4 is not met within 100 iterations, the last iterate is returned without a failure flag.}
\end{procedure}

\begin{procedure}{Residual-checked entropic proximal solve}{alg:jko-prox-checked}
\alginput{the inputs of \cref{alg:jko-prox} and positive entrywise residual tolerances $\varepsilon_i^{\rm res}$.}
\algoutput{$(q,v)$, stored $s$, and inner diagnostics, or an explicit failure.}
\algstep{Set the bracket and clip $s$ as in step 1 of \cref{alg:jko-prox}.}
\algstep{\textbf{for} $n=1,\ldots,160$ \textbf{do}}
\algstep[1]{Set $q=e^s$ and evaluate $g$ and $dg/ds$ as in step 3 of \cref{alg:jko-prox}.
Stop updating the entries with $|g|\le\varepsilon_i^{\rm res}$.
For each remaining entry, replace $s_{\rm lo}$ by $s$ if $g<0$ and $s_{\rm hi}$ by $s$ otherwise.
Propose $s-g/(dg/ds)$; use the updated bracket midpoint if the proposal is outside the open bracket or is nonfinite.}
\algstep[1]{\textbf{if} every entry passes its residual test, \textbf{break}.}
\algstep{\textbf{end for}}
\algstep{After every entry passes the residual test, return $q=e^s$ and $v=v_{\rm in}q/(q+b)$.
Store $s=\log q$ and the inner diagnostics.
Raise a failure if a bracket or output is nonfinite, if a positive density underflows to zero, or if the iteration cap is reached without passing the test.}
\end{procedure}
In \cref{alg:jko-dr}, the entropic trajectories of \cref{sec:flowtests} use \cref{alg:jko-prox}, and the entropic notebook uses \cref{alg:jko-prox-checked}.
The single steps of \cref{tab:common-entropic} iterate steps 3 and 4 of \cref{alg:jko-dr} with \cref{alg:jko-prox-checked} and $\omega=c_\omega\sqrt{\min\nu^{\rm old}\max\nu^{\rm old}}$, but stop by the gap test of \cref{sec:flowcost} instead of the test of step 4.

\begin{procedure}{One entropic DR step}{alg:jko-dr}
\alginput{$\nu^{\rm old}>0$, $\tau$, the Fourier symbol supplied to the solver (\cref{app:records}), tolerance $\varepsilon_{\rm DR}$, iteration cap $K$, and a proximal solver: \cref{alg:jko-prox}, or \cref{alg:jko-prox-checked} with the tolerances $\varepsilon_i^{\rm res}$ of \cref{app:prox-error}.}
\algoutput{the normalized proximal density, outer residual history, and stopping flag.}
\algstep{Set $\beta=1/\tau$, $\omega=\sqrt{\min\nu^{\rm old}\max\nu^{\rm old}}$, and $b=\omega/(\tau\beta)$.
Write $v=\sqrt\beta\,m$, $z=(z_\nu,z_v)$, $x=(x_\nu,x_v)$, and $y=(y_\nu,y_v)$.
Initialize $z_\nu=\nu^{\rm old}$, $z_v=0$, and
$s=\log\max(\nu^{\rm old},10^{-300})$.
Precompute $c_k=\sqrt\beta/\widehat L(k)$ for $k\ne0$, $c_0=0$, and $1+|c_k|^2$.}
\algstep{\textbf{for} $\ell=1,\ldots,K$ \textbf{do}}
\algstep[1]{Compute $x$ by the proximal solver with input $z$, scales $b,\omega$, and the stored $s$.
Set $(\nu_{\rm in},v_{\rm in})=2x-z$ and transform in space.
Project by
\[
 \widehat y_{\nu,k}=\frac{\widehat\nu_{{\rm in},k}+
 \overline{c_k}(\widehat v_{{\rm in},k}+c_k\widehat\nu_k^{\rm old})}{1+|c_k|^2},
 \qquad
 \widehat y_{v,k}=c_k(\widehat y_{\nu,k}-\widehat\nu_k^{\rm old}).
\]
Override $\widehat y_{\nu,0}=\widehat\nu_0^{\rm old}$ and
$\widehat y_{v,0}=\widehat v_{{\rm in},0}$; invert the transforms.}
\algstep[1]{Set $\Delta_{\rm DR}=y-x$, $r=\|\Delta_{\rm DR}\|_2$, and $z\leftarrow z+\Delta_{\rm DR}$.
Record $r/\|z\|_2$, using the updated state.
\textbf{if} $r\le\varepsilon_{\rm DR}\|z\|_2$, \textbf{break}.
The implementation does not insert a denominator floor in this stopping ratio.}
\algstep{\textbf{end for}}
\algstep{Apply the proximal solver once more to the final updated $z$.
Return $\nu=x_\nu/\avg{x_\nu}$.
Determine the outer stopping flag from the last recorded residual; do not infer it from the extra proximal evaluation.}
\end{procedure}

\begin{procedure}{Entropic DR dual check}{alg:jko-dual}
\alginput{the normalized DR density, $\nu^{\rm old}$, $\tau$, and the common operator.}
\algoutput{a floating-point dual lower-bound evaluation; the density is unchanged.}
\algstep{Compute the starting point $\psi$ of the dual check: set $q_i=\max(\nu_i,10^{-14})$, $m_\perp=\C(\nu-\nu^{\rm old})$, and
\[
 m=m_\perp-\frac{\sum_k(m_\perp)_k/q_k}{\sum_k1/q_k}\one,\qquad
 \psi_i=\frac{m_i}{\tau q_i}.
\]
Then $\sum_im_i/q_i=0$, so $\avg\psi=0$, and $\chi=-\C\psi$ gives $L\chi=-\psi$, that is $m=-\tau q\,L\chi$.
Thus the starting point satisfies the optimality relation $m^*=-\tau\nu^*L\chi^*$ of \cref{app:gap}, with $(q,m)$ in place of $(\nu^*,m^*)$.}
\algstep{For each dual evaluation, center $\psi$ by subtracting $\avg\psi\one$, set $\chi=-\C\psi$, and evaluate the dual formula \eqref{eq:dual}.
Use log-sum-exp for the logarithm of the exponential sum.}
\algstep{Maximize $D_\tau(-\C\psi)$ over $\psi$ by L-BFGS applied to $-D_\tau$, with at most 600 iterations, memory 60, strong-Wolfe line search, gradient tolerance $10^{-16}$, and change tolerance $10^{-20}$.
Return the dual value at the final iterate.}
\end{procedure}
The cutoff $q_i=\max(\nu_i,10^{-14})$ enters only the starting point of \cref{alg:jko-dual}; it does not enter the proximal objectives of \cref{alg:jko-prox,alg:jko-prox-checked} or the objective $V$ of \cref{alg:flow}.
The dual check is an accuracy check, not a stopping test for DR, and does not change the accepted density.

\subsection{Entropic IPM: assembly, fallback and validation}\label{app:algorithm-jkoipm}
Let $\sigma>0$ scale the objective and $d_i>0$ scale the variables.
Use $d_i=1$ unless a variant specifies $d_i=\sqrt{\nu_i^{\rm old}}$.
The relative form of the entropy writes $\nu\log\nu-\nu$ as $\nu\log(\nu/\nu^{\rm old})+(\log\nu^{\rm old}-1)\nu$ and sets $z_i=\nu_i^{\rm old}$; the absolute form keeps $\nu\log\nu-\nu$ and sets $z_i=1$.

\Needspace{22\baselineskip}
\begin{procedure}{One entropic conic attempt}{alg:jko-attempt}
\alginput{$\nu^{\rm old}$, $\tau$, operator, and one variant from \cref{tab:ipm-fallback}.}
\algoutput{raw and normalized densities, multiplier, floating-point gap, and solver diagnostics.}
\algstep{Introduce $r_i=\nu_i/d_i$, $\mu_i=v_i/d_i$ with $v=m/\sqrt\tau$ as in \eqref{eq:conic}, action epigraph $t_i$, and entropy epigraph $u_i$.
Minimize
\[
 \sigma\dx\sum_i d_i\bigl[t_i+u_i+(\log z_i-1)r_i\bigr].
\]
The physical objective is the conic objective divided by $\sigma$, plus 1.
For each node impose
\[
 (r_i+t_i,r_i-t_i,\sqrt2\mu_i)\in\mathcal Q,
 \qquad (-u_i,r_i,z_i/d_i)\in K_{\exp}.
\]}
\algstep{In the $L$ form, impose
$\sqrt\tau L\diag(d)\mu-\diag(d)r=-\nu^{\rm old}$.
If row scaling is enabled, divide row $i$ and its right-hand side by $d_i$.
In the $\C$ form, use $d_i=1$, introduce one free scalar $c$, and impose
\[
 \sqrt\tau\mu-\C r-\sqrt\tau c\one=-\C\nu^{\rm old},
 \qquad \avg r=\avg{\nu^{\rm old}}.
\]}
\algstep{Assemble the equality, action-cone, and entropy-cone rows in that order.
Call a fresh Clarabel solver with zero quadratic objective, the chosen linear solver, tolerances, and cap.
Keep all raw arrays and the status, even when the status is not \texttt{Solved}.}
\algstep{Undo the variable scaling: $\nu_{\rm raw}=d\,r$, $m_{\rm raw}=\sqrt\tau d\,\mu$.
For equality multipliers $\psi$, recover
\[
 \chi=\psi/(\sigma\dx)\quad(L\text{ form}),\qquad
 \chi=\C\psi_{1:M}/(\sigma\dx)\quad(\C\text{ form}).
\]
For the row-scaled $L$ form, divide $\chi_i$ by $d_i$ as well.}
\algstep{Set $\nu_{{\rm cand},i}=\max(\nu_{{\rm raw},i},10^{-300})$ and
$\nu=\nu_{\rm cand}/\avg{\nu_{\rm cand}}$.
Recompute $V$ by the direct \cref{alg:momentum-value}, rather than using the raw conic momentum or epigraph value.}
\algstep{Evaluate $D_\tau(\chi)$ and $D_\tau(-\chi)$ by \eqref{eq:dual}.
If $D_\tau(-\chi)>D_\tau(\chi)$, replace $\chi$ by $-\chi$ and record the reversal; by \eqref{eq:gap}, $D_\tau(\chi)\le H_\tau^*$ for every $\chi$, so either sign gives a valid lower bound.
Return the signed gap $g=V-D_\tau(\chi)$ computed with the retained $\chi$, together with the conic residuals, the scaling parameters $\sigma$ and $d$, the raw mass defect $\avg{\nu_{\rm raw}}-\avg{\nu^{\rm old}}$, and the raw minimum $\min_i\nu_{{\rm raw},i}$.}
\end{procedure}
Step 5 normalizes $\nu_{\rm cand}$ to unit mass in floating-point arithmetic, whereas \cref{alg:jko-interval} normalizes it to the exact mass $\varsigma=\avg{\nu^{\rm old}}$ of the stored $\nu^{\rm old}$.
The two normalized densities differ by $|\varsigma-1|$ in $\normone{\cdot}$.
\Cref{app:mass} keeps the mass $\varsigma$ in the dual value \eqref{eq:dualmass} and in the error bound; this is necessary when the gap bound $\eta$ of \cref{alg:jko-interval} is of the order of $|\varsigma-1|$.

\begin{table}[htbp]
\centering\small
\begin{tabular}{@{}cllllr@{}}
\toprule
Order & Constraint & Entropy & $\sigma$ & Linear solver / variable scale $d$ & Cap\\
\midrule
1 & $L$ & relative & $M$ & QDLDL / $\sqrt{\nu^{\rm old}}$ & 150\\
2 & $L$ & relative & $M$ & QDLDL / none & 150\\
3 & $L$ & relative & $1$ & QDLDL / none & 400\\
4 & $L$ & relative & $M$ & faer / $\sqrt{\nu^{\rm old}}$ & 400\\
5 & $\C$ & relative & $1$ & faer / none & 400\\
6 & $\C$ & relative & $1$ & QDLDL / none & 400\\
7 & $\C$ & absolute & $M$ & QDLDL / none & 400\\
\bottomrule
\end{tabular}
\caption{Ordered IPM fallback variants for one entropic step. Constraint: $L$ or $\C$ form of step 2 of \cref{alg:jko-attempt}; $\sigma$: objective scale; variable scale $d$: none means $d=1$; Cap: iteration limit. QDLDL and faer are direct sparse linear-system solvers selectable in Clarabel. All variants use momentum scale $\sqrt\tau$, three solver tolerances $10^{-12}$, and one thread. Variants 1 and 4 also divide conservation row $i$ by $\sqrt{\nu_i^{\rm old}}$. Relative and absolute refer to the entropy encodings in \cref{alg:jko-attempt}; their physical objectives agree. Each attempt starts a new solver. The variant order and shorter first two caps are empirical choices.}
\label{tab:ipm-fallback}
\end{table}

\begin{procedure}{Entropic IPM fallback and acceptance}{alg:jko-robust}
\alginput{one entropic step, acceptance threshold $\eta_{\rm acc}=10^{-14}$, and the bound $\eta_{\rm win}=10^{-12}$; gaps in $(\eta_{\rm acc},\eta_{\rm win}]$ are re-evaluated by \cref{alg:jko-interval}.}
\algoutput{a selected attempt, its acceptance flag, and the complete attempt log.}
\algstep{Initialize empty records $\mathsf{best}$ and $\mathsf{first}$, an empty attempt log, and zero accumulated times.
The record $\mathsf{first}$ stores an accepted candidate awaiting one additional attempt.}
\algstep{\textbf{for} variant $i=1,\ldots,7$ in \cref{tab:ipm-fallback} \textbf{do}}
\algstep[1]{Call \cref{alg:jko-attempt} and denote its returned record by $\mathsf{rec}$.
Build the dense $L$ matrix at its first use and reuse it for subsequent $L$ variants.}
\algstep[1]{Set $g$ to the signed floating-point gap, replacing a nonfinite value by $+\infty$.
\textbf{if} $\eta_{\rm acc}<g\le\eta_{\rm win}$, set the decision value $\eta(\mathsf{rec})$ to the interval upper gap bound from \cref{alg:jko-interval}; \textbf{else} set $\eta(\mathsf{rec})=g$.
Mark $\mathsf{rec}$ accepted exactly when $\eta(\mathsf{rec})\le\eta_{\rm acc}$.}
\algstep[1]{Append the variant, status, iteration count, gap, decision value, and acceptance flag to the log.}
\algstep[1]{\textbf{if} $\mathsf{first}$ is nonempty \textbf{then}}
\algstep[2]{Set $\mathsf{best}=\mathsf{first}$.
\textbf{if} $\mathsf{rec}$ is accepted, evaluate or reuse interval upper gap bounds for both candidates.
Choose $\mathsf{rec}$ only if its bound is strictly smaller; otherwise keep $\mathsf{first}$.
Store both bounds in the log and recompute the selected candidate's acceptance from its bound.}
\algstep[2]{\textbf{break}, including when the interval comparison rejects both candidates.
Do not restart the fallback chain.}
\algstep[1]{\textbf{end if}}
\algstep[1]{\textbf{if} $\mathsf{best}$ is empty or $\eta(\mathsf{rec})<\eta(\mathsf{best})$, set $\mathsf{best}=\mathsf{rec}$.}
\algstep[1]{\textbf{if} $\mathsf{rec}$ is rejected, \textbf{continue} to the next variant.}
\algstep[1]{\textbf{if} the status is \texttt{Solved} or $i=7$, \textbf{break}.}
\algstep[1]{Set $\mathsf{first}=\mathsf{rec}$ and try the next variant once.}
\algstep{\textbf{end for}}
\algstep{\textbf{return} $\mathsf{best}$ and its acceptance flag.
If no attempt passed, this is the candidate with the smallest decision value and a failed flag.
Return every attempt record and the selected attempt's iteration count.
Report accumulated solve and wall times over all attempts, including interval checks.}
\end{procedure}
An attempt can pass the floating-point screen with a small negative gap; this is not treated as a proof of zero error.
The acceptance flag, stored in the field \texttt{certified} of the saved records, is true when $\eta(\mathsf{rec})\le\eta_{\rm acc}$.
The decision value $\eta(\mathsf{rec})$ is an interval upper bound only when it was computed by \cref{alg:jko-interval}; otherwise it is the floating-point gap.
When the interval check of \cref{alg:flow} is selected, it is applied to the attempt selected here.

\begin{procedure}{Interval validation of one entropic candidate}{alg:jko-interval}
\alginput{stored binary64 arrays $(\nu_{\rm cand},\chi,\nu^{\rm old})$, time increment $\tau$, and the Fourier symbol supplied to the solver, as specified in \cref{app:records}.}
\algoutput{outward-rounded bounds for the primal value, dual value, gap, and density $L^1$ error.}
\algstep{Treat all supplied binary64 numbers as exact data.
Use 133-bit interval arithmetic and construct interval circulant matrices for $L$ and $\C$ from the given symbol, including the single Nyquist contribution when $M$ is even.}
\algstep{Compute $\varsigma=\avg{\nu^{\rm old}}$ and
$\nu=\varsigma\nu_{\rm cand}/\avg{\nu_{\rm cand}}$ in interval arithmetic.
Require positive entries, then compute the feasible momentum and primal objective $V$ by \eqref{eq:elimination}.}
\algstep{Evaluate $D_{\tau,\varsigma}(\chi)$ by \eqref{eq:dualmass}.
For its logarithmic exponential sum, subtract a finite common shift before exponentiation and add it after taking the logarithm.
All sums, exponentials, logarithms, and divisions are interval operations.}
\algstep{Set $\eta=\sup(V-D_{\tau,\varsigma})$.
If $\eta\ge0$, form
\[
 \sqrt{2\varsigma\eta}+|\avg{\nu_{\rm cand}}-\varsigma|,
\]
an interval upper bound for $\normone{\nu_{\rm cand}-\nu_\varsigma^*}$, where $\nu_\varsigma^*$ is the successor density of the step with input mass $\varsigma$ (\cref{app:mass}); the second term is $\normone{\nu_{\rm cand}-\nu}$.
If $\eta<0$, no error bound is formed, and the stored error bound is NaN.}
\algstep{Export upper endpoints for $V$, the gap, and the error bound, and a lower endpoint for $D$.
After conversion to binary64, move each bound one representable number outward.
Save the interval width, input mass defect, mass correction, and precision.}
\end{procedure}

\subsection{Nonnegative-diffusivity experiment}\label{app:algorithm-positive}
The driver for \cref{sec:nonnegative-diffusivity} is separate from the notebooks.
\Cref{alg:positive-experiment} lists how this driver differs from the notebook procedures \cref{alg:geo-prox-checked,alg:geo-dr,alg:geo-ipm,alg:geo-candidates} in the sign convention, the proximal map, the candidate selection, and the IPM settings.

\begin{procedure}{Signed and nonnegative geodesic comparison}{alg:positive-experiment}
\alginput{diffuse, convex-order or reversed endpoints, $M,N$, a flag imposing $m\ge0$, and a flag disabling equilibration.}
\algoutput{solver arrays, action-evaluation arrays, diffusivity, and separate interval action bounds.}
\algstep{For DR, pass the symbol $1/\widehat L(k)=-1/(4\pi^2k^2)$ of $\C$ for $k\ne0$ (zero at $k=0$) to the geodesic DR module, which then stores $v=\sqrt\beta\,m$, not the $-\sqrt\beta\,m$ of the notebooks (\cref{app:notebook-conventions}).
Initialize the linear density path and its least-action momentum, using the nonnegative branch of \cref{alg:momentum-value} when requested.}
\algstep{Before each perspective proximal call, replace $v_{\rm in}$ by $\max(v_{\rm in},0)$ only in the nonnegative case.
Evaluate the perspective proximal map by \cref{alg:geo-prox}, then set both returned entries to zero where $2bq_{\rm in}+v_{\rm in}^2\le0$, in the signed and in the nonnegative runs.
Keep the affine projection unchanged.}
\algstep{Compute the DR residual $r_\ell$ of \cref{alg:geo-dr} (relative to the initial residual) at iteration 1 and every 100 iterations, and test stopping only at these iterations.
Record action diagnostics at iteration 1, every 1000 iterations, and a successful stopping check.
Form the action-evaluation path by simplex projection and a $10^{-12}$ uniform mixture as in \cref{alg:positive-path}.
Save the unmodified proximal path, momentum, and midpoint arrays separately.}
\algstep{For IPM, use \cref{alg:geo-ipm}'s conic assembly.
In the nonnegative case append $w\ge0$.
The baseline settings are QDLDL, one thread, tolerance $10^{-11}$, cap 300, and a 180-second limit.
The modified settings are the baseline settings with equilibration enabled, static regularization $10^{-7}$, and both iterative-refinement tolerances $10^{-14}$ with refinement cap 50.
For the diffuse pair, the signed run at $(M,N)=(63,64)$ uses the baseline settings; the signed run at $(M,N)=(127,128)$ and the nonnegative runs at both grids use the modified settings.
For the convex-order and reversed pairs, all runs use the modified settings, and the nonnegative runs at $(127,128)$ are repeated with equilibration disabled.}
\algstep{Clip and normalize the raw IPM interior density.
Test uniform mixtures with weights $0,10^{-14},10^{-12},10^{-10},10^{-8}$, retain those with positive midpoint weights, and select the smallest recomputed action.
This driver does not use the slack candidates or linear-base mixtures of \cref{alg:geo-candidates}.
For both solvers, compute $m=m_\perp-\min_i(m_\perp)_i\one$ with $m_\perp=\C(\nu^{j+1}-\nu^j)$ in the nonnegative case and return $\kappa=m/(hq)$.}
\algstep{For the dual, center the conic seed $\psi$, set $\chi=-\C\psi$, and use \cref{app:positive-dual} with $(-L\chi)_+^2=\psi_+^2$ in the nonnegative case.
Use \eqref{eq:geodual} in the signed case.}
\algstep{Independently normalize all saved density rows, including endpoints, to exact unit mass in 40-decimal interval arithmetic.
Use the exact symbol $-4\pi^2k^2$ and enclose the action and dual as in \cref{app:geointerval,app:positive-dual}.
For an interval minimum use $[\min_i\underline y_i,\min_i\overline y_i]$ when each entry lies in $[\underline y_i,\overline y_i]$; apply the positive-part operation to both interval endpoints.
Export all bounds with outward rounding.}
\end{procedure}
The DR parameters and stopping thresholds of the reported runs are in \cref{tab:positive-agreement}.
Neither the primal--dual gap nor the DR--IPM density difference bounds the error of $\kappa=m/(hq)$: an error $\delta m$ in $m$ changes $\kappa$ by $\delta m/(hq)$, which has no bound as $q\to0$.

\subsection{Conic solver settings and implementation records}\label{app:algorithm-settings}
Clarabel 0.11.1 supplies the interior-point solver.
The following Clarabel settings are active unless a procedure overrides them.

\begin{center}\small
\begin{tabular}{@{}ll@{}}
\toprule
Setting & Value\\
\midrule
Static regularization & Enabled; constant $10^{-8}$, proportional $4.93\cdot10^{-32}$\\
Dynamic regularization & Enabled; threshold $10^{-13}$, shift $2\cdot10^{-7}$\\
Iterative refinement & Enabled; relative $10^{-13}$, absolute $10^{-12}$\\
Refinement cap / stop ratio & 10 / 5\\
Maximum step fraction & $0.99$\\
Entropic-step equilibration & Enabled; 10 iterations, scales in $[10^{-4},10^4]$\\
Baseline geodesic equilibration & Disabled\\
\bottomrule
\end{tabular}
\end{center}
Regularization changes the linear systems used for search directions; it does not add a term to the density objective.
Iterative refinement corrects those linear solves, and equilibration rescales the conic problem.
Raw solver status, residuals, and objective values remain distinct from the independently evaluated primal--dual gap.

The notebook dispatch and final evaluations are in \nolinkurl{code/notebook/paper_notebook_tools.py}.
The DR modules are \nolinkurl{proximal_dr_baseline.py}, \nolinkurl{proximal_dr_scaled.py}, and \nolinkurl{jko_dr_scaled.py} in that directory.
Conic assembly, candidate selection, fallback, and interval validation are in
\nolinkurl{code/ipm_clarabel/geodesic_clarabel.py} and \nolinkurl{code/ipm_clarabel/jko_clarabel.py}.
The nonnegative comparison is in \nolinkurl{paper/dr-ipm-evidence/nonnegative-diffusivity/experiment.py} and \nolinkurl{validate.py} in the same directory; \nolinkurl{convex_order.py} in that directory checks the convex order in exact rational arithmetic and computes the ratio enclosures of the convex-order and reversed pairs.
Each saved run records its parameters, the selected candidate and the modifications applied to it, all IPM attempts, and the stopping status.
\Cref{app:notebooks} states the timing conventions of the notebook and benchmark runs.

\clearpage

\section{Proofs of Eulerian geodesic convergence}\label{app:merged-geodesics}
This appendix proves \cref{geo:thm:time-discrete-convergence,geo:thm:joint-grid,geo:cor:grid-endpoints}.
The liminf inequality of the $\Gamma$-convergence uses weak compactness and the lower semicontinuity of $\Ac$ (\cref{geo:lem:space-time-action}).
The recovery sequence uses time reparametrization and averaging
(\cref{geo:lem:time-average-recovery}). It preserves the endpoints and
any nonnegativity constraint on momentum without spatial smoothing.
The same construction is used in the real-line positive-momentum
corollary, whose additional tightness checks are in
\cref{mart:app:positive-midpoint}.
The heat construction in the proof of
\cref{geo:prop:heat-semigroup-hypotheses} is used separately to prove
finite signed action for arbitrary periodic endpoints.
The final subsections treat the joint Fourier-grid limit and the
change to prescribed positive endpoint approximations.

\subsection{Recovery by time averaging}
\label{geo:sec:time-average-recovery}
A recovery sequence for an admissible continuous pair $(\mu,m)$ is
an admissible discrete approximation with the same endpoints, whose
marginals converge uniformly weakly to $\mu$, whose space--time
momenta converge weakly to $\dd t\otimes m_t$, and whose limiting
upper action is at most $\int_0^1\Ac(\mu_t,m_t)\dd t$.
Such a construction supplies the upper bound needed for
$\Gamma$-convergence. We first recall the standard quadratic duality
formula; for lower semicontinuity and disintegration, see
\cite[Sections~5.3--5.4]{AGS}. We omit the standard proof.

\begin{lemma}[Perspective action]\label{geo:lem:space-time-action}
Let \(Q\) be a Polish space, \(R\) a finite positive Borel measure and \(\mathsf W\) a finite signed Borel measure on \(Q\), and define \(\Ac(R,\mathsf W)\) on \(Q\) by \eqref{geo:eq:space-time-action}.
Let \(C_b(Q)\) be the space of bounded continuous functions on \(Q\). Then
\begin{equation}\label{geo:eq:space-time-duality}
  \Ac(R,\mathsf W)
  =
  \sup_{\varphi\in C_b(Q)}
  \left\{
    \int_Q \varphi\,\dd \mathsf W-\frac12\int_Q \varphi^2\,\dd R
  \right\}.
\end{equation}
Consequently \(\Ac\) is jointly convex and weakly lower
semicontinuous.  Let \(Q=[0,1]\times\T\) or \(Q=[0,1]\times\R\), and
let \(R=\dd t\otimes\mu_t\). If \(\Ac(R,\mathsf W)<\infty\), then
\[
  \mathsf W=\dd t\otimes m_t,
  \qquad
  \Ac(R,\mathsf W)=\int_0^1\Ac(\mu_t,m_t)\,\dd t
\]
for a weakly measurable family \(m_t\ll\mu_t\). Conversely, if
\(\mathsf W=\dd t\otimes m_t\) for a weakly measurable family \(m_t\),
then \(\Ac(R,\mathsf W)=\int_0^1\Ac(\mu_t,m_t)\,\dd t\), also when one
side is infinite.
\end{lemma}

For \(Q=\T\) or \(Q=\R\), use the action \(\Ac\) of
\eqref{geo:eq:action-density} on \(Q\). Let \(\gamma>0\), let \((\Omega,\vartheta)\) be a
finite positive measure space, and let \(\upsilon\mapsto(\mu_\upsilon,m_\upsilon)\)
be weakly measurable, with \(\mu_\upsilon\ge0\) and \(\int(\mu_\upsilon(Q)+|m_\upsilon|(Q))\,\dd\vartheta(\upsilon)<\infty\).  Then
\begin{align}
  \Ac(\gamma\mu,m)&=\gamma^{-1}\Ac(\mu,m),
  &
  \Ac(\mu,\gamma m)&=\gamma^2\Ac(\mu,m),
  \label{geo:eq:action-homogeneity}\\
  \Ac\left(\int\mu_\upsilon\,\dd\vartheta,
           \int m_\upsilon\,\dd\vartheta\right)
  &\le
  \int\Ac(\mu_\upsilon,m_\upsilon)\,\dd\vartheta .
  \label{geo:eq:action-mixture}
\end{align}
These properties follow from \eqref{geo:eq:space-time-duality} on
\(Q\): the substitutions \(\varphi\mapsto\varphi/\gamma\) and
\(\varphi\mapsto \gamma\varphi\) give \eqref{geo:eq:action-homogeneity},
and the supremum of an integral is at most
the integral of the suprema, which gives \eqref{geo:eq:action-mixture}.

\begin{lemma}[Recovery by time averaging]
\label{geo:lem:time-average-recovery}
Let $Q=\T$ with $\Lc$ as above, or $Q=\R$ with $\Lc=\partial_{xx}$,
and use the perspective action \eqref{geo:eq:action-density} on $Q$.
Let $(\mu_t,m_t)_{0\le t\le1}$ be an admissible pair of probability
measures and signed momenta on $Q$, in the sense of
\cref{geo:def:admissible-pair} with $\T$ replaced by $Q$, with finite action
$I=\int_0^1\Ac(\mu_t,m_t)\dd t$.
For $0<\varepsilon<1/8$ and $h=1/N<2\varepsilon$, there is a discrete
pair with the same endpoints and exact constraints
$\Lc m_h^k=\mu_h^{k+1}-\mu_h^k$ such that
\begin{equation}\label{geo:eq:time-average-bound}
 \sum_{k=0}^{N-1}\frac1h
 \Ac\left(\frac{\mu_h^k+\mu_h^{k+1}}2,m_h^k\right)
 \le \frac{1+2h/\varepsilon}{1-4\varepsilon}\,I.
\end{equation}
The construction preserves $m_t\ge0$ if this constraint is imposed.
On $Q=\R$, it preserves finite second moments whenever
$\sup_t\int x^2\mu_t(\dd x)<\infty$.
If $\varepsilon=\varepsilon_h\to0$ and $h/\varepsilon_h\to0$, the affine
marginal interpolations converge uniformly narrowly to $\mu_t$, and
the space--time momenta $\dd t\otimes(m_h^k/h)$ on the time cells
converge weakly to $\dd t\otimes m_t$.
On $\T$, the marginal convergence is convergence in $d_{\mathcal X}$.
\end{lemma}

\begin{proof}
\emph{Stationary endpoints.}
Define the nondecreasing piecewise affine map
\[
 \theta_\varepsilon(t)=
 \begin{cases}
  0,&0\le t\le2\varepsilon,\\
  (t-2\varepsilon)/(1-4\varepsilon),
       &2\varepsilon<t<1-2\varepsilon,\\
  1,&1-2\varepsilon\le t\le1.
 \end{cases}
\]
Set $\widetilde\mu_t=\mu_{\theta_\varepsilon(t)}$ and
$\widetilde m_t=\theta_\varepsilon'(t)m_{\theta_\varepsilon(t)}$,
with zero momentum on the two stationary intervals. This preserves
the equation and endpoints, and preserves nonnegativity when it is imposed. Quadratic homogeneity and
a change of variables give
\begin{equation}\label{geo:eq:recovery-time-change}
 \int_0^1\Ac(\widetilde\mu_t,\widetilde m_t)\dd t
 =\frac{I}{1-4\varepsilon}.
\end{equation}
Extend $\widetilde\mu_t$ outside $[0,1]$ by the constant endpoint
distributions, and extend $\widetilde m_t$ by zero. The distributional equation
then holds on the whole time line.

\emph{Averaged nodes and exact constraints.}
Let $\mathsf k_\varepsilon(t)=(2\varepsilon)^{-1}
\ind_{[-\varepsilon,\varepsilon]}(t)$ and define
\[
 \mu^\varepsilon_t=\int_\R \mathsf k_\varepsilon(t-s)
                         \widetilde\mu_s\dd s,\qquad
 m^\varepsilon_t=\int_\R \mathsf k_\varepsilon(t-s)
                         \widetilde m_s\dd s.
\]
These are weak integrals of measures, and the averaged marginals are
probability measures. On $Q=\R$, a uniform bound on the second moments
of $\mu_t$ gives the same bound for these averages and their affine
interpolations. Time convolution commutes with $\partial_t$ and the
time-independent linear operator $\Lc$. It preserves nonnegative
momentum when this constraint is imposed, and the stationary intervals give
$\mu^\varepsilon_0=\mu_0$ and $\mu^\varepsilon_1=\mu_1$ exactly.
For $t_k=kh$, take
\begin{equation}\label{geo:eq:time-average-nodes}
 \mu_h^k=\mu^\varepsilon_{t_k},\qquad
 m_h^k=\int_{t_k}^{t_{k+1}}m^\varepsilon_t\dd t.
\end{equation}
Consequently
$\Lc m_h^k=\mu_h^{k+1}-\mu_h^k$. If $m_t\ge0$ almost everywhere,
then also $m_h^k\ge0$.

\emph{Comparison with the midpoint weight.}
For each $s\in\R$ put
\[
 b_k(s)=\frac12\bigl(\mathsf k_\varepsilon(t_k-s)
                         +\mathsf k_\varepsilon(t_{k+1}-s)\bigr),\qquad
 \omega_k(s)=\int_{t_k}^{t_{k+1}}\mathsf k_\varepsilon(t-s)\dd t.
\]
Then
\[
 \frac{\mu_h^k+\mu_h^{k+1}}2
       =\int_\R b_k(s)\widetilde\mu_s\dd s,\qquad
 m_h^k=\int_\R\omega_k(s)\widetilde m_s\dd s.
\]
Since $h<2\varepsilon$, $\omega_k(s)>0$ implies $b_k(s)>0$,
up to irrelevant boundary values. Joint convexity and homogeneity of
the perspective action give the weighted mixture inequality
\begin{equation}\label{geo:eq:weighted-time-mixture}
 \Ac\left(\int b_k\widetilde\mu_s\dd s,
          \int\omega_k\widetilde m_s\dd s\right)
 \le\int_\R\frac{\omega_k(s)^2}{b_k(s)}
                  \Ac(\widetilde\mu_s,\widetilde m_s)\dd s.
\end{equation}
Here $0/0$ is defined as zero. To see the inequality, average the pairs
$(\widetilde\mu_s,(\omega_k/b_k)\widetilde m_s)$ with weight $b_k(s)\dd s$;
this weight has mass one, and the action is quadratic in momentum.

For fixed $s$, compare grid cells with the interval
$[s-\varepsilon,s+\varepsilon]$. A cell completely inside this interval
has $\omega_k=h/(2\varepsilon)$ and $b_k=1/(2\varepsilon)$, so it
contributes $h/(2\varepsilon)$ to
$\sum_k\omega_k^2/(h b_k)$. The disjoint full cells have total length
at most $2\varepsilon$, hence their contribution is at most one.
At most two other cells meet the interval in positive length.
Each has $b_k=1/(4\varepsilon)$ and
$0\le\omega_k\le h/(2\varepsilon)$, so contributes at most
$h/\varepsilon$. Therefore, for almost every $s$,
\begin{equation}\label{geo:eq:time-kernel-bound}
 \sum_{k=0}^{N-1}\frac{\omega_k(s)^2}{h b_k(s)}
 \le1+\frac{2h}{\varepsilon}.
\end{equation}
Sum \eqref{geo:eq:weighted-time-mixture} over $k$, divide by $h$,
and use \eqref{geo:eq:time-kernel-bound} and
\eqref{geo:eq:recovery-time-change} to obtain
\eqref{geo:eq:time-average-bound}.

\emph{Convergence of the recovery pairs.}
The maps $\theta_\varepsilon$ converge uniformly to the identity.
Uniform narrow continuity of $\mu_t$, followed by time averaging over
an interval of radius $\varepsilon$ and affine interpolation on
cells of length $h$, gives uniform narrow convergence of the marginals.
For momentum, change variables
$t=2\varepsilon+(1-4\varepsilon)s$ on the nonstationary interval.
Against a continuous compactly supported space--time test function,
the factor $\theta_\varepsilon'$ cancels with $\dd t$.
The resulting time argument tends uniformly to $s$, while
$\int_0^1|m_s|(Q)\dd s<\infty$. Dominated convergence proves convergence
of the reparametrized space--time momenta. Time averaging and then
replacement by their cell averages change the time argument of the
test function by at most $\varepsilon+h$; their total variations stay
bounded. These operations do not change the spatial coordinate, so the
spatial tails of their total variations are bounded by those of the
finite measure $\dd t\otimes|m_t|$. This proves the asserted weak
convergence, also when $Q=\R$.
\end{proof}

\subsection{Proof of \texorpdfstring{\cref{geo:thm:time-discrete-convergence}}{the convergence theorem}}
\begin{proof}
\emph{Step 1: compactness and the liminf inequality.}
Let \(\mu_N\in\mathcal X\) satisfy
\(\liminf_N\mathcal J_{\mathrm g,N}(\mu_N)<\infty\).  The compactness
argument uses no convergence of \(\mu_N\).  For the \(\Gamma\)-liminf
inequality, assume in addition \(\mu_N\to\mu\) in \(\mathcal X\).
Hypothesis \ref{geo:hyp:order} gives
\begin{equation}\label{geo:eq:Lstar-bound}
  \|\Lc^*\varphi\|_\infty
  \le C_\Lc\|\varphi\|_{C^{r+2}},
  \qquad
  C_\Lc:=C\sum_{k\ne0}(1+|k|)^r(2\pi|k|)^{-(r+2)}<\infty.
\end{equation}
Indeed, \(\mathfrak p(0)=0\) and \(|\mathfrak p(-k)|=|\mathfrak p(k)|\) by
\ref{geo:hyp:real}, and
\(|\widehat\varphi(k)|\le\|\varphi\|_{C^{r+2}}(2\pi|k|)^{-(r+2)}\) for
\(k\ne0\), so the Fourier series of \(\Lc^*\varphi\) converges absolutely
for every \(\varphi\in C^{r+2}(\T)\).  Both sides of
\(\int\varphi\,\dd(\Lc m)=\int \Lc^*\varphi\,\dd m\) are
\(\|\cdot\|_{C^{r+2}}\)-continuous in \(\varphi\), so this identity extends
from \(C^\infty(\T)\) to every \(\varphi\in C^{r+2}(\T)\), in particular to
the test functions with \(\|\varphi\|_{C^{r+2}}\le1\) in the definition of
\(d_{-(r+2)}\).
Pass to a subsequence realizing \(\liminf_N\mathcal J_{\mathrm g,N}(\mu_N)\),
and choose momenta whose costs are within \(N^{-1}\) of
\(\mathcal J_{\mathrm g,N}(\mu_N)\).  Their total actions are bounded by some
\(\bar A\).
Write
\[
  \mathsf a_j=\frac1h\Ac(\bar\mu_h^j,m^j),
  \qquad
  \sum_j\,\mathsf a_j\le\bar A.
\]
Since \(\bar\mu_h^j\) is a probability measure,
\begin{equation}\label{geo:eq:action-cs}
  |m^j|(\T)\le\sqrt{2\Ac(\bar\mu_h^j,m^j)}
  =\sqrt{2h\,\mathsf a_j}.
\end{equation}
For \(s<t\), the interval \([s,t]\) meets at most \(n\le(t-s)/h+2\) cells.
For each such cell \(j\) and each \(\varphi\) with
\(\|\varphi\|_{C^{r+2}}\le1\),
\(|\int\varphi\,\dd(\mu_h^{j+1}-\mu_h^j)|=|\int\Lc^*\varphi\,\dd m^j|
\le C_\Lc|m^j|(\T)\le C_\Lc\sqrt{2h\,\mathsf a_j}\) by \eqref{geo:eq:Lstar-bound} and
\eqref{geo:eq:action-cs}.  By Cauchy--Schwarz over these cells,
\(\sum_j\sqrt{2h\,\mathsf a_j}\le\sqrt{2hn}\sqrt{\sum_j\,\mathsf a_j}
\le\sqrt{2\bar A}\sqrt{|t-s|+2h}\).  Summing the cell bounds gives
\begin{equation}\label{geo:eq:holder-estimate}
  d_{-(r+2)}(\mu_N(t),\mu_N(s))
  \le C_\Lc\sqrt{2\bar A}\sqrt{|t-s|+2h}.
\end{equation}
Let \(\delta'>0\) and \(N_0\ge2/\delta'\).  By \eqref{geo:eq:holder-estimate},
\(d_{-(r+2)}(\mu_N(t),\mu_N(s))\le2C_\Lc\sqrt{\bar A\delta'}\) whenever
\(|t-s|\le\delta'\) and \(N\ge N_0\).  For each \(N<N_0\), the affine
interpolants of \(N+1\) nodes form the continuous image of the compact set
\(\Pp(\T)^{N+1}\), hence a compact and therefore equicontinuous subset
of \(\mathcal X\).  Since \(2C_\Lc\sqrt{\bar A\delta'}\to0\) as
\(\delta'\to0\), and only finitely many indices satisfy \(N<N_0\), the
family \(\{\mu_N\}\) is equicontinuous. Since \(\Pp(\T)\) is compact,
the Arzel\`a--Ascoli theorem applies to it.
Thus the family $(\mathcal J_{\mathrm g,N})_{N\ge1}$ is equicoercive:
for every $B<\infty$, the union
$\bigcup_{N\ge1}\{\lambda\in\mathcal X:
\mathcal J_{\mathrm g,N}(\lambda)\le B\}$ is relatively compact in
$\mathcal X$. In particular, every sequence of paths with uniformly
bounded discrete actions has a uniformly weakly convergent subsequence.
Set \(\widetilde m_N(t)=m^j/h\) and
\(\bar\mu_h(t)=\bar\mu_h^j\) on the \(j\)-th cell.  Estimate
\eqref{geo:eq:action-cs} gives
\[
  \int_0^1|\widetilde m_N(t)|(\T)\,\dd t
  =\sum_j|m^j|(\T)\le\sqrt{2\bar A}.
\]
Thus the measures \(\dd t\otimes\widetilde m_N(t)\) on \([0,1]\times\T\) have
total variation at most \(\sqrt{2\bar A}\) and are relatively compact for
weak-star convergence.  After extraction,
\(\dd t\otimes\widetilde m_N(t)\rightharpoonup\mathsf W\) for a signed measure
\(\mathsf W\).  On the
\(j\)-th cell the midpoint and the affine density interpolations differ by
at most \(C_\Lc|m^j|(\T)\le C_\Lc\sqrt{2h\bar A}\) in \(d_{-(r+2)}\); since
\(\mu_N\to\mu\) in \(\mathcal X\), they have the same limit and no further
extraction is needed for
\[
  \dd t\otimes\bar\mu_h(t)\rightharpoonup
  \dd t\otimes\mu_t.
\]
On the \(j\)-th cell, \eqref{geo:eq:action-homogeneity} gives
\(\int_{t_j}^{t_{j+1}}\Ac(\bar\mu_h^j,m^j/h)\,\dd t=h^{-1}\Ac(\bar\mu_h^j,m^j)\).
Summing over \(j\) gives
\begin{equation}\label{geo:eq:discrete-space-time-action}
  \Ac(
    \dd t\otimes\bar\mu_h,
    \dd t\otimes\widetilde m_N)
  =
  \sum_{j=0}^{N-1}\frac1h\Ac(\bar\mu_h^j,m^j).
\end{equation}
The pair \((\mu_N,\widetilde m_N)\) satisfies \eqref{geo:eq:continuity-weak}
exactly.  For \(\varphi\in C_c^\infty((0,1)\times\T)\),
\eqref{geo:eq:Lstar-bound} gives \(\Lc^*\varphi\in C([0,1]\times\T)\).
Uniform convergence of \(\mu_N\) and weak-star convergence of
\(\dd t\otimes\widetilde m_N\) therefore pass both terms of
\eqref{geo:eq:continuity-weak} to the limit, so \((\mu,\mathsf W)\) satisfies
\eqref{geo:eq:continuity-weak}.  Weak lower semicontinuity of \(\Ac\)
(Lemma~\ref{geo:lem:space-time-action}), applied along
\((\dd t\otimes\bar\mu_h,\dd t\otimes\widetilde m_N)
\rightharpoonup(\dd t\otimes\mu_t,\mathsf W)\), together with the disintegration
\(\mathsf W=\dd t\otimes m_t\) of the same lemma, gives
\[
  \mathcal J_{\mathrm g}(\mu)
  \le
  \int_0^1\Ac(\mu_t,m_t)\,\dd t
  \le
  \liminf_{N\to\infty}\mathcal J_{\mathrm g,N}(\mu_N).
\]
Uniform convergence in time preserves both endpoints.  This proves the
\(\Gamma\)-liminf inequality as well as compactness of bounded sublevels.
\emph{Step 2: construction of a recovery sequence.}
Fix $\mu\in\mathcal X$ with $\mathcal J_{\mathrm g}(\mu)<\infty$.
For every sufficiently small $h=1/N$, apply
\cref{geo:lem:time-average-recovery} with $\varepsilon_h=\sqrt h$ to an
admissible pair $(\mu,m)$, and denote the recovered affine curve by
$\mu_h$. The recovered marginal nodes depend only on $\mu$, not on
the choice of momentum realizing this curve. Thus the bound
\eqref{geo:eq:time-average-bound} holds for the same $\mu_h$ for every
such momentum, and taking their infimum gives
\[
 \mathcal J_{\mathrm g,N}(\mu_h)
 \le \frac{1+2\sqrt h}{1-4\sqrt h}
             \mathcal J_{\mathrm g}(\mu).
\]
Since $\mu_h\to\mu$ in $\mathcal X$ and the prefactor tends to one,
this proves the curve-wise $\Gamma$-limsup inequality. If
$\mathcal J_{\mathrm g}(\mu)=+\infty$, the limsup requirement is
automatic. No density lower bound or spatial regularization is used.

\emph{Step 3: minima and minimizers.}
Step~1 gives convergent subsequences for paths of bounded discrete
action. Its lower bound, together with the matching discrete
approximations of Step~2, gives \eqref{geo:eq:value-convergence}.
The continuum infimum is attained by the direct method.  Indeed, every admissible pair \((\mu,m)\) with finite action
satisfies
\[
  d_{-(r+2)}(\mu_t,\mu_s)
  \le
  C_\Lc\sqrt{2|t-s|
    \int_s^t\Ac(\mu_u,m_u)\,\dd u}.
\]
Arzel\`a--Ascoli, weak compactness of the space--time momenta, and
Lemma~\ref{geo:lem:space-time-action} then give a minimizing pair.
The recovery construction supplies a finite discrete competitor for all
sufficiently large \(N\).  Fix such \(N\) and a minimizing sequence
\((\mu_{h,(n)}^j,m_{(n)}^j)_n\).  The nodes lie in the compact set
\(\Pp(\T)^{N-1}\), and \eqref{geo:eq:action-cs} bounds \(|m_{(n)}^j|(\T)\)
uniformly in \(n\).  After extraction, both converge weakly.  For smooth
\(\varphi\), \(\Lc^*\varphi\in C(\T)\) by \ref{geo:hyp:order}, so the
constraints \(\Lc m^j=\mu_h^{j+1}-\mu_h^j\) pass to the limit.  Lower
semicontinuity of \(\Ac\) then shows that the limit is a discrete minimizer.
Finally, Step~1 supplies convergent subsequences of asymptotically
minimizing discrete paths; its lower bound and
\eqref{geo:eq:value-convergence} show that every such limit is
a minimizer of the continuous problem.
\end{proof}
\subsection{Proof of \texorpdfstring{\cref{geo:prop:heat-semigroup-hypotheses}}{the diffusive-operator proposition}}
\begin{proof}
We verify \ref{geo:hyp:real}--\ref{geo:hyp:order} and prove
$\mathsf G_{\partial_{xx}}(\mu^0,\mu^1)<\infty$ for arbitrary
$\mu^0,\mu^1\in\Pp(\T)$.
The symbol \(\mathfrak p(k)=-4\pi^2k^2\) is real and even with
\(\mathfrak p(0)=0\), so \ref{geo:hyp:real} holds, and \ref{geo:hyp:order}
holds with \(r=2\) and \(C=4\pi^2\), since
\(|\mathfrak p(k)|\le4\pi^2(1+|k|)^2\).
It remains to prove \(\mathsf G_{\partial_{xx}}(\mu^0,\mu^1)<\infty\) for
arbitrary \(\mu^0,\mu^1\in\Pp(\T)\). Let $\pi_s$ be
the periodic heat kernel
\[
 \pi_s(x)=\sum_{n\in\mathbb Z}\frac1{\sqrt{4\pi s}}
                 \exp\left(-\frac{|x+n|^2}{4s}\right).
\]
For $s>0$ this is a smooth strictly positive probability density.
For $\lambda\in\Pp(\T)$, the convolution $\pi_s*\lambda$ satisfies
$\pi_s*\lambda\rightharpoonup\lambda$ as $s\downarrow0$, and
$\partial_s(\pi_s*\lambda)=\partial_{xx}(\pi_s*\lambda)$.
Fix \(\delta\in(0,1/2)\).  On \([0,\delta]\), the pair
\[
  \mu_t=\pi_t*\mu^0,
  \qquad
  m_t=\pi_t*\mu^0
\]
joins \(\mu^0\) to \(\pi_\delta*\mu^0\) with action \(\delta/2\).  Use the reversed
connector
\[
  \mu_t=\pi_{1-t}*\mu^1,
  \qquad
  m_t=-\pi_{1-t}*\mu^1
\]
on \([1-\delta,1]\).  Write
\(\pi_\delta*\mu^i=\rho_i\,\dd x\).  The functions \(\rho_i\) are smooth and
strictly positive, so
\[
  \rho_s=(1-s)\rho_0+s\rho_1
\]
has a uniform positive lower bound.  The smooth mean-zero function
\(\rho_1-\rho_0\) lies in the range of \(\partial_{xx}\); let \(w\) be the
smooth periodic mean-zero solution of \(\partial_{xx}w=\rho_1-\rho_0\); then
\[
  \partial_s\rho_s=\partial_{xx}w,
  \qquad
  \int_0^1\Ac(\rho_s\,\dd x,w\,\dd x)\,\dd s<\infty.
\]
On \([\delta,1-\delta]\), set
\[
  s=\frac{t-\delta}{1-2\delta},
  \qquad
  \mu_t=\rho_s\,\dd x,
  \qquad
  m_t=\frac{w}{1-2\delta}\,\dd x.
\]
This middle pair has finite action and matches the two connectors, proving
finiteness.
\end{proof}
\subsection{Positive Fourier averaging}
\label{geo:app:joint-grid}
We prove the joint result for the exact Fourier symbol specified in
\cref{geo:sec:joint-grid}. Throughout this subsection and the following
ones, $\Lc=\partial_{xx}$. Recall the operators of
\cref{geo:def:fourier-grid-operators}: with $\Delta x=1/M$, $x_i=i/M$
and $K_M=\lfloor(M-1)/4\rfloor$,
\[
 \begin{aligned}
 H_M&=(1-M^{-1})F_{K_M}+M^{-1},
 & (P_M\lambda)_i&=\int_\T H_M(x_i-y)\dd\lambda(y),\\
 \widehat L_M(k)&=-4\pi^2k^2\quad(k\in\mathcal I_M),
 & \iota_Mu&=\Delta x\sum_i u_i\delta_{x_i}.
 \end{aligned}
\]
Here $F_K$ is the Fej\'er kernel in \eqref{geo:eq:fejer-kernel},
$L_M$ acts on grid vectors, and $\iota_M$ reconstructs an atomic measure.

\begin{lemma}[Positive Fourier averaging]
\label{geo:lem:positive-fourier-map}
The map $P_M$ of \eqref{geo:eq:fejer-kernel} preserves mass and maps
probability measures to positive grid vectors. For a finite positive
measure $\mu$ and a finite signed measure $m$,
\begin{align}
 \Delta x\sum_i H_M(x_i-y)&=1\quad(y\in\T),
 \label{geo:eq:kernel-partition}\\
 L_MP_Mm&=P_M(\partial_{xx}m),
 \label{geo:eq:kernel-commutation}\\
 \A_M(P_M\mu,P_Mm)&\le\Ac(\mu,m).
 \label{geo:eq:kernel-action-contraction}
\end{align}
The middle identity is interpreted distributionally if necessary.
For every $\varphi\in C(\T)$, the adjoint averaging operator
\begin{equation}\label{geo:eq:kernel-adjoint}
 (P_M^*\varphi)(y):=\Delta x\sum_i\varphi(x_i)H_M(x_i-y)
\end{equation}
satisfies $\|P_M^*\varphi-\varphi\|_\infty\to0$.
Consequently $\iota_MP_M\lambda\rightharpoonup\lambda$, uniformly
when tested against a fixed continuous function over signed measures
with bounded total variation.
\end{lemma}
\begin{proof}
The Fourier support of $H_M$ is contained in $|k|\le K_M<M/2$,
its zero Fourier coefficient is one, and $H_M\ge M^{-1}$.
Summing each Fourier mode over the grid gives
\eqref{geo:eq:kernel-partition}. Differentiation of the finite Fourier
series and the exact symbol of $L_M$ give
\eqref{geo:eq:kernel-commutation}, including the constant mode.

For \eqref{geo:eq:kernel-action-contraction}, the case $\mu=0$ is
immediate. Otherwise it suffices to consider $m=v\mu$ with
$v\in L^2(\mu)$. Weighted Cauchy--Schwarz gives
\[
 \frac{|(P_Mm)_i|^2}{(P_M\mu)_i}
 \le\int H_M(x_i-y)|v(y)|^2\dd\mu(y).
\]
Multiply by $\Delta x/2$, sum over $i$, and use
\eqref{geo:eq:kernel-partition}.

For a fixed Fourier character $e_\ell(x)=e^{2\pi\mathrm i\ell x}$,
once $K_M+|\ell|<M$, exact grid summation gives
$P_M^*e_\ell=\widehat H_M(\ell)e_\ell$.
For $\ell=0$ the coefficient is one; for each fixed nonzero $\ell$
and all sufficiently large $M$, it is
$(1-M^{-1})(1-|\ell|/(K_M+1))$, which tends to one.
Thus $P_M^*\varphi\to\varphi$ uniformly for trigonometric polynomials.
Positivity and \eqref{geo:eq:kernel-partition} give
$\|P_M^*\varphi\|_\infty\le\|\varphi\|_\infty$, so density extends
the convergence to every continuous $\varphi$. Finally,
\[
 \left|\int\varphi\dd(\iota_MP_M\lambda-\lambda)\right|
 \le |\lambda|(\T)\,\|P_M^*\varphi-\varphi\|_\infty.
\]
\end{proof}

\subsection{Joint compactness and the lower bound}
\begin{lemma}[Grid compactness and lower bound]
\label{geo:lem:joint-grid-liminf}
Let $M_n,N_n\to\infty$, and let feasible grid paths have actions
$I_n$ given by \eqref{geo:eq:joint-grid-functional}, with
$\sup_n I_n\le C<\infty$. Their atomic affine marginal
reconstructions have a uniformly weakly convergent subsequence.
Along a further subsequence the momenta converge as signed
space--time measures. Every resulting pair $(\mu,m)$ is admissible,
and its action is at most the lower limit of the corresponding grid
actions. Weak limits of the prescribed endpoint reconstructions are
preserved.
\end{lemma}
\begin{proof}
Let $h_n=1/N_n$ and let $\widetilde\mu_n$ interpolate
$\iota_{M_n}\nu_n^j$ affinely. On each time cell set
\begin{equation}\label{geo:eq:joint-grid-measures}
 \bar\mu_n(t)=\iota_{M_n}q_n^j,\qquad
 w_n(t)=\iota_{M_n}(m_n^j/h_n),\qquad
 R_n=\dd t\otimes\bar\mu_n(t),\quad W_n=\dd t\otimes w_n(t).
\end{equation}
The action identity and Cauchy--Schwarz give
\begin{equation}\label{geo:eq:joint-grid-action-identity}
 \Ac(R_n,W_n)=I_n,\qquad
 \int_0^1|w_n(t)|(\T)^2\dd t\le2I_n\le2C.
\end{equation}
Hence $|W_n|\le\sqrt{2C}$ and, after extraction, $W_n\rightharpoonup W$.

For the character $e_k(x)=e^{2\pi\mathrm i kx}$, let
$k'\in\mathcal I_M$ satisfy $k'\equiv k\pmod M$. The nodal vector of
$e_k$ is an eigenvector of $L_M$ with eigenvalue $-4\pi^2k'^2$,
and $|k'|\le k$ for $k\ge1$. Conservation and
\eqref{geo:eq:joint-grid-action-identity} therefore give
\begin{equation}\label{geo:eq:joint-grid-fourier-modulus}
 |\langle\widetilde\mu_n(t)-\widetilde\mu_n(s),e_k\rangle|
 \le4\pi^2k^2\sqrt{2C|t-s|}.
\end{equation}
The metric
\[
 d_F(\lambda,\eta)=\sum_{k\ge1}\frac{2^{-k}}{k^2}
                  |\langle\lambda-\eta,e_k\rangle|
\]
metrizes weak convergence on $\Pp(\T)$. Arzel\`a--Ascoli yields
$\widetilde\mu_n\to\mu$ uniformly weakly, preserving the endpoint
limits. The midpoint and affine reconstructions differ by at most
$4\pi^2k^2\sqrt{2Ch_n}$ on $e_k$, so
$R_n\rightharpoonup R:=\dd t\otimes\mu_t$ as well.

For a fixed trigonometric polynomial $\varphi$, eventually
$L_{M_n}\varphi=\varphi''$ at the nodes. Thus, for
$\eta\in C_c^\infty((0,1))$,
\[
 \int_0^1\eta'(t)\langle\widetilde\mu_n(t),\varphi\rangle\dd t
 +\int\eta(t)\varphi''(x)\dd W_n(t,x)=0.
\]
Passing to the limit and approximating smooth tests by finite Fourier
sums proves $\partial_t\mu=\partial_{xx}W$.
Finally, \cref{geo:lem:space-time-action} gives
$W=\dd t\otimes m_t$ and
\begin{equation}\label{geo:eq:joint-grid-lower-bound}
 \int_0^1\Ac(\mu_t,m_t)\dd t
 =\Ac(R,W)\le\liminf_n I_n.
\end{equation}
The argument applies along any jointly convergent subsequence,
with its corresponding lower limit of the actions.
\end{proof}

\subsection{Recovery and the joint geodesic theorem}
\begin{proof}[Proof of \cref{geo:thm:joint-grid}]
\Cref{geo:lem:positive-fourier-map} gives positive endpoint vectors
converging weakly to the prescribed marginals.
\Cref{geo:lem:joint-grid-liminf} gives compactness of paths with
bounded action and the $\Gamma$-lower bound, by choosing momenta
within $1/n$ of each curve functional.

For the upper bound, take $\lambda\in\mathcal X$ with finite action
and an admissible momentum $m$ of action $I$. Apply
\cref{geo:lem:time-average-recovery} with $h_n=1/N_n$ and
$\varepsilon_n=\sqrt{h_n}$, and set
\[
 \nu_n^j=P_{M_n}\mu_{h_n}^j,\qquad
 p_n^j=P_{M_n}m_{h_n}^j.
\]
No spatial regularity of $\lambda$ is needed: $P_M\lambda$ evaluates
the smooth convolution $H_M*\lambda$, including for singular measures.
The marginal nodes are taken after time averaging, and the cell
momenta are time integrals. Thus neither construction depends on
choices of representatives on null sets.
Commutation gives the exact grid constraints and endpoints, while
contraction of the action gives
\begin{equation}\label{geo:eq:joint-grid-recovery-bound}
 \sum_j\frac1{h_n}\A_{M_n}
       \left(\frac{\nu_n^j+\nu_n^{j+1}}2,p_n^j\right)
 \le\frac{1+2\sqrt{h_n}}{1-4\sqrt{h_n}}\,I.
\end{equation}
For every node and every $\varphi\in C(\T)$,
$\int\varphi\,\dd(\iota_{M_n}\nu_n^j-\mu_{h_n}^j)
=\int(P_{M_n}^*\varphi-\varphi)\,\dd\mu_{h_n}^j$.
Since $\|P_M^*\varphi\|_\infty\le\|\varphi\|_\infty$ and
$P_M^*\varphi\to\varphi$ for every $\varphi\in C(\T)$
(\cref{geo:lem:positive-fourier-map}), this convergence is uniform on
the unit ball of $C^4(\T)$, which is relatively compact in $C(\T)$.
Hence, at every time, the reconstructed curve and the affine
interpolation of the nodes $\mu_{h_n}^j$ differ by at most
$\sup_{\|\varphi\|_{C^4}\le1}\|P_{M_n}^*\varphi-\varphi\|_\infty\to0$
in $d_{-4}$, the distance $d_{-(r+2)}$ with $r=2$.
This bound and the time-averaging lemma imply convergence of the reconstructed
curves to $\lambda$ in $\mathcal X$, independently of the relative
sizes of $M_n$ and $N_n$. Their nodes depend only on $\lambda$;
taking the infimum over its compatible momenta proves the
$\Gamma$-upper bound for the curve functional.

For fixed $M,N$, positive endpoints admit the affine grid path with
$m^j=C_M(\nu^N-\nu^0)/N$, where $C_M$ has symbol
$-(4\pi^2k^2)^{-1}$ for $k\ne0$ and zero for $k=0$.
This competitor has finite action. Compactness of the density simplex,
the action bound on momenta, and lower semicontinuity give a grid
minimizer. The continuous minimum is finite and attained by
\cref{geo:thm:time-discrete-convergence,geo:prop:heat-semigroup-hypotheses}.
The lower and upper bounds, together with compactness of bounded-action
paths, give convergence of the minimum values and show that every
limit of asymptotically minimizing paths minimizes the continuous problem.
\end{proof}

\subsection{Changing the prescribed grid endpoints}
\begin{proof}[Proof of \cref{geo:cor:grid-endpoints}]
Write $a_M^0=\nu_M^0$, $a_M^1=\nu_M^N$ and
$b_M^\ell=P_M\mu^\ell$. Passing to the limit in
$\iota_Ma_M^\ell\ge c\Delta x\sum_i\delta_{x_i}$ gives
$\mu^\ell\ge c\dd x$, hence $b_M^\ell\ge c\one$.
The lower bound already follows from \cref{geo:lem:joint-grid-liminf};
we construct the upper approximation.

The differences $f_M^\ell=a_M^\ell-b_M^\ell$ have grid mean zero,
$\Delta x\sum_i|f_{M,i}^\ell|\le2$ and
$\iota_Mf_M^\ell\rightharpoonup0$. Let $C_M$ be the mean-zero
inverse of $L_M$. Its one-dimensional Fourier symbol is summable, so
\begin{equation}\label{geo:eq:grid-endpoint-inverse-error}
 \varepsilon_M:=\max_{\ell=0,1}\|C_Mf_M^\ell\|_\infty\longrightarrow0.
\end{equation}
Indeed, fixed Fourier coefficients tend to zero, while for every
$k_0\ge1$ the modes $|k|>k_0$ contribute at most
$(2/(4\pi^2))\sum_{|k|>k_0}k^{-2}$, uniformly in $M$.

For $h=1/N$, choose
\begin{equation}\label{geo:eq:grid-endpoint-bridge-time}
 J=\left\lceil\frac{\varepsilon_M+\sqrt h}{h}\right\rceil,
 \qquad \delta=Jh,\qquad N'=N-2J.
\end{equation}
Then $\delta\to0$ and $N'\to\infty$. Interpolate from $a_M^0$ to
$b_M^0$ over the first $J$ cells, using momentum
$C_M(b_M^0-a_M^0)/J$ in each cell. All midpoint weights are at least
$c$, so the bridge costs at most
\begin{equation}\label{geo:eq:grid-endpoint-bridge-cost}
 \sum_{j<J}\frac1h\A_M(q^j,m^j)\le\frac{\varepsilon_M^2}{2c\delta}.
\end{equation}
The last $J$ cells join $b_M^1$ to $a_M^1$ at the same cost bound;
signed momenta allow both corrections.

On the middle $N'$ cells use the approximation from
\cref{geo:thm:joint-grid}, compressed into $[\delta,1-\delta]$.
For a continuous pair of action $I$, the total action is at most
\begin{equation}\label{geo:eq:grid-endpoint-total-recovery}
 \frac{1}{1-2\delta}
 \frac{1+2/\sqrt{N'}}{1-4/\sqrt{N'}}\,I
 +\frac{\varepsilon_M^2}{c\delta}.
\end{equation}
This tends to $I$, since $\varepsilon_M^2/\delta\le \varepsilon_M\to0$.
The middle time change tends uniformly to the identity; each bridge
joins two endpoint reconstructions with the same weak limit. The
whole path therefore converges uniformly weakly to the prescribed
continuous curve. Its nodes are independent of the compatible
momentum, so taking the infimum over momenta gives the upper bound
for the curve functional. The remaining conclusions follow as in
\cref{geo:thm:joint-grid}.
\end{proof}

\section{Proofs for Eulerian minimizing movements}\label{app:merged-jko}
We first prove the signed limit \cref{jko:thm:main} under its energy
assumptions \textup{(E1)--(E3)} and density floor
\eqref{jko:eq:floor}. Each preliminary lemma uses its stated hypotheses.
The final subsection, \cref{jko:app:positive-compactness}, proves
\cref{jko:prop:positive-compactness} with only \textup{(E1)--(E2)}
and finite initial energy, without a density floor.
\subsection{The frozen cost and its first variation}

The perspective action admits the dual representation, obtained as in
\cref{geo:lem:space-time-action} by completing the square, with matrix-valued
test functions $\mathsf K$,
\begin{equation}\label{jko:eq:perspective-dual}
  \Ac(\sigma,\mathsf m)
  =\sup_{\mathsf K\in C(\Torus;\R^{d\times d})}
  \left\{\int_{\Torus}\mathsf K:\dd \mathsf m
  -\frac12\int_{\Torus}|\mathsf K|^2\dd\sigma\right\}.
\end{equation}
Consequently, it is jointly convex and weak-star lower semicontinuous in
$(\sigma,\mathsf m)$.

Let $\sigma=\rho\dd x$ with $\rho\geq c>0$.  Denote by
$\dot H^2_\rho$ the completion of smooth functions modulo constants for the
norm
\begin{equation}\label{jko:eq:weighted-H2}
  \|\varphi\|_{\dot H^2_\rho}^2
  :=\int_{\Torus}\rho|\nabla^2\varphi|^2\dd x.
\end{equation}
For mean-zero $\varphi$ on $\Torus$, the second-order Poincar\'e inequality
gives
\begin{equation}\label{jko:eq:poincare}
  \|\varphi\|_{H^2}\leq C_P\|\nabla^2\varphi\|_{L^2};
\end{equation}
hence, since $\rho\geq c$,
$\|\varphi\|_{H^2}\leq C_Pc^{-1/2}\|\varphi\|_{\dot H^2_\rho}$.
Whenever $\mathsf G_{\fr}(\mu,\sigma)<+\infty$, the pairing
$\langle\sigma-\mu,\eta\rangle$ below denotes the unique continuous extension
to $\dot H^2_\rho$ of the distributional pairing on smooth functions.  A
matrix-valued measure $\mathsf m$ with $\DivDiv\mathsf m=\sigma-\mu$ is called a
feasible momentum for $(\mu,\sigma)$.  A
finite-action feasible momentum supplies such an extension; any two choices agree
on smooth functions and hence, by density, on the completion.  For $d\geq4$,
elements of $H^2(\Torus)$ need not be continuous, so
$\langle\sigma-\mu,\eta\rangle$ cannot be defined by integrating $\eta$
against $\mu$.

\begin{lemma}[Weighted biharmonic duality]\label{jko:lem:duality}
Let $\mu\in\Pcal(\Torus)$ and let $\sigma=\rho\dd x\in\Pcal(\Torus)$ satisfy
$\rho\geq c>0$ and $\mathsf G_{\fr}(\mu,\sigma)<+\infty$.  There is a unique
mean-zero $\Phi\in\dot H^2_\rho$ such that
\begin{equation}\label{jko:eq:weighted-equation}
  \int_{\Torus}\rho\nabla^2\Phi:\nabla^2\eta\dd x
  =\langle\sigma-\mu,\eta\rangle
  \qquad\text{for all }\eta\in\dot H^2_\rho.
\end{equation}
It satisfies
\begin{equation}\label{jko:eq:frozen-duality}
  \mathsf G_{\fr}(\mu,\sigma)
  =\sup_{\eta\in\dot H^2_\rho}
  \left\{
    \langle\sigma-\mu,\eta\rangle
    -\frac12\int_{\Torus}\rho|\nabla^2\eta|^2\dd x
  \right\},
\end{equation}
and the unique optimal momentum in \eqref{jko:eq:frozen-cost} is
\begin{equation}\label{jko:eq:optimal-M}
  \mathsf m=\rho\nabla^2\Phi\dd x.
\end{equation}
In particular,
\begin{equation}\label{jko:eq:cost-potential-identity}
  2\mathsf G_{\fr}(\mu,\sigma)
  =\int_{\Torus}\rho|\nabla^2\Phi|^2\dd x.
\end{equation}
\end{lemma}

\begin{proof}
Choose one finite-action feasible momentum $\mathsf m=\mathsf v\sigma$.  For smooth
$\eta$,
\[
  |\langle\sigma-\mu,\eta\rangle|
  =\left|\int_{\Torus}\mathsf v:\nabla^2\eta\dd\sigma\right|
  \leq \bigl(2\Ac(\sigma,\mathsf m)\bigr)^{1/2}
       \|\eta\|_{\dot H^2_\rho}.
\]
Thus the right-hand side of \eqref{jko:eq:weighted-equation} is a continuous
linear functional on $\dot H^2_\rho$.  Riesz's theorem gives a unique
mean-zero $\Phi$ satisfying \eqref{jko:eq:weighted-equation}.

For every feasible momentum $\mathsf m=\mathsf v\sigma$ and every $\eta\in\dot H^2_\rho$, Young's
inequality gives
\[
  \langle\sigma-\mu,\eta\rangle
  =\int \mathsf v:\nabla^2\eta\dd\sigma
  \leq \Ac(\sigma,\mathsf m)
       +\frac12\int|\nabla^2\eta|^2\dd\sigma.
\]
This proves that the right-hand side of \eqref{jko:eq:frozen-duality} is a lower
bound for $\mathsf G_{\fr}$.  On the other hand,
$\mathsf m_\Phi:=\rho\nabla^2\Phi\dd x$ is feasible by
\eqref{jko:eq:weighted-equation}, and testing that equation with $\Phi$ yields
\[
  \langle\sigma-\mu,\Phi\rangle
  =\int\rho|\nabla^2\Phi|^2\dd x.
\]
Hence the functional maximized in \eqref{jko:eq:frozen-duality} takes at
$\eta=\Phi$ the value
$\frac12\int\rho|\nabla^2\Phi|^2\dd x=\Ac(\sigma,\mathsf m_\Phi)
\geq\mathsf G_{\fr}(\mu,\sigma)$.  With the opposite inequality
above, \eqref{jko:eq:frozen-duality} holds, $\Phi$ is a maximizer, and
$\mathsf m_\Phi$ is optimal in \eqref{jko:eq:frozen-cost}.  Equality in Young's
inequality forces $\mathsf v=\nabla^2\Phi$, which also proves uniqueness of the
optimal momentum and \eqref{jko:eq:cost-potential-identity}.
\end{proof}

\begin{lemma}[First variation in the final density]\label{jko:lem:first-variation}
Under the assumptions of \cref{jko:lem:duality}, let
$h\in C^\infty(\Torus)$ have zero mean.  For all sufficiently small $|s|$,
put $\sigma_s=(\rho+sh)\dd x$.  Then
\begin{equation}\label{jko:eq:cost-first-variation}
  \left.\frac{\dd}{\dd s}\right|_{s=0}
  \mathsf G_{\fr}(\mu,\sigma_s)
  =\int_{\Torus}
     \left(\Phi-\frac12|\nabla^2\Phi|^2\right)h\dd x,
\end{equation}
where $\Phi$ is the mean-zero potential from \cref{jko:lem:duality}.
\end{lemma}

\begin{proof}
Since $\rho\geq c$ and $h$ is bounded, for
$|s|\leq s_0:=c/(2\|h\|_\infty)$,
$\frac12\rho\leq\rho+sh\leq\frac32\rho$.  The perturbed cost is finite as well.
Indeed, choose a smooth matrix field $\mathsf n$ with $\DivDiv\mathsf n=h$ (solve
$\Delta^2u=h$ with zero mean and take $\mathsf n=\nabla^2u$), and add
$s\mathsf n\dd x$ to any finite-action feasible momentum for $(\mu,\sigma)$.
Comparability of the weights makes the action of this competitor uniformly
bounded for $|s|\leq s_0$.

Let $\Phi_s$ be the mean-zero potential of \cref{jko:lem:duality} for
$(\mu,\sigma_s)$.  By \eqref{jko:eq:cost-potential-identity},
$\int(\rho+sh)|\nabla^2\Phi_s|^2\dd x=2\mathsf G_{\fr}(\mu,\sigma_s)$.  The
competitor above bounds the right side uniformly for $|s|\leq s_0$, and
$\rho\leq2(\rho+sh)$ gives
\[
 \sup_{|s|\leq s_0}\|\Phi_s\|_{\dot H^2_\rho}<+\infty
\]
with this $s_0$.  Subtracting \eqref{jko:eq:weighted-equation} for $\Phi$ from
the corresponding equation for $\Phi_s$, which has weight $\rho+sh$, gives
\[
 \int\rho\nabla^2(\Phi_s-\Phi):\nabla^2\eta\dd x
 =s\int h\eta\dd x-s\int h\nabla^2\Phi_s:\nabla^2\eta\dd x.
\]
Testing with $\eta=\Phi_s-\Phi$, and using the uniform bound on $\Phi_s$ and
the second-order Poincar\'e inequality \eqref{jko:eq:poincare}, gives
$\|\Phi_s-\Phi\|_{\dot H^2_\rho}\leq C|s|$, hence
\begin{equation}\label{jko:eq:potential-stability}
  \|\Phi_s-\Phi\|_{\dot H^2_\rho}\longrightarrow0
  \qquad(s\to0).
\end{equation}

Write $G(s)=\mathsf G_{\fr}(\mu,\sigma_s)$.  The dual problem at $s$
is the supremum in \eqref{jko:eq:frozen-duality} with $\sigma$ replaced by
$\sigma_s$.  Testing it with $\Phi$ gives
\[
  G(s)-G(0)\geq
  s\int h\left(\Phi-\frac12|\nabla^2\Phi|^2\right)\dd x.
\]
Conversely, testing the dual problem at $0$ with $\Phi_s$ gives
\[
  G(s)-G(0)\leq
  s\int h\left(\Phi_s-\frac12|\nabla^2\Phi_s|^2\right)\dd x.
\]
Both inequalities hold for $|s|\leq s_0$.  For $s>0$, dividing by $s$ gives
\[
  \int h\left(\Phi-\frac12|\nabla^2\Phi|^2\right)\dd x
  \leq\frac{G(s)-G(0)}s
  \leq\int h\left(\Phi_s-\frac12|\nabla^2\Phi_s|^2\right)\dd x;
\]
for $s<0$, both inequalities are reversed.
The stability \eqref{jko:eq:potential-stability} now proves
\eqref{jko:eq:cost-first-variation}.  This is the envelope argument; no formal
differentiation of the optimizer is needed.
\end{proof}

\subsection{Existence and discrete Euler equations}

\begin{proposition}[Existence of one step]\label{jko:prop:existence}
If $\Ec(\mu)<+\infty$, then the minimization problem in
\eqref{jko:eq:scheme} has at least one solution.
\end{proposition}

\begin{proof}
By \eqref{jko:eq:frozen-cost}, \eqref{jko:eq:scheme} has the same infimum as
the joint problem: minimize $\frac1\tau\Ac(\rho\,\dd x,\mathsf m)+\Ec(\rho)$
over pairs $(\rho,\mathsf m)$ with $\rho\in\Pcal(\Torus)$ and
$\DivDiv\mathsf m=\rho-\mu$.  The competitor $\rho=\mu$, $\mathsf m=0$ has
finite value.  Let $(\rho_n,\mathsf m_n)$ be a minimizing sequence of this
joint problem.
Because $f\geq0$ and $V_{\rm ext},W$ are bounded, $\Ec$ is bounded below; hence the
action is bounded.  Since the action is finite,
$\mathsf m_n=\mathsf v_n\rho_n\,\dd x$, and the Cauchy--Schwarz inequality
gives
\[
  |\mathsf m_n|(\Torus)^2\leq
  2\Ac(\rho_n\dd x,\mathsf m_n)\,\rho_n(\Torus)
  =2\Ac(\rho_n\dd x,\mathsf m_n).
\]
Compactness of probability measures and weak-star compactness of bounded
sets of matrix-valued measures provide a subsequence along which
$\rho_n\dd x$ converges weakly to a probability measure $\sigma_\infty$ and
$\mathsf m_n$ converges weak-star to a matrix-valued measure $\mathsf m$.  The
affine constraint passes to the limit.  The action is lower semicontinuous by
\eqref{jko:eq:perspective-dual}.  The potential and interaction terms are
continuous under weak convergence because $V_{\rm ext},W$ are continuous.  Let
$\sigma_\infty=\rho^{\rm ac}\dd x+\sigma_\infty^\perp$ be the Lebesgue
decomposition of $\sigma_\infty$, with $\sigma_\infty^\perp$ singular.  Since
$f$ is convex and superlinear by \textup{(E2)}, the weak-star lower
semicontinuous envelope of $\int_{\Torus}f(\rho)\dd x$ is
\[
 \mathcal U(\sigma_\infty)=\int_{\Torus}f(\rho^{\rm ac})\dd x
 +f^\infty\sigma_\infty^\perp(\Torus),
 \qquad
 f^\infty:=\lim_{r\to\infty}\frac{f(r)}r=+\infty.
\]
The joint objective, with $\int_{\Torus}f(\rho)\dd x$ replaced by $\mathcal U$, is
therefore weak-star lower semicontinuous, and its value at
$(\sigma_\infty,\mathsf m)$ is at most the infimum, which is finite.  Hence $\sigma_\infty^\perp=0$, so
$\sigma_\infty=\rho^{\rm ac}\dd x$, and the limit pair
$(\rho^{\rm ac},\mathsf m)$ is a minimizer.
\end{proof}

From now on, $(\rho_k^\tau)_k$ is the sequence chosen in \eqref{jko:eq:scheme},
which satisfies the positivity floor \eqref{jko:eq:floor} by assumption.  Let $\Phi_{k+1}^\tau$ be the mean-zero potential of
\cref{jko:lem:duality} associated with
$(\rho_k^\tau,\rho_{k+1}^\tau)$, and introduce the scaled potential
\begin{equation}\label{jko:eq:scaled-potential}
  \varphi_{k+1}^\tau:=\frac{1}{\tau}\Phi_{k+1}^\tau.
\end{equation}

\begin{proposition}[Euler system]\label{jko:prop:euler-system}
Let $\tfin>0$ and let $\tau$ be small enough for \eqref{jko:eq:floor} to hold
on $[0,\tfin]$.  For every $k$ with $(k+1)\tau\leq\tfin$, there is a constant
$\lambda_{k+1}^\tau\in\R$ such that
\begin{subequations}\label{jko:eq:euler-system}
\begin{align}
  \frac{\rho_{k+1}^\tau-\rho_k^\tau}{\tau}
  &=\DivDiv\left(
      \rho_{k+1}^\tau\nabla^2\varphi_{k+1}^\tau
    \right),\label{jko:eq:euler-continuity}\\
  \varphi_{k+1}^\tau
  -\frac{\tau}{2}|\nabla^2\varphi_{k+1}^\tau|^2
  +V_{\rm ext}+W*\rho_{k+1}^\tau+f'(\rho_{k+1}^\tau)
  &=\lambda_{k+1}^\tau
  \quad\text{a.e.}\label{jko:eq:euler-HJ}
\end{align}
\end{subequations}
\end{proposition}

\begin{proof}
Equation \eqref{jko:eq:euler-continuity} is
\eqref{jko:eq:weighted-equation} divided by $\tau$.  Since
$\rho_{k+1}^\tau\geq c_{\tfin}$, for every smooth zero-mean $h$,
$\rho_{k+1}^\tau+sh$ is a probability density bounded below by $c_{\tfin}/2$
whenever $|s|\leq c_{\tfin}/(2\|h\|_\infty)$.  The growth assumption
\eqref{jko:eq:fprime-growth} guarantees both that
$f'(\rho_{k+1}^\tau)$ is integrable and that differentiation under the
integral is legitimate.  Indeed, let $q=q_{c_{\tfin}/2}$ and
$\alpha=1-1/q$.  On $[c_{\tfin}/2,+\infty)$,
\[
 \left|\bigl[(1+f)^\alpha\bigr]'\right|
 \leq\alpha C_{c_{\tfin}/2}^{1/q}.
\]
Consequently, whenever $r\geq c_{\tfin}$ and $|u|\leq c_{\tfin}/2$,
$1+f(r+u)\leq C_{\tfin}(1+f(r))$ and
$|f'(r+u)|\leq C_{\tfin}(1+f(r))^{1/q}$.  By the mean value theorem, for
$|s|\,\|h\|_\infty\leq c_{\tfin}/2$,
$|f(\rho_{k+1}^\tau+sh)-f(\rho_{k+1}^\tau)|/|s|
\leq C_{\tfin}\|h\|_\infty(1+f(\rho_{k+1}^\tau))^{1/q}$.  The right side is
integrable because $f(\rho_{k+1}^\tau)\in L^1$ and $1/q<1$, so dominated
convergence applies.
Using \cref{jko:lem:first-variation} and the evenness of $W$, the first variation
of the one-step objective in the direction $h$ is
\[
  \int_{\Torus}\left[
    \frac1\tau\left(
      \Phi_{k+1}^\tau
      -\frac12|\nabla^2\Phi_{k+1}^\tau|^2
    \right)
    +V_{\rm ext}+W*\rho_{k+1}^\tau+f'(\rho_{k+1}^\tau)
  \right]h\dd x=0.
\]
An integrable function annihilating every smooth zero-mean test function is
constant.  Substitution of \eqref{jko:eq:scaled-potential} gives
\eqref{jko:eq:euler-HJ}.
\end{proof}

\subsection{A priori estimates and time compactness}

On each interval $(k\tau,(k+1)\tau]$, define
\begin{equation}\label{jko:eq:all-interpolants}
 \underline\rho_\tau(t):=\rho_k^\tau,
 \qquad
 \varphi_\tau(t):=\varphi_{k+1}^\tau,
 \qquad
 Q_\tau(t):=\bar\rho_\tau(t)\nabla^2\varphi_\tau(t),
 \qquad
 \lambda_\tau(t):=\lambda_{k+1}^\tau,
\end{equation}
and set
\begin{equation}\label{jko:eq:defect}
  r_\tau(t,x):=\frac\tau2|\nabla^2\varphi_\tau(t,x)|^2.
\end{equation}
Then, for every $\tfin>0$ and every $\tau\in(0,\hat\tau_{\tfin}]$, on $(0,\tfin)\times\Torus$,
\begin{equation}\label{jko:eq:discrete-continuity}
  \frac{\bar\rho_\tau-\underline\rho_\tau}{\tau}
  =\DivDiv Q_\tau
\end{equation}
in distributions, and
\begin{equation}\label{jko:eq:interpolated-HJ}
  \varphi_\tau-r_\tau
  +V_{\rm ext}+W*\bar\rho_\tau+f'(\bar\rho_\tau)
  =\lambda_\tau(t).
\end{equation}

\begin{lemma}[Action estimate and vanishing defect]\label{jko:lem:action-estimate}
For every $\tfin>0$, there is $C_{\tfin}<+\infty$ such that, for every
$\tau\in(0,\hat\tau_{\tfin}]$,
\begin{align}
  \int_0^{\tfin}\!\!\int_{\Torus}
     \bar\rho_\tau|\nabla^2\varphi_\tau|^2\dd x\dd t
  &\leq 2\bigl(\Ec(\rho_0)-\inf\Ec\bigr),
    \label{jko:eq:action-bound}\\
  \|\varphi_\tau\|_{L^2(0,\tfin;H^2(\Torus))}&\leq C_{\tfin},
    \label{jko:eq:phi-H2-bound}\\
  \|r_\tau\|_{L^1((0,\tfin)\times\Torus)}&\leq C_{\tfin}\tau.
    \label{jko:eq:defect-bound}
\end{align}
Moreover,
\begin{equation}\label{jko:eq:internal-energy-bound}
  \sup_{t\in[0,\tfin]}\int_{\Torus}f(\bar\rho_\tau(t))\dd x\leq C_{\tfin}.
\end{equation}
\end{lemma}

\begin{proof}
Using $\rho_k^\tau$ as a competitor in the $(k+1)$st minimization gives
\[
  \frac{1}{\tau}\mathsf G_{\fr}(\rho_k^\tau,\rho_{k+1}^\tau)
  +\Ec(\rho_{k+1}^\tau)
  \leq\Ec(\rho_k^\tau).
\]
By \eqref{jko:eq:cost-potential-identity} and
\eqref{jko:eq:scaled-potential},
\[
  \frac{2}{\tau}\mathsf G_{\fr}(\rho_k^\tau,\rho_{k+1}^\tau)
  =\tau\int_{\Torus}
     \rho_{k+1}^\tau|\nabla^2\varphi_{k+1}^\tau|^2\dd x.
\]
Summing through the last cell intersecting $(0,\tfin)$ proves
\eqref{jko:eq:action-bound}; a partial last-cell integral is bounded by its full
cell.  The last cell meeting $(0,\tfin)$ ends before $\tfin+\tau\leq\tfin+1$.
In all estimates on $(0,\tfin)$ we therefore use the floor $\hat c_{\tfin}$ of
\eqref{jko:eq:floor} on $[0,\tfin+1]$, and from here on $C_{\tfin}$ denotes a
constant independent of $\tau\in(0,\hat\tau_{\tfin}]$, which may change from one
estimate to the next.  The positivity floor
and the second-order Poincar\'e inequality \eqref{jko:eq:poincare}, with
$\int\varphi_\tau=0$, imply \eqref{jko:eq:phi-H2-bound}.  The positivity floor
and \eqref{jko:eq:action-bound} also give
$\int_0^{\tfin}\!\int_{\Torus}|\nabla^2\varphi_\tau|^2\dd x\dd t\leq C_{\tfin}$.
Multiplying this bound by $\tau/2$ proves
\eqref{jko:eq:defect-bound}.  Finally, the potential and interaction energies are
bounded in absolute value by constants depending only on $V_{\rm ext},W$ and the unit
mass.  Energy monotonicity and $f\geq0$ therefore give
\eqref{jko:eq:internal-energy-bound}.
\end{proof}

To state the time estimate, define, for $\mu,\mu'\in\Pcal(\Torus)$,
\begin{equation}\label{jko:eq:negative-norm}
  \|\mu-\mu'\|_{-2,\infty}
  :=\sup\left\{
    |\langle\mu-\mu',\zeta\rangle|:
    \zeta\in C^\infty(\Torus),\ 
    \|\nabla^2\zeta\|_{L^\infty}\leq1
  \right\}.
\end{equation}
Constants in $\zeta$ are immaterial because $(\mu-\mu')(\Torus)=0$.  This norm
metrizes weak convergence on $\Pcal(\Torus)$.  To see this, normalize the
test functions $\zeta$ to have zero spatial mean.  The unit ball in
\eqref{jko:eq:negative-norm} is then uniformly bounded and equicontinuous in
$C^0$ by the periodic second-order Poincar\'e inequality, hence precompact in
$C^0$.  Weak convergence of probability measures is uniform on this class of
test functions.  The converse follows first for smooth test functions and then,
by density, for all continuous test functions.

\begin{lemma}[Time equicontinuity]\label{jko:lem:time-compactness}
For every $\tfin>0$, every $s,t\in[0,\tfin]$, and every $\tau\in(0,\hat\tau_{\tfin}]$,
\begin{equation}\label{jko:eq:time-modulus}
  \|\bar\rho_\tau(t)-\bar\rho_\tau(s)\|_{-2,\infty}
  \leq C_{\tfin}\sqrt{|t-s|+\tau}.
\end{equation}
\end{lemma}

\begin{proof}
It suffices first to take $s=n\tau$ and $t=m\tau$, with $m>n$.  For a test
function $\zeta$, \eqref{jko:eq:euler-continuity} and Cauchy--Schwarz give
\begin{align*}
 |\langle\rho_m^\tau-\rho_n^\tau,\zeta\rangle|
 &\leq\sum_{k=n}^{m-1}\tau
   \int_{\Torus}\rho_{k+1}^\tau
       |\nabla^2\varphi_{k+1}^\tau|\,|\nabla^2\zeta|\dd x\\
 &\leq \|\nabla^2\zeta\|_\infty
   \left(\sum_{k=n}^{m-1}\tau\right)^{1/2}
   \left(\sum_{k=n}^{m-1}\tau
      \int\rho_{k+1}^\tau
       |\nabla^2\varphi_{k+1}^\tau|^2\dd x
   \right)^{1/2}\\
 &\leq C_{\tfin}\sqrt{(m-n)\tau}\,
       \|\nabla^2\zeta\|_\infty,
\end{align*}
where we also used $\int\rho_{k+1}^\tau=1$.
Moving arbitrary $s,t$ to neighboring grid points adds at most one time cell,
which proves \eqref{jko:eq:time-modulus}.
\end{proof}

\subsection{Strong compactness of the densities}

The bound \eqref{jko:eq:phi-H2-bound} on $\varphi_\tau$ does not by itself
give compactness of $\bar\rho_\tau$ in $L^1$.  Equation
\eqref{jko:eq:interpolated-HJ} transfers the spatial regularity of
$\varphi_\tau$ to $f'(\bar\rho_\tau)$ (\cref{jko:lem:fprime-translations}),
and the strict convexity of $f$ transfers it to $\bar\rho_\tau$
(\cref{jko:lem:rho-space-translations}).
For $y\in\Torus$, write $(\mathsf S_yu)(x)=u(x+y)$.

\begin{lemma}[Spatial translations of $f'(\bar\rho_\tau)$]
\label{jko:lem:fprime-translations}
For every $\tfin>0$ and every $\tau\in(0,\hat\tau_{\tfin}]$,
\begin{equation}\label{jko:eq:fprime-translation-bound}
 \|\mathsf S_yf'(\bar\rho_\tau)-f'(\bar\rho_\tau)\|_{L^1((0,\tfin)\times\Torus)}
 \leq C_{\tfin}|y|+C_{\tfin}\tau.
\end{equation}
\end{lemma}

\begin{proof}
Subtract \eqref{jko:eq:interpolated-HJ} at $(t,x)$ from the same equation at
$(t,x+y)$.
The time-dependent constant $\lambda_\tau$ cancels.  By
\eqref{jko:eq:phi-H2-bound},
\[
 \|\mathsf S_y\varphi_\tau-\varphi_\tau\|_{L^1((0,\tfin)\times\Torus)}
 \leq |y|\,\|\nabla\varphi_\tau\|_{L^1((0,\tfin)\times\Torus)}
 \leq C_{\tfin}|y|.
\]
Similarly,
\[
 \|\mathsf S_yV_{\rm ext}-V_{\rm ext}\|_{L^1((0,\tfin)\times\Torus)}
 +\|\mathsf S_y(W*\bar\rho_\tau)-W*\bar\rho_\tau\|_{L^1((0,\tfin)\times\Torus)}
 \leq C_{\tfin}|y|,
\]
because $V_{\rm ext},W\in C^1$ and $\bar\rho_\tau$ has unit mass at every time.
Finally,
\[
  \|\mathsf S_yr_\tau-r_\tau\|_{L^1}
  \leq2\|r_\tau\|_{L^1}\leq C_{\tfin}\tau.
\]
Combining these three bounds proves \eqref{jko:eq:fprime-translation-bound}.
\end{proof}

\begin{lemma}[Spatial compactness of the densities]\label{jko:lem:rho-space-translations}
For every $\tfin>0$, with $\tau$ restricted to $(0,\hat\tau_{\tfin}]$,
\begin{equation}\label{jko:eq:rho-translation-limit}
  \lim_{|y|\to0}\ \limsup_{\tau\to0}
  \|\mathsf S_y\bar\rho_\tau-\bar\rho_\tau\|_{L^1((0,\tfin)\times\Torus)}=0.
\end{equation}
The family $(\bar\rho_\tau)_{\tau\in(0,\hat\tau_{\tfin}]}$ is also uniformly integrable on
$(0,\tfin)\times\Torus$.
\end{lemma}

\begin{proof}
Superlinearity and \eqref{jko:eq:internal-energy-bound} give uniform
integrability.  More explicitly, define
\[
  \omega_f(R):=\inf_{r\geq R}\frac{f(r)}{r};
  \qquad \omega_f(R)\longrightarrow+\infty.
\]
Then, uniformly in $t$ and $\tau$,
\begin{equation}\label{jko:eq:tail-bound}
  \int_{\{\bar\rho_\tau(t)>R\}}\bar\rho_\tau(t)\dd x
  \leq\frac{C_{\tfin}}{\omega_f(R)}.
\end{equation}

Fix $R>\hat c_{\tfin}$ and set
\[
  m_{\tfin,R}:=\min_{r\in[\hat c_{\tfin},R]}f''(r)>0.
\]
If $r,r'\in[\hat c_{\tfin},R]$, the mean value theorem and $f''\geq m_{\tfin,R}$ on
$[\hat c_{\tfin},R]$ give
\[
  |r-r'|\leq m_{\tfin,R}^{-1}|f'(r)-f'(r')|.
\]
For arbitrary $r,r'\geq \hat c_{\tfin}$, apply this inequality where both are at most
$R$.  On the set where $\max(r,r')>R$, use instead the pointwise bound
\[
 |r-r'|\ind_{\{\max(r,r')>R\}}
 \leq r\ind_{\{r>R\}}+r'\ind_{\{r'>R\}}.
\]
After integration,
\eqref{jko:eq:fprime-translation-bound} and \eqref{jko:eq:tail-bound} give
\begin{equation}\label{jko:eq:rho-translation-estimate}
 \|\mathsf S_y\bar\rho_\tau-\bar\rho_\tau\|_{L^1((0,\tfin)\times\Torus)}
 \leq m_{\tfin,R}^{-1}(C_{\tfin}|y|+C_{\tfin}\tau)
      +\frac{2\tfin C_{\tfin}}{\omega_f(R)}.
\end{equation}
First let $\tau\to0$ and $y\to0$ with $R$ fixed, and then let
$R\to+\infty$.  This proves \eqref{jko:eq:rho-translation-limit}.
\end{proof}

\begin{proposition}[Space-time compactness]\label{jko:prop:strong-compactness}
Every sequence $\tau_n\downarrow0$ has a subsequence and a density $\rho$ such
that, for every $\tfin>0$,
\begin{equation}\label{jko:eq:L1-compactness}
  \bar\rho_{\tau_n}\longrightarrow\rho
  \qquad\text{strongly in }L^1((0,\tfin)\times\Torus)
\end{equation}
and $\bar\rho_{\tau_n}(t)\to\rho(t)$ weakly in $\Pcal(\Torus)$, uniformly in
$t\in[0,\tfin]$.
\end{proposition}

\begin{proof}
Spatial mollification converts the estimate \eqref{jko:eq:time-modulus},
stated in the norm \eqref{jko:eq:negative-norm}, into an $L^1$
time-translation estimate.  Let $\eta_\varepsilon$ be a smooth periodic
spatial mollifier and set
$\bar\rho_\tau^\varepsilon=\eta_\varepsilon*\bar\rho_\tau$.  By
\cref{jko:lem:rho-space-translations},
\begin{equation}\label{jko:eq:mollification-error}
  \lim_{\varepsilon\downarrow0}\limsup_{n\to\infty}
  \|\bar\rho_{\tau_n}^\varepsilon-\bar\rho_{\tau_n}\|_{L^1((0,\tfin)\times\Torus)}
  =0.
\end{equation}
For a time shift $\Delta t>0$, convolution and \cref{jko:lem:time-compactness}
give, for $t,t+\Delta t\in[0,\tfin]$,
\begin{equation}\label{jko:eq:mollified-time-translate}
 \|\bar\rho_{\tau_n}^\varepsilon(t+\Delta t)
      -\bar\rho_{\tau_n}^\varepsilon(t)\|_{L^1(\Torus)}
 \leq C_{\varepsilon,\tfin}\sqrt{\Delta t+\tau_n}.
\end{equation}
Indeed, at each spatial point the convolution is the pairing of the measure
difference with a translate of $\eta_\varepsilon$, whose $C^2$ norm is
bounded by a constant depending only on $\varepsilon$.  Consequently,
\begin{align*}
 &\int_0^{\tfin-\Delta t}
 \|\bar\rho_{\tau_n}(t+\Delta t)-\bar\rho_{\tau_n}(t)\|_{L^1(\Torus)}\dd t\\
 &\qquad\leq
 2\|\bar\rho_{\tau_n}^\varepsilon-\bar\rho_{\tau_n}\|_
       {L^1((0,\tfin)\times\Torus)}
 +C_{\varepsilon,\tfin}\sqrt{\Delta t+\tau_n}.
\end{align*}
Together with \eqref{jko:eq:mollification-error}, letting first
$n\to\infty$, then $\Delta t\downarrow0$, then $\varepsilon\downarrow0$, this
gives
\[
  \lim_{\Delta t\downarrow0}\limsup_{n\to\infty}
  \int_0^{\tfin-\Delta t}
  \|\bar\rho_{\tau_n}(t+\Delta t)-\bar\rho_{\tau_n}(t)\|_{L^1(\Torus)}\dd t
  =0.
\]

Since each $\bar\rho_{\tau_n}\in L^1$, the asymptotic bounds
\eqref{jko:eq:rho-translation-limit} and the time-translation bound above give
translation estimates uniform in $n$.  With the uniform integrability of
\cref{jko:lem:rho-space-translations}, the Fr\'echet--Kolmogorov criterion
gives relative compactness in $L^1((0,\tfin)\times\Torus)$.  This proves \eqref{jko:eq:L1-compactness}.
For the uniform-in-time statement, let $\widetilde\rho_{\tau_n}$ be the
interpolant of the grid values that is affine in $t$ on each cell, in the
vector space of signed measures.  The interpolant $\widetilde\rho_{\tau_n}(t)$
remains in $\Pcal(\Torus)$, differs from $\bar\rho_{\tau_n}(t)$ by
$O(\sqrt{\tau_n})$ in \eqref{jko:eq:negative-norm}, and the proof of
\eqref{jko:eq:time-modulus}, applied also to fractional end cells, gives a
uniform one-half H\"older modulus.  Compactness of $\Pcal(\Torus)$ in that
metric and Arzel\`a--Ascoli now give uniform-in-time weak convergence of a
subsequence of $\widetilde\rho_{\tau_n}$.  The uniform weak limit equals
$\rho(t)$ for almost every $t$ by \eqref{jko:eq:L1-compactness}, and for every
$t$ by weak continuity.  Since $\widetilde\rho_{\tau_n}(t)$ and
$\bar\rho_{\tau_n}(t)$ differ by $O(\sqrt{\tau_n})$ in
\eqref{jko:eq:negative-norm}, uniformly in $t$, also
$\bar\rho_{\tau_n}(t)\to\rho(t)$ weakly, uniformly in $t\in[0,\tfin]$.  A diagonal extraction over $\tfin\in\mathbb N$ gives
one subsequence for all $\tfin>0$.
\end{proof}

\subsection{Identification of the potential and the momentum}

Fix the subsequence from \cref{jko:prop:strong-compactness} and extract a
further subsequence such that $\bar\rho_{\tau_n}\to\rho$ almost everywhere on
$(0,\tfin)\times\Torus$.  Since $\bar\rho_{\tau_n}\geq \hat c_{\tfin}$ almost
everywhere by the positivity floor \eqref{jko:eq:floor}, also
$\rho\geq \hat c_{\tfin}$ almost everywhere, and \eqref{jko:eq:fprime-growth} of
assumption \textup{(E3)} with $c=\hat c_{\tfin}$ applies to $r=\bar\rho_{\tau_n}$.
With \eqref{jko:eq:internal-energy-bound}, it gives
\begin{equation}\label{jko:eq:fprime-Lq}
  \sup_n\|f'(\bar\rho_{\tau_n})\|_{L^{q_{\hat c_{\tfin}}}((0,\tfin)\times\Torus)}<+\infty.
\end{equation}
Because $q_{\hat c_{\tfin}}>1$ and $(0,\tfin)\times\Torus$ has finite measure,
\eqref{jko:eq:fprime-Lq} shows that $(f'(\bar\rho_{\tau_n}))_n$ is uniformly
integrable.  Because $f'$ is continuous on $[\hat c_{\tfin},+\infty)$ by
\textup{(E2)}, $f'(\bar\rho_{\tau_n})\to f'(\rho)$ almost everywhere.
Vitali's theorem then gives
\begin{equation}\label{jko:eq:fprime-strong}
  f'(\bar\rho_{\tau_n})\longrightarrow f'(\rho)
  \qquad\text{strongly in }L^1((0,\tfin)\times\Torus).
\end{equation}
By \textup{(E1)}, $W\in C^2(\Torus)$.  Hence, for almost every
$t\in(0,\tfin)$,
\[
  \|W*(\bar\rho_{\tau_n}(t)-\rho(t))\|_{C^2(\Torus)}
  \leq\|W\|_{C^2(\Torus)}\,
      \|\bar\rho_{\tau_n}(t)-\rho(t)\|_{L^1(\Torus)}.
\]
Integration in $t$ and \eqref{jko:eq:L1-compactness} give
\begin{equation}\label{jko:eq:convolution-strong}
 W*\bar\rho_{\tau_n}\longrightarrow W*\rho
 \qquad\text{strongly in }L^1(0,\tfin;C^2(\Torus)).
\end{equation}

By \eqref{jko:eq:phi-H2-bound}, after extraction,
\begin{equation}\label{jko:eq:phi-weak}
  \varphi_{\tau_n}\weakto\varphi
  \qquad\text{weakly in }L^2(0,\tfin;H^2(\Torus)).
\end{equation}
Test \eqref{jko:eq:interpolated-HJ} against any smooth
$\zeta(t,x)$ such that $\int_{\Torus}\zeta(t,x)\dd x=0$ for every $t$.  For
such $\zeta$, the term
$\int_0^{\tfin}\lambda_{\tau_n}(t)\int_{\Torus}\zeta(t,x)\dd x\dd t$ is zero.  Passing to the limit using
\eqref{jko:eq:defect-bound}, \eqref{jko:eq:fprime-strong},
\eqref{jko:eq:convolution-strong}, and \eqref{jko:eq:phi-weak} gives
\[
  \int_0^{\tfin}\!\!\int_{\Torus}
  \bigl(\varphi+V_{\rm ext}+W*\rho+f'(\rho)\bigr)\zeta\dd x\dd t=0.
\]
By \eqref{jko:eq:chemical-potential}, the integrand is $g\,\zeta$ with
$g:=\varphi+\xi\in L^1((0,\tfin)\times\Torus)$.  Since this identity holds for
every smooth $\zeta$ with zero spatial mean for each $t$,
$g(t,x)=\int_{\Torus}g(t,y)\dd y$ for almost every $(t,x)$.  Hence
$\varphi+\xi$ depends only on $t$.  In particular,
\begin{equation}\label{jko:eq:hessian-identification}
  \nabla^2\varphi=-\nabla^2\xi
  \qquad\text{in distributions}.
\end{equation}
By \eqref{jko:eq:phi-weak},
$\nabla^2\varphi\in L^2((0,\tfin)\times\Torus;\R^{d\times d})$.  Hence
\eqref{jko:eq:hessian-identification} shows that $\nabla^2\xi$ is a function
in $L^2((0,\tfin)\times\Torus;\R^{d\times d})$.  Since $\rho$ is not proved to be
bounded, this does not give the weighted bound
\eqref{jko:eq:limit-regularity}; that bound is proved by
\eqref{jko:eq:B-identification} below.

It remains to identify the limit of
$Q_{\tau_n}=\bar\rho_{\tau_n}\nabla^2\varphi_{\tau_n}$ of
\eqref{jko:eq:all-interpolants}.  In this product, $\bar\rho_{\tau_n}$
converges strongly only in $L^1((0,\tfin)\times\Torus)$,
$\nabla^2\varphi_{\tau_n}$ converges only weakly in
$L^2((0,\tfin)\times\Torus;\R^{d\times d})$, and no upper bound on
$\bar\rho_{\tau_n}$ uniform in $n$ is proved.  Define
\begin{equation}\label{jko:eq:Btau}
  \mathsf J_{\tau_n}:=\sqrt{\bar\rho_{\tau_n}}\,
      \nabla^2\varphi_{\tau_n}.
\end{equation}
By \eqref{jko:eq:Btau} and \eqref{jko:eq:action-bound},
\[
  \|\mathsf J_{\tau_n}\|_{L^2((0,\tfin)\times\Torus;\R^{d\times d})}^2
  =\int_0^{\tfin}\!\!\int_{\Torus}
     \bar\rho_{\tau_n}|\nabla^2\varphi_{\tau_n}|^2\dd x\dd t
  \leq2\bigl(\Ec(\rho_0)-\inf\Ec\bigr).
\]
Hence, along a subsequence, $\mathsf J_{\tau_n}\weakto \mathsf J$ weakly in
$L^2((0,\tfin)\times\Torus;\R^{d\times d})$.  Strong $L^1$ convergence of the densities
implies
\begin{equation}\label{jko:eq:sqrt-strong}
  \sqrt{\bar\rho_{\tau_n}}\longrightarrow\sqrt\rho
  \qquad\text{strongly in }L^2((0,\tfin)\times\Torus),
\end{equation}
because $|\sqrt r-\sqrt{r'}|^2\leq|r-r'|$ for all $r,r'\geq0$.
For every $\Psi\in C([0,\tfin]\times\Torus;\R^{d\times d})$,
\begin{align*}
 \int_0^{\tfin}\!\!\int_{\Torus}\mathsf J_{\tau_n}:\Psi\dd x\dd t
 &=\int_0^{\tfin}\!\!\int_{\Torus}\nabla^2\varphi_{\tau_n}:
       \bigl(\sqrt{\bar\rho_{\tau_n}}\Psi\bigr)\dd x\dd t\\
 &\longrightarrow
   \int_0^{\tfin}\!\!\int_{\Torus}
       \nabla^2\varphi:\bigl(\sqrt\rho\Psi\bigr)\dd x\dd t,
\end{align*}
because $\nabla^2\varphi_{\tau_n}\weakto\nabla^2\varphi$ weakly in
$L^2((0,\tfin)\times\Torus;\R^{d\times d})$ by \eqref{jko:eq:phi-weak} and
$\sqrt{\bar\rho_{\tau_n}}\Psi\to\sqrt\rho\Psi$ strongly in
$L^2((0,\tfin)\times\Torus;\R^{d\times d})$ by \eqref{jko:eq:sqrt-strong}.
Since also $\mathsf J_{\tau_n}\weakto \mathsf J$, the two limits of
$\int_0^{\tfin}\!\!\int_{\Torus}\mathsf J_{\tau_n}:\Psi\dd x\dd t$ coincide for every
such $\Psi$, which gives
\begin{equation}\label{jko:eq:B-identification}
  \mathsf J=\sqrt\rho\,\nabla^2\varphi
   =-\sqrt\rho\,\nabla^2\xi.
\end{equation}
As $\mathsf J\in L^2((0,\tfin)\times\Torus;\R^{d\times d})$, this proves
\eqref{jko:eq:limit-regularity}.

Since $Q_{\tau_n}=\sqrt{\bar\rho_{\tau_n}}\mathsf J_{\tau_n}$ by
\eqref{jko:eq:all-interpolants} and \eqref{jko:eq:Btau}, every
$\Psi\in L^\infty((0,\tfin)\times\Torus;\R^{d\times d})$ satisfies
\begin{align*}
 \int_0^{\tfin}\!\!\int_{\Torus}(Q_{\tau_n}-\sqrt\rho\,\mathsf J):\Psi\dd x\dd t
 &=\int_0^{\tfin}\!\!\int_{\Torus}
    (\sqrt{\bar\rho_{\tau_n}}-\sqrt\rho)\,\mathsf J_{\tau_n}:\Psi\dd x\dd t\\
 &\quad+\int_0^{\tfin}\!\!\int_{\Torus}
    \sqrt\rho\,(\mathsf J_{\tau_n}-\mathsf J):\Psi\dd x\dd t.
\end{align*}
The absolute value of the first term on the right is at most
$\|\sqrt{\bar\rho_{\tau_n}}-\sqrt\rho\|_{L^2}\|\mathsf J_{\tau_n}\|_{L^2}
\|\Psi\|_{L^\infty}$, which tends to zero by \eqref{jko:eq:sqrt-strong}
because $\|\mathsf J_{\tau_n}\|_{L^2}$ is bounded.  The second term tends to zero
because $\sqrt\rho\,\Psi\in L^2((0,\tfin)\times\Torus;\R^{d\times d})$ and
$\mathsf J_{\tau_n}\weakto \mathsf J$.  By \eqref{jko:eq:B-identification},
$\sqrt\rho\,\mathsf J=\rho\nabla^2\varphi$.  Hence
\begin{equation}\label{jko:eq:flux-identification}
  Q_{\tau_n}
  =\sqrt{\bar\rho_{\tau_n}}\mathsf J_{\tau_n}
  \weakto
  \rho\nabla^2\varphi
  =-\rho\nabla^2\xi
  \qquad\text{weakly in }L^1((0,\tfin)\times\Torus;\R^{d\times d}).
\end{equation}

\subsection{Energy inequality}

We prove \eqref{jko:eq:limit-energy-inequality}.  Let $t>0$ and let $m\geq1$
be the integer such that $t\in((m-1)\tau_n,m\tau_n]$.  Then
$\bar\rho_{\tau_n}(t)=\rho_m^{\tau_n}$ and
\[
  \int_0^t\!\!\int_{\Torus}|\mathsf J_{\tau_n}|^2\dd x\dd s
  \leq\int_0^{m\tau_n}\!\!\int_{\Torus}|\mathsf J_{\tau_n}|^2\dd x\dd s.
\]
The proof of \cref{jko:lem:action-estimate} gives, for every $k\geq0$, the
one-step inequality
$\frac1{\tau_n}\mathsf G_{\fr}(\rho_k^{\tau_n},\rho_{k+1}^{\tau_n})
+\Ec(\rho_{k+1}^{\tau_n})\leq\Ec(\rho_k^{\tau_n})$ and the identity
\[
  \frac1{\tau_n}\mathsf G_{\fr}(\rho_k^{\tau_n},\rho_{k+1}^{\tau_n})
  =\frac{\tau_n}2\int_{\Torus}\rho_{k+1}^{\tau_n}
     |\nabla^2\varphi_{k+1}^{\tau_n}|^2\dd x
  =\frac12\int_{k\tau_n}^{(k+1)\tau_n}\!\!\int_{\Torus}
     |\mathsf J_{\tau_n}|^2\dd x\dd s,
\]
where the last equality follows from \eqref{jko:eq:interpolation},
\eqref{jko:eq:all-interpolants} and \eqref{jko:eq:Btau}.  Summing the one-step inequality for $k=0,\dots,m-1$
gives
\[
  \Ec(\rho_m^{\tau_n})
  +\frac12\int_0^{m\tau_n}\!\!\int_{\Torus}|\mathsf J_{\tau_n}|^2\dd x\dd s
  \leq\Ec(\rho_0).
\]
Hence, for every $t>0$,
\[
 \Ec(\bar\rho_{\tau_n}(t))
 +\frac12\int_0^t\!\!\int_{\Torus}|\mathsf J_{\tau_n}|^2\dd x\dd s
 \leq\Ec(\rho_0).
\]
After a further extraction, \eqref{jko:eq:L1-compactness} gives
$\bar\rho_{\tau_n}(t)\to\rho(t)$ in $L^1(\Torus)$ for almost every
$t\in(0,\tfin)$.  Fix such a $t$.
\begin{enumerate}
\item[\textup{(i)}] $\int_{\Torus}f(\rho(t))\dd x\leq\liminf_n
  \int_{\Torus}f(\bar\rho_{\tau_n}(t))\dd x$, since $f$ is convex and lower
  semicontinuous by \textup{(E2)}.
\item[\textup{(ii)}] $\int_{\Torus}V_{\rm ext}\bar\rho_{\tau_n}(t)\dd x
  +\frac12\int_{\Torus}(W*\bar\rho_{\tau_n}(t))\bar\rho_{\tau_n}(t)\dd x$
  converges to the same expression at $\rho(t)$, since $V_{\rm ext}$ and $W$
  are bounded by \textup{(E1)} and $\bar\rho_{\tau_n}(t)\to\rho(t)$ in
  $L^1(\Torus)$.
\item[\textup{(iii)}] The restrictions of $\mathsf J_{\tau_n}$ to $(0,t)\times\Torus$
  converge weakly in $L^2((0,t)\times\Torus;\R^{d\times d})$ to the
  restriction of $\mathsf J$.  Hence
  $\int_0^t\!\!\int_{\Torus}|\mathsf J|^2\dd x\dd s\leq\liminf_n
  \int_0^t\!\!\int_{\Torus}|\mathsf J_{\tau_n}|^2\dd x\dd s$, and
  $|\mathsf J|^2=\rho|\nabla^2\xi|^2$ by \eqref{jko:eq:B-identification}.
\end{enumerate}
Passing to the limit $n\to\infty$ with \textup{(i)--(iii)} in the inequality
for $\Ec(\bar\rho_{\tau_n}(t))$ above gives
\eqref{jko:eq:limit-energy-inequality} for almost every $t\in(0,\tfin)$.
Since $\tfin>0$ is arbitrary, \eqref{jko:eq:limit-energy-inequality} holds for
almost every $t>0$.

\subsection{Passage to the limit and initial trace}

Let $\zeta\in C_c^\infty([0,+\infty)\times\Torus)$ and choose $\tfin>0$ with
$\supp\zeta\subset[0,\tfin)\times\Torus$.  Set
$\bar\rho_{\tau_n}(t):=\rho_0$ for $t\leq0$, so that
$\underline\rho_{\tau_n}(t)=\bar\rho_{\tau_n}(t-\tau_n)$ for every $t>0$.  The
change of variables $t\mapsto t+\tau_n$ in the term that contains
$\underline\rho_{\tau_n}$ gives
\begin{align*}
 \int_0^\infty\!\!\int
   \frac{\bar\rho_{\tau_n}-\underline\rho_{\tau_n}}{\tau_n}\,
       \zeta\dd x\dd t
 &=-\int_0^\infty\!\!\int\bar\rho_{\tau_n}(t)\,
     \frac{\zeta(t+\tau_n,\cdot)-\zeta(t,\cdot)}{\tau_n}\dd x\dd t\\
 &\quad-\frac1{\tau_n}\int_0^{\tau_n}\!\!\int_{\Torus}
     \rho_0\,\zeta(t,\cdot)\dd x\dd t.
\end{align*}
The difference quotient converges to $\partial_t\zeta$ uniformly, and
$\bar\rho_{\tau_n}(t)$ has unit mass for every $t$.  With
\eqref{jko:eq:L1-compactness} on $(0,\tfin)\times\Torus$, the first term on
the right therefore tends to
$-\int_0^\infty\!\!\int\rho\,\partial_t\zeta\dd x\dd t$.  The last term tends
to $-\int_{\Torus}\rho_0\zeta(0,\cdot)\dd x$ because $\zeta$ is continuous.
Hence
\begin{equation}\label{jko:eq:discrete-weak-time}
 \int_0^\infty\!\!\int
   \frac{\bar\rho_{\tau_n}-\underline\rho_{\tau_n}}{\tau_n}\,
       \zeta\dd x\dd t
 \longrightarrow
 -\int_0^\infty\!\!\int\rho\,\partial_t\zeta\dd x\dd t
 -\int_{\Torus}\rho_0\zeta(0,\cdot)\dd x.
\end{equation}
On the other hand, \eqref{jko:eq:flux-identification} with
$\Psi=\nabla^2\zeta$ on $(0,\tfin)\times\Torus$ gives
\[
 \int_0^\infty\!\!\int Q_{\tau_n}:\nabla^2\zeta\dd x\dd t
 \longrightarrow
 -\int_0^\infty\!\!\int
      \rho\nabla^2\xi:\nabla^2\zeta\dd x\dd t.
\]
For every sufficiently large $n$, the left-hand side of \eqref{jko:eq:discrete-weak-time} and
$\int_0^\infty\!\!\int Q_{\tau_n}:\nabla^2\zeta\dd x\dd t$ are equal, by
\eqref{jko:eq:discrete-continuity} tested with $\zeta$ and two integrations
by parts in $x$ on $\Torus$.  Equating their limits gives
\[
  \int_0^\infty\!\!\int_{\Torus}\rho\,\partial_t\zeta\dd x\dd t
  +\int_{\Torus}\rho_0\zeta(0,\cdot)\dd x
  =\int_0^\infty\!\!\int_{\Torus}\rho\nabla^2\xi:\nabla^2\zeta\dd x\dd t,
\]
which is \eqref{jko:eq:weak-formulation}.

The time modulus \eqref{jko:eq:time-modulus} with $s=0$ also gives, for
$t\in[0,\tfin]$,
\[
  \|\bar\rho_{\tau_n}(t)-\rho_0\|_{-2,\infty}
  \leq C_{\tfin}\sqrt{t+\tau_n}.
\]
Letting $n\to\infty$ with \eqref{jko:eq:weak-uniform-main}, and using that
$\|\cdot\|_{-2,\infty}$ is a supremum of weakly continuous linear functionals
and hence weakly lower semicontinuous, gives
$\|\rho(t)-\rho_0\|_{-2,\infty}\leq C_{\tfin}\sqrt t$ for $t\in[0,\tfin]$.
Since this norm metrizes weak convergence on $\Pcal(\Torus)$, the curve
$t\mapsto\rho(t)$ of \eqref{jko:eq:weak-uniform-main} satisfies
$\rho(t)\to\rho_0$ weakly as $t\downarrow0$.

Suppose that \eqref{jko:eq:limit-pde} has at most one weak solution, in the
sense of \eqref{jko:eq:weak-formulation}, that satisfies
$f'(\rho)\in L^1((0,\tfin)\times\Torus)$ and
\eqref{jko:eq:limit-regularity} for every $\tfin>0$, and
\eqref{jko:eq:limit-energy-inequality}; this is the class of
\cref{jko:thm:main}.  Every sequence $\tau_n\downarrow0$ has a subsequence
along which \eqref{jko:eq:strong-main} and \eqref{jko:eq:weak-uniform-main}
hold.  By the proof above, the limit of this subsequence satisfies
$f'(\rho)\in L^1((0,\tfin)\times\Torus)$ for every $\tfin>0$ by
\eqref{jko:eq:fprime-strong},
\eqref{jko:eq:limit-regularity}, \eqref{jko:eq:weak-formulation} and
\eqref{jko:eq:limit-energy-inequality}, so it is this unique solution $\rho$.  A sequence in a metric space whose every subsequence has a further
subsequence converging to the same limit converges to that limit.  Hence
\eqref{jko:eq:strong-main} and \eqref{jko:eq:weak-uniform-main} hold for the
whole family as $\tau\to0$.  This completes the proof of
\cref{jko:thm:main}.

\subsection{The estimates preserved by nonnegative momentum}
\label{jko:app:positive-compactness}
\begin{proof}[Proof of \cref{jko:prop:positive-compactness}]
The proof of \cref{jko:prop:existence} applies to the joint problem
\eqref{jko:eq:positive-scheme}: the competitor $(\rho_k^{\tau,+},0)$
has finite value, the action controls the total variation of momentum,
and the affine conservation constraint is weakly closed. The additional
cone of positive-semidefinite measures is also weakly closed, since
$\int\chi\,e^T\dd\mathsf m\,e\ge0$ passes to the limit for every
continuous $\chi\ge0$ and every $e\in\R^d$. The energy lower
semicontinuity and superlinearity argument in that proof therefore
gives an attained step with an absolutely continuous optimized distribution.

Comparison with this zero-momentum competitor gives
\eqref{jko:eq:positive-energy-bound}. Write
$C_0=2(\Ec(\rho_0)-\inf\Ec)/\kappa$. Set $B_\tau:=B_{k+1}^\tau$ on
$(k\tau,(k+1)\tau]$. Then, on any finite interval $(0,T)$,
\[
 \int_0^T\!\!\int|B_\tau|_F^2\bar\rho_\tau^+\dd x\dd t\le C_0,
 \qquad
 \int_0^T|Q_\tau^+(t)|(\Torus)\dd t\le\sqrt{TC_0}.
\]
The proof of \cref{jko:lem:time-compactness} uses only these two
bounds and the conservation equation. It gives
\[
 \|\bar\rho_\tau^+(t)-\bar\rho_\tau^+(s)\|_{-2,\infty}
 \le\sqrt{C_0}\sqrt{|t-s|+2\tau}.
\]
For the affine marginal interpolation one can replace the last factor
by $\sqrt{|t-s|}$: integrate its exact conservation equation against
a test function and apply Cauchy--Schwarz to the momentum rate on
$(s,t)$. Its distance to $\bar\rho_\tau^+$ is $O(\sqrt\tau)$.
Arzel\`a--Ascoli on the compact space $\Pcal(\Torus)$ gives a
uniformly weakly convergent subsequence, and weak compactness of
$\dd t\otimes Q_\tau^+$ gives a limiting matrix measure $Q$.
A diagonal extraction gives one subsequence on all finite intervals.

Set $R_\tau=\dd t\otimes\bar\rho_\tau^+\dd x$ and
$R=\dd t\otimes\rho_t\dd x$. At this stage the limiting time slices
are probability measures; energy lower semicontinuity and
\eqref{jko:eq:positive-energy-bound} show at every time that their
internal energy is finite, hence they have densities by superlinearity
of $f$. Uniform weak convergence implies $R_\tau\rightharpoonup R$.
The matrix version of \cref{geo:lem:space-time-action}, equivalently
the dual formula \eqref{jko:eq:perspective-dual}, implies
\[
 Q=BR,\qquad B\in L^2(R;\R^{d\times d}),\qquad B\succeq0.
\]
Here $R=\rho_t(x)\dd t\dd x$ is the space--time measure, and
$L^2(R;\R^{d\times d})$ means that, on every finite interval $(0,T)$,
\[
 \int_0^T\!\!\int_{\Torus}|B(t,x)|_F^2\rho_t(x)\dd x\dd t<\infty.
\]
The positivity assertion holds $R$-almost everywhere and follows from
closure of the positive cone. Passing
to the limit in the affine interpolation's exact conservation equation
gives $\partial_t\rho=\DivDiv(\rho B)$ and the prescribed initial
trace, exactly as in the signed proof.

Finally, for fixed $t$, sum the discrete energy inequalities through
the cell containing $t$, and retain only the action on $(0,t)$.
The endpoint energy is lower semicontinuous under the uniform weak
convergence. Lower semicontinuity of the perspective action on
$(0,t)\times\Torus$ gives the other term of
\eqref{jko:eq:positive-limit-energy}. This can be seen first with
continuous tests compactly supported in this open time interval and
then by exhaustion. The limiting action measure is absolutely continuous in time
through $R$, so no boundary atom is lost. This proves the inequality
for every $t$, without the density floor or a dual-potential estimate.
\end{proof}

\begin{remark}[Why the signed $H^2$ duality does not extend]
\label{jko:rem:positive-dual-obstruction}
On $\T$, take the initial density $\mu=1$ and the final density
$\rho=1+\varepsilon\cos(2\pi x)$, $0<\varepsilon<1$.
Every momentum solving $m''=\rho-\mu$ has density
$m(x)=C-a\cos(2\pi x)$, where $a=\varepsilon/(4\pi^2)$.
Nonnegativity requires $C\ge a$, and the quadratic action increases
with $C$ on this interval. Its unique positive minimizer is therefore
$m=a(1-\cos(2\pi x))$, which vanishes only at $x=0$.
If an $H^2(\T)$ dual maximizer satisfied
$m=\rho(\Phi'')_+$, then $\Phi''=m/\rho>0$ almost everywhere,
contradicting $\int_\T\Phi''=0$.
In the relaxed formulation the second derivative is instead
\[
 \Phi''=\frac m\rho\dd x-
          \left(\int_\T\frac m\rho\dd x\right)\delta_0.
\]
This has zero total mass, its positive part is $m/\rho$, and its
negative part is supported on the contact set $\{m=0\}$. It explains
why a strictly positive density does not supply the signed proof's
$H^2$ bound.
\end{remark}

\section{Proofs of martingale convergence}\label{mart:app:proofs}
This appendix proves the results of
\cref{mart:sec:setting,mart:sec:limits,mart:g:sec:m4-theorem}.
\subsection{The positive Eulerian action and the coupling projection}
\begin{proof}[Proof of \cref{mart:prop:positive-eulerian}]
Finite action gives $m_t=v_t\mu_t$ with $v_t\ge0$. Setting $a_t=2v_t$
gives \eqref{mart:eq:momentum-covariance} and the weak Fokker--Planck
equation. The converse follows by the same substitution.
Moreover,
\[
 \int_0^{\tfin}\!\int a_t\dd\mu_t\dd t
 \le\tfin^{1/2}
       \left(\int_0^{\tfin}\!\int a_t^2\dd\mu_t\dd t\right)^{1/2}<\infty.
\]
The superposition principle \cite[Theorem~2.5]{Tre16} yields a continuous
solution of the martingale problem with marginals $\mu_t$ and bracket
density $a_t(X_t)$. Its integrable expected bracket and the initial
second moment make it a square-integrable martingale.
Conversely, if $\dd\langle X\rangle_t=\alpha_t\dd t$ with
$\mathbb E\int\alpha_t^2\dd t<\infty$, It\^o's formula and conditional
projection give the Fokker--Planck equation with
$a_t(x)=\mathbb E[\alpha_t\mid X_t=x]$. Conditional Jensen gives
\[
 \int_0^{\tfin}\!\int a_t^2\dd\mu_t\dd t
 \le\mathbb E\int_0^{\tfin}\alpha_t^2\dd t.
\]
Superposition attains equality for a prescribed Eulerian pair, proving
the equality of the three infima.
\end{proof}

\begin{proof}[Proof of \cref{mart:prop:coupling-momentum}]
The measure in \eqref{mart:eq:segment-momentum} is nonnegative and has
mass $\mathbb E\xi^2/2<\infty$. For $\varphi\in C_c^\infty(\R)$,
Taylor's formula with integral remainder gives
\begin{align*}
 \mathbb E[\varphi(Y)-\varphi(X)]
 &=\mathbb E[\varphi'(X)\xi]
  +\mathbb E\left[\xi^2\int_0^1(1-s)\varphi''(X+s\xi)\dd s\right]\\
 &=\int\varphi''\dd m_\pi.
\end{align*}
The first term vanishes by the martingale constraint. This proves the
distributional equation. The difference of two finite-measure
solutions has second derivative zero, hence is an affine distribution;
the only such finite measure on $\R$ is zero. This proves uniqueness.
If a finite nonnegative solution $m$ exists, test the equation with
smooth convex functions of bounded derivative, using cutoffs. The
errors tend to zero because $m$ is finite and the endpoints have finite
second moments. Thus $\int\varphi\dd\mu\le\int\varphi\dd\sigma$ for
these convex functions. Approximation gives convex order, and
Strassen's theorem \cite{Str65} gives a martingale coupling. The other
direction follows from the projection just constructed.
For the example in \cref{mart:rem:discrete-variables}, direct evaluation
of \eqref{mart:eq:segment-momentum} gives the triangular density
$(b-|z|)_+/2$, which is not absolutely continuous with respect to the
atomic midpoint measure.
For the action comparison following \cref{mart:prop:coupling-momentum},
consider the probability measure
$\widehat\pi(\dd x,\dd y,\dd s)=2(1-s)\pi(\dd x,\dd y)\dd s$
on $\R^2\times[0,1]$ and set $Z=X+s\xi$.
Then $\lambda_\pi$ is the distribution of $Z$ and
\[
 \frac{\dd m_\pi}{\dd\lambda_\pi}(z)
 =\frac12\mathbb E_{\widehat\pi}[\xi^2\mid Z=z].
\]
Conditional Jensen gives
\[
 \Ac(\lambda_\pi,m_\pi)
 =\frac18\mathbb E_{\widehat\pi}
       \left[\left(\mathbb E_{\widehat\pi}[\xi^2\mid Z]\right)^2\right]
 \le\frac18\mathbb E_\pi\xi^4.
\]
The bound is trivial when the right-hand side is infinite.
\end{proof}

\subsection{Convergence of the midpoint scheme with nonnegative momentum}
\label{mart:app:positive-midpoint}
We deduce \cref{mart:cor:positive-midpoint} from the common geodesic
convergence argument and its sign-preserving recovery
(\cref{geo:lem:time-average-recovery}). The state space is now $\R$,
so we verify marginal and momentum tightness explicitly. Both endpoint
laws are fixed; no fourth-moment coupling estimate is used.

\begin{proof}[Proof of \cref{mart:cor:positive-midpoint}]
\emph{Compactness of the marginal laws.}
Consider discrete admissible pairs with actions bounded by $C$, and
write $h=1/N$,
\[
 q_h^k=\frac{\mu_h^k+\mu_h^{k+1}}2,\qquad
 e_h^k=\frac1h\Ac(q_h^k,m_h^k),\qquad
 \sum_ke_h^k\le C.
\]
By \cref{mart:prop:coupling-momentum}, positivity of each momentum
implies
$\mu\le_c\mu_h^k\le_c\sigma$. The same holds for the affine marginal
interpolation $\mu_h(t)$ and for the midpoint weights.
This convex-order interval is narrowly compact in $\Ptwo(\R)$.
For completeness, the convex function $(|x|-r)_+^2$ gives
\begin{equation}\label{mart:eq:positive-uniform-tails}
 \sup_{\lambda\le_c\sigma}
 \int_{\{|x|>2r\}}x^2\lambda(\dd x)
 \le4\int (|x|-r)_+^2\sigma(\dd x)\longrightarrow0.
\end{equation}
Thus the family is tight, its squared positions are uniformly
integrable, and convex order passes to narrow limits by approximation
with convex functions of at most quadratic growth.

The Cauchy--Schwarz inequality gives
$m_h^k(\R)\le\sqrt{2he_h^k}$. For
$\varphi\in C_b^2(\R)$, summing the conservation equations over the
cells met by $[s,t]$ yields
\begin{equation}\label{mart:eq:positive-marginal-modulus}
 |\langle\varphi,\mu_h(t)-\mu_h(s)\rangle|
 \le\|\varphi''\|_\infty
             \sqrt{2C(|t-s|+2h)}.
\end{equation}
On the uniformly tight set above, the metric obtained by taking the
supremum over $\|\varphi\|_{C_b^2}\le1$ metrizes narrow convergence.
This follows by smooth approximation of bounded Lipschitz tests and
compactness on bounded spatial intervals. Arzel\`a--Ascoli therefore
gives a subsequence converging uniformly narrowly to a continuous
curve $\mu_t\in\Ptwo(\R)$, with the prescribed endpoints.
Every $\mu_t$ remains in the same convex-order interval.

\emph{Momentum compactness and the lower bound.}
Define space--time measures
\[
 R_h=\dd t\otimes q_h^k,\qquad
 W_h=\dd t\otimes(m_h^k/h)
 \quad\text{on }[t_k,t_{k+1}).
\]
Their masses satisfy $R_h([0,1]\times\R)=1$ and
$W_h([0,1]\times\R)=\sum_km_h^k(\R)\le\sqrt{2C}$.
For every spatial Borel set $E$, Cauchy--Schwarz also gives
\[
 W_h([0,1]\times E)
 \le\left(2C\,R_h([0,1]\times E)\right)^{1/2}.
\]
The marginal tightness thus implies tightness of both space--time
measures. After extraction, $W_h$ converges weakly to a finite
nonnegative measure $W$. The midpoint and affine marginal
interpolations differ by $O(\sqrt h)$ against the tests in
\eqref{mart:eq:positive-marginal-modulus}; hence
$R_h\rightharpoonup R=\dd t\otimes\mu_t$.

The perspective dual formula of
\cref{geo:lem:space-time-action} applies on $[0,1]\times\R$
with continuous compactly supported tests. The same completing-square
argument and density in $L^2(R)$ prove it on this noncompact space.
It gives lower semicontinuity and
\[
 \Ac(R,W)
 \le\liminf_{h\to0}\Ac(R_h,W_h)
 =\liminf_{h\to0}\sum_ke_h^k.
\]
Thus $W=vR$, with $v\ge0$ and $v\in L^2(R)$.
Set $m_t=v(t,\cdot)\mu_t$. The affine marginal interpolation and the
cellwise momentum rate satisfy the conservation equation exactly.
Passing to the limit against compactly supported smooth tests gives
$\partial_t\mu_t=\partial_{xx}m_t$. Consequently
\begin{equation}\label{mart:eq:positive-liminf}
 \mathsf G^+_{[0,1]}(\mu,\sigma)
 \le\int_0^1\Ac(\mu_t,m_t)\dd t
 \le\liminf_{h\to0}\sum_ke_h^k.
\end{equation}
This proves compactness and the lower bound. If the continuous value
is infinite, it also rules out any bounded subsequence of discrete
minimum values.

\emph{Recovery and minimum values.}
Every finite-action positive Eulerian competitor has a
square-integrable martingale representation by
\cref{mart:prop:positive-eulerian}. Its second moments are bounded by
the terminal second moment. The real-line version of
\cref{geo:lem:time-average-recovery} therefore gives admissible
$\Ptwo(\R)$ recovery nodes with nonnegative momenta and exact
endpoints. Apply it with
$\varepsilon_h=\sqrt h$ for all sufficiently small $h$.
The overhead in \eqref{geo:eq:time-average-bound} tends to one.
Taking the infimum over continuum competitors gives
$\limsup_N\mathsf G_N^+\le\mathsf G^+_{[0,1]}$.
Together with \eqref{mart:eq:positive-liminf}, this proves convergence
of the minimum values, including infinite values.
The equality with $\MB_1^2/8$ is
\cref{mart:prop:positive-eulerian}.

\emph{Attainment and minimizing pairs.}
If the continuous value is finite, recovery provides finite discrete
competitors for every sufficiently large $N$. For a fixed such $N$,
take a minimizing sequence. Its nodes lie in the compact convex-order
interval, and its finitely many momentum measures have uniformly
bounded mass and are tight by the same Cauchy--Schwarz estimate.
Weak compactness, closure of the linear equations and of the positive
cone, and lower semicontinuity give a discrete minimizer.
Apply the preceding compactness argument to asymptotically minimizing
discrete pairs as $N\to\infty$. Equality of the limiting minimum
values forces every extracted limiting pair to minimize the continuous
action. This also proves continuous attainment.
Finally \cref{mart:prop:positive-eulerian} supplies a continuous
square-integrable martingale representation of each limiting
minimizing pair, with the stated endpoint laws and action.
\end{proof}

\subsection{Two dimension-free fourth-moment estimates}\label{mart:app:moments}
\label{mart:c:sec:bdg}
These are explicit fourth-moment versions of the martingale inequalities
of \cite{Bur73}.

\begin{lemma}[Discrete estimate]
\label{mart:c:lem:discrete-bdg}
Let $(M_j)_{0\le j\le n}$ be an $\R^d$-valued martingale with $M_0=0$ and
increments ${B}_j:=M_j-M_{j-1}\in L^4$, and let
$S:=\sqrt{\sum_{j=1}^n|{B}_j|^2}$. Then
\begin{equation}
  \mathbb E|M_n|^4\le C_*\,\mathbb E S^4,
  \qquad
  C_*:=\Bigl(\frac{b+\sqrt{b^2+4}}2\Bigr)^2,
  \quad
  b:=\frac{224}9 .
  \label{mart:c:eq:discrete-bdg}
\end{equation}
\end{lemma}

\begin{proof}
Expanding $|M_j|^2=|M_{j-1}|^2+2\langle M_{j-1},{B}_j\rangle+|{B}_j|^2$ and
using $\mathbb E\bigl(|M_{j-1}|^2\langle M_{j-1},{B}_j\rangle\bigr)=0$, we obtain
\begin{equation}
  \mathbb E|M_j|^4-\mathbb E|M_{j-1}|^4
  \le\mathbb E|{B}_j|^4+6\,\mathbb E\bigl(|M_{j-1}|^2|{B}_j|^2\bigr)
  +4\,\mathbb E\bigl(|M_{j-1}||{B}_j|^3\bigr).
\end{equation}
Let $M^*:=\max_j|M_j|$. Since $|{B}_j|\le2M^*$, summing over $j$ gives
$\mathbb E|M_n|^4\le\mathbb E S^4+14\,\mathbb E(M^{*2}S^2)$. Cauchy--Schwarz and Doob's
inequality $\mathbb E M^{*4}\le(4/3)^4\mathbb E|M_n|^4$ yield
\begin{equation}
  \mathbb E|M_n|^4\le\mathbb E S^4+b\,\bigl(\mathbb E|M_n|^4\bigr)^{1/2}\bigl(\mathbb E S^4\bigr)^{1/2},
\end{equation}
which gives \eqref{mart:c:eq:discrete-bdg} since $\mathbb E|M_n|^4<\infty$.
\end{proof}

\begin{lemma}[Continuous estimate]
\label{mart:c:lem:continuous-bdg}
Let $M$ be an $\R^d$-valued continuous martingale for $(\cG_t)$, let
$\sigma\le\sigma'$ be bounded stopping times and $B\in\cG_\sigma$. Then
\begin{equation}
  \mathbb E\bigl(\one_B|M_{\sigma'}-M_\sigma|^4\bigr)
  \le C_4\,\mathbb E\bigl(\one_B(\tr\langle M\rangle_{\sigma'}
  -\tr\langle M\rangle_\sigma)^2\bigr),
  \qquad
  C_4:=\Bigl(\frac{32}3\Bigr)^2 .
  \label{mart:c:eq:continuous-bdg}
\end{equation}
\end{lemma}

\begin{proof}
Consider the stopped process
$t\mapsto\one_B(M_{(\sigma+t)\wedge\sigma'}-M_\sigma)$ with its stopped
filtration. Its quadratic variation is
$\one_B(\langle M\rangle_{(\sigma+t)\wedge\sigma'}-
\langle M\rangle_\sigma)$. Thus it is enough to treat the case $B=\Omega$,
$\sigma=0$ and $M_0=0$. Let $T$ be a stopping time such that $M^T$ is
bounded. It\^o's formula gives
\begin{equation}
  \dd|M|^4
  =4|M|^2\langle M,\dd M\rangle
  +2|M|^2\tr\dd\langle M\rangle
  +4\langle M,\dd\langle M\rangle M\rangle .
\end{equation}
With $M^*:=\sup_{t\le\sigma'\wedge T}|M_t|$, Cauchy--Schwarz and Doob's
inequality yield
\begin{equation}
  \mathbb E|M_{\sigma'\wedge T}|^4
  \le6\,\mathbb E\bigl(M^{*2}\tr\langle M\rangle_{\sigma'}\bigr)
  \le6\Bigl(\frac43\Bigr)^2
  \bigl(\mathbb E|M_{\sigma'\wedge T}|^4\bigr)^{1/2}
  \bigl(\mathbb E(\tr\langle M\rangle_{\sigma'})^2\bigr)^{1/2}.
\end{equation}
Dividing, letting $T\uparrow\infty$ and using Fatou's lemma gives
\eqref{mart:c:eq:continuous-bdg}.
\end{proof}

\begin{corollary}[Block estimate]
\label{mart:c:cor:block}
Let $X$ be a discrete martingale and $0\le i<j\le N$. Then
\begin{equation}
  \mathbb E|X_j-X_i|^4\le C_*\,(j-i)\sum_{k=i}^{j-1}\mathbb E|\xi_k|^4 .
  \label{mart:c:eq:block}
\end{equation}
\end{corollary}

\begin{proof}
Apply \cref{mart:c:lem:discrete-bdg} to $M_m:=X_{i+m}-X_i$ and use
$\bigl(\sum_{k=i}^{j-1}|\xi_k|^2\bigr)^2\le(j-i)\sum_{k=i}^{j-1}|\xi_k|^4$.
\end{proof}

\subsection{Proof of \texorpdfstring{\cref{mart:c:thm:compactness}}{the martingale result}}
\begin{proof}

Fix $n$ and write $N=N_n$, $\tau=\tau_n$, $X=X^n$ and $\fA=\fA_N(X)$, so that
$\sum_k\mathbb E|\xi_k|^4=\tau\fA$.

(i) First let $2^m\le N$. For $0\le j<2^m$ set $a_j:=\lceil jN2^{-m}\rceil$
and $b_j:=\lfloor(j+1)N2^{-m}\rfloor$. Then $a_j\le b_j$, the blocks
$\{a_j,\dots,b_j-1\}$ are disjoint and contain at most $N2^{-m}$ indices,
and the increment of $\widetilde X$ over $[j2^{-m},(j+1)2^{-m}]$ differs from
$X_{b_j}-X_{a_j}$ by at most $2\max_k|\xi_k|$. Since the fourth power of a
maximum is bounded by the sum of fourth powers, \eqref{mart:c:eq:block} gives
\begin{equation}
  \mathbb E\Delta_m(\widetilde X)^4
  \le8\sum_j\mathbb E|X_{b_j}-X_{a_j}|^4+128\sum_k\mathbb E|\xi_k|^4
  \le(8C_*+128)\,2^{-m}\fA .
  \label{mart:c:eq:coarse-dyadic}
\end{equation}
If $2^m>N$, every level-$m$ dyadic interval meets at most two grid cells, on
which $\widetilde X$ is affine with slope at most $N\max_k|\xi_k|$. Hence
$\Delta_m(\widetilde X)\le2^{-m}N\max_k|\xi_k|$ and
$\mathbb E\Delta_m(\widetilde X)^4\le2^{-4m}N^3\fA\le2^{-m}\fA$. Minkowski's
inequality now gives
\begin{equation}
  \Bigl\|\sum_{m\ge0}2^{\alpha m}\Delta_m(\widetilde X)\Bigr\|_{L^4}
  \le\bigl((8C_*+128)\Lambda\bigr)^{1/4}
  \sum_{m\ge0}2^{-(1/4-\alpha)m}<\infty .
\end{equation}
The dyadic chaining inequality bounds the $C^\alpha$ seminorm by the sum in
\eqref{mart:c:eq:holder}, up to a constant depending only on $\alpha$. Tightness
follows from Arzel\`a--Ascoli and \eqref{mart:c:eq:uniform-action}.

(ii) By orthogonality of martingale increments and Cauchy--Schwarz,
\begin{equation}
  \mathbb E|X_k|^2=\mathbb E|X_0|^2+\sum_{j<k}\mathbb E|\xi_j|^2\le\mathbb E|X_0|^2+\fA^{1/2},
\end{equation}
so the variables $\widetilde X^n_t$ are uniformly integrable. Moreover,
\begin{equation}
  \sup_{t\in[0,1]}\bigl|\widetilde X_t-X_{\lfloor t/\tau\rfloor}\bigr|
  \le\max_k|\xi_k|,
  \qquad
  \mathbb E\max_k|\xi_k|^4\le\tau\fA .
  \label{mart:c:eq:step-linear}
\end{equation}
Let $0\le s_1\le\dots\le s_p\le s<t\le1$ and let $g$ be bounded and
Lipschitz. The martingale property at grid times gives
\begin{equation}
  \mathbb E\bigl[(X_{\lfloor t/\tau\rfloor}-X_{\lfloor s/\tau\rfloor})\,
  g(X_{\lfloor s_1/\tau\rfloor},\dots,X_{\lfloor s_p/\tau\rfloor})\bigr]=0 .
\end{equation}
By \eqref{mart:c:eq:step-linear} and uniform integrability, every limit satisfies
$\mathbb E\bigl[(X_t-X_s)\,g(X_{s_1},\dots,X_{s_p})\bigr]=0$. The same estimate
gives $\bar\mu^{\,n_i}_t\to\Law(X_t)$.

(iii) By conditional Jensen, $(\tr Q_k)^2\le\mathbb E(|\xi_k|^4\mid X_k)$, hence
\begin{equation}
  \int_0^1\!\int(\tr A^N_t)^2\dd\bar\mu_t\dd t
  =\frac1\tau\sum_k\int(\tr Q_k)^2\dd\mu_k
  \le\fA_N(X).
  \label{mart:c:eq:discrete-action-bound}
\end{equation}
Set $\lambda^N:=\bar\mu_t\dd t$. Along the subsequence under consideration,
$\lambda^N$ converges narrowly to $\lambda:=\mu_t\dd t$. By
Cauchy--Schwarz and the tightness of $(\lambda^N)$, the positive
semidefinite matrix measures $B^N:=A^N\lambda^N$ are tight on
$[0,1]\times\R^d$, have total variation at most $\Lambda^{1/2}$ and give
mass at most $|I|^{1/2}\Lambda^{1/2}$ to $I\times\R^d$ for every interval
$I$. After extracting a subsequence realizing the $\liminf$ in
\eqref{mart:c:eq:action-lsc}, they converge narrowly to a positive semidefinite
matrix measure $B$. Set $q^N:=\tr B^N$ and $q:=\tr B$. The dual formula
\begin{equation}
  \int\left|\frac{\dd q}{\dd\lambda}\right|^2\dd\lambda
  =\sup_{f\in C_c}\left(2\int f\dd q-\int f^2\dd\lambda\right)
\end{equation}
(with the left-hand side equal to $+\infty$ when $q\not\ll\lambda$) shows
that $q\ll\lambda$ and gives the required lower bound. Since $B$ is
positive semidefinite and $|B|\le q$, it follows that $B=A\lambda$ for a
positive semidefinite Borel map $A$, with $\tr A=\dd q/\dd\lambda$. This
proves \eqref{mart:c:eq:action-lsc}. For $\varphi\in C_c^3(\R^d)$, Taylor's formula
and the martingale property give
\begin{equation}
  \int\varphi\dd\mu_{k+1}-\int\varphi\dd\mu_k
  =\frac12\int D^2\varphi:Q_k\dd\mu_k+r_k,
  \qquad
  |r_k|\le\frac{\|D^3\varphi\|_\infty}6\,\mathbb E|\xi_k|^3,
  \label{mart:c:eq:taylor}
\end{equation}
and $\sum_k\mathbb E|\xi_k|^3\le N^{1/4}(\tau\fA)^{3/4}=\tau^{1/2}\fA^{3/4}$.
Multiplying \eqref{mart:c:eq:taylor} by the cell averages of $\zeta\in C^1([0,1])$,
summing by parts and passing to the limit gives \eqref{mart:g:eq:weak-fp} for
$\zeta(t)\varphi(x)$. Approximation in the $C^{1,2}$ norm and finite sums of
such tensor products give the identity for all test functions in
$C_c^{1,2}([0,1]\times\R^d)$. The curve $(\mu_t)$ is narrowly continuous
because $X$ has continuous paths, and $\mu_t\in\Ptwo$ by the bound on
second moments.

\end{proof}

\subsection{Proof of \texorpdfstring{\cref{mart:c:prop:upper-bound}}{the martingale result}}
\begin{proof}

Let $X$ be as in \eqref{mart:c:eq:rate} and $Y_k:=X_{t_k}$. \Cref{mart:c:lem:continuous-bdg}
and Cauchy--Schwarz give
\begin{equation}
  \mathbb E|Y_{k+1}-Y_k|^4
  \le C_4\,\mathbb E\Bigl(\int_{t_k}^{t_{k+1}}R_t\dd t\Bigr)^2
  \le C_4\,\tau\,\mathbb E\int_{t_k}^{t_{k+1}}R_t^2\dd t .
\end{equation}
Summing over $k$ gives $\fA_N(Y)\le C_4\Am(\mu,A)$, and
\eqref{mart:c:eq:discrete-bb-martingale} concludes.

\end{proof}

\subsection{Proof of \texorpdfstring{\cref{mart:c:thm:bb-compactness}}{the martingale result}}
\begin{proof}

By \cref{mart:c:prop:upper-bound}, $\fA_N(X^N)\le C_4\MB^2_1(\mu,\sigma)+1$, and
$\mathbb E|X^N_0|^2$ does not depend on $N$, so \cref{mart:c:thm:compactness} applies.
Since $\Law(\widetilde X^N_0)=\mu$ and $\Law(\widetilde X^N_1)=\sigma$ for all $N$,
every limit curve joins $\mu$ to $\sigma$, and \eqref{mart:c:eq:action-lsc} gives
\eqref{mart:c:eq:lower-bound}. For the last claim, let $(N_i)$ realize the
$\liminf$. If it is finite, the same argument along $(N_i)$ produces an
admissible curve from $\mu$ to $\sigma$ with action at most
$\liminf_N3\fB_N(\mu,\sigma)$.

\end{proof}

\subsection{Proof of \texorpdfstring{\cref{mart:c:thm:sharp-upper-bound}}{the martingale result}}
\begin{proof}
Assume that the reference action is finite, and let $X$ be the continuous
martingale in \eqref{mart:c:eq:rate}. Write
\[
 R_t=\tr A_t(X_t),\qquad
 \Lambda=\mathbb E\int_0^1R_t^2\dd t<\infty.
\]
We sample $X$ at exits from small balls. This gives the sharp cost on
most steps; a small reserve of steps guarantees the exact endpoint $X_1$.

\emph{A regular approximation of the rate.}
Fix $\varepsilon,\delta\in(0,1)$ and set $\kappa=\delta/4$.
Extend $R$ by zero to negative times. By truncation and backward
averaging, choose $K,\gamma>0$ so that
\[
 u_t=\frac1\gamma\int_{t-\gamma}^t\min\{R_s,K\}\dd s,
 \qquad \|u-R\|_{L^2(\mathbb P\otimes\dd t)}\le\kappa\delta.
\]
Then $\widehat R_t:=\kappa+u_t$ is adapted, takes values in
$[\kappa,\kappa+K]$, and has a deterministic Lipschitz constant $L\le2K/\gamma$. Moreover,
\begin{equation}\label{mart:c:eq:rate-estimate}
 \|\widehat R-R\|_{L^2(\mathbb P\otimes\dd t)}\le\delta,
 \qquad
 \frac1\kappa\mathbb E\int_0^1(R_t-\widehat R_t)_+\dd t\le\delta.
\end{equation}
All these parameters remain fixed while $N\to\infty$.

\emph{Sampling and control of the random mesh.}
Set $\tau=1/N$, $\theta=\varepsilon/4$, $n=\lfloor(1-\theta)N\rfloor$
and $r=N-n\ge\theta N$. Starting from $S_0=0$, define
\begin{equation}\label{mart:c:eq:exit-times}
 h_k^2=(1+\varepsilon)\tau\widehat R_{S_k},\qquad
 S_{k+1}=\inf\{t\ge S_k:|X_t-X_{S_k}|\ge h_k\}\wedge1,
 \quad k<n.
\end{equation}
Use the remaining $r$ steps to reach time $1$:
\[
 S_{n+j}=S_n+\frac j r(1-S_n),\qquad 0\le j\le r.
\]
These are bounded stopping times, so $Y_k=X_{S_k}$ is a discrete
martingale with the prescribed endpoint marginals. Its discrete labels
are $k\tau$; the physical sampling times $S_k$ may be random.
Write
\[
 \xi_k=Y_{k+1}-Y_k,\qquad
 V_k=\int_{S_k}^{S_{k+1}}R_s\dd s.
\]
For $k<n$, optional stopping gives
$\mathbb E(|\xi_k|^2\mid\cG_{S_k})
=\mathbb E(V_k\mid\cG_{S_k})$.
The continuous martingale $t\mapsto X_{(S_k+t)\wedge S_{k+1}}-X_{S_k}$,
for the filtration $(\cG_{S_k+t})_{t\ge0}$, stays inside the ball of
radius $h_k$, and the trace of its quadratic variation at $t=\infty$ is
$V_k$. Hence the reverse Burkholder--Davis--Gundy inequality
\cite[Theorem~1.1]{MR16}, with exponent $4$, gives, with constants
independent of $N$,
\begin{equation}\label{mart:c:eq:step-tails}
 \mathbb E V_k^2
 \le C_{\mathrm{BDG}}\mathbb E h_k^4\le C\tau^2,
 \qquad
 \mathbb E\max_{k<n}V_k
 \le\Bigl(\sum_{k<n}\mathbb E V_k^2\Bigr)^{1/2}
 \le C\sqrt\tau.
\end{equation}
Let $\underline s=S_k$ for $s\in[S_k,S_{k+1})$, $k<n$, and
$\underline s=s$ for $s\ge S_n$. The Lipschitz bound on $\widehat R$
now controls freezing it at sampling times directly:
\begin{equation}\label{mart:c:eq:pi-error}
 F_N:=\int_0^1|\widehat R_{\underline s}-\widehat R_s|R_s\dd s
 \le L\sum_{k<n}(S_{k+1}-S_k)V_k
 \le L\max_{k<n}V_k,
 \qquad \mathbb E F_N\longrightarrow0.
\end{equation}

\emph{The reserve is rarely needed.}
Let $\Omega_N=\{S_n<1\}\in\cG_{S_n}$. Define the centered sum
\[
 Z_N=\sum_{k<n}\frac{|\xi_k|^2-V_k}{\widehat R_{S_k}}.
\]
Its summands are martingale differences. Since $\widehat R\ge\kappa$,
$|\xi_k|\le h_k$ and \eqref{mart:c:eq:step-tails} imply
$\mathbb E Z_N^2\le C\tau\to0$.
On $\Omega_N$, every exit is complete, so
\[
 1+\frac\varepsilon2-2\tau
 \le(1+\varepsilon)n\tau
 =Z_N+\int_0^{S_n}\frac{R_s}{\widehat R_{\underline s}}\dd s
 \le Z_N+1+\frac1\kappa\int_0^1(R_s-\widehat R_s)_+\dd s
          +\frac{F_N}{\kappa^2}.
\]
Taking expectations on this event and using
\eqref{mart:c:eq:rate-estimate} and \eqref{mart:c:eq:pi-error} yields
\begin{equation}\label{mart:c:eq:failure}
 \limsup_{N\to\infty}\mathbb P(\Omega_N)\le\frac{2\delta}{\varepsilon}.
\end{equation}

\emph{Action bound.}
For the exit steps, $|\xi_k|^4\le h_k^2|\xi_k|^2$ and optional
stopping give
\begin{align}
 \frac1\tau\sum_{k<n}\mathbb E|\xi_k|^4
 &\le(1+\varepsilon)\mathbb E\int_0^{S_n}\widehat R_{\underline s}R_s\dd s\notag\\
 &\le(1+\varepsilon)\bigl(\Lambda+\delta\sqrt\Lambda+\mathbb E F_N\bigr).
 \label{mart:c:eq:exit-cost}
\end{align}
The last inequality uses Cauchy--Schwarz and
$\|\widehat R-R\|_{L^2}\le\delta$.
The reserve increments vanish outside $\Omega_N$, and their physical
intervals have length at most $1/r$. Thus
\cref{mart:c:lem:continuous-bdg} and Cauchy--Schwarz give
\begin{equation}\label{mart:c:eq:reserve-cost}
 \frac1\tau\sum_{k\ge n}\mathbb E|\xi_k|^4
 \le\frac{C_4}{\tau r}\,
      \mathbb E\Bigl(\one_{\Omega_N}\int_{S_n}^1R_s^2\dd s\Bigr)
 \le\frac{C_4}{\theta}\,
      \mathbb E\Bigl(\one_{\Omega_N}\int_0^1R_s^2\dd s\Bigr).
\end{equation}
Since $U:=\int_0^1R_s^2\dd s$ is integrable,
$\omega(q):=\sup_{\mathbb P(B)\le q}\mathbb E(\one_B U)$ tends to zero
as $q\downarrow0$. Combining the two action bounds with
\eqref{mart:c:eq:discrete-bb-martingale}, for every
$q>2\delta/\varepsilon$ we obtain
\[
 \limsup_{N\to\infty}3\fB_N(\mu,\sigma)
 \le(1+\varepsilon)(\Lambda+\delta\sqrt\Lambda)
       +\frac{C_4}{\theta}\omega(q).
\]
First let $\delta\downarrow0$ and $q\downarrow0$ with
$q>2\delta/\varepsilon$, then let $\varepsilon\downarrow0$.
The resulting bound is $\limsup_N3\fB_N(\mu,\sigma)\le\Lambda$.
\end{proof}

\subsection{Proof of \texorpdfstring{\cref{mart:c:cor:bb-convergence}}{the martingale result}}
\begin{proof}
We now combine the lower bound for arbitrary discrete competitors with
the upper bound constructed from a continuous competitor. Denote the
two optimal values by
\[
 V_N:=3\fB_N(\mu,\sigma),\qquad V:=\MB_1^2(\mu,\sigma).
\]
The endpoint measures $\mu$ and $\sigma$ remain fixed throughout.
By \cref{mart:c:thm:bb-compactness}, $V\le\liminf_N V_N$.
Conversely, \cref{mart:c:thm:sharp-upper-bound} gives
$\limsup_N V_N\le\Am((\mu_t),A)$ for every continuous admissible
pair joining these endpoints. Taking the infimum over such pairs
gives $\limsup_N V_N\le V$. Thus $V_N\to V$, including the case
$V=+\infty$, for which the lower bound alone suffices.

If $V<\infty$, consider a limit pair $((\mu_t),A)$ obtained from
asymptotically minimizing chains. Its admissibility and
\eqref{mart:c:eq:lower-bound} give
\[
 V\le\Am((\mu_t),A)\le\liminf_i V_{N_i}=V.
\]
Hence the limiting pair attains the continuous minimum, proving the
second assertion.
\end{proof}

\subsection{The common one-dimensional limit}
\begin{proof}[Proof of \cref{mart:cor:common-limit}]
There are two sequences of discrete optimal values to compare:
$\mathsf G_N^+(\mu,\sigma)$ for the midpoint scheme with nonnegative
momentum, and $\fB_N(\mu,\sigma)$ for the fourth-moment coupling
scheme. Their limits have already been computed separately in
\cref{mart:cor:positive-midpoint,mart:c:cor:bb-convergence}.
The continuous identity $\MB_1^2=8\mathsf G^+_{[0,1]}$ from
\cref{mart:prop:positive-eulerian} therefore gives both equalities
in \eqref{mart:eq:common-discrete-limit}.

For attainment, suppose this common value is finite and take
asymptotically minimizing $M_4$ chains. Compactness gives a limiting
marginal curve $(\mu_t)$ and a nonnegative scalar covariance field
$a_t(x)$. The lower bound and convergence of the optimal values imply
that this pair minimizes $\MB_1^2$. To express the same pair in the
Eulerian variables, define the momentum measure at time $t$ by
\[
 m_t(\dd x):=\frac12a_t(x)\mu_t(\dd x).
\]
Then $\partial_t\mu=\partial_{xx}m$, and
\cref{mart:prop:positive-eulerian} identifies its action with one
eighth of the martingale action. This proves attainment of the
positive Eulerian minimum and the claim about limiting pairs.
The same proposition supplies a continuous martingale with these
marginals through superposition.
\end{proof}

\subsection{Proof of \texorpdfstring{\cref{mart:prop:tensor}}{the martingale result}}
\begin{proof}
This is a pointwise minimization over positive-semidefinite matrices.
The symmetric matrix $H$ is fixed, $A\succeq0$ is the variable, and
$H:A=\tr(HA)$ is their Frobenius product. We first fix the total
covariance $s:=\tr A\ge0$, optimize its directions, and then optimize
$s$.

Write $\lambda=\lambda_{\min}(H)$. Since $H-\lambda I\succeq0$,
\[
 H:A\ge\lambda\tr A=\lambda s.
\]
Equality holds precisely when
$\operatorname{Ran}(A)\subseteq E_{\min}(H)$, the eigenspace of the
smallest eigenvalue. Every $s\ge0$ admits such a matrix, for example
$A=s\,e\otimes e$ with a unit vector $e\in E_{\min}(H)$.
Consequently the original matrix problem reduces exactly to
\[
 \min_{s\ge0}\left\{\frac12s^2+\frac12\lambda s\right\}.
\]
Its minimizer is $s=(-\lambda)_+/2$. If $\lambda\ge0$, then
$s=0$ forces $A=0$; if $\lambda<0$, equality in the spectral bound
requires the stated range condition. This proves
\eqref{mart:g:eq:tensor-min-characterization}.
\end{proof}

\subsection{Proof of \texorpdfstring{\cref{mart:c:cor:gradient-flow}}{the martingale result}}
\begin{proof}
Our goal here is the action bound needed for the compactness theorem.
For each $N$, let $\tau=1/N$, let $\mu_k^N$ be the successive
minimizers, and let $X^N=(X_k^N)_{k=0}^N$ be their Markov gluing:
the transition from $\mu_k^N$ to $\mu_{k+1}^N$ is an optimal $M_4$
martingale coupling. Thus its fourth-moment action is exactly the
sum of the optimized coupling costs in
\eqref{mart:c:eq:action-cost}.

At step $k$, the unchanged marginal $\mu=\mu_k^N$ is admissible and
has zero transport cost. Minimality in \eqref{mart:g:eq:m4-jko} gives
\[
 \frac1{2\tau}M_4(\mu_k^N,\mu_{k+1}^N)
 \le\Em(\mu_k^N)-\Em(\mu_{k+1}^N).
\]
Summing makes the energy differences telescope, so
\begin{equation}
 \fA_N(X^N)=\frac3\tau\sum_{k=0}^{N-1}M_4(\mu_k^N,\mu_{k+1}^N)
 \le6\bigl(\Em(\mu_{\mathrm{in}})-\Em(\mu_N^N)\bigr)\le6\Phi.
\end{equation}
Here $\Phi$ is the uniform bound on energy decrease in
\eqref{mart:c:eq:bounded-below}. The initial distribution is always
$\mu_{\mathrm{in}}$, so its second moment is also uniformly bounded.
These are exactly the two hypotheses of
\cref{mart:c:thm:compactness}. It supplies continuous martingale
limits, their marginal curves and the lower bound
$\Am(\mu,A)\le6\Phi$ asserted here. Identification of a particular
limiting evolution is treated in the potential-energy proofs below.
\end{proof}

\subsection{Proof of \texorpdfstring{\cref{mart:prop:pointwise}}{the martingale result}}
\begin{proof}
The proof first separates the optimization at each starting position,
then explains why at most two increments are needed for a pointwise
minimizer. The distinction between positions and increments is useful
here: after disintegrating a coupling as
$\pi(\dd x,\dd y)=\mu_k(\dd x)\pi_x(\dd y)$, let $K_x$ be the
conditional distribution of $u=y-x$. Thus $K_x$ is centered, whereas
$\pi_x$ has mean $x$; the endpoint kernel $\pi_x$ is obtained by
translating $K_x$ by $x$.

\emph{Separation at the starting position.}
The terminal marginal is free, and the martingale condition is exactly
$\int uK_x(\dd u)=0$ for $\mu_k$-almost every $x$. Consequently the
objective is $\int\mathcal J_{\tau,x}(K_x)\mu_k(\dd x)$, and each
conditional objective is bounded below by its own infimum. This proves
the lower bound in \eqref{mart:g:eq:pointwise-separation}. There is no
constraint coupling the choices of $K_x$ at different starting positions:
the next marginal is the distribution of $x+u$ resulting from these
choices.

\emph{Existence and the two-point description.}
For a fixed $x$, write
$f_{\tau,x}(u)=\eta|u|^4/\tau+\Vm(x+u)$. The quadratic lower bound
on $\Vm$ gives
$f_{\tau,x}(u)\ge c_{\tau}|u|^4-C_{\tau,x}$, with $c_\tau>0$.
A minimizing sequence is therefore tight with bounded fourth moments
and uniformly integrable first moments. The mean-zero condition passes
to its weak limit, which has finite fourth moment by lower
semicontinuity. Lower semicontinuity of the objective then gives a
minimizer.

Fix $(\tau,x)$ and abbreviate $f=f_{\tau,x}$. Its convex envelope
$f^{**}$, the largest lower semicontinuous convex function below $f$, satisfies
\[
 f^{**}(0)
 =\inf_{\substack{K\in\mathcal P(\mathbb R)\\
                  \int u\,K(\dd u)=0}}
   \int f(u)\,K(\dd u).
\]
To prove this identity, let $g(s)$ be the infimum of $\int f\,\dd K$
over probability measures $K$ with $\int u\,K(\dd u)=s$. Mixing two
such measures shows that $g$ is convex. The bound $f\ge-C_{\tau,x}$,
Jensen's inequality for $f^{**}$ and the choice $K=\delta_s$ give
$-C_{\tau,x}\le f^{**}\le g\le f<\infty$. A finite convex function on
$\R$ is continuous, so $g$ is a lower semicontinuous convex function
below $f$. Hence $g\le f^{**}$, and therefore $g=f^{**}$.
Since $f(u)\ge c_\tau|u|^4-C_{\tau,x}$, every $K$ with infinite fourth
moment gives $\int f\,\dd K=+\infty$. The infimum over
$\mathcal P(\mathbb R)$ therefore equals the infimum over $\Pcal_4(\R)$
in \eqref{mart:g:eq:pointwise-separation}.
Thus $f^{**}(0)$ is the smallest average of values of $f$ among
probability distributions with mean zero. Since $f^{**}$ is convex and
finite on $\R$, it has a supporting affine function at zero:
\[
 \ell(u)=f^{**}(0)+pu,
 \qquad \ell(u)\le f^{**}(u)\le f(u).
\]
Its contact set with $f$ is
\[
 C=\{u\in\mathbb R:f(u)=\ell(u)\}.
\]
Since $f-\ell$ is lower semicontinuous and tends to $+\infty$ as
$|u|\to\infty$, this set is compact. For any minimizing centered
distribution $K_*$,
\[
 \int\bigl(f(u)-\ell(u)\bigr)\,K_*(\dd u)
 =f^{**}(0)-f^{**}(0)=0.
\]
The integrand is nonnegative, so $K_*(C)=1$.

If $0\in C$, then $f(0)=f^{**}(0)$, and $\delta_0$ is a minimizer.
Otherwise, because $K_*$ is centered and supported on $C$, there exist
$a,b\in C$ with $a<0<b$. Set
\[
 K_{a,b}=\frac{b}{b-a}\delta_a+\frac{-a}{b-a}\delta_b.
\]
These weights are positive, sum to one, and give mean zero:
\[
 \frac{b}{b-a}a+\frac{-a}{b-a}b=0.
\]
Moreover, $f=\ell$ at both atoms, and therefore
\[
 \int f\,\dd K_{a,b}
 =\int\ell\,\dd K_{a,b}
 =\ell(0)=f^{**}(0).
\]
Hence $K_{a,b}$ is a minimizer with two atoms.

Finally, a probability measure supported on $\{a,b\}$ is centered only
if its weights solve
\[
 \alpha+\beta=1,\qquad \alpha a+\beta b=0,
\]
whose unique solution is the one above. Any convex decomposition of
$K_{a,b}$ into centered probability measures must use measures supported
on $\{a,b\}$, so both must equal $K_{a,b}$. Thus $K_{a,b}$ is an
extreme point of the set of centered probability measures, in agreement
with the one-moment case of \cite{Pin12}.

Finally, the map $(x,K)\mapsto\mathcal J_{\tau,x}(K)$ is lower
semicontinuous, and the quadratic lower bound on $\Vm$ allows the
choice $C_{\tau,x}\le C'_\tau(1+x^2)$ for a constant $C'_\tau$. The
tightness argument above, applied for $x$ in a compact set, therefore
shows that $x\mapsto\inf\mathcal J_{\tau,x}$ is lower semicontinuous.
This function lies between $-C'_\tau(1+x^2)$ and $\Vm(x)$, so its
integral against $\mu_k\in\Ptwo(\R)$ is defined in
$(-\infty,+\infty]$. Measurable selection gives a family
$x\mapsto K_x$ of minimizers, measurable for the completion of
$\mu_k$. Translating these increment distributions gives a martingale
coupling with first marginal $\mu_k$. Integrating their pointwise
minimum values proves the reverse inequality in
\eqref{mart:g:eq:pointwise-separation}. This coupling therefore
attains the infimum on the left, which completes the proof.
\end{proof}

\subsection{Proof of \texorpdfstring{\cref{mart:thm:oned}}{the martingale result}}
\begin{proof}
We first bound the second and fourth moments of an optimal increment,
then determine its variance to leading order. These estimates identify
the generator of every possible continuous-path limit; uniqueness of
the resulting diffusion then gives convergence of the whole family.
Throughout the first two steps, $K=K_{\tau,x}$ is a minimizing
increment distribution. Write $v=v_\tau(x)$ and $w=w_\tau(x)$.
Since $K$ is centered, $v$ is its variance, and $w$ is its fourth moment.

\emph{Moment bounds.}
Uniform concavity and centeredness give
\begin{equation}
  -\frac L2v
  \le\int\bigl(\Vm(x+u)-\Vm(x)-\Vm'(x)u\bigr)K(\dd u)
  \le-\frac{L_0}2v.
  \label{mart:g:eq:concavity-bounds}
\end{equation}
The integral of the linear term is zero. Optimality against
$K=\delta_0$ gives $\mathcal J_{\tau,x}(K)\le\Vm(x)$.
Combining this with the lower bound above and Jensen's inequality
$w\ge v^2$ yields
\begin{equation}
  \frac\eta\tau w\le\frac L2v,
  \qquad
  v\le\frac{L}{2\eta}\tau,
  \qquad
  w\le\frac{L^2}{4\eta^2}\tau^2.
  \label{mart:g:eq:cell-moment-bound}
\end{equation}
Thus optimal increments have variance of order $\tau$, while their
fourth moments are of order $\tau^2$.

\emph{Identification of the variance.}
Taylor's formula gives, uniformly in $x$,
\begin{equation}
  \mathcal J_{\tau,x}(K)-\Vm(x)
  =\frac\eta\tau w+\frac12\Vm''(x)v+O(w^{3/4}).
  \label{mart:g:eq:cell-taylor}
\end{equation}
Indeed, the remainder is bounded by a constant times
$\int|u|^3K(\dd u)\le w^{3/4}$. By
\eqref{mart:g:eq:cell-moment-bound}, it is $O(\tau^{3/2})$.
Put $c(x)=-\Vm''(x)>0$. Using $w\ge v^2$ and completing the square
in the variance gives
\begin{equation}
  \mathcal J_{\tau,x}(K)-\Vm(x)
  \ge-\frac{c(x)^2}{16\eta}\tau
  +\eta\tau\left(\frac v\tau-\frac{c(x)}{4\eta}\right)^2
  -C\tau^{3/2}.
  \label{mart:g:eq:cell-coercive-comparison}
\end{equation}
The square measures the difference between the optimal variance rate
$v/\tau$ and its proposed limit $c(x)/(4\eta)$.
To show that this difference tends to zero, compare the minimizer with
a centered two-point increment distribution. For every $r\ge0$, equal
weights on $\{\pm\sqrt{r\tau}\}$ give $v=r\tau$ and
$w=r^2\tau^2$, hence equality in $w\ge v^2$.
Taking $r=c(x)/(4\eta)$ gives
\begin{equation}
  \mathcal J_{\tau,x}(K_{\tau,x})-\Vm(x)
  \le-\frac{c(x)^2}{16\eta}\tau+C\tau^{3/2}.
\end{equation}
The lower and upper bounds force
\begin{equation}
  \sup_x\left|
  \frac{v_\tau(x)}\tau+\frac{\Vm''(x)}{4\eta}
  \right|\le C\tau^{1/4}.
\end{equation}
This proves the variance assertion in
\eqref{mart:g:eq:variance-asymptotic}; its fourth-moment assertion
already follows from \eqref{mart:g:eq:cell-moment-bound}. The argument
applies to every minimizer, without requiring that the chosen
minimizer itself have only two atoms. The constants depend only on
$\eta$, $L$ and $\|\Vm'''\|_\infty$, so these bounds hold uniformly in
$x$ and in the choice of the minimizers, also when the family of
minimizers changes with the step.

\emph{From increment moments to the limiting diffusion.}
Given $X_k=x$, the next position is $x+u$ with
$u\sim K_{\tau,x}$. Its discrete generator is the expected change
of a test function during one step, divided by $\tau$.
For $\varphi\in C_c^3(\R)$, centeredness and the same Taylor estimate
give
\begin{equation}
  \frac1\tau\int\bigl(\varphi(x+u)-\varphi(x)\bigr)K_{\tau,x}(\dd u)
  =\frac{v_\tau(x)}{2\tau}\varphi''(x)+O(\tau^{1/2}).
  \label{mart:g:eq:generator-convergence}
\end{equation}
The factor $1/2$ comes from the second-order Taylor term. The discrete
generators therefore converge locally uniformly to
$\mathcal L\varphi=-\Vm''\varphi''/(8\eta)$.
For tightness, let $N=\lceil\tfin/\tau\rceil$ and $\tau\le1$, so that
$\tfin\le N\tau\le\tfin+1$. The time change $t\mapsto t/(N\tau)$ maps
the grid $(k\tau)_{0\le k\le N}$ onto $(k/N)_{0\le k\le N}$. By
\eqref{mart:g:eq:cell-moment-bound}, the action
\eqref{mart:c:eq:discrete-action} of the time-changed chain is
$N\sum_{k<N}\mathbb E|X_{k+1}-X_k|^4\le N^2L^2\tau^2/(4\eta^2)
\le(\tfin+1)^2L^2/(4\eta^2)$. Thus \cref{mart:c:thm:compactness} gives
tightness of the time-changed chains in $C([0,1];\R)$. Undoing the
time change and restricting to $[0,\tfin]$ gives tightness of the
linearly interpolated chains in $C([0,\tfin];\R)$.

For completeness, generator convergence identifies the process as
follows. At grid times, subtracting from $\varphi(X_k)$ the accumulated
conditional expected changes gives a martingale. The test function
and its discrete generator are uniformly bounded. Passing this
identity to a subsequential path limit shows that
$\varphi(X_t)-\varphi(X_0)-\int_0^t\mathcal L\varphi(X_s)\dd s$
is a martingale. This is the martingale problem for the operator
$\mathcal L$, equivalently the diffusion
\eqref{mart:g:eq:limit-sde-beta}. Its variance rate is
$-\Vm''/(4\eta)$, and the Fokker--Planck equation has half that
coefficient, explaining $8\eta$ in
\eqref{mart:g:eq:limit-pde-beta}.

Finally, \eqref{mart:g:eq:F-assumptions} makes the diffusion
coefficient bounded, uniformly positive, and Lipschitz. The limiting
martingale problem is therefore unique. Every subsequential path
limit has this same distribution and initial marginal $\mu_0$;
consequently the entire family converges, with the stated limiting
marginals. The superposition principle \cite[Theorem~2.5]{Tre16}
applies because the diffusion coefficient is bounded: every narrowly
continuous probability solution of \eqref{mart:g:eq:limit-pde-beta} is
the marginal curve of a solution of this martingale problem with
initial distribution $\mu_0$. Uniqueness of the martingale problem
therefore gives uniqueness of this solution.
\end{proof}

\subsection{The modified JKO scheme with nonnegative momentum}
\label{mart:app:positive-potential}
\begin{proof}[Proof of \cref{mart:thm:positive-potential}]
The proof separates the estimates for a single minimizing step from
passage to a continuous evolution. We first prove that every step
exists and control the increase of its second moment. The main
identification then compares a minimizing step with an implicit step
of the candidate diffusion equation. This comparison shows that the
discrete momentum per unit time approaches $D\mu$. Compactness,
uniqueness of that diffusion, and the energy identity complete the proof.

For the single-step arguments, $\mu$ is the prescribed old marginal,
$\nu$ is the new marginal being optimized, and $m$ is the finite
nonnegative momentum measure integrated over a time interval of
length $\tau$. Thus $m''=\nu-\mu$. Write
$\Em(\nu)=\int\Vm\dd\nu$, $d=\|D\|_\infty$, and
$S=(\Vm'')_+/\kappa$. Here $D=(-\Vm'')_+/\kappa$ is the candidate
momentum rate determined by the negative part of the second derivative,
whereas $S$ records its positive part. Both $D$ and $S$ are bounded and
Lipschitz, $DS=0$, and $\Vm''=\kappa(S-D)$. The potential has at most
quadratic growth, so $\Em$ is finite and continuous on $\Ptwo(\R)$.

\emph{Integration by parts and existence of each step.}
For a feasible pair $(\nu,m)$ with old marginal $\mu\in\Ptwo(\R)$,
testing the conservation equation with $x^2$ measures the increase
of the second moment, while testing it with $\Vm$ measures the
change in potential energy. Cutoff approximations of these tests give
\begin{equation}\label{mart:eq:positive-potential-identities}
 M_2(\nu)-M_2(\mu)=2m(\R),\qquad
 \Em(\nu)-\Em(\mu)=\int\Vm''\dd m,
 \qquad M_2(\lambda):=\int x^2\lambda(\dd x).
\end{equation}
The cutoff errors vanish because the marginals have finite second moments,
$m$ is finite and $\Vm''$ is bounded. Set
$K=2\Ac(\nu,m)=\int|\dd m/\dd\nu|^2\dd\nu$ for a finite-action
pair. This nonnegative scalar $K$ measures the size of the step momentum.
Since $\nu$ has mass one, Cauchy--Schwarz gives
$m(\R)\le\sqrt K$ and
\[
 \frac\kappa\tau\Ac(\nu,m)+\Em(\nu)
 \ge\Em(\mu)+\frac\kappa{2\tau}K-\kappa d\sqrt K.
\]
Consequently every minimizing sequence has bounded $K$, hence bounded
second moments by \eqref{mart:eq:positive-potential-identities}.
For a spatial Borel set $E$,
\begin{equation}\label{mart:eq:positive-potential-tail-cs}
 m(E)\le\sqrt{K\nu(E)}.
\end{equation}
This also gives tightness of its momentum measures.

A bound on second moments alone does not give compactness in $W_2$;
we also need a uniform bound on their tails. For a spatial cutoff radius $R>0$, use
$\psi_R(x)=(|x|-R)_+^2$, which vanishes on $[-R,R]$ and tests the
quadratic tails outside that interval. A finite measure $m$ with $m''=\nu-\mu$
has a locally continuous density, by integrating this distributional
identity twice. Smooth approximation and cutoff therefore give
\begin{equation}\label{mart:eq:positive-potential-tail-identity}
 \int\psi_R\dd\nu=\int\psi_R\dd\mu+2m(\{|x|>R\}).
\end{equation}
If $M_2(\nu)\le M$ and $K\le K_*$, then
\[
 \int_{|x|>2R}x^2\dd\nu
 \le4\int\psi_R\dd\mu+\frac{8\sqrt{K_*M}}R
 \longrightarrow0
\]
uniformly along the sequence. Thus it is precompact in $W_2$, and
\eqref{mart:eq:positive-potential-tail-cs} gives weak compactness of
the momenta with their total masses preserved. The constraint and
nonnegativity are closed, the perspective action is lower
semicontinuous, and $\Em$ is $W_2$-continuous. The direct method proves
attainment. The zero-momentum competitor makes every subsequent step
finite, so the recursion is defined for all $k$.

\emph{Bounds inherited from the Eulerian conservation law.}
For a minimizing pair write $m=\tau b\nu$, that is,
$b=\tau^{-1}\dd m/\dd\nu\ge0$. The function $b$ is the discrete
momentum rate relative to the new marginal. In particular
$K=\tau^2\int b^2\dd\nu$. Comparison with $(\mu,0)$
and \eqref{mart:eq:positive-potential-identities} imply
\[
 \frac\kappa2\int b^2\dd\nu
 \le-\int\Vm''b\dd\nu\le\kappa d\int b\dd\nu.
\]
Hence
\begin{equation}\label{mart:eq:positive-potential-step-bounds}
 \|b\|_{L^2(\nu)}\le2d,\qquad m(\R)\le2d\tau,
 \qquad M_2(\nu)-M_2(\mu)\le4d\tau.
\end{equation}
By \cref{mart:prop:coupling-momentum}, $\mu\le_c\nu$, so these
marginals admit a martingale coupling. For any such coupling, its
expected squared increment is $M_2(\nu)-M_2(\mu)$. This supplies a
competitor for the Wasserstein cost,
so
\begin{equation}\label{mart:eq:positive-potential-w2-step}
 W_2^2(\mu,\nu)\le M_2(\nu)-M_2(\mu)\le4d\tau.
\end{equation}
If $d=0$, these estimates force every step to be stationary and all
assertions follow. Assume below that $d>0$.

\emph{The resolvent comparison.}
We now construct a feasible comparison pair whose momentum rate is
exactly the candidate $D$. This is an auxiliary competitor; the
minimizing step itself has not yet been identified. Since $D$ is
bounded, nonnegative and Lipschitz,
\[
 |\sqrt{2D(x)}-\sqrt{2D(y)}|^2
 \le2\operatorname{Lip}(D)|x-y|.
\]
Weak existence for a bounded continuous diffusion coefficient and
the one-dimensional Yamada--Watanabe uniqueness criterion
\cite[Theorem~1 and Proposition~1]{YW71} give a unique diffusion distribution
for $\dd X=\sqrt{2D(X)}\dd B$. It is a square-integrable martingale.
Its covariance rate is $2D$, so its marginal equation is
$\partial_t\lambda=\partial_{xx}(D\lambda)$. Starting it with
distribution $\mu$, let $P_t\mu=\Law(X_t)$ denote the marginal at
time $t$. Thus $(P_t)_{t\ge0}$ is the evolution operator of this
diffusion. Define
\[
 \nu_*:=\frac1\tau\int_0^\infty e^{-t/\tau}P_t\mu\dd t,
 \qquad m_*:=\tau D\nu_*.
\]
The measure $\nu_*$ averages this evolution over an independent
exponential time of mean $\tau$, and $m_*$ is its integrated
momentum. Integration of the diffusion equation against the exponential
weight gives the exact resolvent identity
\begin{equation}\label{mart:eq:positive-potential-resolvent}
 \nu_*-\tau\partial_{xx}(D\nu_*)=\mu.
\end{equation}
This equation is the implicit Euler step, or resolvent, of the
candidate diffusion equation. In particular $(\nu_*,m_*)$ is feasible,
including when $D$ vanishes, and
$M_2(\nu_*)-M_2(\mu)\le2d\tau$.

For any finite-action feasible $m=\tau b\nu$, the energy identity
and $\Vm''=\kappa(S-D)$ express the step objective as a quadratic
function of $b$. Completion of squares gives the
exact identity
\begin{equation}\label{mart:eq:positive-potential-square}
 \frac\kappa\tau\Ac(\nu,m)+\Em(\nu)-\Em(\mu)
 =\frac{\kappa\tau}2\int
       \bigl[(b-D)^2-D^2+2Sb\bigr]\dd\nu.
\end{equation}
The term $(b-D)^2$ measures the discrepancy from the candidate
momentum rate. The additional term $Sb$ is nonnegative and penalizes
momentum where the potential is convex. The resolvent competitor has
$b_*=D$ and $Sb_*=0$. For a minimizing pair, comparison with this
competitor therefore gives
\begin{align}
 \int(b-D)^2\dd\nu+2\int Sb\dd\nu
 &\le\int D^2\dd(\nu-\nu_*)\notag\\
 &\le\operatorname{Lip}(D^2)
        \bigl(W_2(\nu,\mu)+W_2(\mu,\nu_*)\bigr)\notag\\
 &\le C\sqrt\tau,
 \qquad C=(2+\sqrt2)\sqrt d\operatorname{Lip}(D^2).
 \label{mart:eq:positive-potential-square-bound}
\end{align}
The second term on the left is nonnegative, so this bounds the
$L^2(\nu)$ discrepancy between $b$ and $D$ for every minimizing step.
Along the recursion, write $b_\tau(t)=b_{k+1}^\tau$ and
$\bar\mu_\tau(t)=\mu_{k+1}^\tau$ on
$(k\tau,(k+1)\tau]$. Cauchy--Schwarz on each time cell and integration
over its actual intersection with $(0,T)$ give
$\int_0^T\int|b_\tau-D|\dd\bar\mu_\tau\dd t
\le T\sqrt C\,\tau^{1/4}$, which is
\eqref{mart:eq:positive-potential-flux-error}.

\emph{Real-line compactness.}
We next combine the single-step bounds over a fixed time interval.
The tail estimate prevents second moments from escaping to infinity,
while convex order gives a time modulus. Together they yield
compactness of the marginal curves in $W_2$.
The step bounds give
$\sup_{k\tau\le T+1}M_2(\mu_k^\tau)\le
M_T:=M_2(\mu_0)+4d(T+1)$. Summing
\eqref{mart:eq:positive-potential-tail-identity}, and using
\eqref{mart:eq:positive-potential-tail-cs} with
$\sqrt K=\tau\|b\|_{L^2(\nu)}\le2d\tau$, yields
\[
 \int\psi_R\dd\mu_k^\tau
 \le\int\psi_R\dd\mu_0+
           \frac{4d(T+1)\sqrt{M_T}}R
 \qquad(k\tau\le T+1,\ \tau\le1).
\]
The second moments are therefore uniformly integrable. All time
marginals lie in a precompact subset of $\Ptwo(\R)$.
Convex order along the chain and the moment identity give
\[
 W_2^2(\mu_k^\tau,\mu_\ell^\tau)
 \le4d(\ell-k)\tau\qquad(k\le\ell).
\]
On $[k\tau,(k+1)\tau]$, the affine interpolation is
$(1-\theta)\mu_k^\tau+\theta\mu_{k+1}^\tau$, with
$\theta=(t-k\tau)/\tau$. It is also increasing in convex order, and
its second moment has slope at most $4d$. It consequently has a uniform
$W_2$ one-half H\"older modulus. Its distance to $\bar\mu_\tau$ is
at most $2\sqrt{d\tau}$. Arzel\`a--Ascoli gives a subsequence
converging uniformly in $W_2$ to a continuous curve $\mu_t$.

\emph{Identification, uniqueness and energy.}
It remains to identify the continuous equation for each subsequential
limit. On each time cell the derivative of the affine interpolation
is $(\mu_{k+1}^\tau-\mu_k^\tau)/\tau$. Thus the discrete constraint
is exactly the conservation equation with momentum rate
$b_\tau\bar\mu_\tau$. By
\eqref{mart:eq:positive-potential-flux-error}, this measure differs
from $D\bar\mu_\tau$ by $O_T(\tau^{1/4})$ in space--time total
variation. Uniform narrow convergence and bounded continuity of $D$
therefore pass the equation to
\eqref{mart:eq:positive-potential-pde}, including its initial trace.
The superposition principle \cite[Theorem~2.5]{Tre16} applies because
$2D$ is bounded. Every narrowly continuous probability solution is
represented by the unique diffusion distribution constructed above.
Thus the process associated with the limiting equation is a continuous
square-integrable martingale. This proves uniqueness
and convergence of the entire family in
\eqref{mart:eq:positive-potential-convergence}.

Finally, telescoping the changes in potential energy at each step
and passing to the limit gives the energy identity. Sum the second identity in
\eqref{mart:eq:positive-potential-identities} over steps. Boundedness
of $\Vm''$ and \eqref{mart:eq:positive-potential-flux-error} identify
its limit as $\int_0^t\int\Vm''D\dd\mu_s\dd s$.
Uniform $W_2$ convergence gives convergence of the endpoint energies;
rounding $t$ to a time node changes the discrete energy by at most
$2\|\Vm''\|_\infty d\tau$. Since $\Vm''D=-\kappa D^2$, this is
\eqref{mart:eq:positive-potential-energy}.
\end{proof}

\begin{proof}[Proof of \cref{mart:cor:positive-potential-normalization}]
The two schemes use different kinetic coefficients, so we compare
the coefficients in their limiting equations. Under
\eqref{mart:g:eq:F-assumptions}, $D=-\Vm''/\kappa$.
The equation \eqref{mart:eq:positive-potential-pde} agrees with
\eqref{mart:g:eq:limit-pde-beta} precisely when $\kappa=8\eta$.
The two convergence theorems and uniqueness identify the same marginal
curve. The stated constants follow from $\eta=1/2$ for $3M_4$ and
$\eta=1/6$ for $M_4$. This comparison concerns the limiting marginal curves; it does not
identify the two finite-step optimizers. For a convex potential,
$D=0$, and
\eqref{mart:eq:positive-potential-step-bounds} gives $m=0$ and
$\nu=\mu$ at every step.
For the quadratic example, $D=a/\kappa$ is constant, so no error
arises from varying $D$ between the optimizing and comparison
marginals. The square
comparison \eqref{mart:eq:positive-potential-square-bound} has zero
right-hand side and forces $b=D$. The constraint is the stated
resolvent equation, whose unit-mass fundamental solution is the
displayed exponential density.
\end{proof}

\subsection{Quadratic fourth-moment example and normalization}
\begin{proof}[Verification of \cref{mart:ex:quadratic}]
Here the Taylor expansion of the potential is exact, so the previous
comparison reduces to a scalar minimization in the variance.
As before, $K$ denotes the centered distribution of the increment
$u$, and $v=\int u^2\dd K$ and $w=\int u^4\dd K$ are its second
and fourth moments. For $\Vm(x)=-ax^2/2$, centeredness cancels the
linear term in $u$, and Jensen's inequality gives
\begin{equation}
  \mathcal J_{\tau,x}(K)-\Vm(x)
  =\frac\eta\tau w-\frac a2v
  \ge\frac\eta\tau v^2-\frac a2v.
\end{equation}
The scalar minimum occurs at $v=a\tau/(4\eta)>0$.
Equality in Jensen's inequality requires $u^2$ to be constant
$K$-almost surely. Thus the only possible increments are
$\pm\sqrt{a\tau/(4\eta)}$, and the mean-zero condition forces
equal weights. This proves the uniqueness claimed in
\cref{mart:ex:quadratic}. Since this increment distribution does not
depend on $x$, its Markov iteration is a symmetric random walk with
variance rate $a/(4\eta)$, whose diffusion limit is Brownian motion
with that variance rate.

The normalization concerns this variance rate and the corresponding
time scale. Substituting $\eta=1/6$ for $M_4$ and $\eta=1/2$ for
$3M_4$ into \eqref{mart:g:eq:limit-pde-beta} gives the two equations
in \cref{mart:cor:normalization}. The constrained momentum equation
has coefficient $1/\kappa$ in place of $1/(8\eta)$, so the two
potential flows agree at the same time when $\kappa=8\eta$;
otherwise they agree after a constant rescaling of time.
\end{proof}

\subsection{The corresponding higher-dimensional problem}
\begin{proof}[Verification of \eqref{mart:g:eq:high-d-cell-bound}]
For $r=\mathbb E|V|^2$, Jensen gives $\mathbb E|V|^4\ge r^2$ and the
spectral bound gives $\mathbb E[V^T\mathcal BV]\le\lambda_{\max}(\mathcal B)r$.
Equality in the first requires constant radius, and equality in the
second requires support in the top eigenspace. The remaining scalar
minimum is attained at $r=\tau\lambda_{\max}(\mathcal B)/(4\eta)$. The centered
two-point law in any top eigendirection realizes both equalities.
An eigenspace of dimension one forces rank-one covariance; if it has
higher dimension, centered laws of constant radius give multiple
minimizing covariance matrices, including the isotropic choice when the
whole space is the top eigenspace.
\end{proof}

\clearpage
\phantomsection

\end{document}